\documentclass[10pt,reqno]{amsart}

\usepackage{graphicx}
\usepackage{times}
\usepackage[colorlinks=true,linkcolor=blue,citecolor=blue]{hyperref}%

\usepackage{xcolor}
\usepackage{stmaryrd}
\usepackage{pdfrender,xcolor}
\usepackage{dsfont}

\newtheorem{theorem}{Theorem}[section]
\newtheorem{lemma}[theorem]{Lemma}

\newtheorem{corollary}[theorem]{Corollary}
\newtheorem{prop}[theorem]{Proposition}

\newtheorem{ass}[theorem]{Assumption}

\theoremstyle{definition}
\newtheorem{definition}[theorem]{Definition}

\theoremstyle{remark}

\newtheorem*{rem}{Remark}
\numberwithin{equation}{section}

\usepackage[top=1.25in, bottom=1.25in, left=1.25in, right=1.25in]{geometry}

\usepackage[scr]{rsfso}
\usepackage[english]{babel}
\usepackage{fancyhdr}
\usepackage{amsmath}
\usepackage{amsthm}
\usepackage{amssymb}
\usepackage{eucal}
\usepackage{comment}
\usepackage{enumitem}
\setlist{leftmargin=*}
\usepackage[integrals]{wasysym}

\makeatletter
\newsavebox{\@brx}
\newcommand{\llangle}[1][]{\savebox{\@brx}{\(\m@th{#1\langle}\)}%
  \mathopen{\copy\@brx\kern-0.5\wd\@brx\usebox{\@brx}}}
\newcommand{\rrangle}[1][]{\savebox{\@brx}{\(\m@th{#1\rangle}\)}%
  \mathclose{\copy\@brx\kern-0.5\wd\@brx\usebox{\@brx}}}
\makeatother

\newcommand\nc{\newcommand}
\nc{\on}{\operatorname}
\nc{\E}{\mathbb{E}}
\nc{\R}{\mathbb R}
\nc{\C}{\mathbb C}
\nc{\Q}{\mathbb Q}
\nc{\Z}{\mathbb Z}
\nc{\N}{\mathrm{N}}
\nc{\F}{\mathbb F}
\nc{\wt}{\widetilde}
\nc{\ol}{\overline}
\nc{\short}[3]{0 \longrightarrow #1 \longrightarrow #2 \longrightarrow #3 \longrightarrow 0}
\nc{\pd}[2]{\frac{\partial #1}{\partial #2}}
\nc{\rnc}{\renewcommand}
\nc{\e}{\varepsilon}
\nc{\DMO}{\DeclareMathOperator}
\nc{\grad}{\nabla}
\nc{\Exp}{\mathrm{Exp}}
\rnc{\t}{\mathrm{t}}
\nc{\s}{\mathrm{s}}
\nc{\x}{\mathrm{x}}
\nc{\y}{\mathrm{y}}
\nc{\w}{\mathrm{w}}
\nc{\z}{\mathrm{z}}
\rnc{\r}{\mathrm{r}}
\nc{\kpz}{$\mathrm{KPZ}(\alpha,\lambda)$ }
\nc{\fsp}{}

\rnc{\leq}{\leqslant}
\rnc{\geq}{\geqslant}
\rnc{\d}{\mathrm{d}}
\rnc{\O}{\mathrm{O}}
\rnc{\exp}{\mathrm{Exp}}
\newenvironment{nouppercase}{%
  \renewcommand{\uppercasenonmath}[1]{}}{}
\title{\fsp\Large Time-inhomogeneous KPZ equation from non-equilibrium Ginzburg-Landau SDEs \vspace{-0.1cm}}

\author{Kevin Yang}
\usepackage{setspace}
\begin{document}
\setstretch{1.0}
\fsp
\raggedbottom
\begin{nouppercase}
\maketitle
\end{nouppercase}
\begin{center}
\today
\end{center}

\begin{abstract}
{We introduce a framework, which is a mesoscopic-fluctuation-scale analog of Yau's method \cite{YauRE} for hydrodynamic limits, for deriving KPZ equations with time-dependent coefficients from time-inhomogeneous interacting particle systems. To our knowledge, this is the first derivation of a time-inhomogeneous KPZ equation whose solution theory has an additional nonlinearity that is absent in the time-homogeneous case. So, we also show global well-posedness for the SPDE. To be concrete, we restrict to time-inhomogeneous Ginzburg-Landau SDEs. The method for deriving KPZ is based on a Cole-Hopf transform, whose analysis is the bulk of this paper. The key ingredient for said analysis is a ``local" second-order Boltzmann-Gibbs principle, which builds on prior work \cite{Y22} of the author. This addresses a ``Big Picture Question" in \cite{KPZAIM} on deriving KPZ equations. It is also, to our knowledge, a first result on KPZ-type limits in a non-equilibrium like that in \cite{CYau}.} 
\end{abstract}

{\hypersetup{linkcolor=blue}
\setcounter{tocdepth}{1}
\tableofcontents}

\fsp

\section{Introduction}
Kardar, Parisi, and Zhang \cite{KPZ} introduced a model for non-equilibrium interface fluctuations that is now ubiquitously called the \emph{Kardar-Parisi-Zhang} (KPZ) equation. This model is the following SPDE (where $\t\geq0$ and $\x\in\mathbb{T}:=\R/\Z$):
\begin{align}
\partial_{\t}\mathfrak{h}(\t,\x) \ = \ \bar{\alpha}\partial_{\x}^{2}\mathfrak{h}(\t,\x)+\bar{\alpha}(\wedge)|\partial_{\x}\mathfrak{h}(\t,\x)|^{2}+\xi(\t,\x). \label{eq:kpzbasic}
\end{align}
Above, $\bar{\alpha}>0$ and $\bar{\alpha}(\wedge)\in\R$ are constants, and $\xi$ is a space-time white noise (i.e. the Gaussian space-time field with covariance kernel $\E\xi(\t,\x)\xi(\s,\y)=\delta_{\t=\s}\delta_{\x=\y}$). Let us give a simple field-theoretic explanation for where \eqref{eq:kpzbasic} comes from.
\begin{enumerate}
\item Because we look at \emph{fluctuations}, we might expect the local behavior (in $\x$) of the interface resembles that of a Gaussian free field (which is just a Brownian bridge in the case of the one-dimensional torus $\mathbb{T}$). Thus, take the corresponding probability measure, which we formally write as $\propto\exp[-\int_{\mathbb{T}}|\grad\varphi|^{2}]\d\varphi$, and its Langevin dynamic. This gives \eqref{eq:kpzbasic} with $\bar{\alpha}(\wedge)=0$.
\item The Langevin dynamic in the previous bullet point is a \emph{reversible} model for interface fluctuations, whereas the main goal of \cite{KPZ} was to write down a model for \emph{non-equilibrium} (or non-reversible) interface models. To this end, \cite{KPZ} allow $\bar{\alpha}(\wedge)\neq0$ in \eqref{eq:kpzbasic}; at the level of physics, the interface now evolves according to its local geometry. The specific choice of a quadratic function of the slope $\partial_{\x}\mathfrak{h}$ is justified with a Taylor expansion heuristic in \cite{KPZ}; see also \cite{HQ} for this heuristic.
\end{enumerate}
Again, \eqref{eq:kpzbasic} is a model for interface \emph{fluctuations}. With this in mind, it is certainly natural to ask about any \emph{universality} of \eqref{eq:kpzbasic}. For example, in point (1) above, we made a choice in modeling local behavior of the interface fluctuations by a Brownian bridge. The question we are most interested in for this paper is what happens to the interface model when we change  the Brownian bridge measure by replacing the quadratic functional $|\grad\varphi|^{2}$ with a general potential $\mathscr{U}(\grad\varphi)$. Moreover, what if $\mathscr{U}$ is time-dependent? (The motivation for this is to add another flavor of non-equilibrium to the interface fluctuations, similar to \cite{CYau}. We note, however, that \cite{CYau} deals with interface fluctuations described by \eqref{eq:kpzbasic} with $\bar{\alpha}(\wedge)=0$.) In this case, one equation, which generalizes \eqref{eq:kpzbasic} to potentials $\mathscr{U}$, that we can look at is the following stochastic PDE (where $\mathscr{U}'$ means derivative in the second input of $\mathscr{U}$):
\begin{align}
\partial_{\t}\wt{\mathfrak{h}}(\t,\x) \ = \ \partial_{\x}\mathscr{U}'(\t,\partial_{\x}\wt{\mathfrak{h}}(\t,\x))+\mathscr{U}'(\t,\partial_{\x}\wt{\mathfrak{h}}(\t,\x))+\xi(\t,\x). \label{eq:glbasic}
\end{align}
(We clarify that the first-order term in \eqref{eq:glbasic} is not the quadratic in \eqref{eq:kpzbasic} for the Gaussian potential $\mathscr{U}(\grad\varphi)=|\grad\varphi|^{2}$. However, because we allow for general classes of nonlinear $\mathscr{U}$, it does not make much difference if we take a squared version of the first-order term in \eqref{eq:glbasic} or not.) The universality claim is now that solutions to \eqref{eq:glbasic} are solutions to \eqref{eq:kpzbasic} (for appropriate coefficients depending on $\mathscr{U}(\t,\cdot)$). A key difficulty in rigorously proving this universality claim is that \eqref{eq:kpzbasic} and \eqref{eq:glbasic} are \emph{singular} SPDEs; in order to solve them, one has to use the probabilistic structure of $\xi$ (see \cite{Hai13, Hai14}). In particular, to make sense of \eqref{eq:glbasic}, we must first discretize (or otherwise regularize) $\xi$ and the equation itself, solve said discretized PDE, and take a limit of the solution as we remove the discretization. The discretized models we are left with are precisely defined as follows. 
\begin{definition}\label{definition:intro4}
 Fix a scaling parameter $\N\geq0$, which is an integer that we eventually take to $\infty$. With notation to be explained shortly, let $\mathbf{J}(\t,\x)$ (for $\t\geq0$ and $\x\in\mathbb{T}(\N):=\Z/\N\Z$) solve the following SDE:
\begin{align}
\d\mathbf{J}(\t,\x) \ = \ \N^{\frac32}\grad^{+}\mathscr{U}'(\t,\mathbf{U}^{\t,\x})\d\t + \N\{\mathscr{U}'(\t,\mathbf{U}^{\t,\x})+\mathscr{U}'(\t,\mathbf{U}^{\t,\x+1})\}\d\t + \sqrt{2}\N^{\frac12}\d\mathbf{b}(\t,\x). \label{eq:hf}
\end{align}
%
\begin{itemize}

\item The process $\mathbf{U}^{\t,\x}:=\N^{1/2}[\mathbf{J}(\t,\x)-\mathbf{J}(\t,\x-1)]$ is a rescaled gradient.
\item The processes $\t\mapsto\mathbf{b}(\t,\x)$ are jointly independent (over $\x\in\mathbb{T}(\N)$) standard Brownian motions.
\item The term $\grad^{+}\mathscr{U}'(\t,\mathbf{U}^{\t,\x}):=\mathscr{U}'(\mathbf{U}^{\t,\x+1})-\mathscr{U}'(\mathbf{U}^{\t,\x})$ is a discretization of the second-order term in \eqref{eq:glbasic}.
\end{itemize}
We isolated the gradient process $\mathbf{U}$ because of its relevance in the literature on interacting particle systems as a \emph{Ginzburg-Landau model}; see \cite{CYau,GPV,FSPDE1}. (Technically, it is a weakly asymmetric Ginzburg-Landau model as in \cite{DGP} because of the first-order term in \eqref{eq:hf}.) As noted in \cite{GPV}, this process has the interpretation of charges moving around a paramagnet. It can be checked that $\mathbf{U}$ itself satisfies the closed system of SDEs below (with notation explained after):
\begin{align}
\d\mathbf{U}^{\t,\x} \ = \ \N^{2}\Delta\mathscr{U}'(\t,\mathbf{U}^{\t,\x})\d\t + \ \N^{\frac32}\grad^{\mathrm{a}}\mathscr{U}'(\t,\mathbf{U}^{\t,\x})\d\t - \sqrt{2}\N\grad^{-}\d\mathbf{b}(\t,\x). \label{eq:glsde}
\end{align}
%
\begin{itemize}

\item The first term $\Delta\mathscr{U}'(\t,\mathbf{U}^{\t,\x}):=\mathscr{U}'(\t,\mathbf{U}^{\t,\x+1})+\mathscr{U}'(\t,\mathbf{U}^{\t,\x-1})-2\mathscr{U}'(\t,\mathbf{U}^{\t,\x})$ is a discrete Laplacian.
\item The second $\grad^{\mathrm{a}}\mathscr{U}'(\t,\mathbf{U}^{\t,\x}):=\mathscr{U}'(\t,\mathbf{U}^{\t,\x+1})-\mathscr{U}'(\t,\mathbf{U}^{\t,\x-1})$ is an asymmetric operator on $\mathscr{L}^{2}(\mathbb{T}(\N))$ acting on $\mathscr{U}'(\t,\mathbf{U}^{\t,\x})$ (hence the superscript ``$\mathrm{a}$").
\item The last term $\grad^{-}\mathbf{b}(\t,\x)=\mathbf{b}(\t,\x-1)-\mathbf{b}(\t,\x)$ is a gradient noise.
\end{itemize}
\end{definition}
\emph{One of our main results (Theorem \ref{theorem:kpz}) says that after shifting $\mathbf{J}$ in space by a ``homogenized" characteristic shift depending only {on} $\mathscr{U}$, then under some assumptions to be explained shortly, we have convergence of $\mathbf{J}$ (after rescaling space $\mathbb{T}(\N)$ to be order $1$), to the following version of \eqref{eq:kpzbasic} but with time-inhomgeneous coefficients}:
\begin{align}
\partial_{\t}\mathbf{h}^{\infty}(\t,\x) \ = \ \bar{\alpha}(\t)\partial_{\x}^{2}\mathbf{h}^{\infty}(\t,\x) + \bar{\alpha}(\t;\wedge)|\partial_{\x}\mathbf{h}^{\infty}(\t,\x)|^{2} + \xi(\t,\x). \label{eq:kpz}
\end{align}
(The $\infty$-superscript suggests this SPDE, which we call the \emph{time-inhomogeneous KPZ}, or TIKPZ, equation, as a scaling limit. {Let us also clarify that $\bar{\alpha}(\t)>0$ and $\bar{\alpha}(\t;\wedge)\in\R$.}) 
\subsubsection{The aforementioned assumptions for convergence to \eqref{eq:kpz}}
Convergence to \eqref{eq:kpz}, even in the time-homogeneous case, is a long-standing open problem; see \cite{KPZAIM}. While a proof in (near) complete generality seems to be out of reach of current methods, our goal is to show convergence at least under three assumptions, the first two of which are loosely stated as follows.
\begin{enumerate}
\item ``Mesoscopic scale" control over the relative entropy of the initial data of \eqref{eq:glsde} with respect to a discretization of the measure $\propto\exp[-\int_{\mathbb{T}}\mathscr{U}(0,\grad\varphi)]\d\varphi$ from $\mathbb{T}$ to $\mathbb{T}(\N)$.
\item A priori {Holder-($1/2-\e$)} spatial regularity of $\mathbf{J}$ on ``mesoscopic scales".
\end{enumerate}
Precise versions of these assumptions unfortunately require quite a bit of setup to state; we instead defer them to the next section (see Definition \ref{definition:entropydata} and Theorem \ref{theorem:kpz}, respectively). Intuitively, they serve the same role as the assumptions from Yau's relative entropy method \cite{YauRE}, which derives the ``hydrodynamic limit", i.e. the leading-order behavior of \eqref{eq:hf} that we study fluctuations about in this paper. (The assumptions that we take in this paper are stronger versions of the ones in Yau's relative entropy method. Indeed, we study finer-scale and more sensitive fluctuations about a hydrodynamic limit that turns out to be constant in space and time, so we need stronger a priori estimates.) We note that the second assumption could be thought of as a version of the ``locally Brownian" ansatz of \cite{KPZ} but for general $\mathscr{U}$, at least heuristically.

Let us now give the third and final (modulo technical a priori bounds) assumption:
\begin{enumerate}
\setcounter{enumi}{2}
\item The potential $\mathscr{U}(\t,\cdot)$ is uniformly convex (uniformly in $\t$).
\end{enumerate}
(For a complete and precise statement, see Assumption \ref{ass:intro8}.) Technically, this is not an immediate analog to any assumptions in Yau's method. However, it is still natural to ask for. Indeed, Yau's method is based on turning entropy bounds into hydrodynamic limits. Our goal is to do the same, but for finer-scale and more sensitive fluctuations. It turns out that this convexity assumption (modulo other technical and less interesting bounds in Assumption \ref{ass:intro8}) is enough, which is in principle not at all clear a priori (at least if one only takes, in addition to convexity, the first two assumptions above). (Actually, we will not need the full strength of convexity, but rather only the log-Sobolev inequalities that are implied by convexity and Bakry-Emery theory.)

To summarize the introduction thus far, one of our main results (Theorem \ref{theorem:kpz}) yields convergence of a spatially shifted version of $\mathbf{J}$ to \eqref{eq:kpz} under mesoscopic scale versions of the assumptions in Yau's relative entropy method. Like this relative entropy method, it provides a derivation of \eqref{eq:kpz} from interacting particle systems like \eqref{eq:hf}-\eqref{eq:glsde} based on entropy, which certainly has its conceptual benefit from the point of view of statistical mechanics. Beyond this, however, Theorem \ref{theorem:kpz} also immediately gives a convergence result (Corollary \ref{corollary:kpz}) in the time-homogeneous case that is an exponential-scale improvement of previous results in \cite{GPV}.
\subsubsection{Solving \eqref{eq:kpz}}
As we mentioned before, \eqref{eq:kpz} is a singular SPDE, so its solution theory requires some care. The approach we take in this paper is the same exponentiation of \eqref{eq:kpz} (i.e. its \emph{Cole-Hopf transform}) as in \cite{BG}. Elementary calculus then shows that the SPDE we must solve to define a solution to \eqref{eq:kpz} is given as follows.
\begin{definition}\label{definition:intro2}
 First, define $\lambda(\t):=\bar{\alpha}(\t;\wedge)/\bar{\alpha}(\t)$. Now, assume that $\bar{\alpha}(\t;\wedge)\neq0$. Set $\mathbf{h}^{\infty}(\t,\x):=\lambda(\t)^{-1}\log\mathbf{Z}^{\infty}(\t,\x)$ where 
\begin{align}
\partial_{\t}\mathbf{Z}^{\infty}(\t,\x) \ = \ \bar{\alpha}(\t)\partial_{\x}^{2}\mathbf{Z}^{\infty}(\t,\x)+\lambda(\t)\mathbf{Z}^{\infty}(\t,\x)\xi(\t,\x)+\partial_{\t}\log|\lambda(\t)|\times\mathbf{Z}^{\infty}(\t,\x)\log\mathbf{Z}^{\infty}(\t,\x). \label{eq:she1}
\end{align}
We will refer to this as the {time-inhomogeneous stochastic heat equation} (TISHE). Rigorously put, $\mathbf{Z}^{\infty}$ is adapted to the filtration generated by $\xi$, and it solves the following, which we obtain formally via the Duhamel formula (with notation explained after):
\begin{align}
\mathbf{Z}^{\infty}(\t,\x) \ &= \ \mathbf{H}(0,\t,\x)(\mathbf{Z}^{\infty}(0,\cdot)) + {\textstyle\int_{0}^{\t}}\mathbf{H}(\s,\t,\x)(\lambda(\s)\mathbf{Z}^{\infty}(\s,\cdot)\xi(\s,\cdot)) \ \d\s \label{eq:she2a}\\
&+ {\textstyle\int_{0}^{\t}}\mathbf{H}(\s,\t,\x)\left(\partial_{\s}\log|\lambda(\s)|\times\mathbf{Z}^{\infty}(\s,\cdot)\log\mathbf{Z}^{\infty}(\s,\cdot)\right)\d\s.\label{eq:she2b}
\end{align}
%
\begin{itemize}

\item For appropriate $\phi:\mathbb{T}\to\R$, let $\mathbf{H}(\s,\t,\x)(\phi)$ be the time-inhomogeneous heat semigroup associated to $\partial_{\t}-\bar{\alpha}(\t)\partial_{\x}^{2}$. More precisely, we have $\mathbf{H}(\s,\t,\x)(\phi(\cdot))=\int_{\mathbb{T}}\mathbf{H}(\s,\t,\x,\y)\phi(\y)\d\y$ where $\partial_{\t}\mathbf{H}(\s,\t,\x,\y)=\bar{\alpha}(\t)\partial_{\x}^{2}\mathbf{H}(\s,\t,\x,\y)$ for any $\s<\t$ and $\mathbf{H}(\s,\s,\x,\y)=\delta_{\x=\y}$. Also, the integral against $\xi$ is in the Ito-Walsh sense \cite{Khosh}.
\end{itemize}
\end{definition}
The main point worth noting about the Cole-Hopf transform and Definition \ref{definition:intro2} is that the exponentiated SPDE is \emph{not} linear in the solution, entirely because of the time-inhomogeneous nature of the interface model it describes. In particular, in the general time-inhomogeneous case, solving the TISHE is itself a nontrivial problem. Proving its almost sure positivity (assuming positive initial data) is also an issue; the work of \cite{Mu} relies on linearity. Ultimately, both of these issues are resolved in our other main result, Theorem \ref{theorem:she}, which requires only some smoothness and strict positivity of $\bar{\alpha}(\t),|\bar{\alpha}(\t;\wedge)|$ (i.e. boundedness of $\bar{\alpha}(\t;\wedge)$ away from $0$).
\subsection{About time-inhomogeneity}
Because the scaling limit for \eqref{eq:hf} is time-inhomogeneous, the interface models we study in this paper are in a non-equilibrium setting similar to \cite{CYau}. (We emphasize \cite{CYau} does not derive KPZ-type scaling limits, but rather scaling limits given by linear SPDEs.) We clarify, however, that the potential for the Ginzburg-Landau model in \cite{CYau} is \emph{not} time-inhomogeneous. In fact, the time-inhomogeneous nature of the limiting SPDEs therein comes from looking at perturbations of a non-constant hydrodynamic limit, so the model is out of thermal equilibrium and thus its macroscopic behavior is dynamic. In principle, we could also try to derive KPZ-type fluctuations about non-constant hydrodynamic limits, but the hydrodynamic limit would have to be the solution to a hyperbolic equation at long times, since KPZ requires a singular asymmetry (even in the weakly asymmetric scale!). Since the long-time behavior of nonlinear hyperbolic equations is not as well understood, we avoid this route. We also emphasize that because we have a time-inhomogeneous potential, the models we consider do \emph{not} have invariant measures! (Indeed, the Langevin part of \eqref{eq:glbasic} is defined with respect to a changing reference measure.) This forces us to develop new techniques for homogenization in interacting particle systems, which we explain shortly.
\subsection{Some of the main innovations}
\subsubsection{Cole-Hopf}
According to what is written in \cite{DGP,GJ15,Hai13,HQ}, there seems to have been a general consensus that Cole-Hopf is not a viable way to get \eqref{eq:kpz} from general Ginzburg-Landau models (i.e. general potentials). Our work \emph{shows} that this is not the case. In a nutshell, Cole-Hopf works for the limiting SPDE \eqref{eq:kpz} because the algebra generates for us a term that cancels the problematic quadratic in \eqref{eq:kpz}. Although it requires heavy calculations to see, the same is ``sufficiently true" for \eqref{eq:hf}. (We do not get exact cancellation of the first-order term in \eqref{eq:hf} via Cole-Hopf like we do for the limit SPDE, but we get cancellation ``in a homogenized sense", i.e. after using hydrodynamic considerations to replace local statistics by their appropriate expectations.)
\subsubsection{Second-order Boltzmann-Gibbs principle}
Let us cite \cite{DGP}, which derives \eqref{eq:kpz} from \eqref{eq:hf} in the time-homogeneous case and under much stronger assumptions than what we need in Theorem \ref{theorem:kpz} and Corollary \ref{corollary:kpz}. The method of \cite{DGP} is based on very precise and sensitive homogenization estimate known as the \emph{second-order Boltzmann-Gibbs principle}. (Roughly speaking, this says that any local statistic, in some averaged sense, is asymptotically computed by a Taylor expansion argument as in \cite{KPZ}.) However, the proof of the second-order Boltzmann-Gibbs principle in \cite{DGP} requires \eqref{eq:glsde} (in the time-homogeneous setting) to be stationary, hence the need for a time-homogeneous potential and the much stronger assumptions. 

In this paper, we are unable to relax the aforementioned assumptions in \cite{DGP} (for proving the second-order Boltzmann-Gibbs principle). On the other hand, we go through the Cole-Hopf map, and thus we only require a \emph{local} version of the second-order Boltzmann-Gibbs principle. Although difficult (and very technical), this local version can be shown under the much more relaxed assumptions that we discussed above. (In words, because we go through Cole-Hopf, the limiting SPDE we derive is much easier to treat. We are never left with the problem of making sense of a singular nonlinearity; we only have to show that certain terms go to zero in the scaling limit, which can be done by averaging over ``more local" scales.) The ingredients for showing this local principle amount to a local equilibrium technique that builds on prior work \cite{Y22}.
\subsubsection{Going beyond \eqref{eq:hf}-\eqref{eq:glsde}}
The two points discussed above (i.e. Cole-Hopf and the second-order Boltzmann-Gibbs principle) should be true for more general interacting particle systems, like exclusions and zero-range processes as in \cite{GJS15}. Showing this is the case certainly requires work, but we do not see any obstructions.
\subsection{Acknowledgements}
We thank Fraydoun Rezakhanlou for discussions and seminars, as well as Herbert Spohn for a brief discussion, pointing out the references \cite{SS,WFS}, and very helpful comments on the paper. {We would also like to give our thanks to the referees who reviewed this paper. Their comments, which significantly improved this paper, are deeply appreciated.} The author was partially supported by a fellowship from the ARCS foundation, the NSF Mathematical Sciences Postdoctoral Fellowship program under Grant. No. DMS-2203075, and the NSF under Grant. No. DMS-1928930 at a program hosted at MSRI in Berkeley, CA during Fall 2021.
%
%
%
\section{Precise statements for the main results}\label{section:main}
\subsection{Time-inhomogeneous KPZ equation}
We first present well-posedness and positivity for TISHE. First, let $\mathscr{C}^{0,\upsilon}(\mathbb{T})$ be the space of Holder continuous functions with exponent $\upsilon$, and set $\mathscr{C}(\mathbb{T}):=\mathscr{C}^{0,0}(\mathbb{T})$ for convenience.
\begin{theorem}\label{theorem:she}
 Let us first recall \eqref{eq:she2a}-\eqref{eq:she2b} and give the setting.
\begin{itemize}

\item Fix a deterministic $\upsilon\in[0,1/2)$ and {fix a possibly random $\mathbf{Z}^{\infty,\mathrm{in}}(\cdot)\in\mathscr{C}^{0,\upsilon}(\mathbb{T})$ that is independent of $\xi$ in \eqref{eq:she2a}}. 
\item Suppose that $\mathbf{Z}^{\infty,\mathrm{in}}$ is strictly positive with probability $1$.
\item {Suppose that $\bar{\alpha}(\t)>0$ and $\bar{\alpha}(\t;\wedge)\in\R$ are smooth in $\t$. Suppose also that $\bar{\alpha}(\t;\wedge)\neq0$ for all $\t\geq0$.}
\end{itemize}
{Then, with probability $1$, there exists a unique random process} $\mathbf{Z}^{\infty}$ {which takes values in $\mathscr{C}([0,\infty)\times\mathbb{T})$} and satisfies \eqref{eq:she2a}-\eqref{eq:she2b} with initial data $\mathbf{Z}^{\infty}(0,\cdot)=\mathbf{Z}^{\infty,\mathrm{in}}(\cdot)$. Moreover, we know $\mathbf{Z}^{\infty}(\t,\cdot)\in\mathscr{C}^{0,\upsilon}(\mathbb{T})$ for all $\t\geq0$, and that $\mathbf{Z}^{\infty}$ is strictly positive, both with probability 1.
\end{theorem}
A proof, modulo easy and standard technical adjustments, is in the appendix. (This paper is already quite long.)
\subsection{Deriving \eqref{eq:kpz} from \eqref{eq:hf}}
Before we present the main results (Theorem \ref{theorem:kpz} and Corollary \ref{corollary:kpz}), let us first introduce the necessary notation and preliminaries. (Because there is quite a bit of notation, we will clarify which points are crucial, and which points are defined the way they are just to make some calculations work.) The organization of what follows is:
\begin{itemize}

\item First, we introduce homogenization measures from the previous section. These will be important for everything else to follow (e.g. computing constants in the limit SPDE, computing the characteristic shift that we must introduce to see a limiting SPDE, etc.).
\item Second, we introduce spatial shifting for \eqref{eq:hf} and other normalizations that are necessary for deriving \eqref{eq:kpz}. We also define a microscopic version of the Cole-Hopf transform which will converge to \eqref{eq:she2a}-\eqref{eq:she2b}.
\item Third, we introduce the three ``Yau-type" assumptions from the previous section. We then state the results.
\end{itemize}
\subsubsection{Homogenization measures}\label{subsubsection:hommeasures}
Recall \eqref{eq:glbasic} is (a perturbation of) the Langevin dynamic for some time-dependent probability measure. The following constructs discretizations for these probability measures (but at the level of the gradient process \eqref{eq:glsde}, not \eqref{eq:hf}). Before we do so, we clarify that because \eqref{eq:glsde} is a gradient process, and because the gradient has a one-dimensional kernel spanned by constant functions, there does not exist one ``homogenization" measure for \eqref{eq:glsde} but a family parameterized by constant shifts (that is then parameterized again by time). Also, because \eqref{eq:glsde} is a gradient process, its sum over $\mathbb{T}(\N)$ is conserved with probability 1, which means that on top of these homogenizations measures, their restrictions onto hyperplanes will also be crucial. Finally, the constants introduced at the end of Definition \ref{definition:intro5} are the ones that show up in the limiting SPDE \eqref{eq:kpz}.
\begin{definition}\label{definition:intro5}
 Throughout this construction, {we assume that $\mathsf{F}\in\mathscr{C}^{\infty}(\R)$ ``does not grow too fast at $\infty$", i.e. that all of the integrals below are well-defined.}
\begin{itemize}

\item Take $\sigma\in\R$ and $\t\geq0$. Set $\E^{\sigma,\t}\mathsf{F}:=\int_{\R}\mathsf{F}(\mathbf{u})\d\mathbb{P}^{\sigma,\t}(\mathbf{u})$, where, with notation explained after, we set the following probability measure for $\mathbf{u}\in\R$:
\begin{align}
\d\mathbb{P}^{\sigma,\t}(\mathbf{u}):=\exp\{\lambda(\sigma,\t)\mathbf{u}-\mathscr{U}(\t,\mathbf{u})+\mathscr{Z}(\sigma,\t)\}\d\mathbf{u}.
\end{align}
Here, the tilt $\lambda(\sigma,\t)$ and normalization constant $\mathscr{Z}(\sigma,\t)$ are chosen so that $\d\mathbb{P}^{\sigma,\t}$ is a probability measure with mean/density $\sigma$, i.e. that $\E^{\sigma,\t}\mathbf{u}=\sigma$ and $\E^{\sigma,\t}1=1$. The probability measure $\d\mathbb{P}^{\sigma,\t}$ (and tensor products with independent copies of itself) is the time-dependent \emph{grand-canonical ensemble} of density $\sigma$.
\item For $\sigma\in\R$ and $\t\geq0$ and $\mathbb{I}\subseteq\mathbb{T}(\N)$, set $\E^{\sigma,\t,\mathbb{I}}\mathsf{F}:=\int_{\R^{\mathbb{I}}}\mathsf{F}(\mathbf{U})\d\mathbb{P}^{\sigma,\t,\mathbb{I}}(\mathbf{U})$, where $\d\mathbb{P}^{\sigma,\t,\mathbb{I}}(\mathbf{U})$ is the probability measure on $\R^{\mathbb{I}}$ obtained by conditioning the product measure $\otimes_{\mathbb{I}}\d\mathbb{P}^{0,\t}$ on $\R^{\mathbb{I}}$ on the following hyperplane of density $\sigma$:
\begin{align}
\mathbb{H}^{\sigma,\mathbb{I}} \ := \ \{(\mathbf{U}(\x))_{\x\in\mathbb{I}}: \ \tfrac{1}{|\mathbb{I}|}{\textstyle\sum_{\x\in\mathbb{I}}}\mathbf{U}(\x)=\sigma\}.
\end{align}
The probability measure $\d\mathbb{P}^{\sigma,\t,\mathbb{I}}$ is the time-dependent \emph{canonical ensemble} of density $\sigma$ on the subset $\mathbb{I}$.
\end{itemize}
Lastly, for any $\t\geq0$, we define $\bar{\alpha}(\t):=\partial_{\sigma}\E^{\sigma,\t}\mathscr{U}'(\t,\cdot)|_{\sigma=0}$ and $\bar{\alpha}(\t;\wedge):=\partial_{\sigma}^{2}\E^{\sigma,\t}\mathscr{U}'(\t,\cdot)|_{\sigma=0}$ and $\lambda(\t):=\bar{\alpha}(\t;\wedge)/\bar{\alpha}(\t)$.
\end{definition}
Let us now introduce the necessary shifts and renormalizations of \eqref{eq:hf} in order to derive {the KPZ equation. In a nutshell, all we do is to shift the discrete torus} $\mathbb{T}(\N)$ by something depending on time only, and then subtract from this shift of \eqref{eq:hf} the usual ``renormalization constant" that shows up in singular SPDEs. Both the shift and this constant are just computed to make sure that the limiting SPDE \eqref{eq:kpz} has neither an infinite speed first-order derivative nor an additional constant term on the RHS. (It is certainly important to know what these constants are as far as proving Theorem \ref{theorem:kpz} is concerned, but the reader is invited to not take them too seriously for now. Their values just fall out of calculations.)

Before we give these constructions, we note that the renormalization constant $\mathscr{R}(\s)$ appearing below (that we subtract from the so-called ``height function" below) does not diverge as $\N\to\infty$, contrary to the usual situation for singular SPDEs. The reason is because the diverging part of this renormalization constant is $\N^{1/2}\E^{0,\s}\mathscr{U}'(\s,\cdot)$. But this is equal to $0$ by calculus. This is one effect of not squaring the first-order term in \eqref{eq:glsde}.
\begin{definition}\label{definition:intro6}
 With notation explained after, we define a height function $\mathbf{h}$, Cole-Hopf (or Gartner \cite{Gar}) transform $\mathbf{Z}$, and auxiliary unshifted Cole-Hopf transform $\mathbf{G}$ as follows:
\begin{align}
\mathbf{h}(\t,\x) \ &:= \ \mathbf{J}(\t,{\y(\x,\t)})-{\textstyle\int_{0}^{\t}}\mathscr{R}(\s)\d\s \\
\mathbf{Z}(\t,\x) \ &:= \ \exp[\lambda(\t)\mathbf{h}(\t,\x)] \\
\mathbf{G}(\t,\x) \ &:= \ \exp[\lambda(\t)\mathbf{J}(\t,\x)-\lambda(\t){\textstyle\int_{0}^{\t}}\mathscr{R}(\s)\d\s] ,
\end{align}
%
\begin{itemize}

\item Recalling $\bar{\alpha}(\cdot)$ from Definition \ref{definition:intro5}, set {$\y(\x,\t):=\x-\lfloor2\N^{3/2}\int_{0}^{\t}\bar{\alpha}(\s)\d\s\rfloor$ (for any $\t\geq0$ and $\x\in\mathbb{T}(\N)$)}.
\item Recalling Definition \ref{definition:intro5}, set $\mathscr{R}(\s)=\frac{1}{12}\lambda(\s)^{3}\E^{0,\s}[\mathscr{U}'(\s,\mathbf{u})\mathbf{u}]+\frac{1}{6}\lambda(\s)^{2}\bar{\alpha}(\s)\E^{0,\s}(\mathbf{u}^{3})$.
\end{itemize}
\end{definition}
\begin{rem}\label{remark:intro7}
 The $1/12$ is exactly the $1/24$ showing up in the renormalization constant for the height function in \cite{BG} (the discrepancy in these factors comes entirely from speeding up our process \eqref{eq:hf} by a factor of $2$ compared to \cite{BG}). Also, the choice of $0$ in the expectations $\E^{0,\s}$ is unimportant; we could replace it with any other $\sigma\in\R$ as long as we subtract from the RHS of \eqref{eq:hf} a $\sigma$-dependent renormalization term (that would then diverge at speed $\N^{1/2}$, contrary to what we said prior to Definition \ref{definition:intro6}).
\end{rem}
\subsubsection{Yau-type assumptions}
We start with a precise version of the relative entropy assumption.
\begin{definition}\label{definition:entropydata}
 Take any probability density $\mathfrak{p}$ with respect to $\mathbb{P}^{0,0,\mathbb{T}(\N)}$. (Recall from Definition \ref{definition:intro5} that $\mathbb{P}^{0,0,\mathbb{T}(\N)}$ is the canonical measure obtained by conditioning the time $0$ grand-canonical product measure on $\R^{\mathbb{T}(\N)}$ to have density $0$.) We say the probability measure $\mathfrak{p}\d\mathbb{P}^{0,0,\mathbb{T}(\N)}$ is \emph{entropy data} if its relative entropy (or \emph{Kullback-Leibler divergence}) with respect to $\mathbb{P}^{0,0,\mathbb{T}(\N)}$ satisfies the following a priori estimate with {small but fixed $\gamma_{\mathrm{KL}}>0$} (which we assume is independent of $\N$):
\begin{align}
\mathfrak{D}_{\mathrm{KL}}(\mathfrak{p}) \ := \ \E^{0,0,\mathbb{T}(\N)}\mathfrak{p}\log\mathfrak{p} \ \leq \ \N^{\frac34-\gamma_{\mathrm{KL}}}.\label{eq:defentropydata}
\end{align}
\end{definition}
\begin{rem}\label{remark:intro12}
 $\mathbb{P}^{0,0,\mathbb{T}(\N)}$ is \emph{not} an invariant measure for \eqref{eq:glsde} if the potential is time-dependent. (However, it can be checked, as in \cite{DGP}, that $\mathbb{P}^{0,\t,\mathbb{T}(\N)}$ is invariant for the time-$\t$ generator of \eqref{eq:glsde}.)
\end{rem}
\begin{rem}\label{remark:intro12b}
 Yau's relative entropy method would have $1$ in place of $3/4$ {in \eqref{eq:defentropydata}}.
\end{rem}
We now introduce a quantitative notion of a priori regularity on mesoscopic scales. Recall $\gamma_{\mathrm{KL}}$ from Definition \ref{definition:entropydata}.
\begin{definition}\label{definition:reg}
 Set $\gamma_{\mathrm{reg}}:={c}\gamma_{\mathrm{KL}}$ {for some small but fixed $c>0$}. Let $|\cdot|_{\mathbb{T}(\N)}$ be geodesic distance on the torus $\mathbb{T}(\N)$. Now, define the stopping time
\begin{align}
\t_{\mathrm{reg}} \ := \ \inf\left\{\t\geq0: \sup_{\x\neq\y\in\mathbb{T}(\N)} \frac{|\mathbf{h}(\t,\x)-\mathbf{h}(\t,\y)|}{\N^{1/2}|\x-\y|_{\mathbb{T}(\N)}^{1/2}} \ \geq \ \N^{\gamma_{\mathrm{reg}}}\right\}\wedge1.
\end{align}
\end{definition}
\begin{rem}\label{remark:intro14}
 {The final time of $1$ on the RHS of the $\t_{\mathrm{reg}}$-formula is not important. It is just a convenient choice of time-horizon for \eqref{eq:hf}-\eqref{eq:glsde} in this paper. Roughly speaking, $\t_{\mathrm{reg}}$ is a blow-up time for the mesoscopic $\mathscr{C}^{0,1/2}(\mathbb{T})$ semi-norm. (By mesoscopic, all we mean is that we allow for a Holder semi-norm to blow-up as $\N^{\gamma_{\mathrm{reg}}}\gg1$. In particular, the reverse of the inequality in the definition of $\t_{\mathrm{reg}}$ is only meaningful for $|\x-\y|\gg\N^{2\gamma_{\mathrm{reg}}}$.) Let us now motivate it. First, we anticipate that $\mathbf{h}(\t,\N\x)\approx\mathrm{TIKPZ}$. By Theorem \ref{theorem:she}, we know $\mathrm{TIKPZ}$ has regularity $\mathscr{C}^{0,\upsilon}(\mathbb{T})$ for any $\upsilon<1/2$. (On the coarse-grained torus $\mathbb{T}(\N)$, this is the same as $\upsilon=1/2$ up to factors of $\log\N\ll\N^{\gamma_{\mathrm{reg}}}$.) Thus, we certainly expect $\t_{\mathrm{reg}}=1$ with high probability (in the sense of Definition \ref{definition:intro15} below). Actually, if $\mathbf{h}(\t,\N\x)$ is (close to) Brownian in $\x$, which is true for the solution to $\mathrm{TIKPZ}$, then because we allow for $\N^{\gamma_{\mathrm{reg}}}$ in $\t_{\mathrm{reg}}$, we expect that $\t_{\mathrm{reg}}=1$ with very high probability. Finally, because $\mathbf{h}$ and $\mathbf{J}$ differ by spatial shift (and constant renormalization), the supremum in the $\t_{\mathrm{reg}}$-formula does not change if we replace $\mathbf{h}$ by $\mathbf{J}$. As $\mathbf{J}(\t,\x)-\mathbf{J}(\t,\x-1)=\N^{1/2}\mathbf{U}^{\t,\x}$ (see Definition \ref{definition:intro4}), for any $\t\leq\t_{\mathrm{reg}}$ and interval $\mathbb{I}\subseteq\mathbb{T}(\N)$ (see Section \ref{subsection:notation} for what ``discrete interval" precisely means),
\begin{align}
|\mathbb{I}|^{-1}|{\textstyle\sum_{\x\in\mathbb{I}}}\mathbf{U}^{\t,\x}| \ = \ \N^{\frac12}|\mathbb{I}|^{-1}|\mathbf{J}(\t,\sup\mathbb{I})-\mathbf{J}(\t,\inf\mathbb{I})| \ \leq \ \N^{\gamma_{\mathrm{reg}}}|\mathbb{I}|^{-\frac12}. 
\end{align}
If $\mathbf{U}^{\t,\x}$ are mean-zero and sub-Gaussian increments of a random walk (bridge), then this bound would hold with exponentially high probability in $\N$ by standard concentration. (This is the $\approx$ Brownian case discussed in the first paragraph.) This remark at least says why $\t_{\mathrm{reg}}=1$ with (very) high probability is a reasonable, and likely true, assumption (in addition to motivation from Yau's method). Actually, in Corollary \ref{corollary:kpz}, we give a nontrivial situation in which it is easy to show $\t_{\mathrm{reg}}=1$ with high probability.}
\end{rem}
Finally, let us introduce a set of assumptions on the potential in \eqref{eq:hf}-\eqref{eq:glsde}. 
\begin{ass}\label{ass:intro8}
 Assume that $\mathscr{U}\in\mathscr{C}^{\infty}(\R\times\R)$ and $\partial_{\mathrm{a}}^{2}\mathscr{U}(\t,\mathrm{a})\in[\mathrm{C}^{-1},\mathrm{C}]$ uniformly in $\t,\mathrm{a}$, for some $\mathrm{C}\geq1$ independent of $\N$. Suppose $\lambda(0,\t)=0$ and $\bar{\alpha}(\t;\wedge)\neq0$ (see {Definition \ref{definition:intro5}}) and $|\partial_{\t}^{\d}\mathscr{U}|+|\partial_{\t}^{\d}\partial_{\mathrm{a}}\mathscr{U}|$ is bounded uniformly in $\t,\mathrm{a}$ and $\d=1,2$.
\end{ass}
Ultimately, all we need from the second-order $\mathrm{a}$-derivative inequalities is the time-$\t$ infinitesimal generator of \eqref{eq:glsde} to admit a log-Sobolev inequality with respect to $\mathbb{P}^{\sigma,\t,\mathbb{T}(\N)}$ that is (sufficiently close to) optimal in the length-scale $|\mathbb{T}(\N)|$. The assumption $\lambda(0,\t)=0$ is for convenience; we can shift $\mathscr{U}(\t,\mathrm{a})$ by a multiple of $\mathrm{a}$ to make it true. (Relevance of $\lambda(\sigma,\t)$ at $\sigma=0$ is explained in Remark \ref{remark:intro7}.) The assumption $\bar{\alpha}(\t;\wedge)\neq0$ guarantees KPZ-type limits. (Note $\mathbf{Z}(\t,\cdot)\equiv1$ if $\bar{\alpha}(\t;\wedge)=0$; see Definition \ref{definition:intro6}.)  The time-derivative bounds in Assumption \ref{ass:intro8} are for convenience. They can be relaxed to polynomial growth. (Smoothness is also a convenience; some finite number of derivatives is sufficient.) Finally, let us clarify that to construct potentials satisfying Assumption \ref{ass:intro8}, one can perturb any time-independent potential satisfying Assumption \ref{ass:intro8} by $\t\mathscr{V}$ with $\mathscr{V}:\R\to\R$ smooth, compactly supported, and sufficiently small in $\mathscr{C}^{2}(\R)$-norm. 
\subsubsection{Main results}
We now state our convergence result. First, a useful convention for the entire paper.
\begin{definition}\label{definition:intro15}
 Say $\mathcal{E}$ is high probability if $\mathbb{P}[\mathcal{E}^{\mathrm{C}}]\to0$ as $\N\to\infty$. Say $\mathcal{E}$ is very high probability if $\mathbb{P}[\mathcal{E}^{\mathrm{C}}]\leq\mathrm{C}(\mathrm{D})\N^{-\mathrm{D}}$ for any $\mathrm{D}>0$.  Say $\mathcal{E}$ is low probability (resp. very low probability) if $\mathcal{E}^{\mathrm{C}}$ is high probability (resp. very high probability). 
\end{definition}
\begin{theorem}\label{theorem:kpz}
 Let us first introduce a set of assumptions.
\begin{itemize}

\item Suppose $\mathbf{h}(0,0)=\mathbf{h}(0,|\mathbb{T}(\N)|)$ (so that $\mathbf{h}(0,\cdot)$ is periodic). Also, suppose $\mathbf{h}(0,\N\cdot)$ converges uniformly on $\mathbb{T}$ in probability to a limit $\mathbf{h}^{\infty}(0,\cdot)$ (where $\mathbf{h}(0,\N\x)$ extends from $\x\in\N^{-1}\mathbb{T}(\N)$ to $\x\in\mathbb{T}$ by linear interpolation). Finally, suppose $\mathbf{h}^{\infty}(0,\cdot)$ is independent of the Brownian motions in \eqref{eq:hf}.
\item Suppose that the law of \eqref{eq:glsde} at $\t=0$ is entropy data {(see Definition \ref{definition:entropydata})}.
\item Suppose $\t_{\mathrm{reg}}=1$ with very high probability, or, if the potential $\mathscr{U}$ is independent of time, that $\t_{\mathrm{reg}}=1$ with high probability.
\end{itemize}
{Then, there exists a coupling between $\mathbf{Z}$ and the solution $\mathbf{Z}^{\infty}$} to \eqref{eq:she2a}-\eqref{eq:she2b} (for an appropriate noise $\xi$) with initial data $\exp[\lambda(0)\mathbf{h}^{\infty}(0,\cdot)]$ so that $\mathbf{Z}(\t,\N\x)-\mathbf{Z}^{\infty}(\t,\x)\to0$ uniformly in $(\t,\x)\in[0,1]\times\mathbb{T}(\N)$ in probability as $\N\to\infty$.
\end{theorem}
The assumption of $\t_{\mathrm{reg}}=1$ with very high probability (e.g. instead of high probability) is just for convenience. We can stop the $(\mathbf{U},\mathbf{J})$ process at $\t_{\mathrm{reg}}$ and replace all our analysis for $(\mathbf{U},\mathbf{J})$ by the same analysis for the stopped process. All this does is modify notation. In particular, parallel to Yau's method \cite{YauRE}, our work holds until a mesoscopic blow-up time $\t_{\mathrm{reg}}$, which should equal the time-horizon 1 with high probability.

Now, assume $\mathscr{U}(\t,\mathrm{a})$ is constant in $\t$. The canonical ensembles in Definition \ref{definition:intro5} become honest invariant measures for \eqref{eq:glsde}; see \cite{DGP}. With invariant measures at our disposal, one can say a lot more. (In what follows, $\e_{\mathrm{Corollary}}$ can be computed and optimized explicitly. The proof of Corollary \ref{corollary:kpz}, given in the final non-appendix section of this paper, shows that we can probably take $\e_{\mathrm{Corollary}}$ to be anything strictly less than $1$, for example, though optimizing it is not as high interest for us as much as its strict positivity is.)
\begin{corollary}\label{corollary:kpz}
 Let us first introduce a set of assumptions.
\begin{itemize}

\item Suppose $\mathbf{h}(0,0)=\mathbf{h}(0,|\mathbb{T}(\N)|)$ (so that $\mathbf{h}(0,\cdot)$ is periodic). Also, suppose $\mathbf{h}(0,\N\cdot)$ converges uniformly on $\mathbb{T}$ in probability to a limit $\mathbf{h}^{\infty}(0,\cdot)$ (where $\mathbf{h}(0,\N\x)$ extends from $\x\in\N^{-1}\mathbb{T}(\N)$ to $\x\in\mathbb{T}$ by linear interpolation). Finally, suppose $\mathbf{h}^{\infty}(0,\cdot)$ is independent of the Brownian motions in \eqref{eq:hf}. (This is the same first bullet point as in Theorem \ref{theorem:kpz}.)
\item Let $\mathfrak{p}$ denote the density of the law of \eqref{eq:glsde} at time $0$ with respect to $\mathbb{P}^{0,0,\mathbb{T}(\N)}$ (see {Definition \ref{definition:intro4}}). Assume $\mathfrak{D}_{\mathrm{KL}}(\mathfrak{p})\leq\N^{\alpha_{\mathrm{KL}}}$, where $\alpha_{\mathrm{KL}}\leq\e_{\mathrm{Corollary}}\gamma_{\mathrm{reg}}$ for $\e_{\mathrm{Corollary}}>0$ sufficiently small but independent of $\N$. (Recall $\gamma_{\mathrm{reg}}$ from {Definition \ref{definition:reg}}.)
\end{itemize}
{Then, there exists a coupling between $\mathbf{Z}$ and the solution $\mathbf{Z}^{\infty}$} to \eqref{eq:she2a}-\eqref{eq:she2b} (for an appropriate noise $\xi$) with initial data $\exp[\lambda(0)\mathbf{h}^{\infty}(0,\cdot)]$ so that $\mathbf{Z}(\t,\N\x)-\mathbf{Z}^{\infty}(\t,\x)\to0$ uniformly in $(\t,\x)\in[0,1]\times\mathbb{T}(\N)$ in probability as $\N\to\infty$.
\end{corollary}
We show Corollary \ref{corollary:kpz} right before the appendix. \cite{DGP} shows Corollary \ref{corollary:kpz} (with weak convergence) if the initial relative entropy is bounded uniformly and independently of $\N$. So, Corollary \ref{corollary:kpz} is an exponential scale improvement of \cite{DGP}. (Indeed, relative entropy is on a log-scale. Moreover, \cite{DGP} only studies time-homogeneous SDEs; it has no version of Theorem \ref{theorem:kpz}.)
\subsection{Previous literature}
We now spend some time putting the introduction and our results in context.
\subsubsection{Derivations of KPZ}
To the author's knowledge, Theorem \ref{theorem:kpz} provides a first (general) derivation of TIKPZ in statistical mechanics processes. There has been tremendous effort in deriving the time-homogeneous KPZ equation, however. For starters, in the seminal paper \cite{KPZ}, Kardar, Parisi, and Zhang justify universality of the KPZ equation by a formal renormalization group heuristic. As for rigorous arguments, in \cite{BG} (see also \cite{CGST,CST,CT}), the time-homogeneous KPZ equation was derived from the height function associated to the ``ASEP". But the models in \cite{BG,CGST,CST,CT} miraculously have enough algebraic structure so that all Boltzmann-Gibbs principles can be completely avoided by coincidence. In particular, the work of \cite{BG,CGST,CST,CT}, in principle, does not extend beyond a very small set of highly special models. In \cite{DT,Y,Y20,Y22}, some perturbations of the \cite{BG} models were successfully treated as well. In this case, only a small (but still challenging) part of a Boltzmann-Gibbs principle is needed.

In \cite{GJ15}, Goncalves-Jara attacked the problem from an entirely different angle that is now known as energy solution theory. It casts the KPZ equation as a nonlinear martingale problem; this makes deriving KPZ a usual issue of convergence of martingale problems. In \cite{DGP}, Diehl-Gubinelli-Perkowski used this idea (and an additional uniqueness result for the martingale problem in \cite{GP}) to prove Corollary \ref{corollary:kpz} under much stronger assumptions. The point is that to run the estimates in the martingale problem theory, one needs a globally stationary process. (In particular, time-inhomogeneous processes are out of the question.) However if one takes this assumption of stationarity, then one can rigorously prove a very strong version of the Boltzmann-Gibbs principle. (Namely, the relevant replacement lemmas in \cite{DGP,GJ15} hold even for small macroscopic length scales.) 

{In \cite{HQ,HS,HX}, regularity structures were used to derive a time-homogeneous KPZ equation from stochastic Hamilton-Jacobi equations (as a toy model for statistical physics processes like Ginzburg-Landau). We now expand on (serious) obstructions in using regularity structures to study the coarse-grained SPDE \eqref{eq:hf}. \cite{HQ,HS,HX} compares (in terms of local regularity) singular SPDEs like KPZ to Gaussian ones like \eqref{eq:kpz} for $\bar{\alpha}(\t;\wedge)=0$. This is why their homogenized coefficients are determined by Gaussian expectations. But for \eqref{eq:hf}-\eqref{eq:glsde}, homogenized coefficients are given by non-Gaussian expectations! This can be felt immediately at a technical level. Indeed, regularity structures depend crucially on optimal regularity for the heat operator in \eqref{eq:kpz}. For \eqref{eq:hf}, even with $\partial_{\mathrm{a}}^{2}\mathscr{U}$ bounds in Assumption \ref{ass:intro8}, the optimal regularity for the fully nonlinear heat operator in \eqref{eq:hf} are the much-too-weak DeGiorgi-Nash-Moser estimates. The homogenized coefficients are intimately connected to said nonlinear operator in \eqref{eq:hf}, so one cannot avoid this issue. (In this way, \eqref{eq:hf}-\eqref{eq:glsde} is very different than models in \cite{HQ,HS,HX}.) Resolving these issues would be very exciting.}

{Besides the aforementioned works on particle systems, we also mention the work \cite{AKQ}, which derives the KPZ equation as a scaling limit for fluctuations of the free energy in a large class of directed polymers. Let us also mention work \cite{CG,DDP2,P2} on deriving the KPZ equation from random walks in random environments, work \cite{DDP1} on deriving the KPZ equation from sticky Brownian motion, as well as work \cite{AC} on deriving the KPZ equation from a class of discrete randomly growing surfaces.}
\subsubsection{Time-inhomogeneous objects of KPZ-type}
The work \cite{BLSZ} solves a time-inhomogeneous generalization of TASEP. Transition kernels are computed by generalizing the method of \cite{MQR}. In \cite{BDR}, the authors compute long-time fluctuations (under a KPZ fixed point scaling) of a KPZ equation with time-inhomogeneous noise. (No log-nonlinearity shows up in its Cole-Hopf map.)
\subsubsection{Well-posedness of TIKPZ and TISHE}
In \cite{DKZ}, the authors prove well-posedness for some nonlinear time-homogeneous stochastic heat equations (SHE). But in \cite{DKZ}, well-posedness means non-explosion. Positivity is not derived in \cite{DKZ}; this is crucial to define TIKPZ via TISHE. We also note that even if regularity structures were applied to TIKPZ itself (which certainly seems possible), we would only get local well-posedness. (Indeed, in \cite{Hai13} global well-posedness for KPZ comes from that for SHE!)
\subsection{A word about the writing of this paper}\label{subsection:writing}
We make a few disclaimers about the writing style of this paper. 
\begin{itemize}

\item The work in this paper is very technical. For this purpose, we always explain proofs in words to supplement any actual mathematics. {In many cases, these intuitive supplements are legitimate proofs, and the ``actual proofs" are just putting the supplements in terms of notation.}
\item In order to write the proof in a more organized fashion, we must unfortunately introduce quite a bit of notation. To remedy this, we always explain things (e.g. calculations) using words and refer back to where pieces of notation are defined. (This way, the reader does not have to remember every piece of notation.) For example, in Section \ref{section:derivesde}, there are necessary and long blocks of calculations. Either immediately before or after each block, we explain the calculation in words and in detail. (With regards to Section \ref{section:derivesde}, in particular, the purpose of this section is to prove a result that in some sense ``has to be true"; see the paragraph right before Proposition \ref{prop:method2}. The appendix has many sections of a similar spirit.) {We also provide a glossary (see Section \ref{section:glossary}) to streamline a lot of the notation in this paper. For example, we clarify that certain types of notation (according to font, input variables, etc.) are meant to designate certain types of objects (e.g. objects related to averaging on different space-time scales, objects related to the potential $\mathscr{U}$, etc.).}
\item {There is a lot of necessary power-counting in the parameter $\N$ in this paper. For this reason, let us clarify that constants $\mathrm{C},\mathrm{D}>0$ that are multiplying exponents in powers of $\N$ will be thought of as big constants, while $\mathrm{c}$ will be thought of as a small constant.}
\end{itemize}
\subsection{Organization and reading this paper}\label{subsection:reading}
Since this paper has a lot of moving parts, we now give an outline for how to read this paper to get the main ideas. Although this paper is long, as we explain in the paragraph below, reading for the main ideas cuts the length of this paper dramatically. Everything else amounts to technical calculations to make the ideas (and precise heuristics for the key estimates) rigorous (i.e. execute all the power-counting, and not much, if anything, else).

Let us be more precise about the writing of this paper (and how to read for the main ideas). Section \ref{section:proofoutline} gives all of the necessary big steps. No proofs will be given in this section, but detailed explanations for each result, why it is true, and what it is saying are given. Everything until (and including) Section \ref{section:bg22proofoutline} spells out the remaining ideas (with precise heuristics, which are themselves honest proofs in many cases, at least modulo power-counting in $\N$). Everything after Section \ref{section:bg22proofoutline}, but before the appendix, executes this power-counting. In the appendix, we prove Theorem \ref{theorem:she}, and we make rigorous heuristics for results that are rather standard to prove. (The rest of the appendix is also either standard or a simple extension of previous results in the literature.) 
\subsection{Some ubiquitous notation}\label{subsection:notation}
Below, we collect notation used throughout this paper (in virtually every section). 
\begin{itemize}

\item For any $\t\geq0$ and space set $\mathbb{K}$, let $\|\|_{\t;\mathbb{K}}$ be the $\mathscr{L}^{\infty}$-norm over $[0,\t]\times\mathbb{K}$. 
\item For any $\mathrm{a}\leq\mathrm{b}\in\R$, set $\llbracket\mathrm{a},\mathrm{b}\rrbracket:=[\mathrm{a},\mathrm{b}]\cap\Z$. Also, by $\Phi\lesssim\Psi$, we mean $|\Phi|=\mathrm{O}[|\Psi|]$, and by $\Phi\gtrsim\Psi$, we mean $|\Psi|\lesssim|\Phi|$. 
\item By a ``discrete interval" in $\mathbb{T}(\N)$, we mean the intersection of an interval in $\N\mathbb{T}$ and $\mathbb{T}(\N)$. By $\inf$ of a discrete interval $\mathbb{I}$, we mean the unique point $\x\in\mathbb{I}$ such that any other point in $\mathbb{I}$ can be reached by traveling from $\x$ in the positive-orientation direction. (Positive orientation means to the right upon identifying $\mathbb{T}(\N)\simeq\llbracket0,\N-1\rrbracket$.) The $\sup$ of $\mathbb{I}$ is the unique point that can be reached in this way from any other point in $\mathbb{I}$. (This is only meant to clarify discrete intervals which loop $\N-1\mapsto0$ because of the periodic boundary of $\mathbb{T}(\N)$.) 
\item Let us extend the gradient notation in Definition \ref{definition:intro4}. Given any function $\phi:\mathbb{T}(\N)\to\R$ and any $\mathfrak{l}\in\Z$, define {$\grad^{\mathfrak{l}}\phi(\x)=\phi(\x+\mathfrak{l})-\phi(\x)$}. Set $\grad^{+}:=\grad^{1}$ and $\grad^{-}:=\grad^{-1}$ and $\grad^{\mathrm{a}}=\grad^{+}-\grad^{-}$. Finally, set $\Delta=\grad^{+}+\grad^{-}$.
\end{itemize}
%
%
%
\section{Outline for the proof of Theorem \ref{theorem:kpz}}\label{section:proofoutline}
The proof of Theorem \ref{theorem:kpz} has many parts. We provide an outline in hopes of clarifying the main points. This section contains precisely stated ingredients. We use them to prove Theorem \ref{theorem:kpz}. We then spend the rest of this paper deriving every ingredient. Before we start, we invite the reader to take a quick look at Section \ref{subsection:notation} (in particular, the gradient notation therein).
\subsection{An SDE for $\mathbf{Z}$}
The SDE for $\mathbf{Z}$ can be computed via its explicit form (as a function of \eqref{eq:hf}-\eqref{eq:glsde}) and the Ito formula. It is given in Proposition \ref{prop:method2}, which we clarify shortly. First, important notation (that we also clarify immediately after stating it); see Section \ref{subsection:notation} for preliminary notation.
\begin{definition}\label{definition:method1}
 First, we let $\mathbb{J}$ denote all jump times $\t\geq0$ of the characteristic-shift $\t\mapsto\lfloor2\N^{3/2}\int_{0}^{\t}\bar{\alpha}(\s)\d\s\rfloor$. Define $\delta(\t\in\mathbb{J})$ to be the measure on $\t\geq0$ given by placing a unit Dirac point mass at every point in $\mathbb{J}$. (Note $\mathbb{J}$ and $\delta(\t\in\mathbb{J})$ are deterministic.) Now, let $\mathbf{H}^{\N}(\s,\t,\x,\y)$ solve $\mathbf{H}^{\N}(\s,\s,\x,\y)=\mathbf{1}(\x=\y)$ and $\partial_{\t}\mathbf{H}^{\N}(\s,\t,\x,\y)=\mathscr{T}(\t)\mathbf{H}^{\N}(\s,\t,\x,\y)$, where $\mathscr{T}(\t)$ is the following time-inhomogeneous discrete differential operator-valued measure acting on $\mathbf{H}^{\N}(\s,\t,\x,\y)$ through the $\x$-variable:
\begin{align}
\mathscr{T}(\t) \ := \ \N^{2}\bar{\alpha}(\t)\Delta + \tfrac14\N\lambda(\t)^{2}\bar{\alpha}(\t)\Delta+\N^{\frac32}\bar{\alpha}(\t)\grad^{+}-\N^{\frac32}\bar{\alpha}(\t)\grad^{-}+\delta(\t\in\mathbb{J})\grad^{-}. \label{eq:method1IT}
\end{align}
For $\phi:\mathbb{T}(\N)\to\R$, let $\mathbf{H}^{\N}(\s,\t,\x)(\phi(\cdot))$ solve 
\begin{align}
\partial_{\t}\mathbf{H}^{\N}(\s,\t,\x)(\phi(\cdot))=\mathscr{T}(\t)\mathbf{H}^{\N}(\s,\t,\x)(\phi(\cdot))
\end{align}
and $\mathbf{H}^{\N}(\s,\s,\x)(\phi(\cdot))=\phi(\x)$. In particular, we have
\begin{align}
\mathbf{H}^{\N}(\s,\t,\x)(\phi(\cdot)) \ := \ \sum_{\y\in\mathbb{T}(\N)}\mathbf{H}^{\N}(\s,\t,\x,\y)\phi(\y).
\end{align}
\end{definition}
To be completely clear, the heat kernel $\mathbf{H}^{\N}(\s,\t,\x,\y)$ is the solution to the integrated equation corresponding to the differential shorthand in Definition \ref{definition:method1}. (In said integrated equation, the $\delta$-function term in $\mathscr{T}(\t)$ becomes a sum over times in $\mathbb{J}$, which is clearly defined.) The $\mathbf{H}^{\N}$ kernel is also the transition density for a random walk, which jumps according to Poisson clocks giving the first four terms in $\mathrm{RHS}\eqref{eq:method1IT}$ as well as at the deterministic set of times $\mathbb{J}$ to the left by 1. (The usual probabilistic intuitions, like non-negativity of $\mathbf{H}^{\N}$ and the fact that it is a probability measure on $\mathbb{T}(\N)$ in its forward variable, therefore hold.) Since the jump speeds of said random walk are independent of its position (see \eqref{eq:method1IT}), the $\mathbf{H}^{\N}$-semigroup factors as the semigroup whose generator equals the first four terms in $\mathrm{RHS}\eqref{eq:method1IT}$ composed with the semigroup corresponding to the deterministic characteristic in Definition \ref{definition:method1}; see Proposition \ref{prop:hke}. (So, to remove the $\delta$-function in \eqref{eq:method1IT}, replace $\x$ by its image under said characteristic.)

The following result is a little involved (in terms of notation), so let us clarify it first. We show that $\mathbf{Z}$ satisfies a discretization of the continuum SHE \eqref{eq:she1}. There is a multiplicative error ${\mathfrak{z}}$, which requires some work to write down, but its key ingredients are comparing $\mathscr{U}'$ to an appropriate ``quadratic" (see \eqref{eq:method2IIb}-\eqref{eq:method2IIc}), and comparing the second-order operator in \eqref{eq:glsde} to a Laplacian (see \eqref{eq:method2IId}). Everything else is more or less technical (and perhaps uninteresting). We emphasize that this SDE for $\mathbf{Z}$ is essentially a consequence of the fact that its scaling limit is \eqref{eq:she1}, i.e. it is equivalent to universality. In particular, in some sense, it \emph{has to be true} (otherwise previous works on universality, e.g. \cite{DGP}, would likely be false). We clarify this point further after Lemma \ref{lemma:method4} (namely the homogenization heuristics that were, in another guise, important in \cite{DGP}).
\begin{prop}\label{prop:method2}
 With notation to be explained afterwards, we have the following SDE:
\begin{align}
\d\mathbf{Z}(\t,\x)&=\mathscr{T}(\t)\mathbf{Z}(\t,\x)\d\t + \sqrt{2}\lambda(\t)\N^{\frac12}\mathbf{Z}(\t,\x)\d\mathbf{b}(\t,{\y(\x,\t)})\nonumber\\
&+ [\partial_{\t}\log|\lambda(\t)|]\mathbf{Z}(\t,\x)\log\mathbf{Z}(\t,\x)\d\t + {\mathfrak{z}}(\t,{\y(\x,\t)})\mathbf{Z}(\t,\x)\d\t. \label{eq:method2I}
\end{align}
\eqref{eq:method2I} is the usual shorthand for the corresponding integrated equation, in which the $\delta$-function term in $\mathscr{T}(\t)\mathbf{Z}(\t,\x)$ makes clear sense as a sum of discrete gradients of $\mathbf{Z}(\t,\x)$ over times $\t\in\mathbb{J}$. 

{The formula for ${\mathfrak{z}}$ is given by the following display, in which
\begin{itemize}
\item we set $\mathscr{W}(\t,\mathbf{u}):=\mathscr{U}(\t,\mathbf{u})-\tfrac12\bar{\alpha}(\t)\mathbf{u}^{2}$;
\item the quantities ${\mathfrak{a}}^{\pm},{\mathfrak{b}},{\mathfrak{c}}$ are possibly random terms such that $|{\mathfrak{a}}^{\pm}(\t,\x)|+|{\mathfrak{b}}(\t,\x)|+|{\mathfrak{c}}(\t,\x)|\lesssim \ \N^{10\gamma_{\mathrm{reg}}}$ for all $\t\leq\t_{\mathrm{reg}}$
\item we define the terms $\mathbf{U}(\t):=\mathbf{U}^{\t,{\y(\x,\t)}}$ and $\mathbf{V}(\t):=\mathbf{U}^{\t,{\y(\x,\t)}+1}$ for convenience, and we define $\E^{0,\t}\mathsf{F}(\mathbf{u}):=\int_{\R}\mathsf{F}(\mathbf{u})\d\mathbb{P}^{0,\t}(\mathbf{u})$: 
\end{itemize}
}%
\begin{align}
&{\mathfrak{z}}(\t,{\y(\x,\t)})\mathbf{Z}(\t,\x)\d\t \nonumber\\
&:= \ \lambda(\t)\N\{\mathscr{U}'(\t,\mathbf{U}(\t))-\bar{\alpha}(\t)\mathbf{U}(\t)-\tfrac12\lambda(\t)[\mathscr{U}'(\t,\mathbf{U}(\t))\mathbf{U}(\t)-1]\}\mathbf{Z}(\t,\x)\d\t \label{eq:method2IIb}\\
&+ \ \lambda(\t)\N\{\mathscr{U}'(\t,\mathbf{V}(\t))-\bar{\alpha}(\t)\mathbf{V}(\t)-\tfrac12\lambda(\t)[\mathscr{U}'(\t,\mathbf{V}(\t))\mathbf{V}(\t)-1]\}\mathbf{Z}(\t,\x)\d\t \label{eq:method2IIc}\\
&+ \ \tfrac12\lambda(\t)\ \N^{\frac32}\grad^{+}\{\mathscr{W}'(\t,\mathbf{U}(\t))\mathbf{Z}(\t,\x)\}\d\t - \tfrac12\lambda(\t)\N^{\frac32}\grad^{-}\{\mathscr{W}'(\t,\mathbf{V}(\t))\mathbf{Z}(\t,\x)\}\d\t \label{eq:method2IId}\\
&+ \ \tfrac14\lambda(\t)^{2}\N\grad^{+}\{[\bar{\alpha}(\t)\mathbf{U}(\t)^{2}-1]\mathbf{Z}(\t,\x)\}\d\t + \tfrac14\lambda(\t)^{2}\N\grad^{-}\{[\bar{\alpha}(\t)\mathbf{V}(\t)^{2}-1]\mathbf{Z}(\t,\x)\}\d\t \label{eq:method2IIe}\\
&+ \ \tfrac{1}{12}\lambda(\t)^{4}\{\E^{0,\t}[\mathscr{U}'(\t,\mathbf{u})\mathbf{u}^{3}]\}\mathbf{Z}(\t,\x)\d\t + \tfrac16\lambda(\t)^{3}\bar{\alpha}(\t)\mathbf{V}(\t)^{3}\mathbf{Z}(\t,\x)\d\t - \lambda(\t)\mathscr{R}(\t)\mathbf{Z}(\t,\x)\d\t \label{eq:method2IIf} \\
&+ \ \tfrac{1}{12}\lambda(\t)^{4}\{\mathscr{U}'(\t,\mathbf{V}(\t))\mathbf{V}(\t)^{3}-\E^{0,\t}[\mathscr{U}'(\t,\mathbf{u})\mathbf{u}^{3}]\}\mathbf{Z}(\t,\x)\d\t \label{eq:method2IIg}\\
&+ \ \N^{\frac12}\grad^{-}\{{\mathfrak{a}}^{-}(\t,{\y(\x,\t)})\mathbf{Z}(\t,\x)\}\d\t + \ \N^{\frac12}\grad^{+}\{{\mathfrak{a}}^{+}(\t,{\y(\x,\t)})\mathbf{Z}(\t,\x)\}\d\t \label{eq:method2IIi} \\
&+ \ \N^{-\frac12}{\mathfrak{b}}(\t,{\y(\x,\t)})\mathbf{Z}(\t,\x)\d\t + \ \N^{-\frac12}{\mathfrak{c}}(\t,{\y(\x,\t)})\mathbf{Z}(\t,\x)\d\t. \label{eq:method2IIj}
\end{align}
Therefore, by the Duhamel principle we obtain the following from \eqref{eq:method2I}:
\begin{align}
\mathbf{Z}(\t,\x) \ &= \ \mathbf{H}^{\N}(0,\t,\x)(\mathbf{Z}(0,\cdot)) + {\textstyle\int_{0}^{\t}}\mathbf{H}^{\N}(\s,\t,\x)(\sqrt{2}\lambda(\s)\N^{\frac12}\mathbf{Z}(\s,\cdot)\d\mathbf{b}(\s,\cdot(\s))) \label{eq:method2III}\\
&+ \ {\textstyle\int_{0}^{\t}}\mathbf{H}^{\N}(\s,\t,\x)(\{\partial_{\s}\log|\lambda(\s)|\}\mathbf{Z}(\s,\cdot)\log\mathbf{Z}(\s,\cdot))\d\s + {\textstyle\int_{0}^{\t}}\mathbf{H}^{\N}(\s,\t,\x)({\mathfrak{z}}(\s,\cdot(\s))\mathbf{Z}(\s,\cdot))\d\s.
\end{align}
\end{prop}
In \eqref{eq:method2IIj} we can certainly combine ${\mathfrak{b}}$ and ${\mathfrak{c}}$. But we have written it in this way to hopefully clarify the proof of Proposition \ref{prop:method2}. It is not hard to see that $\mathscr{T}(\t)\approx\N^{2}\bar{\alpha}(\t)\Delta$ (in a heat kernel sense); see Proposition \ref{prop:hkecont}. So, as we explained earlier, Proposition \ref{prop:method2} says $\mathbf{Z}$ solves a discrete \eqref{eq:she1} with error ${\mathfrak{z}}$. To explain it better, we first make precise what it means to say ${\mathfrak{z}}$ is an error. First, we give a definition. It essentially classifies local statistics according to ``germs" of their homogenized versions. As we explain after Lemma \ref{lemma:method4}, this construction determines the leading-order asymptotics of any local statistic.
\begin{definition}\label{definition:method3}
 First, recall the notation of Definition \ref{definition:intro5}.
\begin{itemize}

\item Define $\mathrm{CT}$ to be the set of all ``centered" terms, which we define to be functions $\mathsf{F}:\R_{\geq0}\times\R\to\R$ for which $\E^{\sigma,\t}\mathsf{F}(\t,\cdot)$ is smooth in $(\sigma,\t)$ and such that $\E^{0,\t}\mathsf{F}(\t,\cdot)=0$ for all $\t\geq0$.
\item Define $\mathrm{LCT}$ to consist of all ``linearly centered/corrected" terms, which we define to be functions $\mathsf{F}:\R_{\geq0}\times\R\to\R$ for which $\E^{\sigma,\t}\mathsf{F}(\t,\cdot)$ is smooth in $(\sigma,\t)$ and such that $\E^{0,\t}\mathsf{F}(\t,\cdot)=\partial_{\sigma}\E^{\sigma,\t}\mathsf{F}(\t,\cdot)|_{\sigma=0}=0$.
\item Define $\mathrm{QCT}$ to consist of all ``quadratically centered/corrected" terms, which we define to be functions $\mathsf{F}:\R_{\geq0}\times\R\to\R$ for which $\E^{\sigma,\t}\mathsf{F}(\t,\cdot)$ is smooth in $(\sigma,\t)$ and such that $\E^{0,\t}\mathsf{F}(\t,\cdot)=\partial_{\sigma}\E^{\sigma,\t}\mathsf{F}(\t,\cdot)|_{\sigma=0}=\partial_{\sigma}^{2}\E^{\sigma,\t}\mathsf{F}(\t,\cdot)|_{\sigma=0}=0$.
\end{itemize}
\end{definition}
\begin{lemma}\label{lemma:method4}
 Recall the notation of {Proposition \ref{prop:method2}}. As functions evaluated at $(\t,\mathbf{U})$:
\begin{itemize}

\item After dividing by $\mathbf{Z}(\t,\x)\d\t$, \eqref{eq:method2IIb} belongs to $\mathrm{QCT}$. The same is true for \eqref{eq:method2IIc}.
\item The function $\mathscr{W}'(\t,\mathbf{U})$ in \eqref{eq:method2IId} belongs to $\mathrm{LCT}$.
\item The term $\bar{\alpha}(\t)\mathbf{U}^{2}-1$ in $\grad$-operators in \eqref{eq:method2IIe} belongs to $\mathrm{CT}$. After dividing by $\mathbf{Z}(\t,\x)$, \eqref{eq:method2IIf} and \eqref{eq:method2IIg} each belong to $\mathrm{CT}$.
\end{itemize}
\end{lemma}
Definition \ref{definition:method3} and Lemma \ref{lemma:method4} are \emph{algebraic} conditions. They depend only on the grand-canonical measures in Definition \ref{definition:intro5}, not dynamics of \eqref{eq:hf}-\eqref{eq:glsde}. In particular, Lemma \ref{lemma:method4} follows from calculus. (Grand-canonical measures for \eqref{eq:glsde} are sufficiently simple to compute expectations appearing {in} Definition \ref{definition:method3}. They are also nontrivial enough to make the computations interesting.) We now explain what Lemma \ref{lemma:method4} says (i.e., why exactly ${\mathfrak{z}}=\mathrm{o}(1)$). The heuristic from \cite{KPZ} says that after homogenizing, say, the local statistic \eqref{eq:method2IIb}, its contribution in the limit TISHE is given by its degree $\d=0,1,2$ derivatives. Because \eqref{eq:method2IIb} is $\mathrm{QCT}$, this means its contribution is zero. We can be more precise about this. Heuristically, $\mathrm{CT}$ statistics are fluctuating. So ``averaging" (or ``homogenizing") something in $\mathrm{CT}$ with respect to length-scale $\mathfrak{l}>0$ gives $\mathrm{O}(\mathfrak{l}^{-1/2})$. This is a CLT-type statement. It turns out that, by formal reasoning via \emph{Boltzmann-Gibbs principles}, or more precisely via local CLTs, the $\mathrm{k}$-th leading-order correction to this CLT-type estimate is given by $\mathrm{k}$-th derivatives in Definition \ref{definition:method3}. (See Section 4 of \cite{GJ15} for more precise versions of this statement.) The $\mathrm{k}$-th order correction is formally $\mathrm{O}(\mathfrak{l}^{-\mathrm{k}/2})$; again, see \cite{GJ15}. In the case of $\mathrm{QCT}$, we know $\mathrm{k}=0,1,2$-corrections are zero. So, $\mathrm{QCT}$-statistics are formally $\mathrm{O}(\mathfrak{l}^{-3/2})$ in a homogenized sense. Assuming that one can homogenize on $\mathfrak{l}\approx\e\N$ (for $\e>0$ fixed), this beats the $\N$-factor in \eqref{eq:method2IIb}-\eqref{eq:method2IIc}. Thus \eqref{eq:method2IIb}-\eqref{eq:method2IIc} are errors in a homogenized sense. Similar formal arguments for $\mathrm{LCT}$ statistics and $\mathrm{CT}$ statistics then show that the rest of ${\mathfrak{z}}$ are also errors in the same homogenized sense. We conclude this by noting that we do not need to homogenize on macroscopic scales $\mathfrak{l}\approx\e\N$ but only mesoscopic scales $\mathfrak{l}\approx\N^{2/3+\e}$. This is why our homogenization analysis can be done on local space-time scales. In particular, this last observation is crucial to our work.
\subsection{Proof of Theorem \ref{theorem:kpz} assuming a priori estimates}
Our goal now is to get that, assuming a priori (stochastic) estimates for ${\mathfrak{z}}\mathbf{Z}$ in \eqref{eq:method2I}, Theorem \ref{theorem:kpz} follows. (This is the ``rough paths" philosophy mentioned near the end of the introduction.)
\subsubsection{Preliminary smoothing}
Since $\mathbf{h}$ has a priori regularity before time $\t_{\mathrm{reg}}$, so does $\mathbf{Z}$. So, we can mollify $\mathbf{Z}$ by mesoscopic smoothing scale while changing it by only $\mathrm{o}(1)$. This ends up being very convenient for technical reasons in this paper (such as dealing with singularities of heat kernels in our analysis of \eqref{eq:method2III}), though it is perhaps unnecessary.
\begin{definition}\label{definition:method5}
 Set $\mathbf{S}(\t,\x):=\mathbf{H}^{\N}(\t,\t(\N),\x)\{\mathbf{Z}(\t,\cdot)\}$, where $\t\mapsto\t(\N):=\t+\ \N^{-100\gamma_{\mathrm{reg}}}$ with $\gamma_{\mathrm{reg}}$ from Definition \ref{definition:reg}.
\end{definition}
\begin{lemma}\label{lemma:method6}
We have the deterministic estimate $\|\mathbf{S}-\mathbf{Z}\|_{\t_{\mathrm{reg}};\mathbb{T}(\N)}\lesssim\N^{-20\gamma_{\mathrm{reg}}}\{\|\mathbf{Z}\|_{\t_{\mathrm{reg}};\mathbb{T}(\N)}\wedge\|\mathbf{S}\|_{\t_{\mathrm{reg}};\mathbb{T}(\N)}\}$.
\end{lemma}
Let us now record the following elementary calculation. It computes the SDE satisfied by $\mathbf{S}$ as a regularization of the SDE for $\mathbf{Z}$. (This SDE is effectively the one satisfied by $\mathbf{Z}$, but everything is regularized by the short-time heat operator.) Then, it rewrites it in Duhamel form, and decomposes the contribution of ${\mathfrak{z}}(\t,{\y(\x,\t)})\mathbf{Z}(\t,\x)$ in \eqref{eq:method2I} based on Definition \ref{definition:method3}.
\begin{lemma}\label{lemma:method7}
 Recall notation of {Definition \ref{definition:method1}} and {Proposition \ref{prop:method2}}. For any $\t\geq0$ and $\x\in\mathbb{T}(\N)$, we have 
\begin{align}
&\d\mathbf{S}(\t,\x) \nonumber\\
&= \ \mathscr{T}(\t(\N))\mathbf{S}(\t,\x)\d\t + \mathbf{H}^{\N}(\t,\t(\N),\x)\{\sqrt{2}\lambda(\t)\N^{\frac12}\mathbf{Z}(\t,\cdot)\d\mathbf{b}(\t,\cdot(\t))\}\label{eq:method7Ia} \\
&+ \mathbf{H}^{\N}(\t,\t(\N),\x)\{\partial_{\t}\log|\lambda(\t)|\times\mathbf{Z}(\t,\cdot)\log\mathbf{Z}(\t,\cdot)\}\d\t + \mathbf{H}^{\N}(\t,\t(\N),\x)({\mathfrak{z}}(\t,\cdot(\t))\mathbf{Z}(\t,\cdot))\d\t \label{eq:method7Ib}\\
&= \ \mathscr{T}(\t(\N))\mathbf{S}(\t,\x)\d\t + \mathbf{H}^{\N}(\t,\t(\N),\x)\{\sqrt{2}\lambda(\t)\N^{\frac12}\mathbf{S}(\t,\cdot)\d\mathbf{b}(\t,\cdot(\t))\}\label{eq:method7Ic}\\
&+ \ \mathbf{H}^{\N}(\t,\t(\N),\x)(\partial_{\t}\log|\lambda(\t)|\times\mathbf{S}(\t,\cdot)\log\mathbf{S}(\t,\cdot))\d\t + \mathbf{H}^{\N}(\t,\t(\N),\x)({\mathfrak{z}}(\t,\cdot(\t))\mathbf{Z}(\t,\cdot))\d\t \label{eq:method7Id}\\
&+ \ \mathbf{H}^{\N}(\t,\t(\N),\x)(\partial_{\t}\log|\lambda(\t)|\times\{\mathbf{Z}(\t,\cdot)\log\mathbf{Z}(\t,\cdot)-\mathbf{S}(\t,\cdot)\log\mathbf{S}(\t,\cdot)\})\d\t \label{eq:method7Ie}\\
&+ \ \mathbf{H}^{\N}(\t,\t(\N),\x)(\sqrt{2}\lambda(\t)\N^{\frac12}\{\mathbf{Z}(\t,\cdot)-\mathbf{S}(\t,\cdot)\}\d\mathbf{b}(\t,\cdot(\t))).\label{eq:method7If}
\end{align}
By the Duhamel principle and the semigroup property for $\mathbf{H}^{\N}$, we ultimately get
\begin{align}
&\mathbf{S}(\t,\x) \nonumber\\
&= \ \mathbf{H}^{\N}(0,\t(\N),\x)(\mathbf{Z}(0,\cdot)) + {\textstyle\int_{0}^{\t}}\mathbf{H}^{\N}(\s,\t(\N),\x)(\sqrt{2}\lambda(\s)\N^{\frac12}\mathbf{S}(\s,\cdot)\d\mathbf{b}(\s,\cdot(\s))) \label{eq:method7IIa}\\
&+ {\textstyle\int_{0}^{\t}}\mathbf{H}^{\N}(\s,\t(\N),\x)(\partial_{\s}\log|\lambda(\s)|\times\mathbf{S}(\s,\cdot)\log\mathbf{S}(\s,\cdot))\d\s \ + \ {\textstyle\int_{0}^{\t}}\mathbf{H}^{\N}(\s,\t(\N),\x)({\mathfrak{z}}(\s,\cdot(\s))\mathbf{Z}(\s,\cdot))\d\s \label{eq:method7IIb}\\
&+ \ {\textstyle\int_{0}^{\t}}\mathbf{H}^{\N}(\s,\t(\N),\x)(\partial_{\s}\log|\lambda(\s)|\times\{\mathbf{Z}(\s,\cdot)\log\mathbf{Z}(\s,\cdot)-\mathbf{S}(\s,\cdot)\log\mathbf{S}(\s,\cdot)\})\d\s \label{eq:method7IIc}\\
&+ \ {\textstyle\int_{0}^{\t}}\mathbf{H}^{\N}(\s,\t(\N),\x)(\sqrt{2}\lambda(\s)\N^{\frac12}\{\mathbf{Z}(\s,\cdot)-\mathbf{S}(\s,\cdot)\}\d\mathbf{b}(\s,\cdot(\s))). \label{eq:method7IId}
\end{align}
Finally, we have the following decomposition with notation explained afterwards:
\begin{align}
{\textstyle\int_{0}^{\t}}\mathbf{H}^{\N}(\s,\t(\N),\x)({\mathfrak{z}}(\s,\cdot(\s))\mathbf{Z}(\s,\cdot))\d\s \ = \ \mathrm{QCT}(\t,\x)+\mathrm{LCT}(\t,\x)+\mathrm{CT}(\t,\x)+\mathrm{An}(\t,\x). \label{eq:method7III}
\end{align}
Above, we have introduced the following terms in which $\mathbf{U}(\s)=\mathbf{U}^{\s,\cdot(\s)}$ and $\mathbf{V}(\s)=\mathbf{U}^{\s,\cdot(\s)+1}$ for convenience:
{\small
\begin{align}
&\mathrm{QCT}(\t,\x)\nonumber\\
&:= {\textstyle\int_{0}^{\t}}\mathbf{H}^{\N}(\s,\t(\N),\x)(\lambda(\s)\N\{\mathscr{U}'(\s,\mathbf{U}(\s))-\bar{\alpha}(\s)\mathbf{U}(\s)-\tfrac12\lambda(\s)[\mathscr{U}'(\s,\mathbf{U}(\s))\mathbf{U}(\s)-1]\}\mathbf{Z}(\s,\cdot))\d\s \label{eq:method7IVa}\\
&+ {\textstyle\int_{0}^{\t}}\mathbf{H}^{\N}(\s,\t(\N),\x)(\lambda(\s)\N\{\mathscr{U}'(\s,\mathbf{V}(\s))-\bar{\alpha}(\s)\mathbf{V}(\s)-\tfrac12\lambda(\s)[\mathscr{U}'(\s,\mathbf{V}(\s))\mathbf{V}(\s)-1]\}\mathbf{Z}(\s,\cdot))\d\s \label{eq:method7IVb}\\
&\mathrm{LCT}(\t,\x) \nonumber\\
&:= \ \tfrac12{\textstyle\int_{0}^{\t}}\mathbf{H}^{\N}(\s,\t(\N),\x)(\lambda(\s)\N^{\frac32}\grad^{+}\{\mathscr{W}'(\s,\mathbf{U}(\s))\mathbf{Z}(\s,\cdot)\}-\lambda(\s)\N^{\frac32}\grad^{-}\{\mathscr{W}'(\s,\mathbf{V}(\s))\mathbf{Z}(\s,\cdot)\})\d\s \label{eq:method7IVc}\\
&\mathrm{CT}(\t,\x) \nonumber\\
&:= \ \tfrac14{\textstyle\int_{0}^{\t}}\mathbf{H}^{\N}(\s,\t(\N),\x)(\lambda(\s)^{2}\N\grad^{+}\{[\bar{\alpha}(\s)\mathbf{U}(\s)^{2}-1]\mathbf{Z}(\s,\cdot)\}+\lambda(\s)^{2}\N\grad^{-}\{[\bar{\alpha}(\s)\mathbf{V}(\s)^{2}-1]\mathbf{Z}(\s,\cdot)\})\d\s \label{eq:method7IVd}\\
&+ \ {\textstyle\int_{0}^{\t}}\mathbf{H}(\s,\t(\N),\x)(\{\tfrac{1}{12}\lambda(\s)^{4}\E^{0,\s}[\mathscr{U}'(\s,\mathbf{u})\mathbf{u}^{3}]+\tfrac16\lambda(\s)^{3}\bar{\alpha}(\s)\mathbf{V}(\s)^{3}-\lambda(\s)\mathscr{R}(\s)\}\mathbf{Z}(\s,\cdot))\d\s \label{eq:method7IVe}\\
&+ \ {\textstyle\int_{0}^{\t}}\mathbf{H}^{\N}(\s,\t(\N),\x)(\{\tfrac{1}{12}\lambda(\s)^{4}\mathscr{U}'(\s,\mathbf{V}(\s))\mathbf{V}(\s)^{3}-\E^{0,\s}[\mathscr{U}'(\s,\mathbf{u})\mathbf{u}^{3})]\}\mathbf{Z}(\s,\cdot))\d\s \label{eq:method7IVf}\\
&\mathrm{An}(\t,\x) \nonumber\\
&:= \ {\textstyle\int_{0}^{\t}}\mathbf{H}^{\N}(\s,\t(\N),\x)\{\N^{-\frac12}{\mathfrak{b}}(\s,\cdot(\s))\mathbf{Z}(\s,\cdot)+\N^{-\frac12}{\mathfrak{c}}(\s,\cdot(\s))\mathbf{Z}(\s,\cdot)\}\d\s \label{eq:method7IVh}\\
&+ \ {\textstyle\int_{0}^{\t}}\mathbf{H}^{\N}(\s,\t(\N),\x)(\N^{\frac12}\grad^{-}\{{\mathfrak{a}}^{-}(\s,\cdot(\s))\mathbf{Z}(\s,\cdot)\}+\N^{\frac12}\grad^{+}\{{\mathfrak{a}}^{+}(\s,\cdot(\s))\mathbf{Z}(\s,\cdot)\})\d\s.\label{eq:method7IVi}
\end{align}
}
\end{lemma}
\subsubsection{The a priori estimates}
Let us now introduce the a priori stochastic estimates needed to prove Theorem \ref{theorem:kpz}. We package them as stopping times. We then use them to basically stop $\mathbf{S}$ and get something we can analyze via standard SPDE tools.
\begin{definition}\label{definition:method8}
 Recall $\gamma_{\mathrm{KL}}$ from Definition \ref{definition:entropydata}, and recall $\gamma_{\mathrm{reg}}$ and $\t_{\mathrm{reg}}$ from Definition \ref{definition:reg}. We define $\gamma_{\mathrm{ap}}={c}\gamma_{\mathrm{reg}}$ and $\beta_{\mathrm{BG}}={c'}\gamma_{\mathrm{KL}}$ {for some fixed $c,c'>0$ chosen such that $\beta_{\mathrm{BG}}$ at least some large but fixed factor times $\gamma_{\mathrm{reg}},\gamma_{\mathrm{ap}}$}. We now define the following stopping times with explanation given after:
{\small
\begin{align}
\t_{\mathrm{ap}}&:=\inf\left\{\t\geq0: \ \log\N \ \leq \ \|\mathbf{Z}\|_{\t;\mathbb{T}(\N)}\vee\|\mathbf{S}\|_{\t;\mathbb{T}(\N)} + \|\mathbf{Z}^{-1}\|_{\t;\mathbb{T}(\N)}\vee\|\mathbf{S}^{-1}\|_{\t;\mathbb{T}(\N)} \ \leq \ \N^{\gamma_{\mathrm{ap}}} \right\} \wedge \t_{\mathrm{reg}} \\
\t_{\mathrm{BG}}&:=\inf\left\{\t\geq0: \ \|\mathrm{QCT}\|_{\t;\mathbb{T}(\N)}+\|\mathrm{LCT}\|_{\t;\mathbb{T}(\N)}+\|\mathrm{CT}\|_{\t;\mathbb{T}(\N)}\geq \N^{-\beta_{\mathrm{BG}}}\|\mathbf{Z}\|_{\t;\mathbb{T}(\N)}\right\}\wedge\t_{\mathrm{reg}} \\
\t_{\mathrm{rest}}&:=\inf\{\t\geq0: \  \|\mathrm{An}\|_{\t;\mathbb{T}(\N)} \geq\N^{-\frac{\gamma_{\mathrm{reg}}}{2}}\|\mathbf{Z}\|_{\t;\mathbb{T}(\N)} \\
&\quad\quad\quad\quad\quad\quad\mathrm{or}\quad \|\eqref{eq:method7IIc}\|_{\t;\mathbb{T}(\N)} + \|\eqref{eq:method7IId}\|_{\t;\mathbb{T}(\N)}\geq\N^{-\frac{\gamma_{\mathrm{reg}}}{2}+\gamma_{\mathrm{ap}}}\}\wedge\t_{\mathrm{reg}}.\nonumber
\end{align}
}(In words, $\t_{\mathrm{ap}}$ gives a priori estimates that should be redundant if $\mathbf{Z}$ has a space-time continuous limit. Next, $\t_{\mathrm{BG}}$ controls the contribution of the first three error terms on the RHS of \eqref{eq:method7III}; the subscript ``BG" refers to ``Boltzmann-Gibbs". Lastly, $\t_{\mathrm{rest}}$ controls the remaining error term in \eqref{eq:method7III}.) Set $\t_{\mathrm{st}}:=\t_{\mathrm{ap}}\wedge\t_{\mathrm{BG}}\wedge\t_{\mathrm{rest}}$. We now define $\mathbf{Y}(\cdot,\cdot)$ to solve the following stochastic integral equation obtained by stopping $\mathbf{S}$:
\begin{align}
\mathbf{Y}(\t,\x) \ &= \ \mathbf{H}^{\N}(0,\t(\N),\x)(\mathbf{Z}(0,\cdot)) + {\textstyle\int_{0}^{\t}}\mathbf{H}^{\N}(\s,\t(\N),\x)(\sqrt{2}\lambda(\s)\N^{\frac12}\mathbf{Y}(\s,\cdot)\d\mathbf{b}(\s,\cdot(\s))) \label{eq:method8Ia}\\
&+ \ {\textstyle\int_{0}^{\t}}\mathbf{H}^{\N}(\s,\t(\N),\x)(\partial_{\s}\log|\lambda(\s)|\times\mathbf{Y}(\s,\cdot)\log\mathbf{Y}(\s,\cdot))\d\s \label{eq:method8Ib}\\
&+ \ \mathbf{1}(\t\leq\t_{\mathrm{st}}){\textstyle\int_{0}^{\t}}\mathbf{H}^{\N}(\s,\t(\N),\x)({\mathfrak{z}}(\s,\cdot(\s))\mathbf{Z}(\s,\cdot))\d\s \label{eq:method8Ic}\\
&+ \ \mathbf{1}(\t\leq\t_{\mathrm{st}}){\textstyle\int_{0}^{\t}}\mathbf{H}^{\N}(\s,\t(\N),\x)(\partial_{\s}\log|\lambda(\s)|\times\{\mathbf{Z}(\s,\cdot)\log\mathbf{Z}(\s,\cdot)-\mathbf{S}(\s,\cdot)\log\mathbf{S}(\s,\cdot)\})\d\s \label{eq:method8Id}\\
&+ \ \mathbf{1}(\t\leq\t_{\mathrm{st}}){\textstyle\int_{0}^{\t}}\mathbf{H}^{\N}(\s,\t(\N),\x)(\sqrt{2}\lambda(\s)\N^{\frac12}\{\mathbf{Z}(\s,\cdot)-\mathbf{S}(\s,\cdot)\}\d\mathbf{b}(\s,\cdot(\s))).\label{eq:method8Ie}
\end{align}
Now, let $\mathbf{W}(\cdot,\cdot)$ solve the same stochastic equation but formally forgetting all error terms \eqref{eq:method8Ic}-\eqref{eq:method8Ie}:
\begin{align}
\mathbf{W}(\t,\x) \ &= \ \mathbf{H}^{\N}(0,\t(\N),\x)(\mathbf{Z}(0,\cdot)) + {\textstyle\int_{0}^{\t}}\mathbf{H}^{\N}(\s,\t(\N),\x)(\sqrt{2}\lambda(\s)\N^{\frac12}\mathbf{W}(\s,\cdot)\d\mathbf{b}(\s,\cdot(\s))) \label{eq:method8IIa}\\
&+ \ {\textstyle\int_{0}^{\t}}\mathbf{H}^{\N}(\s,\t(\N),\x)(\partial_{\s}\log|\lambda(\s)|\times\mathbf{W}(\s,\cdot)\log\mathbf{W}(\s,\cdot))\d\s. \label{eq:method8IIb}
\end{align}
Lastly, for any $\t\geq0$, we extend $\mathbf{Y}(\t,\N\x)$ and $\mathbf{W}(\t,\N\x)$ from $\x\in\N^{-1}\mathbb{T}(\N)$ to $\x\in\mathbb{T}$ by linear interpolation.
\end{definition}
\begin{rem}\label{remark:method9}
 Because $\mathbf{Z}$ and $\mathbf{S}$ (and their inverses) are continuous in time with probability 1, for sufficiently large $\N$ (depending on $\gamma_{\mathrm{ap}}$), as soon as one of $\mathbf{Z}$ or $\mathbf{S}$ or their inverses exceeds $\log\N$, it is also $\leq\N^{\gamma_{\mathrm{ap}}}$. So, the constraint $\leq\N^{\gamma_{\mathrm{ap}}}$ in the definition of $\t_{\mathrm{ap}}$ is redundant. We have included it, however, because we will basically only use $\t_{\mathrm{ap}}$ to make sure that $\mathbf{Z}$, $\mathbf{S}$, their inverses, and similar exponentials are $\lesssim\N^{\gamma_{\mathrm{ap}}}$ before time $\t_{\mathrm{ap}}$. (It is only for technical, uninteresting reasons that we need the upper bound of $\log\N$ instead of $\N^{\gamma_{\mathrm{ap}}}$.) Next, we note that it is easy to see that $\mathbf{Y}$ and $\mathbf{W}$ are adapted to the filtration generated by the Brownian motions, so products with Brownian motions are well-defined. Indeed, for $\mathbf{Y}$, we have the following SDE (with a jump at $\t_{\mathrm{st}}$):
\begin{align}
\d\mathbf{Y}(\t,\x) \ &= \ \mathscr{T}(\t(\N))\mathbf{Y}(\t,\x)\d\t + \mathbf{H}^{\N}(\t,\t(\N),\x)\{\sqrt{2}\lambda(\t)\N^{\frac12}\mathbf{Y}(\t,\cdot)\d\mathbf{b}(\t,\cdot(\t))\} \\
&+ \ \mathbf{H}^{\N}(\t,\t(\N),\x)\{\partial_{\t}\log|\lambda(\t)|\times\mathbf{Y}(\t,\cdot)\log\mathbf{Y}(\t,\cdot)\}\d\t \\
&+ \ \mathbf{1}(\t\leq\t_{\mathrm{st}})\mathbf{H}^{\N}(\t,\t(\N),\x)\left({\mathfrak{z}}(\t,\cdot(\t))\mathbf{Z}(\t,\cdot)\right)\d\t \\
&+ \ \mathbf{1}(\t\leq\t_{\mathrm{st}})\mathbf{H}^{\N}(\t,\t(\N),\x)(\partial_{\t}\log|\lambda(\t)|\times\{\mathbf{Z}(\t,\cdot)\log\mathbf{Z}(\t,\cdot)-\mathbf{S}(\t,\cdot)\log\mathbf{S}(\t,\cdot)\})\d\t \\
&+ \ \mathbf{1}(\t\leq\t_{\mathrm{st}})\mathbf{H}^{\N}(\t,\t(\N),\x)(\sqrt{2}\lambda(\t)\N^{\frac12}\{\mathbf{Z}(\t,\cdot)-\mathbf{S}(\t,\cdot)\}\d\mathbf{b}(\t,\cdot(\t))).
\end{align}
Similarly for $\mathbf{W}$, we have the following discretized TISHE with additional smoothing via $\mathbf{H}^{\N}(\t,\t(\N),\x)$:
\begin{align}
\d\mathbf{W}(\t,\x) \ &= \ \mathscr{T}(\t(\N))\mathbf{W}(\t,\x)\d\t + \mathbf{H}^{\N}(\t,\t(\N),\x)\{\sqrt{2}\lambda(\t)\N^{\frac12}\mathbf{W}(\t,\cdot)\d\mathbf{b}(\t,\cdot(\t))\} \\
&+ \ \mathbf{H}^{\N}(\t,\t(\N),\x)\{\partial_{\t}\log|\lambda(\t)|\times\mathbf{W}(\t,\cdot)\log\mathbf{W}(\t,\cdot)\}\d\t. 
\end{align}
These SDEs and \eqref{eq:method8Ia}-\eqref{eq:method8Ie} and \eqref{eq:method8IIa}-\eqref{eq:method8IIb} are really defined only until their respective explosion times. We will eventually show these explosion times to be independent of $\N$ so this point is ultimately unimportant. But technically, it must be made.
\end{rem}
For $\t\leq\t_{\mathrm{st}}$, we can drop indicators $\mathbf{1}(\t\leq\t_{\mathrm{st}})$ in \eqref{eq:method8Ia}-\eqref{eq:method8Ie}, giving back the {$\mathbf{S}$} SDE. By Ito theory, this ultimately yields:
\begin{lemma}\label{lemma:method10}
 With probability 1, we have $\mathbf{Y}(\t,\x)=\mathbf{S}(\t,\x)$ for all $0\leq\t\leq\t_{\mathrm{st}}$ and $\x\in\mathbb{T}(\N)$.
\end{lemma}
\subsubsection{Comparing $\mathbf{Y}$ and $\mathbf{W}$}
We now show why the a priori estimates defining $\mathbf{Y}$ are useful (they let us directly compare $\mathbf{Y}$ to $\mathbf{W}$). This is made precise in the following, which also states the desired convergence in Theorem \ref{theorem:kpz} for $\mathbf{W}$ instead of $\mathbf{Z}$.
\begin{prop}\label{prop:method11}
 With high probability, we know that $\|\mathbf{Y}-\mathbf{W}\|_{1;\mathbb{T}(\N)}\lesssim\N^{-\frac{\gamma_{\mathrm{reg}}}{100}}\{1+\|\mathbf{W}\|_{1;\mathbb{T}(\N)}+\|\mathbf{W}^{-1}\|_{1;\mathbb{T}(\N)}\}$. Moreover, there exists a coupling between $\mathbf{W}$ and the solution $\mathbf{Z}^{\infty}$ to \eqref{eq:she2a}-\eqref{eq:she2b} with initial data $\mathbf{Z}^{\infty,\mathrm{in}}(0,\cdot)$ from the statement of {Theorem \ref{theorem:kpz}} such that $|\mathbf{W}(\t,\N\x)-\mathbf{Z}^{\infty}(\t,\x)|\to0$ uniformly over $0\leq\t\leq1$ and $\x\in\mathbb{T}(\N)$ in probability as $\N\to\infty$.
\end{prop}
According to Theorem \ref{theorem:she}, $\mathrm{TISHE}$ is a well-posed $(1+1)$-dimensional stochastic heat PDE with locally smooth nonlinearity. (By well-posed, we also mean strictly positive with probability 1, so that its inverse is also well-defined and continuous.) Thus, convergence of $\mathbf{W}$ to $\mathrm{TISHE}$ basically follows by stability of well-posed (stochastic) PDEs under discretizations. On the other hand, the bound on $\mathbf{Y}-\mathbf{W}$  follows by similar ideas, except one now uses the a priori estimates from $\t_{\mathrm{st}}$; see Definition \ref{definition:method8}. (The difference is controlled by $\mathbf{W}$ and $\mathbf{W}^{-1}$ because the coefficients in the $\mathbf{W}$-equation are smooth if $\mathbf{W},\mathbf{W}^{-1}$ are away from $\infty$.) We give a proof in the appendix modulo easy and elementary technicalities. Lastly, as we explain in the proof of Theorem \ref{theorem:kpz} below, for our purposes, it would be enough to just state $\mathbf{Y}\to\mathrm{TISHE}$ in Proposition \ref{prop:method11}. We stated Proposition \ref{prop:method11} with more detail to {express} the ``rough paths" nature of our argument. Also, we want to highlight the benefit of a priori estimates in $\t_{\mathrm{st}}$.
\subsubsection{The a priori estimates in $\t_{\mathrm{st}}$ hold with high probability}
Lemma \ref{lemma:method10} identifies $\mathbf{S}$ and $\mathbf{Y}$ until $\t_{\mathrm{st}}$. So, we (ultimately) get $\mathbf{S}\approx\mathrm{TISHE}$ until $\t_{\mathrm{st}}$. To obtain it until time 1, or more precisely until the regularity blow-up time $\t_{\mathrm{reg}}$ from Definition \ref{definition:reg}, we need to prove that $\t_{\mathrm{st}}=\t_{\mathrm{reg}}$ with high probability. (We could settle for deriving KPZ until $\t_{\mathrm{st}}$ instead of $\t_{\mathrm{reg}}$. But this result would follow by Proposition \ref{prop:method11}. In particular, all our stochastic homogenization analysis is in proving the following.)
\begin{prop}\label{prop:method12}
 With high probability, we have $\t_{\mathrm{st}}=\t_{\mathrm{reg}}$.
\end{prop}
\begin{proof}[Proof of {Theorem \ref{theorem:kpz}}]
Fix any $\e>0$ independent of $\N$. We first claim the following (to be justified afterwards):
\begin{align}
&\mathbb{P}[\|\mathbf{S}(\cdot,\N\cdot)-\mathbf{Z}^{\infty}(\cdot,\cdot)>\e\|_{1;\mathbb{T}}] \ \leq \ \mathbb{P}[\|\mathbf{Y}(\cdot,\N\cdot)-\mathbf{Z}^{\infty}(\cdot,\cdot)\|_{\t_{\mathrm{st}};\mathbb{T}}>\e]+\mathbb{P}[\t_{\mathrm{st}}\neq1] \label{eq:kpz1a}\\
&\leq \ \mathbb{P}[\|\mathbf{Y}-\mathbf{W}\|_{\t_{\mathrm{st}};\mathbb{T}(\N)}>\tfrac{\e}{2}]+\mathbb{P}[\|\mathbf{W}(\cdot,\N\cdot)-\mathbf{Z}^{\infty}(\cdot,\cdot)\|_{\t_{\mathrm{st}};\mathbb{T}}>\tfrac{\e}{2}]+\mathbb{P}[\t_{\mathrm{st}}\neq1] \label{eq:kpz1b}\\
&= \ \mathrm{o}(1)+\mathbb{P}[\t_{\mathrm{st}}\neq1] \ \leq \ \mathrm{o}(1)+\mathbb{P}[\t_{\mathrm{reg}}\neq1] \ \leq \ \mathrm{o}(1). \label{eq:kpz1c}
\end{align}
\eqref{eq:kpz1a} follows from first replacing $1;\mathbb{T}$ on the LHS by $\t_{\mathrm{st}};\mathbb{T}$. By the union bound, the cost is the last term in \eqref{eq:kpz1a}. We then use Lemma \ref{lemma:method10} to replace $\mathbf{S}(\cdot,\N\cdot)$ by $\mathbf{Y}(\cdot,\N\cdot)$ before time $\t_{\mathrm{st}}$. \eqref{eq:kpz1b} follows by triangle inequality. The identity in \eqref{eq:kpz1c} follows by Proposition \ref{prop:method11}. The inequalities in \eqref{eq:kpz1c} follow by Proposition \ref{prop:method12} and by $\t_{\mathrm{reg}}=1$ with high probability. Now, we claim
\begin{align}
&\mathbb{P}[\|\mathbf{Z}(\cdot,\N\cdot)-\mathbf{S}(\cdot,\N\cdot)\|_{1;\mathbb{T}}>\e] \nonumber\\
&\leq \ \mathbb{P}[\|\mathbf{Z}-\mathbf{S}\|_{1;\mathbb{T}(\N)}>\e] \ \leq \ \mathbb{P}[\|\mathbf{Z}-\mathbf{S}\|_{\t_{\mathrm{reg}};\mathbb{T}(\N)}>\e]+\mathbb{P}[\t_{\mathrm{reg}}\neq1] \label{eq:kpz2a}\\
&\leq \ {\mathbb{P}[\|\mathbf{Z}-\mathbf{S}\|_{\t_{\mathrm{reg}};\mathbb{T}(\N)}\gtrsim \N^{-\gamma_{\mathrm{ap}}}\e\|\mathbf{Z}\|_{\t_{\mathrm{reg}};\mathbb{T}(\N)}]}+\mathbb{P}[\t_{\mathrm{reg}}\neq1] \ = \ \mathrm{o}(1). \label{eq:kpz2b}
\end{align}
The first bound in \eqref{eq:kpz2a} follows since $\|\|_{1;\N\mathbb{T}}\geq\|\|_{1;\mathbb{T}(\N)}$. (Indeed, we know $\mathbb{T}(\N)\subseteq\N\mathbb{T}$.) The second bound in \eqref{eq:kpz2a} follows by union bound. The first bound in \eqref{eq:kpz2b} follows because of {$\|\mathbf{Z}\|_{\t_{\mathrm{reg}};\mathbb{T}(\N)}\lesssim \N^{\gamma_{\mathrm{ap}}}$}. The last estimate in \eqref{eq:kpz2b} follows by Lemma \ref{lemma:method6}. Note $\e>0$ was independent of $\N$ but otherwise arbitrary in \eqref{eq:kpz1a}-\eqref{eq:kpz2b}. By triangle inequality, union bound, and \eqref{eq:kpz1a}-\eqref{eq:kpz1c} and \eqref{eq:kpz2a}-\eqref{eq:kpz2b}, we get $\mathbf{Z}(\cdot,\N\cdot)-\mathbf{Z}^{\infty}(\cdot,\cdot)\to0$ in probability as claimed. This finishes the proof.
\end{proof}
\subsubsection{What is left to do}
We proved Theorem \ref{theorem:kpz} assuming Proposition \ref{prop:method2}, Lemmas \ref{lemma:method4}, \ref{lemma:method6}, \ref{lemma:method7}, and \ref{lemma:method10}, and Proposition \ref{prop:method12}. It remains to prove these. (These types of remarks will be made for organizational clarity.)
\subsection{Proof outline of Proposition \ref{prop:method12}}
The purpose of this brief subsection is to break Proposition \ref{prop:method12} into two sets of bounds. The first is an analytic lemma, whose proof is standard (S)PDE procedure. The second is the set of Boltzmann-Gibbs principles that we discussed after Lemma \ref{lemma:method4}. {Our stochastic homogenization is only used to prove the latter result.}
\begin{lemma}\label{lemma:method13}
 We have the following. The first is deterministic and the second is with high probability:
\begin{align}
\|\mathrm{An}\|_{\t_{\mathrm{st}};\mathbb{T}(\N)} \ \lesssim \ \ \N^{-\frac13}\|\mathbf{Z}\|_{\t_{\mathrm{st}};\mathbb{T}(\N)} \quad\mathrm{and}\quad \|\eqref{eq:method7IIc}\|_{\t_{\mathrm{st}};\mathbb{T}(\N)} + \|\eqref{eq:method7IId}\|_{\t_{\mathrm{st}};\mathbb{T}(\N)}\ \lesssim \ \ \N^{-\frac34\gamma_{\mathrm{reg}}+\gamma_{\mathrm{ap}}}. \label{eq:method13I}
\end{align}
\end{lemma}
\begin{theorem}\label{theorem:bg}
 Recall the constant $\beta_{\mathrm{BG}}$ from {Definition \ref{definition:method8}}. With high probability, we have
\begin{align}
\|\mathrm{QCT}\|_{\t_{\mathrm{st}};\mathbb{T}(\N)} \ &\lesssim \ \ \N^{-2\beta_{\mathrm{BG}}}\|\mathbf{Z}\|_{\t_{\mathrm{st}};\mathbb{T}(\N)} \label{eq:bg2} \\
\|\mathrm{LCT}\|_{\t_{\mathrm{st}};\mathbb{T}(\N)} \ &\lesssim \ \ \N^{-2\beta_{\mathrm{BG}}}\|\mathbf{Z}\|_{\t_{\mathrm{st}};\mathbb{T}(\N)} \label{eq:bg1} \\
\|\mathrm{CT}\|_{\t_{\mathrm{st}};\mathbb{T}(\N)} \ &\lesssim \ \ \N^{-2\beta_{\mathrm{BG}}}\|\mathbf{Z}\|_{\t_{\mathrm{st}};\mathbb{T}(\N)} \label{eq:hl}.
\end{align}
\end{theorem}
\begin{rem}\label{remark:method15}
 The proof of Theorem \ref{theorem:bg} is the ``stochastic homogenization" heart of this paper. (The proofs of Proposition \ref{prop:method2} and Lemma \ref{lemma:method4} are the ``algebraic" heart.) For now, let us just emphasize that \eqref{eq:bg2}-\eqref{eq:hl} are better than what $\t_{\mathrm{BG}}=\t_{\mathrm{reg}}$ asks for. In other words, \eqref{eq:bg2}-\eqref{eq:hl} is ``self-propagating". This is crucial. (Indeed, if we change $2\beta_{\mathrm{BG}}$ to $\beta_{\mathrm{BG}}$ in \eqref{eq:bg2}-\eqref{eq:hl}, the resulting estimates are trivial by definition of $\t_{\mathrm{st}}$. In this case, Theorem \ref{theorem:bg} would be pointless, and there would be no homogenization analysis necessary.)
\end{rem}
Let us now deduce Proposition \ref{prop:method12} from Lemma \ref{lemma:method13} and Theorem \ref{theorem:bg}. In a nutshell, these estimates imply that $\t_{\mathrm{BG}}$ and $\t_{\mathrm{rest}}$ (see Definition \ref{definition:method8}) equal $\t_{\mathrm{reg}}$ with high probability. We are left to show the same for $\t_{\mathrm{ap}}$. For this, use that $\t_{\mathrm{BG}},\t_{\mathrm{rest}}=\t_{\mathrm{reg}}$ with high probability to show that $\mathbf{Z}$ is essentially a perturbation of $\mathbf{W}$. But $\mathbf{W}$ is a discretization of a well-posed SPDE, so it is both well-behaved and stable under perturbations. In particular, $\mathbf{Z}$ is well-behaved, and thus $\t_{\mathrm{ap}}=\t_{\mathrm{reg}}$ with high probability (recall from Definition \ref{definition:method8} that $\t_{\mathrm{ap}}$ is a stopping time used to control $\mathbf{Z}$ from above and below).
\begin{proof}[Proof of {Proposition \ref{prop:method12}}]
Recall {that} $\t_{\mathrm{st}}=\t_{\mathrm{ap}}\wedge\t_{\mathrm{BG}}\wedge\t_{\mathrm{rest}}\wedge\t_{\mathrm{reg}}$; see Definition \ref{definition:method8}. By $\t_{\mathrm{st}}\leq\t_{\mathrm{reg}}$ and union bound, if $\t_{\mathrm{st}}\neq\t_{\mathrm{reg}}$ (i.e. $\t_{\mathrm{st}}<\t_{\mathrm{reg}}$), then one of the stopping times defining $\t_{\mathrm{st}}$ must both equal $\t_{\mathrm{st}}$ and be strictly less than $\t_{\mathrm{reg}}$, so
\begin{align}
\mathbb{P}[\t_{\mathrm{st}}\neq\t_{\mathrm{reg}}] \ &= \ \mathbb{P}[\t_{\mathrm{st}}<\t_{\mathrm{reg}}] \nonumber\\
&\leq \ \mathbb{P}[\t_{\mathrm{st}}=\t_{\mathrm{ap}}<\t_{\mathrm{reg}}]+\mathbb{P}[\t_{\mathrm{st}}=\t_{\mathrm{rest}}<\t_{\mathrm{reg}}]+\mathbb{P}[\t_{\mathrm{st}}=\t_{\mathrm{BG}}<\t_{\mathrm{reg}}]. \label{eq:method121}
\end{align}
Take the last term in \eqref{eq:method121}. Because {$\t_{\mathrm{BG}}\leq\t_{\mathrm{reg}}$}, on the event $\t_{\mathrm{BG}}<\t_{\mathrm{reg}}$, the lower bound defining $\t_{\mathrm{BG}}$ (see Definition \ref{definition:method8}) must be realized at time $\t_{\mathrm{st}}=\t_{\mathrm{BG}}$. (This also requires continuity of $\mathrm{QCT},\mathrm{LCT},\mathrm{CT}$; this holds with probability 1 because these are Riemann integrals.) But realizing the lower bound defining $\t_{\mathrm{BG}}$ at time $\t_{\mathrm{st}}$ means \eqref{eq:bg2}-\eqref{eq:hl} fails. To summarize, the last term in \eqref{eq:method121} is bounded by the probability that one of \eqref{eq:bg2}-\eqref{eq:hl} fails. Thus, the last term in \eqref{eq:method121} is $\mathrm{o}(1)$ by Theorem \ref{theorem:bg}. The same argument, but replacing $\t_{\mathrm{BG}}$ by $\t_{\mathrm{rest}}$ and \eqref{eq:bg2}-\eqref{eq:hl} by \eqref{eq:method13I}, shows the second-to-last term in \eqref{eq:method121} is $\mathrm{o}(1)$. We are left to control the first term on the far RHS of \eqref{eq:method121}. To this end, we claim the following (that we justify afterwards):
\begin{align}
\mathbb{P}[\t_{\mathrm{st}}=\t_{\mathrm{ap}}<\t_{\mathrm{reg}}] \ &\leq \ \mathbb{P}[\|\mathbf{Z}\|_{\t_{\mathrm{st}};\mathbb{T}(\N)}\vee\|\mathbf{S}\|_{\t_{\mathrm{st}};\mathbb{T}(\N)}\vee\|\mathbf{Z}^{-1}\|_{\t_{\mathrm{st}};\mathbb{T}(\N)}\vee\|\mathbf{S}^{-1}\|_{\t_{\mathrm{st}};\mathbb{T}(\N)}\geq\log\N] \label{eq:method122a}\\
&\leq \ \mathbb{P}[\|\mathbf{S}\|_{\t_{\mathrm{st}};\mathbb{T}(\N)}\vee\|\mathbf{S}^{-1}\|_{\t_{\mathrm{st}};\mathbb{T}(\N)}\gtrsim\log\N] \label{eq:method122b}\\
&= \ \mathbb{P}[\|\mathbf{Y}\|_{\t_{\mathrm{st}};\mathbb{T}(\N)}\vee\|\mathbf{Y}^{-1}\|_{\t_{\mathrm{st}};\mathbb{T}(\N)}\gtrsim\log\N] \label{eq:method122c}\\
&\leq \ \mathbb{P}[\|\mathbf{W}\|_{\t_{\mathrm{st}};\mathbb{T}(\N)}\vee\|\mathbf{W}^{-1}\|_{\t_{\mathrm{st}};\mathbb{T}(\N)}\gtrsim\log\N]+\mathrm{o}(1). \label{eq:method122d}
\end{align}
To get \eqref{eq:method122a}, we note $\t_{\mathrm{ap}}=\t_{\mathrm{ap}}\wedge\t_{\mathrm{reg}}$, so $\t_{\mathrm{ap}}<\t_{\mathrm{reg}}$ means the lower bound defining $\t_{\mathrm{ap}}$ is realized at time $\t_{\mathrm{ap}}=\t_{\mathrm{st}}$. \eqref{eq:method122b} follows from Lemma \ref{lemma:method6}. Indeed, Lemma \ref{lemma:method6} implies $\|\mathbf{Z}-\mathbf{S}\|=\mathrm{o}(1)\|\mathbf{S}\|$. In particular, $\|\mathbf{Z}\|\lesssim\|\mathbf{S}\|$. This lets us remove $\|\mathbf{Z}\|$ inside the probability on the RHS of \eqref{eq:method122a}. It also implies $\|\mathbf{Z}-\mathbf{S}\|=\mathrm{o}(1)\|\mathbf{Z}\|$, which, by ordinary calculus, gives the estimates $\|\mathbf{Z}^{-1}-\mathbf{S}^{-1}\|\lesssim\|\mathbf{Z}^{-1}\|\|\mathbf{S}^{-1}\|\|\mathbf{Z}-\mathbf{S}\|\lesssim\mathrm{o}(1)\|\mathbf{S}^{-1}\|$. From this last pair of bounds, we have $\|\mathbf{Z}^{-1}\|\lesssim\|\mathbf{S}^{-1}\|$, so we can also drop $\|\mathbf{Z}^{-1}\|$ on the RHS of \eqref{eq:method122a}. This justifies \eqref{eq:method122b}. \eqref{eq:method122c} follows by Lemma \ref{lemma:method10}. \eqref{eq:method122d} follows from the same reasoning as \eqref{eq:method122b}. (Except, we swap $(\mathbf{Z},\mathbf{S})$ for $(\mathbf{Y},\mathbf{W})$ and use Proposition \ref{prop:method11} instead of Lemma \ref{lemma:method6}. Because the $\mathbf{Y}-\mathbf{W}$ estimate in Proposition \ref{prop:method11} {holds with high probability}, we get the extra $\mathrm{o}(1)$ in \eqref{eq:method122d}.) Now, because $\N^{-1}\mathbb{T}(\N)\subseteq\mathbb{T}$ and $\t_{\mathrm{st}}\leq1$,
\begin{align}
\mathbb{P}[\|\mathbf{W}\|_{\t_{\mathrm{st}};\mathbb{T}(\N)}\vee\|\mathbf{W}^{-1}\|_{\t_{\mathrm{st}};\mathbb{T}(\N)}\gtrsim\log\N] \ \leq \ \mathbb{P}[\|\mathbf{W}(\cdot,\N\cdot)\|_{1;\mathbb{T}}\vee\|\mathbf{W}(\cdot,\N\cdot)^{-1}\|_{1;\mathbb{T}}\gtrsim\log\N]. \label{eq:method122e}
\end{align}
By Proposition \ref{prop:method11}, we know $\mathbf{W}(\cdot,\N\cdot)$ converges uniformly (in probability) on $[0,1]\times\mathbb{T}$ to the solution of TISHE with strictly positive and continuous initial data. Theorem \ref{theorem:she} says this solution is continuous and positive on $[0,1]\times\mathbb{T}$ with probability 1. Thus, $\mathbf{W}(\cdot,\N\cdot)^{-1}$ converges uniformly (in probability) on $[0,1]\times\mathbb{T}$ to the inverse of said TISHE solution; this TISHE-inverse is continuous on $[0,1]\times\mathbb{T}$ with probability 1. Ultimately, $\mathrm{RHS}\eqref{eq:method122e}=\mathrm{o}(1)$. This finishes the proof.
\end{proof}
\subsection{What is left}
Proposition \ref{prop:method2}, Lemmas \ref{lemma:method4}, \ref{lemma:method6}, \ref{lemma:method7}, \ref{lemma:method10}, and \ref{lemma:method13}, and Theorem \ref{theorem:bg}. Having said this, let us explain the organization for the rest of the paper. Proposition \ref{prop:method2} and Lemma \ref{lemma:method4} are shown by long calculations. They are also perhaps not so interesting. In fact, they must be true for universality to be true, assumptions aside. Lemmas \ref{lemma:method6}, \ref{lemma:method7}, \ref{lemma:method10}, and \ref{lemma:method13} are more or less standard. Thus, we give proofs of all these results in the appendix. 

On the other hand, Theorem \ref{theorem:bg} is interesting; it is a stochastic estimate whose proof is the technical heart of this paper, so we focus on it first. It has a number of separate steps, some of which are complicated in their own right, so we spend the rest of the non-appendix sections on this (and Corollary \ref{corollary:kpz}). 

In the appendix, besides the proofs of Lemmas \ref{lemma:method6}, \ref{lemma:method7}, \ref{lemma:method10}, and \ref{lemma:method13}, we present auxiliary results used throughout. These include estimates for the heat kernel in Definition \ref{definition:method1}.
%
%
%
\section{Outline for proof of \eqref{eq:bg2}: second-order Boltzmann-Gibbs principle}\label{section:bg2}
We break the proof of \eqref{eq:bg2} into two steps: homogenization and Taylor expansion (like in \cite{KPZ}) for a homogenized statistic.
\subsection{Main homogenization}
\eqref{eq:bg2} asks to control a time-integrated heat operator acting on $\mathfrak{q}\mathbf{Z}$, where $\mathfrak{q}$ is $\mathrm{QCT}$. The key step is to replace $\mathfrak{q}$ by one of the following \emph{local} (mesoscopic) equilibrium expectations, which are much smoother objects (as a function of space) than $\mathfrak{q}$ itself. (Namely, homogenization at \emph{local} scales. The smoothness of these expectations is important when we Taylor expand them.) For an intuitive description of the following construction, see immediately after Definition \ref{definition:bg21}.
\begin{definition}\label{definition:bg21}
 First, fix any integer $\mathfrak{l}\geq1$, and consider the following construction.
\begin{itemize}
\item Define $\mathbb{I}(\mathfrak{l},+):=\llbracket1,\mathfrak{l}\rrbracket$ and $\mathbb{I}(\mathfrak{l},-):=\llbracket-\mathfrak{l}+1,0\rrbracket$. (These are intervals of length $\mathfrak{l}$ pointing to the right and left, respectively, centered at the origin.) Now, we define the following $\mathbf{U}$ density on these intervals (shifted to be centered at $\y$) at time $\s$:
\begin{align}
\sigma(\s,\y;\mathfrak{l},\pm) \ := \ \mathfrak{l}^{-1}{\textstyle\sum_{\mathrm{j}\in\mathbb{I}(\mathfrak{l},\pm)}}\mathbf{U}^{\s,\y+\mathrm{j}}. \label{eq:bg21I}
\end{align}
\item We now introduce expectations with respect to canonical measures (see Definition \ref{definition:intro5}) with densities in \eqref{eq:bg21I}. Take any test {function} $\mathsf{F}^{\pm}\in\mathscr{C}^{\infty}(\R^{\mathbb{I}(\mathfrak{l},\pm)})$, so that $\mathsf{F}^{\pm}(\mathbf{U})$ depends only on $\mathbf{U}(\x)$ for $\x\in\mathbb{I}(\mathfrak{l},\pm)$. (In particular, it makes sense to take its canonical measure expectation.) Also, recall canonical measure expectations $\E^{\sigma,\t,\mathbb{I}}$ from Definition \ref{definition:intro5}. Set
\begin{align}
\E^{\mathfrak{l},\pm}[\mathsf{F}^{\pm};\s,\y] \ := \ \E^{\sigma(\s,\y;\mathfrak{l},\pm),\s,\mathbb{I}(\mathfrak{l},\pm)}[\mathsf{F}^{\pm}]. \label{eq:bg21II}
\end{align}
\end{itemize}
Note that \eqref{eq:bg21I}-\eqref{eq:bg21II} are functionals $\R^{\y+\mathbb{I}(\mathfrak{l},\pm)}\to\R$. For a matter of notational convention, when we evaluate these functionals at $\mathbf{U}^{\s,\cdot}\in\R^{\mathbb{T}(\N)}$, we actually evaluate them at the projection of $\mathbf{U}^{\s,\cdot}\in\R^{\mathbb{T}(\N)}$ onto its marginals $\R^{\y+\mathbb{I}(\mathfrak{l},\pm)}$.
\end{definition}
In words, $\sigma(\s,\y;\mathfrak{l},\pm)$ is the average charge at time $\s$ on the block $\y+\mathbb{I}(\mathfrak{l},\pm)$, which comes with an orientation $\pm$. On the other hand, $\E^{\mathfrak{l},\pm}(\cdot;\s,\y)$ is the associated canonical measure expectation on a domain of length $\mathfrak{l}$ with orientation $\pm$ centered at the space-time point $(\s,\y)$. (All of this extends to any real $\mathfrak{l}\geq0$ upon replacing $\mathfrak{l}\mapsto\lfloor\mathfrak{l}\rfloor$.) 

Let us clarify the orientation $\pm$. We eventually want to multiply $\E^{\mathfrak{l},\pm}[\mathsf{F}^{\pm};\s,\y(\s)]$-terms by the Gartner map $\mathbf{Z}$; see Lemma \ref{lemma:method7}. By Definition \ref{definition:intro6} this is the same as multiplying $\E^{\mathfrak{l},\pm}[\mathsf{F}^{\pm};\s,\y]$ by $\mathbf{G}$ and shifting $\y\mapsto\y(\s)$. Fix the orientation $+$. In this case $\E^{\mathfrak{l},\pm}[\mathsf{F}^{\pm};\s,\y(\s)]$ depends only on $\mathbf{U}^{\s,\y+\mathrm{j}}$ for $\mathrm{j}>0$. Also, $\mathbf{G}(\s,\y)$ depends only on $\mathbf{U}^{\s,\x}$ for $\x\leq\y$ (modulo the height shift $\mathbf{G}(\s,0)$ that does not affect the local dynamics of $\mathbf{G}(\s,\y)$, except for the trivial case $\y=0$). Therefore both factors $\E^{\mathfrak{l},\pm}[\mathsf{F}^{\pm};\s,\y]$ and $\mathbf{G}(\s,\y)$ are, in some sense, ``decoupled" as functionals of $\mathbf{U}^{\s,\cdot}$. For the orientation $-$, note $\E^{\mathfrak{l},\pm}[\mathsf{F}^{\pm};\s,\y]$ depends only on $\mathbf{U}^{\s,\y-\mathrm{j}}$ for $\mathrm{j}\geq0$. On the other hand, $\mathbf{G}(\s,\y)$ depends on $\mathbf{U}^{\s,\x}$ for $\x\leq\y$ only through their sum. As the sum of $\mathbf{U}^{\s,\x}$ over all $\x\in\mathbb{T}(\N)$ is constant, this means $\mathbf{G}(\s,\y)$ depends only on the average $\mathbf{U}^{\s,\x}$ over $\x\not\leq\y$. ``Decoupling" therefore also holds in this case.

Although the terms in \eqref{eq:bg21I} are functionals on $\R^{\y+\mathbb{I}(\mathfrak{l},\pm)}$ that we evaluate at $\mathbf{U}^{\s,\cdot}$, we have chosen notation that instead portrays them as space-time functions on $\R\times\mathbb{T}(\N)$. (This is also the use of the notation \eqref{eq:bg21II}.) This is a matter of convenience; we will think of and use \eqref{eq:bg21I}-\eqref{eq:bg21II} as space-time functions with a priori stochastic estimates that give us ``analytic" bounds after space-time integration. (See, for example, Proposition \ref{prop:bg22}.) It is only in proving said a priori stochastic estimates when we look into the structure of \eqref{eq:bg21I}-\eqref{eq:bg21II} as functions on $\R^{\y+\mathbb{I}(\mathfrak{l},\pm)}$. (We will do this for various space-time dependent functionals of \eqref{eq:glsde} for the same reason, such as the ones to be introduced in Section \ref{section:bg22proofoutline}.)
\subsubsection{The main result for this step}
Let us now explain the following result intuitively, before giving its precise statement. Let $\mathfrak{q}$ by a local $\mathrm{QCT}$ statistic. We expect it to locally homogenize/average out, so we should be able to replace it by its homogenized version (i.e. \eqref{eq:bg21II} with $\mathsf{F}=\mathfrak{q}$). Below, we take $\mathfrak{l}$ of order $\gg\N^{2/3}$; the reason for this was explained after Lemma \ref{lemma:method4}. Proposition \ref{prop:bg22} estimates the error in this replacement. (The only other point that maybe asks for clarification is the orientation $\pm$ appearing below. The error term ${\mathfrak{z}}$ in Proposition \ref{prop:method2} has a $\mathrm{QCT}$ term both with and without a shift in space by $+1$. Now, see the paragraph after Definition \ref{definition:bg21}.)
\begin{prop}\label{prop:bg22}
 Set $\mathfrak{l}(\mathrm{hom}):=\lfloor\N^{2/3+\gamma_{\mathrm{KL}}}\rfloor$. Take $\mathfrak{q}(\t,\cdot)\in\mathrm{QCT}$ jointly smooth and satisfying $|\partial_{\mathbf{u}}^{\d}\mathfrak{q}(\t,\mathbf{u})|\lesssim_{\d}1+\mathbf{u}^{10}$. Define the following ``replacement functionals" (that account for the error in replacing $\mathfrak{q}$ by its length-$\mathfrak{l}(\mathrm{hom})$-local average):
\begin{align}
\mathds{R}^{\mathfrak{q},-}(\s,\y) \ &:= \ \mathfrak{q}(\s,\mathbf{U}^{\s,\y})-\E^{\mathfrak{l}(\mathrm{hom}),-}[\mathfrak{q}(\s,\cdot);\s,\y]\nonumber\\
\mathds{R}^{\mathfrak{q},+}(\s,\y) \ &:= \ \mathfrak{q}(\s,\mathbf{U}^{\s,\y+1})-\E^{\mathfrak{l}(\mathrm{hom}),+}[\mathfrak{q}(\s,\cdot);\s,\y]. \nonumber
\end{align}
Let us also define $\mathscr{R}^{\mathfrak{q}}(\t,\x):=\mathscr{R}^{\mathfrak{q},-}(\t,\x)+\mathscr{R}^{\mathfrak{q},+}(\t,\x)$, where $\mathscr{R}^{\mathfrak{q},\pm}$ are the following integrated versions of $\mathds{R}^{\mathfrak{q},\pm}$:
\begin{align}
\mathscr{R}^{\mathfrak{q},\pm}(\t,\x)\ := \ {\textstyle\int_{0}^{\t}}\mathbf{H}^{\N}(\s,\t(\N),\x)\{\N\mathds{R}^{\mathfrak{q},\pm}(\s,\cdot(\s))\mathbf{Z}(\s,\cdot)\}\d\s.
\end{align}
If we recall $\beta_{\mathrm{BG}}$ from {Definition \ref{definition:method8}}, then with high probability, we have
\begin{align}
\|\mathscr{R}^{\mathfrak{q}}\|_{\t_{\mathrm{st}};\mathbb{T}(\N)} \ \leq \ \|\mathscr{R}^{\mathfrak{q},-}\|_{\t_{\mathrm{st}};\mathbb{T}(\N)}+\|\mathscr{R}^{\mathfrak{q},+}\|_{\t_{\mathrm{st}};\mathbb{T}(\N)} \ \lesssim \  \N^{-3\beta_{\mathrm{BG}}}\|\mathbf{Z}\|_{\t_{\mathrm{st}};\mathbb{T}(\N)}. \label{eq:bg22I}
\end{align}
\end{prop}
%
\subsection{Taylor expansion}
Following the discussion after Lemma \ref{lemma:method4}, we now take advantage of homogenization (Proposition \ref{prop:bg22}) to show asymptotic vanishing of $\E^{\mathfrak{l}(\mathrm{hom}),\pm}(\mathfrak{q}(\s,\cdot);\s,\y)$. As noted there and as in \cite{KPZ}, this is a Taylor expansion-type calculation.
\begin{lemma}\label{lemma:bg23}
 Recall $\gamma_{\mathrm{reg}}$ in {Definition \ref{definition:reg}} and $\beta_{\mathrm{BG}}$ in {Definition \ref{definition:method8}}. Take $\mathfrak{l}\geq  \N^{\beta_{\mathrm{BG}}}$. With probability 1, we have
\begin{align}
\|\E^{\mathfrak{l},\pm}[\mathfrak{q}(\t,\cdot);\t,\x]\|_{\t_{\mathrm{st}};\mathbb{T}(\N)}\lesssim  \N^{10\gamma_{\mathrm{reg}}}\mathfrak{l}^{-\frac32}. \label{eq:bg23I}
\end{align}
\end{lemma}
The norm on the LHS of \eqref{eq:bg23I} is with respect to $(\t,\x)$ on the LHS. \eqref{eq:bg23I} will follow by a Taylor expansion (see after Lemma \ref{lemma:method4}). Indeed, if we replace canonical expectation by a grand-canonical one in $\mathrm{LHS}\eqref{eq:bg23I}$, Taylor expansion around $\sigma=0$ implies $\mathrm{LHS}\eqref{eq:bg23I}\lesssim|\sigma(\t,\x;\mathfrak{l},\pm)|^{3}$, since $\mathfrak{q}\in\mathrm{QCT}$. Square-root cancellation (see Remark \ref{remark:intro14}) gives $|\sigma(\t,\x;\mathfrak{l},\pm)|^{3}\lesssim\N^{\mathrm{O}(\gamma_{\mathrm{reg}})}\mathfrak{l}^{-3/2}$. So, we are left to replace canonical expectations by grand-canonical ones in $\mathrm{LHS}\eqref{eq:bg23I}$. This is the usual \emph{equivalence of ensembles}.
\subsection{Proof of \eqref{eq:bg2} given Proposition \ref{prop:bg22} and Lemma \ref{lemma:bg23}}
This argument is essentially ``replace a local $\mathrm{QCT}$ statistic by its homogenized version via Proposition \ref{prop:bg22}, and bound the homogenized version by Lemma \ref{lemma:bg23}". Let us make this precise.

We claim {that} $\mathfrak{q}(\t,\mathbf{u})={\lambda(\t)(\mathscr{U}'(\t,\mathbf{u})-\bar{\alpha}(\t)\mathbf{u}-2^{-1}\lambda(\t)[\mathscr{U}'(\t,\mathbf{u})\mathbf{u}-1])}$ satisfies constraints of Proposition \ref{prop:bg22}. (This is by Lemma \ref{lemma:method4} and Assumption \ref{ass:intro8}.) Next, we claim that with high probability,
\begin{align}
\|\mathrm{QCT}\|_{\t_{\mathrm{st}};\mathbb{T}(\N)} \ &\lesssim \ \|\mathscr{R}^{\mathfrak{q},-}\|_{\t_{\mathrm{st}};\mathbb{T}(\N)}+\|\mathscr{R}^{\mathfrak{q},+}\|_{\t_{\mathrm{st}};\mathbb{T}(\N)} + \|\mathrm{Fin}^{\mathfrak{q},-}\|_{\t_{\mathrm{st}};\mathbb{T}(\N)}+\|\mathrm{Fin}^{\mathfrak{q},+}\|_{\t_{\mathrm{st}};\mathbb{T}(\N)} \label{eq:bg21a} \\
&\lesssim \ \N^{-3\beta_{\mathrm{BG}}}\|\mathbf{Z}\|_{\t_{\mathrm{st}};\mathbb{T}(\N)}+\|\mathrm{Fin}^{\mathfrak{q},-}\|_{\t_{\mathrm{st}};\mathbb{T}(\N)}+\|\mathrm{Fin}^{\mathfrak{q},+}\|_{\t_{\mathrm{st}};\mathbb{T}(\N)}, \label{eq:bg21b}
\end{align}
where we have introduced $\mathrm{Fin}^{\mathfrak{q},\pm}$, which are the space-time integrals below (for the homogenized $\mathfrak{q}$):
\begin{align}
\mathrm{Fin}^{\mathfrak{q},\pm}(\t,\x) \ := \ {\textstyle\int_{0}^{\t}}\mathbf{H}^{\N}(\s,\t(\N),\x)\{\N\E^{\mathfrak{l}(\mathrm{hom}),\pm}[\mathfrak{q}(\s,\mathbf{u});\s,\cdot(\s)]\mathbf{Z}(\s,\cdot)\}\d\s.
\end{align}
(In the definition of $\mathrm{Fin}^{\mathfrak{q},\pm}$, $\mathbf{u}$ is the expectation-dummy-variable.) \eqref{eq:bg21a} is by triangle inequality. \eqref{eq:bg21b} is by Proposition \ref{prop:bg22}. By contractivity of the $\mathbf{H}^{\N}$ operator ({see \eqref{eq:hke3} with $\mathrm{m}=0$}), by $\t_{\mathrm{st}}\leq1$, and by Lemma \ref{lemma:bg23}, we have the deterministic estimates
\begin{align}
\|\mathrm{Fin}^{\mathfrak{q},\pm}\| &\lesssim \N\|\t_{\mathrm{st}}\times|\E^{\mathfrak{l}(\mathrm{hom}),\pm}(\mathfrak{q}(\t,\cdot);\t,\x)\mathbf{Z}|\| \nonumber\\
&\lesssim \N^{1+10\gamma_{\mathrm{reg}}}\mathfrak{l}(\mathrm{hom})^{-\frac32}\|\mathbf{Z}\| \lesssim \N^{-\frac32\gamma_{\mathrm{KL}}+10\gamma_{\mathrm{reg}}}\|\mathbf{Z}\| \lesssim \N^{-3\beta_{\mathrm{BG}}}\|\mathbf{Z}\|, \label{eq:bg21c}
\end{align}
where $\|\|=\|\|_{\t_{\mathrm{st}}}$. (The last bound follows because $\gamma_{\mathrm{reg}},\beta_{\mathrm{BG}}\leq{c}\gamma_{\mathrm{KL}}$ {for some small but fixed $c>0$}.) \eqref{eq:bg21a}-\eqref{eq:bg21c} complete the proof. \qed
\subsubsection{A word about {Proposition \ref{prop:bg22}} (and {Lemma \ref{lemma:bg23}})}\label{subsection:bg2word}
Let us make a short clarification that will be useful for proving \eqref{eq:bg1} and \eqref{eq:hl}. For the proof of Proposition \ref{prop:bg22}, we only use the assumption of $\mathfrak{q}(\t,\cdot)\in\mathrm{QCT}$ in the capacity that Lemma \ref{lemma:bg23} holds for $\mathfrak{q}(\t,\cdot)$. Moreover, instead of evaluating the heat operator in $\mathscr{R}^{\mathfrak{q},\pm}$ in Proposition \ref{prop:bg22} at forward time $\t(\N)=\t+\N^{-100\gamma_{\mathrm{reg}}}$, it suffices to evaluate it at the forward time $\t+2^{-1}\N^{-100\gamma_{\mathrm{reg}}}$. All we need is that this forward time, which depends on $\t$, is separated from $[0,\t]$ by $\gtrsim\N^{-100\gamma_{\mathrm{reg}}}$. (This regularizes the $\mathbf{H}^{\N}$ heat kernel.)
%
%
%
\section{Proofs of \eqref{eq:bg1},\eqref{eq:hl}: first-order Boltzmann-Gibbs principle and hydrodynamic limit estimate}\label{section:bg1hlproofoutline}
This section has the same ideas and architecture as Section \ref{section:bg2}. (It turns out to be easier for technical reasons that are perhaps irrelevant.) So, the reader is invited to skip it in a first reading. For the same reason, in this section, we give analogs of Proposition \ref{prop:bg22} and Lemma \ref{lemma:bg23} then show \eqref{eq:bg1}, \eqref{eq:hl}.

The following result is an estimate for space-time integrals against functionals that belong to $\mathrm{LCT}$. In particular, the result below will be used with the special choice of $\mathfrak{d}(\t,\mathbf{U}):=\mathscr{W}'(\t,\mathbf{U})=\mathscr{U}'(\t,\mathbf{u})-\bar{\alpha}(\t)\mathbf{u}$; see Proposition \ref{prop:method2} for $\mathscr{W}$, and see Lemma \ref{lemma:method4} for the fact that this choice of functional belongs to $\mathrm{LCT}$.
\begin{prop}\label{prop:bg1hl1}
 Suppose $\mathfrak{d}\in\mathrm{LCT}$ and $\mathfrak{w}\in\mathrm{CT}$ satisfy $|\partial_{\mathbf{u}}^{\d}\mathfrak{d}(\t,\mathbf{u})|+|\partial_{\mathbf{u}}^{\d}\mathfrak{w}(\t,\mathbf{u})|\lesssim_{\d}1+\mathbf{u}^{10}$. Set
\begin{align}
\mathds{R}^{\mathfrak{d},-}(\s,\y) \ &:= \ \mathfrak{d}(\s,\mathbf{U}^{\s,\y})-\E^{\mathfrak{l}(\mathrm{hom}),-}(\mathfrak{d}(\s,\cdot);\s,\y) \nonumber\\
\mathds{R}^{\mathfrak{d},+}(\s,\y) \ &:= \ \mathfrak{d}(\s,\mathbf{U}^{\s,\y+1})-\E^{\mathfrak{l}(\mathrm{hom}),+}(\mathfrak{d}(\s,\cdot);\s,\y) \nonumber \\
\mathds{R}^{\mathfrak{w},-}(\s,\y) \ &:= \ \mathfrak{w}(\s,\mathbf{U}^{\s,\y})-\E^{\mathfrak{l}(\mathrm{hom}),-}(\mathfrak{w}(\s,\cdot);\s,\y) \nonumber\\
\mathds{R}^{\mathfrak{w},+}(\s,\y) \ &:= \ \mathfrak{w}(\s,\mathbf{U}^{\s,\y+1})-\E^{\mathfrak{l}(\mathrm{hom}),+}(\mathfrak{w}(\s,\cdot);\s,\y). \nonumber
\end{align}
Let us also define $\mathscr{R}^{\mathfrak{d}}(\t,\x):=\mathscr{R}^{\mathfrak{d},-}(\t,\x)+\mathscr{R}^{\mathfrak{d},+}(\t,\x)$, where $\mathscr{R}^{\mathfrak{d},\pm}$ are the following integrated versions of $\mathds{R}^{\mathfrak{d},\pm}$:
\begin{align}
\mathscr{R}^{\mathfrak{d},?}(\t,\x)\ := \ {\textstyle\int_{0}^{\t}}\mathbf{H}^{\N}(\s,\t(\N),\x)\{\N^{\frac32}\grad^{?}[\mathds{R}^{\mathfrak{d},?}(\s,\cdot(\s))\mathbf{Z}(\s,\cdot)]\}\d\s.
\end{align}
Finally, define the following set of space-time integrals but catered to $\mathfrak{w}$:
\begin{align}
\mathscr{R}^{\mathfrak{w},1,\pm}(\t,\x)\ &:= \ {\textstyle\int_{0}^{\t}}\mathbf{H}^{\N}(\s,\t(\N),\x)\{\mathds{R}^{\mathfrak{w},\pm}(\s,\cdot(\s))\mathbf{Z}(\s,\cdot)\}\d\s \\
\mathscr{R}^{\mathfrak{w},2,+}(\t,\x)\ &:= \ {\textstyle\int_{0}^{\t}}\mathbf{H}^{\N}(\s,\t(\N),\x)\{\N\grad^{-}[\mathds{R}^{\mathfrak{w},+}(\s,\cdot(\s))\mathbf{Z}(\s,\cdot)]\}\d\s \\
\mathscr{R}^{\mathfrak{w},2,-}(\t,\x)\ &:= \ {\textstyle\int_{0}^{\t}}\mathbf{H}^{\N}(\s,\t(\N),\x)\{\N\grad^{+}[\mathds{R}^{\mathfrak{w},-}(\s,\cdot(\s))\mathbf{Z}(\s,\cdot)]\}\d\s.
\end{align}
Recall $\beta_{\mathrm{BG}}$ from {Definition \ref{definition:method8}}, and set $\|\|:=\|\|_{\t_{\mathrm{st}};\mathbb{T}(\N)}$ for convenience. With high probability, we have 
\begin{align}
\|\mathscr{R}^{\mathfrak{d}}\| + \sum_{?=\pm}\|\mathscr{R}^{\mathfrak{w},1,?}\| + \sum_{?=\pm}\|\mathscr{R}^{\mathfrak{w},2,?}\| \ &\leq \ \|\mathscr{R}^{\mathfrak{d},+}\| + \|\mathscr{R}^{\mathfrak{d},-}\| + \sum_{?=\pm}\|\mathscr{R}^{\mathfrak{w},1,?}\| + \sum_{?=\pm}\|\mathscr{R}^{\mathfrak{w},2,?}\| \nonumber\\
&\lesssim \  \N^{-2\beta_{\mathrm{BG}}}\|\mathbf{Z}\|. \label{eq:bg1hl1I}
\end{align}
\end{prop}
\begin{lemma}\label{lemma:bg1hl2}
 Recall $\gamma_{\mathrm{reg}}$ in {Definition \ref{definition:reg}} and $\beta_{\mathrm{BG}}$ in {Definition \ref{definition:method8}}. Recall {Definition \ref{definition:bg21}}. Take $\mathfrak{l}\geq  \N^{\beta_{\mathrm{BG}}}$. With probability 1,
\begin{align}
\|\E^{\mathfrak{l},\pm}(\mathfrak{d}(\t,\cdot);\t,\x)\|_{\t_{\mathrm{st}};\mathbb{T}(\N)} \ &\lesssim \  \N^{10\gamma_{\mathrm{reg}}}\mathfrak{l}^{-1} \label{eq:bg1hl2Ia}\\
\|\E^{\mathfrak{l},\pm}(\mathfrak{w}(\t,\cdot);\t,\x)\|_{\t_{\mathrm{st}};\mathbb{T}(\N)} \ &\lesssim \  \N^{10\gamma_{\mathrm{reg}}}\mathfrak{l}^{-\frac12}. \label{eq:bg1hl2Ib}
\end{align}
\end{lemma}
(For heuristic proof of Lemma \ref{lemma:bg1hl2}, see the paragraph after Lemma \ref{lemma:bg23}. The only difference is that $\mathrm{LCT}$ functionals subtract only linear projections. So $\mathrm{LHS}\eqref{eq:bg1hl2Ia}\lesssim|\sigma(\t,\x;\mathfrak{l},\pm)|^{2}$. Similarly, $\mathrm{CT}$ functionals only de-mean, so $\mathrm{LHS}\eqref{eq:bg1hl2Ib}\lesssim|\sigma(\t,\x;\mathfrak{l},\pm)|$.)
\subsection{Proof of \eqref{eq:bg1},\eqref{eq:hl} given Proposition \ref{prop:bg1hl1} and Lemma \ref{lemma:bg1hl2}}
We first show \eqref{eq:bg1}. Note $\mathfrak{d}(\t,\mathbf{u}):=\mathscr{U}'(\t,\mathbf{u})-\bar{\alpha}(\t)\mathbf{u}$ is $\mathrm{LCT}$ and $|\partial_{\mathbf{u}}^{\d}\mathfrak{d}(\t,\mathbf{u})|\lesssim_{\d}1+\mathbf{u}^{10}$; see Lemma \ref{lemma:method4} and Assumption \ref{ass:intro8}. By Proposition \ref{prop:bg1hl1} and explanation of \eqref{eq:bg21a}-\eqref{eq:bg21b},
\begin{align}
\|\mathrm{LCT}\| \ \lesssim \ \|\mathscr{R}^{\mathfrak{d},+}\|+\|\mathscr{R}^{\mathfrak{d},-}\| +\|\mathrm{Fin}^{\mathfrak{d},+}\|+\|\mathrm{Fin}^{\mathfrak{d},-}\| \ \lesssim \  \N^{-2\beta_{\mathrm{BG}}}\|\mathbf{Z}\|+\|\mathrm{Fin}^{\mathfrak{d},+}\|+\|\mathrm{Fin}^{\mathfrak{d},-}\| \label{eq:bg11}
\end{align}
with high probability, where $\|\|=\|\|_{\t_{\mathrm{st}};\mathbb{T}(\N)}$, and where $\mathrm{Fin}^{\mathfrak{d},\pm}$ is defined by
\begin{align}
\mathrm{Fin}^{\mathfrak{d},?}(\t,\x) \ := \ {\textstyle\int_{0}^{\t}}\mathbf{H}^{\N}(\s,\t(\N),\x)\{\N^{\frac32}\grad^{?}[\E^{\mathfrak{l}(\mathrm{hom}),\pm}(\mathfrak{d}(\s,\mathbf{u});\s,\cdot(\s))\mathbf{Z}(\s,\cdot)]\}\d\s.
\end{align}
As the operator $\mathbf{H}^{\N}$ is a convolution operator ({since it is the semigroup for a spatially homogeneous infinitesimal generator}), any constant-coefficient discrete gradient $\grad^{?}$ commutes with $\mathbf{H}^{\N}$. By the operator bound for $\N\grad^{?}\mathbf{H}^{\N}$ in {\eqref{eq:hke3}} and by Lemma \ref{lemma:bg1hl2}, with probability 1, we have
\begin{align}
\|\mathrm{Fin}^{\mathfrak{d},\pm}\| \ &\lesssim \ \sup_{0\leq\t\leq\t_{\mathrm{st}}}{\textstyle\int_{0}^{\t}}|\t(\N)-\s|^{-\frac12}\d\s \times \|\N^{\frac12}\E^{\mathfrak{l}(\mathrm{hom}),\pm}(\mathfrak{d}(\s,\mathbf{u});\s,\cdot(\s))\mathbf{Z}(\s,\cdot)\| \nonumber\\
&\lesssim \ \t_{\mathrm{st}}^{\frac12}\N^{\frac12+10\gamma_{\mathrm{reg}}}\mathfrak{l}(\mathrm{hom})^{-1}\|\mathbf{Z}\|, \label{eq:bg11b}
\end{align}
(The second bound requires $\t(\N)\geq\t$, which follows from construction in Definition \ref{definition:method5}, and doing the time-integral in \eqref{eq:bg11b}.) Because $\mathfrak{l}(\mathrm{hom})\geq\N^{2/3}$ (see Proposition \ref{prop:bg22}) and ${\gamma_{\mathrm{reg}}}>0$ is small (see Definitions \ref{definition:entropydata}, \ref{definition:reg}), the far RHS of \eqref{eq:bg11b} is $\lesssim\N^{-2\beta_{\mathrm{BG}}}\|\mathbf{Z}\|$. Using this with \eqref{eq:bg11} gives \eqref{eq:bg1} with high probability. Now, we prove \eqref{eq:hl}. We set $\mathfrak{w}(\t,\mathbf{u};1)=\bar{\alpha}(\t)\mathbf{u}^{2}-1$ and $\mathfrak{w}(\t,\mathbf{u};2)=2^{-1}\lambda(\t)^{4}\E^{0,\t}[\mathscr{U}'(\t,\mathbf{u})\mathbf{u}^{3}]+6^{-1}\lambda(\t)^{3}\bar{\alpha}(\t)\mathbf{u}^{3}-\lambda(\t)\mathscr{R}(\t)$ and $\mathfrak{w}(\t,\mathbf{u};3)=2^{-1}\lambda(\t)^{4}\{\mathscr{U}'(\t,\mathbf{u})\mathbf{u}^{3}-\E^{0,\t}[\mathscr{U}'(\t,\mathbf{u})\mathbf{u}^{3}]\}$. Now, note that $\mathfrak{w}(\t,\cdot;\mathrm{k})$ is $\mathrm{CT}$ and $|\partial_{\mathbf{u}}^{\d}\mathfrak{w}(\t,\mathbf{u};\mathrm{k})|\lesssim_{\d}1+\mathbf{u}^{10}$ for $\mathrm{k}=1,2,3$. This follows by Lemma \ref{lemma:method4} and Assumption \ref{ass:intro8}. Thus, Proposition \ref{prop:bg1hl1} plus the explanation of \eqref{eq:bg21a}-\eqref{eq:bg21b} and \eqref{eq:bg11} give the following high probability estimate:
\begin{align}
\|\mathrm{CT}\| \ &\lesssim \ \sum_{\mathrm{k}=1,2,3}\sum_{\mathrm{n}=1,2}\sum_{?=\pm}\{\|\mathscr{R}^{\mathfrak{w}(\cdot,\cdot;\mathrm{k}),\mathrm{n},?}\|+\|\mathrm{Fin}^{\mathfrak{w}(\cdot,\cdot;\mathrm{k}),\mathrm{n},?}\|\} \nonumber\\
&\lesssim \  \N^{-2\beta_{\mathrm{BG}}}\|\mathbf{Z}\|+\sum_{\mathrm{k}=1,2,3}\sum_{\mathrm{n}=1,2}\sum_{?=\pm}\|\mathrm{Fin}^{\mathfrak{w}(\cdot,\cdot;\mathrm{k}),\mathrm{n},?}\|\label{eq:hl1a}
\end{align}
where
\begin{align}
\mathrm{Fin}^{\mathfrak{w}(\cdot,\cdot;\mathrm{k}),1,\pm}(\t,\x)\ &:= \ {\textstyle\int_{0}^{\t}}\mathbf{H}^{\N}(\s,\t(\N),\x)\{\E^{\mathfrak{l}(\mathrm{hom}),\pm}(\mathfrak{w}(\s,\mathbf{u};\mathrm{k});\s,\cdot(\s))\mathbf{Z}(\s,\cdot)\}\d\s \label{eq:hl1b}\\
\mathrm{Fin}^{\mathfrak{w}(\cdot,\cdot;\mathrm{k}),2,+}(\t,\x)\ &:= \ {\textstyle\int_{0}^{\t}}\mathbf{H}^{\N}(\s,\t(\N),\x)\{\N\grad^{-}[\E^{\mathfrak{l}(\mathrm{hom}),+}((\mathfrak{w}(\s,\mathbf{u};\mathrm{k});\s,\cdot(\s))\mathbf{Z}(\s,\cdot)]\}\d\s\label{eq:hl1c}\\
\mathrm{Fin}^{\mathfrak{w}(\cdot,\cdot;\mathrm{k}),2,-}(\t,\x)\ &:= \ {\textstyle\int_{0}^{\t}}\mathbf{H}^{\N}(\s,\t(\N),\x)\{\N\grad^{+}[\E^{\mathfrak{l}(\mathrm{hom}),-}((\mathfrak{w}(\s,\mathbf{u};\mathrm{k});\s,\cdot(\s))\mathbf{Z}(\s,\cdot)]\}\d\s.\label{eq:hl1d}
\end{align}
Lemma \ref{lemma:bg1hl2} implies each expectation in \eqref{eq:hl1b}-\eqref{eq:hl1d} is $\lesssim\N^{10\gamma_{\mathrm{reg}}}\mathfrak{l}(\mathrm{hom})^{-1/2}$. Thus, similar to \eqref{eq:bg21c} and \eqref{eq:bg11b}, we get
\begin{align}
\|\mathrm{Fin}^{\mathfrak{w}(\cdot,\cdot;\mathrm{k}),\mathrm{n},?}\|  \ \lesssim \ \N^{10\gamma_{\mathrm{reg}}}\mathfrak{l}(\mathrm{hom})^{-\frac12}\|\mathbf{Z}\| \ \lesssim \ \N^{-2\beta_{\mathrm{BG}}}\|\mathbf{Z}\|
\end{align}
with probability 1, where the last bound follows by $\mathfrak{l}(\mathrm{hom})\geq\N^{2/3}$. Plugging this into \eqref{eq:hl1a} gives \eqref{eq:hl}, so we are done. \qed
\subsection{What is left}
As we explained right before Section \ref{section:bg2}, proofs of Proposition \ref{prop:method2}, and Lemmas \ref{lemma:method4}, \ref{lemma:method6}, \ref{lemma:method7}, \ref{lemma:method10}, \ref{lemma:method13} will be given in the appendix. We must now show Proposition \ref{prop:bg22}, Lemma \ref{lemma:bg23}, Proposition \ref{prop:bg1hl1}, and Lemma \ref{lemma:bg1hl2}. Lemmas \ref{lemma:bg23} and \ref{lemma:bg1hl2} are not as interesting, since slightly weaker versions are shown in \cite{DGP}, for example. For this reason, we defer their proofs to the appendix as well. The proof of Proposition \ref{prop:bg1hl1} will ultimately follow from the same argument as Proposition \ref{prop:bg22}. (In a nutshell, Propositions \ref{prop:bg1hl1} and \ref{prop:bg22} ask for the same thing, but the former asks for a weaker estimate. We make this precise in the appendix since it is elementary.) So, the proof of Proposition \ref{prop:bg22} (and of Corollary \ref{corollary:kpz}) is all we have left before the appendix.
%
%
%
\section{Outline for proof of Proposition \ref{prop:bg22}}\label{section:bg22proofoutline}
Again, we give here ingredients for the proof of Proposition \ref{prop:bg22}. We defer their proofs to forthcoming sections. While doing so, we give intuitive descriptions of what each ingredient says, and a more-than-intuitive explanation for why it is true and how it is proved.

{Before we begin, let us emphasize that for the rest of this paper, any estimates concerning functions $\mathfrak{q}\in\mathrm{QCT}$ depend only smoothness of $\mathfrak{q}$ and sub-polynomial bounds on both it and its derivatives.}
\subsection{Multiscale I}
Proposition \ref{prop:bg22} asks to estimate the cost in replacing $\mathfrak{q}$ by its homogenized expectation with respect to length-scale $\mathfrak{l}(\mathrm{hom})$. It turns out to be much more effective to first replace $\mathfrak{q}$ by homogenized expectation on some length-scale $\mathfrak{l}$. Then, we interpolate between $\mathfrak{l}$ and $\mathfrak{l}(\mathrm{hom})$. To see why, all of our homogenization is done by first comparing the law of the processes to local equilibrium measures from Definition \ref{definition:intro5} on mesoscopic space-time scales. (This is the folklore ``local equilibration" of many-body processes.) The larger (or ``less local") the scale, the harder it is to compare to any local equilibrium. This makes our analysis, when interpolating $\mathfrak{l}\mapsto\mathfrak{l}(\mathrm{hom})$, to deteriorate at larger scales. However, by Lemma \ref{lemma:bg23}, homogenized expectations are better controlled a priori at larger scales. These competing factors ultimately cancel out.
\begin{definition}\label{definition:bg24}
 Recall $\beta_{\mathrm{BG}}$ from Definition \ref{definition:method8}. Fix $2\beta_{\mathrm{BG}}\leq\delta_{\mathrm{BG}}\leq3\beta_{\mathrm{BG}}$ so that $\N^{\delta_{\mathrm{BG}}}$ is an integer and $\N^{\mathrm{j}(\infty)\delta_{\mathrm{BG}}}=\mathfrak{l}(\mathrm{hom})$ from Proposition \ref{prop:bg22} for some positive integer $\mathrm{j}(\infty)$. Set $\mathfrak{l}(\mathrm{j}):=\N^{\mathrm{j}\delta_{\mathrm{BG}}}\wedge\mathfrak{l}(\mathrm{hom})$ for $\mathrm{j}\geq1$. Also, for $\mathrm{j}>1$, we define the ``renormalization" functional (see \cite{GJ15} for the naming) below, interpolating the length-scale (in the canonical measure expectation of our $\mathrm{QCT}$-functional $\mathfrak{q}$) from $\mathfrak{l}(\mathrm{j}-1)$ to $\mathfrak{l}(\mathrm{j})$:
\begin{align}
\mathds{R}^{\mathfrak{q},\pm,\mathrm{j}}(\s,\y) \ := \ \E^{\mathfrak{l}(\mathrm{j}-1),\pm}[\mathfrak{q}(\s,\cdot);\s,\y]-\E^{\mathfrak{l}(\mathrm{j}),\pm}[\mathfrak{q}(\s,\cdot);\s,\y]. \nonumber
\end{align}
On the other hand, for the case $\mathrm{j}=1$, let us instead define the following ``renormalization" functionals:
\begin{align}
\mathds{R}^{\mathfrak{q},-,1}(\s,\y) \ &:= \ \mathfrak{q}(\s,\mathbf{U}^{\s,\y})-\E^{\mathfrak{l}(1),-}[\mathfrak{q}(\s,\cdot);\s,\y] \nonumber\\
\mathds{R}^{\mathfrak{q},+,1}(\s,\y) \ &:= \ \mathfrak{q}(\s,\mathbf{U}^{\s,\y+1})-\E^{\mathfrak{l}(1),+}[\mathfrak{q}(\s,\cdot);\s,\y]. \nonumber
\end{align}
These are functions $\mathds{R}^{\mathfrak{q},\pm,\mathrm{j}}(\s,\mathbf{U})$ for $\mathbf{U}\in\R^{\y+\mathbb{I}(\mathfrak{l}(\mathrm{j}),\pm)}$ that we evaluate at the projection of $\mathbf{U}^{\s,\cdot}$ onto its $\R^{\y+\mathbb{I}(\mathfrak{l}(\mathrm{j}),\pm)}$-marginals. Next, let us define the following space-time (heat-operator) integrated version of the above $\mathds{R}^{\mathfrak{q},\pm,\mathrm{j}}$-functionals:
\begin{align}
\mathscr{R}^{\mathfrak{q},\pm,\mathrm{j}}(\t,\x)\ := \ {\textstyle\int_{0}^{\t}}\mathbf{H}^{\N}(\s,\t(\N),\x)(\N\mathds{R}^{\mathfrak{q},\pm,\mathrm{j}}(\s,\cdot(\s))\mathbf{Z}(\s,\cdot))\d\s. \nonumber
\end{align}
\end{definition}
Throughout this section (and related contexts), the reader should think of $\mathrm{j}$ as length-scale parameters for homogenization.

The following result controls the error in upgrading the length-scale in the canonical expectation. As shown immediately afterwards, said result, {plugged with} the triangle inequality and a bound on the number of interpolation steps needed to get from $\mathfrak{l}(1)$ to $\mathfrak{l}(\mathrm{hom})$, yields Proposition \ref{prop:bg22} after fairly elementary considerations.
\begin{prop}\label{prop:bg25}
 With high probability, we have have the following estimate:
\begin{align}
{\textstyle\sup_{\mathrm{j}}}\|\mathscr{R}^{\mathfrak{q},\pm,\mathrm{j}}\|_{\t_{\mathrm{st}};\mathbb{T}(\N)} \ \lesssim \  \N^{-3\beta_{\mathrm{BG}}}\|\mathbf{Z}\|_{\t_{\mathrm{st}};\mathbb{T}(\N)}. \label{eq:bg25I}
\end{align}
\end{prop}
\begin{proof}[Proof of {Proposition \ref{prop:bg22}}]
It suffices to show that with high probability, we have
\begin{align}
\|\mathscr{R}^{\mathfrak{q},\pm}\|_{\t_{\mathrm{st}};\mathbb{T}(\N)} \ \leq \ {\textstyle\sum_{\mathrm{j}}}\|\mathscr{R}^{\mathfrak{q},\pm,\mathrm{j}}\|_{\t_{\mathrm{st}};\mathbb{T}(\N)} \ \lesssim \ {\textstyle\sup_{\mathrm{j}}}\|\mathscr{R}^{\mathfrak{q},\pm,\mathrm{j}}\|_{\t_{\mathrm{st}};\mathbb{T}(\N)} \ \lesssim \  \N^{-3\beta_{\mathrm{BG}}}\|\mathbf{Z}\|_{\t_{\mathrm{st}};\mathbb{T}(\N)}.
\end{align}
The first estimate is triangle inequality. (Indeed, $\mathscr{R}^{\mathfrak{q},\pm,\mathrm{j}}$ sum over $\mathrm{j}$ to $\mathscr{R}^{\mathfrak{q},\pm}$; see Proposition \ref{prop:bg22} and Definition \ref{definition:bg24}.) The second follows by $\mathrm{j}(\infty)\lesssim1$. (Indeed, observe $\mathscr{R}^{\mathfrak{q},\pm,\mathrm{j}}=0$ for all $\mathrm{j}>\mathrm{j}(\infty)$; see Definition \ref{definition:bg24}. To prove $\mathrm{j}(\infty)\lesssim1$, it is enough to note $\mathrm{j}(\infty)$ is at most the number of steps of size $\gtrsim1$ to get from 0 to 1; see Definition \ref{definition:bg24}.) The last bound is Proposition \ref{prop:bg25}.
\end{proof}
\subsection{A priori (technical) cutoff}
We are now left to show Proposition \ref{prop:bg25}. Recall from the paragraph before Definition \ref{definition:bg24} that $\mathds{R}^{\mathfrak{q},\pm,\mathrm{j}}$-terms in Definition \ref{definition:bg24} should have improving a priori estimates as $\mathrm{j}$ increases. We introduce said a priori bounds below, in a way that preserves another important ``fluctuation" property of $\mathds{R}^{\mathfrak{q},\pm,\mathrm{j}}$. (We clarify Definition \ref{definition:bg26} and all of this shortly.)
\begin{definition}\label{definition:bg26}
 Retain the notation of Definition \ref{definition:bg24}. In what follows, we will also use the functions $\chi(\mathrm{a};\upsilon):=\chi(\upsilon^{-1}\mathrm{a})$, where $\upsilon>0$ and $\mathrm{a}\in\R$. We assume $\chi$ is smooth with support contained in $[-2,2]$ and that $\chi(\mathrm{a})=1$ for all $\mathrm{a}\in[-1,1]$. In particular, $\chi(\cdot;\upsilon)$ has support contained in $[-2\upsilon,2\upsilon]$ and is identically 1 on $[-\upsilon,\upsilon]$. Also, we have $|\chi'(\mathrm{a};\upsilon)|\lesssim\upsilon^{-1}$. Now, for $\mathrm{j}\geq1$, set
\begin{align}
\mathds{R}^{\chi,\mathfrak{q},\pm,\mathrm{j}}(\s,\y) &:= \{\mathds{R}^{\mathfrak{q},\pm,\mathrm{j}}(\s,\y)\cdot\chi[\mathds{R}^{\mathfrak{q},\pm,\mathrm{j}}(\s,\y);\upsilon_{\mathrm{j}-1}]\}\nonumber\\
&- \E^{\mathfrak{l}(\mathrm{j}),\pm}\{\mathds{R}^{\mathfrak{q},\pm,\mathrm{j}}(\s,\mathbf{U})\cdot\chi[\mathds{R}^{\mathfrak{q},\pm,\mathrm{j}}(\s,\mathbf{U});\upsilon_{\mathrm{j}-1}]\}. \label{eq:bg26I}
\end{align}
Above, we defined {$\upsilon_{0}= \N^{20\gamma_{\mathrm{reg}}}$} and $\upsilon_{\mathrm{j}}= \N^{20\gamma_{\mathrm{reg}}}\mathfrak{l}(\mathrm{j})^{-3/2}$ for $\mathrm{j}>1$. We briefly explain what \eqref{eq:bg26I} is. The first term in $\mathrm{RHS}\eqref{eq:bg26I}$ is determined by Definition \ref{definition:bg24}. In particular, it is a function on $\R^{\y+\mathbb{I}(\mathfrak{l}(\mathrm{j}),\pm)}$ evaluated at the time-$\s$ data $\mathbf{U}^{\s,\cdot}$ of the process \eqref{eq:glsde}. The final term in \eqref{eq:bg26I} is the $\E^{\mathfrak{l}(\mathrm{j}),\pm}$-expectation over the ``dummy" variable $\mathbf{U}$ (see Definition \ref{definition:bg21}) of this function. In particular, if we write anything inside $\E^{\mathfrak{l}(\mathrm{j}),\pm}$-expectation, we never evaluate it at $\mathbf{U}^{\s,\cdot}$. We always evaluate it at the expectation ``dummy" variable $\mathbf{U}$ in $\E^{\mathfrak{l}(\mathrm{j}),\pm}$. This is the case for anything in canonical ensemble expectations unless otherwise noted. Now, we define
\begin{align}
\mathscr{R}^{\chi,\mathfrak{q},\pm,\mathrm{j}}(\t,\x)\ := \ {\textstyle\int_{0}^{\t}}\mathbf{H}^{\N}(\s,\t(\N),\x)\{\N\mathds{R}^{\chi,\mathfrak{q},\pm,\mathrm{j}}(\s,\cdot(\s))\mathbf{Z}(\s,\cdot)\}\d\s. \nonumber
\end{align}
\end{definition}
Throughout this section (and related contexts), the reader should think of $\chi$ as denoting cutoff a la Definition \ref{definition:bg26}. It is easy to see that by construction, we have $\mathds{R}^{\chi,\mathfrak{q},\pm,\mathrm{j}}\lesssim\upsilon_{\mathrm{j}-1}$ with probability 1. This implements the improving a priori bounds mentioned before Definition \ref{definition:bg26}. We clarify two important properties of $\mathds{R}^{\chi,\mathfrak{q},\pm,\mathrm{j}}(\s,\y)$. First, it depends only on $\mathbf{U}^{\s,\x}$ for $\x\in\y+\mathbb{I}(\mathfrak{l}(\mathrm{j}),\pm)$. This follows by construction of $\mathds{R}^{\mathfrak{q},\pm,\mathrm{j}}$ in Definition \ref{definition:bg24} and of $\E^{\mathfrak{l}(\mathrm{j}),\pm}$ in Definition \ref{definition:bg21}. Moreover, if we take canonical ensemble expectation of $\mathds{R}^{\chi,\mathfrak{q},\pm,\mathrm{j}}(\s,\y)$ (at time $\s$, given any $\sigma$ and any superset of $\y+\mathbb{I}(\mathfrak{l}(\mathrm{j}),\pm)$), we get zero. (Indeed, $\mathds{R}^{\chi,\mathfrak{q},\pm,\mathrm{j}}$ is centered with respect to canonical measure; see Lemma \ref{lemma:vanishcanonical}.) The same is true of $\mathds{R}^{\mathfrak{q},\pm,\mathrm{j}}$ of Definition \ref{definition:bg24}; see \cite{GJ15}, proof of Lemma 2. 

The following says that the a priori estimate via the $\chi$-cutoff can be introduced with small cost. We explain it further below.
\begin{lemma}\label{lemma:bg27}
 Recall $\beta_{\mathrm{BG}}$ from {Definition \ref{definition:method8}} and the notation of {Definitions \ref{definition:bg24} and \ref{definition:bg26}}. With high probability, we have 
\begin{align}
{\textstyle\sup_{\mathrm{j}}}\|\mathscr{R}^{\mathfrak{q},\pm,\mathrm{j}}-\mathscr{R}^{\chi,\mathfrak{q},\pm,\mathrm{j}}\|_{\t_{\mathrm{st}};\mathbb{T}(\N)} \ \lesssim \  \N^{-3\beta_{\mathrm{BG}}}\|\mathbf{Z}\|_{\t_{\mathrm{st}};\mathbb{T}(\N)}. \label{eq:bg27I}
\end{align}
\end{lemma}
Lemma \ref{lemma:bg27} bounds the error in swapping $\mathds{R}^{\mathfrak{q},\pm,\mathrm{j}}(\s,\y)$ by $\mathds{R}^{\chi,\mathfrak{q},\pm,\mathrm{j}}(\s,\y)$ (upon integrating against heat kernel in space-time). First, note their difference has two terms. (Recall notation from Definitions \ref{definition:bg24} and \ref{definition:bg26}.) The first is given by $\mathds{R}^{\mathfrak{q},\pm,\mathrm{j}}(\s,\y)$ times $1-\chi$, where $\chi$ is short-hand for the first $\chi$-factor in \eqref{eq:bg26I}. Note $1-\chi$ is zero unless $|\mathds{R}^{\mathfrak{q},\pm,\mathrm{j}}(\s,\y)|\gtrsim\upsilon_{\mathrm{j}-1}=\N^{20\gamma_{\mathrm{reg}}}\mathfrak{l}(\mathrm{j}-1)^{-1/2}$. (See Definition \ref{definition:bg26}.) But $\mathds{R}^{\mathfrak{q},\pm,\mathrm{j}}$ is a difference of $\E^{\mathfrak{l}(\mathrm{k}),\pm}$-expectations of $\mathfrak{q}\in\mathrm{QCT}$ for $\mathrm{k}=\mathrm{j}-1,\mathrm{j}$. Therefore, by Lemma \ref{lemma:bg23}, this lower bound for $|\mathds{R}^{\mathfrak{q},\pm,\mathrm{j}}(\s,\y)|$ never happens before time $\t_{\mathrm{st}}$. Thus ${\mathds{R}^{\mathfrak{q},\pm,\mathrm{j}}(\s,\y)}(1-\chi)=0$, and the first term in $\mathds{R}^{\mathfrak{q},\pm,\mathrm{j}}-\mathds{R}^{\chi,\mathfrak{q},\pm,\mathrm{j}}$ is zero. 

The second (and remaining) term in this difference is the final $\E^{\mathfrak{l}(\mathrm{j}),\pm}$-term in \eqref{eq:bg26I}. We would like to, again, drop the $\chi$-cutoff in the $\E^{\mathfrak{l}(\mathrm{j}),\pm}$-term in \eqref{eq:bg26I} to get canonical ensemble expectation of $\mathds{R}^{\mathfrak{q},\pm,\mathrm{j}}$. As we noted right before Lemma \ref{lemma:bg27}, said expectation is zero! This would treat the remaining term in the difference $\mathds{R}^{\mathfrak{q},\pm,\mathrm{j}}-\mathds{R}^{\chi,\mathfrak{q},\pm,\mathrm{j}}$. But, to drop said $\chi$-cutoff, we cannot apply the same argument as before. Indeed, inside $\E^{\mathfrak{l}(\mathrm{j}),\pm}$ in \eqref{eq:bg26I}, the $\mathds{R}^{\mathfrak{q},\pm,\mathrm{j}}\chi(\mathds{R}^{\mathfrak{q},\pm,\mathrm{j}})$-functional is not being evaluated at the time-$\s$ data of the process $\mathbf{U}^{\s,\cdot}$. As clarified in Definition \ref{definition:bg26}, it is being evaluated with respect to the expectation ``dummy" variable $\mathbf{U}\sim\mathbb{P}^{\sigma,\s,\y+\mathbb{I}(\mathfrak{l}(\mathrm{j}),\pm)}$, in which $\sigma$ is the average of $\mathbf{U}^{\s,\z}$ over $\z\in\y+\mathbb{I}(\mathfrak{l}(\mathrm{j}),\pm)$. 

However, the earlier argument  that controls $|\mathds{R}^{\mathfrak{q},\pm,\mathrm{j}}|$ from above, can be salvaged as follows. Recall $\mathds{R}^{\mathfrak{q},\pm,\mathrm{j}}$ is a difference of canonical expectations; see Definition \ref{definition:bg24}. The charge densities of these canonical ensemble expectations are now being sampled according to $\mathbf{U}\sim\mathbb{P}^{\sigma,\s,\y+\mathbb{I}(\mathfrak{l}(\mathrm{j}),\pm)}$, where, again, $\sigma$ is the average of $\mathbf{U}^{\s,\z}$ for all $\z\in\y+\mathbb{I}(\mathfrak{l}(\mathrm{j}),\pm)$. The earlier argument estimates $\mathds{R}^{\mathfrak{q},\pm,\mathrm{j}}$ by controlling these canonical ensemble expectations. In the case where the charge densities for these canonical expectations are given by averages of $\mathbf{U}^{\s,\x}$ (for $\s\leq\t_{\mathrm{st}}\leq\t_{\mathrm{reg}}$ and over some set of $\x$), one can just use Lemma \ref{lemma:bg23}. (This is exactly what we did.) But, by going through the proof of Lemma \ref{lemma:bg23} (and of Lemma \ref{lemma:ee1}), it is easy to see that all we need is the charge density for the canonical ensemble expectation to satisfy a natural Brownian-type bound (see \eqref{eq:ee1I1}). Thus, we need to show that when $\mathbf{U}\sim\mathbb{P}^{\sigma,\s,\y+\mathbb{I}(\mathfrak{l}(\mathrm{j}),\pm)}$, where, again, $\sigma$ is the average of $\mathbf{U}^{\s,\z}$ over $\z\in\y+\mathbb{I}(\mathfrak{l}(\mathrm{j}),\pm)$, the average of $\mathbf{U}(\z)$ over $\z\in\y+\mathbb{I}(\mathfrak{l}(\mathrm{j}),\pm)$ and over $\z\in\y+\mathbb{I}(\mathfrak{l}(\mathrm{j}-1),\pm)$ both admit Brownian-type bounds like \eqref{eq:ee1I1}. (Indeed, these two spatial sets are those that the canonical expectations in $\mathds{R}^{\mathfrak{q},\pm,\mathrm{j}}$ are defined on; see Definition \ref{definition:bg24}.) 

For the average over the bigger set $\y+\mathbb{I}(\mathfrak{l}(\mathrm{j}),\pm)$, by construction, this is simply $\sigma$, the average of $\mathbf{U}^{\s,\z}$ over $\z\in\y+\mathbb{I}(\mathfrak{l}(\mathrm{j}),\pm)$. Thus, by Remark \ref{remark:intro14} we get a Brownian-type bound on the average of $\mathbf{U}(\z)$ over $\z\in\y+\mathbb{I}(\mathfrak{l}(\mathrm{j}),\pm)$. It remains to derive the Brownian-type bound for the average of $\mathbf{U}(\z)$ over $\z\in\y+\mathbb{I}(\mathfrak{l}(\mathrm{j}-1),\pm)$. Again, we have said Brownian estimate for the average over $\z\in\y+\mathbb{I}(\mathfrak{l}(\mathrm{j}),\pm)$. Now note that if $\mathbf{U}\sim\mathbb{P}^{\sigma,\s,\y+\mathbb{I}(\mathfrak{l}(\mathrm{j}),\pm)}$, then $\mathbf{U}(\z)$ are random walk bridge steps. So, it suffices to note that from standard random walk bridge estimates, any increment of said random walk bridge is a martingale (which certainly satisfies Brownian-type bounds with very high probability by the Azuma inequality){,} plus another random variable that is controlled by the average drift of $\mathbf{U}(\z)$. But the average drift of $\mathbf{U}(\z)$, namely the average on the larger set $\z\in\y+\mathbb{I}(\mathfrak{l}(\mathrm{j}),\pm)$, satisfies a Brownian estimate as we already showed. So the average of $\mathbf{U}(\z)$ over $\z\in\y+\mathbb{I}(\mathfrak{l}(\mathrm{j}-1),\pm)$ satisfies a natural Brownian-type bound. As noted several sentences ago, we may now conclude this (formal) proof of Lemma \ref{lemma:bg27}. In the rigorous proof, we just make this argument quantitative.
\subsection{Replacement-by-spatial-average}
Unsurprisingly, in this paper homogenization is done through averaging in space-time. Provided the ``local" spirit of our homogenization analysis, let us adjust $\mathscr{R}^{\chi,\mathfrak{q},\pm,\mathrm{j}}$ in Definition \ref{definition:bg26} by introducing local averaging in the heat operator $\mathbf{H}^{\N}$. We start with spatial averaging. To this end, {we present some convenient notation}.
\begin{definition}\label{definition:bg28}
 Fix $\mathsf{F}\in\mathscr{C}^{\infty}(\R^{\mathbb{T}(\N)})$. Let $\mathfrak{l}(\mathsf{F})$ be the smallest positive integer such that for any $\mathbf{U}\in\R^{\mathbb{T}(\N)}$, $\mathsf{F}(\mathbf{U})$ depends only on $\mathbf{U}(\x)$ for $\x$ in a discrete interval of length $\mathfrak{l}(\mathsf{F})$. We call such $\mathfrak{l}(\mathsf{F})$ the \emph{support length} of $\mathsf{F}$. 

For any integer ${\mathfrak{m}}\geq1$ and $\mathsf{F}(\t,\cdot)\in\mathscr{C}^{\infty}(\R^{\mathbb{T}(\N)})$, set
\begin{align}
{\mathds{A}_{\mathbf{X}}^{{\mathfrak{m}},\pm}}(\mathsf{F}\mathbf{Z};\s,\y(\s)) \ &:= \ {\mathfrak{m}^{-1}}{\textstyle{\sum}_{\mathrm{j}=0}^{{\mathfrak{m}}-1}}\mathsf{F}(\s,\mathbf{U}^{\s,\y(\s)\pm2\mathrm{j}\mathfrak{l}(\mathsf{F})})\mathbf{Z}(\s,\y\pm2\mathrm{j}\mathfrak{l}(\mathsf{F})).
\end{align}
In words, ${\mathds{A}_{\mathbf{X}}^{{\mathfrak{m}},\pm}}(\mathsf{F}\mathbf{Z};\s,\y)$ introduces the spatial-average of {$\mathfrak{m}$}-many shifts of $\mathsf{F}\mathbf{Z}$, where each of these shifts are ``centered" at the space-time point $(\s,\y)$. The $\pm$ denotes ``orientation" of shifting. Lastly, shifting in space by multiples of $2\mathfrak{l}(\mathsf{F})$ is to ensure that no two copies of $\mathsf{F}$ that appear in said average both depend on $\mathbf{U}^{\s,\x}$ for the same $\x$. (So in some sense, the ``supports" of these shifts are mutually disjoint.) 
\end{definition}
We now apply Definition \ref{definition:bg28} to our setting, which we can only clarify after the lemma of this step (Lemma \ref{lemma:bg210}).
\begin{definition}\label{definition:bg29}
 Take ${c}\gamma_{\mathrm{KL}}\leq\alpha(\mathrm{j})\leq2{c}\gamma_{\mathrm{KL}}$ {for some $c>0$ small but fixed} so that $ \N^{3/4+\alpha(\mathrm{j})}$ is a positive integer multiple of $\mathfrak{l}(\mathrm{j})\geq1$ in Definition \ref{definition:bg24}. Next, for $\mathrm{j}\geq1$, set $\mathfrak{m}(\mathrm{j}):= \N^{3/4+\alpha(\mathrm{j})}\mathfrak{l}(\mathrm{j})^{-1}$. Now, define the following time-integrated heat operator actions:
\begin{align}
\mathscr{A}^{\mathbf{X}}\mathscr{R}^{\chi,\mathfrak{q},\pm,\mathrm{j}}(\t,\x) \ &:= \ {\textstyle\int_{0}^{\t}}\mathbf{H}^{\N}(\s,\t(\N),\x)\{\N{\mathds{A}_{\mathbf{X}}^{\mathfrak{m}(\mathrm{j}),\pm}}[\mathds{R}^{\chi,\mathfrak{q},\pm,\mathrm{j}}\mathbf{Z};\s,\cdot(\s)]\}\d\s.
\end{align}
Let us also define $\mathscr{T}^{\pm,\mathrm{j}}$ as the following discrete spatial differential operator (for which we allow only two choices of signs):
\begin{align}
\mathscr{T}^{\pm,\mathrm{j}} \ := \ -\mathfrak{m}(\mathrm{j})^{-1}\{\grad^{\pm2\mathfrak{l}(\mathrm{j})}+\ldots+\grad^{\pm2\mathfrak{m}(\mathrm{j})\mathfrak{l}(\mathrm{j})}\}.
\end{align}
\end{definition}
\begin{lemma}\label{lemma:bg210}
 We have the following deterministic estimate for any $\mathrm{j}=1,\ldots,\mathrm{j}(\infty)$:
\begin{align}
\|\mathscr{R}^{\chi,\mathfrak{q},\pm,\mathrm{j}}-\sum_{0\leq\d\leq4}(\mathscr{T}^{\pm,\mathrm{j}})^{\d}\mathscr{A}^{\mathbf{X}}\mathscr{R}^{\chi,\mathfrak{q},\pm,\mathrm{j}}\|_{\t_{\mathrm{st}};\mathbb{T}(\N)} \ \lesssim \  \N^{-3\beta_{\mathrm{BG}}}\|\mathbf{Z}\|_{\t_{\mathrm{st}};\mathbb{T}(\N)}. \label{eq:bg210I}
\end{align}
\end{lemma}
Let us now briefly explain Definition \ref{definition:bg29} and Lemma \ref{lemma:bg210}. The only difference between $\mathscr{A}^{\mathbf{X}}\mathscr{R}^{\chi,\mathfrak{q},\pm,\mathrm{j}}$ from Definition \ref{definition:bg29} and $\mathscr{R}^{\chi,\mathfrak{q},\pm,\mathrm{j}}$ in Definition \ref{definition:bg26} is that the latter integrates $\N\mathds{R}^{\chi,\mathfrak{q},\pm,\mathrm{j}}\mathbf{Z}$ in space-time against the heat kernel, whereas the former first averages $\N\mathds{R}^{\chi,\mathfrak{q},\pm,\mathrm{j}}\mathbf{Z}$ in space with respect to length-scale $\mathfrak{m}(\mathrm{j})\mathfrak{l}(\mathrm{j})\lesssim\N^{3/4+\e}$ (for some small $\e>0$). (This length-scale comes from the observation that the aforementioned spatial average is of $\mathfrak{m}(\mathrm{j})$-many copies of $\N\mathds{R}^{\chi,\mathfrak{q},\pm,\mathrm{j}}\mathbf{Z}$, each copy being shifted by its support length $\mathfrak{l}(\mathrm{j})$; see Definition \ref{definition:bg26}.) Thus the cost in replacing $\mathscr{R}^{\chi,\mathfrak{q},\pm,\mathrm{j}}$ by $\mathscr{A}^{\mathbf{X}}\mathscr{R}^{\chi,\mathfrak{q},\pm,\mathrm{j}}$ is the error in freezing the $\mathbf{H}^{\N}$ heat kernel on blocks of length $\lesssim\N^{3/4+\e}$. This yields the $\d=1$ term (but without $\mathscr{A}^{\mathbf{X}}$) on the LHS of \eqref{eq:bg210I}. Now, we introduce $\mathscr{A}^{\mathbf{X}}$ into this $\d=1$ term, which yields the actual $\d=1$ term in \eqref{eq:bg210I}. By the same token, the cost is the $\d=2$ term (but without $\mathscr{A}^{\mathbf{X}}$) on the LHS of \eqref{eq:bg210I}. We then iterate until we derive an error given by $(\mathscr{T}^{\pm,\mathrm{j}})^{\d}\mathscr{R}^{\chi,\mathfrak{q},\pm,\mathrm{j}}$ for $\d=5$. By {\eqref{eq:hke2}} (regularity of the heat kernel), every power of $\mathscr{T}^{\pm,\mathrm{j}}$ introduces a factor of $\N^{-1}\N^{3/4}=\N^{-1/4}$ (up to small powers of $\N$). Five of these factors beat the $\N$ factor in $\mathscr{R}^{\chi,\mathfrak{q},\pm,\mathrm{j}}$ in Definition \ref{definition:bg26}, thereby concluding the proof of \eqref{eq:bg210I}.
\subsection{Replacement-by-time-average}
We now do the same but for time-averages. We clarify everything at the end of this step. (See after Proposition \ref{prop:bg212}.)
\begin{definition}\label{definition:bg211}
 First, we set some convenient notation. For any $\s\geq0$, define the following modified Gartner transform, in which we fix the time-parameter for the coupling constants $\lambda(\cdot)$ from Definitions \ref{definition:intro5} and \ref{definition:intro6} to be $\s$:
\begin{align}
\mathbf{G}^{\s}(\t,\x) \ := \ \exp[\lambda(\s)\mathbf{J}(\t,\x)-\lambda(\s){\textstyle\int_{0}^{\t}}\mathscr{R}(\r)\d\r]
\end{align}
Now, fix any $\tau>0$, integer ${\mathfrak{m}}\geq1$, and $\mathsf{F}(\t,\cdot)\in\mathscr{C}^{\infty}(\R^{\mathbb{T}(\N)})$. Define the following mixed space-time averages:
\begin{align}
&\mathds{A}^{{\mathfrak{m}},\tau,\pm}(\mathsf{F}\mathbf{Z};\s,\y(\s)) \nonumber\\
&:= \tau^{-1}{\textstyle\int_{0}^{\tau}}\{{\mathfrak{m}}^{-1}{\textstyle{\sum}_{\mathrm{j}=0}^{{\mathfrak{m}}-1}}\mathsf{F}(\s-\r,\mathbf{U}^{\s-\r,\y(\s-\r)\pm2\mathrm{j}\mathfrak{l}(\mathsf{F})})\mathbf{G}^{\s}(\s-\r,\y(\s-\r)\pm2\mathrm{j}\mathfrak{l}(\mathsf{F}))\}\d\r. \label{eq:bg211I}
\end{align}
Let us now specialize to the current situation. First we set $\tau(\mathrm{j}):= \N^{-3/2} \N^{3/4+\alpha(\mathrm{j})}$, where $\alpha(\mathrm{j})$ is from Definition \ref{definition:bg29}. Then, define the space-time integrated version
\begin{align}
\mathscr{A}^{\mathbf{X},\mathbf{T}}\mathscr{R}^{\chi,\mathfrak{q},\pm,\mathrm{j}}(\t,\x) \ := \ {\textstyle\int_{\tau(\mathrm{j})}^{\t}}\mathbf{H}^{\N}(\s,\t(\N),\x)\{\N\mathds{A}^{\mathfrak{m}(\mathrm{j}),\tau(\mathrm{j}),\pm}[\mathds{R}^{\chi,\mathfrak{q},\pm,\mathrm{j}}\mathbf{Z};\s,\cdot(\s)]\}\d\s. \label{eq:bg211II}
\end{align}
\end{definition}
A quick word about the use of $\mathbf{G}^{\s}$ instead of the original Gartner transform $\mathbf{G}$ in Definition \ref{definition:bg211}. We are essentially freezing the coupling constant $\lambda(\s)$ in \eqref{eq:bg211I}. We always take $\tau$ quite small (since we only look at local space-time scales), so by smoothness of the coupling constant, this freezing ultimately has a cost that can be controlled fairly directly. The reason why we need a frozen coupling constant is technical. (In a nutshell, \eqref{eq:hf} is easier to work with than the SDE for $\lambda(\t)\mathbf{J}(\t,\x)$, as far as our analysis of local averages is concerned. Indeed, by the product rule, the SDE for $\lambda(\t)\mathbf{J}(\t,\x)$ depends on $\mathbf{J}(\t,\x)$, whereas \eqref{eq:hf} depends only on \eqref{eq:glsde}. This means invariant measures for \eqref{eq:hf} are easy to compute, and this becomes hard once we add a $\mathbf{J}$-dependent drift.)
\begin{prop}\label{prop:bg212}
 With high probability, the following holds simultaneously for all $0\leq\d\leq4$ and $1\leq\mathrm{j}\leq\mathrm{j}(\infty)$:
\begin{align}
\|(\mathscr{T}^{\pm,\mathrm{j}})^{\d}\mathscr{A}^{\mathbf{X}}\mathscr{R}^{\chi,\mathfrak{q},\pm,\mathrm{j}}-(\mathscr{T}^{\pm,\mathrm{j}})^{\d}\mathscr{A}^{\mathbf{X},\mathbf{T}}\mathscr{R}^{\chi,\mathfrak{q},\pm,\mathrm{j}}\|_{\t_{\mathrm{st}};\mathbb{T}(\N)} \ \lesssim \  \N^{-3\beta_{\mathrm{BG}}}\|\mathbf{Z}\|_{\t_{\mathrm{st}};\mathbb{T}(\N)}. \label{eq:bg212I}
\end{align}
\end{prop}
The objects in Definition \ref{definition:bg211} average space-time shifts of $\mathsf{F}\mathbf{G}^{\s}$ on the space-time block $\{\s+[-\tau,0]\}\times\{\y(\cdot)\pm\llbracket0,2\mathfrak{l}\mathfrak{l}(\mathsf{F})\rrbracket\}$, where $\cdot$ denotes the time-variable for the Galilean shift in Definition \ref{definition:intro6}. The specialization in Definition \ref{definition:bg211} is for $\mathfrak{l}=\mathfrak{m}(\mathrm{j})$. So, the aforementioned block has \emph{spatial} scale $\lesssim\mathfrak{l}\mathfrak{l}(\mathsf{F})=\mathfrak{m}(\mathrm{j})\mathfrak{l}(\mathrm{j})$ (for $\mathsf{F}=\mathds{R}^{\chi,\mathfrak{q},\pm,\mathrm{j}}$, whose support length is $\mathfrak{l}(\mathrm{j})$; see Definition \ref{definition:bg26}). Observe that $\mathfrak{m}(\mathrm{j})\mathfrak{l}(\mathrm{j})\approx\N^{3/4}$ up to small powers of $\N$. Also, for time-scale $\tau(\mathrm{j})$, the relevant length-scale for \eqref{eq:hf}-\eqref{eq:glsde} is $\N^{3/2}\tau(\mathrm{j})\approx\N^{3/4}$. (Indeed, \eqref{eq:glsde} has an asymmetric drift with speed $\mathrm{O}(\N^{3/2})$.) Thus, our choices of $\tau(\mathrm{j})$ and $\mathfrak{m}(\mathrm{j})$ in Definitions \ref{definition:bg29}, \ref{definition:bg211} are ``naturally compatible". 

We now explain Proposition \ref{prop:bg212}. It asks to estimate the cost in replacing local spatial-average by its time-average on scale $\tau(\mathrm{j})$. (The change in integration-lower-limit from 0 to $\tau(\mathrm{j})$ in \eqref{eq:bg211II} is simply to guarantee we are not looking at \eqref{eq:glsde} for negative times.) Similar to Lemma \ref{lemma:bg210} (but in the time-direction), by time-regularity for the $\mathbf{H}^{\N}$ heat kernel (see {\eqref{eq:hke4}-\eqref{eq:hke5}}), this cost is $\lesssim\tau(\mathrm{j})\lesssim\N^{-3/4}$, at least modulo small powers of $\N$. This is not enough to beat the $\N$-factor in $\mathscr{A}^{\mathbf{X}}\mathscr{R}^{\chi,\mathfrak{q},\pm,\mathrm{j}}$! (See Definition \ref{definition:bg29}.) What saves us is the spatial-average introduced in Lemma \ref{lemma:bg210}. More precisely, by square-root cancellation and the fluctuating property of $\mathds{R}^{\chi,\mathfrak{q},\pm,\mathrm{j}}$, we expect that the average over spatial-shifts of $\mathds{R}^{\chi,\mathfrak{q},\pm,\mathrm{j}}\mathbf{G}^{\s}$ on a space-set of length $\N^{3/4}$ is $\lesssim\N^{-3/8+\e}$. (The fluctuating property of $\mathds{R}^{\chi,\mathfrak{q},\pm,\mathrm{j}}$ extends to $\mathds{R}^{\chi,\mathfrak{q},\pm,\mathrm{j}}\mathbf{G}^{\s}$ because of the ``decoupling" remark we made after {Definition \ref{definition:bg21}}.) Including this factor to the time-regularity upper bound $\N^{-3/4}$ beats the aforementioned $\N$ factor in $\mathscr{A}^{\mathbf{X}}\mathscr{R}^{\chi,\mathfrak{q},\pm,\mathrm{j}}$. This gives \eqref{eq:bg212I}. (Note that we only have to average on length $\N^{1/2+\e}$-blocks. Let us also note that technically, we average only $\mathfrak{m}(\mathrm{j})$-many shifts of $\mathds{R}^{\chi,\mathfrak{q},\pm,\mathrm{j}}\mathbf{G}^{\s}$, not $\mathfrak{m}(\mathrm{j})\mathfrak{l}(\mathrm{j})$-many. Therefore spatial-averaging, in principle, does not gain us $\mathfrak{m}(\mathrm{j})^{-1/2}\mathfrak{l}(\mathrm{j})^{-1/2}$ by square-root cancellations, but instead $\mathfrak{m}(\mathrm{j})^{-1/2}$. However, by construction in Definition \ref{definition:bg26}, we also gain $\mathds{R}^{\chi,\mathfrak{q},\pm,\mathrm{j}}\lesssim\mathfrak{l}(\mathrm{j})^{-3/2}$ modulo small powers of $\N$, so the above formal argument follows.)
\subsection{Finishing touches}
Having introduced mesoscopic space-time averaging in Lemma \ref{lemma:bg210} and Proposition \ref{prop:bg212}, we now leverage this to do homogenization. Ultimately, we get the following that we explain after proving Proposition \ref{prop:bg25}.
\begin{prop}\label{prop:bg213}
 With high probability, the following holds simultaneously for all $0\leq\d\leq4$ and $1\leq\mathrm{j}\leq\mathrm{j}(\infty)$:
\begin{align}
\|(\mathscr{T}^{\pm,\mathrm{j}})^{\d}\mathscr{A}^{\mathbf{X},\mathbf{T}}\mathscr{R}^{\chi,\mathfrak{q},\pm,\mathrm{j}}\|_{\t_{\mathrm{st}};\mathbb{T}(\N)} \ \lesssim \  \N^{-3\beta_{\mathrm{BG}}}\|\mathbf{Z}\|_{\t_{\mathrm{st}};\mathbb{T}(\N)}. \label{eq:bg213I}
\end{align}
\end{prop}
\begin{proof}[Proof of {Proposition \ref{prop:bg25}} assuming every result stated afterwards]
First, note that the intersection of $\mathrm{O}(1)$-many high probability events is also high probability. (Use union bound for their complements.) By this, the triangle inequality, and high probability estimates of Lemmas \ref{lemma:bg27}, \ref{lemma:bg210} and Propositions \ref{prop:bg212}, \ref{prop:bg213}, we obtain the following with high probability in which $\|\|=\|\|_{\t_{\mathrm{st}};\mathbb{T}(\N)}$:
\begin{align}
&{\textstyle\sup_{\mathrm{j}}}\|\mathscr{R}^{\mathfrak{q},\pm,\mathrm{j}}\| \nonumber\\
&\lesssim \ {\textstyle\sup_{\mathrm{j}}}\left(\|\mathscr{R}^{\mathfrak{q},\pm,\mathrm{j}}-\mathscr{R}^{\chi,\mathfrak{q},\pm,\mathrm{j}}\| + \|\mathscr{R}^{\chi,\mathfrak{q},\pm,\mathrm{j}}-{\textstyle\sum_{\d}}(\mathscr{T}^{\pm,\mathrm{j}})^{\d}\mathscr{A}^{\mathbf{X}}\mathscr{R}^{\chi,\mathfrak{q},\pm,\mathrm{j}}\|\right)\nonumber \\
&+ \ {\textstyle\sup_{\mathrm{d}}}{\textstyle\sup_{\mathrm{j}}}\left(\|(\mathscr{T}^{\pm,\mathrm{j}})^{\d}\mathscr{A}^{\mathbf{X}}\mathscr{R}^{\chi,\mathfrak{q},\pm,\mathrm{j}}-(\mathscr{T}^{\pm,\mathrm{j}})^{\d}\mathscr{A}^{\mathbf{X},\mathbf{T}}\mathscr{R}^{\chi,\mathfrak{q},\pm,\mathrm{j}}\|+\|(\mathscr{T}^{\pm,\mathrm{j}})^{\d}\mathscr{A}^{\mathbf{X},\mathbf{T}}\mathscr{R}^{\chi,\mathfrak{q},\pm,\mathrm{j}}\|\right) \nonumber\\
&\lesssim \  \N^{-3\beta_{\mathrm{BG}}}\|\mathbf{Z}\|. \nonumber
\end{align}
This is exactly the desired estimate \eqref{eq:bg25I}, so we are done.
\end{proof}
\subsubsection{Explanation of {Proposition \ref{prop:bg213}}: Multiscale II}\label{section:msII}
The term $\mathds{A}^{\mathfrak{l},\tau,\pm}(\mathsf{F}\mathbf{Z};\s,\y(\s))$ is an average of shifts of $\mathsf{F}\mathbf{G}$ over a space-time block of size $\tau\times\mathfrak{l}$. Thus, square-root cancellation implies this term is roughly $\lesssim\N^{-1}\tau^{-1/2}\mathfrak{l}^{-1/2}$. (The extra factor of $\N^{-1}$ comes from the $\N^{2}$ speed of \eqref{eq:glsde}. This implies faster homogenization, thus better estimates.) Choose $\tau=\tau(\mathrm{j})$ and $\mathfrak{l}=\mathfrak{m}(\mathrm{j})\mathfrak{l}(\mathrm{j})$, which is basically the case we are in (see the paragraph after Proposition \ref{prop:bg212}). In this case, we have $\N^{-1}\tau^{-1/2}\mathfrak{l}^{-1/2}\ll\N^{-1}$; see Definitions \ref{definition:bg29}, \ref{definition:bg211}. Thus, we expect $\mathds{A}^{\mathfrak{l},\tau,\pm}(\mathsf{F}\mathbf{Z};\s,\y(\s))$ to beat the $\N$-factor in \eqref{eq:bg211II}, roughly yielding the desired bound \eqref{eq:bg213I}. Unfortunately, this type of argument, thus far, can only be made rigorous for a stationary SDE. However \eqref{eq:glsde} does not even have a stationary measure! We discuss how to establish square-root cancellations in space-time using only a local equilibrium shortly. For now, let us take it for granted. Even in this case, to make rigorous the formal argument from the beginning of this paragraph, one needs to establish that the law of the process \eqref{eq:glsde} is sufficiently close to some local equilibrium on a space-time block of size $\tau(\mathrm{j})\times\mathfrak{m}(\mathrm{j})\mathfrak{l}(\mathrm{j})$. In our analysis, we compare to local equilibrium via the classical relative entropy inequality; see Appendix 1.8 of \cite{KL}. This lets us compare to local equilibrium on larger scales if we have stronger a priori estimates over the space-time average $\mathds{A}^{\mathfrak{l},\tau,\pm}(\mathsf{F}\mathbf{Z};\s,\y(\s))$ that hold on exponential large-deviations scale. But, based on what can be done for general stationary SDEs, square-root cancellation holds only in second moment. So we need to turn $\mathscr{L}^{2}$ into large deviations. To this end, {we present an algorithm}.
\begin{itemize}

\item Decompose the space-time average $\mathds{A}^{\mathfrak{l},\tau,\pm}(\mathsf{F}\mathbf{Z};\s,\y(\s))$ (for $\tau=\tau(\mathrm{j})$ and $\mathfrak{l}=\mathfrak{m}(\mathrm{j})\mathfrak{l}(\mathrm{j})$) as an average of $\mathds{A}^{\mathfrak{l},\tau,\pm}(\mathsf{F}\mathbf{Z};\s,\y(\s))$ for much smaller $\tau,\mathfrak{l}$. (In doing so we assume the space-time sets on which we average have mutually disjoint interiors.) Intuitively, this is morally the same as breaking a square with area equal to 1 square meter into 10000 mutually disjoint squares with area 1 square centimeter. Averaging over the big square is the same as averaging over each smaller square and then averaging the averages.
\item For each smaller-scale average $\mathds{A}^{\mathfrak{l},\tau,\pm}(\mathsf{F}\mathbf{Z};\s,\y(\s))$, artificially introduce a sub-optimal but nontrivial upper bound cutoff. As these cutoffs are sub-optimal, the square-root cancellations at local equilibrium (that we have taken for granted for now) show that the error term behind introducing said cutoffs is negligible in the large-$\N$ limit \emph{if} we can use said local equilibrium bounds. But because we are at much smaller space-time scales, we can directly compare to local equilibrium.
\item Next, we glue together smaller-scale averages $\mathds{A}^{\mathfrak{l},\tau,\pm}(\mathsf{F}\mathbf{Z};\s,\y(\s))$, \emph{now with a deterministic upper bound cutoff}, into averages on slightly larger space-time scales. In doing so, the a priori upper bounds glue as well. Thus the slightly larger-scale averages have the same deterministic upper bound cutoff. Next, artificially upgrade this cutoff into a stronger one that is still sub-optimal according to square-root cancellation heuristics \emph{on this larger space-time scale}. (Because we are now averaging over the larger space-time scales, we have better square-root cancellation.) The new cutoff is still sub-optimal. So, local equilibrium analysis shows the error behind introducing this new cutoff vanishes in the large-$\N$ limit. But, to use this, we need to compare the law of \eqref{eq:glsde} to local equilibrium on a larger scale. \emph{The point is that the deterministic upper bound obtained from the smaller-scale analysis lets us compare to local equilibrium on a larger scale by relative entropy inequality}.
\item Now, iterate. In particular, we again glue together the second-scale averages with improved a priori bounds into averages on slightly larger scales with the same improved a priori upper bounds. Now, with larger averaging scale we can further improve the a priori bounds. We continue until we reach the original scales $\tau=\tau(\mathrm{j})$ and $\mathfrak{l}=\mathfrak{m}(\mathrm{j})\mathfrak{l}(\mathrm{j})$. Again, the point is that improved a priori upper bounds let us compare to local equilibrium on larger scales. Averaging on said scales improves a priori bounds. 
\item A quick word to provide context: this strategy is exactly the one used in Section 7 of \cite{Y}, but with important modifications. First, our local averaging takes into account the Gartner transform $\mathbf{Z}$, and the space-time scales are chosen more optimally here than in \cite{Y}. Otherwise, the method works the same. (One technical difference is that in \cite{Y}, the largest space-time scale needed is $\tau\approx\N^{-1}$ and $\mathfrak{l}\approx\N^{1/2}$, but here we need $\tau\approx\N^{-3/4}$ and $\mathfrak{l}\approx\N^{2/3}$. The role of the entropy data assumption, i.e. Definition \ref{definition:entropydata}, is to access these slightly larger but still local scales.)
\end{itemize}
\subsubsection{Explanation of {Proposition \ref{prop:bg213}}: square-root cancellations using local equilibrium}\label{section:sqle}
The setting of this step is as follows. Take a (possibly) small-scale space-time average $\mathds{A}^{\mathfrak{l},\tau,\pm}(\mathsf{F}\mathbf{Z};\s,\y(\s))$; for example, take {$\mathsf{F}=\mathds{R}^{\chi,\mathfrak{q},\pm,\mathrm{j}}$}. This depends on the value of the $\mathbf{U}$ process on some space-time block of dimensions $\tau\times\mathfrak{l}$. If we take $\tau\approx\N^{-3/2}\mathfrak{l}+\N^{-2}\mathfrak{l}^{2}$, then we can replace $\mathbf{U}$, which satisfies the SDE \eqref{eq:glsde}, by a ``localized SDE". Precisely, take a neighborhood of the $\tau\times\mathfrak{l}$ block with dimensions slightly $\gg$ than $\tau\times\mathfrak{l}$. Because \eqref{eq:glsde} propagates (in space) according to speed $\N^{3/2}$ asymmetry and speed $\N^{2}$ symmetry, information outside this neighborhood does not affect $\mathds{A}^{\mathfrak{l},\tau,\pm}(\mathsf{F}\mathbf{Z};\s,\y(\s))$ with very high probability. This is by standard speed of propagation estimates for random walks. So, instead of averaging the value of $\mathsf{F}\mathbf{Z}$ along the path of \eqref{eq:glsde}, we instead average along the trajectory of the SDE \eqref{eq:glsde} but on the aforementioned $\tau\times\mathfrak{l}$ neighborhood (instead of $\mathbb{T}(\N)$) with periodic boundary. This is our localization step. In view of this, throughout the rest of this step, we will only ever refer to \eqref{eq:glsde} but on a spatial torus of size $\mathfrak{l}$ and for times until $\tau$. 

Assume that we have reduced to local equilibrium initial data for the local $\mathbf{U}$ dynamics. (Justifying the reduction is explained in the previous step.) In this case, where these $\mathbf{U}$ dynamics are time-homogeneous, they are also stationary; see Remark \ref{remark:intro12}. To derive square-root cancellations of space-time averages then amounts to the Kipnis-Varadhan estimate; see \cite{CLO,DGP}. This needs two main ingredients. First is knowing that the law of local $\mathbf{U}$ dynamics at any time is (sufficiently close to) a local equilibrium. (Again, in the time-homogeneous case, local equilibrium is invariant measure.) Second is a ``good" SDE representation for the time-reversal of the local $\mathbf{U}$ dynamics. For stationary SDEs, which we have in the time-homogeneous case, both are guaranteed. Let us describe these ingredients in the case of time-inhomogeneous dynamics with local equilibrium initial data.
\begin{enumerate}
\item (\emph{Closeness to local equilibrium}). The time $\t$ law of the local $\mathbf{U}$ dynamics with local equilibrium initial data can be computed via the Kolmogorov forward equation. If $\mathfrak{p}(\t)$ denotes the density of $\mathbf{U}$ at time $\t$ with respect to time $\t$ local equilibrium, then this Kolmogorov equation reads $\partial_{\t}\mathfrak{p}(\t)=\mathscr{L}(\t)\mathfrak{p}(\t)+\partial_{\t}\mathscr{V}(\t)\mathfrak{p}(\t)$. Here, $\mathscr{L}(\t)$ is a time-inhomogeneous Markov generator that, thus, has a maximum principle. Also, $\mathscr{V}(\t)$ is a ``total potential" or ``Hamiltonian". It is $\log$ of the Lebesgue density of local equilibrium measure at time $\t$, up to a sign; see Definition \ref{definition:intro5}. Indeed, this Kolmogorov PDE is classical for time-independent reference measure (for $\mathfrak{p}(\t)$). The additional multiplicative term comes from time-differentiating the local equilibrium itself. We must do this to compute total time-derivative of $\mathfrak{p}(\t)$. Because we always take $\tau\mathfrak{l}\leq\tau(\mathrm{j})\mathfrak{m}(\mathrm{j})\mathfrak{l}(\mathrm{j})\lesssim\N^{\e}$ (for very small $\e>0$), the maximum principle implies $\mathfrak{p}(\t)\lesssim\exp[\mathrm{O}(\N^{\e})]$. This is not so good. However, note that $\partial_{\t}\mathscr{V}(\t)$ is fluctuating (mean-zero with respect to time $\t$ local equilibrium). This follows from differentiating $\E^{\sigma,\t,\mathbb{I}}1=1$ (for any $\sigma,\t,\mathbb{I}$) in $\t$. Thus, we have cancellation in $\partial_{\t}\mathscr{V}(\t)$. This gives power-saving that beats $\N^{\e}$, ultimately giving an $\mathrm{O}(1)$ estimate for $\mathfrak{p}(\t)$. (When we make this precise, said $\mathrm{O}(1)$ estimate will be in some mixed $\mathscr{L}^{1},\mathscr{L}^{\infty}$ sense. Indeed, the cancellations in $\partial_{\t}\mathscr{V}(\t)$ are very high probability, not deterministic. Thus $\mathscr{L}^{\infty}$ alone is impossible. Also, $\mathscr{L}^{1}$ lets us rule out the very low probability event in which we do not see cancellations in $\partial_{\t}\mathscr{V}(\t)$.) 
\item (\emph{Time-reversed process}). Consider a localization of \eqref{eq:glsde} with local equilibrium initial data. If we reverse it, we get the same SDE, but with some adjustments. First, replace all times with their reversals and replace $\grad^{\mathrm{a}}\mathscr{U}'$ in \eqref{eq:glsde} by its negative. (This is just the statement that adjoint of the time-$\t$ generator of \eqref{eq:glsde} with respect to time-$\t$ local equilibrium is the original generator but the term giving $\grad^{\mathrm{a}}\mathscr{U}'$ picks up a sign because it is anti-symmetric. See also Section 2 of \cite{DGP}.) We must also add additional drift. Up to $\mathrm{O}(1)$ powers of the density $\mathfrak{p}(\t)$ and $\mathfrak{p}(\t)^{-1}$, this drift is the square root of the Dirichlet form of $\mathfrak{p}(\t)$. (The powers of $\mathfrak{p}(\t)$ and $\mathfrak{p}(\t)^{-1}$ are harmless by the above bullet point. Also, in the stationary time-homogeneous case, the Dirichlet form is zero as $\mathfrak{p}(\t)\equiv1$. Thus, what we said is clearly true for the stationary case.) Now, we want to remove this additional drift. (This would give us the situation one directly ends up in for the time-homogeneous stationary case.) Call the time-reversed process $\mathbf{V}$. Denote by $\mathbf{X}$ the time-reversal without additional Dirichlet form drift. By Girsanov, the law of the time-reversed \emph{process} $\mathbf{V}$ has Radon-Nikodym density with respect to the \emph{process} $\mathbf{X}$ given by an exponential martingale. Its log is equal to a martingale plus the time-integrated Dirichlet form. (Indeed, in general its log is a martingale plus a time-integrated square of the additional drift.) So, as estimating relative entropy requires estimating the log-Radon-Nikodym density, we can again use the relative entropy inequality to replace $\mathbf{V}$ by $\mathbf{X}$, if we can estimate the time-integrated Dirichlet form of $\mathfrak{p}(\t)$. But this can be done via standard diffusion theory, namely by energy dissipation (differentiate the $\mathscr{L}^{2}$-norm of $\mathfrak{p}(\t)$). This is another heat kernel estimate for the local $\mathbf{U}$ dynamics. It is treated like in the first bullet point. We conclude this point by clarifying one interesting technical detail. Given sufficiently small times $\tau$ (so small lengths $\mathfrak{l}$), estimating the time-integrated Dirichlet form of $\mathfrak{p}(\t)$ is rather straightforward. (Indeed, the integration domain in said Dirichlet form is small.) For larger times $\tau$, as in the previous point, we need sharper a priori estimates for time-averages $\mathds{A}^{\mathfrak{l},\tau,\pm}(\mathsf{F}\mathbf{Z};\s,\y(\s))$ if we want to use the entropy inequality as explained in this paragraph. But these come from the bootstrapping algorithm from the previous step! In particular, the ``Multiscale II" and the current step work in unison.
\end{enumerate}
Of course, we must make all these estimates precise. The main obstruction to this end is the fact that local space-time averages in Definition \ref{definition:bg211} are of functions of \eqref{eq:hf}-\eqref{eq:glsde}, not just \eqref{eq:glsde}. First, observe that \eqref{eq:hf} is an SDE that is (spatially) local in the solution to \eqref{eq:glsde}. Thus, the speed of propagation in the first paragraph of Section \ref{section:sqle} lets us replace \eqref{eq:hf} by a ``localized" version as well. Second, observe that \eqref{eq:hf} (and its localizations given by replacing $\mathbb{T}(\N)$ with a discrete interval $\mathbb{K}\subseteq\mathbb{T}(\N)$) at any fixed spatial point $\x$ has invariant measure equal to Lebesgue measure on $\R$. (Indeed, \eqref{eq:hf} is an $\R$-valued SDE with coefficients that are independent of the solution.) So ultimately, we only have to adjust the above bullet points by replacing local equilibrium measures therein by its tensor product with $\mathrm{Leb}(\R)$. (There is the technical issue that $\mathrm{Leb}(\R)$ is not a probability measure. Thus we have to introduce some type of a priori cutoff to the SDE \eqref{eq:hf}. This will ultimately come from working before time $\t_{\mathrm{st}}$, where we have a priori estimates for the Gartner transform $\mathbf{Z}$.)
\subsection{What is left}
For the non-appendix sections of this paper, all we have to prove are the results in this section (Lemmas \ref{lemma:bg27}, \ref{lemma:bg210}, Propositions \ref{prop:bg212}, \ref{prop:bg213}) and Corollary \ref{corollary:kpz}. As we explained in Section \ref{subsection:reading}, everything from here until the proof of Corollary \ref{corollary:kpz} executes the power-counting behind the proof heuristics given in this section. In particular, what follows is fairly technical, though we try to explain every step and notation for clarity's purpose. But, every idea in what follows has already been explained in this section or an earlier one (so the role of what follows is just to make this paper into rigorous mathematics).
%
%
%
\section{Proofs of Lemmas \ref{lemma:bg27}, \ref{lemma:bg210}}\label{section:proofsstart}
The point of this section is just to turn the proof heuristics given after the statements of Lemmas \ref{lemma:bg27}, \ref{lemma:bg210} into mathematics. (There are no additional subtleties.) Since these arguments are fairly elementary and somewhat besides the most interesting result {(Proposition \ref{prop:bg213})} from the previous section, especially given their explanations in the previous section, the reader is invited to skip this section in a first reading.
\subsection{Proof of Lemma \ref{lemma:bg27}}
Recall notation from Definitions \ref{definition:bg24} and \ref{definition:bg26}. Define $\mathds{R}^{\dagger,\mathfrak{q},\pm,\mathrm{j}}:=\mathds{R}^{\mathfrak{q},\pm,\mathrm{j}}-\mathds{R}^{\chi,\mathfrak{q},\pm,\mathrm{j}}$. By contractivity of $\mathbf{H}^{\N}$ (see {\eqref{eq:hke3} with $\mathrm{m}=0$}) and triangle inequality, we have the following deterministic bound, in which $\|\|=\|\|_{\t_{\mathrm{st}};\mathbb{T}(\N)}$:
\begin{align}
\|\mathscr{R}^{\mathfrak{q},\pm,\mathrm{j}}-\mathscr{R}^{\chi,\mathfrak{q},\pm,\mathrm{j}}\| \ &\leq \ \|{\textstyle\int_{0}^{\t}}\mathbf{H}^{\N}(\s,\t(\N),\x)\{\N|\mathds{R}^{\dagger,\mathfrak{q},\pm,\mathrm{j}}(\s,\cdot(\s))|\mathbf{Z}(\s,\cdot)\}\d\s\| \nonumber\\
&\lesssim \ \N\|\mathds{R}^{\dagger,\mathfrak{q},\pm,\mathrm{j}}\mathbf{Z}\| \ \lesssim \ \|\N\mathds{R}^{\dagger,\mathfrak{q},\pm,\mathrm{j}}\|\|\mathbf{Z}\|. \label{eq:bg27I1}
\end{align}
(In \eqref{eq:bg27I1}, we use $\t_{\mathrm{st}}\leq1$; see Definitions \ref{definition:method8}, \ref{definition:reg}. Also, the second norm in \eqref{eq:bg27I1} is with respect to $\t,\x$.) Now, it suffices to show
\begin{align}
\|\mathds{R}^{\dagger,\mathfrak{q},\pm,\mathrm{j}}\|_{\t_{\mathrm{st}};\mathbb{T}(\N)} \ \lesssim \  \N^{-10} \label{eq:bg27I2}
\end{align}
with high probability. (It actually holds almost surely.) Take $\t\leq\t_{\mathrm{st}}$. We compute $\mathds{R}^{\dagger,\mathfrak{q},\pm,\mathrm{j}}$ by recalling Definition \ref{definition:bg26}. We claim the following holds, where we implicitly evaluate $\E^{\mathfrak{l}(\mathrm{j}),\pm}$ at $(\t,\x)$. (We do this implicit evaluating at $(\t,\x)$ throughout this proof.)
\begin{align}
&\mathds{R}^{\dagger,\mathfrak{q},\pm,\mathrm{j}}(\t,\x) \nonumber\\
&= \ \mathds{R}^{\mathfrak{q},\pm,\mathrm{j}}(\t,\x)\times\{1-\chi[\mathds{R}^{\mathfrak{q},\pm,\mathrm{j}}(\t,\x);\upsilon_{\mathrm{j}-1}]\}+\E^{\mathfrak{l}(\mathrm{j}),\pm}\{\mathds{R}^{\mathfrak{q},\pm,\mathrm{j}}(\t,\mathbf{U})\cdot\chi[\mathds{R}^{\mathfrak{q},\pm,\mathrm{j}}(\t,\mathbf{U});\upsilon_{\mathrm{j}-1}]\} \\
&= \ \E^{\mathfrak{l}(\mathrm{j}),\pm}\{\mathds{R}^{\mathfrak{q},\pm,\mathrm{j}}(\t,\mathbf{U})\cdot\chi[\mathds{R}^{\mathfrak{q},\pm,\mathrm{j}}(\t,\mathbf{U});\upsilon_{\mathrm{j}-1}]\}. \label{eq:bg27I3}
\end{align}
Indeed, the first term on the RHS of the first identity is supported on the event $|\mathds{R}^{\mathfrak{q},\pm,\mathrm{j}}(\t,\x)|\geq\upsilon_{\mathrm{j}-1}=\N^{20\gamma_{\mathrm{reg}}}\mathfrak{l}(\mathrm{j}-1)^{-1/2}$. However, by construction of $\mathds{R}^{\mathfrak{q},\pm,\mathrm{j}}(\t,\x)$ as a difference of $\E^{\mathfrak{l}(\mathrm{k}),\pm}$-terms (for $\mathrm{k}=\mathrm{j}-1,\mathrm{j}$) and by Lemma \ref{lemma:bg23}, we deduce that this event is empty because $\t\leq\t_{\mathrm{st}}\leq\t_{\mathrm{reg}}$. So \eqref{eq:bg27I3} is true. Now, by Remark \ref{remark:intro14} the charge density in $\E^{\mathfrak{l}(\mathrm{j}),\pm}$ in \eqref{eq:bg27I3} is $\lesssim\N^{\gamma_{\mathrm{reg}}}\mathfrak{l}(\mathrm{j})^{-1/2}$ in absolute value. Let $\mathcal{E}$ be the event (in the probability space that $\E^{\mathfrak{l}(\mathrm{j}),\pm}$ is defined on) where the charge density on $\x+\mathbb{I}(\mathfrak{l}(\mathrm{j}-1),\pm)$ is $\lesssim\N^{2\gamma_{\mathrm{reg}}}\mathfrak{l}(\mathrm{j}-1)^{-1/2}$ in absolute value. (In particular, $\mathbf{1}[\mathcal{E}]$ and $\mathbf{1}[\mathcal{E}^{\mathrm{C}}]$ will never be evaluated at $\mathbf{U}^{\t,\cdot}$. Like with the clarification after \eqref{eq:bg26I}, it will always be a function of the expectation ``dummy" variable $\mathbf{U}$ in $\E^{\mathfrak{l}(\mathrm{j}),\pm}$.) Now, we claim that
\begin{align}
\eqref{eq:bg27I3} \ &= \ \E^{\mathfrak{l}(\mathrm{j}),\pm}\{\mathds{R}^{\mathfrak{q},\pm,\mathrm{j}}(\t,\mathbf{U})\}+\E^{\mathfrak{l}(\mathrm{j}),\pm}\{\mathds{R}^{\mathfrak{q},\pm,\mathrm{j}}(\t,\mathbf{U})\cdot(1-\chi[\mathds{R}^{\mathfrak{q},\pm,\mathrm{j}}(\t,\mathbf{U});\upsilon_{\mathrm{j}-1}])\}\label{eq:bg27I4}\\
&= \ \E^{\mathfrak{l}(\mathrm{j}),\pm}\{\mathbf{1}[\mathcal{E}]\mathds{R}^{\mathfrak{q},\pm,\mathrm{j}}(\t,\mathbf{U})\cdot(1-\chi[\mathds{R}^{\mathfrak{q},\pm,\mathrm{j}}(\t,\mathbf{U});\upsilon_{\mathrm{j}-1}])\}\nonumber\\
&\quad+\E^{\mathfrak{l}(\mathrm{j}),\pm}\{\mathbf{1}[\mathcal{E}^{\mathrm{C}}]\mathds{R}^{\mathfrak{q},\pm,\mathrm{j}}(\t,\mathbf{U})\cdot(1-\chi[\mathds{R}^{\mathfrak{q},\pm,\mathrm{j}}(\t,\mathbf{U});\upsilon_{\mathrm{j}-1}])\}. \nonumber
\end{align}
Indeed, the first term on the RHS of \eqref{eq:bg27I4} is zero by Lemma \ref{lemma:vanishcanonical} as noted after Definition \ref{definition:bg26}. We claim $\E^{\mathfrak{l}(\mathrm{j}),\pm}\mathbf{1}[\mathcal{E}^{\mathrm{C}}]\lesssim\N^{-{\mathrm{D}}}$ {for some large but fixed $\mathrm{D}>0$}; we take this for granted for now. By this and Cauchy-Schwarz and $|\chi|\lesssim1$, we can control the second term in the last line above:
\begin{align}
|\E^{\mathfrak{l}(\mathrm{j}),\pm}\{\mathbf{1}[\mathcal{E}^{\mathrm{C}}]\mathds{R}^{\mathfrak{q},\pm,\mathrm{j}}(\t,\mathbf{U})\cdot(1-\chi[\mathds{R}^{\mathfrak{q},\pm,\mathrm{j}}(\t,\mathbf{U});\upsilon_{\mathrm{j}-1}])\}| \ \lesssim \ \N^{-100}\{\E^{\mathfrak{l}(\mathrm{j}),\pm}\mathds{R}^{\mathfrak{q},\pm,\mathrm{j}}(\t,\mathbf{U})^{2}\}^{\frac12}.
\end{align}
Recall from Definition \ref{definition:bg24} that $\mathds{R}^{\mathfrak{q},\pm,\mathrm{j}}(\t,\mathbf{U})$ is the difference of canonical ensemble expectations of $\mathfrak{q}$. In this context, again, the canonical measure expectations are not with respect to charge densities determined by the process at time $\t$, namely $\mathbf{U}^{\t,\cdot}$. Rather, these charge densities are sampled according to $\mathbf{U}$, the expectation dummy variable in $\E^{\mathfrak{l}(\mathrm{j}),\pm}$. Therefore, via the tower property of expectation and Cauchy-Schwarz, when computing the expectation from the RHS of the previous display, we can collapse all iterated expectations into a single $\E^{\mathfrak{l}(\mathrm{j}),\pm}$ and get the following, which is deterministic like every other display in this proof:{
\begin{align}
&\E^{\mathfrak{l}(\mathrm{j}),\pm}\mathds{R}^{\mathfrak{q},\pm,\mathrm{j}}(\t,\mathbf{U})^{2} \\\
&\lesssim \ \E^{\mathfrak{l}(\mathrm{j}),\pm}\{|\mathfrak{q}(\t,\cdot)|^{2};\t,\x\}\nonumber\\
&= \ \E^{\mathfrak{l}(\mathrm{j}),\pm}\{[|\mathfrak{q}(\t,\cdot)|^{2}-\E^{0,\t}|\mathfrak{q}(\t,\cdot)|^{2}];\t,\x\} + \E^{0,\t}|\mathfrak{q}(\t,\cdot)|^{2} \ = \ \mathrm{O}(1) \quad\mathrm{for}\quad \t\leq\t_{\mathrm{st}}\leq\t_{\mathrm{reg}},\nonumber
\end{align}
where the last bound in the previous display can be justified as follows. Since $|\mathfrak{q}(\t,\mathbf{u})|\lesssim1+|\mathbf{u}|^{C}$ for some $C=\mathrm{O}(1)$, and because the probability measure in $\E^{0,\t}$ is sub-Gaussian (see Assumption \ref{ass:intro8}), we deduce that $\E^{0,\t}|\mathfrak{q}(\t,\cdot)|^{2}\lesssim1$. On the other hand, we note that $|\mathfrak{q}(\t,\cdot)|^{2}-\E^{0,\t}|\mathfrak{q}(\t,\cdot)|^{2}$ belongs to the $\mathrm{CT}$ set by construction, so we can apply Lemma \ref{lemma:ee1} and \eqref{eq:eep1} (the latter of which holds for any $\mathfrak{a}(\t,\mathbf{u})$ with polynomial growth in $\mathbf{u}$). This implies that $\E^{\mathfrak{l}(\mathrm{j}),\pm}\{[|\mathfrak{q}(\t,\cdot)|^{2}-\E^{0,\t}|\mathfrak{q}(\t,\cdot)|^{2}];\t,\x\}=\mathrm{O}(\N^{10\gamma_{\mathrm{reg}}}|\mathfrak{l}(\mathrm{j})|^{-1/2})=\mathrm{O}(1)$, where the last bound is because $\mathfrak{l}(\mathrm{j})=N^{\delta_{\mathrm{BG}}}$ with $\delta_{\mathrm{BG}}$ equal to a big constant times $\gamma_{\mathrm{reg}}$. (See Definition \ref{definition:method8} and Definition \ref{definition:bg24}.) Therefore, the last bound in the previous display holds.} From the previous two displays, we get
\begin{align}
|\E^{\mathfrak{l}(\mathrm{j}),\pm}\{\mathbf{1}[\mathcal{E}^{\mathrm{C}}]\mathds{R}^{\mathfrak{q},\pm,\mathrm{j}}(\t,\mathbf{U})\cdot(1-\chi[\mathds{R}^{\mathfrak{q},\pm,\mathrm{j}}(\t,\mathbf{U});\upsilon_{\mathrm{j}-1}])\}| \ \lesssim \ \N^{-100}.
\end{align}
{Let us study the first term in the line after \eqref{eq:bg27I4}. For this, we clarify two points. First, whatever is inside the expectation $\E^{\mathfrak{l}(\mathrm{j}),\pm}$ is no longer a function of the process $\mathbf{U}^{\t,\cdot}$ at time $\t$. Instead, they are each functions of the expectation (dummy) variable in $\E^{\mathfrak{l}(\mathrm{j}),\pm}$. In particular, the first term in the line after \eqref{eq:bg27I4} is a function of $\mathbf{U}^{\t,\cdot}$ only through the average of $\mathbf{U}^{\t,\z}$ for $\z\in\x+\mathbb{I}(\mathfrak{l}(\mathrm{j}),\pm)$ that determines only the charge density for $\E^{\mathfrak{l}(\mathrm{j}),\pm}$; see Definition \ref{definition:bg21}. Second, the factor $1-\chi[\mathds{R}^{\mathfrak{q},\pm,\mathrm{j}}(\t,\mathbf{U});\upsilon_{\mathrm{j}-1}]$, which, again, is a function of the expectation (dummy) variable, is supported on the event where $|\mathds{R}^{\mathfrak{q},\pm,\mathrm{j}}(\t,\mathbf{U})|\geq\upsilon_{\mathrm{j}-1}=\N^{20\gamma_{\mathrm{reg}}}\mathfrak{l}(\mathrm{j}-1)^{-1/2}$; see Definition \ref{definition:bg26}. Also, recall $\mathcal{E}$ from before \eqref{eq:bg27I4}. On $\mathcal{E}$, we may follow the proof of Lemma \ref{lemma:ee1} to get $1-\chi[\mathds{R}^{\mathfrak{q},\pm,\mathrm{j}}(\t,\mathbf{U});\upsilon_{\mathrm{j}-1}]=0$. We clarify that by $1-\chi[\mathds{R}^{\mathfrak{q},\pm,\mathrm{j}}(\t,\mathbf{U});\upsilon_{\mathrm{j}-1}]$ we no longer mean a function of the process $\mathbf{U}^{\t,\cdot}$ but of the expectation variable in $\E^{\mathfrak{l}(\mathrm{j}),\pm}$. Let us now clarify what we mean by following the proof of Lemma \ref{lemma:ee1} to establish $1-\chi[\mathds{R}^{\mathfrak{q},\pm,\mathrm{j}}(\t,\mathbf{U});\upsilon_{\mathrm{j}-1}]=0$. We recall in Definition \ref{definition:bg24} that $\mathds{R}^{\mathfrak{q},\pm,\mathrm{j}}$ is a difference between two canonical expectations of $\mathfrak{q}$. One is on the set $\x+\mathbb{I}(\mathfrak{l}(\mathrm{j}-1),\pm)$, with charge density sampled by the measure in $\E^{\mathfrak{l}(\mathrm{j}),\pm}$, and another is on the set $\x+\mathbb{I}(\mathfrak{l}(\mathrm{j}),\pm)$. On the event $\mathcal{E}$, the charge densities are $\lesssim\N^{2\gamma_{\mathrm{reg}}}\mathfrak{l}(\mathrm{j}-1)^{-1/2}$ and $\lesssim\N^{2\gamma_{\mathrm{reg}}}\mathfrak{l}(\mathrm{j})^{-1/2}$ in absolute value, respectively. (The latter bound follows from Remark \ref{remark:intro14}; see the paragraph after \eqref{eq:bg27I3}.) Now use the proof of Lemma \ref{lemma:ee1} but replace $\sigma(\t,\x;\mathfrak{l},\pm)$ by charge densities on $\x+\mathbb{I}(\mathfrak{l}(\mathrm{j}-1),\pm)$ and $\x+\mathbb{I}(\mathfrak{l}(\mathrm{j}),\pm)$. This means both canonical ensemble expectations in $\mathds{R}^{\mathfrak{q},\pm,\mathrm{j}}$ are $\lesssim\N^{13\gamma_{\mathrm{reg}}}\mathfrak{l}(\mathrm{j}-1)^{-1/2}$. (Technically, Lemma \ref{lemma:ee1} requires a bound on charge densities of $\lesssim\N^{\gamma_{\mathrm{reg}}}|\mathfrak{l}(\mathrm{j}-1)|^{-1/2}$, not the bound of $\N^{2\gamma_{\mathrm{reg}}}|\mathfrak{l}(\mathrm{j}-1)|^{-1/2}$ that we get from $\mathcal{E}$. But the dependence of the bound in Lemma \ref{lemma:ee1} on the charge density is at worst cubic. Thus, the slightly worse charge density estimate that we have on $\mathcal{E}$ forces us to multiply the $\mathrm{RHS}\eqref{eq:ee1I}$ by $\N^{3\gamma_{\mathrm{reg}}}$. This explains $\N^{13\gamma_{\mathrm{reg}}}$, instead of $\N^{10\gamma_{\mathrm{reg}}}$ in \eqref{eq:ee1I}.) In any case, $|\mathds{R}^{\mathfrak{q},\pm,\mathrm{j}}(\t,\x)|\lesssim\N^{13\gamma_{\mathrm{reg}}}\mathfrak{l}(\mathrm{j}-1)^{-1/2}\leq\N^{20\gamma_{\mathrm{reg}}}\mathfrak{l}(\mathrm{j}-1)^{-1/2}=\upsilon_{\mathrm{j}-1}$ holds on $\mathcal{E}$. Thus $1-\chi[\mathds{R}^{\mathfrak{q},\pm,\mathrm{j}}(\t,\mathbf{U});\upsilon_{\mathrm{j}-1}]=0$ on $\mathcal{E}$; see Definition \ref{definition:bg26} (in particular, the definition of $\chi$ therein). Ultimately, we have}
\begin{align}
\E^{\mathfrak{l}(\mathrm{j}),\pm}\{\mathbf{1}[\mathcal{E}]\mathds{R}^{\mathfrak{q},\pm,\mathrm{j}}(\t,\mathbf{U})\cdot(1-\chi[\mathds{R}^{\mathfrak{q},\pm,\mathrm{j}}(\t,\mathbf{U});\upsilon_{\mathrm{j}-1}])\} \ = \ 0. \label{eq:bg27I6}
\end{align}
Combining the previous six displays proves \eqref{eq:bg27I2}. So, we are left to show that $\mathcal{E}$ (from before \eqref{eq:bg27I4}) is very high probability (with respect to the law of $\E^{\mathfrak{l}(\mathrm{j}),\pm}$), {i.e. that the following holds with very high probability with respect to the law of $\E^{\mathfrak{l}(\mathrm{j}),\pm}$:
\begin{align*}
\mathfrak{l}(\mathrm{j}-1)^{-1}\left|\sum_{\z\in\x+\mathds{I}(\mathfrak{l}(\mathrm{j}-1),\pm)}\mathbf{U}(\z)\right|\ \lesssim\ \N^{2\gamma_{\mathrm{reg}}}\mathfrak{l}(\mathrm{j}-1)^{-\frac12}.
\end{align*}
To this end, we first recall that $a:=\mathfrak{l}(\mathrm{j})^{-1}\sum_{\z\in\x+\mathds{I}(\mathfrak{l}(\mathrm{j}),\pm)}\mathbf{U}(\z)$ (we use this notation for this proof only) satisfies $|a|\lesssim\N^{2\gamma_{\mathrm{reg}}}\mathfrak{l}(\mathrm{j})^{-1/2}$; see the paragraph before \eqref{eq:bg27I4}. Next, we use this bound to obtain
\begin{align*}
\mathfrak{l}(\mathrm{j}-1)^{-1}\left|\sum_{\z\in\x+\mathds{I}(\mathfrak{l}(\mathrm{j}-1),\pm)}\mathbf{U}(\z)\right|&\lesssim\mathfrak{l}(\mathrm{j}-1)^{-1}\left|\sum_{\z\in\x+\mathds{I}(\mathfrak{l}(\mathrm{j}-1),\pm)}(\mathbf{U}(\z)-a)\right|+\N^{2\gamma_{\mathrm{reg}}}\mathfrak{l}(\mathrm{j})^{-1/2}.
\end{align*}
If $\mathbf{U}(\z)$ is distributed via the law of $\E^{\mathfrak{l}(\mathrm{j}),\pm}$, then $\mathbf{U}(\z)-a$ are the steps in a random walk bridge with zero drift of length $\mathfrak{l}(\mathrm{j}-1)$. Moreover, the law of $\mathbf{U}(\z)-a$ is sub-Gaussian; see Assumption \ref{ass:intro8}. In particular, by standard random walk bridge theory, we know that $\sum_{\z\in\x+\mathds{I}(\mathfrak{l}(\mathrm{j}-1),\pm)}\mathbf{U}(\z)$ is, up to an error of size $\lesssim|a|\lesssim\N^{2\gamma_{\mathrm{reg}}}\mathfrak{l}(\mathrm{j})^{-1/2}$, equal to a sum of $\mathfrak{l}(\mathrm{j}-1)$-many sub-Gaussian martingale increments. (The aforementioned ``standard random walk bridge theory" is just a random walk version of the fact that a Brownian bridge $b_{t}$ on $t\in[0,1]$ has the representation $w_{t}-tw_{1}$, where $w_{t}$ is a standard Brownian motion.) The desired high probability bound now follows by applying Azuma's martingale inequality to the first term on the RHS of the previous display. As noted earlier, we are now done.} \qed
\subsection{Proof of Lemma \ref{lemma:bg210}}
First, {we make} some observations. Note {that} ${\mathds{A}_{\mathbf{X}}^{\mathfrak{m}(\mathrm{j}),\pm}}\{\mathds{R}^{\chi,\mathfrak{q},\pm,\mathrm{j}}\mathbf{Z};\s,\y(\s)\}$ is the average over $\mathrm{k}=0,\ldots,\mathfrak{m}(\mathrm{j})-1$ of $\mathds{R}^{\chi,\mathfrak{q},\pm,\mathrm{j}}[\s,\y(\s)\pm2\mathrm{k}\mathfrak{l}(\mathds{R}^{\chi,\mathfrak{q},\pm,\mathrm{j}})]\mathbf{Z}[\s,\y\pm2\mathrm{k}\mathfrak{l}(\mathds{R}^{\chi,\mathfrak{q},\pm,\mathrm{j}})]$, where $\mathfrak{l}(\mathds{R}^{\chi,\mathfrak{q},\pm,\mathrm{j}})$ is the support length of $\mathds{R}^{\chi,\mathfrak{q},\pm,\mathrm{j}}$. By Definition \ref{definition:bg26} (see also the paragraph following it), we get $\mathfrak{l}(\mathds{R}^{\chi,\mathfrak{q},\pm,\mathrm{j}})=\mathfrak{l}(\mathrm{j})$; see Definition \ref{definition:bg24} for $\mathfrak{l}(\mathrm{j})$. So, ${\mathds{A}_{\mathbf{X}}^{\mathfrak{m}(\mathrm{j}),\pm}}\{\mathds{R}^{\chi,\mathfrak{q},\pm,\mathrm{j}}\mathbf{Z};\s,\y(\s)\}$ is the average over $\mathrm{k}=0,\ldots,\mathfrak{m}(\mathrm{j})-1$ of gradients of $\mathds{R}^{\chi,\mathfrak{q},\pm,\mathrm{j}}[\s,\y(\s)]\mathbf{Z}[\s,\y]$ (with length-scale $\pm2\mathrm{k}\mathfrak{l}(\mathrm{j})$). By Definition \ref{definition:bg29},
\begin{align}
\mathds{R}^{\chi,\mathfrak{q},\pm,\mathrm{j}}(\s,\y(\s))\mathbf{Z}(\s,\y)-{\mathds{A}_{\mathbf{X}}^{\mathfrak{m}(\mathrm{j}),\pm}}\{\mathds{R}^{\chi,\mathfrak{q},\pm,\mathrm{j}}\mathbf{Z};\s,\y(\s)\} \ = \ \mathscr{T}^{\pm,\mathrm{j}}\{\mathds{R}^{\chi,\mathfrak{q},\pm,\mathrm{j}}(\s,\y(\s))\mathbf{Z}(\s,\y)\}. \label{eq:bg2101}
\end{align}
We now observe that convolution on $\mathbb{T}(\N)$, viewed as an operator on functions $\mathbb{T}(\N)\to\R$, commutes with any discrete gradient. Because $\mathscr{T}^{\pm,\mathrm{j}}$ is an average of discrete gradients, we deduce $\mathscr{T}^{\pm,\mathrm{j}}$ commutes with convolution. Using \eqref{eq:bg2101} with this, we get the following. Take $\Phi,\Psi:[0,\infty)\times\mathbb{T}(\N)\to\R$. Let {$\Phi(\x)\star_{\y}\Psi(\y):=\sum_{\y\in\mathbb{T}(\N)}\Phi(\x-\y)\Psi(\y)$} be convolution of $\Phi,\Psi$ on $\mathbb{T}(\N)$ evaluated at $\x$. Now, we get
\begin{align}
&{\textstyle\int_{0}^{\t}}\Phi(\s,\x)\star_{\y}\N\{\mathds{R}^{\chi,\mathfrak{q},\pm,\mathrm{j}}(\s,\y(\s))\mathbf{Z}(\s,\y)-{\mathds{A}_{\mathbf{X}}^{\mathfrak{m}(\mathrm{j}),\pm}}[\mathds{R}^{\chi,\mathfrak{q},\pm,\mathrm{j}}\mathbf{Z};\s,\y(\s)]\}\d\s \label{eq:bg2102a}\\
&= \ {\textstyle\int_{0}^{\t}}\Phi(\s,\x)\star_{\y}\mathscr{T}^{\pm,\mathrm{j}}\{\N\mathds{R}^{\chi,\mathfrak{q},\pm,\mathrm{j}}(\s,\y(\s))\mathbf{Z}(\s,\y)\}\d\s \\
&= \ {\textstyle\int_{0}^{\t}}[\mathscr{T}^{\pm,\mathrm{j}}\Phi(\s,\x)]\star_{\y}[\N\mathds{R}^{\chi,\mathfrak{q},\pm,\mathrm{j}}(\s,\y(\s))\mathbf{Z}(\s,\y)]\d\s.\label{eq:bg2102b}
\end{align}
(Again, the first identity in \eqref{eq:bg2102a}-\eqref{eq:bg2102b} follows from \eqref{eq:bg2101}, and the second follows from the commutativity of convolution and discrete gradient.) Recall $\mathscr{R}^{\chi,\mathfrak{q},\pm,\mathrm{j}}$ and $\mathscr{A}^{\mathbf{X}}\mathscr{R}^{\chi,\mathfrak{q},\pm,\mathrm{j}}$ from Definitions \ref{definition:bg26} and \ref{definition:bg29}. We now apply \eqref{eq:bg2102a}-\eqref{eq:bg2102b} iteratively, so that $\mathscr{R}^{\chi,\mathfrak{q},\pm,\mathrm{j}}$ equals
\begin{align}
&\mathscr{A}^{\mathbf{X}}\mathscr{R}^{\chi,\mathfrak{q},\pm,\mathrm{j}} + \mathscr{T}^{\pm,\mathrm{j}}\mathscr{R}^{\chi,\mathfrak{q},\pm,\mathrm{j}} \ = \ \mathscr{A}^{\mathbf{X}}\mathscr{R}^{\chi,\mathfrak{q},\pm,\mathrm{j}}+\mathscr{T}^{\pm,\mathrm{j}}\mathscr{A}^{\mathbf{X}}\mathscr{R}^{\chi,\mathfrak{q},\pm,\mathrm{j}} + (\mathscr{T}^{\pm,\mathrm{j}})^{2}\mathscr{R}^{\chi,\mathfrak{q},\pm,\mathrm{j}} \\
&= \ \ldots \ = \ \mathscr{A}^{\mathbf{X}}\mathscr{R}^{\chi,\mathfrak{q},\pm,\mathrm{j}} + \mathscr{T}^{\pm,\mathrm{j}}\mathscr{A}^{\mathbf{X}}\mathscr{R}^{\chi,\mathfrak{q},\pm,\mathrm{j}}+(\mathscr{T}^{\pm,\mathrm{j}})^{2}\mathscr{A}^{\mathbf{X}}\mathscr{R}^{\chi,\mathfrak{q},\pm,\mathrm{j}} +\ldots+ (\mathscr{T}^{\pm,\mathrm{j}})^{5}\mathscr{R}^{\chi,\mathfrak{q},\pm,\mathrm{j}}. \label{eq:bg2103c}
\end{align}
The first identity holds as the difference between the LHS and the first term on the RHS is \eqref{eq:bg2102a} with $\Phi(\s,\x):=\mathbf{H}^{\N}(\s,\t(\N),\x)$. (Indeed, by Definitions \ref{definition:bg26} and \ref{definition:bg29}, the difference between $\mathscr{R}^{\chi,\mathfrak{q},\pm,\mathrm{j}}(\t,\x)$ and $\mathscr{A}^{\mathbf{X}}\mathscr{R}^{\chi,\mathfrak{q},\pm,\mathrm{j}}(\t,\x)$ is a spatial-averaging that turns $\mathds{R}^{\chi,\mathfrak{q},\pm,\mathrm{j}}(\s,\y(\s))\mathbf{Z}(\s,\y)$ to ${\mathds{A}_{\mathbf{X}}^{\mathfrak{m}(\mathrm{j}),\pm}}\{\mathds{R}^{\chi,\mathfrak{q},\pm,\mathrm{j}}\mathbf{Z};\s,\y(\s)\}$.) The second identity holds by using the same reasoning. This turns $\mathscr{T}^{\pm,\mathrm{j}}\mathscr{R}^{\chi,\mathfrak{q},\pm,\mathrm{j}}$ into the last two terms in the first line. \eqref{eq:bg2103c} follows by iterating until we get $(\mathscr{T}^{\pm,\mathrm{j}})^{5}\mathscr{R}^{\chi,\mathfrak{q},\pm,\mathrm{j}}$. Given the previous display, to show the desired claim \eqref{eq:bg210I}, it suffices to show the following estimate for the last term in \eqref{eq:bg2103c}:
\begin{align}
\|(\mathscr{T}^{\pm,\mathrm{j}})^{5}\mathscr{R}^{\chi,\mathfrak{q},\pm,\mathrm{j}}\|_{\t_{\mathrm{st}};\mathbb{T}(\N)} \ \lesssim \  \N^{-3\beta_{\mathrm{BG}}}\|\mathbf{Z}\|_{\t_{\mathrm{st}};\mathbb{T}(\N)}. \label{eq:bg2104}
\end{align}
Again, recall $\mathscr{R}^{\chi,\mathfrak{q},\pm,\mathrm{j}}$ in Definition \ref{definition:bg26}. (It is the time-integral of heat operators acting on the space-time function $\N\mathds{R}^{\chi,\mathfrak{q},\pm,\mathrm{j}}(\s,\y(\s))\mathbf{Z}(\s,\y)$.) Thus, we can compute $(\mathscr{T}^{\pm,\mathrm{j}})^{5}\mathscr{R}^{\chi,\mathfrak{q},\pm,\mathrm{j}}$ by moving $(\mathscr{T}^{\pm,\mathrm{j}})^{5}$ onto the heat kernel. Now recall that $|\mathds{R}^{\chi,\mathfrak{q},\pm,\mathrm{j}}(\s,\cdot)|\lesssim  \N^{20\gamma_{\mathrm{reg}}}\mathfrak{l}(\mathrm{j}-1)^{-3/2}\leq  \N^{20\gamma_{\mathrm{reg}}}$; this follows from construction in Definition \ref{definition:bg26}. So, for any $\t\leq\t_{\mathrm{st}}$ and $\x\in\mathbb{T}(\N)$, we get
\begin{align}
&|(\mathscr{T}^{\pm,\mathrm{j}})^{5}\mathscr{R}^{\chi,\mathfrak{q},\pm,\mathrm{j}}(\t,\x)| \nonumber\\
&= \ |{\textstyle\int_{0}^{\t}{\sum}_{\y}}(\mathscr{T}^{\pm,\mathrm{j}})^{5}\mathbf{H}^{\N}(\s,\t(\N),\x-\y)\cdot\N\mathds{R}^{\chi,\mathfrak{q},\pm,\mathrm{j}}(\s,\y(\s))\mathbf{Z}(\s,\y)\d\s| \label{eq:bg2105a}\\
&\leq \ \|\N\mathds{R}^{\chi,\mathfrak{q},\pm,\mathrm{j}}\|\|\mathbf{Z}\|{\textstyle\int_{0}^{\t}{\sum}_{\y}}|(\mathscr{T}^{\pm,\mathrm{j}})^{5}\mathbf{H}^{\N}(\s,\t(\N),\x-\y)|\d\s \ \lesssim \  \N^{1+20\gamma_{\mathrm{reg}}}\|\mathbf{Z}\|\Upsilon, \label{eq:bg2105b}
\end{align}
where $\|\|=\|\|_{\t_{\mathrm{st}};\mathbb{T}(\N)}$ and $\Upsilon$ is defined and estimated in the following calculation (that we explain afterwards):
\begin{align}
\Upsilon \ &:= \ {\int_{0}^{\t}\sum_{\y\in\mathbb{T}(\N)}}|(\mathscr{T}^{\pm,\mathrm{j}})^{5}\mathbf{H}^{\N}(\s,\t(\N),\x-\y)|\d\s \nonumber\\
&\lesssim \ \sup_{{1}\leq|\mathrm{k}_{\mathrm{i}}|\lesssim\mathfrak{m}(\mathrm{j})\mathfrak{l}(\mathrm{j})}{\int_{0}^{\t}\sum_{\y\in\mathbb{T}(\N)}}|\grad^{\mathrm{k}_{1}}\ldots\grad^{\mathrm{k}_{5}}\mathbf{H}^{\N}(\s,\t(\N),\x-\y)|\d\s. \label{eq:bg2106}
\end{align}
Indeed, the operator $\mathscr{T}^{\pm,\mathrm{j}}$ averages discrete gradients on length-scales $\lesssim\mathfrak{m}(\mathrm{j})\mathfrak{l}(\mathrm{j})$; we refer to Definition \ref{definition:bg29}. Therefore, its fifth-power is an average of compositions of five discrete gradients whose length-scales are all $\mathrm{O}(\mathfrak{m}(\mathrm{j})\mathfrak{l}(\mathrm{j}))$. Bounding the average by the supremum gives \eqref{eq:bg2106}. Recall from Definition \ref{definition:bg29} that $|\mathfrak{m}(\mathrm{j})\mathfrak{l}(\mathrm{j})|\lesssim\N^{3/4+\alpha}(\mathrm{j})$ with $\alpha(\mathrm{j})\leq 2{c}\gamma_{\mathrm{KL}}$ {for some small but fixed $c>0$}. We now observe that the heat kernel $\mathbf{H}^{\N}$ is smooth on macroscopic scales. Thus, each $\grad^{\mathrm{k}}$ yields a factor of $\lesssim  \N^{-1}|\mathrm{k}||\t(\N)-\s|^{-1/2}$. (The extra factor $|\t(\N)-\s|^{-1/2}$ is the usual heat kernel singularity.) Precisely, by {\eqref{eq:hke3}} and \eqref{eq:bg2106} and $\t\leq\t_{\mathrm{st}}\leq1$, we claim
\begin{align}
\Upsilon \ &\lesssim \ \sup_{1\leq|\mathrm{k}_{\mathrm{i}}|\lesssim\mathfrak{m}(\mathrm{j})\mathfrak{l}(\mathrm{j})} \N^{-5}{\textstyle\prod_{\mathrm{i}=1}^{5}}|\mathrm{k}_{\mathrm{i}}|{\textstyle\int_{0}^{\t}}|\t(\N)-\s|^{-\frac52}\d\s \nonumber\\
&\lesssim \ \N^{-5}|\mathfrak{m}(\mathrm{j})\mathfrak{l}(\mathrm{j})|^{5}|\t(\N)-\t|^{-\frac52} \ \lesssim \ \N^{-\frac54+5\alpha(\mathrm{j})+500\gamma_{\mathrm{reg}}}. \label{eq:bg2107}
\end{align}
Indeed, the first bound in \eqref{eq:bg2107} comes from the aforementioned heat kernel spatial regularity estimates in {\eqref{eq:hke3}}. The second bound comes from first replacing each $|\mathrm{k}_{\mathrm{i}}|$ with its maximal value $|\mathfrak{m}(\mathrm{j})\mathfrak{l}(\mathrm{j})|$. Then, we note that for $0\leq\s\leq\t$, the factor $|\t(\N)-\s|^{-5/2}$ is maximized at $\s=\t$. The last inequality in \eqref{eq:bg2107} now follows from $|\mathfrak{m}(\mathrm{j})\mathfrak{l}(\mathrm{j})|\lesssim  \N^{3/4+\alpha(\mathrm{j})}$ (see the previous paragraph) and by definition of $\t(\N)=\t+ \N^{-100\gamma_{\mathrm{reg}}}$ in Definition \ref{definition:method5}. By \eqref{eq:bg2105a}-\eqref{eq:bg2105b}, \eqref{eq:bg2106}, and \eqref{eq:bg2107}, with probability 1, 
\begin{align}
|(\mathscr{T}^{\pm,\mathrm{j}})^{5}\mathscr{R}^{\chi,\mathfrak{q},\pm,\mathrm{j}}(\t,\x)| \ \lesssim \  \N^{-\frac14+20\gamma_{\mathrm{reg}}+5\alpha(\mathrm{j})+500\gamma_{\mathrm{reg}}}\|\mathbf{Z}\|_{\t_{\mathrm{st}};\mathbb{T}(\N)}. \label{eq:bg2108}
\end{align}
Recall that $\alpha(\mathrm{j})\leq 2{c}\gamma_{\mathrm{KL}}$ {for some small but fixed $c>0$} from the paragraph before \eqref{eq:bg2107}. Recall also that $\gamma_{\mathrm{KL}}\leq1$; see Definition \ref{definition:entropydata}. Lastly, recall $\gamma_{\mathrm{reg}}={c}\gamma_{\mathrm{KL}}$ from Definition \ref{definition:reg}, and thus $\gamma_{\mathrm{reg}}>0$ is small. Hence, the exponent on the RHS of \eqref{eq:bg2108} is at most $-1/3$. But $1/3$ is bigger than $3\beta_{\mathrm{BG}}$, because $\beta_{\mathrm{BG}}$ is small (see Definition \ref{definition:method8}). Combining this paragraph with \eqref{eq:bg2108}, we get \eqref{eq:bg2104}. As noted right before \eqref{eq:bg2104}, this completes the proof of the desired estimate \eqref{eq:bg210I}, so we are done. \qed
%
%
%
\section{Local equilibrium estimates}\label{section:le}
We now prepare ingredients for proofs of Propositions \ref{prop:bg212}, \ref{prop:bg213}. First, we gather estimates concerning (local comparison to) local equilibrium measures; see Definition \ref{definition:intro5}. This is the purpose of this section. Besides one result on entropy production, the analysis in this section does not force us to deal with the time-inhomogeneous nature of \eqref{eq:hf}-\eqref{eq:glsde}. (Even for entropy production, the actual estimates do not become more difficult due to time-inhomogeneity.) We deal with stochastic homogenization problems and obstructions due to time-inhomogeneity in the next section.

Moreover, because this section does not quite see any real problems from time-inhomogeneity, the results here are somewhat standard. In particular, in a first reading the reader is invited to skip proofs (which we include for the sake of being complete, since there are a few details that prevent us from directly citing previous work). In any case, we will give a parallel to every result in this section for the sake of providing context.
\subsection{Square-root cancellations in space}
The following is a large-deviations estimate for $\mathbb{P}^{\sigma,\t,\mathbb{I}}$-measures. It is essentially a general martingale inequality, and it does not depend on the measure being $\mathbb{P}^{\sigma,\t,\mathbb{I}}$. First, some notation.
\begin{definition}\label{definition:le1}
 Take any $\mathfrak{a}:\R^{\mathbb{T}(\N)}\to\R$. The \emph{support} of $\mathfrak{a}$ is the smallest discrete interval $\mathbb{I}\subseteq\mathbb{T}(\N)$ such that $\mathfrak{a}$ depends only on $\mathbf{U}(\x)$ for $\x\in\mathbb{I}$. (By ``discrete interval", we mean the intersection of an interval in $\N\mathbb{T}$ and the lattice $\mathbb{T}(\N)\subseteq\N\mathbb{T}$.)
\end{definition}
\begin{lemma}\label{lemma:le2}
 Fix $\mathbb{I}\subseteq\mathbb{T}(\N)$, {$\mathrm{m}\in\mathbb{N}$}, $\sigma\in\R$, and $\t\geq0$. Take functionals $\mathfrak{a}(\cdot;1),\ldots,\mathfrak{a}(\cdot,;\mathrm{m})$ and discrete intervals $\mathbb{J}(1),\ldots,\mathbb{J}(\mathrm{m})\subseteq\mathbb{J}$. Now, we define $\mathscr{F}(\mathrm{j})$ to be the sigma-algebra generated by $\mathbf{U}(\z)$ for $\z\in\mathbb{J}(1)\cup\ldots\cup\mathbb{J}(\mathrm{j}-1)$, where $\mathbf{U}\in\R^{\mathbb{J}}$ denotes our ``dummy" variable. Suppose that $\E^{\sigma,\t,\mathbb{J}}\{\mathfrak{a}(\mathbf{U};\mathrm{j})|\mathscr{F}(\mathrm{j})\}=0$, where $\E^{\sigma,\t,\mathbb{J}}\{\cdot|\mathscr{F}(\mathrm{j})\}$ is $\E^{\sigma,\t,\mathbb{J}}$ but conditioning on $\mathscr{F}(\mathrm{j})$. Letting $\|\|_{\infty}$ be sup-norm, we have the following estimate for any constant $\mathrm{C}>0$:
\begin{align}
\mathbb{P}^{\sigma,\t,\mathbb{J}}[\mathrm{m}^{-1}|\mathfrak{a}(\mathbf{U};1)+\ldots+\mathfrak{a}(\mathbf{U};\mathrm{m})|\gtrsim\mathrm{C}\mathrm{m}^{-\frac12}{\sup}_{\mathrm{j}}\|\mathfrak{a}(\cdot,\mathrm{j})\|_{\infty}] \ \lesssim \ \exp\{-\mathrm{C}^{2}\}. \label{eq:le2I}
\end{align}
If $\mathfrak{a}(\cdot,;\mathrm{j})$ are each sub-Gaussian with variance parameter $\mathfrak{v}(\mathrm{j})^{2}$, then \eqref{eq:le2I} holds upon replacing $\|\mathfrak{a}(\cdot,\mathrm{j})\|_{\infty}\mapsto\mathfrak{v}(\mathrm{j})$.
\end{lemma}
\begin{proof}
The sequence $\mathrm{j}\mapsto\mathscr{F}(\mathrm{j})$ is a filtration. By $\E^{\sigma,\t,\mathbb{J}}\{\mathfrak{a}(\mathbf{U};\mathrm{j})|\mathscr{F}(\mathrm{j})\}=0$, the discrete-time process $\mathrm{j}\mapsto\mathrm{m}^{-1}\{\mathfrak{a}(\mathbf{U};1)+\ldots+\mathfrak{a}(\mathbf{U};\mathrm{j})\}$ is a martingale with respect to this filtration. So, \eqref{eq:le2I} follows by {the Azuma inequality} (even after $\|\mathfrak{a}(\cdot,\mathrm{j})\|_{\infty}\mapsto\mathfrak{v}(\mathrm{j})$).
\end{proof}
\subsection{Functional inequalities}
We now present some estimates with respect to $\mathbb{P}^{\sigma,\t,\mathbb{I}}$-measures related to the generator of \eqref{eq:glsde} (or ``localized" versions of \eqref{eq:glsde} on subsets $\mathbb{I}$). We start with a log-Sobolev inequality. This has two consequences. First, it helps compare any probability measure with local equilibrium at the cost of relative entropy plus large deviations. Second, it implies a spectral gap for the (localized) generator of \eqref{eq:glsde}. This spectral gap will be clarified later in this subsection. Lastly, we emphasize that the results of this subsection do not reflect {the} time-inhomogeneity of the dynamics \eqref{eq:glsde}, except for the fact that our generators depend on $\t$. (In particular, fix $\t$, and pretend \eqref{eq:glsde} is time-homogeneous with potential $\mathscr{U}(\t,\cdot)$. This subsection shows estimates for generators of such processes that {are} uniform in the family $\{\mathscr{U}(\t,\cdot)\}_{\t\lesssim1}$.) Again, we will denote dummy variables by $\mathbf{U}$.
\begin{definition}\label{definition:le3}
 Fix any subset $\mathbb{I}\subseteq\mathbb{T}(\N)$, any time $\t\geq0$, and any charge density $\sigma\in\R$. We let $\mathfrak{p}$ denote a probability density with respect to $\mathbb{P}^{\sigma,\t,\mathbb{I}}$. We define the \emph{Fisher information} of $\mathfrak{p}$ (or $\mathfrak{p}\d\mathbb{P}^{\sigma,\t,\mathbb{I}}$) with respect to $\mathbb{P}^{\sigma,\t,\mathbb{I}}$ as
\begin{align}
\mathfrak{D}_{\mathrm{FI}}^{\sigma,\t,\mathbb{I}}(\mathfrak{p}) \ = \ \mathfrak{D}_{\mathrm{FI}}^{\sigma,\t,\mathbb{I}}(\mathfrak{p}\d\mathbb{P}^{\sigma,\t,\mathbb{I}}) \ = \ \sum_{\x,\x+1\in\mathbb{I}}\E^{\sigma,\t,\mathbb{I}}|\mathrm{D}_{\x}\sqrt{\mathfrak{p}(\mathbf{U})}|^{2} \quad\mathrm{where}\quad \mathrm{D}_{\x}=\partial_{\mathbf{U}(\x+1)}-\partial_{\mathbf{U}(\x)}. \label{eq:le3doperators}
\end{align}
We call the quantity $\mathfrak{D}_{\mathrm{FI}}^{\sigma,\t,\mathbb{I}}(\mathfrak{p}^{2})$ the \emph{Dirichlet form} of $\mathfrak{p}$ with respect to $\mathbb{P}^{\sigma,\t,\mathbb{I}}$. Next, define the following $\mathbb{I}$-local relative entropy:
\begin{align}
\mathfrak{D}_{\mathrm{KL}}^{\sigma,\t,\mathbb{I}}(\mathfrak{p}) \ = \ \mathfrak{D}_{\mathrm{KL}}^{\sigma,\t,\mathbb{I}}(\mathfrak{p}\d\mathbb{P}^{\sigma,\t,\mathbb{I}}) \ = \ \E^{\sigma,\t,\mathbb{I}}\mathfrak{p}\log\mathfrak{p}.
\end{align}
For convenience, in the case $\mathbb{I}=\mathbb{T}(\N)$, we also set $\mathfrak{D}^{\sigma,\t}_{\mathrm{FI}}:=\mathfrak{D}^{\sigma,\t,\mathbb{T}(\N)}_{\mathrm{FI}}$ and $\mathfrak{D}^{\sigma,\t}_{\mathrm{KL}}:=\mathfrak{D}^{\sigma,\t,\mathbb{T}(\N)}_{\mathrm{KL}}$.
\end{definition}
\begin{lemma}\label{lemma:le4}
 Fix any discrete interval $\mathbb{I}\subseteq\mathbb{T}(\N)$, time $\t\lesssim1$, and charge density $\sigma\in\R$. Let $\mathfrak{p}$ denote a probability density with respect to $\mathbb{P}^{\sigma,\t,\mathbb{I}}$. We first claim the following logarithmic Sobolev inequality, which depends diffusively in the length-scale $|\mathbb{I}|$:
\begin{align}
\mathfrak{D}_{\mathrm{KL}}^{\sigma,\t,\mathbb{I}}(\mathfrak{p}) \ \lesssim \ |\mathbb{I}|^{2}\mathfrak{D}_{\mathrm{FI}}^{\sigma,\t,\mathbb{I}}(\mathfrak{p}). \label{eq:le4I}
\end{align}
We {stress that} the implied constant is independent of $\t$ and $\sigma$. For any constant $\kappa>0$ and function $\mathfrak{a}:\R^{\mathbb{I}}\to\R$, we also have
\begin{align}
\E^{\sigma,\t,\mathbb{I}}\mathfrak{p}\mathfrak{a} \ \lesssim \ \tfrac{1}{\kappa}\mathfrak{D}_{\mathrm{KL}}^{\sigma,\t,\mathbb{I}}(\mathfrak{p}) + \tfrac{1}{\kappa}\log\E^{\sigma,\t,\mathbb{I}}\exp\left(\kappa|\mathfrak{a}|\right) \lesssim \ \tfrac{|\mathbb{I}|^{2}}{\kappa}\mathfrak{D}_{\mathrm{FI}}^{\sigma,\t,\mathbb{I}}(\mathfrak{p}) + \tfrac{1}{\kappa}\log\E^{\sigma,\t,\mathbb{I}}\exp\left(\kappa|\mathfrak{a}|\right). \label{eq:le4II}
\end{align}
\end{lemma}
\begin{proof}
\eqref{eq:le4I} follows by second derivative bounds on $\mathscr{U}(\t,\cdot)$ in Assumption \ref{ass:intro8} and Bakry-Emery estimates; see (3.20)-(3.22) in \cite{CYau}. The first upper bound in \eqref{eq:le4II} is {the} classical duality between relative entropy and large deviations; see the inequality after Proposition 8.1 in Appendix 1 of \cite{KL}. The second bound in \eqref{eq:le4II} follows by \eqref{eq:le4I}.
\end{proof}
We now introduce notation for the generator of the joint process \eqref{eq:hf}-\eqref{eq:glsde}, but localized (in space) from the global torus $\mathbb{T}(\N)$ to any discrete interval $\mathbb{I}$. Its necessity comes from the local aspect of our analysis. We will give an intuitive explanation for the following construction in Remark \ref{remark:le6}. 
\begin{definition}\label{definition:le5}
 Take a discrete interval $\mathbb{I}\subseteq\mathbb{T}(\N)$ and $\t\geq0$. First, we let $\mathbb{S}(\N)$ be the one-dimensional torus of length $2\N^{20\gamma_{\mathrm{reg}}}$. We choose coordinates $\mathbb{S}(\N)\simeq[-\N^{20\gamma_{\mathrm{reg}}},\N^{20\gamma_{\mathrm{reg}}}]$ equipped with a periodic boundary. 
Set $\mathscr{L}^{\mathrm{tot}}(\t,\mathbb{I}):=\mathscr{L}^{\mathrm{curr}}(\t,\mathbb{I})+\mathscr{L}(\t,\mathbb{I})$, which acts on sufficiently smooth functions $\mathsf{F}:\mathbb{S}(\N)\times\R^{\mathbb{I}}\to\R$, as follows. Set $\mathscr{L}(\t,\mathbb{I}):=\mathscr{L}^{\mathrm{S}}(\t,\mathbb{I})+\mathscr{L}^{\mathrm{A}}(\t,\mathbb{I})$. Here, $\mathscr{L}^{?}(\t,\mathbb{I})$ are the following differential operators, in which $\mathrm{a}\in\mathbb{S}(\N)$ and $\mathbf{U}\in\R^{\mathbb{I}}$ are dummy variables, and $\grad^{\mathbb{I},?}$ is just $\grad^{?}$ from Definition \ref{definition:intro4} but with respect to periodic boundary conditions on the discrete interval $\mathbb{I}\subseteq\mathbb{T}(\N)$:
\begin{align}
(\mathscr{L}^{\mathrm{S}}(\t,\mathbb{I})\mathsf{F})(\mathrm{a},\mathbf{U}) \ &:= \  \N^{2}\sum_{\x\in\mathbb{I}}(\mathrm{D}_{\x}^{2}\mathsf{F}-\grad^{\mathbb{I},+}\mathscr{U}'[\t,\mathbf{U}(\x)]\mathrm{D}_{\x}\mathsf{F})(\mathrm{a},\mathbf{U}) \\
(\mathscr{L}^{\mathrm{A}}(\t,\mathbb{I})\mathsf{F})(\mathrm{a},\mathbf{U}) \ &:= \  \N^{\frac32}\sum_{\x\in\mathbb{I}}\{\grad^{\mathbb{I},\pm}\mathscr{U}'[\t,\mathbf{U}(\x)]\}\partial_{\mathbf{U}(\x)}\mathsf{F}(\mathrm{a},\mathbf{U}).
\end{align}
The $\mathrm{D}$-operators are from {\eqref{eq:le3doperators} in} Definition \ref{definition:le3}; they act on $\mathsf{F}$ through the second variable $\mathbf{U}$. Let us now define $\mathscr{L}^{\mathrm{curr}}(\t,\mathbb{I})$ as follows. In what follows, $\mathrm{a}$-derivatives are with respect to the periodic boundary conditions on $\mathbb{S}(\N)$, which is the state space for $\mathrm{a}$ (and completely unrelated to $\mathbb{I}$ or $\mathbb{T}(\N)$). {We also introduce the coefficient ${b(\t,\mathbf{U}):=\N^{\frac32}\grad^{+}\mathscr{U}'(\t,\mathbf{U}(0))\d\t + \N\{\mathscr{U}'(\t,\mathbf{U}(0))+\mathscr{U}'(\t,\mathbf{U}(1))\}}$ given by the drift in \eqref{eq:hf} at $\x=0$}. {With this notation, we define}
\begin{align}
(\mathscr{L}^{\mathrm{curr}}(\t,\mathbb{I})\mathsf{F})(\mathrm{a},\mathbf{U}) \ {:=} \  \N\partial_{\mathrm{a}}^{2}\mathsf{F}(\mathrm{a},\mathbf{U}) + {b(\t,\mathbf{U})}\partial_{\mathrm{a}}\mathsf{F}(\mathrm{a},\mathbf{U}).
\end{align}
Set $\mathbb{P}^{\mathrm{Leb},\sigma,\t,\mathbb{I}}:=\mathrm{Leb}[\mathbb{S}(\N)]\times\mathbb{P}^{\sigma,\t,\mathbb{I}}$, where $\mathrm{Leb}[\mathbb{S}(\N)]$ is {the} uniform measure on $\mathbb{S}(\N)$; see Definition \ref{definition:intro5} for $\mathbb{P}^{\sigma,\t,\mathbb{I}}$. (Note {that} $\mathrm{Leb}[\mathbb{S}(\N)]$ and $\mathbb{P}^{\mathrm{Leb},\sigma,\t,\mathbb{I}}$ are probability measures, as $\mathbb{S}(\N)$ is finite volume.) Let $\E^{\mathrm{Leb},\sigma,\t,\mathbb{I}}$ be {the} expectation with respect to $\mathbb{P}^{\mathrm{Leb},\sigma,\t,\mathbb{I}}$. Also, let $\E^{\mathrm{Leb}}$ be {the} expectation with respect to $\mathrm{Leb}[\mathbb{S}(\N)]$. (In particular, $\E^{\mathrm{Leb}}$ just denotes integration on the continuum torus $\mathbb{S}(\N)$ with respect to normalized Lebesgue measure. For clarity, we note that by Fubini, we have $\E^{\mathrm{Leb},\sigma,\t,\mathbb{I}}=\E^{\mathrm{Leb}}\E^{\sigma,\t,\mathbb{I}}=\E^{\sigma,\t,\mathbb{I}}\E^{\mathrm{Leb}}$.) 
\end{definition}
\begin{rem}\label{remark:le6}
 Suppose $\mathbb{S}(\N)$ were replaced by $\R$ for now. Also, take $\mathbb{I}=\mathbb{T}(\N)$ for simplicity. In this case, $\mathscr{L}^{\mathrm{tot}}(\t,\mathbb{I})$ equals the infinitesimal generator for the joint process $\t\mapsto(\mathbf{J}(\t,0),\mathbf{U}^{\t,\cdot})$; see \eqref{eq:hf}-\eqref{eq:glsde}. We refer to Section 2 of \cite{DGP} for why $\mathscr{L}(\t,\mathbb{I})$ is the generator for $\mathbf{U}^{\t,\cdot}$. (In a nutshell, if we take $\mathsf{F}(\mathrm{a},\mathbf{U})=\mathbf{U}(\x)$ for a fixed $\x\in\mathbb{T}(\N)$ and apply $\mathscr{L}(\t,\mathbb{I})$ to this $\mathsf{F}$, then $\mathscr{L}^{\mathrm{S}}(\t,\mathbb{I})$ yields the $\Delta$-term in \eqref{eq:glsde} and $\mathscr{L}^{\mathrm{A}}(\t,\mathbb{I})$ yields the $\grad^{\mathrm{a}}$-term.) Moreover, the operator $\mathscr{L}^{\mathrm{curr}}(\t,\mathbb{I})$ describes the evolution of \eqref{eq:hf} at $\x=0=\inf\mathbb{T}(\N)$ as a Brownian motion plus drifts depending only on $\mathbf{U}^{\t,\x}$ for $\x$ within 1 of $0$. (We have allowed ourselves a general choice of said drift; it turns out to be quite important to include in the evolution of \eqref{eq:hf} the renormalization term $\mathscr{R}(\t)$ from Definition \ref{definition:intro6}.) For the case of any general discrete interval $\mathbb{I}\subseteq\mathbb{T}(\N)$, the same picture holds, but \eqref{eq:hf}-\eqref{eq:glsde} are ``localized in space" from $\mathbb{T}(\N)$ to $\mathbb{I}$. Also, the use of $\mathbb{S}(\N)$ instead of $\R$ as the state space for $\mathbf{J}(\t,0)$ is entirely technical. (The benefit is that uniform measure on $\mathbb{S}(\N)$ is a well-defined probability measure; Lebesgue on $\R$ is not.) If we exclusively work before time $\t_{\mathrm{st}}$ from Definition \ref{definition:method8}, this reduction $\R\mapsto\mathbb{S}(\N)$ basically comes for free. (This is by a priori estimates on $\mathbf{Z}$ before $\t_{\mathrm{st}}$.) Now, note $\mathbf{J}(\t,\cdot)$ can be recovered by $(\mathbf{J}(\t,0),\mathbf{U}^{\t,\cdot})$. Indeed, any function on a torus is recovered by its value at one point and its gradients. Finally, summation-by-parts lets us move $\grad^{\mathbb{I},\pm}$-gradients in $\mathscr{L}^{\mathrm{A}}(\t,\mathbb{I})$ onto $\mathbf{U}(\x)$-partials therein. So $\mathscr{L}^{\mathrm{A}}(\t,\mathbb{I})$ is a sum of $\mathrm{D}_{\x}$.
\end{rem}
We now present important facts about $\mathscr{L}^{\mathrm{tot}}(\t,\mathbb{I})$ and its resolvents. First, recall $\E^{\mathrm{Leb},\sigma,\t,\mathbb{I}}$ from Definition \ref{definition:le5}. The following result (Lemma \ref{lemma:le7}) is essentially a spectral gap estimate, like that in Section 3 of \cite{DGP}. It effectively says that the total generator $\mathscr{L}^{\mathrm{tot},\mathrm{sym}}(\t,\mathbb{I})$ can be inverted when acting on functions orthogonal to its kernel, and said inverse has nice bounds. The only difference between the following and Section 3 of \cite{DGP} is that the generator $\mathscr{L}^{\mathrm{tot},\mathrm{sym}}(\t,\mathbb{I})$ includes the dynamics of the height function at $\inf\mathbb{I}$, not just the $\mathbf{U}$-process localized to $\mathbb{I}$. Since we still know explicitly the invariant measure of the fixed-time operator $\mathscr{L}^{\mathrm{tot},\mathrm{sym}}(\t,\mathbb{I})$, essentially nothing changes.
\begin{lemma}\label{lemma:le7}
 Fix a discrete interval $\mathbb{I}\subseteq\mathbb{T}(\N)$, time $\t\geq0$, and charge density $\sigma\in\R$. We first claim the infinitesimal generator $\mathscr{L}^{\mathrm{tot}}(\t,\mathbb{I})$ has $\E^{\mathrm{Leb},\sigma,\t,\mathbb{I}}$ as an invariant measure. Second, define $\mathscr{L}^{\mathrm{tot},\mathrm{sym}}(\t,\mathbb{I})$ as the symmetric part of $\mathscr{L}^{\mathrm{tot}}(\t,\mathbb{I})$ with respect to $\E^{\mathrm{Leb},\sigma,\t,\mathbb{I}}$. Now, set $\mathscr{G}^{\mathrm{tot},\mathrm{sym}}(\t,\mathbb{I}):=(-\mathscr{L}^{\mathrm{tot},\mathrm{sym}}(\t,\mathbb{I}))^{-1}$. Take $\mathfrak{a}(\cdot;\mathrm{j})$ in {Lemma \ref{lemma:le2}}, and set $\mathfrak{a}=\mathrm{m}^{-1}(\mathfrak{a}(\cdot;1)+\ldots+\mathfrak{a}(\cdot;\mathrm{m}))$. (Also take the subsets $\mathbb{J}(\mathrm{j})$ in {Lemma \ref{lemma:le2}}.) Fix $\varphi:\mathbb{S}(\N)\to\R$. For the function $\mathsf{F}(\mathrm{a},\mathbf{U}):=\varphi(\mathrm{a})\mathfrak{a}(\mathbf{U}):\mathbb{S}(\N)\times\R^{\mathbb{I}}\to\R$, we have
\begin{align}
&\E^{\mathrm{Leb},\sigma,\t,\mathbb{I}}\left[\mathsf{F}\times\mathscr{G}^{\mathrm{tot},\mathrm{sym}}(\t,\mathbb{I})\mathsf{F}\right] \ \lesssim \  \N^{-2}\mathrm{m}^{-1}\{\E^{\mathrm{Leb}}|\varphi|^{2}\}\times{\sup}_{\mathrm{j}}\{|\mathbb{J}(\mathrm{j})|^{2}\|\mathfrak{a}(\cdot;\mathrm{j})\|_{\infty}^{2}\}. \label{eq:le7II}
\end{align}
Lastly, \eqref{eq:le7II} is true if we swap $\mathscr{G}^{\mathrm{tot},\mathrm{sym}}(\t,\mathbb{I})$ for $[\lambda-\mathscr{L}^{\mathrm{tot},\mathrm{sym}}(\t,\mathbb{I})]^{-1}$ given any $\lambda\geq0$.
\end{lemma}
(Let us intuitively describe the bound \eqref{eq:le7II} before we start the proof. The $\N^{-2}$-factor comes from inverting the operator $\mathscr{L}^{\mathrm{tot},\mathrm{sym}}(\t,\mathbb{I})$, which comes with a speed $\N^{2}$. The $\mathrm{m}^{-1}$-factor comes from orthogonality, which results from the assumptions in Lemma \ref{lemma:le2}. Everything else on the RHS of \eqref{eq:le7II} is just a second moment, since we are working in Hilbert spaces. We note that the $|\mathbb{I}(\mathrm{j})|^{2}$-factor on the RHS of \eqref{eq:le7II} is there simply to reflect that the spectral gap, i.e. relaxation of the $\mathbf{U}$ process, decays quadratically in the length-scale, just like a simple random walk.)
\begin{proof}
Our strategy is to show that including the dynamics of the height function at $\inf\mathbb{I}$ introduces no difficulty and can be ``factored out", since the invariant measures in Definition \ref{definition:le5} factor. Then, we appeal to the work in Section 3 of \cite{DGP}. Fix any test function $\mathsf{H}:\mathbb{S}(\N)\times\R^{\mathbb{I}}\to\R$ that is smooth with compact support. We claim
\begin{align}
&\E^{\mathrm{Leb},\sigma,\t,\mathbb{I}}\mathscr{L}^{\mathrm{tot}}(\t,\mathbb{I})\mathsf{H} \nonumber\\
&= \ \E^{\mathrm{Leb}}\E^{\sigma,\t,\mathbb{I}}\mathscr{L}^{\mathrm{tot}}(\t,\mathbb{I})\mathsf{H}(\mathrm{a},\mathbf{U}) \ = \ \E^{\mathrm{Leb}}\E^{\sigma,\t,\mathbb{I}}\mathscr{L}^{\mathrm{curr}}(\t,\mathbb{I})\mathsf{H}(\mathrm{a},\mathbf{U})+\E^{\mathrm{Leb}}\E^{\sigma,\t,\mathbb{I}}\mathscr{L}(\t,\mathbb{I})\mathsf{H}(\mathrm{a},\mathbf{U}) \nonumber \\
&= \ \E^{\mathrm{Leb}}\E^{\sigma,\t,\mathbb{I}}\{\N\partial_{\mathrm{a}}^{2}\mathsf{H}(\mathrm{a},\mathbf{U}) + {b(\t,\mathbf{U})}\partial_{\mathrm{a}}\mathsf{H}(\mathrm{a},\mathbf{U})\} + \E^{\mathrm{Leb}}\E^{\sigma,\t,\mathbb{I}}\mathscr{L}(\t,\mathbb{I})\mathsf{H}(\mathrm{a},\mathbf{U}) \\
&= \ \E^{\mathrm{Leb}}\E^{\sigma,\t,\mathbb{I}}\{\N\partial_{\mathrm{a}}^{2}\mathsf{H}(\mathrm{a},\mathbf{U}) + {b(\t,\mathbf{U})}\partial_{\mathrm{a}}\mathsf{H}(\mathrm{a},\mathbf{U})\} \ = \ 0.
\end{align}
The first line follows by definition of $\E^{\mathrm{Leb},\sigma,\t,\mathbb{I}}$ and $\mathscr{L}^{\mathrm{tot}}(\t,\mathbb{I})$; see Definition \ref{definition:le5}. The second line is by definition of $\mathscr{L}^{\mathrm{curr}}$; see Definition \ref{definition:le5}. The first identity in the last line follows because $\mathscr{L}(\t,\mathbb{I})$ has $\E^{\sigma,\t,\mathbb{I}}$ as an invariant measure for any $\sigma$; see Section 2 of \cite{DGP} and Remark \ref{remark:intro12}. The last line follows by swapping $\E^{\mathrm{Leb}}$ and $\E^{\sigma,\t,\mathbb{I}}$, then integrating-by-parts in $\E^{\mathrm{Leb}}$. (Recall $\E^{\mathrm{Leb}}$ is just integration on the continuum torus $\mathbb{S}(\N)$ with respect to Lebesgue measure. Thus $\E^{\mathrm{Leb}}\partial_{\mathrm{a}}=0$.) So $\E^{\mathrm{Leb},\sigma,\t,\mathbb{I}}$ is invariant for $\mathscr{L}^{\mathrm{tot}}(\t,\mathbb{I})$. We now compute the symmetric part $\mathscr{L}^{\mathrm{tot},\mathrm{sym}}(\t,\mathbb{I})$. Take smooth, compactly supported $\mathsf{H}^{(1)},\mathsf{H}^{(2)}:\mathbb{S}(\N)\times\R^{\mathbb{I}}\to\R$. Following the previous display, we claim $\E^{\mathrm{Leb},\sigma,\t,\mathbb{I}}\mathsf{H}^{(1)}(\mathrm{a},\mathbf{U})\mathscr{L}^{\mathrm{tot}}(\t,\mathbb{I})\mathsf{H}^{(2)}(\mathrm{a},\mathbf{U})$ {can be written as}
{
\begin{align*}
&\E^{\mathrm{Leb}}\E^{\sigma,\t,\mathbb{I}}\mathsf{H}^{(1)}(\mathrm{a},\mathbf{U})\{\N\partial_{\mathrm{a}}^{2}\mathsf{H}^{(2)}(\mathrm{a},\mathbf{U}) + {b(\t,\mathbf{U})}\partial_{\mathrm{a}}\mathsf{H}^{(2)}(\mathrm{a},\mathbf{U})\}\\
&+ \E^{\mathrm{Leb}}\E^{\sigma,\t,\mathbb{I}}\mathsf{H}^{(1)}(\mathrm{a},\mathbf{U})\mathscr{L}(\t,\mathbb{I})\mathsf{H}^{(2)}(\mathrm{a},\mathbf{U}) \\
&= \ \E^{\mathrm{Leb}}\E^{\sigma,\t,\mathbb{I}}\mathsf{H}^{(2)}(\mathrm{a},\mathbf{U})\{\N\partial_{\mathrm{a}}^{2}\mathsf{H}^{(1)}(\mathrm{a},\mathbf{U})-{b(\t,\mathbf{U})}\partial_{\mathrm{a}}\mathsf{H}^{(1)}(\mathrm{a},\mathbf{U})\}\\
&+ \E^{\mathrm{Leb}}\E^{\sigma,\t,\mathbb{I}}\mathsf{H}^{(2)}(\mathrm{a},\mathbf{U})\mathscr{L}(\t,\mathbb{I})^{\ast}\mathsf{H}^{(1)}(\mathrm{a},\mathbf{U}).
\end{align*}
}The second identity follows by integration-by-parts in the $\d\mathrm{a}$ integral and by definition (upon setting $\mathscr{L}(\t,\mathbb{I})^{\ast}$ as the adjoint of $\mathscr{L}(\t,\mathbb{I})$ with respect to $\E^{\sigma,\t,\mathbb{I}}$). (The $\d\mathrm{a}$-integration-by-parts depends heavily on the fact that $\mathscr{B}$ depends only on $\mathbf{U}$, not on $\mathrm{a}$. The physics behind this statement is that the growth of KPZ-type height functions depends only on the local slope, not the height itself.)  By Section 2 of \cite{DGP} and Remark \ref{remark:intro12}, we get $\mathscr{L}(\t,\mathbb{I})^{\ast}=\mathscr{L}^{\mathrm{S}}(\t,\mathbb{I})-\mathscr{L}^{\mathrm{A}}(\t,\mathbb{I})$, so $\mathscr{L}^{?}(\t,\mathbb{I})$ is symmetric if $?=\mathrm{S}$ and anti-symmetric for $?=\mathrm{A}$. By the above display, the adjoint $\mathscr{L}^{\mathrm{tot}}(\t,\mathbb{I})^{\ast}$ with respect to $\E^{\mathrm{Leb},\sigma,\t,\mathbb{I}}$ is just $\mathscr{L}^{\mathrm{tot}}(\t,\mathbb{I})$ itself, with two adjustments. First, in $\mathscr{L}^{\mathrm{curr}}(\t,\mathbb{I})$, change $\mathscr{B}\mapsto-\mathscr{B}$. Second, in $\mathscr{L}(\t,\mathbb{I})$, change $\mathscr{L}^{\mathrm{A}}\mapsto-\mathscr{L}^{\mathrm{A}}$. So, we deduce by this and $\mathscr{L}^{\mathrm{tot},\mathrm{sym}}(\t,\mathbb{I}) := \tfrac12\{\mathscr{L}^{\mathrm{tot}}(\t,\mathbb{I})+\mathscr{L}^{\mathrm{tot}}(\t,\mathbb{I})^{\ast}\}$ that
{
\begin{align}
\mathscr{L}^{\mathrm{tot},\mathrm{sym}}(\t,\mathbb{I}) = \tfrac12\{\N\partial_{\mathrm{a}}^{2}+\N\partial_{\mathrm{a}}^{2}+\mathscr{L}^{\mathrm{S}}(\t,\mathbb{I})+\mathscr{L}^{\mathrm{S}}(\t,\mathbb{I})\} = \N\partial_{\mathrm{a}}^{2}+\mathscr{L}^{\mathrm{S}}(\t,\mathbb{I}). \label{eq:le7II0}
\end{align}
}We now (formally) establish \eqref{eq:le7II}; we make it rigorous at the end. We use Sobolev duality to compute $\mathrm{LHS}\eqref{eq:le7II}$. In particular, $\mathrm{LHS}\eqref{eq:le7II}$ is a negative Sobolev norm of degree $-1$ as $\mathscr{G}^{\mathrm{tot},\mathrm{sym}}(\t,\mathbb{I})$ is an inverse-Laplacian-type operator. Precisely, by Section 6 in Appendix 1 of \cite{KL}, we have the following dual form of the degree $-1$ Sobolev norm in terms of a degree $+1$ norm:
\begin{align}
\mathrm{LHS}\eqref{eq:le7II} \ &\lesssim \ {\textstyle\sup_{\mathsf{H}}}\left(2\E^{\mathrm{Leb},\sigma,\t,\mathbb{I}}\mathsf{F}\mathsf{H}+\E^{\mathrm{Leb},\sigma,\t,\mathbb{I}}\mathsf{H}\cdot\mathscr{L}^{\mathrm{tot},\mathrm{sym}}(\t,\mathbb{I})\mathsf{H}\right). \label{eq:le7II3a}
\end{align}
To be clear, the supremum on the RHS of \eqref{eq:le7II3a} is over compactly supported and smooth functions $\mathsf{H}:\mathbb{S}(\N)\times\R^{\mathbb{I}}\to\R$, which are dense in positive-degree Sobolev spaces. To estimate the RHS of \eqref{eq:le7II3a}, we claim
\begin{align}
\E^{\mathrm{Leb},\sigma,\t,\mathbb{I}}\mathsf{H}\cdot\mathscr{L}^{\mathrm{tot},\mathrm{sym}}(\t,\mathbb{I})\mathsf{H} \ &= \ \E^{\mathrm{Leb},\sigma,\t,\mathbb{I}}\mathsf{H}(\mathrm{a},\mathbf{U})\cdot\N\partial_{\mathrm{a}}^{2}\mathsf{H}(\mathrm{a},\mathbf{U}) + \E^{\mathrm{Leb},\sigma,\t,\mathbb{I}}\mathsf{H}\cdot\mathscr{L}^{\mathrm{S}}(\t,\mathbb{I})\mathsf{H} \\
&= \ \E^{\sigma,\t,\mathbb{I}}\E^{\mathrm{Leb}}\mathsf{H}(\mathrm{a},\mathbf{U})\cdot\N\partial_{\mathrm{a}}^{2}\mathsf{H}(\mathrm{a},\mathbf{U})+ \E^{\mathrm{Leb},\sigma,\t,\mathbb{I}}\mathsf{H}\cdot\mathscr{L}^{\mathrm{S}}(\t,\mathbb{I})\mathsf{H} \\
&= \ -\N\E^{\sigma,\t,\mathbb{I}}\E^{\mathrm{Leb}}|\partial_{\mathrm{a}}\mathsf{H}(\mathrm{a},\mathbf{U})|^{2}+ \E^{\mathrm{Leb},\sigma,\t,\mathbb{I}}\mathsf{H}\cdot\mathscr{L}^{\mathrm{S}}(\t,\mathbb{I})\mathsf{H}.
\end{align}
The first line follows by \eqref{eq:le7II0}. The second line follows by $\E^{\mathrm{Leb},\sigma,\t,\mathbb{I}}=\E^{\sigma,\t,\mathbb{I}}\E^{\mathrm{Leb}}$ (see Definition \ref{definition:le5}). The third line follows from integrating-by-parts in $\E^{\mathrm{Leb}}$. (Again, $\E^{\mathrm{Leb}}$ is just Lebesgue-measure integration on the continuum torus $\mathbb{S}(\N)$. So, the third line is the usual Laplacian-to-Dirichlet-energy calculation.) So, for an upper bound on $\mathrm{RHS}\eqref{eq:le7II3a}$, we can forget $\mathrm{a}$-differentials to swap $\mathscr{L}^{\mathrm{tot},\mathrm{sym}}$ for $\mathscr{L}^{\mathrm{S}}$. This gives the first line below; we explain the rest after:
\begin{align}
\mathrm{RHS}\eqref{eq:le7II3a} \ &\leq \ {\textstyle\sup_{\mathsf{H}}}\{2\E^{\mathrm{Leb},\sigma,\t,\mathbb{I}}\mathsf{F}\mathsf{H}+\E^{\mathrm{Leb},\sigma,\t,\mathbb{I}}\mathsf{H}\cdot\mathscr{L}^{\mathrm{S}}(\t,\mathbb{I})\mathsf{H}\} \label{eq:le7II3b} \\
&= \ {\textstyle\sup_{\mathsf{H}}}\E^{\mathrm{Leb}}\{2\E^{\sigma,\t,\mathbb{I}}\mathsf{F}(\mathrm{a},\mathbf{U})\mathsf{H}(\mathrm{a},\mathbf{U})+\E^{\sigma,\t,\mathbb{I}}\mathsf{H}(\mathrm{a},\mathbf{U})\mathscr{L}^{\mathrm{S}}(\t,\mathbb{I})\mathsf{H}(\mathrm{a},\mathbf{U})\} \label{eq:le7II3c}\\
&\leq \ \E^{\mathrm{Leb}}{\textstyle\sup_{\mathsf{G}}}\{2\E^{\sigma,\t,\mathbb{I}}\mathsf{F}(\mathrm{a},\mathbf{U})\mathsf{G}(\mathbf{U})+\E^{\sigma,\t,\mathbb{I}}\mathsf{G}(\mathbf{U})\mathscr{L}^{\mathrm{S}}(\t,\mathbb{I})\mathsf{G}(\mathbf{U})\}\label{eq:le7II3d} \\
&= \ \E^{\mathrm{Leb}}{\textstyle\sup_{\mathsf{G}}}\{2\E^{\sigma,\t,\mathbb{I}}\varphi(\mathrm{a})\mathfrak{a}(\mathbf{U})\mathsf{G}(\mathbf{U})+\E^{\sigma,\t,\mathbb{I}}\mathsf{G}(\mathbf{U})\mathscr{L}^{\mathrm{S}}(\t,\mathbb{I})\mathsf{G}(\mathbf{U})\} \label{eq:le7II3e}\\
&= \ \E^{\mathrm{Leb}}|\varphi(\mathrm{a})|^{2}{\textstyle\sup_{\mathsf{G}}}\{2\E^{\sigma,\t,\mathbb{I}}\mathfrak{a}(\mathbf{U})\mathsf{G}(\mathbf{U})+\E^{\sigma,\t,\mathbb{I}}\mathsf{G}(\mathbf{U})\mathscr{L}^{\mathrm{S}}(\t,\mathbb{I})\mathsf{G}(\mathbf{U})\}.\label{eq:le7II3f}
\end{align}
The first line is explained in the previous paragraph. The second line follows by $\E^{\mathrm{Leb},\sigma,\t,\mathbb{I}}=\E^{\mathrm{Leb}}\E^{\sigma,\t,\mathbb{I}}$. The third line follows by moving the supremum inside $\E^{\mathrm{Leb}}$; this gives an upper bound, because we get to optimize the test function we take supremum over per the dummy variable $\mathrm{a}$ in $\E^{\mathrm{Leb}}$. (In \eqref{eq:le7II3d}-\eqref{eq:le7II3f}, the sup is over smooth, compactly supported functions $\mathsf{G}$ in just the $\mathbf{U}$ variable.) The fourth line follows by definition of $\mathsf{F}(\mathrm{a},\mathbf{U})=\varphi(\mathrm{a})\mathfrak{a}(\mathbf{U})$. The final line follows by reparameterizing the sup in \eqref{eq:le7II3e} by $\mathsf{G}\mapsto\varphi(\mathrm{a})\mathsf{G}$ per $\mathrm{a}$. We now claim the following estimate, which we justify shortly:
\begin{align}
{\textstyle\sup_{\mathsf{G}}}\{2\E^{\sigma,\t,\mathbb{I}}\mathfrak{a}(\mathbf{U})\mathsf{G}(\mathbf{U})+\E^{\sigma,\t,\mathbb{I}}\mathsf{G}(\mathbf{U})\mathscr{L}^{\mathrm{S}}(\t,\mathbb{I})\mathsf{G}(\mathbf{U})\} \ \lesssim \  \N^{-2}\mathrm{m}^{-1} {\sup}_{\mathrm{j}}\{|\mathbb{J}(\mathrm{j})|^{2}\|\mathfrak{a}(\cdot;\mathrm{j})\|_{\infty}^{2}\}. \label{eq:le7II3g}
\end{align}
\eqref{eq:le7II3g}, \eqref{eq:le7II3a}, and \eqref{eq:le7II3b}-\eqref{eq:le7II3f} imply \eqref{eq:le7II} and thus complete the proof. It remains to explain \eqref{eq:le7II3g}. To this end, follow Section 3 of \cite{DGP}. This, in turn, follows Proposition 7 of \cite{GJ15}. To be self-contained, we give a brief description of the argument. First, recall $\mathfrak{a}$ is the average of $\mathfrak{a}(\cdot;1),\ldots,\mathfrak{a}(\cdot;\mathrm{m})$. By assumption in Lemma \ref{lemma:le2}, these functionals are orthogonal with respect to $\E^{\sigma,\t,\mathbb{I}}$. (The proof of) Proposition 7 in \cite{GJ15} shows the orthogonality holds in negative-Sobolev-norm (of degree $-1$) as well. This explains the factor of $\mathrm{m}^{-1}$ in $\mathrm{RHS}\eqref{eq:le7II3g}$. We are left to bound degree $-1$ Sobolev norms of $\mathfrak{a}(\cdot;\mathrm{j})$ for all $\mathrm{j}$. Precisely, we must derive $\E^{\sigma,\t,\mathbb{I}}\mathfrak{a}(\mathbf{U};\mathrm{j})\{[-\mathscr{L}^{\mathrm{S}}(\t,\mathbb{I})]^{-1}\mathfrak{a}(\mathbf{U};\mathrm{j})\}\lesssim\N^{-2}|\mathbb{J}(\mathrm{j})|^{2}\|\mathfrak{a}(\cdot;\mathrm{j})\|_{\infty}^{2}$. This bound would follow immediately by a spectral gap for $-\mathscr{L}^{\mathrm{S}}(\t,\mathbb{I})$ that we get from the LSI \eqref{eq:le4I}, if we replace $|\mathbb{J}(\mathrm{j})|$ in the desired resolvent bound with $|\mathbb{I}|$. However, the support of $\mathfrak{a}(\mathbf{U};\mathrm{j})$ is in $\mathbb{J}(\mathrm{j})$ by assumption in Lemma \ref{lemma:le2}, and therefore $\E^{\sigma,\t,\mathbb{I}}\mathfrak{a}(\mathbf{U};\mathrm{j})\{[-\mathscr{L}^{\mathrm{S}}(\t,\mathbb{I})]^{-1}\mathfrak{a}(\mathbf{U};\mathrm{j})\}\lesssim\N^{-2}|\mathbb{J}(\mathrm{j})|^{2}\|\mathfrak{a}(\cdot;\mathrm{j})\|_{\infty}^{2}$ still follows from spectral gap considerations. (The point is that every nearest-neighbor bond in $\mathbb{I}$ corresponds to a non-negative term in the dual form $\E^{\sigma,\t,\mathbb{I}}\mathfrak{a}(\mathbf{U};\mathrm{j})\{[-\mathscr{L}^{\mathrm{S}}(\t,\mathbb{I})]\mathfrak{a}(\mathbf{U};\mathrm{j})\}$. So we drop all bonds whose vertices are not both in $\mathbb{J}(\mathrm{j})$ for the sake of an upper bound on the $-\mathscr{L}^{\mathrm{S}}(\t,\mathbb{I})^{-1}$-form. We are then left with the $-\mathscr{L}^{\mathrm{S}}(\t,\mathbb{I}(\mathrm{j}))^{-1}$-form, not the $-\mathscr{L}^{\mathrm{S}}(\t,\mathbb{I})^{-1}$-form.) Also, we clarify that Proposition 7 in \cite{GJ15} does not give an $\N^{-2}$-factor that we say appears on the RHS of \eqref{eq:le7II3g}. This is just a matter of convention; the $\N^{-2}$ factor in \cite{GJ15} is delegated to the estimate immediately prior to Remark 8 in \cite{GJ15}. This finishes our formal argument for \eqref{eq:le7II}. The reason why it is formal is that not every $\mathscr{L}$-operator in this proof is bijective, thus its inverse is not well-defined. Instead, we should regularize the resolvent $-\mathscr{L}\mapsto\lambda-\mathscr{L}$, prove estimates independent of $\lambda$, and then take $\lambda\to0$. This is standard, so we do not do it. It remains to prove \eqref{eq:le7II} but replace $\mathscr{G}^{\mathrm{tot},\mathrm{sym}}(\t,\mathbb{I})$ by $[\lambda-\mathscr{L}^{\mathrm{tot},\mathrm{sym}}(\t,\mathbb{I})]^{-1}$ given any $\lambda\geq0$. For this, it suffices to note that as quadratic forms, we have $[\lambda-\mathscr{L}^{\mathrm{tot},\mathrm{sym}}(\t,\mathbb{I})]^{-1}\leq[-\mathscr{L}^{\mathrm{tot},\mathrm{sym}}(\t,\mathbb{I})]^{-1}=\mathscr{G}^{\mathrm{tot},\mathrm{sym}}(\t,\mathbb{I})$.
\end{proof}
\subsection{Entropy production}
We now present a classical entropy production bound (see \cite{GPV}). This says the evolution of relative entropy of \eqref{eq:glsde} with respect to local equilibrium dissipates via Fisher information and grows according to time-evolution of the local equilibrium reference measure, or equivalently, the time-derivative of $\mathscr{U}(\t,\cdot)$ (see \cite{YauRE}). Thus, to bound entropy production, which is crucial to compare local statistics to local equilibrium, we need control on $\partial_{\t}\mathscr{U}(\t,\cdot)$; see the first bullet in Section \ref{section:sqle}. We start with an auxiliary estimate {that} will be important for the proof of entropy production (Lemma \ref{lemma:le8}). Then, we use Lemma \ref{lemma:le8} and improve Lemma \ref{lemma:le8a} to Lemma \ref{lemma:le9}. (Lemma \ref{lemma:le9} compares local statistics to local equilibrium.)

The proof of the following auxiliary estimate would essentially follow by the standard entropy inequality (Lemma \ref{lemma:le4}) if we had $\mathfrak{D}_{\mathrm{KL}}$ on the RHS instead of $\mathfrak{D}_{\mathrm{FI}}$. The log-Sobolev inequality, which is guaranteed by convexity (Assumption \ref{ass:intro8}), is then crucial to completing the proof of \eqref{eq:le8aI}.
\begin{lemma}\label{lemma:le8a}
 Fix any $\mathbb{I}\subseteq\mathbb{T}(\N)$, $\mathfrak{a}:[0,\infty)\times\R^{\mathbb{I}}\to\R$, and any $\t\lesssim1$. For any $\s,\y$, we let $\mathbf{U}^{\s,\y+\cdot}$ be the configuration in $\R^{\mathbb{T}(\N)}$ obtained after shifting the configuration $\z\mapsto\mathbf{U}^{\s,\z}$ by $\y$ in the spatial variable $\z\in\mathbb{T}(\N)$. Now, take a probability density $\mathfrak{p}$ with respect to $\mathbb{P}^{0,0,\mathbb{T}(\N)}$. Let $\mathfrak{p}(\t)$ be the density with respect to $\mathbb{P}^{0,\t,\mathbb{T}(\N)}$ for the measure $\mathbb{P}(\t)$ obtained after time-$\t$ evolution of \eqref{eq:glsde} with initial law $\mathfrak{p}\d\mathbb{P}^{0,0,\mathbb{T}(\N)}$. For any $\kappa>0$, we have the following in which for the last term, $\mathbf{U}$ is dummy variable for $\E^{\sigma,\s,\mathbb{I}}$:
\begin{align}
{\int_{0}^{\t}|\mathbb{T}(\N)|^{-1}{\sum_{\y\in\mathbb{T}(\N)}}}\E|\mathfrak{a}(\s,\mathbf{U}^{\s,\y+\cdot})|\d\s \ &\lesssim \ \tfrac{1}{\kappa}\N^{-1}|\mathbb{I}|^{3}{\textstyle\int_{0}^{\t}}\mathfrak{D}_{\mathrm{FI}}^{0,\s}(\mathfrak{p}(\s))\d\s \label{eq:le8aI}\\
&+\tfrac{1}{\kappa}\sup_{\substack{\sigma\in\R\\0\leq\s\leq\t}}\log\E^{\sigma,\s,\mathbb{I}}\exp\{\kappa|\mathfrak{a}(\s,\mathbf{U})|\}.\nonumber
\end{align}
Upon relabeling spatial variables, the same estimate holds if we replace $\y$ by any time-dependent shift $\y+\z(\s)$ on $\mathrm{LHS}\eqref{eq:le8aI}$.
\end{lemma}
\begin{proof}
Fix any $\s\geq0$ and $\y\in\mathbb{T}(\N)$. Because $\mathfrak{a}(\s,\mathbf{U})$ depends only on $\mathbf{U}(\x)$ for $\x\in\mathbb{I}$, we know that $\mathfrak{a}(\s,\mathbf{U}^{\s,\y+\cdot})$ depends only on $\mathbf{U}^{\s,\x}$ for $\x\in-\y+\mathbb{I}=:\mathbb{I}(\y)$. Now, some constructions. In what follows, $\E^{\s}$ always denotes expectation with respect to the law of $\mathbf{U}^{\s,\cdot}$. Let $\Pi^{\mathbb{I}(\y)}\E^{\s}$ be {the} expectation with respect to the marginal onto $\R^{\mathbb{I}(\y)}$. For any $\sigma\in\R$, we also set $\Pi^{\mathbb{I}(\y),\sigma}\E^{\s}$ as $\Pi^{\mathbb{I}(\y)}\E^{\s}$ after further conditioning on the set of all $\mathbf{U}\in\R^{\mathbb{I}(\y)}$ so that the average of $\mathbf{U}(\x)$ over $\x\in\mathbb{I}(\y)$ equals $\sigma$. Lastly, let $\mathfrak{p}^{\s,\y,\sigma}$ be the probability density of $\Pi^{\mathbb{I}(\y),\sigma}\E^{\s}$ with respect to $\mathbb{P}^{\sigma,\s,\mathbb{I}(\y)}$. We now claim that for some probability measure $\mathcal{Q}$ on $\R$,
\begin{align}
\E|\mathfrak{a}(\s,\mathbf{U}^{\s,\y+\cdot})| \ = \ \Pi^{\mathbb{I}(\y)}\E^{\s}|\mathfrak{a}(\s,\mathbf{U}^{\s,\y+\cdot})| \ &= \ {\textstyle\int_{\R}}\Pi^{\mathbb{I}(\y),\sigma}\E^{\s}|\mathfrak{a}(\s,\mathbf{U}^{\s,\y+\cdot})|\d\mathcal{Q}(\sigma) \nonumber\\
&= \ {\textstyle\int_{\R}}\E^{\sigma,\s,\mathbb{I}(\y)}\mathfrak{p}^{\s,\y,\sigma}|\mathfrak{a}(\s,\mathbf{U})|\d\mathcal{Q}(\sigma). \label{eq:le8aI0}
\end{align}
The first identity follows because $\mathfrak{a}(\s,\mathbf{U}^{\s,\y+\cdot})$ has support $\mathbb{I}(\y)$ (see the previous paragraph). The second follows by conditioning on the average of $\mathbf{U}(\x)$ over $\x\in\mathbb{I}(\y)$. (Technically, $\mathcal{Q}$ may depend on $\s,\y$, but this is not important.) The final identity follows by definition of $\mathfrak{p}^{\s,\y,\sigma}$. (We note $\mathbf{U}\in\R^{\mathbb{I}(\y)}$ on the far RHS is the expectation dummy variable; it is unrelated to $\mathbf{U}^{\s,\cdot}$.) By \eqref{eq:le4II},
\begin{align}
&{\textstyle\int_{\R}}\E^{\sigma,\s,\mathbb{I}(\y)}\mathfrak{p}^{\s,\y,\sigma}|\mathfrak{a}(\s,\mathbf{U})|\d\mathcal{Q}(\sigma) \nonumber\\
&\lesssim \ {\textstyle\int_{\R}}\tfrac{1}{\kappa}|\mathbb{I}(\y)|^{2}\mathfrak{D}^{\sigma,\s,\mathbb{I}(\y)}_{\mathrm{FI}}(\mathfrak{p}^{\s,\y,\sigma})\d\mathcal{Q}(\sigma)+{\textstyle\int_{\R}}\tfrac{1}{\kappa}\log\E^{\sigma,\s,\mathbb{I}(\y)}\exp\{\kappa|\mathfrak{a}(\s,\mathbf{U})|\}\d\mathcal{Q}(\sigma) \label{eq:le8aI1a} \\
&\lesssim \ \tfrac{1}{\kappa}|\mathbb{I}|^{2}{\textstyle\int_{\R}}\mathfrak{D}^{\sigma,\s,\mathbb{I}(\y)}_{\mathrm{FI}}(\mathfrak{p}^{\s,\y,\sigma})\d\mathcal{Q}(\sigma)+{\textstyle\sup_{\sigma}}\tfrac{1}{\kappa}\log\E^{\sigma,\s,\mathbb{I}(\y)}\exp\{\kappa|\mathfrak{a}(\s,\mathbf{U})|\}. \label{eq:le8aI1b}
\end{align}
\eqref{eq:le8aI1b} follows because $\mathbb{I}$ and $\mathbb{I}(\y)$ are shifts of each other, so their sizes are the same. (Also, we recall $\mathcal{Q}$ is a probability measure, so integrating against $\mathcal{Q}$ is bounded by a supremum over $\sigma$.) Now, observe that canonical ensembles are invariant under shifts in space, so for the second term in \eqref{eq:le8aI1b}, we can replace $\mathbb{I}(\y)$ by $\mathbb{I}(0)$. (Indeed, by Definition \ref{definition:intro5}, $\mathbb{P}^{\sigma,\s,\mathbb{I}(\y)}$ is the law of a random walk bridge whose steps are distributed via something independent of $\y$; they just are indexed by $\mathbb{I}(\y)$. When replacing $\mathbb{I}(\y)\mapsto\mathbb{I}$, all we do is change step indices.) So, as $\t\lesssim1$, space-time averaging the last term in \eqref{eq:le8aI1b} is bounded by its time-sup:
\begin{align}
\frac{1}{\kappa}{\int_{0}^{\t}}|\mathbb{T}(\N)|^{-1}{\sum_{\y}}{\sup_{\sigma}}\log\E^{\sigma,\s,\mathbb{I}(\y)}\exp\{\kappa|\mathfrak{a}(\s,\mathbf{U})|\}\d\s \ \lesssim \ \frac{1}{\kappa}{\sup_{\sigma}\sup_{\s\leq\t}}\log\E^{\sigma,\s,\mathbb{I}}\exp\{\kappa|\mathfrak{a}(\s,\mathbf{U})|\}. \label{eq:le8aI2}
\end{align}
On the other hand, we also have the following estimate, which we explain afterwards:
\begin{align}
{\textstyle\int_{0}^{\t}}{\textstyle\sum_{\y}}{\textstyle\int_{\R}}\mathfrak{D}^{\sigma,\s,\mathbb{I}(\y)}_{\mathrm{FI}}(\mathfrak{p}^{\s,\y,\sigma})\d\mathcal{Q}(\sigma)\d\s \ \lesssim \ |\mathbb{I}|{\textstyle\int_{0}^{\t}}\mathfrak{D}_{\mathrm{FI}}^{0,\s}(\mathfrak{p}(\s))\d\s. \label{eq:le8aI3}
\end{align}
\eqref{eq:le8aI3} is a classical convexity estimate for entropy production in hydrodynamic limits; see Lemma 2.2 in \cite{GPV}, for example. (In a nutshell, the Fisher information on the RHS of \eqref{eq:le8aI3} is summing the local energy attached to each nearest-neighbor bond in $\mathbb{T}(\N)$. The LHS of \eqref{eq:le8aI3} is summing the local energy attached to each nearest-neighbor bond in $\mathbb{I}$ and then sliding $\mathbb{I}$ across until it covers $\mathbb{T}(\N)$. The factor $|\mathbb{I}|$ is exactly the redundancy factor; when sliding $\mathbb{I}$, every bond is accounted for in $|\mathbb{I}|$-many different shifts. Technically, $\mathrm{LHS}\eqref{eq:le8aI3}$ takes $\R^{\mathbb{I}(\mathrm{y})}$-marginals before taking expectations/local energies; see Definition \ref{definition:le3}. Convexity of Fisher information says this only makes things smaller; this is why \eqref{eq:le8aI3} is not equality. By law of total expectation, \eqref{eq:le8aI3} writes $\Pi^{\mathbb{I}(\y)}\E^{\s}$, from the definition of Fisher information, as an average of $\Pi^{\mathbb{I}(\y),\sigma}\E^{\s}$ over all $\sigma\in\R$. This explains the structure of $\mathrm{LHS}\eqref{eq:le8aI3}$.) To get the desired bound \eqref{eq:le8aI}, it now suffices to combine \eqref{eq:le8aI0}, \eqref{eq:le8aI1a}-\eqref{eq:le8aI1b}, \eqref{eq:le8aI2}, and \eqref{eq:le8aI3}. 
\end{proof}
The following is a standard entropy production estimate; see \cite{GPV}, for example. As in \cite{GPV}, the time-integrated Fisher information is controlled by the initial relative entropy times the relaxation-speed factor $\N^{-2}$; from this, \eqref{eq:le8I} would follow immediately. However, in the time-inhomogeneous case, the reference measure for the Fisher information changes in time, so we need a ``relative entropy production" estimate as in \cite{YauRE}. More precisely, the time-derivative of $\mathfrak{D}_{\mathrm{KL}}$ functionals contains a term whose main factor is the time derivative of the reference measure, i.e. terms of the form $\partial_{\t}\mathscr{U}$. This turns out to be a fluctuating functional (as explained in point (1) in Section \ref{section:sqle}). In particular, we have a sum of $|\mathbb{T}(\N)|$-many fluctuating terms, which should have contribution $\lesssim\N^{1/2}$ by square-root cancellations. Multiplying this by the same $\N^{-2}$ speed-factor gives something much smaller than the RHS of \eqref{eq:le8I}. All the work in the proof of this is in making precise the ``fluctuating" notion. To this end, we use the preliminary estimate in Lemma \ref{lemma:le7} to reduce to showing that $\partial_{\t}\mathscr{U}$ is mean-zero with respect to the explicit measures in Definition \ref{definition:intro5}. But this is just calculus. Otherwise, the proof of Lemma \ref{lemma:le8} has nothing more to it, so the reader is again invited to skip the proof in a first reading.

(This is only for the reader interested in reading the proof of Lemma \ref{lemma:le8}. The way we show that a sum of $|\mathbb{T}(\N)|$-many fluctuating terms has square-root cancellation is by a standard one-block, two-blocks argument. Gather the fluctuating terms into groups of a slowly diverging size. Use ergodic theory to replace the average of each group by a ``local expectation". Use ergodic theory, again, to show that this ``local expectation" is stable under slowly increasing the size of the group over which we average, until we hit a sufficiently large size.) 
\begin{lemma}\label{lemma:le8}
 Suppose $\mathfrak{p}$ is the probability density with respect to $\mathbb{P}^{0,0,\mathbb{T}(\N)}$ for a probability measure that belongs to the class of entropy data (see {Definition \ref{definition:entropydata}} for the definition of entropy data). Again, we let $\mathfrak{p}(\t)$ denote the density with respect to $\mathbb{P}^{0,\t,\mathbb{T}(\N)}$ for the measure $\mathbb{P}(\t)$ obtained after time-$\t$ evolution of \eqref{eq:glsde} with initial law $\mathfrak{p}\d\mathbb{P}^{0,0,\mathbb{T}(\N)}$. For any $\t\leq1$, we have the bound
\begin{align}
{\textstyle\int_{0}^{\t}}\mathfrak{D}_{\mathrm{FI}}^{0,\s}(\mathfrak{p}(\s))\d\s \ \lesssim \  \N^{-\frac54-\gamma_{\mathrm{KL}}}. \label{eq:le8I}
\end{align}
\end{lemma}
\begin{proof}
Let $\mathbb{H}$ be the hyperplane of $\mathbf{U}\in\R^{\mathbb{T}(\N)}$ such that the average of $\mathbf{U}(\x)$ for $\x\in\mathbb{T}(\N)$ is zero. (Note $\mathbb{H}$ is the support of the canonical measure $\mathbb{P}^{0,\t,\mathbb{T}(\N)}$ for any $\t\geq0$.) We also clarify Lebesgue measure on $\mathbb{H}$ is the measure induced by Euclidean metric on $\mathbb{H}$. This metric is induced from $\mathbb{H}\subseteq\R^{\mathbb{T}(\N)}$ and standard Euclidean metric on $\R^{\mathbb{T}(\N)}$. Now, we let $\mathscr{L}(\t,\mathbb{T}(\N))^{\dagger}$ be adjoint of $\mathscr{L}(\t,\mathbb{T}(\N))$ with respect to Lebesgue measure on $\mathbb{H}$; see Definition \ref{definition:le5} for $\mathscr{L}(\t,\mathbb{T}(\N))$. We also let $\mathfrak{p}(\t)^{\dagger}$ be the density for the law of $\mathbf{U}^{\t,\cdot}$ with respect to Lebesgue measure on $\mathbb{H}$. By the Kolmogorov forward equation, we have
\begin{align}
\partial_{\t}\mathfrak{p}(\t)^{\dagger} \ = \ \mathscr{L}(\t,\mathbb{T}(\N))^{\dagger}\mathfrak{p}(\t)^{\dagger}, \label{eq:le8I1}
\end{align}
where $\mathscr{L}(\t,\mathbb{T}(\N))^{\dagger}$ acts on the implicit $\mathbf{U}$-variable that $\mathfrak{p}(\t)^{\dagger}$ depends on (and that we have omitted). Now, let $\mathfrak{p}(0,\t,\mathbb{T}(\N))$ be the density of $\mathbb{P}^{0,\t,\mathbb{T}(\N)}$ with respect to Lebesgue measure on $\mathbb{H}$. Thus, we know $\mathfrak{p}(\t)^{\dagger}=\mathfrak{p}(\t)\mathfrak{p}(0,\t,\mathbb{T}(\N))$. We now claim
{\small
\begin{align}
&\partial_{\t}\mathfrak{p}(\t) \nonumber\\
&= [\mathscr{L}(\t,\mathbb{T}(\N))^{\dagger}\mathfrak{p}(\t)^{\dagger}]\mathfrak{p}(0,\t,\mathbb{T}(\N))^{-1} - \{\partial_{\t}\mathfrak{p}(0,\t,\mathbb{T}(\N))\}\cdot\mathfrak{p}(0,\t,\mathbb{T}(\N))^{-2}\mathfrak{p}(\t)^{\dagger}\label{eq:le8II2a}\\
&= [\mathscr{L}(\t,\mathbb{T}(\N))^{\dagger}\mathfrak{p}(\t)^{\dagger}]\mathfrak{p}(0,\t,\mathbb{T}(\N))^{-1} - \{[\partial_{\t}\mathfrak{p}(0,\t,\mathbb{T}(\N))]\cdot\mathfrak{p}(0,\t,\mathbb{T}(\N))^{-1}\}\cdot\mathfrak{p}(0,\t,\mathbb{T}(\N))^{-1}\mathfrak{p}(\t)^{\dagger}\label{eq:le8II2aa}\\
&= [\mathscr{L}(\t,\mathbb{T}(\N))^{\dagger}\mathfrak{p}(\t)^{\dagger}]\mathfrak{p}(0,\t,\mathbb{T}(\N))^{-1} - \partial_{\t}\log\mathfrak{p}(0,\t,\mathbb{T}(\N))\cdot\mathfrak{p}(\t)\label{eq:le8II2aaa}\\
&= \mathscr{L}(\t,\mathbb{T}(\N))^{\ast}\mathfrak{p}(\t)- \partial_{\t}\log\mathfrak{p}(0,\t,\mathbb{T}(\N))\cdot\mathfrak{p}(\t). \label{eq:le8II2b}
\end{align}
}\eqref{eq:le8II2a} follows from $\mathfrak{p}(\t)^{\dagger}=\mathfrak{p}(\t)\mathfrak{p}(0,\t,\mathbb{T}(\N))$ and Leibniz rule. \eqref{eq:le8II2aa}-\eqref{eq:le8II2aaa} follow from elementary manipulations. \eqref{eq:le8II2b} follows because acting on the Lebesgue density $\mathfrak{p}(\t)^{\dagger}$ by $\mathscr{L}(\t,\mathbb{T}(\N))^{\dagger}$ and then changing measure via $\mathfrak{p}(0,\t,\mathbb{T}(\N))^{-1}$ is the same as just acting on the $\mathbb{P}^{0,\t,\mathbb{T}(\N)}$ density $\mathfrak{p}(\t)$ by the $\mathbb{P}^{0,\t,\mathbb{T}(\N)}$-adjoint $\mathscr{L}(\t,\mathbb{T}(\N))^{\ast}$. (This is the usual calculus that implies the adjoint in the Kolmogorov forward equation is always with respect to the reference measure that the density is defined on.) Let us now compute $\mathfrak{p}(0,\t,\mathbb{T}(\N))$. As a function of $\mathbf{U}\in\R^{\mathbb{T}(\N)}$, define the ``Hamiltonian" $\mathscr{H}(\t,\mathbf{U})$ as the sum over $\x\in\mathbb{T}(\N)$ of $\mathscr{U}(\t,\mathbf{U}(\x))$. We claim
\begin{align}
\mathfrak{p}(0,\t,\mathbb{T}(\N)) \ = \ \exp\{-\mathscr{H}(\t,\mathbf{U})+\mathscr{P}(\t)\} \ =: \ \exp\{-\mathscr{HP}(\t,\mathbf{U})\}, \label{eq:le8II3a}
\end{align}
where $\mathscr{HP}(\t,\mathbf{U})$ is ``Hamiltonian plus pressure", i.e. the sum over $\x\in\mathbb{T}(\N)$ of $\mathscr{U}(\t,\mathbf{U}(\x))-|\mathbb{T}(\N)|^{-1}\mathscr{P}(\t)$, and, if $\d\mathrm{Leb}(\cdot;\mathbb{H})$ is Lebesgue measure on $\mathbb{H}$,
\begin{align}
\mathscr{P}(\t) \ := \ -\log{\textstyle\int_{\mathbb{H}}}\exp(-\mathscr{H}(\t,\mathbf{U}))\d\mathrm{Leb}(\mathbf{U};\mathbb{H}). \label{eq:le8II3b}
\end{align}
\eqref{eq:le8II3a}-\eqref{eq:le8II3b} holds for the following reasons. It is proportional to $\exp[-\mathscr{H}(\t,\mathbf{U})]$, which, up to {a} constant factor, is the density of the grand-canonical measure $\mathbb{P}^{0,\t}$ with respect to {the} Lebesgue measure on $\R^{\mathbb{T}(\N)}$. (Indeed, by Definition \ref{definition:intro5}, $\mathbb{P}^{0,\t,\mathbb{T}(\N)}$ is just $\mathbb{P}^{0,\t}$ conditioned on the subset $\mathbb{H}$. So, up to constants, its Lebesgue density on the set $\mathbb{H}$ we condition on is the same as the Lebesgue density of $\mathbb{P}^{0,\t}$ with respect to {the} Lebesgue measure on $\mathbb{R}^{\mathbb{T}(\N)}$. There is also the indicator of $\mathbb{H}$ that we must multiply $\mathfrak{p}(0,\t,\mathbb{T}(\N))$ by. But this is redundant, as the reference measure $\mathrm{Leb}(\cdot;\mathbb{H})$ has this factor.) It now suffices to note $\mathscr{P}(\t)$ is exactly the constant that makes $\mathfrak{p}(0,\t,\mathbb{T}(\N))$ a probability density on $\mathbb{H}$. Now recall the notation of Definition \ref{definition:le3}. For some $\upsilon\gtrsim1$, we claim
\begin{align}
&\partial_{\t}\mathfrak{D}_{\mathrm{KL}}^{0,\t}(\mathfrak{p}(\t)) \ = \ \E^{0,\t,\mathbb{T}(\N)}\mathfrak{p}(\t)\mathscr{L}(\t,\mathbb{T}(\N))\log\mathfrak{p}(\t) + \E^{0,\t,\mathbb{T}(\N)}\mathfrak{p}(\t)\partial_{\t}\log\mathfrak{p}(\t) \label{eq:le8II4a}\\
&= \ \E^{0,\t,\mathbb{T}(\N)}\mathfrak{p}(\t)\mathscr{L}(\t,\mathbb{T}(\N))\log\mathfrak{p}(\t) + \E^{0,\t,\mathbb{T}(\N)}\partial_{\t}\mathfrak{p}(\t)\label{eq:le8II4b}\\
&= \ \E^{0,\t,\mathbb{T}(\N)}\mathfrak{p}(\t)\mathscr{L}(\t,\mathbb{T}(\N))\log\mathfrak{p}(\t)+\E^{0,\t,\mathbb{T}(\N)}\mathscr{L}(\t,\mathbb{T}(\N))^{\ast}\mathfrak{p}(\t)+\E^{0,\t,\mathbb{T}(\N)}\mathfrak{p}(\t)\partial_{\t}\mathscr{HP}(\t,\mathbf{U})\label{eq:le8II4c}\\
&= \ \E^{0,\t,\mathbb{T}(\N)}\mathfrak{p}(\t)\mathscr{L}^{\mathrm{S}}(\t,\mathbb{T}(\N))\log\mathfrak{p}(\t)+\E^{0,\t,\mathbb{T}(\N)}\mathfrak{p}(\t)\mathscr{L}^{\mathrm{A}}(\t,\mathbb{T}(\N))\log\mathfrak{p}(\t)+\label{eq:le8II4d}\\
&+ \ \E^{0,\t,\mathbb{T}(\N)}\mathfrak{p}(\t)\partial_{\t}\mathscr{HP}(\t,\mathbf{U})\nonumber\\
&\leq \ -\upsilon\N^{2}\mathfrak{D}_{\mathrm{FI}}^{0,\t}(\mathfrak{p}(\t))+\E^{0,\t,\mathbb{T}(\N)}\mathfrak{p}(\t)\partial_{\t}\mathscr{HP}(\t,\mathbf{U}).\label{eq:le8II4e}
\end{align}
To show \eqref{eq:le8II4a}, we know the LHS equals $\partial_{\t}\E^{0,\t,\mathbb{T}(\N)}\{\mathfrak{p}(\t)\log\mathfrak{p}(\t)\}=\partial_{\t}\E\log\mathfrak{p}(\t,\mathbf{U}^{\t,\cdot})$ in which the latter expectation is with respect to the law of $\mathbf{U}^{\t,\cdot}$. Leibniz rule produces two terms. Kolmogorov backward equation says the first is $\E\mathscr{L}(\t,\mathbb{T}(\N))\log\mathfrak{p}(\t)$, which is the first term on the RHS of \eqref{eq:le8II4a}. However, the test function that we apply Kolmogorov to is time-dependent as well. So we need to take its time-derivative. Thus, the second term is $\E\partial_{\t}\log\mathfrak{p}(\t)$, which is the last term in \eqref{eq:le8II4a}. This gives \eqref{eq:le8II4a}. \eqref{eq:le8II4b} is calculus. \eqref{eq:le8II4c} follows by \eqref{eq:le8II2a}-\eqref{eq:le8II2b}, \eqref{eq:le8II3a}. To explain \eqref{eq:le8II4d}, we note the second term in \eqref{eq:le8II4c} is zero. (Indeed, replace $\mathscr{L}(\t,\mathbb{T}(\N))^{\ast}$ by $\mathscr{L}(\t,\mathbb{T}(\N))$, and instead of having it act on $\mathfrak{p}(\t)$, let it act on $1$. But $\mathscr{L}(\t,\mathbb{T}(\N))$ is a differential, so its action on 1 vanishes.) \eqref{eq:le8II4e} follows by classical, explicit calculation for the first term in \eqref{eq:le8II4d}; see Section 2 of \cite{GPV}. It then suffices to note the second term in \eqref{eq:le8II4d} is zero. (Indeed, $\mathscr{L}^{\mathrm{A}}(\t,\mathbb{T}(\N))$ is a first-order differential. So, the second term in \eqref{eq:le8II4d} is $\E^{0,\t,\mathbb{T}(\N)}\mathscr{L}^{\mathrm{A}}(\t,\mathbb{T}(\N))\mathfrak{p}(\t)$ by calculus. Note the $\E^{0,\t,\mathbb{T}(\N)}$-adjoint of $\mathscr{L}^{\mathrm{A}}(\t,\mathbb{T}(\N))$ is $-\mathscr{L}^{\mathrm{A}}(\t,\mathbb{T}(\N))$; see right before \eqref{eq:le7II0}. Now, follow the justification for \eqref{eq:le8II4c}.) Integrate the above differential inequality in time over $\s\in[0,\t]$. This gives
\begin{align}
\N^{2}{\textstyle\int_{0}^{\t}}\mathfrak{D}_{\mathrm{FI}}^{0,\s}(\mathfrak{p}(\s))\d\s \ &\lesssim \ \mathfrak{D}_{\mathrm{KL}}^{0,0}(\mathfrak{p}(0))-\mathfrak{D}_{\mathrm{KL}}^{0,\t}(\mathfrak{p}(\t))+{\textstyle\int_{0}^{\t}}\E^{0,\s,\mathbb{T}(\N)}\mathfrak{p}(\s)\partial_{\s}\mathscr{HP}(\s,\mathbf{U})\d\s \\
&\leq \ \mathfrak{D}_{\mathrm{KL}}^{0,0}(\mathfrak{p}(0))+{\textstyle\int_{0}^{\t}}\E^{0,\s,\mathbb{T}(\N)}\mathfrak{p}(\s)(\mathbf{1}[\s\leq\t_{\mathrm{reg}}]+\mathbf{1}[\t_{\mathrm{reg}}<\s])\partial_{\s}\mathscr{HP}(\s,\mathbf{U})\d\s. \label{eq:le8II5}
\end{align}
(Note \eqref{eq:le8II5} follows because relative entropy is non-negative, so we can drop the second term on the RHS of the first line.) We now give properties of $\mathscr{HP}$ (i.e. ``Hamiltonian plus pressure") before bounding \eqref{eq:le8II5}. First, we compute $\partial_{\s}\mathscr{P}(\s)$ (see \eqref{eq:le8II3b}). For any $\x\in\mathbb{T}(\N)$, we claim
\begin{align}
\partial_{\s}\mathscr{P}(\s) \ &= \ -\partial_{\s}\log{\textstyle\int_{\mathbb{H}}}\exp\{-\mathscr{H}(\s,\mathbf{U})\} \ = \ -[{\textstyle\int_{\mathbb{H}}}\exp\{-\mathscr{H}(\s,\mathbf{U})\}]^{-1}\partial_{\s}{\textstyle\int_{\mathbb{H}}}\exp\{-\mathscr{H}(\s,\mathbf{U})\} \label{eq:le8II7a}\\
&= \ [{\textstyle\int_{\mathbb{H}}}\exp\{-\mathscr{H}(\s,\mathbf{U})\}]^{-1}{\textstyle\int_{\mathbb{H}}}\partial_{\s}\mathscr{H}(\s,\mathbf{U})\exp\{-\mathscr{H}(\s,\mathbf{U})\} \ = \ \E^{0,\s,\mathbb{T}(\N)}\partial_{\s}\mathscr{H}(\s,\mathbf{U})  \label{eq:le8II7b}\\
&= \ |\mathbb{T}(\N)|\E^{0,\s,\mathbb{T}(\N)}\partial_{\s}\mathscr{U}(\s,\mathbf{U}(\x)). \nonumber
\end{align}
\eqref{eq:le8II7a} follows by calculus. The first identity in \eqref{eq:le8II7b} follows by chain rule. The second follows by the definition of $\E^{0,\s,\mathbb{T}(\N)}$. (Indeed, we can replace the $[\cdot]^{-1}$-factor in \eqref{eq:le8II7b} by putting its negative-log in $\exp\{-\mathscr{H}(\s,\mathbf{U})\}$. This gives the integral over $\mathbb{H}$ of $\partial_{\s}\mathscr{H}(\s,\mathbf{U})$ against \eqref{eq:le8II3a} at time $\t=\s$. This is $\E^{0,\s,\mathbb{T}(\N)}$.) The last identity holds because $\mathscr{H}(\s,\mathbf{U})$ is the sum over $\x\in\mathbb{T}(\N)$ of $\mathscr{U}(\s,\mathbf{U}(\x))$. But, under $\E^{0,\s,\mathbb{T}(\N)}$, the $\mathbf{U}(\x)$ are exchangeable. (This is because the grand-canonical measure, that we condition to get $\E^{0,\s,\mathbb{T}(\N)}$, is invariant under permutations on $\mathbb{T}(\N)$. So is the set $\mathbb{H}$ we condition on.) Now, by definition, $\mathscr{HP}$ from after \eqref{eq:le8II3a} is the sum over $\x\in\mathbb{T}(\N)$ of $\mathscr{UP}(\s,\mathbf{U}(\x)):=\mathscr{U}(\s,\mathbf{U}(\x))-|\mathbb{T}(\N)|^{-1}\mathscr{P}(\s)$ ("potential $\mathscr{U}$ plus pressure"). So, \eqref{eq:le8II7a}-\eqref{eq:le8II7b} gives
\begin{align}
\partial_{\s}\mathscr{HP}(\s,\mathbf{U}) &= {\textstyle\sum_{\x}}\partial_{\s}\mathscr{UP}(\s,\mathbf{U}(\x)) \quad\mathrm{and}\quad \E^{0,\s,\mathbb{T}(\N)}\partial_{\s}\mathscr{UP}(\s,\mathbf{U}(\x))=0. \label{eq:le8II6c}
\end{align}
Moreover, because $|\partial_{\s}\mathscr{U}(\s,\cdot)|\lesssim1$ (see Assumption \ref{ass:intro8}), we also know $|\partial_{\s}\mathscr{UP}(\s,\cdot)|\leq|\partial_{\s}\mathscr{U}(\s,\cdot)|+|\mathbb{T}(\N)|^{-1}|\partial_{\s}\mathscr{P}(\s)|\lesssim1$. (Indeed, $|\mathbb{T}(\N)|^{-1}|\partial_{\s}\mathscr{P}(\s)|\lesssim\E^{0,\s,\mathbb{T}(\N)}|\partial_{\s}\mathscr{U}(\s,\mathbf{U}(\x))|\lesssim1$ by \eqref{eq:le8II7a}-\eqref{eq:le8II7b}.) Thus, by \eqref{eq:le8II6c}, we also get $|\partial_{\s}\mathscr{HP}(\s,\mathbf{U})|\lesssim|\mathbb{T}(\N)|$. We now estimate the time-integral in \eqref{eq:le8II5}. First, recall $\t_{\mathrm{reg}}=1$ with very high probability (see Theorem \ref{theorem:kpz}). We assumed $\t\leq1$ in this lemma, so in \eqref{eq:le8II5}, we always have $\s\leq1$. This means $\t_{\mathrm{reg}}<\s$ is very low probability in \eqref{eq:le8II5}. Thus, {for any large but finite $\mathrm{D}>0$}, we have
\begin{align}
{\textstyle\int_{0}^{\t}}\E^{0,\s,\mathbb{T}(\N)}\mathfrak{p}(\s)\mathbf{1}[\t_{\mathrm{reg}}<\s]\partial_{\s}\mathscr{HP}(\s,\mathbf{U})\d\s \ \lesssim \ |\mathbb{T}(\N)|{\textstyle\int_{0}^{\t}}\E^{0,\s,\mathbb{T}(\N)}\mathfrak{p}(\s)\mathbf{1}[\t_{\mathrm{reg}}<\s]\d\s \ \lesssim \ \N^{-{\mathrm{D}}}. \label{eq:le8II10}
\end{align}
Take the $\s\leq\t_{\mathrm{reg}}$ term in \eqref{eq:le8II5}. To this end, {we give} some more notation. Set $\alpha(\mathscr{U},\s):=\partial_{\sigma}\E^{\sigma,\s}\partial_{\s}\mathscr{UP}(\s,\mathbf{u})|_{\sigma=0}$, where $\mathbf{u}$ is the expectation dummy variable, and $\mathscr{UP}$ is the ``$\mathscr{U}$ plus pressure" from right before \eqref{eq:le8II6c}. We also define the ``centered" term $\mathscr{CP}(\s,\mathbf{u}):=\partial_{\s}\mathscr{UP}(\s,\mathbf{u})-\alpha(\mathscr{U},\s)\mathbf{u}$. Next, using notation from Definition \ref{definition:bg21}, we define the following centering (via local equilibrium expectation) of the $\mathscr{UP}$-term:
\begin{align}
\mathscr{UP}(\s,\y;\mathfrak{l}(\N)) \ := \ \partial_{\s}\mathscr{UP}(\s,\mathbf{U}(\y))-\E^{\mathfrak{l}(\N),+}[\partial_{\s}\mathscr{UP}(\s,\cdot);\s,\y] \quad\mathrm{where}\quad \mathfrak{l}(\N) \ := \ \N^{\frac14+100\gamma_{\mathrm{KL}}}. \label{eq:le8II11}
\end{align}
Lastly, let $\mathscr{AP}(\s,\y;\mathfrak{l}(\N))$ be the average of $\mathscr{UP}(\s,\y+2\mathrm{k}\mathfrak{l}(\N);\mathfrak{l}(\N))$ over $\mathrm{k}=0,\ldots,\lfloor\mathfrak{l}(\N)\rfloor$. We claim the following, in which $\Phi$ is the time-integral in \eqref{eq:le8II5} but without the $\mathbf{1}[\t_{\mathrm{reg}}<\s]$ term therein and, for convenience, we set $\mathbf{1}[\mathcal{E}(\s)]:=\mathbf{1}[\s\leq\t_{\mathrm{reg}}]$:
\begin{align}
\Phi \ &= \ {\textstyle\int_{0}^{\t}}{\textstyle\sum_{\y}}\E^{0,\s,\mathbb{T}(\N)}\{\mathfrak{p}(\s)\mathbf{1}(\mathcal{E}(\s))\partial_{\s}\mathscr{UP}(\s,\mathbf{U}(\y))\}\d\s \label{eq:le8I9a}\\
&= \ {\textstyle\int_{0}^{\t}}{\textstyle\sum_{\y}}\E^{0,\s,\mathbb{T}(\N)}\{\mathfrak{p}(\s)\mathbf{1}(\mathcal{E}(\s))\mathscr{UP}(\s,\y;\mathfrak{l}(\N))\}\d\s\label{eq:le8I9b}\\
&+ \ {\textstyle\int_{0}^{\t}}{\textstyle\sum_{\y}}\E^{0,\s,\mathbb{T}(\N)}\{\mathfrak{p}(\s)\mathbf{1}(\mathcal{E}(\s))\E^{\mathfrak{l}(\N),+}[\partial_{\s}\mathscr{UP}(\s,\cdot);\s,\y]\}\d\s \nonumber\\
&= \ {\textstyle\int_{0}^{\t}}{\textstyle\sum_{\y}}\E^{0,\s,\mathbb{T}(\N)}\mathfrak{p}(\s)\mathbf{1}(\mathcal{E}(\s))\mathscr{AP}(\s,\y;\mathfrak{l}(\N))\d\s \label{eq:le8I9c}\\
&+ \ {\textstyle\int_{0}^{\t}}{\textstyle\sum_{\y}}\E^{0,\s,\mathbb{T}(\N)}\{\mathfrak{p}(\s)\mathbf{1}(\mathcal{E}(\s))\E^{\mathfrak{l}(\N),+}[\partial_{\s}\mathscr{UP}(\s,\cdot);\s,\y]\}\d\s.\nonumber
\end{align}
\eqref{eq:le8I9a} is by \eqref{eq:le8II6c}. \eqref{eq:le8I9b} is by \eqref{eq:le8II11}. To establish \eqref{eq:le8I9c}, we first leave the last term in \eqref{eq:le8I9b} alone. For the first term in \eqref{eq:le8I9b}, note the only thing that depends on the $\y$-sum variable is $\mathscr{UP}(\s,\y;\mathfrak{l}(\N))$. Pull the $\y$-sum through to hit just this term. If we replace $\mathscr{UP}(\s,\y;\mathfrak{l}(\N))$ by the average $\mathscr{AP}(\s,\y;\mathfrak{l}(\N))$ of its spatial shifts, the error we get is an average of discrete gradients of $\mathscr{UP}(\s,\y;\mathfrak{l}(\N))$. Summing these discrete gradients over all $\y\in\mathbb{T}(\N)$ makes them vanish. This gives \eqref{eq:le8I9c}. We now control the second integral in \eqref{eq:le8I9c}. We first claim the following identity, which replaces $\mathscr{UP}$ by $\mathscr{CP}$ from right before \eqref{eq:le8II11}:
{\fontsize{9.5}{12}
\begin{align}
&{\textstyle\int_{0}^{\t}}{\textstyle\sum_{\y}}\E^{0,\s,\mathbb{T}(\N)}\{\mathfrak{p}(\s)\mathbf{1}(\mathcal{E}(\s))\E^{\mathfrak{l}(\N),+}[\partial_{\s}\mathscr{UP}(\s,\cdot);\s,\y]\}\d\s \nonumber\\
&= \ {\textstyle\int_{0}^{\t}}{\textstyle\sum_{\y}}\E^{0,\s,\mathbb{T}(\N)}\{\mathfrak{p}(\s)\mathbf{1}(\mathcal{E}(\s))\E^{\mathfrak{l}(\N),+}[\mathscr{CP}(\s,\cdot);\s,\y]\}\d\s. \label{eq:le8I9d}
\end{align}
}By construction of $\mathscr{CP}$, the difference $\partial_{\s}\mathscr{UP}(\s,\mathbf{u})-\mathscr{CP}(\s,\mathbf{u})$ equals $\alpha(\mathscr{U},\s)\mathbf{u}$. Taking $\E^{\mathfrak{l}(\N),+}$ of this term gives the average of $\mathbf{U}^{\s,\y+\mathrm{k}}$ over some set of $\mathrm{k}$. (Indeed, for any $\mathbb{I}$ and $\x\in\mathbb{I}$ and $\sigma$ and $\t$, we know $\E^{\sigma,\t,\mathbb{I}}\mathbf{U}(\x)=\sigma$. To derive this, note that $\mathbf{U}(\x)$ are exchangeable in $\x$ with respect to $\mathbb{P}^{\sigma,\t,\mathbb{I}}$; see after \eqref{eq:le8II7b}. Thus, $\E^{\sigma,\t,\mathbb{I}}\mathbf{U}(\x)$ is unchanged if we replace $\mathbf{U}(\x)$ by its average over $\x\in\mathbb{I}$. This is $\sigma$. Therefore, by Definition \ref{definition:bg21}, we have $\E^{\mathfrak{l}(\N),+}[\mathbf{U}(\x);\s,\y]=\sigma(\s,\y;\mathfrak{l}(\N),+)$, where $\mathbf{U}$ denotes expectation dummy-variable. It now suffices to recall, from Definition \ref{definition:bg21}, that $\sigma(\s,\y;\mathfrak{l}(\N),+)$ is an average of $\mathbf{U}^{\s,\y+\mathrm{k}}$ over some set of $\mathrm{k}$.) Because the expectation $\E^{0,\s,\mathbb{T}(\N)}$ is supported on the hyperplane $\mathbb{H}$ from the beginning of this argument, when we sum $\mathbf{U}^{\s,\y+\mathrm{k}}$ first over all $\y\in\mathbb{T}(\N)$, we get zero with probability 1. \eqref{eq:le8I9d} therefore follows. We now claim
\begin{align}
\mathscr{CP}(\s,\mathbf{U}(\y)) \ &= \ \partial_{\s}\mathscr{UP}(\s,\mathbf{U}(\y))-\alpha(\mathscr{U},\s)\mathbf{U}(\y) \nonumber\\
&= \ \partial_{\s}\mathscr{UP}(\s,\mathbf{U}(\y))-\E^{0,\s,\mathbb{T}(\N)}\partial_{\s}\mathscr{UP}(\s,\cdot)-\alpha(\mathscr{U},\s)\mathbf{U}(\y) \nonumber\\
&= \ \{\partial_{\s}\mathscr{UP}(\s,\mathbf{U}(\y))-\E^{0,\s}\partial_{\s}\mathscr{UP}(\s,\cdot)-\alpha(\mathscr{U},\s)\mathbf{U}(\y)\} \label{eq:le8I12}\\
&+ \ \{\E^{0,\s}\partial_{\s}\mathscr{UP}(\s,\cdot)-\E^{0,\s,\mathbb{T}(\N)}\partial_{\s}\mathscr{UP}(\s,\cdot)\}. \nonumber
\end{align}
The last identity is easy to check. The first two are by definition and then \eqref{eq:le8II6c}. We first estimate the last $\{\}$-term in \eqref{eq:le8I12}{; to} this end, we use the equivalence of ensembles. In particular, we apply Corollary B.3 in \cite{DGP}. (We did this in the proof of Lemma \ref{lemma:ee1}.) In the notation of Corollary B.3 in \cite{DGP}, we take $F=\partial_{\s}\mathscr{UP}(\s,\cdot)$. Its support length $\ell$ is 1 (since $\partial_{\s}\mathscr{UP}(\s,\cdot)$ depends only on a real-valued input $\cdot$). Again, in Corollary B.3 from \cite{DGP}, the $N$-factor therein is the length of the canonical ensemble domain. It is equal to $|\mathbb{T}(\N)|=\N$ in our case. Therefore, we get the following deterministic bound, which is basically the statement that if one takes Brownian bridge and Brownian motion (with the same average drift), their local increments are very close in law:
\begin{align}
|\E^{0,\s}\partial_{\s}\mathscr{UP}(\s,\cdot)-\E^{0,\s,\mathbb{T}(\N)}\partial_{\s}\mathscr{UP}(\s,\cdot)| \ \lesssim \ |\mathbb{T}(\N)|^{-1}. \label{eq:le8I13}
\end{align}
(Technically, the implied constant in \eqref{eq:le8I13} depends continuously on $\E^{0,\s}\mathbf{u}^{2}\lesssim1$, where the bound comes by uniform convexity in Assumption \ref{ass:intro8}. It also depends continuously on $\partial_{\sigma}^{\d}\E^{\sigma,\s}\partial_{\s}\mathscr{UP}(\s,\cdot)|_{\sigma=0}$ for $\d\lesssim1$. But because $\partial_{\s}\mathscr{UP}\lesssim1$, which we noted after \eqref{eq:le8II6c}, $\partial_{\sigma}^{\d}\E^{\sigma,\s}\partial_{\s}\mathscr{UP}(\s,\cdot)|_{\sigma=0}\lesssim1$ by the last paragraph in the proof of Lemmas \ref{lemma:bg23}, \ref{lemma:bg1hl2}. Therefore, the implied constant in \eqref{eq:le8I13} is $\lesssim1$.) We clarify \eqref{eq:le8I13} will be used to bound the last term in \eqref{eq:le8I9c} (see \eqref{eq:le8I12}). In this spirit, we now estimate the the first $\{\}$-term in \eqref{eq:le8I12}. We claim that this term, as a function of $\mathbf{U}(\y)\in\R$, is $\mathrm{LCT}$ (see Definition \ref{definition:method3}). Indeed, we have subtracted from $\partial_{\s}\mathscr{UP}(\s,\cdot)$ its $\E^{0,\s}$-mean and something that vanishes in $\E^{0,\s}$. On the other hand, 
\begin{align}
&\partial_{\sigma}\E^{\sigma,\s}[\partial_{\s}\mathscr{UP}(\s,\mathbf{U}(\y))-\E^{0,\s}\partial_{\s}\mathscr{UP}(\s,\mathbf{U}(\y))- \alpha(\mathscr{U},\s)\mathbf{U}(\y)]|_{\sigma=0} \\
&= \ \partial_{\sigma}\E^{\sigma,\s}[\partial_{\s}\mathscr{UP}(\s,\mathbf{U}(\y))]|_{\sigma=0}-\partial_{\sigma}\E^{0,\s}\partial_{\s}\mathscr{UP}(\s,\mathbf{U}(\y))|_{\sigma=0}-\alpha(\mathscr{U},\s)\partial_{\sigma}\E^{\sigma,\s}\mathbf{U}(\y)|_{\sigma=0} \\
&= \ \partial_{\sigma}\E^{\sigma,\s}[\partial_{\s}\mathscr{UP}(\s,\mathbf{U}(\y))]|_{\sigma=0} - \alpha(\mathscr{U},\s)\partial_{\sigma}\sigma|_{\sigma=0} \ = \  \partial_{\sigma}\E^{\sigma,\s}[\partial_{\s}\mathscr{UP}(\s,\mathbf{U}(\y))]|_{\sigma=0} - \alpha(\mathscr{U},\s) \\
&= \ 0.
\end{align}
(The first identity follows from linearity of expectation and derivative. The second identity holds for the following reasons. First, the second term in the second line is zero, since it is the $\sigma$-derivative of something independent of $\sigma$. We also have $\E^{\sigma,\s}\mathbf{U}(\y)=\sigma$ by construction; see Definition \ref{definition:intro5}. This gives the second identity. The final identity follows from definition of $\alpha(\mathscr{U},\s)$; see after \eqref{eq:le8II10}. So the first $\{\}$-term in \eqref{eq:le8I12} is $\mathrm{LCT}$. Lemma \ref{lemma:ee1} for $\mathrm{j}=2$ gives the following (where $\mathbf{1}[\mathcal{E}(\s)]$ just means $\s\leq\t_{\mathrm{reg}}$):
\begin{align}
\mathbf{1}[\mathcal{E}(\s)]|\E^{\mathfrak{l}(\N),+}\{\partial_{\s}\mathscr{UP}(\s,\mathbf{U}(\y))-\E^{0,\s}\partial_{\s}\mathscr{UP}(\s,\cdot)-\alpha(\mathscr{U},\s)\mathbf{U}(\y)\}| \ &\lesssim \ \N^{10\gamma_{\mathrm{reg}}}\mathfrak{l}(\N)^{-1} \nonumber\\
&\lesssim \N^{-\frac14-99\gamma_{\mathrm{KL}}}. \label{eq:le8I14}
\end{align}
The last bound above follows by \eqref{eq:le8II11} and the fact that $\gamma_{\mathrm{reg}}$ is {at most a small constant times} $\gamma_{\mathrm{KL}}$ (see Definition \ref{definition:reg}). (Technically, we should also include the norm from Lemma \ref{lemma:ee1} of the $\{\}$-term in \eqref{eq:le8I14}. Because $|\partial_{\s}\mathscr{UP}(\s,\cdot)|\lesssim1$ as noted right after \eqref{eq:le8II6c}, this norm is $\lesssim1$ by the reasoning in the last paragraph of the proof of Lemmas \ref{lemma:bg23},\ref{lemma:bg1hl2}.) Via \eqref{eq:le8I13}, \eqref{eq:le8I14}, the $\E^{\mathfrak{l}(\N),+}$-term in the second integral of \eqref{eq:le8I9c} is $\lesssim\N^{-1}+\N^{-1/4-99\gamma_{\mathrm{KL}}}$. Note $\mathbf{1}(\mathcal{E}(\s))\leq1$ and $\mathfrak{p}(\s)$ is a probability density with respect to $\mathbb{P}^{0,\s,\mathbb{T}(\N)}$. So
\begin{align}
{\textstyle\int_{0}^{\t}}{\textstyle\sum_{\y}}\E^{0,\s,\mathbb{T}(\N)}\{\mathfrak{p}(\s)\mathbf{1}(\mathcal{E}(\s))\E^{\mathfrak{l}(\N),+}[\mathscr{CP}(\s,\cdot);\s,\y]\}\d\s \ &\lesssim \  \N^{-\frac14-99\gamma_{\mathrm{KL}}}|\mathbb{T}(\N)| + |\mathbb{T}(\N)||\mathbb{T}(\N)|^{-1} \nonumber\\
&\lesssim \  \N^{\frac34-99\gamma_{\mathrm{KL}}}. \label{eq:le8I14b}
\end{align}
where the last bound is because $\t\leq1$. By \eqref{eq:le8I9d}, this controls the second integral in \eqref{eq:le8I9c}. (We record the exact bound when relevant.) Let us now control the first integral in \eqref{eq:le8I9c}. We apply Lemma \ref{lemma:le8a} with the following choices. First, take $\kappa=1$. Take $\mathfrak{a}(\s,\mathbf{U}^{\s,\y+\cdot})=\mathscr{AP}(\s,\y;\mathfrak{l}(\N))$. Recalling $\mathscr{AP}(\s,\y;\mathfrak{l}(\N))$ from right after \eqref{eq:le8II11}, its support $\mathbb{I}$ (in the sense of Definition \ref{definition:le1}) is a discrete interval with length $\lesssim\mathfrak{l}(\N)^{2}$. Indeed, it is the average of $\mathscr{UP}(\s,\y+2\mathrm{k}\mathfrak{l}(\N);\mathfrak{l}(\N))$ whose supports are length $\mathfrak{l}(\N)$. (Indeed, see \eqref{eq:le8II11}, and see Definition \ref{definition:bg21}, from which it is clear that $\E^{\mathfrak{l},\pm}$-terms have support of length $\mathfrak{l}$). Moreover, we always take shift-indices $|\mathrm{k}|\lesssim\mathfrak{l}(\N)$. Therefore, Lemma \ref{lemma:le8a} gives
{
\begin{align}
&{\textstyle\int_{0}^{\t}}{\textstyle\sum_{\y}}\E^{0,\s,\mathbb{T}(\N)}\mathfrak{p}(\s)\mathbf{1}(\mathcal{E}(\s))\mathscr{AP}(\s,\y;\mathfrak{l}(\N))\d\s \nonumber\\
&\lesssim \ |\mathbb{I}|^{3}{\textstyle\int_{0}^{\t}}\mathfrak{D}_{\mathrm{FI}}^{0,\s}(\mathfrak{p}(\s))\d\s+|\mathbb{T}(\N)|\sup_{\substack{\sigma\in\R\\0\leq\s\leq\t}}\log\E^{\sigma,\s,\mathbb{I}}\exp\{|\mathfrak{a}(\s,\mathbf{U})|\}. \label{eq:le8I15}
\end{align}
}(We have multiplied everything in the upper bound coming from Lemma \ref{lemma:le8a} by $|\mathbb{T}(\N)|=\N$, because the LHS of \eqref{eq:le8I15} sums over $\mathbb{T}(\N)$; it does not average over $\mathbb{T}(\N)$.) We will now estimate the last term in \eqref{eq:le8I15}. Recall that we chose $\mathfrak{a}(\s,\mathbf{U}^{\s,\y+\cdot})=\mathscr{AP}(\s,\y;\mathfrak{l}(\N))$. By definition of the latter (see after \eqref{eq:le8II11}), we have that $\mathfrak{a}(\s,\mathbf{U})$ is an average of $\mathfrak{a}(\s,\mathbf{U};\mathrm{k})$ for $\mathrm{k}=0,\ldots,\mathfrak{l}(\N)$ with $\mathfrak{a}(\s,\mathbf{U}^{\s,\y+\cdot};\mathrm{k})=\mathscr{UP}(\s,\y+2\mathrm{k}\mathfrak{l}(\N);\mathfrak{l}(\N))$. We now claim the following about $\mathfrak{a}(\s,\mathbf{U};\mathrm{k})$. 
\begin{itemize}
\item First, the supports of $\mathfrak{a}(\s,\mathbf{U};\mathrm{k})$, which we denote with $\mathbb{I}(\mathrm{k})$, are mutually disjoint. (Indeed, by \eqref{eq:le8II11}, the support of $\mathfrak{a}(\s,\mathbf{U}^{\s,\y+\cdot};\mathrm{k})=\mathscr{UP}(\s,\y+2\mathrm{k}\mathfrak{l}(\N);\mathfrak{l}(\N))$ is a discrete interval length of $\mathfrak{l}(\N)$; see Definition \ref{definition:bg21}, from which it is clear that the support length of the $\E^{\mathfrak{l}(\N),+}$-term in \eqref{eq:le8II11} is $\mathfrak{l}(\N)$. It now suffices to note that we spatially shift $\mathfrak{a}(\s,\mathbf{U}^{\s,\y+\cdot};\mathrm{k})=\mathscr{UP}(\s,\y+2\mathrm{k}\mathfrak{l}(\N);\mathfrak{l}(\N))$ by multiples of $2\mathfrak{l}(\N)$, and that shifts of intervals of length $\mathfrak{l}(\N)$ by multiples of $2\mathfrak{l}(\N)$ are mutually disjoint.) 
\item Second, we claim $|\mathfrak{a}(\s,\mathbf{U};\mathrm{k})|\lesssim1$ with probability 1. (Indeed, this follows from $\mathfrak{a}(\s,\mathbf{U}^{\s,\y+\cdot};\mathrm{k})=\mathscr{UP}(\s,\y+2\mathrm{k}\mathfrak{l}(\N);\mathfrak{l}(\N))$, by \eqref{eq:le8II11}, and by $|\partial_{\s}\mathscr{UP}(\s,\cdot)|\lesssim1$; see right after \eqref{eq:le8II6c}.) 
\item Third, we claim that expectation of $\mathfrak{a}(\s,\mathbf{U};\mathrm{k})$ with respect to any canonical ensemble on its support $\mathbb{I}(\mathrm{k})$ is zero. (This follows from $\mathfrak{a}(\s,\mathbf{U}^{\s,\y+\cdot};\mathrm{k})=\mathscr{UP}(\s,\y+2\mathrm{k}\mathfrak{l}(\N);\mathfrak{l}(\N))$, by \eqref{eq:le8II11}, and by Lemma \ref{lemma:vanishcanonical}.) 
\end{itemize}
From these three claims, we deduce that $\mathfrak{a}(\s,\mathbf{U})$, which averages $\mathfrak{a}(\s,\mathbf{U};\mathrm{k})$ over $\mathrm{k}=0,\ldots,\mathfrak{l}(\N)$, is sub-Gaussian with variance $\lesssim\mathfrak{l}(\N)^{-1}$. (Indeed, we use Lemma \ref{lemma:le2} with our choice of $\mathbb{J}(\mathrm{k})$ as the support of $\mathfrak{a}(\s,\mathbf{U};\mathrm{k})$. The constraints in Lemma \ref{lemma:le2} are satisfied for the following reason. If one conditions on $\mathbf{U}(\x)$ for $\x\in\mathbb{J}(\mathrm{j})$ and $\mathrm{j}\neq\mathrm{k}$, the law of any canonical ensemble on $\mathbb{J}(1)\cup\ldots\cup\mathbb{J}(\mathfrak{l}(\N))$ projects to a canonical ensemble on $\mathbb{J}(\mathrm{k})$ when taking $\R^{\mathbb{J}(\mathrm{k})}$-marginals. This is just the statement that if one takes a random walk bridge and conditions one some increments, the remaining increments are distributed via random walk bridges. It now suffices to recall that $\mathfrak{a}(\s,\mathbf{U};\mathrm{k})$ vanishes with respect to \emph{any} canonical ensemble expectation on its support $\mathbb{J}(\mathrm{k})$.) Thus, by sub-Gaussianity of $\mathfrak{a}(\s,\mathbf{U})$,
\begin{align}
\log\E^{\sigma,\s,\mathbb{I}}\exp\{|\mathfrak{a}(\s,\mathbf{U})|\} \ \leq \ \log\exp\{\mathrm{O}(\mathfrak{l}(\N)^{-1})\} \ \lesssim \ \mathfrak{l}(\N)^{-1} \ = \ \N^{-\frac14-100\gamma_{\mathrm{KL}}}. \label{eq:le8I16}
\end{align}
This is uniform over all $\sigma$ and $0\leq\s\leq\t\leq1$ (everything is smooth in $\s$). By \eqref{eq:le8II5}, \eqref{eq:le8II10}, \eqref{eq:le8I9a}-\eqref{eq:le8I9c}, \eqref{eq:le8I9d}, \eqref{eq:le8I14b}, \eqref{eq:le8I15}, and \eqref{eq:le8I16}, we deduce the following estimate:
\begin{align}
\N^{2}{\textstyle\int_{0}^{\t}}\mathfrak{D}_{\mathrm{FI}}^{0,\s}(\mathfrak{p}(\s))\d\s \ &\lesssim \ \mathfrak{D}_{\mathrm{KL}}^{0,0}(\mathfrak{p}(0))+|\mathbb{I}|^{3}{\textstyle\int_{0}^{\t}}\mathfrak{D}_{\mathrm{FI}}^{0,\s}(\mathfrak{p}(\s))\d\s+\N^{\frac34-99\gamma_{\mathrm{KL}}}+|\mathbb{T}(\N)|\N^{-\frac14-100\gamma_{\mathrm{KL}}} \\
&\lesssim \ \mathfrak{D}_{\mathrm{KL}}^{0,0}(\mathfrak{p}(0))+|\mathbb{I}|^{3}{\textstyle\int_{0}^{\t}}\mathfrak{D}_{\mathrm{FI}}^{0,\s}(\mathfrak{p}(\s))\d\s+\N^{\frac34-99\gamma_{\mathrm{KL}}}. \label{eq:le8I17}
\end{align}
Recall that $|\mathbb{I}|\lesssim\mathfrak{l}(\N)^{2}=\N^{1/2+200\gamma_{\mathrm{KL}}}$; see \eqref{eq:le8II11} and the paragraph after \eqref{eq:le8I14b}. Thus, $|\mathbb{I}|^{3}\ll\N^{2}$. This lets us move the first term in \eqref{eq:le8I17} to the LHS of the first line of the previous display and deduce the following estimate: 
\begin{align}
\{1-\mathrm{o}(1)\}\N^{2}{\textstyle\int_{0}^{\t}}\mathfrak{D}_{\mathrm{FI}}^{0,\s}(\mathfrak{p}(\s))\d\s\lesssim\mathfrak{D}_{\mathrm{KL}}^{0,0}(\mathfrak{p}(0))+\N^{\frac34-99\gamma_{\mathrm{KL}}}. \label{eq:le8I18} 
\end{align}
\eqref{eq:le8I} follows as $\mathfrak{p}(0)$ is entropy data (see Definition \ref{definition:entropydata}). Observe that this argument uses $\t_{\mathrm{reg}}=1$ with very high probability, and not just with high probability. We only have the latter in the case where $\mathscr{U}(\t,\cdot)$ is independent of $\t$ (see Theorem \ref{theorem:kpz}). Thus suppose that $\mathscr{U}(\t,\cdot)$ is independent of $\t$ and avoid using $\t_{\mathrm{reg}}=1$ with very high probability. In this case, \eqref{eq:le8II5} still holds. But the integral over $[0,\t]$ therein is zero as $\mathscr{U}(\t,\cdot)$ is independent of $\t$. (Indeed, by \eqref{eq:le8II3b} and right before it, if $\mathscr{U}$ is independent of $\t$, so is the $\mathscr{HP}$-term in \eqref{eq:le8II5}.) So \eqref{eq:le8I18} follows directly and gives \eqref{eq:le8I}. This finishes the proof.
\end{proof}
We now use Lemma \ref{lemma:le8} to upgrade Lemma \ref{lemma:le7}.
\begin{lemma}\label{lemma:le9}
 Retain the settings of {Lemmas \ref{lemma:le8a}, \ref{lemma:le8}}. We have the following for any $\kappa>0$:
\begin{align}
\mathrm{LHS}\eqref{eq:le8aI} \ \lesssim \ \tfrac{1}{\kappa}\N^{-\frac94-\gamma_{\mathrm{KL}}}|\mathbb{I}|^{3}+\tfrac{1}{\kappa}\sup_{\sigma\in\R}\sup_{0\leq\s\leq\t}\log\E^{\sigma,\s,\mathbb{I}}\exp\{\kappa|\mathfrak{a}(\s,\mathbf{U})|\}. \label{eq:le9I}
\end{align}
\end{lemma}
\begin{proof}
Combine \eqref{eq:le8aI} and \eqref{eq:le8I}.
\end{proof}
\subsection{Localizing \eqref{eq:hf}-\eqref{eq:glsde} via speed of propagation}
The reduction-to-local-equilibrium estimate \eqref{eq:le9I} deteriorates quite badly as the support-length $|\mathbb{I}|$ grows. In particular, it is effective for \emph{local} statistics. We eventually want to use \eqref{eq:le9I} to estimate averages of local fluctuations (like $\mathds{R}^{\chi,\mathfrak{q},\pm,\mathrm{j}}$-terms in Definition \ref{definition:bg26}) in space \emph{and time}. But time-averages are, in principle, global statistics. Indeed, the evolution \eqref{eq:glsde} is basically a heat flow on $\mathbb{T}(\N)$ with diffusion-speed $\N^{2}$ and asymmetry-speed $\N^{3/2}$. But, modulo exponentially small errors, such heat flows propagate only distance $\approx\N\t^{1/2}+\N^{3/2}\t$ in space by time $\t$, up to arbitrarily small powers of $\N$. So, for mesoscopic $\t$, time-averages are also basically local statistics. The point of this subsection is to make this random walk heuristic precise (basically by linearizing \eqref{eq:glsde} to get a discrete parabolic equation, with a good random walk kernel, on $\mathbb{T}(\N)$ and localizations to sub-intervals $\mathbb{I}\subseteq\mathbb{T}(\N)$). In what follows, see the end of Section \ref{section:main} for what $\inf$ and $\sup$ of a discrete interval means.
\begin{definition}\label{definition:le10}
 Take any $\mathfrak{t}\geq0$ and discrete interval $\mathbb{I}\subseteq\mathbb{T}(\N)$. Define $\mathfrak{l}(\mathfrak{t},\mathbb{I})=\lfloor\N^{\gamma_{\mathrm{ap}}}\{\N\mathfrak{t}^{1/2}+\N^{3/2}\mathfrak{t}+|\mathbb{I}|\}\rfloor$. We now define $\mathbb{I}(\mathfrak{t}):=\llbracket\inf\mathbb{I}-\mathfrak{l}(\mathfrak{t},\mathbb{I}),\sup\mathbb{I}+\mathfrak{l}(\mathfrak{t},\mathbb{I})\rrbracket$. Let $\mathbf{U}^{\t,\cdot}[\mathbb{I}(\mathfrak{t})]$ solve the following for $\t\geq0$ and $\x\in\mathbb{I}(\mathfrak{t})$ (with notation explained after):
\begin{align}
\d\mathbf{U}^{\t,\x}[\mathbb{I}(\mathfrak{t})] = \N^{2}\Delta^{\mathbb{I}(\mathfrak{t})}\mathscr{U}'(\t,\mathbf{U}^{\t,\x}[\mathbb{I}(\mathfrak{t})])\d\t+\N^{\frac32}\grad^{\mathbb{I}(\mathfrak{t}),\mathrm{a}}\mathscr{U}'(\t,\mathbf{U}^{\t,\x}[\mathbb{I}(\mathfrak{t})])\d\t - \sqrt{2}\N\grad^{\mathbb{I}(\mathfrak{t}),-}\d\mathbf{b}(\t,\x). \label{eq:glsdeloc}
\end{align}
Recall $\grad^{\mathbb{I}(\mathfrak{t}),?}$ from Definition \ref{definition:le5}. We define $\Delta^{\mathbb{I}(\mathfrak{t})}=\grad^{\mathbb{I}(\mathfrak{t}),+}+\grad^{\mathbb{I}(\mathfrak{t}),-}$. (Thus, it is just $\Delta$ from Definition \ref{definition:intro4} but with respect to periodic boundary data on $\mathbb{I}(\mathfrak{t})$.) We write initial data (as an element in $\R^{\mathbb{I}(\mathfrak{t})}$) for \eqref{eq:glsdeloc} on a case-by-case basis. Now, define $\mathbf{J}(\t,\x;\mathbb{I}(\mathfrak{t}))$ by the gradient relation $-\N^{1/2}\grad^{\mathbb{I}(\mathfrak{t}),-}\mathbf{J}(\t,\x;\mathbb{I}(\mathfrak{t})):=\mathbf{U}^{\t,\x}[\mathbb{I}(\mathfrak{t})]$ and $\mathbf{J}(\t,\inf\mathbb{I}(\mathfrak{t});\mathbb{I}(\mathfrak{t})):=\Pi^{\mathbb{S}(\N)}\mathbf{J}(\t;\mathbb{I}(\mathfrak{t}))$, in which $\Pi^{\mathbb{S}(\N)}$ is the canonical quotient map $\R\to\mathbb{S}(\N)\simeq[-\N^{20\gamma_{\mathrm{reg}}},\N^{20\gamma_{\mathrm{reg}}}]$ (see Definition \ref{definition:le5}), and where $\mathbf{J}(\t;\mathbb{I}(\mathfrak{t}))$ solves
\begin{align}
\d\mathbf{J}(\t;\mathbb{I}(\mathfrak{t})) \ := \ &\N^{\frac32}\grad^{\mathbb{I}(\mathfrak{t}),+}\mathscr{U}'(\t,\mathbf{U}^{\t,\inf\mathbb{I}(\mathfrak{t})}[\mathbb{I}(\mathfrak{t})])\d\t-\mathscr{R}(\t)\d\t+\sqrt{2}\N^{\frac12}\d\mathbf{b}(\t,\inf\mathbb{I}(\mathfrak{t})) \label{eq:hfloc}\\
+ \ &\N\{\mathscr{U}'(\t,\mathbf{U}^{\t,\inf\mathbb{I}(\mathfrak{t})}[\mathbb{I}(\mathfrak{t})])+\mathscr{U}'(\t,\mathbf{U}^{\t,\inf\mathbb{I}(\mathfrak{t})+1}[\mathbb{I}(\mathfrak{t})])\}\d\t.\nonumber
\end{align}
(We used $\mathscr{R}(\t)$ from Definition \ref{definition:intro6}. Again, we specify initial data for \eqref{eq:hfloc} on a case-by-case basis when relevant.)
\end{definition}
For now, let us defer discussing the motivation behind constructing $\mathbf{J}(\t,\x;\mathbb{I}(\mathfrak{t}))$ the way that we did in Definition \ref{definition:le10}. (This motivation is given after Definition \ref{definition:le14}.) First, we make precise the heuristic before Definition \ref{definition:le10}. (Namely, we compare \eqref{eq:glsde}, \eqref{eq:glsdeloc} on local space-time scales.)
\begin{lemma}\label{lemma:le12}
 Fix a discrete interval $\mathbb{I}\subseteq\mathbb{T}(\N)$ and $\mathfrak{t}\lesssim1$. Suppose $\mathbf{U}^{\t,\cdot}$ solves \eqref{eq:glsde}, and suppose $\mathbf{U}^{\t,\cdot}[\mathbb{I}(\mathfrak{t})]$ solves \eqref{eq:glsdeloc}. Assume that the Brownian motions $\mathbf{b}(\t,\x)$ (for $\x\in\mathbb{I}(\mathfrak{t})$) are the same for \eqref{eq:glsde}, \eqref{eq:glsdeloc}. (Brownian motions $\mathbf{b}(\t,\x)$ for $\x\not\in\mathbb{I}(\mathfrak{t})$ are chosen independently.) Now, fix any $\mathfrak{t}(\mathrm{in})\geq0$ and assume $\mathbf{U}^{\mathfrak{t}(\mathrm{in}),\x}=\mathbf{U}^{\mathfrak{t}(\mathrm{in}),\x}[\mathbb{I}(\mathfrak{t})]$ for all $\x\in\mathbb{I}(\mathfrak{t})$. Assume $|\mathbf{U}^{\mathfrak{t}(\mathrm{in}),\x}|\lesssim\N$ for all $\x\in\mathbb{T}(\N)$. Lastly, define $\mathbb{I}(+):=\llbracket\inf\mathbb{I}-10,\sup\mathbb{I}+10\rrbracket$. With very high probability, {we have the following for any large but fixed $D>0$:}
\begin{align}
\sup_{\mathfrak{t}(\mathrm{in})\leq\t\leq\mathfrak{t}(\mathrm{in})+\mathfrak{t}}\sup_{\x\in\mathbb{I}(+)}|\mathbf{U}^{\t,\x}-\mathbf{U}^{\t,\x}[\mathbb{I}(\mathfrak{t})]| \ \lesssim \ \N^{-D}. \label{eq:le12I}
\end{align}
\end{lemma}
\begin{rem}\label{eq:le13}
 Our proof for Lemma \ref{lemma:le12} will use convexity of $\mathscr{U}$; see Assumption \ref{ass:intro8}. But a totally adequate substitute for Lemma \ref{lemma:le12}, which is a total variation bound as opposed to Wasserstein-type bound, can be obtained via heat kernel estimates for SDEs \eqref{eq:glsde}, \eqref{eq:glsdeloc}. Said heat kernel estimates would give control for speed of propagation (in space) for the SDEs. Proving them requires only log-Sobolev inequality \eqref{eq:le4I}. (This log-Sobolev method is, by now, standard. See \cite{Davies} for a general picture, and see \cite{L} for a specialization to speed of propagation for interacting particle systems.) The log-Sobolev method, however, is much more complicated than the following proof via convexity, as evidenced by \cite{L}. Given the length of this paper, we give the proof via convexity. (We make this remark to show log-Sobolev is enough. Convexity only makes things easy to write.)
\end{rem}
\begin{proof}
First, a preliminary estimate. Set $\tau\eqref{eq:glsde}$ to be the first time $\mathfrak{t}(\mathrm{in})\leq\t\leq\mathfrak{t}(\mathrm{in})+\mathfrak{t}$ such that the supremum of $|\mathbf{U}^{\t,\x}|$ over $\x\in\mathbb{T}(\N)$ exceeds $\N^{{\mathrm{D}}}$. Let $\tau\eqref{eq:glsdeloc}$ be the same but for $\mathbf{U}^{\t,\x}[\mathbb{I}(\mathfrak{t})]$ in place of $\mathbf{U}^{\t,\x}$ and $\mathbb{I}(\mathfrak{t})$ in place of $\mathbb{T}(\N)$. We claim
\begin{align}
\mathbb{P}[\tau\eqref{eq:glsde}\neq\mathfrak{t}(\mathrm{in})+\mathfrak{t}]+\mathbb{P}[\tau\eqref{eq:glsdeloc}\neq\mathfrak{t}(\mathrm{in})+\mathfrak{t}] \ \lesssim \  \exp[-\N]. \label{eq:le12I1}
\end{align}
Take \eqref{eq:le12I1} for now. Set $\mathbb{I}(\circ):=\llbracket\inf\mathbb{I}-2^{-1}\mathfrak{l}(\mathfrak{t},\mathbb{I}),\sup\mathbb{I}+2^{-1}\mathfrak{l}(\mathfrak{t},\mathbb{I})\rrbracket\supseteq\mathbb{I}(+)$. Set $\mathbf{D}^{\t,\x}:=\mathbf{U}^{\t,\x}-\mathbf{U}^{\t,\x}[\mathbb{I}(\mathfrak{t})]$. The Brownian motions in \eqref{eq:glsde} and \eqref{eq:glsdeloc} are equal for $\x$ of distance 1 from $\mathbb{I}(\circ)$. (Indeed, $\mathbb{I}(\circ)$ is a subset of $\mathbb{I}(\mathfrak{t})$ that is $2^{-1}\mathfrak{l}(\mathfrak{t},\mathbb{I})$-away from the boundary of $\mathbb{I}(\mathfrak{t})$. It now suffices to recall that the Brownian motions are the same for all $\x\in\mathbb{I}(\mathfrak{t})$ and to note $2^{-1}\mathfrak{l}(\mathfrak{t},\mathbb{I})\gg1$; see Definition \ref{definition:le10}.) We claim that if $\x\in\mathbb{I}(\circ)$, then the following SDE holds for some process $\mathbf{M}^{\t,\x}$:
\begin{align}
&\d\mathbf{D}^{\t,\x} \nonumber\\
&= \ \N^{2}\Delta^{\mathbb{I}(\mathfrak{t})}\{\mathscr{U}'(\t,\mathbf{U}^{\t,\x})-\mathscr{U}'(\t,\mathbf{U}^{\t,\x}[\mathbb{I}(\mathfrak{t})])\}\d\t + \N^{\frac32}\grad^{\mathbb{I}(\mathfrak{t}),\mathrm{a}}\{\mathscr{U}'(\t,\mathbf{U}^{\t,\x})-\mathscr{U}'(\t,\mathbf{U}^{\t,\x}[\mathbb{I}(\mathfrak{t})])\}\d\t \label{eq:le12I2a} \\
&= \ \N^{2}\Delta^{\mathbb{I}(\mathfrak{t})}\{\mathscr{U}''(\t,\mathbf{M}^{\t,\x})\mathbf{D}^{\t,\x}\}\d\t + \N^{\frac32}\grad^{\mathbb{I}(\mathfrak{t}),\mathrm{a}}\{\mathscr{U}''(\t,\mathbf{M}^{\t,\x})\mathbf{D}^{\t,\x}\}\d\t. \label{eq:le12I2b}
\end{align}
Indeed, the Brownian motions in $\mathbf{U}^{\t,\x}-\mathbf{U}^{\t,\x}[\mathbb{I}(\mathfrak{t})]$ cancel if $\x\in\mathbb{I}(\circ)$ as explained in the previous paragraph. Also, for $\x\in\mathbb{I}(\circ)$, which is of distance more than 1 from the boundary of $\mathbb{I}(\mathfrak{t})$, the gradients $\grad^{\mathbb{K},\mathrm{a}}$ for $\mathbb{K}=\mathbb{T}(\N),\mathbb{I}(\mathfrak{t})$ are the same when acting on functions of $\mathbf{U}^{\t,\x}$. (For $\x\in\mathbb{I}(\circ)$, we also get $\Delta^{\mathbb{I}(\mathfrak{t})}=\Delta$ when acting on functions of $\mathbf{U}^{\t,\x}$.) This proves \eqref{eq:le12I2a}. \eqref{eq:le12I2b} follows by the mean-value theorem. Now, set $\mathbf{I}^{\t,\x}=\mathbf{D}^{\t,\x}\mathbf{1}[\x\in\mathbb{I}(\circ)]=:\mathbf{D}^{\t,\x}\mathbf{1}[\mathbb{I}(\circ)]$. We claim, with notation explained after,
\begin{align}
\d\mathbf{I}^{\t,\x} \ &= \ \mathbf{1}[\mathbb{I}(\circ)]\times\N^{2}\Delta^{\mathbb{I}(\mathfrak{t})}\{\mathscr{U}''(\t,\mathbf{M}^{\t,\x})\mathbf{D}^{\t,\x}\}\d\t \label{eq:le12I3a}\\
&+ \  \mathbf{1}[\mathbb{I}(\circ)]\times\N^{\frac32}\grad^{\mathbb{I}(\mathfrak{t}),\mathrm{a}}\{\mathscr{U}''(\t,\mathbf{M}^{\t,\x})\mathbf{D}^{\t,\x}\}\d\t\nonumber\\
&= \ \N^{2}\Delta^{\mathbb{I}(\mathfrak{t})}\{\mathscr{U}''(\t,\mathbf{M}^{\t,\x})\mathbf{I}^{\t,\x}\}\d\t+\N^{\frac32}\grad^{\mathbb{I}(\mathfrak{t}),\mathrm{a}}\{\mathscr{U}''(\t,\mathbf{M}^{\t,\x})\mathbf{I}^{\t,\x}\}\d\t \label{eq:le12I3b}\\
&+ \ [\mathbf{1}[\mathbb{I}(\circ)],\mathscr{O}]\{\mathscr{U}''(\t,\mathbf{M}^{\t,\x})\mathbf{D}^{\t,\x}\}\d\t. \nonumber
\end{align}
\eqref{eq:le12I3a} is by $\d\mathbf{I}^{\t,\x}=\mathbf{1}[\mathbb{I}(\circ)]\times\d\mathbf{D}^{\t,\x}$ and \eqref{eq:le12I2a}-\eqref{eq:le12I2b}. In \eqref{eq:le12I3b}, $\mathscr{O}:=\N^{2}\Delta^{\mathbb{I}(\mathfrak{t})}+\N^{3/2}\grad^{\mathbb{I}(\mathfrak{t}),\mathrm{a}}$. Also, $[,]$ is the commutator for operators, where $\mathbf{1}[\mathbb{I}(\circ)]$ is identified with the operator given by multiplication by $\mathbf{1}[\mathbb{I}(\circ)]$. Thus, to get \eqref{eq:le12I3b}, we first move $\mathbf{1}[\mathbb{I}(\circ)]$ in $\Delta^{\mathbb{I}(\mathfrak{t})}$ and $\grad^{\mathbb{I}(\mathfrak{t}),\mathrm{a}}$ in \eqref{eq:le12I3a}. This turns $\mathbf{D}^{\t,\x}$ into $\mathbf{I}^{\t,\x}$, giving the first two terms in \eqref{eq:le12I3b}. The cost is the commutator $[\mathbf{1}[\mathbb{I}(\circ)],\N^{2}\Delta^{\mathbb{I}(\mathfrak{t})}]+[\mathbf{1}[\mathbb{I}(\circ)],\N^{3/2}\grad^{\mathbb{I}(\mathfrak{t}),\mathrm{a}}]=[\mathbf{1}[\mathbb{I}(\circ)],\mathscr{O}]$. Now, we let $\Gamma[\s,\t,\x,\y]$ solve the PDE $\partial_{\t}\Gamma[\s,\t,\x,\y]=\mathscr{Y}\Gamma[\s,\t,\x,\y]$ and $\Gamma[\s,\s,\x,\y]=\mathbf{1}[\x=\y]$. Here, we defined $\mathscr{Y}\phi(\x):=\N^{2}\Delta^{\mathbb{I}(\mathfrak{t})}\{\mathscr{U}''(\t,\mathbf{M}^{\t,\x})\phi(\x)\}+\N^{3/2}\grad^{\mathbb{I}(\mathfrak{t}),\mathrm{a}}\{\mathscr{U}''(\t,\mathbf{M}^{\t,\x})\phi(\x)\}$ that acts on $\Gamma[\s,\t,\x,\y]$ in the $\x$-variable. ($\mathscr{Y}$ and this PDE are on the space-time $[0,\infty)\times\mathbb{I}(\mathfrak{t})$.) By \eqref{eq:le12I3a}-\eqref{eq:le12I3b} and Duhamel, 
\begin{align}
\mathbf{I}^{\t,\x} \ = \ \int_{\mathfrak{t}(\mathrm{in})}^{\t}\sum_{\y\in\mathbb{I}(\mathfrak{t})}\Gamma[\s,\t,\x,\y]\cdot[\mathbf{1}[\mathbb{I}(\circ)],\mathscr{O}]\{\mathscr{U}''(\s,\mathbf{M}^{\s,\y})\mathbf{D}^{\s,\y}\}\d\s. \label{eq:le12I3c}
\end{align}
(Technically, \eqref{eq:le12I3c} should have an initial data term given by integrating-in-space $\Gamma[0,\t,\x,\y]\mathbf{I}^{\mathfrak{t}(\mathrm{in}),\y}$. However, $\mathbf{I}^{\mathfrak{t}(\mathrm{in}),\y}=0$ by construction.) We now claim the inequality $|[\mathbf{1}[\mathbb{I}(\circ)],\mathscr{O}]\{\mathscr{U}''(\s,\mathbf{M}^{\s,\y})\mathbf{D}^{\s,\y}\}|\lesssim\N^{{\mathrm{O}(1)}}\mathbf{1}[\mathrm{dist}(\y,\mathbb{I}(+))\gtrsim\mathfrak{l}(\mathfrak{t},\mathbb{I})]$ for $\s\leq\tau\eqref{eq:glsde}\wedge\tau\eqref{eq:glsdeloc}$. (Recall these two stopping times from the beginning of this proof.) To prove this bound, {we perform} two steps. First, by Assumption \ref{ass:intro8}, we know $\mathscr{U}''=\mathrm{O}(1)$. We also know that $|\mathbf{D}^{\s,\y}|\leq|\mathbf{U}^{\s,\y}|+|\mathbf{U}^{\s,\y}[\mathbb{I}(\mathfrak{t})]|\lesssim\N^{{\mathrm{O}(1)}}$ by triangle inequality and then the definitions of $\tau\eqref{eq:glsde},\tau\eqref{eq:glsdeloc}$. (This is if $\s\leq\tau\eqref{eq:glsde}\wedge\tau\eqref{eq:glsdeloc}$.) Next, we study the commutator $[\mathbf{1}[\mathbb{I}(\circ)],\mathscr{O}]\phi=\mathbf{1}[\mathbb{I}(\circ)]\mathscr{O}\phi-\mathscr{O}\{\mathbf{1}[\mathbb{I}(\circ)]\phi\}$ for any $\phi:\mathbb{I}(\mathfrak{t})\to\R$. By the discrete Leibniz rule (see the proof of Lemma \ref{lemma:method22}), an explicit calculations implies $\mathscr{O}\{\mathbf{1}[\mathbb{I}(\circ)]\phi\}=\mathbf{1}[\mathbb{I}(\circ)]\mathscr{O}\phi+\mathrm{error}$, in which $\mathrm{error}$ is given by a linear combination of discrete gradients of $\mathbf{1}[\mathbb{I}(\circ)]$. Thus, $\mathrm{error}$ is supported in a neighborhood of radius $10$ centered at the boundary points of $\mathbb{I}(\circ)$. But the distance between the boundary of $\mathbb{I}(\circ)$ and $\mathbb{I}(+)$ is $\gtrsim\mathfrak{l}(\mathfrak{t},\mathbb{I})$. (Indeed, $\mathbb{I}(\circ)$ is a radius $2^{-1}\mathfrak{l}(\mathfrak{t},\mathbb{I})$ neighborhood of $\mathbb{I}$, and $\mathbb{I}(+)$ is a radius 10 neighborhood of $\mathbb{I}$.) This proves the claim in the second line of this paragraph. Now, if we evaluate \eqref{eq:le12I3c} for any $\x\in\mathbb{I}(+)$, because $\tau\eqref{eq:glsde}\wedge\tau\eqref{eq:glsdeloc}=\mathfrak{t}(\mathrm{in})+\mathfrak{t}$ with very high probability, we get the following with very high probability for $\mathfrak{t}(\mathrm{in})\leq\t\leq\mathfrak{t}(\mathrm{in})+\mathfrak{t}$:
\begin{align}
|\mathbf{I}^{\t,\x}| \ \lesssim \ \mathfrak{t}|\mathbb{I}(\mathfrak{t})|\N^{{\mathrm{O}(1)}}\sup_{\mathfrak{t}(\mathrm{in})\leq\s\leq\t}\sup_{|\y-\x|\gtrsim\mathfrak{l}(\mathfrak{t},\mathbb{I})}|\Gamma[\s,\t,\x,\y]|. \label{eq:le12I4a}
\end{align}
Now, recall $\partial_{\t}\Gamma[\s,\t,\x,\y]=\mathscr{Y}\Gamma[\s,\t,\x,\y]$ and $\Gamma[\s,\s,\x,\y]=\mathbf{1}[\x=\y]$. Via duality between Kolmogorov forward and backward PDEs, we get $\partial_{\s}\Gamma[\s,\t,\y,\x]=-\mathscr{Y}^{\ast}\Gamma[\s,\t,\y,\x]$, where $\mathscr{Y}^{\ast}$ is adjoint with respect to uniform measure on $\mathbb{I}(\mathfrak{t})$. (It acts on $\y$. The additional negative sign is because going forward in $\s$ means going backwards in the time-parameter for the adjoint process.) As $1\lesssim\mathscr{U}''\lesssim1$, $\mathscr{Y}^{\ast}$ is the infinitesimal generator for a random walk whose symmetric jump speed is $\lesssim\N^{2}$ and whose asymmetric jump rate is $\lesssim\N^{3/2}$. (Indeed, recall $\mathscr{Y}$ from right before \eqref{eq:le12I3c}. To compute its adjoint, it is enough to know $\Delta^{\mathbb{I}(\mathfrak{t})}$ is symmetric with respect to the uniform measure $\mathbb{I}(\mathfrak{t})$, and $\grad^{\mathbb{I}(\mathfrak{t}),\mathrm{a}}$ is asymmetric. This can be checked by noting $[\grad^{\mathbb{I}(\mathfrak{t}),\mathfrak{l}}]^{\ast}=\grad^{\mathbb{I}(\mathfrak{t}),-\mathfrak{l}}$, which is a discrete integration-by-parts, and construction of $\Delta^{\mathbb{I}(\mathfrak{t})}$ and $\grad^{\mathbb{I}(\mathfrak{t}),\mathrm{a}}$ in Definitions \ref{definition:le5} and \ref{definition:le10}.) In particular, $\Gamma[\s,\t,\y,\x]$ is the transition probability for said random walk. Standard concentration bounds for random walks (like Azuma) then give that $\Gamma[\s,\t,\y,\x]\lesssim\exp[-\N^{\gamma_{\mathrm{ap}}}]$ if $|\y-\x|\geq\N^{\gamma_{\mathrm{ap}}}\{\N[\t-\s]^{1/2}+\N^{3/2}[\t-\s]\}$. By construction of $\mathfrak{l}(\mathfrak{t},\mathbb{I})$ in Definition \ref{definition:le10}, this is always the case if $|\y-\x|\gtrsim\mathfrak{l}(\mathfrak{t},\mathbb{I})$ and $\mathfrak{t}(\mathrm{in})\leq\s\leq\t\leq\mathfrak{t}(\mathrm{in})+\mathfrak{t}$. Since we are taking a sup over a set of all pairs $(\y,\x)$ defined by a constraint that is invariant under swapping $\y,\x$, we can swap $\Gamma[\s,\t,\x,\y]$ with $\Gamma[\s,\t,\y,\x]$ on the RHS of \eqref{eq:le12I4a}. So
\begin{align}
|\mathbf{I}^{\t,\x}| \ \lesssim \ \mathfrak{t}|\mathbb{I}(\mathfrak{t})|\N^{{\mathrm{O}(1)}}\sup_{\mathfrak{t}(\mathrm{in})\leq\s\leq\t}\sup_{|\y-\x|\gtrsim\mathfrak{l}(\mathfrak{t},\mathbb{I})}|\Gamma[\s,\t,\y,\x]| \lesssim \ \mathfrak{t}|\mathbb{I}(\mathfrak{t})|\N^{{\mathrm{O}(1)}}\exp[-\N^{\gamma_{\mathrm{ap}}}] \ \lesssim_{\mathrm{D}} \ \N^{-\mathrm{D}}, \label{eq:le12I4b}
\end{align}
where the last bound follows because $\mathfrak{t},|\mathbb{I}(\mathfrak{t})|\lesssim\N$ (by assumption and by $\mathbb{I}(\mathfrak{t})\subseteq\mathbb{T}(\N)$) and because polynomials are always beat by $\exp[-\N^{\gamma_{\mathrm{ap}}}]$. We clarify \eqref{eq:le12I4b} holds with very high probability simultaneously for $\x\in\mathbb{I}(\mathfrak{t})$. \eqref{eq:le12I} now follows since $\mathbf{I}^{\t,\x}=\mathbf{U}^{\t,\x}-\mathbf{U}^{\t,\x}[\mathbb{I}(\mathfrak{t})]$ for all $\x\in\mathbb{I}(\circ)\supseteq\mathbb{I}(+)$ (see right before \eqref{eq:le12I3a}). It now remains to prove the estimate \eqref{eq:le12I1} that we took for granted. We prove \eqref{eq:le12I1} for $\tau\eqref{eq:glsde}$; the proof for $\tau\eqref{eq:glsdeloc}$ is the same. (Just replace $\mathbb{T}(\N)$ by $\mathbb{I}(\mathfrak{t})$.) By using Lemma \ref{lemma:ste}, we claim the following. It reduces to estimating a probability for a supremum over a finite set that we can apply union bound to (if we have exponential-scale estimates). Below, we {let $\mathrm{C}>0$ be a large but fixed constant, and we} set $\mathbb{X}:=\{[\mathfrak{t}(\mathrm{in}),\mathfrak{t}(\mathrm{in})+\mathfrak{t}]\cap\N^{-\mathrm{D}}\Z\}\cup\{\mathfrak{t}(\mathrm{in})\}$ (for large, fixed $\mathrm{D}\geq0$):
\begin{align}
\mathbb{P}[\tau\eqref{eq:glsde}\neq\mathfrak{t}(\mathrm{in})+\mathfrak{t}] \ \lesssim \ \mathbb{P}\left\{\sup_{\t\in\mathbb{X}}\sup_{\x\in\mathbb{T}(\N)}|\mathbf{U}^{\t,\x}|\gtrsim\N^{{\mathrm{C}}}\right\} + \exp[-\N]. \label{eq:le12I5}
\end{align}
Indeed, Lemma \ref{lemma:ste} says that on very short time-increments, \eqref{eq:glsde} evolves very little. Precisely, since the coefficients in \eqref{eq:glsde} are uniformly Lipschitz in the solution by Assumption \ref{ass:intro8}, $\mathbf{U}^{\t,\x}$ is controlled by its values over the very fine time-discretization $\mathbb{X}$. (We clarify the LHS of \eqref{eq:le12I5} is the same as the first term on the RHS, but with $\mathfrak{t}(\mathrm{in})+[0,\mathfrak{t}]$ instead of its fine discretization $\mathbb{X}$.) Now, let $\mathbf{A}^{\t,\cdot}$ solve the following time-homogeneous SDE with the same initial data $\mathbf{A}^{\mathfrak{t}(\mathrm{in}),\cdot}=\mathbf{U}^{\mathfrak{t}(\mathrm{in}),\cdot}$
\begin{align}
\d\mathbf{A}^{\t,\x} \ = \ \N^{2}\Delta\mathscr{U}'(\mathfrak{t}(\mathrm{in}),\mathbf{A}^{\t,\x})\d\t + \N^{2}\grad^{\mathrm{a}}\mathscr{U}'(\mathfrak{t}(\mathrm{in}),\mathbf{A}^{\t,\x})\d\t - \sqrt{2}\N\grad^{-}\d\mathbf{b}(\t,\x). \label{eq:le12I6}
\end{align}
(The benefit of \eqref{eq:le12I6} is that it has invariant measures. We use the invariant measures to compute one-point statistics explicitly for \eqref{eq:le12I6} with stationary initial data. We then compare to stationary data via the parabolic maximum principle that is satisfied by a linearization of \eqref{eq:le12I6}. This is like how \eqref{eq:le12I3c} comes from the linearization of an SDE with analogous features as \eqref{eq:le12I6}.) By the Girsanov theorem, the Radon-Nikodym derivative of the process \eqref{eq:glsde} with respect to the process \eqref{eq:le12I6}, both restricted to times in $\mathfrak{t}(\mathrm{in})+[0,\mathfrak{t}]$, is the exponential martingale $\exp\{\mathscr{N}(\mathfrak{t}(\mathrm{in}),\mathfrak{t})+\mathscr{Q}(\mathfrak{t}(\mathrm{in}),\mathfrak{t})\}$. Here, $\mathscr{N}$ is a martingale, and $\mathscr{Q}$ satisfies the following estimate, in which we view $\mathscr{Q}$ as a function of an arbitrary continuous path $\mathbf{x}(\s,\z)$ for $(\s,\z)\in[\mathfrak{t}(\mathrm{in}),\mathfrak{t}(\mathrm{in})+\mathfrak{t}]\times\mathbb{I}(\mathfrak{t})$ for the display below:
{\small
\begin{align}
&|\mathscr{Q}(\mathfrak{t}(\mathrm{in}),\mathfrak{t})| \nonumber\\
&\lesssim \ \N^{2}{\textstyle\int_{\mathfrak{t}(\mathrm{in})}^{\mathfrak{t}(\mathrm{in})+\mathfrak{t}}}{\textstyle\sum_{\y}}|[\grad^{-}]^{-1}\{\Delta[\mathscr{U}'(\mathfrak{t}(\mathrm{in}),\mathbf{x}(\s,\z))-\mathscr{U}'(\s,\mathbf{x}(\s,\z))]+\grad^{\mathrm{a}}[\mathscr{U}'(\mathfrak{t}(\mathrm{in}),\mathbf{x}(\s,\z))-\mathscr{U}'(\s,\mathbf{x}(\s,\z))]\}|^{2}\d\s \nonumber \\
&\lesssim \ \N^{2}|\mathbb{T}(\N)|\mathfrak{t} \ \lesssim \ \N^{3}. \label{eq:le12I7}
\end{align}
}The first line follows by Girsanov. In particular, to remove a drift (at the level of laws of SDEs), one has to apply the inverse of the operator hitting the Brownian motion to the drift that is being removed, square it per $\y\in\mathbb{T}(\N)$, and then integrate in time. (The drift we are removing is the error obtained when we freeze the potential in \eqref{eq:glsde} at time $\mathfrak{t}(\mathrm{in})$, thereby giving \eqref{eq:le12I6}. Also, the factor $\N^{2}$ appearing in both lines above comes from the fact that the drift we are removing is $\mathrm{O}(\N^{2})$, and the operator hitting the Brownian motion in \eqref{eq:glsde}, \eqref{eq:le12I6} is $\mathrm{O}(\N)$. Taking the inverse of $\mathrm{O}(\N)$, hitting $\mathrm{O}(\N^{2})$, and squaring gives $\N^{2}$.) The second line follows by $\partial_{\s}\mathscr{U}'(\s,\cdot)\lesssim1$, by $[\grad^{-}]^{-1}\Delta=-\grad^{+}$ (because $\Delta$ is a second-order gradient), and that $[\grad^{-}]^{-1}\grad^{\mathrm{a}}$ is the identity plus a shift operator on $\mathbb{T}(\N)$; see Definition \ref{definition:intro4}. (We also use the assumption $\mathfrak{t}\lesssim1$.) Now, by Cauchy-Schwarz,
\begin{align}
&\mathbb{P}\left\{\sup_{\t\in\mathbb{X}}\sup_{\x\in\mathbb{T}(\N)}|\mathbf{U}^{\t,\x}|\gtrsim\N^{{\mathrm{C}}}\right\} \nonumber\\
&\lesssim \ [\E\exp\{2\mathscr{N}(\mathfrak{t}(\mathrm{in}),\mathfrak{t})+2\mathscr{Q}(\mathfrak{t}(\mathrm{in}),\mathfrak{t})\}]^{\frac12}\mathbb{P}\left\{\sup_{\t\in\mathbb{X}}\sup_{\x\in\mathbb{T}(\N)}|\mathbf{A}^{\t,\x}|\gtrsim\N^{{\mathrm{C}}}\right\}^{\frac12}. \label{eq:le12I8}
\end{align}
Standard exponential martingale moment upper bounds show that the first factor on the RHS of \eqref{eq:le12I8} is $\lesssim\{\E\exp[\mathrm{O}(\mathscr{Q}(\mathfrak{t}(\mathrm{in}),\mathfrak{t})]\}^{1/2}\lesssim\exp[\mathrm{O}(\N^{{\mathrm{O}(1)}})]$; see \eqref{eq:le12I7}. {We think of the exponent $\mathrm{O}(1)$ as being much smaller than $\mathrm{C}>0$.} We now control the second term on the RHS of \eqref{eq:le12I8}. Let $\mathbf{B}^{\t,\x}$ solve \eqref{eq:le12I6} but with initial data $\mathbf{B}^{\mathfrak{t}(\mathrm{in}),\cdot}$ distributed as $\mathbb{P}^{0,\mathfrak{t}(\mathrm{in})}$. Again, by the mean-value theorem as in \eqref{eq:le12I2b}, we know $\mathbf{C}^{\t,\x}:=\mathbf{A}^{\t,\x}-\mathbf{B}^{\t,\x}$ solves the SDE
\begin{align}
\d\mathbf{C}^{\t,\x} \ = \ \N^{2}\Delta\{\mathscr{U}''(\mathfrak{t}(\mathrm{in}),\mathbf{N}^{\t,\x})\mathbf{C}^{\t,\x}\}\d\t+\N^{\frac32}\grad^{\mathrm{a}}\{\mathscr{U}''(\mathfrak{t}(\mathrm{in}),\mathbf{N}^{\t,\x})\mathbf{C}^{\t,\x}\}\d\t \ =: \ \mathscr{E}\mathbf{C}^{\t,\x}\d\t, \label{eq:le12I9}
\end{align}
for some process $\mathbf{N}^{\t,\x}$. Let $\Gamma^{\mathfrak{t}(\mathrm{in})}[\s,\t,\x,\y]$ solve $\partial_{\t}\Gamma^{\mathfrak{t}(\mathrm{in})}[\s,\t,\x,\y]=\mathscr{E}\Gamma^{\mathfrak{t}(\mathrm{in})}[\s,\t,\x,\y]$ and $\Gamma^{\mathfrak{t}(\mathrm{in})}[\s,\s,\x,\y]=\mathbf{1}[\x=\y]$, where $\mathscr{E}$ in \eqref{eq:le12I9} acts on $\x$. So $\Gamma^{\mathfrak{t}(\mathrm{in})}$ is the fundamental solution for \eqref{eq:le12I9}. Using the reasoning after \eqref{eq:le12I4a}, we know $\Gamma^{\mathfrak{t}(\mathrm{in})}[\s,\t,\x,\y]$ is a probability measure on $\mathbb{T}(\N)$ in the $\x$-variable if $\s\leq\t$, and thus for $\s\leq\t$, we know $\Gamma^{\mathfrak{t}(\mathrm{in})}[\s,\t,\x,\y]\leq1$. Using this with the assumption $|\mathbf{A}^{\mathfrak{t}(\mathrm{in}),\cdot}|=|\mathbf{U}^{\mathfrak{t}(\mathrm{in}),\cdot}|\lesssim\N$, we have the following deterministic bound for all $\t\geq\mathfrak{t}(\mathrm{in})$ and $\x\in\mathbb{T}(\N)$:
\begin{align}
|\mathbf{C}^{\t,\x}| \ \leq \ {\textstyle\sum_{\y}}\Gamma^{\mathfrak{t}(\mathrm{in})}[\mathfrak{t}(\mathrm{in}),\t,\x,\y]\times|\mathbf{C}^{\mathfrak{t}(\mathrm{in}),\y}| \ \lesssim \ |\mathbb{T}(\N)|\{\N+{\textstyle\sup_{\y}}|\mathbf{B}^{\mathfrak{t}(\mathrm{in}),\y}|\}. \label{eq:le12I10}
\end{align}
Now, by the triangle inequality, if $|\mathbf{A}^{\t,\x}|\gtrsim\N^{{\mathrm{C}}}$ for some $\t\in\mathbb{X}$ and $\x\in\mathbb{T}(\N)$, then $\N^{{\mathrm{C}}}\lesssim|\mathbf{A}^{\t,\x}|\leq|\mathbf{C}^{\t,\x}|+|\mathbf{B}^{\t,\x}|\lesssim\N^{2}+\N\|\mathbf{B}\|$, where $\|\mathbf{B}\|$ is the supremum of $|\mathbf{B}^{\t,\x}|$ over $\t\in\mathbb{X}$ and $\x\in\mathbb{T}(\N)$. (This last bound uses \eqref{eq:le12I10}.) Thus,
\begin{align}
\mathbb{P}\left\{\sup_{\t\in\mathbb{X}}\sup_{\x\in\mathbb{T}(\N)}|\mathbf{A}^{\t,\x}|\gtrsim\N^{{\mathrm{C}}}\right\} \ \lesssim \ \mathbb{P}\left\{\sup_{\t\in\mathbb{X}}\sup_{\x\in\mathbb{T}(\N)}|\mathbf{B}^{\t,\x}|\gtrsim\N^{{\mathrm{C}}}-\N^{2}\right\}. \label{eq:le12I11}
\end{align}
By a union bound, it suffices to pull the double supremum outside the probability on the RHS of \eqref{eq:le12I11}, if we insert a factor of $|\mathbb{X}||\mathbb{T}(\N)|\lesssim\N^{\mathrm{D}}$ for some $\mathrm{D}\geq0$ (see immediately before \eqref{eq:le12I5}). Now, because $\mathbf{B}^{\mathfrak{t}(\mathrm{in}),\cdot}\sim\mathbb{P}^{0,\mathfrak{t}(\mathrm{in})}$ and $\mathbf{B}^{\t,\cdot}$ solves \eqref{eq:le12I6}, we know $\mathbf{B}^{\t,\cdot}\sim\mathbb{P}^{0,\mathfrak{t}(\mathrm{in})}$ for all $\t\geq\mathfrak{t}(\mathrm{in})$. (Indeed, $\mathbb{P}^{0,\mathfrak{t}(\mathrm{in})}$ is invariant for \eqref{eq:le12I6}; see Section 2 of \cite{DGP}.) Because $\mathbb{P}^{0,\mathfrak{t}(\mathrm{in})}$ has sub-Gaussian tails and $\mathrm{O}(1)$ variance (by Assumption \ref{ass:intro8}), we deduce $\mathbb{P}\{|\mathbf{B}^{\t,\x}|\gtrsim\N^{{\mathrm{C}}}-\N^{2}\}\lesssim\exp[-{c\N^{2{\mathrm{C}}}}]$. Thus, the RHS of \eqref{eq:le12I11} is $\lesssim\exp[-{c\N^{2{\mathrm{C}}}}]$ {for a possibly different $c>0$}, as exponentials beat polynomials. {Now, if we choose $\mathrm{C}>0$ large enough, then} using this with \eqref{eq:le12I5}, \eqref{eq:le12I8}, \eqref{eq:le12I11} provides $\mathbb{P}[\tau\eqref{eq:glsde}\neq\mathfrak{t}(\mathrm{in})+\mathfrak{t}]\lesssim\exp[-\N]$. This is the desired estimate \eqref{eq:le12I1} (for $\tau\eqref{eq:glsde}$). As noted prior to \eqref{eq:le12I5}, the proof of the $\tau\eqref{eq:glsdeloc}$ estimate is the same. This finishes the proof.
\end{proof}
\subsubsection{Speed of propagation for the current processes, i.e. \eqref{eq:hf} and $\mathbf{J}(\t,\cdot;\mathbb{I}(\mathfrak{t}))$ in Definition \ref{definition:le10}}
In a nutshell, we do the following in this part.
\begin{itemize}
\item First, we show that the SDE \eqref{eq:hf} can be approximated by a localization of it to the discrete interval $\mathbb{I}(\mathfrak{t})$ (in the same sense as Lemma \ref{lemma:le12}). This is the content of Lemma \ref{lemma:le15}.
\item Next, we show that for the aforementioned localization of \eqref{eq:hf}, we can project its values to the torus $\mathbb{S}(\N)\simeq[-\N^{20\gamma_{\mathrm{reg}}},\N^{20\gamma_{\mathrm{reg}}}]$ for free with sufficiently high probability. The resulting process then has generator given in Definition \ref{definition:le5}, letting us then use inequalities like that in Lemma \ref{lemma:le7} to study this process later in this paper. This is the content of Lemma \ref{lemma:le16}.
\end{itemize}
We first introduce notation for an auxiliary process, which is basically $\mathbf{J}(\t,\cdot;\mathbb{I}(\mathfrak{t}))$ but without projecting to $\mathbb{S}(\N)$ as we did in Definition \ref{definition:le10}.
\begin{definition}\label{definition:le14}
 First, take the setting of Definition \ref{definition:le10}.  Now, define $\mathtt{J}(\t,\x;\mathbb{I}(\mathfrak{t}))$ to satisfy the gradient equation $-\N^{1/2}\grad^{\mathbb{I}(\mathfrak{t}),-}\mathtt{J}(\t,\x;\mathbb{I}(\mathfrak{t}))=\mathbf{U}^{\t,\x}[\mathbb{I}(\mathfrak{t})]$ and $\mathtt{J}(\t,\inf\mathbb{I}(\mathfrak{t});\mathbb{I}(\mathfrak{t}))=\mathbf{J}(\t;\mathbb{I}(\mathfrak{t}))$ with $\mathbf{J}(\t;\mathbb{I}(\mathfrak{t}))$ defined by \eqref{eq:hfloc}. (We specify initial data on a case-by-case basis.)
\end{definition}
The reason for Definition \ref{definition:le14} is (somewhat) explained via the following two results. The first compares $\mathtt{J}(\t,\x;\mathbb{I}(\mathfrak{t}))$ to \eqref{eq:hf} in the same setting of Lemma \ref{lemma:le12}. (Roughly, $\mathtt{J}(\t,\x;\mathbb{I}(\mathfrak{t}))$ satisfies \eqref{eq:hf} but $\mathbb{T}(\N)\mapsto\mathbb{I}(\mathfrak{t})$ and $\mathbf{U}^{\t,\cdot}\mapsto\mathbf{U}^{\t,\cdot}[\mathbb{I}(\mathfrak{t})]$. Because \eqref{eq:hf} is a local evolution equation, the aforementioned comparison follows from Lemma \ref{lemma:le12}.) The second compares $\mathtt{J}(\t,\x;\mathbb{I}(\mathfrak{t}))$ to $\mathbf{J}(\t,\x;\mathbb{I}(\mathfrak{t}))$ by proving the projection to $\mathbb{S}(\N)$ does nothing with extremely high probability, which should be true because of our a priori estimates from the stopping times in Definition \ref{definition:method8}. To be rigorous, however, for completely technical (and perhaps uninteresting) reasons, we need other assumptions that go beyond the setting of Lemma \ref{lemma:le12}, for example that the process is locally (close to) a canonical ensemble in Definition \ref{definition:intro5}. The fact that these two comparisons hold in different settings is why we need to introduce Definition \ref{definition:le14}.
\begin{lemma}\label{lemma:le15}
 Assume the setting of {Lemma \ref{lemma:le12}}, and set $\mathbb{I}(\sim):=\llbracket\inf\mathbb{I}-9,\sup\mathbb{I}+9\rrbracket$. Now, suppose $\mathtt{J}(\mathfrak{t}(\mathrm{in}),\inf\mathbb{I}(\mathfrak{t});\mathbb{I}(\mathfrak{t}))=\mathbf{J}(\mathfrak{t}(\mathrm{in}),\inf\mathbb{I}(\mathfrak{t}))$, where $\mathbf{J}(\t,\x)$ solves \eqref{eq:hf}. Next, recall $\mathscr{R}(\s)$ from {Definition \ref{definition:intro6}}. With very high probability, we have
\begin{align}
\sup_{\mathfrak{t}(\mathrm{in})\leq\t\leq\mathfrak{t}(\mathrm{in})+\mathfrak{t}}\sup_{\x\in\mathbb{I}(\sim)}|\{\mathbf{J}(\t,\x)-{\textstyle\int_{\mathfrak{t}(\mathrm{in})}^{\t}}\mathscr{R}(\s)\d\s\}-\mathtt{J}(\t,\x;\mathbb{I}(\mathfrak{t}))| \ \lesssim \ \N^{-500}. \label{eq:le15I}
\end{align}
\end{lemma}
\begin{proof}
Just for convenience, set $\mathbf{J}^{!}(\t,\x):=\mathbf{J}(\t,\x)-${\small$\int_{\mathfrak{t}(\mathrm{in})}^{\t}\mathscr{R}(\s)\d\s$}. We claim the following for all $\t\geq\mathfrak{t}(\mathrm{in})$ and $\x\in\mathbb{I}(\mathfrak{t})$:
{
\begin{align}
\d\mathtt{J}(\t,\x;\mathbb{I}(\mathfrak{t})) \ &:= \ \N^{\frac32}\grad^{\mathbb{I}(\mathfrak{t}),+}\mathscr{U}'(\t,\mathbf{U}^{\t,\x}[\mathbb{I}(\mathfrak{t})])\d\t+\N\{\mathscr{U}'(\t,\mathbf{U}^{\t,\x}[\mathbb{I}(\mathfrak{t})]) \label{eq:le15I1}\\
&+ \ \mathscr{U}'(\t,\mathbf{U}^{\t,\x+1}[\mathbb{I}(\mathfrak{t})])\}\d\t-\mathscr{R}(\t)\d\t+\sqrt{2}\N^{\frac12}\d\mathbf{b}(\t,\x). \nonumber
\end{align}
}(Thus, $\mathtt{J}(\t,\cdot;\mathbb{I}(\mathfrak{t}))$ solves a version of \eqref{eq:hf} but localized $\mathbb{T}(\N)\mapsto\mathbb{I}(\mathfrak{t})$. To verify \eqref{eq:le15I1}, first apply $-\N^{1/2}\grad^{\mathbb{I}(\mathfrak{t}),-}$ to both sides. The RHS of \eqref{eq:le15I1} turns into the RHS of \eqref{eq:glsdeloc}. By the uniqueness of solutions to the SDE \eqref{eq:glsdeloc}, this means \eqref{eq:le15I1} satisfies the gradient relation for $\mathtt{J}$ in Definition \ref{definition:le14}. Moreover, \eqref{eq:le15I1} holds for $\x=\inf\mathbb{I}(\mathfrak{t})$ by Definitions \ref{definition:le10}, \ref{definition:le14}. Therefore, the solution to \eqref{eq:le15I1} evaluated at $\x=\inf\mathbb{I}(\mathfrak{t})$ equals $\mathtt{J}(\t,\inf\mathbb{I}(\mathfrak{t});\mathbb{I}(\mathfrak{t}))$ for all times $\t$. In particular, the solution to \eqref{eq:le15I1} has the same value at $\x=\inf\mathbb{I}(\mathfrak{t})$ and discrete gradient as $\mathtt{J}(\t,\x;\mathbb{I}(\mathfrak{t}))$. Because any function on the torus $\mathbb{I}(\mathfrak{t})$ is determined uniquely by these two pieces of data, \eqref{eq:le15I1} follows.) Set $\mathds{J}(\t,\x):=\mathbf{J}^{!}(\t,\x)-\mathtt{J}(\t,\x;\mathbb{I}(\mathfrak{t}))$. We claim that, for $\t\geq\mathfrak{t}(\mathrm{in})$,
\begin{align}
\d\mathds{J}(\t,\x) \ &= \ \N^{\frac32}\grad^{\mathbb{I}(\mathfrak{t}),+}\{\mathscr{U}'(\t,\mathbf{U}^{\t,\x})-\mathscr{U}'(\t,\mathbf{U}^{\t,\x}[\mathbb{I}(\mathfrak{t})])\}\d\t \label{eq:le15I2a}\\
&+ \ \N\{\mathscr{U}'(\t,\mathbf{U}^{\t,\x})-\mathscr{U}'(\t,\mathbf{U}^{\t,\x}[\mathbb{I}(\mathfrak{t})])\}\d\t \label{eq:le15I2b}\\
&+ \ \N\{\mathscr{U}'(\t,\mathbf{U}^{\t,\x+1})-\mathscr{U}'(\t,\mathbf{U}^{\t,\x+1}[\mathbb{I}(\mathfrak{t})])\}\d\t. \label{eq:le15I2c}
\end{align}
Indeed, Brownian motions in \eqref{eq:hf}, \eqref{eq:le15I1} cancel if $\x\in\mathbb{I}(\sim)$, because, by construction, we have coupled the Brownian motions for all $\x\in\mathbb{I}(\mathfrak{t})\supseteq\mathbb{I}(\sim)$. (This containment follows as $\mathbb{I}(\mathfrak{t})$ is a neighborhood of $\mathbb{I}$ of radius $\mathfrak{l}(\mathfrak{t},\mathbb{I})\gg1$, and $\mathbb{I}(\sim)$ is a neighborhood of radius $9\lesssim1$.) The $-\mathscr{R}(\t)$-drift in \eqref{eq:le15I1} is cancelled because, in $\mathbf{J}^{!}$, we subtracted from $\mathbf{J}(\t,\x)$ the renormalization of speed $\mathscr{R}(\t)$. Finally, the gradients $\grad^{+}$ and $\grad^{\mathbb{I}(\mathfrak{t}),+}$ are the same when acting at $\x\in\mathbb{I}(\sim)$, since $\mathbb{I}(\sim)$ is distance $\gg1$ from the boundary of $\mathbb{I}(\mathfrak{t})$, which is where the gradients start to disagree. See after \eqref{eq:le12I2b} for a similar discussion. Now, we claim $\mathds{J}(\mathfrak{t}(\mathrm{in}),\x)\equiv0$ with probability 1 given any $\x\in\mathbb{I}(\mathfrak{t})$. Indeed, by assumption, we know $-\N^{1/2}\grad^{\mathbb{I}(\mathfrak{t}),-}\mathds{J}(\mathfrak{t}(\mathrm{in}),\x)=\mathbf{U}^{\mathfrak{t}(\mathrm{in}),\x}-\mathbf{U}^{\mathfrak{t}(\mathrm{in}),\x}[\mathbb{I}(\mathfrak{t})]=0$. (This last identity is assumed in the statement of the lemma.) Thus, $\mathds{J}(\mathfrak{t}(\mathrm{in}),\x)$ is constant in $\x\in\mathbb{I}(\mathfrak{t})$. So, to justify this claim, it suffices to show $\mathds{J}(\mathfrak{t}(\mathrm{in}),\inf\mathbb{I}(\mathfrak{t}))=0$; this holds by assumption. We now claim integrating \eqref{eq:le15I2a}-\eqref{eq:le15I2c} gives
\begin{align}
\mathrm{LHS}\eqref{eq:le15I} \ \lesssim \ {\textstyle\int_{\mathfrak{t}(\mathrm{in})}^{\mathfrak{t}(\mathrm{in})+\mathfrak{t}}}\sup_{\x\in\mathbb{I}(\sim)}\{|\eqref{eq:le15I2a}+\eqref{eq:le15I2b}+\eqref{eq:le15I2c}|\}\d\t \ \lesssim \ \N^{300}\mathrm{LHS}\eqref{eq:le12I}. \label{eq:le15I3}
\end{align}
(To establish the last estimate in \eqref{eq:le15I3}, we use $|\mathscr{U}''|\lesssim1$ and $\mathfrak{t}\leq\N^{300}$.) Observe \eqref{eq:le15I3} is deterministic. Using the very high probability estimate \eqref{eq:le12I} gives \eqref{eq:le15I} with very high probability. This completes the proof.
\end{proof}
\begin{lemma}\label{lemma:le16}
 Take the setting of {Lemma \ref{lemma:le15}}. We now introduce more assumptions. Suppose $|\mathbb{I}(\mathfrak{t})|\lesssim\N^{5/6}$ and $\mathfrak{t}|\mathbb{I}(\mathfrak{t})|\lesssim\N^{\gamma}$ for $\gamma\leq{c}\gamma_{\mathrm{KL}}$, {where $c>0$ is a small but fixed constant}. Next, we assume $\mathbf{U}^{\mathfrak{t}(\mathrm{in}),\cdot}[\mathbb{I}(\mathfrak{t})]$ has law given by $\mathbb{P}^{\sigma,\mathfrak{t}(\mathrm{in}),\mathbb{I}(\mathfrak{t})}$. We also assume $|\sigma|\lesssim\N^{\gamma_{\mathrm{reg}}}|\mathbb{I}(\mathfrak{t})|^{-1/2}$. Suppose $\mathbf{J}(\mathfrak{t}(\mathrm{in}),\inf\mathbb{I}(\mathfrak{t});\mathbb{I}(\mathfrak{t}))=\mathtt{J}(\mathfrak{t}(\mathrm{in}),\inf\mathbb{I}(\mathfrak{t});\mathbb{I}(\mathfrak{t}))$, where $\mathbf{J}(\t,\cdot;\mathbb{I}(\mathfrak{t}))$ is from {Definition \ref{definition:le10}}. Now, define the stopping time
\begin{align}
\tau[\mathtt{J}] \ := \ \inf\{\t\in\mathfrak{t}(\mathrm{in})+[0,\mathfrak{t}]: \ {\sup_{\x\in\mathbb{I}}}|\mathtt{J}(\t,\x;\mathbb{I}(\mathfrak{t}))|\geq\N^{15\gamma_{\mathrm{reg}}}\}\wedge[\mathfrak{t}(\mathrm{in})+\mathfrak{t}].
\end{align}
With very high probability, we have the following:
\begin{align}
\sup_{\mathfrak{t}(\mathrm{in})\leq\t\leq\tau[\mathtt{J}]}\sup_{\x\in\mathbb{I}(\mathfrak{t})}|\mathbf{J}(\t,\x;\mathbb{I}(\mathfrak{t}))-\mathtt{J}(\t,\x;\mathbb{I}(\mathfrak{t}))| \ = \ 0. \label{eq:le16I}
\end{align}
\end{lemma}
\begin{proof}
We note $-\N^{1/2}\grad^{\mathbb{I}(\mathfrak{t}),-}\phi(\t,\x)=\mathbf{U}^{\t,\x}[\mathbb{I}(\mathfrak{t})]$ if $\phi(\t,\x)=\mathbf{J}(\t,\x;\mathbb{I}(\mathfrak{t})),\mathtt{J}(\t,\x;\mathbb{I}(\mathfrak{t}))$; see Definitions \ref{definition:le10} and \ref{definition:le14}. Thus, $\mathbf{J}(\t,\x;\mathbb{I}(\mathfrak{t}))-\mathtt{J}(\t,\x;\mathbb{I}(\mathfrak{t}))$ is constant in $\x$ for all $\t\geq\mathfrak{t}(\mathrm{in})$ and $\x\in\mathbb{I}(\mathfrak{t})$. (Indeed, its discrete gradient on $\mathbb{I}(\mathfrak{t})$ is zero.) So,
\begin{align}
\sup_{\substack{\mathfrak{t}(\mathrm{in})\leq\t\leq\tau[\mathtt{J}]\\ \x\in\mathbb{I}(\mathfrak{t})}}|\mathbf{J}(\t,\x;\mathbb{I}(\mathfrak{t}))-\mathtt{J}(\t,\x;\mathbb{I}(\mathfrak{t}))| \ = \ \sup_{\mathfrak{t}(\mathrm{in})\leq\t\leq\tau[\mathtt{J}]}|\mathbf{J}(\t,\inf\mathbb{I}(\mathfrak{t});\mathbb{I}(\mathfrak{t}))-\mathtt{J}(\t,\inf\mathbb{I}(\mathfrak{t});\mathbb{I}(\mathfrak{t}))|. \label{eq:le16I1}
\end{align}
The two functions of $\t$ in $\mathrm{RHS}\eqref{eq:le16I1}$ are the solution to \eqref{eq:hfloc} and the solution to \eqref{eq:hfloc} then projected to $\mathbb{S}(\N)$. We assumed that their data at $\t=\mathfrak{t}(\mathrm{in})$ are the same; see the statement of this lemma. So the term in the supremum in $\mathrm{RHS}\eqref{eq:le16I1}$ is zero for all $\t\leq\tau[\mathbb{I}(\mathfrak{t})]$, where $\tau[\mathbb{I}(\mathfrak{t})]$ is the first time that $\mathtt{J}(\t,\inf\mathbb{I}(\mathfrak{t});\mathbb{I}(\mathfrak{t}))$ leaves our choice of coordinates $\mathbb{S}(\N)\simeq[-\N^{20\gamma_{\mathrm{reg}}},\N^{20\gamma_{\mathrm{reg}}}]$. (This follows by uniqueness of the Ito SDE \eqref{eq:hfloc} and that $\tau[\mathtt{J}],\tau[\mathbb{I}(\mathfrak{t})]$ are stopping times.) We deduce
\begin{align}
\mathbb{P}[\mathrm{RHS}\eqref{eq:le16I1}\neq0] \ \leq \ \mathbb{P}[\tau[\mathbb{I}(\mathfrak{t})]\leq\tau[\mathtt{J}]]. \label{eq:le16I2}
\end{align}
Write $\mathtt{J}(\t,\inf\mathbb{I}(\mathfrak{t});\mathbb{I}(\mathfrak{t}))=\mathtt{J}(\t,\inf\mathbb{I}(\mathfrak{t});\mathbb{I}(\mathfrak{t}))-\mathtt{J}(\t,\w;\mathbb{I}(\mathfrak{t}))+\mathtt{J}(\t,\w;\mathbb{I}(\mathfrak{t}))$, where $\w\in\mathbb{I}$ is fixed. Given any $\t\leq\tau[\mathtt{J}]$, we know $|\mathtt{J}(\t,\w;\mathbb{I}(\mathfrak{t}))|\leq\N^{15\gamma_{\mathrm{reg}}}$ by construction of $\tau[\mathtt{J}]$. Thus, by the triangle inequality, if $|\mathtt{J}(\t,\inf\mathbb{I}(\mathfrak{t});\mathbb{I}(\mathfrak{t}))|\geq\N^{20\gamma_{\mathrm{reg}}}$, we know that $|\mathtt{J}(\t,\inf\mathbb{I}(\mathfrak{t});\mathbb{I}(\mathfrak{t}))-\mathtt{J}(\t,\w;\mathbb{I}(\mathfrak{t}))|\gtrsim\N^{15\gamma_{\mathrm{reg}}}$. We obtain the following, where $\mathcal{E}$ is the event $|\mathtt{J}(\t,\inf\mathbb{I}(\mathfrak{t});\mathbb{I}(\mathfrak{t}))-\mathtt{J}(\t,\w;\mathbb{I}(\mathfrak{t}))|\gtrsim1$ for some time $\t\in\mathfrak{t}(\mathrm{in})+[0,\mathfrak{t}]$:
\begin{align}
\mathbb{P}[\tau[\mathbb{I}(\mathfrak{t})]\leq\tau[\mathtt{J}]] \ \leq \ \mathbb{P}[\mathcal{E}]. \label{eq:le16I3}
\end{align}
By construction, $\mathtt{J}(\t,\w;\mathbb{I}(\mathfrak{t}))-\mathtt{J}(\t,\inf\mathbb{I}(\mathfrak{t});\mathbb{I}(\mathfrak{t}))$ is $\N^{-1/2}$ times a random walk with steps $\mathbf{U}^{\t,\y}[\mathbb{I}(\mathfrak{t})]$:
\begin{align}
\mathtt{J}(\t,\w;\mathbb{I}(\mathfrak{t}))-\mathtt{J}(\t,\inf\mathbb{I}(\mathfrak{t});\mathbb{I}(\mathfrak{t})) \ = \ {\fontsize{9}{12}\N^{-\frac12}\sum_{\y=\inf\mathbb{I}(\mathfrak{t})+1}^{\w}\mathbf{U}^{\t,\y}[\mathbb{I}(\mathfrak{t})]}. \label{eq:le16I4}
\end{align}
{Take $\t\in\mathfrak{t}(\mathrm{in})+[0,\mathfrak{t}]$. For now, \emph{assume} that $\mathbf{U}^{\t,\cdot}[\mathbb{I}(\mathfrak{t})]\sim\mathbb{P}^{\sigma,\t,\mathbb{I}(\mathfrak{t})}$ with $\sigma$ from the statement of the lemma. (This is not necessarily true; we remedy this shortly.) In this case, $\mathrm{RHS}\eqref{eq:le16I4}$ is $\N^{-1/2}$ times a random walk bridge increment with length at most $|\mathbb{I}(\mathfrak{t})|$. The average drift of the random walk bridge is $\sigma\lesssim\N^{\gamma_{\mathrm{reg}}}|\mathbb{I}(\mathfrak{t})|^{-1/2}$ by assumption. Following the end of the proof of Lemma \ref{lemma:bg27}, we deduce $\mathrm{RHS}\eqref{eq:le16I4}$ is $\lesssim$ the sum of the following three quantities. The first is $\N^{-1/2}|\mathbb{I}(\mathfrak{t})||\sigma|\lesssim\N^{-1/2+\gamma_{\mathrm{reg}}}|\mathbb{I}(\mathfrak{t})|^{1/2}$. This is the contribution from the average drift. The second quantity is a sum of $\lesssim|\mathbb{I}(\mathfrak{t})|$-many terms $\N^{-1/2}\mathbf{V}$, where $\mathbf{V}+\sigma\sim\mathbb{P}^{\sigma,\t}$ are independent. The third has the same form as the second. (Again, the second and third quantities that we just explained come from a random walk analog of the representation of Brownian bridge as $\mathrm{B}(\t)-\t\mathrm{B}(1)$.) Thus, by the triangle inequality, we deduce that on $\mathcal{E}$, the sum of $\lesssim|\mathbb{I}(\mathfrak{t})|$-many sub-Gaussian martingale increments with variance parameter $\lesssim\N^{-1}$ must, in absolute value, exceed $1-\N^{-1/2+\gamma_{\mathrm{reg}}}|\mathbb{I}(\mathfrak{t})|^{1/2}\gtrsim1$; this last bound follows by $|\mathbb{I}(\mathfrak{t})|\lesssim\N^{5/6}$. If $\mathrm{Mart}$ is said martingale, then
\begin{align}
\mathbb{P}[\mathcal{E}] \ \leq \ \mathbb{P}[|\mathrm{Mart}|\gtrsim1] \ \leq \ \mathbb{P}[|\mathrm{Mart}|\gtrsim\N^{-\frac12+\frac{1}{99}}|\mathbb{I}(\mathfrak{t})|^{\frac12}] \ \lesssim \ \exp[-\N^{\frac{1}{100}}]. \label{eq:le16I5}
\end{align}
The second inequality follows by the assumption $|\mathbb{I}(\mathfrak{t})|\lesssim\N^{5/6}$. The last bound follows by Azuma. (The exponent $1/99$ pushes beyond the natural martingale scale by $\gtrsim\N^{1/100}$, which gives the exponential decay in the far RHS of \eqref{eq:le16I5}.) Recall that this is all under the assumption that $\mathbf{U}^{\t,\cdot}[\mathbb{I}(\mathfrak{t})]\sim\mathbb{P}^{\sigma,\t,\mathbb{I}(\mathfrak{t})}$ with $\sigma$ from the statement of the lemma. Again, this is not necessarily true. But assumed $\mathfrak{t}|\mathbb{I}(\mathfrak{t})|\lesssim\N^{\gamma}$. So Lemma \ref{lemma:kv8} implies that \eqref{eq:le16I5} still holds if we multiply by $\exp[\mathrm{O}(\N^{\gamma})]$ everywhere after the first $\leq$-sign in \eqref{eq:le16I5}. (We clarify Lemma \ref{lemma:kv8} does not need the current lemma, so there is no circular reasoning.) This exponential factor is overwhelmed by the exponential decay in the far RHS of \eqref{eq:le16I5}. Thus, we deduce $\mathbb{P}[\mathcal{E}]$ is still exponentially small in $\N$. Combining this with \eqref{eq:le16I1}, \eqref{eq:le16I2}, \eqref{eq:le16I3}, \eqref{eq:le16I4} finishes the proof.}
\end{proof}
%
%
%
\section{Non-equilibrium Kipnis-Varadhan inequality}\label{section:kv}
The goal of this section is to make rigorous a version of the following. Space-time averages of local fluctuations have square-root cancellation. See Section \ref{section:sqle} for discussion of the non-equilibrium-type challenges (coming from time-inhomogeneity of \eqref{eq:hf}-\eqref{eq:glsde}) that we must resolve to this end. Ultimately, we prove the following estimate; see after its statement for an intuitive explanation of what it is actually saying. (The reader is invited to go directly there before reading the statement of Proposition \ref{prop:kv1} to get a clearer handle on the statement.) First, recall notation of Definition \ref{definition:le10}.
\begin{prop}\label{prop:kv1}
 Take $\mathfrak{t}(\mathrm{in}),\mathfrak{t}\geq0$ and any non-empty discrete interval $\mathbb{I}\subseteq\mathbb{T}(\N)$. Now, we define $\gamma_{\mathrm{KV}}={c}\gamma_{\mathrm{KL}}$ with $\gamma_{\mathrm{KL}}$ from {Definition \ref{definition:entropydata}}, {where $c>0$ is a small but fixed constant}. Suppose that $\mathfrak{t}|\mathbb{I}(\mathfrak{t})|\lesssim\N^{\gamma_{\mathrm{av}}}$ and $|\mathbb{I}(\mathfrak{t})|\gtrsim\N^{1/10}$, in which $\gamma_{\mathrm{av}}\leq{c}\gamma_{\mathrm{KV}}$. 

Now consider the joint process $\t\mapsto(\mathbf{J}(\t,\inf\mathbb{I}(\mathfrak{t});\mathbb{I}(\mathfrak{t})),\mathbf{U}^{\t,\cdot}[\mathbb{I}(\mathfrak{t})])$. Assume that its time-$\mathfrak{t}(\mathrm{in})$ data $(\mathbf{J}(\mathfrak{t}(\mathrm{in}),\inf\mathbb{I}(\mathfrak{t});\mathbb{I}(\mathfrak{t})),\mathbf{U}^{\mathfrak{t}(\mathrm{in}),\cdot}[\mathbb{I}(\mathfrak{t})])$ is distributed according to the product measure $\mathbb{P}^{\delta[0],\sigma,\mathfrak{t}(\mathrm{in}),\mathbb{I}(\mathfrak{t})}=\delta[0]\otimes\mathbb{P}^{\sigma,\mathfrak{t}(\mathrm{in}),\mathbb{I}(\mathfrak{t})}$ on $\mathbb{S}(\N)\times\R^{\mathbb{I}(\mathfrak{t})}$, where $\delta[0]$ is the Dirac mass at $0\in\mathbb{S}(\N)$ and $|\sigma|\lesssim1$. 

Fix $\varphi\in\mathscr{L}^{2}(\mathbb{S}(\N))\cap\mathscr{L}^{\infty}(\mathbb{S}(\N))${, and fix any $\mathrm{m}\in\mathbb{N}$}. Consider $\mathfrak{a}(\t,\cdot;1),\ldots,\mathfrak{a}(\t,\cdot:\mathrm{m}):\R^{\mathbb{I}(\mathfrak{t})}\to\R$ such that for $\t\leq\mathfrak{t}(\mathrm{in})+\mathfrak{t}$, the functions $\mathfrak{a}(\t,\cdot;\mathrm{k})$ are smooth outside of a deterministic and finite set of jump times. Assume $\mathfrak{a}(\t,\cdot;\mathrm{k})$ satisfy assumptions in {Lemma \ref{lemma:le2}} with sets $\mathbb{J}(\t,\mathrm{k})$. Assume that $\mathbb{J}(\t,\mathrm{k})$ are $\t$-dependent shifts of $\mathbb{J}(\mathrm{k})$. We let $\mathfrak{a}(\t,\cdot)$ denote the average of $\mathfrak{a}(\t,\cdot;1),\ldots,\mathfrak{a}(\t,\cdot{;}\mathrm{m})$. Let $\mathrm{U}$ be uniform $[-1,1]$. Assume that it is independent of everything else. Now, set the following for $\t\in\mathfrak{t}(\mathrm{in})+[0,\mathfrak{t}]$, in which $\mathrm{U}+\mathbf{J}(\s,\inf\mathbb{I}(\mathfrak{t});\mathbb{I}(\mathfrak{t}))$ is with respect to periodic boundary on $\mathbb{S}(\N)$:
\begin{align}
\mathscr{A}(\t) \ := \ [\t-\mathfrak{t}(\mathrm{in})]^{-1}{\textstyle\int_{\mathfrak{t}(\mathrm{in})}^{\t}}\{\varphi(\mathrm{U}+\mathbf{J}(\s,\inf\mathbb{I}(\mathfrak{t});\mathbb{I}(\mathfrak{t})))\mathfrak{a}(\s,\mathbf{U}^{\s,\cdot}[\mathbb{I}(\mathfrak{t})])\}\d\s. \label{eq:kv1I}
\end{align}
(Let $\mathscr{A}(\mathfrak{t}(\mathrm{in}))$ be the integrand in $\mathrm{RHS}\eqref{eq:kv1I}$ at $\s=\mathfrak{t}(\mathrm{in})$.) {Then, with notation explained after, for any $\t\in\mathfrak{t}(\mathrm{in})+[0,\mathfrak{t}]$ and $\mathscr{B}\geq0$, we have the following for any large but fixed $\mathrm{D}>0$:}
\begin{align}
&\N^{-20\gamma_{\mathrm{reg}}}\E\{|\mathscr{A}(\t)|^{2}\mathbf{1}[|\mathscr{A}(\t)|\lesssim\mathscr{B}]\} \label{eq:kv1II}\\
&\lesssim \ \N^{\gamma_{\mathrm{KV}}}\mathfrak{t}|\mathbb{I}(\mathfrak{t})|^{\frac12}\mathscr{B}^{2} + \N^{-{\mathrm{D}}}\mathscr{B}^{2}\nonumber\\
&+ \ \N^{-2}\mathrm{m}^{-1}\|\varphi\|^{2}\times{\sup}_{\mathrm{k}}\{|\mathbb{J}(\mathrm{k})|^{2}\|\mathfrak{a}(\cdot,\cdot;\mathrm{k})\|_{\infty}^{2}\}\times\{[\t-\mathfrak{t}(\mathrm{in})]^{-1}+1\}. \nonumber
\end{align}
{The expectation} in $\mathrm{LHS}\eqref{eq:kv1II}$ is with respect to the law of $\t\mapsto(\mathbf{J}(\t,\inf\mathbb{I}(\mathfrak{t});\mathbb{I}(\mathfrak{t})),\mathbf{U}^{\t,\cdot}[\mathbb{I}(\mathfrak{t})])$ for $\t\in\mathfrak{t}(\mathrm{in})+[0,\mathfrak{t}]$. In $\mathrm{RHS}\eqref{eq:kv1II}$, we used $\|\varphi\|^{2}:=\E^{\mathrm{Leb}}|\varphi|^{2}+\sup_{\mathrm{a}\in\R}|\varphi(\mathrm{a})|^{2}$. In the same setting, we also have the following estimate:
\begin{align}
\N^{-20\gamma_{\mathrm{reg}}}\E\{\mathscr{B}^{2}\mathbf{1}[|\mathscr{A}(\t)|\geq\mathscr{B}]\} \ \lesssim \ \mathrm{RHS}\eqref{eq:kv1II}. \label{eq:kv1III}
\end{align}
\end{prop}
Let us now explain what the bounds \eqref{eq:kv1II} and \eqref{eq:kv1III} say. We want to bound the time-average of $\varphi\times\mathfrak{a}$, where $\varphi$ is a function of the current $\mathbf{J}(\cdot,\inf\mathbb{I}(\mathfrak{t});\mathbb{I}(\mathfrak{t}))$, and $\mathfrak{a}$ is an average of ``orthogonal", fluctuating, and local functionals of $\mathbf{U}^{\cdot,\cdot}[\mathbb{I}(\mathfrak{t})]$. As in Lemma \ref{lemma:le7}, the $\varphi$-factor is harmless and contributes just its norm in \eqref{eq:kv1II}. For the $\mathfrak{a}$-average, like with the usual Kipnis-Varadhan inequality (see Appendix 1.6 of \cite{KL}, for example), we can control its time-average by resolvent estimates as in Lemma \ref{lemma:le7}; this explains the last term in \eqref{eq:kv1II} (the factor $[\t-\mathfrak{t}(\mathrm{in})]^{-1}+1$ is there just to control the time-scale on which we integrate). The first two terms in the bound in \eqref{eq:kv1II} come from controlling the time-inhomogeneity of the dynamics. Indeed, the $\mathscr{B}^{2}$-factor is an a priori estimate for the square of $\mathscr{A}(\t)$ on the event where $|\mathscr{A}(\t)|\lesssim\mathscr{B}$. The factor $\mathfrak{t}$ comes from the time-scale on which we integrate, and the $|\mathbb{I}(\mathfrak{t})|^{1/2}$-factor comes from the fluctuating property explained in point (1) of Section \ref{section:sqle}. (The extra factors of $\N^{\gamma_{\mathrm{KV}}}$ and $\N^{20\gamma_{\mathrm{reg}}}$ are harmless; the latter, for example, comes from changing the initial law of the current process to uniform on its state space $\mathbb{S}(\N)\simeq[-\N^{20\gamma_{\mathrm{reg}}},\N^{20\gamma_{\mathrm{reg}}}]$.) The point of this section is to make this paragraph rigorous, so the reader is invited to skip in a first reading. (The final estimate \eqref{eq:kv1III} is basically the same. Modulo indicators, we could control {it} by the LHS of \eqref{eq:kv1II} by Chebyshev. But the role of the indicator in \eqref{eq:kv1II} is to give an a priori upper bound of $\mathscr{B}^{2}$, which the LHS of \eqref{eq:kv1III} clearly has anyway.)

Let us now be a little more precise about the previous heuristic. The initial data of $\mathbf{U}^{\t,\cdot}[\mathbb{I}(\mathfrak{t})]$ is local equilibrium $\mathbb{P}^{\sigma,\mathfrak{t}(\mathrm{in}),\mathbb{I}(\mathfrak{t})}$. In particular, Proposition \ref{prop:kv1} addresses non-equilibrium aspects arising from time-inhomogeneity of SDEs, not from our choice of initial data for any SDEs. (For example, it addresses the issue that there is no notion of stationarity for \eqref{eq:hf}-\eqref{eq:glsde}.) Note that if $\mathbf{U}^{\t,\cdot}[\mathbb{I}(\mathfrak{t})]$ was stationary, \eqref{eq:kv1II} would be the Kipnis-Varadhan bound; see Section 4 of \cite{CLO}. The problem is to somehow derive it in the non-equilibrium and time-inhomogeneous case. To this end, we now intuitively quantify the discussion in Section \ref{section:sqle}. As we explained there, we compare the SDE \eqref{eq:glsdeloc} to another SDE (for which the proof of Kipnis-Varadhan is accessible). The cost of comparison is $\lesssim$ the first two terms in $\mathrm{RHS}\eqref{eq:kv1II}$. (The second term in $\mathrm{RHS}\eqref{eq:kv1II}$ is completely harmless.) Let us briefly explain why. Clearly, the term inside $\E$ in $\mathrm{LHS}\eqref{eq:kv1II}$ is $\lesssim\mathscr{B}^{2}$. This gives the $\mathscr{B}^{2}$-factor. Now, recall more precisely from Section \ref{section:sqle} that the cost of comparison depends on a time-integrated Dirichlet form on $\mathbb{I}(\mathfrak{t})$. We integrate-in-time over $\mathfrak{t}(\mathrm{in})+[0,\mathfrak{t}]$, which is of length $\mathfrak{t}$. This explains $\mathfrak{t}$ in $\mathrm{RHS}\eqref{eq:kv1II}$. {The} Dirichlet form sums energies per bond in $\mathbb{I}(\mathfrak{t})$. Thus we expect $|\mathbb{I}(\mathfrak{t})|$ in $\mathrm{RHS}\eqref{eq:kv1II}$. However, we have $|\mathbb{I}(\mathfrak{t})|^{1/2}$! This power-saving comes from the fluctuating property of the ``total potential" in the first bullet point in Section \ref{section:sqle}. (The assumption $|\mathbb{I}(\mathfrak{t})|\gtrsim\N^{1/10}$ is just to ensure enough spatial fluctuations of the total potential; the exponent $1/10$ is just something noticeably bigger than $\gamma_{\mathrm{KL}}$.)

Before we proceed, unless otherwise mentioned, we emphasize that $|\sigma|\lesssim1$ throughout this section.
\subsection{Preliminary constructions}
We need to modify $\mathbf{J}(\t,\cdot;\mathbb{I}(\mathfrak{t}))$ and $\mathbf{U}^{\t,\cdot}[\mathbb{I}(\mathfrak{t})]$ for entirely technical reasons. This will help us take advantage of fluctuations in the total potential from Section \ref{section:sqle}. First, some other constructions.
\begin{definition}\label{definition:kv2}
 Fix $\sigma\in\R$ and a discrete interval $\mathbb{K}\subseteq\mathbb{T}(\N)$. Let $\mathbb{H}^{\sigma,\mathbb{K}}$ be the hyperplane in $\R^{\mathbb{K}}$ consisting of all $\mathbf{U}$ such that the average of $\mathbf{U}(\x)$ over $\x\in\mathbb{K}$ is $\sigma$. Let $\d^{\sigma,\mathbb{K}}$ be Lebesgue measure on $\mathbb{H}^{\sigma,\mathbb{K}}$. (It is induced by Euclidean metric on $\mathbb{H}^{\sigma,\mathbb{K}}\subseteq\R^{\mathbb{K}}$, which is the metric on $\mathbb{H}^{\sigma,\mathbb{K}}$ determined by this containment and Euclidean metric on $\R^{\mathbb{K}}$.)
\end{definition}
\begin{lemma}\label{lemma:kv3}
 Let $\mathfrak{p}[\sigma,\t,\mathbb{I}(\mathfrak{t})]$ be the density of $\mathbb{P}^{\sigma,\t,\mathbb{I}(\mathfrak{t})}$ with respect to $\d^{\sigma,\mathbb{I}(\mathfrak{t})}$. {Then, as functions of $\mathbf{U}\in\R^{\mathbb{I}(\mathfrak{t})}$, we have (with notation explained and intuitively clarified after)}
\begin{align}
\mathfrak{p}[\sigma,\t,\mathbb{I}(\mathfrak{t})] \ = \ \exp\{-\mathscr{HP}(\t,\mathbf{U};\sigma)\} \ := \ \exp\{-{\textstyle\sum_{\x}}\mathscr{UP}(\t,\mathbf{U}(\x);\sigma)\}, \label{eq:kv3I}
\end{align}
where the sum is over $\x\in\mathbb{I}(\mathfrak{t})$, and where $\mathscr{UP}(\t,\mathbf{U}(\x);\sigma):=\mathscr{U}(\t,\mathbf{U}(\x))-|\mathbb{I}(\mathfrak{t})|^{-1}\mathscr{P}[\sigma,\t,\mathbb{I}(\mathfrak{t})]$ with
\begin{align}
\mathscr{P}[\sigma,\t,\mathbb{I}(\mathfrak{t})] \ := \ -\log{\textstyle\int_{\mathbb{H}^{\sigma,\mathbb{I}(\mathfrak{t})}}}\exp\{-{\textstyle\sum_{\x}}\mathscr{U}(\t,\mathbf{U}(\x))\}\d^{\sigma,\mathbb{I}(\mathfrak{t})}(\mathbf{U}). \label{eq:kv3II}
\end{align}
(Intuitively, $\mathscr{UP}$ is ``$\mathscr{U}$ plus pressure", and $\mathscr{HP}$ is ``Hamiltonian plus pressure".) We also claim the inequalities $|\partial_{\t}\mathscr{UP}(\t,\cdot;\sigma)|\lesssim1$ and $|\partial_{\t}^{2}\mathscr{UP}(\t,\cdot;\sigma)|\lesssim\N^{100}$ uniformly in $\t\geq0$. Lastly, for any $\t,\x,\sigma$, we have the following fluctuation property from Section \ref{section:sqle}:
\begin{align}
\E^{\sigma,\t,\mathbb{I}(\mathfrak{t})}\partial_{\t}\mathscr{UP}(\t,\mathbf{U}(\x);\sigma) \ = \ 0. \label{eq:kv3III}
\end{align}
\end{lemma}
\begin{proof}
Everything except the second-time-derivative bound was already proved for $\sigma=0$ and $\mathbb{I}(\mathfrak{t})\mapsto\mathbb{T}(\N)$; see \eqref{eq:le8II3a}, \eqref{eq:le8II3b}, \eqref{eq:le8II6c}, and right after \eqref{eq:le8II6c}. The proof works for any $\sigma$ and interval in $\mathbb{T}(\N)$. We are left to obtain the second-time-derivative bound. By Assumption \ref{ass:intro8}, we know that $|\partial_{\t}^{2}\mathscr{U}(\t,\cdot)|\lesssim1$. So, because $\mathscr{UP}(\t,\mathbf{U}(\x);\sigma):=\mathscr{U}(\t,\mathbf{U}(\x))-|\mathbb{I}(\mathfrak{t})|^{-1}\mathscr{P}[\sigma,\t,\mathbb{I}(\mathfrak{t})]$, it suffices to show $|\partial_{\t}^{2}\mathscr{P}[\sigma,\t,\mathbb{I}(\mathfrak{t})]|\lesssim\N^{100}$. We proceed with a direct calculation. We already know $\partial_{\t}\mathscr{P}[\sigma,\t,\mathbb{I}(\mathfrak{t})]$; use \eqref{eq:le8II7a}-\eqref{eq:le8II7b} but replace $(\s,0,\mathbb{T}(\N))\mapsto(\t,\sigma,\mathbb{I}(\mathfrak{t}))$. Ultimately, we claim the following calculation:
\begin{align}
&\partial_{\t}^{2}\mathscr{P}[\sigma,\t,\mathbb{I}(\mathfrak{t})] \nonumber\\
&= \ |\mathbb{T}(\N)|\E^{\sigma,\t,\mathbb{I}(\mathfrak{t})}\partial_{\t}^{2}\mathscr{U}(\t,\mathbf{U}(\x))+|\mathbb{T}(\N)|{\textstyle\int_{\mathbb{H}^{\sigma,\mathbb{I}(\mathfrak{t})}}}\partial_{\t}\mathscr{U}(\t,\mathbf{U}(\x))\partial_{\t}\mathfrak{p}[\sigma,\t,\mathbb{I}(\mathfrak{t})]\d^{\sigma,\mathbb{I}(\mathfrak{t})}(\mathbf{U}). \label{eq:kv31}
\end{align}
By \eqref{eq:le8II7a}-\eqref{eq:le8II7b}, we know $\mathrm{LHS}\eqref{eq:kv31}$ is equal to $\partial_{\t}$ of the far RHS of \eqref{eq:le8II7b}. (Again, this is with $(\s,0,\mathbb{T}(\N))\mapsto(\t,\sigma,\mathbb{I}(\mathfrak{t}))$.) By Leibniz rule, $\partial_{\t}$ hits both $\partial_{\t}\mathscr{U}(\t,\mathbf{U}(\x))$ and the measure for $\E^{\sigma,\t,\mathbb{I}(\mathfrak{t})}$. The former gives the first term in $\mathrm{RHS}\eqref{eq:kv31}$. The latter gives the second term therein, because $\E^{\sigma,\t,\mathbb{I}(\mathfrak{t})}$ is integration against $\mathbb{P}^{\sigma,\t,\mathbb{I}(\mathfrak{t})}$, whose $\d^{\sigma,\mathbb{I}(\mathfrak{t})}$ density is $\mathfrak{p}[\sigma,\t,\mathbb{I}(\mathfrak{t})]$ by construction in the lemma. By \eqref{eq:kv3I}, we know $\partial_{\t}\mathfrak{p}[\sigma,\t,\mathbb{I}(\mathfrak{t})]$ is equal to $\mathfrak{p}[\sigma,\t,\mathbb{I}(\mathfrak{t})]$ times the sum of $\partial_{\t}\mathscr{UP}(\t,\mathbf{U}(\x);\sigma)\lesssim1$ terms (over all $\x\in\mathbb{I}(\mathfrak{t})$). So the last term in \eqref{eq:kv31} is $\lesssim\E^{\sigma,\t,\mathbb{I}(\mathfrak{t})}\mathfrak{p}[\sigma,\t,\mathbb{I}(\mathfrak{t})]|\mathbb{T}(\N)|^{2}\lesssim|\mathbb{T}(\N)|^{2}$. Assumption \ref{ass:intro8} also states that the first term on the RHS of \eqref{eq:kv31} is $\lesssim|\mathbb{T}(\N)|$. Thus, $\partial_{\t}^{2}\mathscr{P}[\sigma,\t,\mathbb{I}(\mathfrak{t})]\lesssim|\mathbb{T}(\N)|^{2}\lesssim\N^{100}$, and we are done.
\end{proof}
We now give the main construction of this subsection. Recall $\mathscr{HP}$ ("Hamiltonian plus pressure") from Lemma \ref{lemma:kv3}. It basically restricts processes to the domain on which a square-root cancellation estimate for the $\mathscr{HP}$-term holds in a way that is friendly for stochastic calculus. See after Remark \ref{remark:kv5} for more explanation.
\begin{definition}\label{definition:kv4}
 Take $\mathcal{E}[\mathfrak{t}(\mathrm{in}),\mathfrak{t};\mathbb{I}(\mathfrak{t})]\subseteq\mathbb{H}^{\sigma,\mathbb{I}(\mathfrak{t})}$ so that if $\mathbf{U}\in\mathcal{E}[\mathfrak{t}(\mathrm{in}),\mathfrak{t};\mathbb{I}(\mathfrak{t})]$ and $\t\in\mathfrak{t}(\mathrm{in})+[0,\mathfrak{t}]$, we have $|\partial_{\t}\mathscr{HP}(\t,\mathbf{U};\sigma)|\lesssim\N^{\gamma_{\mathrm{KV}}}|\mathbb{I}(\mathfrak{t})|^{1/2}$. Also, assume that if $\mathbf{U}\not\in\mathcal{E}[\mathfrak{t}(\mathrm{in}),\mathfrak{t};\mathbb{I}(\mathfrak{t})]$, we know $|\partial_{\t}\mathscr{HP}(\t,\mathbf{U};\sigma)|\gtrsim\N^{\gamma_{\mathrm{KV}}}|\mathbb{I}(\mathfrak{t})|^{1/2}$ for some $\t\in\mathfrak{t}(\mathrm{in})+[0,\mathfrak{t}]$. (See Proposition \ref{prop:kv1} for $\gamma_{\mathrm{KV}}$. Also, assume $\mathcal{E}[\mathfrak{t}(\mathrm{in}),\mathfrak{t};\mathbb{I}(\mathfrak{t})]\subseteq\mathbb{H}^{\sigma,\mathbb{I}(\mathfrak{t})}$ has smooth boundary.) 

For the same set of $\t$ and for $\tau[\mathcal{E}]$ defined to be the first time $\t\in\mathfrak{t}(\mathrm{in})+[0,\mathfrak{t}]$ where $\mathbf{U}^{\t,\cdot}[\mathbb{I}(\mathfrak{t})]\not\in\mathcal{E}[\mathfrak{t}(\mathrm{in}),\mathfrak{t};\mathbb{I}(\mathfrak{t})]$, we let $\mathbf{U}^{\t,\cdot}[\mathbb{I}(\mathfrak{t}),\mathcal{E}]$ solve \eqref{eq:glsdeloc} until $\tau[\mathcal{E}]$. Also, set $\mathbf{J}(\t;\mathbb{I}(\mathfrak{t}),\mathcal{E})=\mathbf{J}(\t\wedge\tau[\mathcal{E}],\inf\mathbb{I}(\mathfrak{t});\mathbb{I}(\mathfrak{t}))$; see Definition \ref{definition:le10}. Now, define $\mathbb{P}^{\mathrm{Leb},\sigma,\t,\mathbb{I}(\mathfrak{t}),\mathcal{E}}:=\mathrm{Leb}(\mathbb{S}(\N))\otimes\mathbb{P}^{\sigma,\t,\mathbb{I}(\mathfrak{t}),\mathcal{E}}$, where $\mathbb{P}^{\sigma,\t,\mathbb{I}(\mathfrak{t}),\mathcal{E}}$ is defined to be the measure $\mathbb{P}^{\sigma,\t,\mathbb{I}(\mathfrak{t})}$ from Definition \ref{definition:intro5} but conditioned on $\mathcal{E}[\mathfrak{t}(\mathrm{in}),\mathfrak{t};\mathbb{I}(\mathfrak{t})]$.
\end{definition}
\begin{rem}\label{remark:kv5}
 A word on building $\mathcal{E}[\mathfrak{t}(\mathrm{in}),\mathfrak{t};\mathbb{I}(\mathfrak{t})]$. By Assumption \ref{ass:intro8}, the function $\mathbf{U}\mapsto\partial_{\t}\mathscr{HP}(\t,\mathbf{U};\sigma)$ is uniformly Lipschitz (also uniformly over $\t$). (Indeed, by Lemma \ref{lemma:kv3}, we know $\mathbf{U}\mapsto\partial_{\t}\mathscr{HP}(\t,\mathbf{U};\sigma)$ is a sum over $\x\in\mathbb{I}(\mathfrak{t})$ of $\mathbf{U}\mapsto\partial_{\t}\mathscr{U}(\t,\mathbf{U}(\x))$ plus something constant in $\mathbf{U}$.) So, $\mathbf{U}\mapsto\sup_{\t}|\partial_{\t}\mathscr{HP}(\t,\mathbf{U};\sigma)|$ is uniformly Lipschitz. Take a smooth function that is uniformly within $1$ of $\mathbf{U}\mapsto\sup_{\t}|\partial_{\t}\mathscr{HP}(\t,\mathbf{U};\sigma)|$. Via Sard's theorem, we know there exists $\N^{\gamma_{\mathrm{KV}}}|\mathbb{I}(\mathfrak{t})|^{1/2}\lesssim\mathrm{a}\lesssim\N^{\gamma_{\mathrm{KV}}}|\mathbb{I}(\mathfrak{t})|^{1/2}$ so the level set of this smooth function (for value $\mathrm{a}$) is smooth. Let $\mathcal{E}[\mathfrak{t}(\mathrm{in}),\mathfrak{t};\mathbb{I}(\mathfrak{t})]$ be the union over all $0\leq\mathrm{a}'\leq\mathrm{a}$ of level sets of the smooth function (for the value $\mathrm{a}'$). One can readily check that this set satisfies {all of the} conditions on Definition \ref{definition:kv4}.
\end{rem}
In Definition \ref{definition:kv4}, we just stop $\t\mapsto(\mathbf{J}(\t,\inf\mathbb{I}(\mathfrak{t});\mathbb{I}(\mathfrak{t})),\mathbf{U}^{\t,\cdot}[\mathbb{I}(\mathfrak{t})])$ when $\mathbf{U}^{\t,\cdot}[\mathbb{I}(\mathfrak{t})]$ exits $\mathcal{E}[\mathfrak{t}(\mathrm{in}),\mathfrak{t};\mathbb{I}(\mathfrak{t})]$. Note $\mathcal{E}[\mathfrak{t}(\mathrm{in}),\mathfrak{t};\mathbb{I}(\mathfrak{t})]$ is when the time-derivative of the total potential $\mathscr{HP}$ has square-root cancellation. Intuitively, because of the extra $\N^{\gamma_{\mathrm{KV}}}$ factor, by \eqref{eq:kv3III}, we expect these cancellations with exponentially high probability. The following result makes this precise. (Again, for the reader interested in reading the proof of \eqref{eq:kv6I}, see before Lemma \ref{lemma:le8} for an intuitive description of how we leverage/make precise ``square-root cancellations".)
\begin{lemma}\label{lemma:kv6}
 Recall $\gamma_{\mathrm{av}},\gamma_{\mathrm{KV}}$ from {Proposition \ref{prop:kv1}}. Suppose $\mathfrak{t}|\mathbb{I}(\mathfrak{t})|\lesssim\N^{\gamma_{\mathrm{av}}}$. For any $\mathfrak{t}(\mathrm{in})\leq\t\leq\mathfrak{t}(\mathrm{in})+\mathfrak{t}$, we have
\begin{align}
\mathbb{P}^{\mathrm{Leb},\sigma,\t,\mathbb{I}(\mathfrak{t})}\{\mathcal{E}[\mathfrak{t}(\mathrm{in}),\mathfrak{t};\mathbb{I}(\mathfrak{t})]^{\mathrm{C}}\} \ = \ \mathbb{P}^{\sigma,\t,\mathbb{I}(\mathfrak{t})}\{\mathcal{E}[\mathfrak{t}(\mathrm{in}),\mathfrak{t};\mathbb{I}(\mathfrak{t})]^{\mathrm{C}}\} \ \lesssim \ \exp\{\N^{-\frac12\gamma_{\mathrm{KV}}}\}. \label{eq:kv6I}
\end{align}
\end{lemma}
\begin{proof}
The identity in \eqref{eq:kv6I} follows because the event $\mathcal{E}[\mathfrak{t}(\mathrm{in}),\mathfrak{t};\mathbb{I}(\mathfrak{t})]^{\mathrm{C}}$ is independent of the $\mathbb{S}(\N)$-variable in $\mathbb{P}^{\mathrm{Leb},\sigma,\t,\mathbb{I}(\mathfrak{t})}$. Now, by Definition \ref{definition:kv4}, we first have the following estimate, where $\mathbf{U}$ denotes the dummy-variable for $\mathbb{P}^{\sigma,\t,\mathbb{I}(\mathfrak{t})}$:
\begin{align}
\mathrm{LHS}\eqref{eq:kv6I} \ \lesssim \ \mathbb{P}^{\sigma,\t,\mathbb{I}(\mathfrak{t})}\{{\textstyle\sup_{\s}}|\partial_{\s}\mathscr{HP}(\s,\mathbf{U};\sigma)| \gtrsim \N^{\gamma_{\mathrm{KV}}}|\mathbb{I}(\mathfrak{t})|^{\frac12}\}, \label{eq:kv6I1}
\end{align}
where the supremum is over $\mathfrak{t}(\mathrm{in})\leq\s\leq\mathfrak{t}(\mathrm{in})+\mathfrak{t}$ and $\mathscr{HP}$ is from Lemma \ref{lemma:kv3}. We now claim that $|\partial_{\s}^{2}\mathscr{HP}(\s,\mathbf{U};\sigma)|\lesssim\N^{200}$. (This follows by definition of $\mathscr{HP}(\s,\mathbf{U};\sigma)$ as a sum of $\mathscr{UP}(\s,\mathbf{U}(\x);\sigma)$ and $|\partial_{\s}^{2}\mathscr{UP}(\s,\mathbf{U}(\x);\sigma)|\lesssim\N^{100}$; see Lemma \ref{lemma:kv3}.) So, by Lemma \ref{lemma:steeasy}, the sup in $\mathbb{P}^{\sigma,\t,\mathbb{I}(\mathfrak{t})}$ in $\mathrm{RHS}\eqref{eq:kv6I1}$ is $\lesssim\N^{-99}$ plus the sup over $\s$ in a very fine discretization of $\mathfrak{t}(\mathrm{in})+[0,\mathfrak{t}]$ of mesh size $\N^{{\mathrm{C}}}$ {for some $\mathrm{C}>0$ large but fixed}. Because $\mathfrak{t}\lesssim\N^{\gamma_{\mathrm{av}}}$, we know said discretization of $\mathfrak{t}(\mathrm{in})+[0,\mathfrak{t}]$ of mesh size $\N^{{\mathrm{C}}}$ has cardinality $\lesssim\N^{{\mathrm{C}}}$. Therefore, by a union bound, we can pull the $\sup_{\s}$ in $\mathrm{RHS}\eqref{eq:kv6I1}$ outside $\mathbb{P}^{\sigma,\t,\mathbb{I}(\mathfrak{t})}$ if we give up a factor of $\N^{{\mathrm{C}}}\lesssim\exp[\N^{\gamma_{\mathrm{av}}}]$. This ultimately turns \eqref{eq:kv6I1} into the following estimate (again, see Lemma \ref{lemma:steeasy}):
\begin{align}
\mathrm{LHS}\eqref{eq:kv6I} \ \lesssim \ \exp[\N^{\gamma_{\mathrm{av}}}]\times{\textstyle\sup_{\s}}\mathbb{P}^{\sigma,\t,\mathbb{I}(\mathfrak{t})}\{|\partial_{\s}\mathscr{HP}(\s,\mathbf{U};\sigma)| \gtrsim \N^{\gamma_{\mathrm{KV}}}|\mathbb{I}(\mathfrak{t})|^{\frac12}\}. \label{eq:kv6I2}
\end{align}
Since $\gamma_{\mathrm{av}}\leq{\mathrm{c}}\gamma_{\mathrm{KV}}$ (see Proposition \ref{prop:kv1}), $\exp[\N^{\gamma_{\mathrm{av}}}]\exp[-\N^{2\gamma_{\mathrm{KV}}/3}]\lesssim\exp[\N^{-\gamma_{\mathrm{KV}}/2}]$. So by \eqref{eq:kv6I2}, it suffices to get
\begin{align}
{\textstyle\sup_{\s}}\mathbb{P}^{\sigma,\t,\mathbb{I}(\mathfrak{t})}\{|\partial_{\s}\mathscr{HP}(\s,\mathbf{U};\sigma)| \gtrsim \N^{\gamma_{\mathrm{KV}}}|\mathbb{I}(\mathfrak{t})|^{\frac12}\} \ \lesssim \ \exp\{\N^{-\frac23\gamma_{\mathrm{KV}}}\}. \label{eq:kv6I3}
\end{align}
By Lemma \ref{lemma:kv3}, the Radon-Nikodym derivative of $\mathbb{P}^{\sigma,\t,\mathbb{I}(\mathfrak{t})}$ with respect to $\mathbb{P}^{\sigma,\s,\mathbb{I}(\mathfrak{t})}$ is $\mathfrak{p}[\sigma,\t,\mathbb{I}(\mathfrak{t})]\mathfrak{p}[\sigma,\s,\mathbb{I}(\mathfrak{t})]^{-1}$. Lemma \ref{lemma:kv3} computes this ratio to be the following exponentiated sum over $\x\in\mathbb{I}(\mathfrak{t})$, which we estimate and then explain:
\begin{align}
\mathfrak{p}[\sigma,\t,\mathbb{I}(\mathfrak{t})]\mathfrak{p}[\sigma,\s,\mathbb{I}(\mathfrak{t})]^{-1} \ &= \ \exp\{{\textstyle\sum_{\x}}[\mathscr{UP}(\s,\mathbf{U}(\x);\sigma)-\mathscr{UP}(\t,\mathbf{U}(\x);\sigma)]\} \ \lesssim \ \exp\{\mathrm{O}(\N^{\gamma_{\mathrm{av}}})\}. \label{eq:kv6I3b}
\end{align}
The identity in \eqref{eq:kv6I3b} is by \eqref{eq:kv3I}. To prove the upper bound in \eqref{eq:kv6I3b}, first note that the middle of \eqref{eq:kv6I3b} is bounded above by $\exp\{\mathrm{O}(|\t-\s|\times|\mathbb{I}(\mathfrak{t})|)\}$. This follows from the time-derivative estimate for $\mathscr{UP}$ in Lemma \ref{lemma:kv3} and the fact that the sum is over $\x\in\mathbb{I}(\mathfrak{t})$. But $\s,\t\in\mathfrak{t}(\mathrm{in})+[0,\mathfrak{t}]$, so $|\t-\s|\lesssim\mathfrak{t}$. It now suffices to recall $\mathfrak{t}|\mathbb{I}(\mathfrak{t})|\lesssim\N^{\gamma_{\mathrm{av}}}$; see the statement of the lemma. So
\begin{align}
&{\textstyle\sup_{\s}}\mathbb{P}^{\sigma,\t,\mathbb{I}(\mathfrak{t})}\{|\partial_{\s}\mathscr{HP}(\s,\mathbf{U};\sigma)| \gtrsim \N^{\gamma_{\mathrm{KV}}}|\mathbb{I}(\mathfrak{t})|^{\frac12}\} \nonumber\\
&\lesssim \ \exp\{\mathrm{O}(\N^{\gamma_{\mathrm{av}}})\}{\textstyle\sup_{\s}}\mathbb{P}^{\sigma,\s,\mathbb{I}(\mathfrak{t})}\{|\partial_{\s}\mathscr{HP}(\s,\mathbf{U};\sigma)| \gtrsim \N^{\gamma_{\mathrm{KV}}}|\mathbb{I}(\mathfrak{t})|^{\frac12}\}. \nonumber
\end{align}
To prove \eqref{eq:kv6I3}, which would complete the proof (as noted right before \eqref{eq:kv6I3}), by the previous display, it suffices to show
\begin{align}
{\textstyle\sup_{\s}}\mathbb{P}^{\sigma,\s,\mathbb{I}(\mathfrak{t})}\{|\partial_{\s}\mathscr{HP}(\s,\mathbf{U};\sigma)| \gtrsim \N^{\gamma_{\mathrm{KV}}}|\mathbb{I}(\mathfrak{t})|^{\frac12}\} \ \lesssim \ \exp\{-\N^{\frac56\gamma_{\mathrm{KV}}}\}. \label{eq:kv6I4}
\end{align}
(Indeed, $\mathrm{RHS}\eqref{eq:kv6I4}$ absorbs the exponential prefactor in the display before \eqref{eq:kv6I4} and therefore yields \eqref{eq:kv6I3}.) We now recall $\mathscr{HP}$ from Lemma \ref{lemma:kv3}. By adding and subtracting $\E^{\sigma,\s}\partial_{\s}\mathscr{UP}(\s,\mathbf{u};\sigma)$ for each $\x\in\mathbb{T}(\N)$ in the definition of $\mathscr{HP}$, in which $\mathbf{u}$ is the dummy-variable for the $\E^{\sigma,\s}$ expectation, we get the following decomposition where the sum is over $\x\in\mathbb{I}(\mathfrak{t})$:
\begin{align}
\partial_{\s}\mathscr{HP}(\s,\mathbf{U};\sigma) \ = \ {\textstyle\sum_{\x}}\{\partial_{\s}\mathscr{UP}(\s,\mathbf{U}(\x);\sigma)-\E^{\sigma,\s}\partial_{\s}\mathscr{UP}(\s,\mathbf{u};\sigma)\}+|\mathbb{I}(\mathfrak{t})|\E^{\sigma,\s}\partial_{\s}\mathscr{UP}(\s,\mathbf{u};\sigma). \label{eq:kv6I5}
\end{align}
We claim the following bound. This first step below follows by \eqref{eq:kv3III}. The second is the equivalence of ensembles bound \eqref{eq:le8I13} but replacing $(0,\mathbb{T}(\N))\mapsto(\sigma,\mathbb{I}(\mathfrak{t}))$. (We claim \eqref{eq:le8I13} with $0$ replaced by $\sigma$ holds. The only role $0$ played is in bounding some moments that we explained after \eqref{eq:le8I13}. But all we need is $|0|\lesssim1$. And, we assumed $|\sigma|\lesssim1$ for the fixed $\sigma$ of this section.)
\begin{align}
|\mathbb{I}(\mathfrak{t})|\E^{\sigma,\s}\partial_{\s}\mathscr{UP}(\s,\mathbf{u};\sigma) \ = \ |\mathbb{I}(\mathfrak{t})|\{\E^{\sigma,\s}\partial_{\s}\mathscr{UP}(\s,\mathbf{u};\sigma)-\E^{\sigma,\s,\mathbb{I}(\mathfrak{t})}\partial_{\s}\mathscr{UP}(\s,\mathbf{u};\sigma)\} \ \lesssim \ 1. \label{eq:kv6I6}
\end{align}
If $|\mathrm{LHS}\eqref{eq:kv6I5}|\gtrsim\N^{\gamma_{\mathrm{KV}}}|\mathbb{I}(\mathfrak{t})|^{1/2}\gg1$, then by the deterministic bound \eqref{eq:kv6I6}, we deduce that the sum over $\x\in\mathbb{I}(\mathfrak{t})$ in \eqref{eq:kv6I5} is $\gtrsim\N^{\gamma_{\mathrm{KV}}}|\mathbb{I}(\mathfrak{t})|^{1/2}$. Thus, to show \eqref{eq:kv6I4}, which would finish this proof, it suffices to show the following instead:
\begin{align}
&{\textstyle\sup_{\s}}\mathbb{P}^{\sigma,\s,\mathbb{I}(\mathfrak{t})}\{|{\textstyle\sum_{\x}}\{\partial_{\s}\mathscr{UP}(\s,\mathbf{U}(\x);\sigma)-\E^{\sigma,\s}\partial_{\s}\mathscr{UP}(\s,\mathbf{u};\sigma)\}| \gtrsim \N^{\gamma_{\mathrm{KV}}}|\mathbb{I}(\mathfrak{t})|^{\frac12}\} \nonumber\\
&\lesssim \ \exp\{-\N^{\frac56\gamma_{\mathrm{KV}}}\}. \label{eq:kv6I7}
\end{align}
For convenience, we define the ``centered" term $\mathscr{CP}(\s,\mathbf{U}(\x);\sigma):=\partial_{\s}\mathscr{UP}(\s,\mathbf{U}(\x);\sigma)-\E^{\sigma,\s}\partial_{\s}\mathscr{UP}(\s,\mathbf{u};\sigma)$. We note $|\mathscr{CP}|\lesssim1$ by Lemma \ref{lemma:kv3}. We also note $\E^{\sigma,\s}\mathscr{CP}(\s,\mathbf{u};\sigma)=0$, where $\mathbf{u}$ is the expectation dummy variable. 

Now comes the ``one-block, two-blocks part" of our proof. This is built on the dyadic sequence $\mathfrak{i}[\mathrm{k}]:=2^{\mathrm{k}}$. With this sequence, define $\mathbf{U}\mapsto\sigma(\mathbf{U},\x,\mathrm{k})$ to be the average of $\mathbf{U}(\y)$ over $\y\in\x+\llbracket0,\mathfrak{i}[\mathrm{k}]-1\rrbracket$. Let $\E^{\sigma(\mathbf{U},\x,\mathrm{k})}$ be {the} canonical ensemble expectation on $\x+\llbracket0,\mathfrak{i}[\mathrm{k}]-1\rrbracket$ with charge density $\sigma(\mathbf{U},\x,\mathrm{k})$ at time $\s$. (We omit $\s$ from notation just for convenience. We take $1\leq\mathrm{k}\leq\mathrm{k}(\infty)$, where $\mathrm{k}(\infty)$ is the smallest integer such that $\mathfrak{i}[\mathrm{k}]=2^{\mathrm{k}}\geq\N^{-\gamma_{\mathrm{KV}}/100}|\mathbb{I}(\mathfrak{t})|$. Indeed, we never want $\mathfrak{i}[\mathrm{k}]$ to exceed $|\mathbb{I}(\mathfrak{t})|$. Otherwise, $\E^{\sigma(\mathbf{U},\x,\mathrm{k})}$ {would} stabilize as a sequence in $\mathrm{k}$.) Next, we define the following terms, in which $\mathrm{k}>1$:
\begin{align}
\mathscr{CP}(\s,\x,\mathbf{U};\sigma,1)&:=\mathscr{CP}(\s,\mathbf{U}(\x);\sigma)-\E^{\sigma(\mathbf{U},\x,1)}\mathscr{CP}(\s,\x,\cdot;\sigma) \label{eq:kv6I8a}\\
\mathscr{CP}(\s,\x,\mathbf{U};\sigma,\mathrm{k})&:=\E^{\sigma(\mathbf{U},\x,\mathrm{k}-1)}\mathscr{CP}(\s,\x,\cdot;\sigma)-\E^{\sigma(\mathbf{U},\x,\mathrm{k})}\mathscr{CP}(\s,\x,\cdot;\sigma). \label{eq:kv6I8b}
\end{align}
A few notes about \eqref{eq:kv6I8a}-\eqref{eq:kv6I8b}. First, on the RHS of \eqref{eq:kv6I8a}-\eqref{eq:kv6I8b}, the $\cdot$ in the expectation denotes an expectation-dummy-variable. (We choose to avoid overloading $\mathbf{U}$ for notation.) Now, we claim the following estimate that we justify afterwards:
\begin{align}
|\E^{\sigma(\mathbf{U},\x,\mathrm{k})}\mathscr{CP}(\s,\x,\cdot;\sigma)| \ \lesssim \ 1\wedge\{\mathfrak{i}[\mathrm{k}]^{-1}+|\sigma(\mathbf{U},\x,\mathrm{k})-\sigma|\}. \label{eq:kv6I8c}
\end{align}
The bound $\mathrm{LHS}\eqref{eq:kv6I8c}\lesssim1$ follows immediately via $|\mathscr{CP}|\lesssim1$; see the paragraph {after \eqref{eq:kv6I7}}. To establish \eqref{eq:kv6I8c}, it suffices to assume that $|\sigma(\mathbf{U},\x,\mathrm{k})-\sigma|\lesssim1$. In this case, the bound \eqref{eq:kv6I8c} (after dropping the $1\wedge$) follows by first replacing $\E^{\sigma(\mathbf{U},\x,\mathrm{k})}$ by the corresponding grand-canonical expectation. (The error we pick up after this is $\lesssim\mathfrak{i}[\mathrm{k}]^{-1}$. This follows just from the bound \eqref{eq:kv6I6} upon replacing $(\mathbb{I}(\mathfrak{t}),\sigma,\partial_{\s}\mathscr{UP})\mapsto(\x+\llbracket0,\mathfrak{i}[\mathrm{k}]-1\rrbracket,\sigma(\mathbf{U},\x,\mathrm{k}),\mathscr{CP})$. Indeed, all \eqref{eq:kv6I6} needs from either $\sigma$ or $\partial_{\s}\mathscr{UP}$ is the bound $|\sigma|+|\partial_{\s}\mathscr{UP}|\lesssim1$. However, this stays true after replacing $(\sigma,\partial_{\s}\mathscr{UP})\mapsto(\sigma(\mathbf{U},\x,\mathrm{k}),\mathscr{CP})$. Indeed, we have assumed $|\sigma(\mathbf{U},\x,\mathrm{k})-\sigma|\lesssim1$ and $|\sigma|\lesssim1$.) It now suffices to Taylor expand the resulting grand-canonical expectation around $\sigma(\mathbf{U},\x,\mathrm{k})\approx\sigma$ as in the proof of Lemma \ref{lemma:ee1}. (See \eqref{eq:ee1I3}, but Taylor expand only up to first order. In doing so, the implied constant in \eqref{eq:kv6I8c} actually picks up the norm in Lemma \ref{lemma:ee1} of $\mathscr{CP}$. But $|\mathscr{CP}|\lesssim1$ and $\sigma(\mathbf{U},\x,\mathrm{k}),\sigma=\mathrm{O}(1)$ by assumption. So, by following the last paragraph in the proofs of Lemmas \ref{lemma:bg23}, \ref{lemma:bg1hl2}, said norm is $\lesssim1$.) We have now explained \eqref{eq:kv6I8c}. 

Now, $|\sigma(\mathbf{U},\x,\mathrm{k})-\sigma|$ is sub-Gaussian with variance parameter $\lesssim\mathfrak{i}[\mathrm{k}]^{-1}$. (Indeed, it is an average of $\mathfrak{i}[\mathrm{k}]$-many random walk bridge steps, where the random walk bridge has drift 0 since we subtracted $\sigma$. By uniform convexity of $\mathscr{U}$ in Assumption \ref{ass:intro8}, the steps have sub-Gaussian distribution. Thus, as in the end of the proof of Lemma \ref{lemma:bg27}, we know that $\sigma(\mathbf{U},\x,\mathrm{k})-\sigma$ is an average of $\mathfrak{i}[\mathrm{k}]$-many sub-Gaussian martingale increments, up to an error of $\mathfrak{i}[\mathrm{k}]^{-1}$ per increment. Now use Azuma to get sub-Gaussianity of $\sigma(\mathbf{U},\x,\mathrm{k})-\sigma$.) From \eqref{eq:kv6I8c} and the previous couple of sentences, we ultimately get that $\mathscr{CP}(\s,\x,\mathbf{U};\sigma,\mathrm{k})$ is sub-Gaussian with variance parameter $\lesssim\mathfrak{i}[\mathrm{k}-1]^{-1}$. (If $\mathrm{k}=1$, it is sub-Gaussian and $\lesssim1$. This is clear because $|\mathscr{CP}(\s,\cdot;\sigma)|\lesssim1$ as we noted after \eqref{eq:kv6I7}.) \eqref{eq:kv6I8c} also gives that $\E^{\sigma(\mathbf{U},\x,\mathrm{k}(\infty))}\mathscr{CP}(\s,\x,\cdot;\sigma)$ is sub-Gaussian with variance parameter $\lesssim\mathfrak{i}[\mathrm{k}(\infty)-1]^{-1}\lesssim\N^{\gamma_{\mathrm{KV}}/100}|\mathbb{I}(\mathfrak{t})|^{-1}$. Finally, by Lemma \ref{lemma:vanishcanonical}, $\mathscr{CP}(\s,\x,\mathbf{U};\sigma,\mathrm{k})$ vanishes with respect to any canonical expectation on its support. 

We now prove \eqref{eq:kv6I7}. By definition of $\mathscr{CP}$ and telescoping sum, the sum inside $\mathbb{P}^{\sigma,\s,\mathbb{I}(\mathfrak{t})}$ in \eqref{eq:kv6I7} is
\begin{align}
= \ {\textstyle\sum_{\x}}\mathscr{CP}(\s,\mathbf{U}(\x);\sigma) \ = \ {\textstyle\sum_{\mathrm{k}=1}^{\mathrm{k}(\infty)}\sum_{\x}}\mathscr{CP}(\s,\x,\mathbf{U};\sigma,\mathrm{k})+{\textstyle\sum_{\x}}\E^{\sigma(\mathbf{U},\x,\mathrm{k}(\infty))}\mathscr{CP}(\s,\x,\cdot;\sigma). \label{eq:kv6I9}
\end{align}
Note $\mathrm{k}(\infty)\lesssim\log\N$, since $\mathrm{k}(\infty)$ is at most the number of dyadic scales needed to hit $|\mathbb{I}(\mathfrak{t})|\leq\N$. So, if $|\eqref{eq:kv6I9}|\gtrsim\N^{\gamma_{\mathrm{KV}}}|\mathbb{I}(\mathfrak{t})|^{1/2}$, then either at least one of the $\mathrm{k}$-summands or the last term in \eqref{eq:kv6I9} must exceed $\gtrsim\N^{\gamma_{\mathrm{KV}}}[\log\N]^{-1}|\mathbb{I}(\mathfrak{t})|^{1/2}\gtrsim\N^{\gamma_{\mathrm{KV}}/2}|\mathbb{I}(\mathfrak{t})|^{1/2}$ in absolute value. By a union bound to account for all $\lesssim\mathrm{k}(\infty)\lesssim\log\N$-many possibilities, we deduce
\begin{align}
&\mathrm{LHS}\eqref{eq:kv6I7} \nonumber\\
&\lesssim \ (\log\N)\times{\textstyle\sup_{\s}\sup_{\mathrm{k}=1}^{\mathrm{k}(\infty)}}\mathbb{P}^{\sigma,\s,\mathbb{I}(\mathfrak{t})}\{|{\textstyle\sum_{\x}}\mathscr{CP}(\s,\x,\mathbf{U};\sigma,\mathrm{k})|\gtrsim\N^{\frac12\gamma_{\mathrm{KV}}}|\mathbb{I}(\mathfrak{t})|^{\frac12}\} \label{eq:kv6I10a} \\
&+ \ (\log\N)\times{\textstyle\sup_{\s}}\mathbb{P}^{\sigma,\s,\mathbb{I}(\mathfrak{t})}\{|{\textstyle\sum_{\x}}\E^{\sigma(\mathbf{U},\x,\mathrm{k}(\infty))}\mathscr{CP}(\s,\x,\cdot;\sigma)|\gtrsim\N^{\frac12\gamma_{\mathrm{KV}}}|\mathbb{I}(\mathfrak{t})|^{\frac12}\}. \label{eq:kv6I10b} \\
&= \ (\log\N)\times{\textstyle\sup_{\s}\sup_{\mathrm{k}=1}^{\mathrm{k}(\infty)}}\mathbb{P}^{\sigma,\s,\mathbb{I}(\mathfrak{t})}\{|\mathbb{I}(\mathfrak{t})|^{-1}|{\textstyle\sum_{\x}}\mathscr{CP}(\s,\x,\mathbf{U};\sigma,\mathrm{k})|\gtrsim\N^{\frac12\gamma_{\mathrm{KV}}}|\mathbb{I}(\mathfrak{t})|^{-\frac12}\} \label{eq:kv6I10c} \\
&+ \ (\log\N)\times{\textstyle\sup_{\s}}\mathbb{P}^{\sigma,\s,\mathbb{I}(\mathfrak{t})}\{|\mathbb{I}(\mathfrak{t})|^{-1}|{\textstyle\sum_{\x}}\E^{\sigma(\mathbf{U},\x,\mathrm{k}(\infty))}\mathscr{CP}(\s,\x,\cdot;\sigma)|\gtrsim\N^{\frac12\gamma_{\mathrm{KV}}}|\mathbb{I}(\mathfrak{t})|^{-\frac12}\}. \label{eq:kv6I10d}
\end{align}
Take any $\mathbb{P}^{\sigma,\s,\mathbb{I}(\mathfrak{t})}$-term in \eqref{eq:kv6I10c}. We will rewrite the average of $\mathscr{CP}(\s,\x,\mathbf{U};\sigma,\mathrm{k})$ over $\x\in\mathbb{I}(\mathfrak{t})$ as follows. By construction in \eqref{eq:kv6I8a}-\eqref{eq:kv6I8b},  the support of $\mathscr{CP}(\s,\x,\mathbf{U};\sigma,\mathrm{k})$ is some discrete interval of length $\mathfrak{i}[\mathrm{k}]$. We first write $\mathbb{I}(\mathfrak{t})$ as a union of ``clusters", where each cluster is of size $\gtrsim|\mathbb{I}(\mathfrak{t})|\mathfrak{i}[\mathrm{k}]^{-1}$, and any pair of points in any common cluster are separated from each other by $2\mathfrak{i}[\mathrm{k}]$. With this construction, the average of $\mathscr{CP}(\s,\x,\mathbf{U};\sigma,\mathrm{k})$ over $\x\in\mathbb{I}(\mathfrak{t})$ is the same thing as an average, over all clusters, of the average of $\mathscr{CP}(\s,\x,\mathbf{U};\sigma,\mathrm{k})$ over all $\x$ in a fixed cluster. As $\mathscr{CP}(\s,\x,\mathbf{U};\sigma,\mathrm{k})$ vanishes with respect to canonical expectations on its support, by Lemma \ref{lemma:le2}, the average of $\mathscr{CP}(\s,\x,\mathbf{U};\sigma,\mathrm{k})$ over any cluster of $\x$ is sub-Gaussian with variance parameter given by multiplying the inverse-cluster-size $|\mathbb{I}(\mathfrak{t})|^{-1}\mathfrak{i}[\mathrm{k}]$ by the variance parameter of $\mathscr{CP}(\s,\x,\mathbf{U};\sigma,\mathrm{k})$, which we recall equals $\lesssim\mathfrak{i}[\mathrm{k}-1]^{-1}$. (More precise explanation can be found in the paragraph prior to \eqref{eq:le8I16}.) As averaging sub-Gaussian variables keeps sub-Gaussian property, we know the average of $\mathscr{CP}(\s,\x,\mathbf{U};\sigma,\mathrm{k})$ over $\x\in\mathbb{I}(\mathfrak{t})$ is sub-Gaussian with variance parameter $\lesssim|\mathbb{I}(\mathfrak{t})|^{-1}\mathfrak{i}[\mathrm{k}]\mathfrak{i}[\mathrm{k}-1]^{-1}=2|\mathbb{I}(\mathfrak{t})|^{-1}$. Thus, $\eqref{eq:kv6I10c}\lesssim\exp[-\upsilon\N^{\gamma_{\mathrm{KV}}}]$ with $\upsilon\gtrsim1$. (In particular, this exponential decay beats the $\log\N$ growth in \eqref{eq:kv6I10c}.) As for \eqref{eq:kv6I10d}, the same argument works. We just use that $\E^{\sigma(\mathbf{U},\x,\mathrm{k}(\infty))}\mathscr{CP}(\s,\x,\cdot;\sigma)$ is sub-Gaussian with variance parameter $\lesssim\N^{\gamma_{\mathrm{KV}}/100}|\mathbb{I}(\mathfrak{t})|^{-1}$. Since the variance parameter has now grown by a factor of $\N^{\gamma_{\mathrm{KV}}/100}$, our estimate for \eqref{eq:kv6I10d} is slightly worse. But at any rate, we can still get $\eqref{eq:kv6I10d}\lesssim\exp[-\upsilon\N^{\gamma_{\mathrm{KV}}}\N^{-\gamma_{\mathrm{KV}}/100}]$. Using the past few sentences and \eqref{eq:kv6I10a}-\eqref{eq:kv6I10d} gives \eqref{eq:kv6I7}. As noted before \eqref{eq:kv6I7}, this gives \eqref{eq:kv6I4} and finishes the proof. 
\end{proof}
We now show the processes in Definition \ref{definition:kv4} are also good proxies for \eqref{eq:glsdeloc}-\eqref{eq:hfloc} in the following ``analytic" sense.
\begin{lemma}\label{lemma:kv7}
 The process $\t\mapsto(\mathbf{J}(\t;\mathbb{I}(\mathfrak{t}),\mathcal{E}),\mathbf{U}^{\t,\cdot}[\mathbb{I}(\mathfrak{t}),\mathcal{E}])$ in {Definition \ref{definition:kv4}} is Markov. Recall $\mathbb{P}^{\mathrm{Leb},\sigma,\t,\mathbb{I}(\mathfrak{t}),\mathcal{E}}=\mathrm{Leb}(\mathbb{S}(\N))\otimes\mathbb{P}^{\sigma,\t,\mathbb{I}(\mathfrak{t}),\mathcal{E}}$ in {Definition \ref{definition:kv4}}. The infinitesimal generator of $\t\mapsto(\mathbf{J}(\t;\mathbb{I}(\mathfrak{t}),\mathcal{E}),\mathbf{U}^{\t,\cdot}[\mathbb{I}(\mathfrak{t}),\mathcal{E}])$ is the time-inhomogeneous operator $\mathscr{L}^{\mathrm{tot}}(\t,\mathbb{I}(\mathfrak{t}))$ from {Definition \ref{definition:le5}} but with (vanishing) Dirichlet boundary conditions on the boundary of $\mathcal{E}[\mathfrak{t}(\mathrm{in}),\mathfrak{t};\mathbb{I}(\mathfrak{t})]$. Finally, the adjoint of this (time $\t$) generator with respect to $\mathbb{P}^{\mathrm{Leb},\sigma,\t,\mathbb{I}(\mathfrak{t}),\mathcal{E}}$ is the differential operator given by the adjoint of $\mathscr{L}^{\mathrm{tot}}(\t,\mathbb{I}(\mathfrak{t}))$ with respect to $\mathbb{P}^{\mathrm{Leb},\sigma,\t,\mathbb{I}(\mathfrak{t})}$, but with (vanishing) Dirichlet boundary conditions on the boundary of $\mathcal{E}[\mathfrak{t}(\mathrm{in}),\mathfrak{t};\mathbb{I}(\mathfrak{t})]$.
\end{lemma}
\begin{proof}
These claims are standard. But we could not find a reference that would make the proof easier. So, we give one complete proof here. The Markov property we claimed follows by the strong Markov property for the SDEs \eqref{eq:hf}-\eqref{eq:glsde}. The infinitesimal Markov generator is computed by using standard Ito theory. (Indeed, we still have an Ito formula for the processes constructed in Definition \ref{definition:kv4}. It comes from the Ito formula for \eqref{eq:glsdeloc}-\eqref{eq:hfloc}.) We now compute the adjoint. Take any smooth functions $\mathsf{H},\mathsf{F}\in\mathscr{C}^{\infty}(\mathbb{S}(\N)\times\mathcal{E}[\mathfrak{t}(\mathrm{in}),\mathfrak{t};\mathbb{I}(\mathfrak{t})])$ with compact support and with vanishing data on the boundary of $\mathcal{E}[\mathfrak{t}(\mathrm{in}),\mathfrak{t};\mathbb{I}(\mathfrak{t})]$. We have 
\begin{align}
\E^{\mathrm{Leb},\sigma,\t,\mathbb{I}(\mathfrak{t}),\mathcal{E}}\mathsf{H}\cdot\mathscr{L}^{\mathrm{tot}}(\t,\mathbb{I}(\mathfrak{t}))\mathsf{F} \ &\propto \ \E^{\mathrm{Leb},\sigma,\t,\mathbb{I}(\mathfrak{t})}(\mathbf{1}\{\mathcal{E}[\mathfrak{t}(\mathrm{in}),\mathfrak{t};\mathbb{I}(\mathfrak{t})]\}\cdot\mathsf{H})\cdot\mathscr{L}^{\mathrm{tot}}(\t,\mathbb{I}(\mathfrak{t}))\mathsf{F} \label{eq:kv7I1a} \\
&= \ \E^{\mathrm{Leb},\sigma,\t,\mathbb{I}(\mathfrak{t})}\mathscr{L}^{\mathrm{tot}}(\t,\mathbb{I}(\mathfrak{t}))^{\ast}(\mathbf{1}\{\mathcal{E}[\mathfrak{t}(\mathrm{in}),\mathfrak{t};\mathbb{I}(\mathfrak{t})]\}\cdot\mathsf{H})\cdot\mathsf{F}. \label{eq:kv7I1b}
\end{align}
\eqref{eq:kv7I1a} follows by the definition of $\E^{\mathrm{Leb},\sigma,\t,\mathbb{I}(\mathfrak{t}),\mathcal{E}}$ as $\E^{\mathrm{Leb},\sigma,\t,\mathbb{I}(\mathfrak{t})}$ but conditioning on $\mathbb{S}(\N)\times\mathcal{E}[\mathfrak{t}(\mathrm{in}),\mathfrak{t};\mathbb{I}(\mathfrak{t})]$; see Definition \ref{definition:kv4}. \eqref{eq:kv7I1b} follows by definition of $\mathscr{L}^{\mathrm{tot}}(\t,\mathbb{I}(\mathfrak{t}))^{\ast}$ as $\E^{\mathrm{Leb},\sigma,\t,\mathbb{I}(\mathfrak{t})}$-adjoint of $\mathscr{L}^{\mathrm{tot}}(\t,\mathbb{I}(\mathfrak{t}))$. We already computed $\mathscr{L}^{\mathrm{tot}}(\t,\mathbb{I}(\mathfrak{t}))^{\ast}$ in the proof of Lemma \ref{lemma:le7}. It is a second-order differential operator in the $\mathbb{S}(\N)$-variable plus a second-order differential operator in the $\mathbb{H}^{\sigma,\mathbb{I}(\mathfrak{t})}$-variable; see right before \eqref{eq:le7II0}. (The support of $\E^{\mathrm{Leb},\sigma,\t,\mathbb{I}(\mathfrak{t})}$ in \eqref{eq:kv7I1b} is $\mathbb{S}(\N)\times\mathbb{H}^{\sigma,\mathbb{I}(\mathfrak{t})}$.) We now claim that
\begin{align}
&\E^{\mathrm{Leb},\sigma,\t,\mathbb{I}(\mathfrak{t})}\mathscr{L}^{\mathrm{tot}}(\t,\mathbb{I}(\mathfrak{t}))^{\ast}(\mathbf{1}\{\mathcal{E}[\mathfrak{t}(\mathrm{in}),\mathfrak{t};\mathbb{I}(\mathfrak{t})]\}\cdot\mathsf{H})\cdot\mathsf{F} \nonumber\\
&= \ \E^{\mathrm{Leb},\sigma,\t,\mathbb{I}(\mathfrak{t})}(\mathbf{1}\{\mathcal{E}[\mathfrak{t}(\mathrm{in}),\mathfrak{t};\mathbb{I}(\mathfrak{t})]\}\cdot\mathscr{L}^{\mathrm{tot}}(\t,\mathbb{I}(\mathfrak{t}))^{\ast}\mathsf{H})\cdot\mathsf{F}. \label{eq:kv7I1c}
\end{align}
Plugging the identity \eqref{eq:kv7I1c} into \eqref{eq:kv7I1a}-\eqref{eq:kv7I1b} finishes the proof (that is, upon rewriting $\E^{\mathrm{Leb},\sigma,\t,\mathbb{I}(\mathfrak{t})}\mathbf{1}\{\mathcal{E}[\mathfrak{t}(\mathrm{in}),\mathfrak{t};\mathbb{I}(\mathfrak{t})]\}\propto\E^{\mathrm{Leb},\sigma,\t,\mathbb{I}(\mathfrak{t}),\mathcal{E}}$, where this proportionality constant cancels the one in \eqref{eq:kv7I1a}). So, we are left to show \eqref{eq:kv7I1c}. By the Leibniz rule (and our representation of $\mathscr{L}^{\mathrm{tot}}(\t,\mathbb{I}(\mathfrak{t}))^{\ast}$ as a second-order operator), \eqref{eq:kv7I1c} holds if we put in $\mathrm{RHS}\eqref{eq:kv7I1c}$ terms of the form
{\fontsize{9.5}{12}
\begin{align}
&\E^{\mathrm{Leb},\sigma,\t,\mathbb{I}(\mathfrak{t})}(\mathsf{H}\mathsf{F}\times\mathscr{O}^{(2)}\mathbf{1}\{\mathcal{E}[\mathfrak{t}(\mathrm{in}),\mathfrak{t};\mathbb{I}(\mathfrak{t})]\}\nonumber\\
&+\mathsf{H}\mathsf{F}\times\mathscr{O}^{(1)}\mathbf{1}\{\mathcal{E}[\mathfrak{t}(\mathrm{in}),\mathfrak{t};\mathbb{I}(\mathfrak{t})]\}\nonumber\\
&+ \mathsf{F}\times\mathscr{O}^{(1),1}\mathsf{H}\times\mathscr{O}^{(1),2}\mathbf{1}\{\mathcal{E}[\mathfrak{t}(\mathrm{in}),\mathfrak{t};\mathbb{I}(\mathfrak{t})]\}), \label{eq:kv7I2}
\end{align}
}{where $\mathscr{O}^{(2)}$ is a second-order differential operator, and where $\mathscr{O}^{(1)},\mathscr{O}^{(1),1},\mathscr{O}^{(1),2}$ are first-order differential operators. Any first-order differential operator acting on $\mathbf{1}\{\mathcal{E}[\mathfrak{t}(\mathrm{in}),\mathfrak{t};\mathbb{I}(\mathfrak{t})]\}$ is a measure supported at the boundary of $\mathcal{E}[\mathfrak{t}(\mathrm{in}),\mathfrak{t};\mathbb{I}(\mathfrak{t})]$. (Said measure can certainly equal zero. For example, it vanishes if we differentiate with respect to the $\mathbb{S}(\N)$-variable. Otherwise, this is true of the indicator function of any set with a smooth boundary.) Because $\mathsf{F},\mathsf{H}$ vanish on such a set, the last two terms in \eqref{eq:kv7I2} vanish. For the first term, integrate-by-parts. This rewrites the first term in \eqref{eq:kv7I2} as expectation of first derivatives of $\mathbf{1}\{\mathcal{E}[\mathfrak{t}(\mathrm{in}),\mathfrak{t};\mathbb{I}(\mathfrak{t})]\}$ times linear combinations of $\mathsf{H}\mathsf{F}$ and its first derivatives. The Leibniz rule implies that a first-order operator acting on $\mathsf{H}\mathsf{F}$ has a factor of $\mathsf{H}$ and/or $\mathsf{F}$. Therefore, when we multiply it by a first-derivative of $\mathbf{1}\{\mathcal{E}[\mathfrak{t}(\mathrm{in}),\mathfrak{t};\mathbb{I}(\mathfrak{t})]\}$, we get zero. This paragraph shows that \eqref{eq:kv7I2} is zero. As noted before \eqref{eq:kv7I2}, we deduce \eqref{eq:kv7I1c}, so we are done.}
\end{proof}
\subsection{Stability of local equilibrium}
We recall the process $\t\mapsto(\mathbf{J}(\t,\inf\mathbb{I}(\mathfrak{t});\mathbb{I}(\mathfrak{t})),\mathbf{U}^{\t,\cdot}[\mathbb{I}(\mathfrak{t})])$ has law at $\t=\mathfrak{t}(\mathrm{in})$ distributed as $\delta[0]\otimes\mathbb{P}^{\sigma,\mathfrak{t}(\mathrm{in}),\mathbb{I}(\mathfrak{t})}$; see Proposition \ref{prop:kv1}. (For Proposition \ref{prop:kv1}, we will anyway change measure for $\mathbf{J}(\t,\inf\mathbb{I}(\mathfrak{t});\mathbb{I}(\mathfrak{t}))$ in order to give it $\mathrm{Leb}(\mathbb{S}(\N))$ initial data; see after Proposition \ref{prop:kv1}.) Now, our goal is the following ``stability". For sake of a discussion, assume that the law of $(\mathbf{J}(\t,\inf\mathbb{I}(\mathfrak{t});\mathbb{I}(\mathfrak{t})),\mathbf{U}^{\t,\cdot}[\mathbb{I}(\mathfrak{t})])$ at $\t=\mathfrak{t}(\mathrm{in})$ actually equals $\mathbb{P}^{\mathrm{Leb},\sigma,\mathfrak{t}(\mathrm{in}),\mathbb{I}(\mathfrak{t})}$; see Definition \ref{definition:le5}. Given any sufficiently small (i.e. ``local") space-time, the local equilibrium $\mathbb{P}^{\mathrm{Leb},\sigma,\t,\mathbb{I}(\mathfrak{t})}$ is approximately constant in time. For sufficiently small space-time, we thus expect that the law of $\t\mapsto(\mathbf{J}(\t,\inf\mathbb{I}(\mathfrak{t});\mathbb{I}(\mathfrak{t})),\mathbf{U}^{\t,\cdot}[\mathbb{I}(\mathfrak{t})])$ is somewhat ``close" to $\mathbb{P}^{\mathrm{Leb},\sigma,\t,\mathbb{I}(\mathfrak{t})}$. (Indeed, if $\t\mapsto\mathscr{U}(\t,\cdot)$ is constant in $\t$ then $\mathbb{P}^{\mathrm{Leb},\sigma,\t,\mathbb{I}(\mathfrak{t})}$ is invariant for $\t\mapsto(\mathbf{J}(\t,\inf\mathbb{I}(\mathfrak{t});\mathbb{I}(\mathfrak{t})),\mathbf{U}^{\t,\cdot}[\mathbb{I}(\mathfrak{t})])$.) We give two versions of this. The first is pointwise bounds on heat kernels for $(\mathbf{J}(\t,\inf\mathbb{I}(\mathfrak{t});\mathbb{I}(\mathfrak{t})),\mathbf{U}^{\t,\cdot}[\mathbb{I}(\mathfrak{t})])$ and the stopped processes in Definition \ref{definition:kv4}. In addition to this, we also derive a collection of energy and gradient bounds for heat kernels. We emphasize that these estimates are PDE bounds. The first bound, namely \eqref{eq:kv8I}, results from a Gronwall argument.
\begin{lemma}\label{lemma:kv8}
 Fix $\mathfrak{t}(\mathrm{in}),\mathfrak{t}\geq0$ and discrete interval $\mathbb{I}\subseteq\mathbb{T}(\N)$. Suppose $\mathfrak{t}|\mathbb{I}(\mathfrak{t})|\lesssim\N^{\gamma_{\mathrm{av}}}$; see {Proposition \ref{prop:kv1}} for $\gamma_{\mathrm{av}}$. Suppose the data of $(\mathbf{J}(\t,\inf\mathbb{I}(\mathfrak{t});\mathbb{I}(\mathfrak{t})),\mathbf{U}^{\t,\cdot}[\mathbb{I}(\mathfrak{t})])$ at $\t=\mathfrak{t}(\mathrm{in})$ is distributed as $\mathbb{P}^{\mathrm{Leb},\sigma,\mathfrak{t}(\mathrm{in}),\mathbb{I}(\mathfrak{t})}$. Take $\t\in\mathfrak{t}(\mathrm{in})+[0,\mathfrak{t}]$. Let $\mathbb{P}[\mathfrak{t}(\mathrm{in}),\t]$ be the law of $(\mathbf{J}(\t,\inf\mathbb{I}(\mathfrak{t});\mathbb{I}(\mathfrak{t})),\mathbf{U}^{\t,\cdot}[\mathbb{I}(\mathfrak{t})])$. Let $\mathfrak{p}[\mathfrak{t}(\mathrm{in}),\t]$ be the Radon-Nikodym derivative of $\mathbb{P}[\mathfrak{t}(\mathrm{in}),\t]$ with respect to $\mathbb{P}^{\mathrm{Leb},\sigma,\t,\mathbb{I}(\mathfrak{t})}$. Let $\|\|_{\infty}$ denote the sup-norm over $\mathbb{S}(\N)\times\mathbb{H}^{\sigma,\mathbb{I}(\mathfrak{t})}$. {Then we} have the following estimate uniformly over $\t\in\mathfrak{t}(\mathrm{in})+[0,\mathfrak{t}]$:
\begin{align}
\|\mathfrak{p}[\mathfrak{t}(\mathrm{in}),\t]\|_{\infty}+\|\mathfrak{p}[\mathfrak{t}(\mathrm{in}),\t]^{-1}\|_{\infty} \ \lesssim \ \exp[\mathrm{O}(\N^{\gamma_{\mathrm{av}}})]. \label{eq:kv8I}
\end{align}
Suppose that the stopped process $\t\mapsto(\mathbf{J}(\t;\mathbb{I}(\mathfrak{t}),\mathcal{E}),\mathbf{U}^{\t,\cdot}[\mathbb{I}(\mathfrak{t}),\mathcal{E}])$ has $\t=\mathfrak{t}(\mathrm{in})$ data distributed according to $\mathbb{P}^{\mathrm{Leb},\sigma,\mathfrak{t}(\mathrm{in}),\mathbb{I}(\mathfrak{t}),\mathcal{E}}$; see {Definition \ref{definition:kv4}}. Let $\mathbb{P}[\mathfrak{t}(\mathrm{in}),\t;\mathcal{E}]$ be the law of this process at time $\t\in\mathfrak{t}(\mathrm{in})+[0,\mathfrak{t}]$. Let $\mathfrak{p}[\mathfrak{t}(\mathrm{in}),\t;\mathcal{E}]$ be the Radon-Nikodym derivative of $\mathbb{P}[\mathfrak{t}(\mathrm{in}),\t;\mathcal{E}]$ with respect to $\mathbb{P}^{\mathrm{Leb},\sigma,\t,\mathbb{I}(\mathfrak{t}),\mathcal{E}}$. Additionally assume $|\mathbb{I}(\mathfrak{t})|\gtrsim\N^{1/10}$. {Then} 
\begin{align}
\|\mathfrak{p}[\mathfrak{t}(\mathrm{in}),\t;\mathcal{E}]\|_{\infty}+\|\mathfrak{p}[\mathfrak{t}(\mathrm{in}),\t;\mathcal{E}]^{-1}\|_{\infty} \ \lesssim \ 1. \label{eq:kv8II}
\end{align}
\end{lemma}
\begin{proof}
We first claim the following Kolmogorov forward equation that we justify after:
\begin{align}
\partial_{\t}\mathfrak{p}[\mathfrak{t}(\mathrm{in}),\t] \ = \ \mathscr{L}^{\mathrm{tot}}(\t,\mathbb{I}(\mathfrak{t}))^{\ast}\mathfrak{p}[\mathfrak{t}(\mathrm{in}),\t]+\mathfrak{p}[\mathfrak{t}(\mathrm{in}),\t]\partial_{\t}\mathscr{HP}(\t,\mathbf{U};\sigma). \label{eq:kv8I1}
\end{align}
Above, $\mathscr{L}^{\mathrm{tot}}(\t,\mathbb{I}(\mathfrak{t}))^{\ast}$ denotes the adjoint of $\mathscr{L}^{\mathrm{tot}}(\t,\mathbb{I}(\mathfrak{t}))$ (see Definition \ref{definition:le5}) with respect to $\mathbb{P}^{\mathrm{Leb},\sigma,\t,\mathbb{I}(\mathfrak{t})}$. See Lemma \ref{lemma:kv3} for $\mathscr{HP}$. Indeed, \eqref{eq:kv8I1} follows by a calculation that is similar to \eqref{eq:le8II2a}-\eqref{eq:le8II2b}. More precisely, the first term in $\mathrm{RHS}\eqref{eq:kv8I1}$ is the Kolmogorov forward equation. To justify the last term in \eqref{eq:kv8I1}, note $\mathfrak{p}[\mathfrak{t}(\mathrm{in}),\t]$ is a density with respect to a time-dependent measure $\mathbb{P}^{\mathrm{Leb},\sigma,\t,\mathbb{I}(\mathfrak{t})}=\mathrm{Leb}(\mathbb{S}(\N))\otimes\mathbb{P}^{\sigma,\t,\mathbb{I}(\mathfrak{t})}$. The Lebesgue density of this measure (with respect to $\mathrm{Leb}(\mathbb{S}(\N))\otimes\mathrm{Leb}(\mathbb{H}^{\sigma,\mathbb{I}(\mathfrak{t})})$) is the density of $\mathbb{P}^{\sigma,\t,\mathbb{I}(\mathfrak{t})}$ with respect to $\mathrm{Leb}(\mathbb{H}^{\sigma,\mathbb{I}(\mathfrak{t})})$. (The $\mathrm{Leb}(\mathbb{S}(\N))$-factors cancel.) But the density of $\mathbb{P}^{\sigma,\t,\mathbb{I}(\mathfrak{t})}$ with respect to $\mathrm{Leb}(\mathbb{H}^{\sigma,\mathbb{I}(\mathfrak{t})})$ is computed in \eqref{eq:kv3I}. Like \eqref{eq:le8II2a}-\eqref{eq:le8II2b}, \eqref{eq:kv8I1} would hold if $\partial_{\t}\mathscr{HP}(\t,\mathbf{U};\sigma)$ equals $-\log\eqref{eq:kv3I}$; this can be checked easily. We now show \eqref{eq:kv8I}. By Lemma \ref{lemma:kv3}, we know $\partial_{\t}\mathscr{HP}(\t,\mathbf{U};\sigma)$ is a sum over $\x\in\mathbb{I}(\mathfrak{t})$ of $\partial_{\t}\mathscr{UP}(\t,\mathbf{U}(\x);\sigma)$, each of which are $\lesssim1$ with probability 1. Thus, $\partial_{\t}\mathscr{HP}(\t,\mathbf{U};\sigma)\lesssim|\mathbb{I}(\mathfrak{t})|$. Parabolic max-min principles (or Feynman-Kac) for \eqref{eq:kv8I1} therefore get the following string of estimates for $\t\in\mathfrak{t}(\mathrm{in})+[0,\mathfrak{t}]$, where the last bound below follows by assumption:
\begin{align}
\mathrm{LHS}\eqref{eq:kv8I} \ \lesssim \ \exp\{{\textstyle\int_{\mathfrak{t}(\mathrm{in})}^{\t}}\mathrm{O}(|\mathbb{I}(\mathfrak{t})|)\} \ \lesssim \ \exp[\mathrm{O}(\mathfrak{t}|\mathbb{I}(\mathfrak{t})|)] \ \lesssim \ \mathrm{RHS}\eqref{eq:kv8I}. \label{eq:kv8I2}
\end{align}
It remains to get \eqref{eq:kv8II}. We first unfold the reference measure $\mathbb{P}^{\mathrm{Leb},\sigma,\t,\mathbb{I}(\mathfrak{t}),\mathcal{E}}$ for $\mathfrak{p}[\mathfrak{t}(\mathrm{in}),\t;\mathcal{E}]$. By Definition \ref{definition:kv4}, this measure is just $\mathbb{P}^{\mathrm{Leb},\sigma,\t,\mathbb{I}(\mathfrak{t}),\mathcal{E}}=\mathrm{Leb}(\mathbb{S}(\N))\otimes\mathbb{P}^{\sigma,\t,\mathbb{I}(\mathfrak{t}),\mathcal{E}}$. So its density with respect to $\mathrm{Leb}(\mathbb{S}(\N))\otimes\mathrm{Leb}(\mathcal{E}[\mathfrak{t}(\mathrm{in}),\mathfrak{t};\mathbb{I}(\mathfrak{t})])$ is the density of $\mathbb{P}^{\sigma,\t,\mathbb{I}(\mathfrak{t}),\mathcal{E}}$ with respect to $\mathrm{Leb}(\mathcal{E}[\mathfrak{t}(\mathrm{in}),\mathfrak{t};\mathbb{I}(\mathfrak{t})])$, where $\mathrm{Leb}(\mathcal{E}[\mathfrak{t}(\mathrm{in}),\mathfrak{t};\mathbb{I}(\mathfrak{t})])$ denotes the Lebesgue measure on $\mathbb{H}^{\sigma,\mathbb{I}(\mathfrak{t})}$ restricted to $\mathcal{E}[\mathfrak{t}(\mathrm{in}),\mathfrak{t};\mathbb{I}(\mathfrak{t})]\subseteq\mathbb{H}^{\sigma,\mathbb{I}(\mathfrak{t})}$. Recall $\mathbb{P}^{\sigma,\t,\mathbb{I}(\mathfrak{t}),\mathcal{E}}$ is $\mathbb{P}^{\sigma,\t,\mathbb{I}(\mathfrak{t})}$ conditioned on $\mathcal{E}[\mathfrak{t}(\mathrm{in}),\mathfrak{t};\mathbb{I}(\mathfrak{t})]\subseteq\mathbb{H}^{\sigma,\mathbb{I}(\mathfrak{t})}$. So on $\mathcal{E}[\mathfrak{t}(\mathrm{in}),\mathfrak{t};\mathbb{I}(\mathfrak{t})]$, the density of $\mathbb{P}^{\sigma,\t,\mathbb{I}(\mathfrak{t}),\mathcal{E}}$ with respect to $\mathrm{Leb}(\mathcal{E}[\mathfrak{t}(\mathrm{in}),\mathfrak{t};\mathbb{I}(\mathfrak{t})])$ is the density of $\mathbb{P}^{\sigma,\t,\mathbb{I}(\mathfrak{t})}$ with respect to $\mathrm{Leb}(\mathbb{H}^{\sigma,\mathbb{I}(\mathfrak{t})})$ divided by the probability under $\mathbb{P}^{\sigma,\t,\mathbb{I}(\mathfrak{t})}$ of $\mathcal{E}[\mathfrak{t}(\mathrm{in}),\mathfrak{t};\mathbb{I}(\mathfrak{t})]$. We ultimately get the following, where $\mathfrak{p}[\sigma,\t,\mathbb{I}(\mathfrak{t}),\mathcal{E}]$ is the density of $\mathbb{P}^{\mathrm{Leb},\sigma,\t,\mathbb{I}(\mathfrak{t}),\mathcal{E}}$ with respect to $\mathrm{Leb}(\mathbb{S}(\N))\otimes\mathrm{Leb}(\mathcal{E}[\mathfrak{t}(\mathrm{in}),\mathfrak{t};\mathbb{I}(\mathfrak{t})])$, and $\mathfrak{p}[\sigma,\t,\mathbb{I}(\mathfrak{t})]$ is {the} density of $\mathbb{P}^{\mathrm{Leb},\sigma,\t,\mathbb{I}(\mathfrak{t})}$ with respect to $\mathrm{Leb}(\mathbb{S}(\N))\otimes\mathrm{Leb}(\mathbb{H}^{\sigma,\mathbb{I}(\mathfrak{t})})$:
\begin{align}
\mathfrak{p}[\sigma,\t,\mathbb{I}(\mathfrak{t}),\mathcal{E}] \ = \ \mathbf{1}\{\mathcal{E}[\mathfrak{t}(\mathrm{in}),\mathfrak{t};\mathbb{I}(\mathfrak{t})]\}\times[\mathbb{P}^{\sigma,\t,\mathbb{I}(\mathfrak{t})}\{\mathcal{E}[\mathfrak{t}(\mathrm{in}),\mathfrak{t};\mathbb{I}(\mathfrak{t})]\}]^{-1}\mathfrak{p}[\sigma,\t,\mathbb{I}(\mathfrak{t})]. \label{eq:kv8II1}
\end{align}
(One can directly check that {the} $\mathrm{Leb}(\mathbb{S}(\N))$-factors in $\mathbb{P}^{\mathrm{Leb},\sigma,\t,\mathbb{I}(\mathfrak{t}),\mathcal{E}}=\mathrm{Leb}(\mathbb{S}(\N))\otimes\mathbb{P}^{\sigma,\t,\mathbb{I}(\mathfrak{t}),\mathcal{E}}$ and $\mathrm{Leb}(\mathbb{S}(\N))\otimes\mathrm{Leb}(\mathcal{E}[\mathfrak{t}(\mathrm{in}),\mathfrak{t};\mathbb{I}(\mathfrak{t})])$ {are} cancelled, and $\mathrm{RHS}\eqref{eq:kv8II1}$ is a probability density with respect to $\mathrm{Leb}(\mathcal{E}[\mathfrak{t}(\mathrm{in}),\mathfrak{t};\mathbb{I}(\mathfrak{t})])$.) Now, by ultimately the same calculation that gave \eqref{eq:kv8I1}, we claim the following (where $\ast$ in \eqref{eq:kv8II2} now denotes $\mathbb{P}^{\mathrm{Leb},\sigma,\t,\mathbb{I}(\mathfrak{t}),\mathcal{E}}$-adjoint):
\begin{align}
\partial_{\t}\mathfrak{p}[\mathfrak{t}(\mathrm{in}),\t;\mathcal{E}] \ = \ \mathscr{L}^{\mathrm{tot}}(\t,\mathbb{I}(\mathfrak{t}))^{\ast}\mathfrak{p}[\mathfrak{t}(\mathrm{in}),\t;\mathcal{E}]-\mathfrak{p}[\mathfrak{t}(\mathrm{in}),\t;\mathcal{E}]\partial_{\t}\log\mathfrak{p}[\sigma,\t,\mathbb{I}(\mathfrak{t}),\mathcal{E}]. \label{eq:kv8II2}
\end{align}
We use \eqref{eq:kv8II1} to compute the log in \eqref{eq:kv8II2}. Note $\mathfrak{p}[\mathfrak{t}(\mathrm{in}),\t;\mathcal{E}]$ is supported on $\mathcal{E}[\mathfrak{t}(\mathrm{in}),\mathfrak{t};\mathbb{I}(\mathfrak{t})]$. Therefore, to compute the last term in \eqref{eq:kv8II2}, it suffices to forget the indicator of $\mathcal{E}[\mathfrak{t}(\mathrm{in}),\mathfrak{t};\mathbb{I}(\mathfrak{t})]$ in $\mathrm{RHS}\eqref{eq:kv8II1}$. We now get
\begin{align}
&\mathfrak{p}[\mathfrak{t}(\mathrm{in}),\t;\mathcal{E}]\partial_{\t}\log\mathfrak{p}[\sigma,\t,\mathbb{I}(\mathfrak{t}),\mathcal{E}] \nonumber\\
&= \ \mathfrak{p}[\mathfrak{t}(\mathrm{in}),\t;\mathcal{E}]\partial_{\t}\log\mathfrak{p}[\sigma,\t,\mathbb{I}(\mathfrak{t})]-\mathfrak{p}[\mathfrak{t}(\mathrm{in}),\t;\mathcal{E}]\partial_{\t}\log\mathbb{P}^{\sigma,\t,\mathbb{I}(\mathfrak{t})}\{\mathcal{E}[\mathfrak{t}(\mathrm{in}),\mathfrak{t};\mathbb{I}(\mathfrak{t})]\} \label{eq:kv8II3a} \\
&= \ -\mathfrak{p}[\mathfrak{t}(\mathrm{in}),\t,\mathcal{E}]\partial_{\t}\mathscr{HP}(\t,\mathbf{U};\sigma)-\mathfrak{p}[\mathfrak{t}(\mathrm{in}),\t;\mathcal{E}]\partial_{\t}\log\mathbb{P}^{\sigma,\t,\mathbb{I}(\mathfrak{t})}\{\mathcal{E}[\mathfrak{t}(\mathrm{in}),\mathfrak{t};\mathbb{I}(\mathfrak{t})]\}. \label{eq:kv8II3b}
\end{align}
\eqref{eq:kv8II3b} follows by \eqref{eq:kv3I}. Since $\mathfrak{p}[\mathfrak{t}(\mathrm{in}),\t;\mathcal{E}]$ is supported on $\mathcal{E}[\mathfrak{t}(\mathrm{in}),\mathfrak{t};\mathbb{I}(\mathfrak{t})]$, by Definition \ref{definition:kv4}, we know the first term in \eqref{eq:kv8II3b} is $\lesssim\N^{\gamma_{\mathrm{KV}}}|\mathbb{I}(\mathfrak{t})|^{1/2}\mathfrak{p}[\mathfrak{t}(\mathrm{in}),\t,\mathcal{E}]$ (in absolute value). On the other hand, we claim the following {(in which $\mathrm{D}>0$ is large but fixed)}: 
\begin{align}
&\partial_{\t}\log\mathbb{P}^{\sigma,\t,\mathbb{I}(\mathfrak{t})}\{\mathcal{E}[\mathfrak{t}(\mathrm{in}),\mathfrak{t};\mathbb{I}(\mathfrak{t})]\} \nonumber\\
&= \ (\mathbb{P}^{\sigma,\t,\mathbb{I}(\mathfrak{t})}\{\mathcal{E}[\mathfrak{t}(\mathrm{in}),\mathfrak{t};\mathbb{I}(\mathfrak{t})]\})^{-1}\partial_{\t}\mathbb{P}^{\sigma,\t,\mathbb{I}(\mathfrak{t})}\{\mathcal{E}[\mathfrak{t}(\mathrm{in}),\mathfrak{t};\mathbb{I}(\mathfrak{t})]\} \label{eq:kv8II4a} \\
&= \ (\mathbb{P}^{\sigma,\t,\mathbb{I}(\mathfrak{t})}\{\mathcal{E}[\mathfrak{t}(\mathrm{in}),\mathfrak{t};\mathbb{I}(\mathfrak{t})]\})^{-1}\partial_{\t}(1-\mathbb{P}^{\sigma,\t,\mathbb{I}(\mathfrak{t})}\{\mathcal{E}[\mathfrak{t}(\mathrm{in}),\mathfrak{t};\mathbb{I}(\mathfrak{t})]^{\mathrm{C}}\})\label{eq:kv8II4b} \\
&= \ -(\mathbb{P}^{\sigma,\t,\mathbb{I}(\mathfrak{t})}\{\mathcal{E}[\mathfrak{t}(\mathrm{in}),\mathfrak{t};\mathbb{I}(\mathfrak{t})]\})^{-1}\partial_{\t}{\textstyle\int_{\mathcal{E}[\mathfrak{t}(\mathrm{in}),\mathfrak{t};\mathbb{I}(\mathfrak{t})]^{\mathrm{C}}}}\mathfrak{p}[\sigma,\t,\mathbb{I}(\mathfrak{t})] \label{eq:kv8II4c} \\
&= \ (\mathbb{P}^{\sigma,\t,\mathbb{I}(\mathfrak{t})}\{\mathcal{E}[\mathfrak{t}(\mathrm{in}),\mathfrak{t};\mathbb{I}(\mathfrak{t})]\})^{-1}{\textstyle\int_{\mathcal{E}[\mathfrak{t}(\mathrm{in}),\mathfrak{t};\mathbb{I}(\mathfrak{t})]^{\mathrm{C}}}}\partial_{\t}\mathscr{HP}(\t,\mathbf{U};\sigma)\mathfrak{p}[\sigma,\t,\mathbb{I}(\mathfrak{t})] \label{eq:kv8II4d} \\
&\lesssim \ \N^{200}(\mathbb{P}^{\sigma,\t,\mathbb{I}(\mathfrak{t})}\{\mathcal{E}[\mathfrak{t}(\mathrm{in}),\mathfrak{t};\mathbb{I}(\mathfrak{t})]\})^{-1}\mathbb{P}^{\sigma,\t,\mathbb{I}(\mathfrak{t})}\{\mathcal{E}[\mathfrak{t}(\mathrm{in}),\mathfrak{t};\mathbb{I}(\mathfrak{t})]^{\mathrm{C}}\} \ \lesssim \ \N^{-{\mathrm{D}}}. \label{eq:kv8II4e} 
\end{align}
\eqref{eq:kv8II4a} follows by calculus. \eqref{eq:kv8II4b} follows by elementary probability. \eqref{eq:kv8II4c} follows since $\mathfrak{p}[\sigma,\t,\mathbb{I}(\mathfrak{t})]$ is the Lebesgue density of $\mathbb{P}^{\sigma,\t,\mathbb{I}(\mathfrak{t})}$; see Lemma \ref{lemma:kv3}. \eqref{eq:kv8II4d} follows by pulling $\partial_{\t}$ inside the integral in \eqref{eq:kv8II4c} and then using \eqref{eq:kv3I}. The first estimate in \eqref{eq:kv8II4e} follows by $|\partial_{\t}\mathscr{HP}|\lesssim\N^{\mathrm{O}(1)}$. (Indeed, it is a sum of $|\mathbb{I}(\mathfrak{t})|$-many $\partial_{\t}\mathscr{UP}$-terms, and $|\partial_{\t}\mathscr{UP}|\lesssim1$ with probability 1; see Lemma \ref{lemma:kv3}. We then use, again, that $\mathfrak{p}[\sigma,\t,\mathbb{I}(\mathfrak{t})]$ is Lebesgue density of $\mathbb{P}^{\sigma,\t,\mathbb{I}(\mathfrak{t})}$.) The last bound in \eqref{eq:kv8II4e} follows via Lemma \ref{lemma:kv6}. By \eqref{eq:kv8II4a}-\eqref{eq:kv8II4e} and the paragraph before it and \eqref{eq:kv8II3a}-\eqref{eq:kv8II3b}, we deduce the last term in \eqref{eq:kv8II2} is, in absolute value, $\lesssim\N^{\gamma_{\mathrm{KV}}}|\mathbb{I}(\mathfrak{t})|^{1/2}\mathfrak{p}[\mathfrak{t}(\mathrm{in}),\t;\mathcal{E}]$. We claim that the parabolic max-min principles that we used to get \eqref{eq:kv8I2} then turn \eqref{eq:kv8II2} into
\begin{align}
\mathrm{LHS}\eqref{eq:kv8II} \ &\lesssim \ \exp\{{\textstyle\int_{\mathfrak{t}(\mathrm{in})}^{\t}}\mathrm{O}(\N^{\gamma_{\mathrm{KV}}}|\mathbb{I}(\mathfrak{t})|^{\frac12})\d\s\} \nonumber\\
&\lesssim \ \exp\{\mathrm{O}(\N^{\gamma_{\mathrm{KV}}}\mathfrak{t}|\mathbb{I}(\mathfrak{t})|^{\frac12})\} \nonumber\\
&\lesssim \ \exp\{\mathrm{O}(\N^{\gamma_{\mathrm{KV}}-\frac{1}{20}}\mathfrak{t}|\mathbb{I}(\mathfrak{t})|)\} \lesssim 1. \label{eq:kv8II5}
\end{align}
Indeed, the third bound in \eqref{eq:kv8II5} follows because $|\mathbb{I}(\mathfrak{t})|\gtrsim\N^{1/10}$ by assumption. The last bound follows because $\mathfrak{t}|\mathbb{I}(\mathfrak{t})|\leq\N^{\gamma_{\mathrm{av}}}$, and because $\gamma_{\mathrm{av}},\gamma_{\mathrm{KV}}$ are small. \eqref{eq:kv8II5} implies the remaining desired estimate \eqref{eq:kv8II}, so we are done.
\end{proof}
\begin{rem}\label{remark:kv9}
 For convenience, in the following result (Lemma \ref{lemma:kv10}), let $\mathfrak{p}[\mathfrak{t}(\mathrm{in}),\t]$ denote the density for the law of $\mathbf{U}^{\t,\cdot}[\mathbb{I}(\mathfrak{t})]$ with respect to $\mathbb{P}^{\sigma,\t,\mathbb{I}(\mathfrak{t})}$. (This contrasts with Lemma \ref{lemma:kv8}.) However, this is technically not {an} abuse of notation; let us explain why. Take $\mathfrak{p}[\mathfrak{t}(\mathrm{in}),\t]$ from Lemma \ref{lemma:kv8}. The object we denote by $\mathfrak{p}[\mathfrak{t}(\mathrm{in}),\t]$ here is $\E^{\mathrm{Leb}}\mathfrak{p}[\mathfrak{t}(\mathrm{in}),\t]$, i.e., the marginal onto the variable corresponding to $\mathbf{U}^{\t,\cdot}[\mathbb{I}(\mathfrak{t})]$. For $\t=\mathfrak{t}(\mathrm{in})$, this is simply $\E^{\mathrm{Leb}}1=1$ by assumption in Lemma \ref{lemma:kv8}. This clearly equals the initial data of the density for the law of $\mathbf{U}^{\t,\cdot}[\mathbb{I}(\mathfrak{t})]$ with respect to $\mathbb{P}^{\sigma,\t,\mathbb{I}(\mathfrak{t})}$. We also know $\E^{\mathrm{Leb}}\mathfrak{p}[\mathfrak{t}(\mathrm{in}),\t]$ solves a Kolmogorov forward equation, since it is the density for the process $\t\mapsto\mathbf{U}^{\t,\cdot}[\mathbb{I}(\mathfrak{t})]$. This Kolmogorov PDE is \eqref{eq:kv8I1}, but replace $\mathscr{L}^{\mathrm{tot}}(\t,\mathbb{I}(\mathfrak{t}))$ with $\mathscr{L}(\t,\mathbb{I}(\mathfrak{t}))$. By Definition \ref{definition:le5} and the paragraph right before \eqref{eq:le7II0}, the only difference between the adjoints of these operators is some derivatives in the $\mathrm{a}$-variable we are taking $\E^{\mathrm{Leb}}$ with respect to. Thus, their action on $\E^{\mathrm{Leb}}\mathfrak{p}[\mathfrak{t}(\mathrm{in}),\t]$ vanishes, since this term integrates out said $\mathrm{a}$-variable. In particular, $\E^{\mathrm{Leb}}\mathfrak{p}[\mathfrak{t}(\mathrm{in}),\t]$ solves \eqref{eq:kv8I1} with the same data at time $\t=\mathfrak{t}(\mathrm{in})$ as $\mathfrak{p}[\mathfrak{t}(\mathrm{in}),\t]$, meaning these two are the same. (This is just a rigorous interpretation of the fact that $\mathbb{P}^{\mathrm{Leb}}$ is invariant for $\mathbf{J}(\t,\inf\mathbb{I}(\mathfrak{t});\mathbb{I}(\mathfrak{t}))$; see Lemma \ref{lemma:le7}.) A similar discussion also shows $\E^{\mathrm{Leb}}\mathfrak{p}[\mathfrak{t}(\mathrm{in}),\t;\mathcal{E}]=\mathfrak{p}[\mathfrak{t}(\mathrm{in}),\t;\mathcal{E}]$ in the setting of Lemma \ref{lemma:le8}. One consequence of this remark is that Lemma \ref{lemma:le8} holds if we replace $(\mathfrak{p}[\mathfrak{t}(\mathrm{in}),\t],\mathfrak{p}[\mathfrak{t}(\mathrm{in}),\t;\mathcal{E}])$ by their respective $\E^{\mathrm{Leb}}$-expectations. (Note this is also obvious by convexity of $\|\|_{\infty}$-norms and of the function $\eta\mapsto\eta^{-1}$ for $\eta>0$.)
\end{rem}
We now give derivative bounds. The first is a classical parabolic energy estimate obtained by differentiating an $\mathscr{L}^{2}$-norm. The second is a pointwise bound that is much more sub-optimal, but suffices for our purposes. (Like \eqref{eq:kv8I}, it follows by a Gronwall argument, albeit a more involved one.)
\begin{lemma}\label{lemma:kv10}
 Take any $\mathfrak{t}(\mathrm{in}),\mathfrak{t}\geq0$ and any discrete interval $\mathbb{I}\subseteq\mathbb{T}(\N)$. Assume that $\mathfrak{t}|\mathbb{I}(\mathfrak{t})|\lesssim\N^{\gamma_{\mathrm{av}}}$ and $|\mathbb{I}(\mathfrak{t})|\gtrsim\N^{1/10}$; see {Proposition \ref{prop:kv1}} for $\gamma_{\mathrm{av}}$. Take $\t\in\mathfrak{t}(\mathrm{in})+[0,\mathfrak{t}]$. Let $\mathfrak{p}[\mathfrak{t}(\mathrm{in}),\t]$ denote the density for the law of $\mathbf{U}^{\t,\cdot}[\mathbb{I}(\mathfrak{t})]$ with respect to $\mathbb{P}^{\sigma,\t,\mathbb{I}(\mathfrak{t})}$, determined by the initial data $\mathfrak{p}[\mathfrak{t}(\mathrm{in}),\mathfrak{t}(\mathrm{in})]\equiv1$. Recall the notation in {Definition \ref{definition:le3}}. For any $\t\in\mathfrak{t}(\mathrm{in})+[0,\mathfrak{t}]$, we have {the following in which $\mathrm{D}>0$ is any large but fixed constant:}
\begin{align}
\N^{2}{\textstyle\int_{\mathfrak{t}(\mathrm{in})}^{\t}}\mathfrak{D}_{\mathrm{FI}}^{\sigma,\s,\mathbb{I}(\mathfrak{t})}(\mathfrak{p}[\mathfrak{t}(\mathrm{in}),\s]^{2})\d\s \ \lesssim \ \N^{\gamma_{\mathrm{KV}}}\mathfrak{t}|\mathbb{I}(\mathfrak{t})|^{\frac12}+\N^{-{\mathrm{D}}}. \label{eq:kv10I}
\end{align}
Let us define $\|\|_{\infty}$ as the sup-norm for functions of $\mathbf{U}\in\R^{\mathbb{I}(\mathfrak{t})}$. We also define $|\mathbf{U}|$ as the standard Euclidean length of $\mathbf{U}\in\R^{\mathbb{I}(\mathfrak{t})}$. Given any time $\t\in\mathfrak{t}(\mathrm{in})+[0,\mathfrak{t}]$ and any point $\x$ such that $\x,\x+1\in\mathbb{I}(\mathfrak{t})$, we have
\begin{align}
\|\{1+|\mathbf{U}|^{100}\}^{-1}\mathrm{D}_{\x}\mathfrak{p}[\mathfrak{t}(\mathrm{in}),\t]\|_{\infty} \ \lesssim \ \exp[\mathrm{O}(\N^{\gamma_{\mathrm{av}}})]. \label{eq:kv10II}
\end{align}
\end{lemma}
\begin{proof}
We start with \eqref{eq:kv10I}. We first claim the following Kolmogorov forward equation:
\begin{align}
\partial_{\s}\mathfrak{p}[\mathfrak{t}(\mathrm{in}),\s] \ = \ \mathscr{L}(\s,\mathbb{I}(\mathfrak{t}))^{\ast}\mathfrak{p}[\mathfrak{t}(\mathrm{in}),\s]+\mathfrak{p}[\mathfrak{t}(\mathrm{in}),\s]\partial_{\s}\mathscr{HP}(\s,\mathbf{U};\sigma). \label{eq:kv10I0a}
\end{align}
The adjoint on the RHS is with respect to $\E^{\sigma,\s,\mathbb{I}(\mathfrak{t})}$. To derive \eqref{eq:kv10I0a}, we refer to \eqref{eq:le8II2a}-\eqref{eq:le8II2b}. (Indeed, this calculation does not use anything about the specific torus $\mathbb{T}(\N)$, and it holds if we replace it by $\mathbb{I}(\mathfrak{t})$ with periodic boundary. Similarly, the specific charge density of $0$ in this calculation plays no role either, and we can swap it with $\sigma$. Finally, we compute the log of the Lebesgue density of $\mathbb{P}^{\sigma,\s,\mathbb{I}(\mathfrak{t})}$ in \eqref{eq:le8II2b} by \eqref{eq:kv3I}.) Now, recall {the} notation from Definition \ref{definition:le5}. We claim that for some $\upsilon\gtrsim1$, we have
\begin{align}
&\partial_{\s}\E^{\sigma,\s,\mathbb{I}(\mathfrak{t})}|\mathfrak{p}[\mathfrak{t}(\mathrm{in}),\s]-1|^{2} \ = \ \partial_{\s}\E^{\sigma,\s,\mathbb{I}(\mathfrak{t})}\{\mathfrak{p}[\mathfrak{t}(\mathrm{in}),\s]^{2}-2\mathfrak{p}[\mathfrak{t}(\mathrm{in}),\s]+1\} \ = \ \partial_{\s}\E^{\sigma,\s,\mathbb{I}(\mathfrak{t})}\mathfrak{p}[\mathfrak{t}(\mathrm{in}),\s]^{2} \nonumber\\
&= \ \E^{\sigma,\s,\mathbb{I}(\mathfrak{t})}\mathfrak{p}[\mathfrak{t}(\mathrm{in}),\s]\partial_{\s}\mathfrak{p}[\mathfrak{t}(\mathrm{in}),\s]+\E^{\sigma,\s,\mathbb{I}(\mathfrak{t})}\mathfrak{p}[\mathfrak{t}(\mathrm{in}),\s]\mathscr{L}(\s,\mathbb{I}(\mathfrak{t}))\mathfrak{p}[\mathfrak{t}(\mathrm{in}),\s] \nonumber\\
&= \ \E^{\sigma,\s,\mathbb{I}(\mathfrak{t})}\mathfrak{p}[\mathfrak{t}(\mathrm{in}),\s]\mathscr{L}(\s,\mathbb{I}(\mathfrak{t}))^{\ast}\mathfrak{p}[\mathfrak{t}(\mathrm{in}),\s]+\E^{\sigma,\s,\mathbb{I}(\mathfrak{t})}\mathfrak{p}[\mathfrak{t}(\mathrm{in}),\s]^{2}\partial_{\s}\mathscr{HP}(\s,\mathbf{U};\sigma) \nonumber\\
&+ \ \E^{\sigma,\s,\mathbb{I}(\mathfrak{t})}\mathfrak{p}[\mathfrak{t}(\mathrm{in}),\s]\mathscr{L}(\s,\mathbb{I}(\mathfrak{t}))\mathfrak{p}[\mathfrak{t}(\mathrm{in}),\s]\nonumber \\
&= \ 2\E^{\sigma,\s,\mathbb{I}(\mathfrak{t})}\mathfrak{p}[\mathfrak{t}(\mathrm{in}),\s]\mathscr{L}^{\mathrm{S}}(\s,\mathbb{I}(\mathfrak{t}))\mathfrak{p}[\mathfrak{t}(\mathrm{in}),\s]+\E^{\sigma,\s,\mathbb{I}(\mathfrak{t})}\mathfrak{p}[\mathfrak{t}(\mathrm{in}),\s]^{2}\partial_{\s}\mathscr{HP}(\s,\mathbf{U};\sigma) \nonumber \\
&\leq \ -\upsilon\N^{2}\mathfrak{D}_{\mathrm{FI}}^{\sigma,\s,\mathbb{I}(\mathfrak{t})}(\mathfrak{p}[\mathfrak{t}(\mathrm{in}),\s]^{2})+\E^{\sigma,\s,\mathbb{I}(\mathfrak{t})}\mathfrak{p}[\mathfrak{t}(\mathrm{in}),\s]^{2}\partial_{\s}\mathscr{HP}(\s,\mathbf{U};\sigma). \label{eq:kv10I1}
\end{align}
The first identity is elementary. The second line is the Kolmogorov backward equation for expectations of observables of Markov processes. (Recall that the generator of $\mathbf{U}^{\s,\cdot}[\mathbb{I}(\mathfrak{t})]$ is $\mathscr{L}(\s,\mathbb{I}(\mathfrak{t}))$; see Remark \ref{remark:le6}.) The third line follows by \eqref{eq:kv10I0a}. The fourth line follows by combining the first and last terms in the line before, and by recalling the symmetric part of $\mathscr{L}(\s,\mathbb{I}(\mathfrak{t}))$ with respect to $\E^{\sigma,\s,\mathbb{I}(\mathfrak{t})}$ equals $\mathscr{L}^{\mathrm{S}}(\s,\mathbb{I}(\mathfrak{t}))$; see the paragraph before \eqref{eq:le7II0}. \eqref{eq:kv10I1} is a standard integration-by-parts. (This implies that the Dirichlet form of a Markov generator is equal to its quadratic form, up to appropriate sign; see the beginning of Section 3 of \cite{DGP}, for example. The extra $\N^{2}$ factor, which is omitted in this part of \cite{DGP}, comes from the $\N^{2}$ speed in $\mathscr{L}^{\mathrm{S}}$; see Definition \ref{definition:le5}.) At this point in the proof, we will use a forthcoming result (Lemma \ref{lemma:kv11}). Its proof does not need the current lemma, only Lemma \ref{lemma:kv8}. Thus, there is no circular reasoning. Let us estimate the last term in \eqref{eq:kv10I1}. First, {we give some extra} notation. Let $\mathfrak{p}[\mathfrak{t}(\mathrm{in}),\s;\mathcal{E}]$ be the density of the law of the stopped process $\mathbf{U}^{\s,\cdot}[\mathbb{I}(\mathfrak{t});\mathcal{E}]$ with respect to $\mathbb{P}^{\sigma,\s,\mathbb{I}(\mathfrak{t}),\mathcal{E}}$; see Definition \ref{definition:kv4}. Use $\mathrm{a}^{2}=\mathrm{b}^{2}+(\mathrm{a}+\mathrm{b})(\mathrm{a}-\mathrm{b})$ for $\mathrm{a}=\mathfrak{p}[\mathfrak{t}(\mathrm{in}),\s]$ and $\mathrm{b}=\mathfrak{p}[\mathfrak{t}(\mathrm{in}),\s;\mathcal{E}]$. This implies the following identity:
\begin{align}
\E^{\sigma,\s,\mathbb{I}(\mathfrak{t})}\mathfrak{p}[\mathfrak{t}(\mathrm{in}),\s]^{2}\partial_{\s}\mathscr{HP}(\s,\mathbf{U};\sigma) \ &= \ \E^{\sigma,\s,\mathbb{I}(\mathfrak{t})}\mathfrak{p}[\mathfrak{t}(\mathrm{in}),\s;\mathcal{E}]^{2}\partial_{\s}\mathscr{HP}(\s,\mathbf{U};\sigma)+\mathrm{Err}, \label{eq:kv10I2a}
\end{align}
where $\mathrm{Err}$ is an error that is computed and estimated as follows (with explanation given afterwards):
\begin{align}
\mathrm{Err} \ &= \ \E^{\sigma,\s,\mathbb{I}(\mathfrak{t})}\{\mathfrak{p}[\mathfrak{t}(\mathrm{in}),\s]+\mathfrak{p}[\mathfrak{t}(\mathrm{in}),\s;\mathcal{E}]\}\{\mathfrak{p}[\mathfrak{t}(\mathrm{in}),\s]-\mathfrak{p}[\mathfrak{t}(\mathrm{in}),\s;\mathcal{E}]\}\partial_{\s}\mathscr{HP}(\s,\mathbf{U};\sigma) \label{eq:kv10I2b}\\
&\lesssim \ \exp[\mathrm{O}(\N^{\gamma_{\mathrm{av}}})]\times\E^{\sigma,\s,\mathbb{I}(\mathfrak{t})}|\mathfrak{p}[\mathfrak{t}(\mathrm{in}),\s]-\mathfrak{p}[\mathfrak{t}(\mathrm{in}),\s;\mathcal{E}]|. \label{eq:kv10I2c}
\end{align}
\eqref{eq:kv10I2c} follows by \eqref{eq:kv8I}-\eqref{eq:kv8II} and because $|\partial_{\s}\mathscr{HP}|\lesssim\N^{200}$; see the paragraph after \eqref{eq:kv8II4e}. Now, for convenience, we set $p=\mathbb{P}^{\sigma,\s,\mathbb{I}(\mathfrak{t})}\{\mathcal{E}[\mathfrak{t}(\mathrm{in}),\mathfrak{t};\mathbb{I}(\mathfrak{t})]\}$. By definition, we get that $\mathfrak{p}[\mathfrak{t}(\mathrm{in}),\s;\mathcal{E}]\times p$ is a probability density with respect to $\mathbb{P}^{\sigma,\s,\mathbb{I}(\mathfrak{t})}$. (The factor $p$ is the change-of-measure when changing the reference measure $\mathbb{P}^{\sigma,\s,\mathbb{I}(\mathfrak{t}),\mathcal{E}}\mapsto\mathbb{P}^{\sigma,\s,\mathbb{I}(\mathfrak{t})}$.) Therefore, if we modify \eqref{eq:kv10I2c} by including a factor $p$ to $\mathfrak{p}[\mathfrak{t}(\mathrm{in}),\s;\mathcal{E}]$, the expectation in \eqref{eq:kv10I2c} is the total variation distance of $\mathbf{U}^{\s,\cdot}[\mathbb{I}(\mathfrak{t})]$ and $\mathbf{U}^{\s,\cdot}[\mathbb{I}(\mathfrak{t}),\mathcal{E}]$. By Lemma \ref{lemma:kv11}, this total variation is $\lesssim\exp[-\N^{\gamma_{\mathrm{KV}}/3}]$. The cost in this modification is $\exp[\mathrm{O}(\N^{\gamma_{\mathrm{av}}})]\times\E^{\sigma,\s,\mathbb{I}(\mathfrak{t})}\mathfrak{p}[\mathfrak{t}(\mathrm{in}),\s;\mathcal{E}]|1-p|$. By \eqref{eq:kv8II} and Lemma \ref{lemma:kv6}, this cost is $\lesssim\exp[\mathrm{O}(\N^{\gamma_{\mathrm{av}}})]\exp[-\N^{\gamma_{\mathrm{KV}}/2}]\lesssim\exp[-\N^{\gamma_{\mathrm{KV}}/3}]$, because $\gamma_{\mathrm{av}}\leq{\mathrm{c}}\gamma_{\mathrm{KV}}$ for some ${\mathrm{c}>0}$ small; see Proposition \ref{prop:kv1}. Using this paragraph and the previous two displays, we deduce
\begin{align}
\mathrm{LHS}\eqref{eq:kv10I2a} \ &\lesssim \ \E^{\sigma,\s,\mathbb{I}(\mathfrak{t})}\mathfrak{p}[\mathfrak{t}(\mathrm{in}),\s;\mathcal{E}]^{2}|\partial_{\s}\mathscr{HP}(\s,\mathbf{U};\sigma)|+\exp[-\N^{\frac13\gamma_{\mathrm{KV}}}] \nonumber\\
&\lesssim \ \N^{\gamma_{\mathrm{KV}}}|\mathbb{I}(\mathfrak{t})|^{\frac12}+\exp[-\N^{\frac13\gamma_{\mathrm{KV}}}], \label{eq:kv10I3}
\end{align}
where the last bound above follows because $\mathfrak{p}[\mathfrak{t}(\mathrm{in}),\s;\mathcal{E}]\lesssim1$ (see \eqref{eq:kv8II}) and {since} $\mathfrak{p}[\mathfrak{t}(\mathrm{in}),\s;\mathcal{E}]$ is supported on a set where $|\partial_{\s}\mathscr{HP}|\lesssim\N^{\gamma_{\mathrm{KV}}}|\mathbb{I}(\mathfrak{t})|^{1/2}$ (see Definition \ref{definition:kv4}). We now integrate \eqref{eq:kv10I1} over $\s\in[\mathfrak{t}(\mathrm{in}),\t]$ and use \eqref{eq:kv10I3}. This gives
\begin{align}
\mathrm{LHS}\eqref{eq:kv10I} \ &\lesssim \ \E^{\sigma,\mathfrak{t}(\mathrm{in}),\mathbb{I}(\mathfrak{t})}|\mathfrak{p}[\mathfrak{t}(\mathrm{in}),\mathfrak{t}(\mathrm{in})]-1|^{2}-\E^{\sigma,\t,\mathbb{I}(\mathfrak{t})}|\mathfrak{p}[\mathfrak{t}(\mathrm{in}),\t]-1|^{2}+{\textstyle\int_{\mathfrak{t}(\mathrm{in})}^{\t}}\mathrm{RHS}\eqref{eq:kv10I3}\d\s \nonumber\\
&\lesssim \ \mathrm{RHS}\eqref{eq:kv10I},
\end{align}
where the last bound follows because $\mathfrak{p}[\mathfrak{t}(\mathrm{in}),\mathfrak{t}(\mathrm{in})]-1=0$, because the second term on the RHS is $\leq0$, and because $\t-\mathfrak{t}(\mathrm{in})\leq\mathfrak{t}\lesssim\N^{\gamma_{\mathrm{av}}}$. This gives \eqref{eq:kv10I}. We now get \eqref{eq:kv10II}. We start with \eqref{eq:kv10I0a}. Recall $\mathscr{L}(\s,\mathbb{I}(\mathfrak{t}))^{\ast}=\mathscr{L}^{\mathrm{S}}(\s,\mathbb{I}(\mathfrak{t}))-\mathscr{L}^{\mathrm{A}}(\s,\mathbb{I}(\mathfrak{t}))$; see the paragraph before \eqref{eq:le7II0}. Now, a couple of observations. First, recall from Remark \ref{remark:le6} that $\mathscr{L}^{\mathrm{A}}(\s,\mathbb{I}(\mathfrak{t}))$ is a linear combination of $\mathrm{D}_{\x}$-operators. More precisely, as functions of $\mathbf{U}\in\R^{\mathbb{I}(\mathfrak{t})}$, \eqref{eq:kv10I0a} implies the following PDE:
\begin{align}
\partial_{\s}\mathfrak{p}[\mathfrak{t}(\mathrm{in}),\s] \ &= \ \N^{2}{\textstyle\sum_{\x}}\mathrm{D}_{\x}^{2}\mathfrak{p}[\mathfrak{t}(\mathrm{in}),\s] + \N^{2}{\textstyle\sum_{\x}}\mathrm{B}(\s,\x,\mathbf{U})\mathrm{D}_{\x}\mathfrak{p}[\mathfrak{t}(\mathrm{in}),\s] + \partial_{\s}\mathscr{HP}(\s,\mathbf{U};\sigma)\mathfrak{p}[\mathfrak{t}(\mathrm{in}),\s]. \label{eq:kv10II1}
\end{align}
All of the sums in \eqref{eq:kv10II1} are over $\x\in\mathbb{I}(\mathfrak{t})$. Also, given any $\x,\z\in\mathbb{I}(\mathfrak{t})$, we have $|\partial_{\s}\mathrm{B}(\s,\x,\mathbf{U})|\lesssim1$, and $|\partial_{\mathbf{U}(\z)}\mathrm{B}(\s,\x,\mathbf{U})|\lesssim1$. (Indeed, the $\mathrm{B}(\s,\x,\mathbf{U})$ are linear combinations of $\mathscr{U}'(\s,\mathbf{U}(\w))$ for $|\w-\x|\lesssim1$; said bounds now follow by Assumption \ref{ass:intro8}.) Let us study the second-order operator in \eqref{eq:kv10II1}. \eqref{eq:kv10II1} is a PDE on $[\mathfrak{t}(\mathrm{in}),\mathfrak{t}(\mathrm{in})+\mathfrak{t}]\times\mathbb{H}^{\sigma,\mathbb{I}(\mathfrak{t})}$. Now, let us put a Riemannian metric on $\mathbb{H}^{\sigma,\mathbb{I}(\mathfrak{t})}$. Let $\vec{\mathrm{e}}(1),\ldots,\vec{\mathrm{e}}(|\mathbb{I}(\mathfrak{t})|)$ be the standard Euclidean basis for $\R^{\mathbb{I}(\mathfrak{t})}$. Set $\vec{\mathrm{h}}(\mathrm{k})=\vec{\mathrm{e}}(\mathrm{k}+1)-\vec{\mathrm{e}}(\mathrm{k})$ given $1\leq\mathrm{k}<|\mathbb{I}(\mathfrak{t})|$. This provides a basis of $\mathbb{H}^{0,\mathbb{I}(\mathfrak{t})}$, namely a linear isomorphism $\R^{|\mathbb{I}(\mathfrak{t})|-1}\simeq\mathbb{H}^{0,\mathbb{I}(\mathfrak{t})}$. The metric we take on $\mathbb{H}^{0,\mathbb{I}(\mathfrak{t})}$ is the one induced by this isomorphism and {the} Euclidean metric on $\R^{|\mathbb{I}(\mathfrak{t})|-1}$. {An affine} shift $\mathbb{H}^{0,\mathbb{I}(\mathfrak{t})}\simeq\mathbb{H}^{\sigma,\mathbb{I}(\mathfrak{t})}$ gives a flat metric on $\mathbb{H}^{\sigma,\mathbb{I}(\mathfrak{t})}$ for any $\sigma$. Under this metric, the tangent space is spanned by mutually orthogonal $\mathrm{D}_{\x}$-operators for $\x\in\mathbb{I}(\mathfrak{t})$. (To be totally clear, flatness follows since the metric on $\mathbb{H}^{\sigma,\mathbb{I}(\mathfrak{t})}$ is determined by affine map on $\R^{|\mathbb{I}(\mathfrak{t})|-1}$. The tangent space is computed by tracking what happens to Euclidean differentials along $\vec{\mathrm{e}}(\mathrm{k})$ under this affine map.) Thus, the second-order operator in \eqref{eq:kv10II1} is a Laplacian (times $\N^{2}$) with respect to this metric. Now, for $\t\geq0$ and $\mathbf{U},\mathbf{V}\in\mathbb{H}^{\sigma,\mathbb{I}(\mathfrak{t})}$, let $\Gamma[\t,\mathbf{U},\mathbf{V}]$ solve $\Gamma[0,\mathbf{U},\mathbf{V}]=\delta[\mathbf{U}=\mathbf{V}]$ and
\begin{align}
\partial_{\t}\Gamma[\t,\mathbf{U},\mathbf{V}] \ = \ \N^{2}{\textstyle\sum_{\x}}\mathrm{D}_{\x}^{2}\Gamma[\t,\mathbf{U},\mathbf{V}], \label{eq:kv10II2}
\end{align}
where $\mathrm{D}_{\x}$ act on $\mathbf{U}$. (It turns out to not matter, as the Laplacian in $\mathrm{RHS}\eqref{eq:kv10II2}$ is self-adjoint with respect to the measure induced by the metric on $\mathbb{H}^{\sigma,\mathbb{I}(\mathfrak{t})}$ defining it.) By Duhamel, we have the following, where we now write $\mathfrak{p}$ as a function on $\mathbb{H}^{\sigma,\mathbb{I}(\mathfrak{t})}$:
\begin{align}
\mathfrak{p}[\mathfrak{t}(\mathrm{in}),\t,\mathbf{U}] \ &= \ {\textstyle\int_{\mathbb{H}^{\sigma,\mathbb{I}(\mathfrak{t})}}}\Gamma[\t-\mathfrak{t}(\mathrm{in}),\mathbf{U},\mathbf{V}]\mathfrak{p}[\mathfrak{t}(\mathrm{in}),\mathfrak{t}(\mathrm{in}),\mathbf{V}]\d^{\sigma,\mathbb{I}(\mathfrak{t})}(\mathbf{V}) \label{eq:kv10II3a}\\
&+ \ {\textstyle\int_{\mathfrak{t}(\mathrm{in})}^{\t}\int_{\mathbb{H}^{\sigma,\mathbb{I}(\mathfrak{t})}}}\Gamma[\t-\s,\mathbf{U},\mathbf{V}]\times\N^{2}{\textstyle\sum_{\x}}\mathrm{B}(\s,\x,\mathbf{V})\mathrm{D}_{\x}\mathfrak{p}[\mathfrak{t}(\mathrm{in}),\s,\mathbf{V}]\d^{\sigma,\mathbb{I}(\mathfrak{t})}(\mathbf{V})\d\s \label{eq:kv10II3b}\\
&+ \ {\textstyle\int_{\mathfrak{t}(\mathrm{in})}^{\t}\int_{\mathbb{H}^{\sigma,\mathbb{I}(\mathfrak{t})}}}\Gamma[\t-\s,\mathbf{U},\mathbf{V}]\times\partial_{\s}\mathscr{HP}(\s,\mathbf{V};\sigma)\mathfrak{p}[\mathfrak{t}(\mathrm{in}),\s,\mathbf{V}]\d^{\sigma,\mathbb{I}(\mathfrak{t})}(\mathbf{V})\d\s. \label{eq:kv10II3c} 
\end{align}
(In words, $\Gamma[\t,\mathbf{U},\mathbf{V}]$ is just the heat kernel for a standard Brownian motion on $\mathbb{H}^{\sigma,\mathbb{I}(\mathfrak{t})}$ with respect to the flat metric from before \eqref{eq:kv10II2}. In particular, it satisfies the usual gradient bounds for Brownian motions, where the usual role of derivative is played by $\mathrm{D}_{\x}$-operators. Indeed, the Laplacian, i.e. the generator of the Brownian motion, squares $\mathrm{D}_{\x}$-operators and sums over $\x$. We give precise bounds when relevant.) Note $\mathrm{RHS}\eqref{eq:kv10II3a}\equiv1$. Indeed, $\mathfrak{p}[\mathfrak{t}(\mathrm{in}),\mathfrak{t}(\mathrm{in})]\equiv1$ by assumption, and $\Gamma$ is a probability density in its forward variable. Thus, we clearly have $\mathrm{D}_{\x}\mathrm{RHS}\eqref{eq:kv10II3a}=0$ and the following Holder norm estimate for any $\upsilon\in[0,1]$:
\begin{align}
\|\{1+|\mathbf{U}|^{10}\}^{-1}\eqref{eq:kv10II3a}\|_{\mathscr{C}^{0,\upsilon}(\mathbb{H}^{\sigma,\mathbb{I}(\mathfrak{t})})} \ \lesssim \ \|\eqref{eq:kv10II3a}\|_{\mathscr{C}^{0,\upsilon}(\mathbb{H}^{\sigma,\mathbb{I}(\mathfrak{t})})} \ \lesssim \ 1. \label{eq:kv10II4aa}
\end{align}
{(To be completely clear, the Holder space $\mathscr{C}^{0,\upsilon}(\mathbb{H}^{\sigma,\mathbb{I}(\mathfrak{t})})$ is defined with respect to the metric on $\mathbb{H}^{\sigma,\mathbb{I}(\mathfrak{t})}$, constructed after \eqref{eq:kv10II1}. In particular, we know {\small$\mathrm{D}_{\x}:\mathscr{C}^{0,1}(\mathbb{H}^{\sigma,\mathbb{I}(\mathfrak{t})})\cap\mathscr{C}^{\infty}(\mathbb{H}^{\sigma,\mathbb{I}(\mathfrak{t})})\to\mathscr{C}^{0}(\mathbb{H}^{\sigma,\mathbb{I}(\mathfrak{t})})$} is uniformly bounded, where the domain is equipped with the $\mathscr{C}^{0,1}(\mathbb{H}^{\sigma,\mathbb{I}(\mathfrak{t})})$-norm. Also, the first estimate in \eqref{eq:kv10II4aa} follows because multiplying by a function in $\mathscr{C}^{0,\upsilon}(\mathbb{H}^{\sigma,\mathbb{I}(\mathfrak{t})})$ is bounded on $\mathscr{C}^{0,\upsilon}(\mathbb{H}^{\sigma,\mathbb{I}(\mathfrak{t})})$ with operator norm given by the $\mathscr{C}^{0,\upsilon}(\mathbb{H}^{\sigma,\mathbb{I}(\mathfrak{t})})$ -norm of said multiplier.) We study \eqref{eq:kv10II3c}. Letting $\mathrm{D}_{\x}$ act only on $\mathbf{U}$-variables in the following display, we claim the calculation below (with explanation given afterwards):}
\begin{align}
|\mathrm{D}_{\x}\eqref{eq:kv10II3c}| \ &\leq \ {\textstyle\int_{\mathfrak{t}(\mathrm{in})}^{\t}\int_{\mathbb{H}^{\sigma,\mathbb{I}(\mathfrak{t})}}}|\mathrm{D}_{\x}\Gamma[\t-\s,\mathbf{U},\mathbf{V}]|\times|\partial_{\s}\mathscr{HP}(\s,\mathbf{V};\sigma)|\mathfrak{p}[\mathfrak{t}(\mathrm{in}),\s,\mathbf{V}]\d^{\sigma,\mathbb{I}(\mathfrak{t})}(\mathbf{V})\d\s \label{eq:kv10II4a}\\
&\lesssim \ \N^{200}\exp[\mathrm{O}(\N^{\gamma_{\mathrm{av}}})]{\textstyle\int_{\mathfrak{t}(\mathrm{in})}^{\t}\int_{\mathbb{H}^{\sigma,\mathbb{I}(\mathfrak{t})}}}|\mathrm{D}_{\x}\Gamma[\t-\s,\mathbf{U},\mathbf{V}]|\d^{\sigma,\mathbb{I}(\mathfrak{t})}(\mathbf{V})\d\s \label{eq:kv10II4b} \\
&\lesssim \ \N^{200}\exp[\mathrm{O}(\N^{\gamma_{\mathrm{av}}})]{\textstyle\int_{\mathfrak{t}(\mathrm{in})}^{\t}}|\t-\s|^{-\frac12}\d\s \ \lesssim \ \N^{200}\exp[\mathrm{O}(\N^{\gamma_{\mathrm{av}}})]\mathfrak{t}^{\frac12} \ \lesssim \ \exp[\mathrm{O}(\N^{\gamma_{\mathrm{av}}})]. \label{eq:kv10II4c}
\end{align}
\eqref{eq:kv10II4a} follows by triangle inequality. \eqref{eq:kv10II4b} follows by \eqref{eq:kv8I} and $|\partial_{\s}\mathscr{HP}|\lesssim\N^{200}$; see right after \eqref{eq:kv10I2c}. The first estimate in \eqref{eq:kv10II4c} follows from standard gradient estimates for Brownian motion heat kernels. (Indeed, one can check this is true for the Gaussian heat kernel on Euclidean {spaces} of any finite dimension. In particular, there is no dimensional prefactor. To see this for the Euclidean case, note the Gaussian heat kernel factorizes into one-dimensional heat kernels. Taking partial in a fixed standard basis direction affects one of the one-dimensional factors, so the first implied constant in \eqref{eq:kv10II4c} is the one in the one-dimensional case. See the beginning of \cite{W}, which holds in a general geometric setting.)  The second bound in \eqref{eq:kv10II4c} follows by integration and $\t\in\mathfrak{t}(\mathrm{in})+[0,\mathfrak{t}]$. The final bound in \eqref{eq:kv10II4c} uses $\mathfrak{t}\leq\mathfrak{t}|\mathbb{I}(\mathfrak{t})|\lesssim\N^{\gamma_{\mathrm{av}}}$ (since $|\mathbb{I}(\mathfrak{t})|$ is always a positive integer). It is otherwise elementary. Since integration against $\Gamma[\t-\s,\mathbf{U},\cdot]$ is a contractive operator $\mathscr{L}^{\infty}(\mathbb{H}^{\sigma,\mathbb{I}(\mathfrak{t})})\to\mathscr{L}^{\infty}(\mathbb{H}^{\sigma,\mathbb{I}(\mathfrak{t})})$ if $\s\leq\t$, we also know $|\eqref{eq:kv10II3c}|\lesssim\exp[\mathrm{O}(\N^{\gamma_{\mathrm{av}}})]$ by basically the same calculation as \eqref{eq:kv10II4a}-\eqref{eq:kv10II4c}. (Just forget the $|\t-\s|^{-1/2}$-factor therein.) Via interpolation of Holder norms, we then deduce the following for any $\upsilon\in[0,1]$:
\begin{align}
\|\{1+|\mathbf{U}|^{10}\}^{-1}\eqref{eq:kv10II3c}\|_{\mathscr{C}^{0,\upsilon}(\mathbb{H}^{\sigma,\mathbb{I}(\mathfrak{t})})} \ \lesssim \ \|\eqref{eq:kv10II3c}\|_{\mathscr{C}^{0,\upsilon}(\mathbb{H}^{\sigma,\mathbb{I}(\mathfrak{t})})} \ \lesssim \ \exp[\mathrm{O}(\N^{\gamma_{\mathrm{av}}})]. \label{eq:kv10II4d}
\end{align}
We now study \eqref{eq:kv10II3b}. Upon relabeling the sum-variable in \eqref{eq:kv10II3b} from $\x$ to $\y$, we have the following with explanation after:
\begin{align}
\eqref{eq:kv10II3b} \ = \ &-\N^{2}{\textstyle\sum_{\y}}{\textstyle\int_{\mathfrak{t}(\mathrm{in})}^{\t}\int_{\mathbb{H}^{\sigma,\mathbb{I}(\mathfrak{t})}}}\mathrm{D}_{\y}\Gamma[\t-\s,\mathbf{U},\mathbf{V}]\times\mathrm{B}(\s,\y,\mathbf{V})\mathfrak{p}[\mathfrak{t}(\mathrm{in}),\s,\mathbf{V}]\d^{\sigma,\mathbb{I}(\mathfrak{t})}(\mathbf{V})\d\s \label{eq:kv10II5aa}\\[-1.2mm]
&-\N^{2}{\textstyle\sum_{\y}}{\textstyle\int_{\mathfrak{t}(\mathrm{in})}^{\t}\int_{\mathbb{H}^{\sigma,\mathbb{I}(\mathfrak{t})}}}\Gamma[\t-\s,\mathbf{U},\mathbf{V}]\times\mathrm{D}_{\y}\mathrm{B}(\s,\y,\mathbf{V})\mathfrak{p}[\mathfrak{t}(\mathrm{in}),\s,\mathbf{V}]\d^{\sigma,\mathbb{I}(\mathfrak{t})}(\mathbf{V})\d\s. \label{eq:kv10II5ab}
\end{align}
\eqref{eq:kv10II5aa}-\eqref{eq:kv10II5ab} follows from integration-by-parts for $\mathrm{D}_{\x}$ in $\mathbb{H}^{\sigma,\mathbb{I}(\mathfrak{t})}$. (Indeed, $\mathbb{H}^{\sigma,\mathbb{I}(\mathfrak{t})}$ is a flat hyperplane whose tangent space is spanned by mutually orthogonal $\mathrm{D}_{\x}$-operators.) We then use the Leibniz rule. Let us now treat the second line \eqref{eq:kv10II5ab}. We claim:
\begin{align}
&|\mathrm{D}_{\x}\eqref{eq:kv10II5ab}| \nonumber\\
&\leq \ \N^{2}|\mathbb{I}(\mathfrak{t})|{\textstyle\sup_{\y}}{\textstyle\int_{\mathfrak{t}(\mathrm{in})}^{\t}\int_{\mathbb{H}^{\sigma,\mathbb{I}(\mathfrak{t})}}}|\mathrm{D}_{\x}\Gamma[\t-\s,\mathbf{U},\mathbf{V}]|\times|\mathrm{D}_{\y}\mathrm{B}(\s,\y,\mathbf{V})|\mathfrak{p}[\mathfrak{t}(\mathrm{in}),\s,\mathbf{V}]\d^{\sigma,\mathbb{I}(\mathfrak{t})}(\mathbf{V})\d\s \label{eq:kv10II5ba}\\
&\lesssim \ \N^{2}|\mathbb{I}(\mathfrak{t})|\exp[\mathrm{O}(\N^{\gamma_{\mathrm{av}}})]{\textstyle\int_{\mathfrak{t}(\mathrm{in})}^{\t}\int_{\mathbb{H}^{\sigma,\mathbb{I}(\mathfrak{t})}}}|\mathrm{D}_{\x}\Gamma[\t-\s,\mathbf{U},\mathbf{V}]|\d^{\sigma,\mathbb{I}(\mathfrak{t})}(\mathbf{V})\d\s \label{eq:kv10II5bb}\\
&\lesssim \ \N^{2}|\mathbb{I}(\mathfrak{t})|\exp[\mathrm{O}(\N^{\gamma_{\mathrm{av}}})]{\textstyle\int_{\mathfrak{t}(\mathrm{in})}^{\t}}|\t-\s|^{-\frac12}\d\s \ \lesssim \ \exp[\mathrm{O}(\N^{\gamma_{\mathrm{av}}})]. \label{eq:kv10II5bc}
\end{align}
\eqref{eq:kv10II5ba} follows by triangle inequality. \eqref{eq:kv10II5bb} follows by \eqref{eq:kv8I} and because $\mathrm{B}$ is uniformly Lipschitz; see after \eqref{eq:kv10II1}. Also, $|\eqref{eq:kv10II5ab}|\lesssim\exp[\mathrm{O}(\N^{\gamma_{\mathrm{av}}})]$ by the same calculation. (Forget $\mathrm{D}_{\x}$ and $|\t-\s|^{-1/2}$.) Again, for $\upsilon\in[0,1]$, interpolation gives
\begin{align}
\|\{1+|\mathbf{U}|^{10}\}^{-1}\eqref{eq:kv10II5ab}\|_{\mathscr{C}^{0,\upsilon}(\mathbb{H}^{\sigma,\mathbb{I}(\mathfrak{t})})} \ \lesssim \ \|\eqref{eq:kv10II5ab}\|_{\mathscr{C}^{0,\upsilon}(\mathbb{H}^{\sigma,\mathbb{I}(\mathfrak{t})})} \ \lesssim \ \exp[\mathrm{O}(\N^{\gamma_{\mathrm{av}}})]. \label{eq:kv10II5ca}
\end{align}
We treat $\mathrm{RHS}\eqref{eq:kv10II5aa}$. Set $\|\|=\|\|_{\mathscr{C}^{0,\upsilon}(\mathbb{H}^{\sigma,\mathbb{I}(\mathfrak{t})})}$. For $\upsilon\in(0,1)$ and $\kappa(\upsilon)\gtrsim_{\upsilon}1$, we claim that the quantity $\|\{1+|\mathbf{U}|^{10}\}^{-1}\mathrm{RHS}\eqref{eq:kv10II5aa}\|$ is
{\small
\begin{align}
&\lesssim \ \N^{2}|\mathbb{I}(\mathfrak{t})|{\textstyle\sup_{\y}}{\textstyle\int_{\mathfrak{t}(\mathrm{in})}^{\t}}\|\{1+|\mathbf{U}|^{10}\}^{-1}{\textstyle\int_{\mathbb{H}^{\sigma,\mathbb{I}(\mathfrak{t})}}}\mathrm{D}_{\y}\Gamma[\t-\s,\mathbf{U},\mathbf{V}]\times\mathrm{B}(\s,\y,\mathbf{V})\mathfrak{p}[\mathfrak{t}(\mathrm{in}),\s,\mathbf{V}]\d^{\sigma,\mathbb{I}(\mathfrak{t})}(\mathbf{V})\|\d\s \label{eq:kv10II5daa}\\
&\lesssim \ \N^{2}|\mathbb{I}(\mathfrak{t})|{\textstyle\sup_{\y}}{\textstyle\int_{\mathfrak{t}(\mathrm{in})}^{\t}}\|{\textstyle\int_{\mathbb{H}^{\sigma,\mathbb{I}(\mathfrak{t})}}}\mathrm{D}_{\y}\Gamma[\t-\s,\mathbf{U},\mathbf{V}]\times\{\mathrm{B}(\s,\y,\mathbf{V})-\mathrm{B}(\s,\y,\mathbf{U})\}\mathfrak{p}[\mathfrak{t}(\mathrm{in}),\s,\mathbf{V}]\d^{\sigma,\mathbb{I}(\mathfrak{t})}(\mathbf{V})\|\d\s \label{eq:kv10II5da}\\
&+ \ \N^{2}|\mathbb{I}(\mathfrak{t})|{\textstyle\sup_{\y}}{\textstyle\int_{\mathfrak{t}(\mathrm{in})}^{\t}}\|\{1+|\mathbf{U}|^{10}\}^{-1}{\textstyle\int_{\mathbb{H}^{\sigma,\mathbb{I}(\mathfrak{t})}}}\mathrm{D}_{\y}\Gamma[\t-\s,\mathbf{U},\mathbf{V}]\times\mathrm{B}(\s,\y,\mathbf{U})\mathfrak{p}[\mathfrak{t}(\mathrm{in}),\s,\mathbf{V}]\d^{\sigma,\mathbb{I}(\mathfrak{t})}(\mathbf{V})\|\d\s\label{eq:kv10II5db}\\
&\lesssim \ \N^{2}|\mathbb{I}(\mathfrak{t})|{\textstyle\sup_{\y}}{\textstyle\int_{\mathfrak{t}(\mathrm{in})}^{\t}}\|{\textstyle\int_{\mathbb{H}^{\sigma,\mathbb{I}(\mathfrak{t})}}}\mathrm{D}_{\y}\Gamma[\t-\s,\mathbf{U},\mathbf{V}]\times\{\mathrm{B}(\s,\y,\mathbf{V})-\mathrm{B}(\s,\y,\mathbf{U})\}\mathfrak{p}[\mathfrak{t}(\mathrm{in}),\s,\mathbf{V}]\d^{\sigma,\mathbb{I}(\mathfrak{t})}(\mathbf{V})\|\d\s \label{eq:kv10II5dc}\\
&+ \ \N^{2}|\mathbb{I}(\mathfrak{t})|{\textstyle\sup_{\y}}{\textstyle\int_{\mathfrak{t}(\mathrm{in})}^{\t}}\|{\textstyle\int_{\mathbb{H}^{\sigma,\mathbb{I}(\mathfrak{t})}}}\mathrm{D}_{\y}\Gamma[\t-\s,\mathbf{U},\mathbf{V}]\times\mathfrak{p}[\mathfrak{t}(\mathrm{in}),\s,\mathbf{V}]\d^{\sigma,\mathbb{I}(\mathfrak{t})}(\mathbf{V})\|\d\s \label{eq:kv10II5dd}\\
&\lesssim \ \N^{2}|\mathbb{I}(\mathfrak{t})|\exp[\mathrm{O}(\N^{\gamma_{\mathrm{av}}})]{\textstyle\int_{\mathfrak{t}(\mathrm{in})}^{\t}}\{1+|\t-\s|^{-1+\kappa(\upsilon)}\}\d\s \ \lesssim_{\upsilon} \ \N^{2}|\mathbb{I}(\mathfrak{t})|\exp[\mathrm{O}(\N^{\gamma_{\mathrm{av}}})] \ \lesssim \ \exp[\mathrm{O}(\N^{\gamma_{\mathrm{av}}})].\label{eq:kv10II5de}
\end{align}
}The first bound is by triangle inequality. \eqref{eq:kv10II5da}-\eqref{eq:kv10II5db} follows by writing $\mathrm{B}(\s,\y,\mathbf{V})=\mathrm{B}(\s,\y,\mathbf{V})-\mathrm{B}(\s,\y,\mathbf{U})+\mathrm{B}(\s,\y,\mathbf{U})$, the triangle inequality and the fact that $\|\{1+|\mathbf{U}|^{10}\}^{-1}\|\lesssim1$. (This bound is the reason why \eqref{eq:kv10II5da} does not have $\{1+|\mathbf{U}|^{10}\}^{-1}$ inside the norm. Indeed, since $\|\|$ is sub-multiplicative, so that Holder spaces are Banach algebras, the cost we pay in removing $\{1+|\mathbf{U}|^{10}\}^{-1}$ from the norm in \eqref{eq:kv10II5da} is $\|\{1+|\mathbf{U}|^{10}\}^{-1}\|\lesssim1$.) \eqref{eq:kv10II5dc}-\eqref{eq:kv10II5dd} follows by first leaving \eqref{eq:kv10II5da} alone. Then, we pull out $\mathrm{B}(\s,\y,\mathbf{U})$ from the $\mathbb{H}^{\sigma,\mathbb{I}(\mathfrak{t})}$-integral in \eqref{eq:kv10II5db}. Afterwards, we remove the resulting factor $\{1+|\mathbf{U}|^{10}\}^{-1}\mathrm{B}(\s,\y,\mathbf{U})$ outside said integral by noting $\|\{1+|\mathbf{U}|^{10}\}^{-1}\mathrm{B}(\s,\y,\mathbf{U})\|\lesssim1$ and the reasoning for \eqref{eq:kv10II5da}-\eqref{eq:kv10II5db}. This $\|\|$-estimate can be checked by the uniform Lipschitz property of $\mathrm{B}(\s,\y,\mathbf{U})$ and elementary calculations; see after \eqref{eq:kv10II1}. Let us now explain how to get \eqref{eq:kv10II5de}. First, the norm in \eqref{eq:kv10II5dd} is $\lesssim\|\mathfrak{p}[\mathfrak{t}(\mathrm{in})\|_{\infty}\lesssim\exp[\mathrm{O}(\N^{\gamma_{\mathrm{av}}})]$ (see \eqref{eq:kv8I}) times $|\t-\s|^{-1+\kappa(\upsilon)}$, where $\kappa(\upsilon)\to0$ as $\upsilon\to1$ from below. This is just the usual operator norm estimate for the $\Gamma$-semigroup $\mathscr{C}^{1,\upsilon}(\mathbb{H}^{\sigma,\mathbb{I}(\mathfrak{t})})\to\mathscr{C}^{0}(\mathbb{H}^{\sigma,\mathbb{I}(\mathfrak{t})})$. It can be shown by interpolating the gradient bounds in \cite{W} (which hold for gradients of any order) with $\mathscr{C}^{0}$-contractivity of $\Gamma$. (Compare with the Euclidean case, similar to the paragraph after \eqref{eq:kv10II4c}. As a reality check, if $\upsilon=1$, then $\kappa(\upsilon)=0$ as we claimed. The resulting bound is just a second-derivative estimate for a Gaussian heat kernel.) Ultimately, we know that $\eqref{eq:kv10II5dd}\lesssim$ the first term in \eqref{eq:kv10II5de}. To estimate \eqref{eq:kv10II5dc}, first note that $|\mathrm{B}(\s,\y,\mathbf{V})-\mathrm{B}(\s,\y,\mathbf{U})|\lesssim|\mathbf{V}-\mathbf{U}|$; see after \eqref{eq:kv10II1}. This factor only dampens the on-diagonal singularity of $\mathrm{D}_{\y}\Gamma[\t-\s,\mathbf{U},\mathbf{V}]$ at $\mathbf{V}=\mathbf{U}$. Thus, this term has a better heat operator estimate than \eqref{eq:kv10II5dd}. This shows the first bound in \eqref{eq:kv10II5de}. The rest of \eqref{eq:kv10II5de} follows by doing the time-integral and proceeding as in \eqref{eq:kv10II5bc}. Let us now combine \eqref{eq:kv10II3a}-\eqref{eq:kv10II3c}, \eqref{eq:kv10II4aa}, \eqref{eq:kv10II4d}, \eqref{eq:kv10II5aa}-\eqref{eq:kv10II5ab}, \eqref{eq:kv10II5ca}, and \eqref{eq:kv10II5daa}-\eqref{eq:kv10II5de}. This shows (for any fixed $\upsilon\in(0,1)$)
\begin{align}
\|\{1+|\mathbf{U}|^{10}\}^{-1}\mathfrak{p}[\mathfrak{t}(\mathrm{in}),\t]\|_{\mathscr{C}^{0,\upsilon}(\mathbb{H}^{\sigma,\mathbb{I}(\mathfrak{t})})} \ \lesssim_{\upsilon} \ \exp[\mathrm{O}(\N^{\gamma_{\mathrm{av}}})]. \label{eq:kv10II5ea}
\end{align}
We would like to now upgrade $\upsilon\mapsto1$. Note that \eqref{eq:kv10II3a} and \eqref{eq:kv10II3c} have $\upsilon=1$ estimates; see \eqref{eq:kv10II4a}-\eqref{eq:kv10II4c} and right before \eqref{eq:kv10II4aa}. As for \eqref{eq:kv10II3b}, we also have $\upsilon=1$ estimates for \eqref{eq:kv10II5ab}; see \eqref{eq:kv10II5ba}-\eqref{eq:kv10II5bc}. Thus, we only have to upgrade $\upsilon\mapsto1$ for $\mathrm{RHS}\eqref{eq:kv10II5aa}$. To this end, we first decompose $\mathrm{RHS}\eqref{eq:kv10II5aa}$ via the following display (which is elementary to check):
{\small
\begin{align}
&\mathrm{RHS}\eqref{eq:kv10II5aa} \nonumber\\
&= \ -\N^{2}{\textstyle\sum_{\y}}{\textstyle\int_{\mathfrak{t}(\mathrm{in})}^{\t}\mathrm{B}(\s,\y,\mathbf{U})\{1+|\mathbf{U}|^{10}\}\int_{\mathbb{H}^{\sigma,\mathbb{I}(\mathfrak{t})}}}\mathrm{D}_{\y}\Gamma[\t-\s,\mathbf{U},\mathbf{V}]\times\tfrac{\mathfrak{p}[\mathfrak{t}(\mathrm{in}),\s,\mathbf{V}]}{1+|\mathbf{V}|^{10}}\d^{\sigma,\mathbb{I}(\mathfrak{t})}(\mathbf{V})\d\s \label{eq:kv10II6a}\\
&-\N^{2}{\textstyle\sum_{\y}}{\textstyle\int_{\mathfrak{t}(\mathrm{in})}^{\t}\int_{\mathbb{H}^{\sigma,\mathbb{I}(\mathfrak{t})}}}\mathrm{D}_{\y}\Gamma[\t-\s,\mathbf{U},\mathbf{V}][\mathrm{B}(\s,\y,\mathbf{V})\{1+|\mathbf{V}|^{10}\}-\mathrm{B}(\s,\y,\mathbf{U})\{1+|\mathbf{U}|^{10}\}] \label{eq:kv10II6b}\\
&\quad\quad\quad\quad\quad\quad\quad\quad\quad\quad\quad\quad\quad\quad\quad\quad\quad\quad\quad\quad\quad\quad\quad\quad\quad\quad\quad\quad\quad\quad\times\tfrac{\mathfrak{p}[\mathfrak{t}(\mathrm{in}),\s,\mathbf{V}]}{1+|\mathbf{V}|^{10}}\d^{\sigma,\mathbb{I}(\mathfrak{t})}(\mathbf{V})\d\s.\nonumber
\end{align}
}We will now control $\eqref{eq:kv10II6a}$. First, let $\llangle\rrangle$ be the $\mathscr{C}^{0,\varrho}(\mathbb{H}^{\sigma,\mathbb{I}(\mathfrak{t})})$-norm for $\varrho<1$. (We eventually take $\varrho\to1$. We cannot take $\varrho=1$ as certain Gaussian heat semigroups needed in the following paragraph only work on Holder spaces for $\varrho\neq1$.) We claim
\begin{align}
\llangle\{1+|\mathbf{U}|^{100}\}^{-1}\eqref{eq:kv10II6a}\rrangle &\lesssim \N^{2}|\mathbb{I}(\mathfrak{t})|{\textstyle\sup_{\y}}{\textstyle\int_{\mathfrak{t}(\mathrm{in})}^{\t}\llangle\int_{\mathbb{H}^{\sigma,\mathbb{I}(\mathfrak{t})}}}\mathrm{D}_{\y}\Gamma[\t-\s,\mathbf{U},\mathbf{V}]\tfrac{\mathfrak{p}[\mathfrak{t}(\mathrm{in}),\s,\mathbf{V}]}{1+|\mathbf{V}|^{10}}\d^{\sigma,\mathbb{I}(\mathfrak{t})}(\mathbf{V})\rrangle\d\s \label{eq:kv10II7a} \\
&\lesssim \N^{2}|\mathbb{I}(\mathfrak{t})|\exp[\mathrm{O}(\N^{\gamma_{\mathrm{av}}})]{\textstyle\int_{\mathfrak{t}(\mathrm{in})}^{\t}}|\t-\s|^{-1+\kappa}\d\s \ \lesssim \ \exp[\mathrm{O}(\N^{\gamma_{\mathrm{av}}})]. \label{eq:kv10II7b}
\end{align}
Here, $\kappa\gtrsim1$. We clarify \eqref{eq:kv10II7a}. Note $\{1+|\mathbf{U}|^{100}\}^{-1}\eqref{eq:kv10II6a}$ is just $\eqref{eq:kv10II6a}$ but replacing $\mathrm{B}(\s,\y,\mathbf{U})\{1+|\mathbf{U}|^{10}\}\mapsto\mathrm{B}(\s,\y,\mathbf{U})\{1+|\mathbf{U}|^{10}\}\{1+|\mathbf{U}|^{100}\}^{-1}$. This is uniformly bounded since $\mathrm{B}(\s,\y,\mathbf{U})$ is uniformly Lipschitz in $\mathbf{U}$; see after \eqref{eq:kv10II1}. For the same reason, it is uniformly Lipschitz in $\mathbf{U}$. Thus, by interpolation, its $\llangle\rrangle$-norm is $\lesssim1$. Because $\llangle\rrangle$ is sub-multiplicative, for the sake of upper bound we can forget the resulting factor $\mathrm{B}(\s,\y,\mathbf{U})\{1+|\mathbf{U}|^{10}\}\{1+|\mathbf{U}|^{100}\}^{-1}$. This shows \eqref{eq:kv10II7a}. Let us now explain \eqref{eq:kv10II7b}. It suffices to show the $\llangle\rrangle$ in $\mathrm{RHS}\eqref{eq:kv10II7a}$ is $\lesssim|\t-\s|^{-1+\kappa}$ for some $\kappa\gtrsim1$. This $\llangle\rrangle$ is $\lesssim$ the $\mathscr{C}^{1,\varrho}(\mathbb{H}^{\sigma,\mathbb{I}(\mathfrak{t})})$-norm of the $\Gamma$-operator acting on $\{1+|\mathbf{V}|^{10}\}^{-1}\mathfrak{p}[\mathfrak{t}(\mathrm{in}),\s,\mathbf{V}]$. If we only knew $\mathscr{C}^{0}$-estimates for $\{1+|\mathbf{V}|^{10}\}^{-1}\mathfrak{p}[\mathfrak{t}(\mathrm{in}),\s,\mathbf{V}]$, standard Gaussian heat kernel estimates show this norm is $\lesssim|\t-\s|^{-[1+\varrho]/2}$ times the sup-norm of $\{1+|\mathbf{V}|^{10}\}^{-1}\mathfrak{p}[\mathfrak{t}(\mathrm{in}),\s,\mathbf{V}]$. However, by \eqref{eq:kv10II5ea}, we know the $\mathscr{C}^{0,\upsilon}(\mathbb{H}^{\sigma,\mathbb{I}(\mathfrak{t})})$-norm (for $\upsilon=2/3$) of $\{1+|\mathbf{V}|^{10}\}^{-1}\mathfrak{p}[\mathfrak{t}(\mathrm{in}),\s,\mathbf{V}]$ is $\lesssim\exp[\mathrm{O}(\N^{\gamma_{\mathrm{av}}})]$ as well. Thus, the $\llangle\rrangle$ in $\mathrm{RHS}\eqref{eq:kv10II7a}$ gets a regularization in its short-time singularity. Namely, it is $\lesssim|\t-\s|^{-[1+\varrho]/2+\upsilon/2}\exp[\mathrm{O}(\N^{\gamma_{\mathrm{av}}})]$. It now suffices to take $\kappa=1-[1+\varrho]/2+\upsilon/2$, which is $\gtrsim1$ if $\varrho\leq1$ and $\upsilon=2/3$. (For these $\Gamma$-semigroup bounds, compare to the standard Gaussian case on Euclidean {spaces}. We can also interpolate gradient bounds in \cite{W} with $\mathscr{C}^{0}(\mathbb{H}^{\sigma,\mathbb{I}(\mathfrak{t})})$-contractivity as we noted after \eqref{eq:kv10II5de}.) The second bound in \eqref{eq:kv10II7b} follows by doing the integral and proceeding as in \eqref{eq:kv10II4c}. Now, we study \eqref{eq:kv10II6b}. Recall $\mathrm{B}(\s,\y,\mathbf{V})$ is uniformly Lipschitz in $\mathbf{V}$; see after \eqref{eq:kv10II1}. By triangle inequality and elementary manipulations,
\begin{align}
&|\mathrm{B}(\s,\y,\mathbf{V})\{1+|\mathbf{V}|^{10}\}-\mathrm{B}(\s,\y,\mathbf{U})\{1+|\mathbf{U}|^{10}\}| \nonumber\\
&\lesssim \ |\mathrm{B}(\s,\y,\mathbf{V})-\mathrm{B}(\s,\y,\mathbf{U})|\{1+|\mathbf{U}|^{10}\}+|\mathrm{B}(\s,\y,\mathbf{V})|\cdot||\mathbf{V}|^{10}-|\mathbf{U}|^{10}| \nonumber\\
&\lesssim \ |\mathbf{V}-\mathbf{U}|\{1+|\mathbf{U}|^{10}\}+|\mathrm{B}(\s,\y,\mathbf{V})-\mathrm{B}(\s,\y,\mathbf{U})|\cdot||\mathbf{V}|^{10}-|\mathbf{U}|^{10}|\label{eq:kv10II8a}\\
&+ \ |\mathrm{B}(\s,\y,\mathbf{U})|\cdot||\mathbf{V}|^{10}-|\mathbf{U}|^{10}| \nonumber\\
&\lesssim \ |\mathbf{V}-\mathbf{U}|\{1+|\mathbf{U}|^{10}\}+|\mathbf{V}-\mathbf{U}|\cdot||\mathbf{V}|^{10}-|\mathbf{U}|^{10}|+|\mathbf{U}|\cdot||\mathbf{V}|^{10}-|\mathbf{U}|^{10}| .\label{eq:kv10II8b}
\end{align}
{The application of a Taylor expansion} for the \emph{analytic} polynomials of $|\mathbf{V}|,|\mathbf{U}|$ shows $\eqref{eq:kv10II8b}\lesssim\{1+|\mathbf{U}|^{20}\}|\mathbf{V}-\mathbf{U}|$. Now, take any $\x\in\mathbb{I}(\mathfrak{t})$. We claim that the following estimate holds (with explanation given after):
\begin{align}
\tfrac{|\mathrm{D}_{\x}\eqref{eq:kv10II6b}|}{1+|\mathbf{U}|^{100}} \ &\lesssim \ \N^{2}|\mathbb{I}(\mathfrak{t})|\exp[\mathrm{O}(\N^{\gamma_{\mathrm{av}}})]{\textstyle\sup_{\y}}{\textstyle\int_{\mathfrak{t}(\mathrm{in})}^{\t}\int_{\mathbb{H}^{\sigma,\mathbb{I}(\mathfrak{t})}}}|\mathrm{D}_{\x}\mathrm{D}_{\y}\Gamma[\t-\s,\mathbf{U},\mathbf{V}]||\mathbf{V}-\mathbf{U}|\d^{\sigma,\mathbb{I}(\mathfrak{t})}(\mathbf{V})\d\s \label{eq:kv10II9a}\\
&\lesssim \ \N^{2}|\mathbb{I}(\mathfrak{t})|\exp[\mathrm{O}(\N^{\gamma_{\mathrm{av}}})]{\textstyle\sup_{\y}}{\textstyle\int_{\mathfrak{t}(\mathrm{in})}^{\t}}|\t-\s|^{-\frac12}\d\s \ \lesssim \ \exp[\mathrm{O}(\N^{\gamma_{\mathrm{av}}})]. \label{eq:kv10II9b}
\end{align}
\eqref{eq:kv10II9a} holds for the following reason. First, the $\mathrm{D}_{\x}$-operator hits the $\Gamma$-kernel. Next, we apply \eqref{eq:kv10II8b} and the sentence after it to bound $|\mathrm{B}(\s,\y,\mathbf{V})\{1+|\mathbf{V}|^{10}\}-\mathrm{B}(\s,\y,\mathbf{U})\{1+|\mathbf{U}|^{10}\}|\lesssim\{1+|\mathbf{U}|^{20}\}|\mathbf{V}-\mathbf{U}|$ in \eqref{eq:kv10II6b}. $\{1+|\mathbf{U}|^{20}\}$ can be removed via the $\{1+|\mathbf{U}|^{100}\}^{-1}$-factor in $\mathrm{LHS}\eqref{eq:kv10II9a}$. Finally, we use \eqref{eq:kv8I} to bound the $\mathfrak{p}[\mathfrak{t}(\mathrm{in}),\s]$-term in \eqref{eq:kv10II6b} by $\exp[\mathrm{O}(\N^{\gamma_{\mathrm{av}}})]$. This gives \eqref{eq:kv10II9a}. Let us explain the first bound in \eqref{eq:kv10II9b}. If we did not have $|\mathbf{V}-\mathbf{U}|$ in $\mathrm{RHS}\eqref{eq:kv10II9a}$, then the first bound in \eqref{eq:kv10II9b} would hold if we replace $|\t-\s|^{-1/2}$ by $|\t-\s|^{-1}$. This is the typical second-derivative bound for standard Gaussian heat kernels. Now, the $|\mathbf{V}-\mathbf{U}|$-factor in $\mathrm{RHS}\eqref{eq:kv10II9a}$ only helps, because it partially regularizes the short-time singularity of the heat kernel near $\mathbf{V}\approx\mathbf{U}$. (Thus, by parabolic scaling, $|\mathbf{V}-\mathbf{U}|$ gives square-root savings $|\t-\s|^{-1}\mapsto|\t-\s|^{-1/2}$. This is exactly the reasoning used to bound \eqref{eq:kv10II5dc}.) The rest of \eqref{eq:kv10II9b} follows from {the} reasoning for \eqref{eq:kv10II5bc}. A calculation like \eqref{eq:kv10II9a}-\eqref{eq:kv10II9b} also shows $\{1+|\mathbf{U}|^{100}\}^{-1}|\eqref{eq:kv10II6b}|\lesssim\exp[\mathrm{O}(\N^{\gamma_{\mathrm{av}}})]$. (Forget $\mathrm{D}_{\x}$ and $|\t-\s|^{-1/2}$ therein.) So, by interpolation, we get
\begin{align}
\llangle\{1+|\mathbf{U}|^{100}\}^{-1}\eqref{eq:kv10II6b}\rrangle \ \lesssim \ \exp[\mathrm{O}(\N^{\gamma_{\mathrm{av}}})]. \label{eq:kv10II10}
\end{align}
(Recall $\llangle\rrangle$ from before \eqref{eq:kv10II7b}.) Combine \eqref{eq:kv10II3a}-\eqref{eq:kv10II3c}, \eqref{eq:kv10II4aa}, \eqref{eq:kv10II4d}, \eqref{eq:kv10II5aa}-\eqref{eq:kv10II5ab}, \eqref{eq:kv10II5ca}, \eqref{eq:kv10II6a}-\eqref{eq:kv10II6b},  \eqref{eq:kv10II7a}-\eqref{eq:kv10II7b}, and \eqref{eq:kv10II10}. This ultimately implies the following Holder estimate uniformly in $\varrho<1$:
\begin{align}
\|\{1+|\mathbf{U}|^{100}\}^{-1}\mathfrak{p}[\mathfrak{t}(\mathrm{in}),\t]\|_{\mathscr{C}^{0,\varrho}(\mathbb{H}^{\sigma,\mathbb{I}(\mathfrak{t})})} \ \lesssim \ \exp[\mathrm{O}(\N^{\gamma_{\mathrm{av}}})]. \label{eq:kv10II8}
\end{align}
(Some displays listed before \eqref{eq:kv10II8} have the weight $\{1+|\mathbf{U}|^{10}\}^{-1}$. As we noted after \eqref{eq:kv10II4aa}, making the weight smaller makes $\mathscr{C}^{0,\varrho}(\mathbb{H}^{\sigma,\mathbb{I}(\mathfrak{t})})$-norms smaller, so \eqref{eq:kv10II8} follows.) Because \eqref{eq:kv10II8} is uniform in $\varrho<1$ and the RHS is independent of $\varrho$, \eqref{eq:kv10II8} holds if we replace $\varrho\mapsto1$ on the LHS. This implies \eqref{eq:kv10II}; see after \eqref{eq:kv10II4aa}. So, we are done.
\end{proof}
\subsection{Comparing \eqref{eq:hf}-\eqref{eq:glsde} with stopped versions from Definition \ref{definition:kv4}}
We continue Lemma \ref{lemma:kv6} and the paragraph before it. In words, stopping \eqref{eq:glsdeloc}-\eqref{eq:hfloc} at $\tau[\mathcal{E}]$ in Definition \ref{definition:kv4} does nothing with extremely high probability.
\begin{lemma}\label{lemma:kv11}
 Suppose $\mathfrak{t}|\mathbb{I}(\mathfrak{t})|\lesssim\N^{\gamma_{\mathrm{av}}}$ (see {Proposition \ref{prop:kv1}} for $\gamma_{\mathrm{av}}$). Recall $\gamma_{\mathrm{KV}}$ in {Proposition \ref{prop:kv1}}. Consider the processes $\t\mapsto(\mathbf{J}(\t,\inf\mathbb{I}(\mathfrak{t});\mathbb{I}(\mathfrak{t})),\mathbf{U}^{\t,\cdot}[\mathbb{I}(\mathfrak{t})])$ and $\t\mapsto(\mathbf{J}(\t;\mathbb{I}(\mathfrak{t}),\mathcal{E}),\mathbf{U}^{\t,\cdot}[\mathbb{I}(\mathfrak{t}),\mathcal{E}])$ from {Definitions \ref{definition:le10}, \ref{definition:kv4}}. We take $\t\in\mathfrak{t}(\mathrm{in})+[0,\mathfrak{t}]$. Suppose the former has law at $\t=\mathfrak{t}(\mathrm{in})$ distributed as $\mathbb{P}^{\mathrm{Leb},\sigma,\mathfrak{t}(\mathrm{in}),\mathbb{I}(\mathfrak{t})}$ in {Definition \ref{definition:le5}}. Suppose the latter has law at $\t=\mathfrak{t}(\mathrm{in})$ distributed as $\mathbb{P}^{\mathrm{Leb},\sigma,\mathfrak{t}(\mathrm{in}),\mathbb{I}(\mathfrak{t}),\mathcal{E}}$ in {Definition \ref{definition:kv4}}. There exists a coupling of the two initial data so that the probability these two processes are not the same for all $\t\in\mathfrak{t}(\mathrm{in})+[0,\mathfrak{t}]$ is $\lesssim\exp[-\N^{\gamma_{\mathrm{KV}}/3}]$.
\end{lemma}
\begin{proof}
If $\mathbf{U}^{\mathfrak{t}(\mathrm{in}),\cdot}[\mathbb{I}(\mathfrak{t})]\in\mathcal{E}[\mathfrak{t}(\mathrm{in}),\mathfrak{t};\mathbb{I}(\mathfrak{t})]$, then take the initial data of the two joint processes to be equal. Otherwise, sample them independently. Since the two joint processes are the same SDE for $\mathfrak{t}(\mathrm{in})\leq\t\leq\tau[\mathcal{E}]$ (see Definition \ref{definition:kv4}), it suffices to show
\begin{align}
\mathbb{P}\{\tau[\mathcal{E}]\wedge(\mathfrak{t}(\mathrm{in})+\mathfrak{t})\neq\mathfrak{t}(\mathrm{in})+\mathfrak{t}\} \ \lesssim \ \exp\{\N^{-\frac13\gamma_{\mathrm{KV}}}\}. \label{eq:kv11I1}
\end{align}
If $\mathbf{U}^{\t,\cdot}[\mathbb{I}(\mathfrak{t})]\not\in\mathcal{E}[\mathfrak{t}(\mathrm{in}),\mathfrak{t};\mathbb{I}(\mathfrak{t})]$, then by Definition \ref{definition:kv4}, we know $|\partial_{\s}\mathscr{HP}(\s,\mathbf{U}^{\t,\cdot}[\mathbb{I}(\mathfrak{t})];\sigma)|\gtrsim\N^{\gamma_{\mathrm{KV}}}|\mathbb{I}(\mathfrak{t})|^{1/2}$ for some $\mathfrak{t}(\mathrm{in})\leq\s\leq\mathfrak{t}(\mathrm{in})+\mathfrak{t}$. By definition of $\tau[\mathcal{E}]$ as an exit time for $\mathcal{E}[\mathfrak{t}(\mathrm{in}),\mathfrak{t};\mathbb{I}(\mathfrak{t})]$, we deduce
\begin{align}
\mathrm{LHS}\eqref{eq:kv11I1} \ \lesssim \ \mathbb{P}\{{\textstyle\sup_{\s,\t}}|\partial_{\s}\mathscr{HP}(\s,\mathbf{U}^{\t,\cdot}[\mathbb{I}(\mathfrak{t})];\sigma)| \ \gtrsim \ \N^{\gamma_{\mathrm{KV}}}|\mathbb{I}(\mathfrak{t})|^{\frac12}\}, \label{eq:kv11I2}
\end{align}
where the supremum is for $\s,\t\in\mathfrak{t}(\mathrm{in})+[0,\mathfrak{t}]$. Recall in Lemma \ref{lemma:kv3} that, as a function of $\mathbf{U}\in\R^{\mathbb{I}(\mathfrak{t})}$, the term $\mathscr{HP}(\s,\mathbf{U};\sigma)$ is a sum over $\x\in\mathbb{I}(\mathfrak{t})$ of $\mathscr{UP}(\s,\mathbf{U}(\x);\sigma)$. By Lemma \ref{lemma:kv3}, we know $|\partial_{\s}^{2}\mathscr{UP}(\s,\cdot;\sigma)|\lesssim\N^{100}$. Thus, $|\partial_{\s}^{2}\mathscr{HP}(\s,\mathbf{U};\sigma)|\lesssim\N^{200}$. We also know, by construction of $\mathscr{UP}(\s,\mathbf{U}(\x);\sigma)$ as $\mathscr{U}(\s,\mathbf{U}(\x);\sigma)$ plus a $\mathbf{U}$-independent constant, that $\partial_{\s}\mathscr{UP}(\s,\mathbf{U}(\x);\sigma)$ is Lipschitz in $\mathbf{U}$ with Lipschitz norm $\lesssim1$ bounded uniformly in $\s,\sigma$; see Assumption \ref{ass:intro8}. Thus, we know that $\partial_{\s}\mathscr{HP}(\s,\mathbf{U};\sigma)$ is Lipschitz in $\mathbf{U}$ with Lipschitz norm $\lesssim\N^{200}$ bounded uniformly in $\s,\sigma$. At least intuitively, the sup from the previous display is then controlled by values on some very fine polynomial-in-$\N$-sized discretization of $(\s,\t)\in[\mathfrak{t}(\mathrm{in}),\mathfrak{t}(\mathrm{in})+\mathfrak{t}]\times[\mathfrak{t}(\mathrm{in}),\mathfrak{t}(\mathrm{in})+\mathfrak{t}]$. By a union bound over the $\lesssim\exp[\N^{\gamma_{\mathrm{av}}}]$-many points in this discretization, we may pull the double-sup outside $\mathbb{P}$ on the RHS of the previous display, if we include a factor of $\lesssim\exp[\N^{\gamma_{\mathrm{av}}}]$. This is done precisely by Lemma \ref{lemma:ste}. It says this is exactly the case if $\mathbf{U}^{\t,\cdot}[\mathbb{I}(\mathfrak{t})]$ is not extremely wild (which we show will not be the case with extremely high probability). So, $\mathrm{RHS}\eqref{eq:kv11I2}$ is {(for some $\mathrm{C}=\mathrm{O}(1)$)}
\begin{align}
&\lesssim \ \exp[\N^{\gamma_{\mathrm{av}}}]{\textstyle\sup_{\s,\t}}\mathbb{P}\{|\partial_{\s}\mathscr{HP}(\s,\mathbf{U}^{\t,\cdot}[\mathbb{I}(\mathfrak{t})];\sigma)| \ \gtrsim \ \N^{\gamma_{\mathrm{KV}}}|\mathbb{I}(\mathfrak{t})|^{\frac12}\}\nonumber\\
&+\exp[\N^{\gamma_{\mathrm{av}}}]{\textstyle\sup_{\t,\x}}\mathbb{P}\{|\mathbf{U}^{\t,\x}[\mathbb{I}(\mathfrak{t})]|\gtrsim\N^{{\mathrm{C}}}\}+\exp[-\N^{99}]. \label{eq:kv11I2new}
\end{align}
The first double-sup on the RHS is over $\s,\t\in\mathfrak{t}(\mathrm{in})+[0,\mathfrak{t}]$. The second is over $\t\in\mathfrak{t}(\mathrm{in})+[0,\mathfrak{t}]$ and $\x\in\mathbb{I}(\mathfrak{t})$. We now claim
\begin{align}
\mathbb{P}\{|\partial_{\s}\mathscr{HP}(\s,\mathbf{U}^{\t,\cdot}[\mathbb{I}(\mathfrak{t})];\sigma)| \ \gtrsim \ \N^{\gamma_{\mathrm{KV}}}|\mathbb{I}(\mathfrak{t})|^{\frac12}\} \ &\lesssim \ \exp[\mathrm{O}(\N^{\gamma_{\mathrm{av}}})]\mathbb{P}^{\sigma,\t,\mathbb{I}(\mathfrak{t})}\{\mathcal{E}[\mathfrak{t}(\mathrm{in}),\mathfrak{t};\mathbb{I}(\mathfrak{t})]^{\mathrm{C}}\} \nonumber\\
&\lesssim \ \exp[-\N^{\frac13\gamma_{\mathrm{KV}}}]. \label{eq:kv11I2new2}
\end{align}
The first bound follows by the change-of-measure estimate \eqref{eq:kv8I} and construction of $\mathcal{E}[\mathfrak{t}(\mathrm{in}),\mathfrak{t};\mathbb{I}(\mathfrak{t})]$; see Definition \ref{definition:kv4}. The second bound follows by Lemma \ref{lemma:kv6}, which gives $\mathbb{P}^{\sigma,\t,\mathbb{I}(\mathfrak{t})}\{\mathcal{E}[\mathfrak{t}(\mathrm{in}),\mathfrak{t};\mathbb{I}(\mathfrak{t})]^{\mathrm{C}}\}\lesssim\exp[-\N^{\gamma_{\mathrm{KV}}/2}]$, and the fact $\N^{\gamma_{\mathrm{av}}}\ll\N^{\gamma_{\mathrm{KV}}/2}$; see Proposition \ref{prop:kv1}. Next, we claim the following estimate, which we explain afterwards:
\begin{align}
\mathbb{P}\{|\mathbf{U}^{\t,\x}[\mathbb{I}(\mathfrak{t})]|\gtrsim\N^{{\mathrm{C}}}\} \ \lesssim \ \exp[\mathrm{O}(\N^{\gamma_{\mathrm{av}}})]\mathbb{P}^{\sigma,\t,\mathbb{I}(\mathfrak{t})}[|\mathbf{U}(\x)|\gtrsim\N^{{\mathrm{C}}}] \ \lesssim \ \exp[-\N^{99}]. \label{eq:kv11I2new3}
\end{align}
The first bound follows from the change-of-measure estimate \eqref{eq:kv8I}. The second bound follows from the fact that if one takes a random walk with sub-Gaussian step distribution (see Assumption \ref{ass:intro8}) and conditions on the average drift to be $\sigma\lesssim1$, then one just gets a random walk bridge whose steps are also sub-Gaussian with mean $\sigma\lesssim1$ and variance parameter $\lesssim1$. (This bounds the middle probability above by an exponential factor that beats $\exp[\mathrm{O}(\N^{\gamma_{\mathrm{av}}})]$.) Combining the last four displays gives \eqref{eq:kv11I1}. (Indeed, {the factor of $\exp[\N^{\gamma_{\mathrm{av}}}]$ on the RHS of \eqref{eq:kv11I2new} is overwhelmed by the exponential decay on the RHS of \eqref{eq:kv11I2new2} and \eqref{eq:kv11I2new3}, respectively}.) This completes the proof (as we noted immediately before \eqref{eq:kv11I1}).
\end{proof}
\subsection{Time-reversing the non-equilibrium SDEs \eqref{eq:glsdeloc}-\eqref{eq:hfloc}}
See Section \ref{section:sqle} for {the} motivation of this subsection.
\begin{definition}\label{definition:kv12}
 Fix $\mathfrak{t}(\mathrm{in}),\mathfrak{t}\geq0$. Consider $\t\mapsto(\mathbf{J}(\t,\inf\mathbb{I}(\mathfrak{t});\mathbb{I}(\mathfrak{t})),\mathbf{U}^{\t,\cdot}[\mathbb{I}(\mathfrak{t})])$ for $\t\in\mathfrak{t}(\mathrm{in})+[0,\mathfrak{t}]$. We assume that its law at $\t=\mathfrak{t}(\mathrm{in})$ is distributed as $\mathbb{P}^{\mathrm{Leb},\sigma,\mathfrak{t}(\mathrm{in}),\mathbb{I}(\mathfrak{t})}$ in Definition \ref{definition:le5}. Define a time-reversal map $\t\mapsto\t[\ast]$ on $\mathfrak{t}(\mathrm{in})+[0,\mathfrak{t}]\to\mathfrak{t}(\mathrm{in})+[0,\mathfrak{t}]$. This is the unique linear function such that $[\mathfrak{t}(\mathrm{in})+\mathfrak{t}][\ast]=\mathfrak{t}(\mathrm{in})$ and $\mathfrak{t}(\mathrm{in})[\ast]=\mathfrak{t}(\mathrm{in})+\mathfrak{t}$. We set $\mathbf{U}^{\t,\cdot,\ast}[\mathbb{I}(\mathfrak{t})]:=\mathbf{U}^{\t[\ast],\cdot}[\mathbb{I}(\mathfrak{t})]$ and $\mathbf{J}^{\ast}(\t,\cdot;\mathbb{I}(\mathfrak{t})):=\mathbf{J}(\t[\ast],\cdot;\mathbb{I}(\mathfrak{t}))$ for our notation for time-reversals of $\t\mapsto(\mathbf{J}(\t,\inf\mathbb{I}(\mathfrak{t});\mathbb{I}(\mathfrak{t})),\mathbf{U}^{\t,\cdot}[\mathbb{I}(\mathfrak{t})])$.
\end{definition}
Time-reversing SDEs is a well understood procedure. We now put it into practice below.
\begin{lemma}\label{lemma:kv13}
 We have the following SDE for $\t\mapsto\mathbf{U}^{\t,\cdot,\ast}[\mathbb{I}(\mathfrak{t})]$ for times $\t\in\mathfrak{t}(\mathrm{in})+[0,\mathfrak{t}]$, which we explain after:
\begin{align}
\d\mathbf{U}^{\t,\x,\ast}[\mathbb{I}(\mathfrak{t})] \ &= \ \N^{2}\Delta^{\mathbb{I}(\mathfrak{t})}\mathscr{U}'(\t[\ast],\mathbf{U}^{\t,\x,\ast}[\mathbb{I}(\mathfrak{t})])\d\t-\N^{\frac32}\grad^{\mathbb{I}(\mathfrak{t}),\mathrm{a}}\mathscr{U}'(\t[\ast],\mathbf{U}^{\t,\x,\ast}[\mathbb{I}(\mathfrak{t})])\d\t\nonumber\\
&+ \ \mathscr{I}(\t[\ast],\x)\d\t- \sqrt{2}\N\grad^{\mathbb{I}(\mathfrak{t}),-}\d\mathbf{b}(\t[\ast],\x). \nonumber
\end{align}
{The} $\Delta^{\mathbb{I}(\mathfrak{t})}$ and $\grad^{\mathbb{I}(\mathfrak{t}),?}$ operators can be found in {Definition \ref{definition:le10}}. To construct $\mathscr{I}$, first let $\mathfrak{p}[\mathfrak{t}(\mathrm{in}),\t]$ be the Radon-Nikodym derivative for the law of $\mathbf{U}^{\t,\cdot}[\mathbb{I}(\mathfrak{t})]$ with respect to $\mathbb{P}^{\sigma,\t,\mathbb{I}(\mathfrak{t})}$. The drift $\mathscr{I}$ above is defined as
\begin{align}
\mathscr{I}(\t,\x) \ := \ 2\N^{2}\mathfrak{p}[\mathfrak{t}(\mathrm{in}),\t]^{-1}\grad^{\mathbb{I}(\mathfrak{t}),-}\grad^{\mathbb{I}(\mathfrak{t}),+}\partial_{\mathbf{U}(\x)}\mathfrak{p}[\mathfrak{t}(\mathrm{in}),\t]. \label{eq:kv13I}
\end{align}
To be clear, we evaluate $\mathrm{RHS}\eqref{eq:kv13I}$ at $\mathbf{U}^{\t,\cdot}[\mathbb{I}(\mathfrak{t})]$, so $\mathscr{I}(\t[\ast],\x)$ is evaluated at $\mathbf{U}^{\t,\cdot,\ast}[\mathbb{I}(\mathfrak{t})]$. We also claim $\mathbf{J}^{\ast}(\t,\inf\mathbb{I}(\mathfrak{t});\mathbb{I}(\mathfrak{t}))=\Pi^{\mathbb{S}(\N)}\mathbf{J}^{\ast}(\t;\mathbb{I}(\mathfrak{t}))$, where $\mathbf{J}^{\ast}(\t;\mathbb{I}(\mathfrak{t}))=\mathbf{J}(\t[\ast];\mathbb{I}(\mathfrak{t}))$; for relevant notation, see {Definition \ref{definition:le10}}. Moreover, we have the SDE
\begin{align}
\d\mathbf{J}^{\ast}(\t;\mathbb{I}(\mathfrak{t})) \ &= \ -\N^{\frac32}\grad^{\mathbb{I}(\mathfrak{t}),+}\mathscr{U}'(\t[\ast],\mathbf{U}^{\t,\inf\mathbb{I}(\mathfrak{t}),\ast}[\mathbb{I}(\mathfrak{t})])\d\t+\mathscr{R}(\t[\ast])\d\t+\sqrt{2}\N^{\frac12}\d\mathbf{b}(\t[\ast],\x) \label{eq:kv13IIa}\\
&- \ \N\{\mathscr{U}'(\t[\ast],\mathbf{U}^{\t,\inf\mathbb{I}(\mathfrak{t}),\ast}[\mathbb{I}(\mathfrak{t})])+\mathscr{U}'(\t[\ast],\mathbf{U}^{\t,\inf\mathbb{I}(\mathfrak{t})+1,\ast}[\mathbb{I}(\mathfrak{t})])\}\d\t. \label{eq:kv13IIb}
\end{align}
\end{lemma}
\begin{proof}
Time-reverse $\t\mapsto\mathrm{Joint}(\t):=(\mathbf{J}(\t,\inf\mathbb{I}(\mathfrak{t});\mathbb{I}(\mathfrak{t})),\mathbf{U}^{\t,\cdot}[\mathbb{I}(\mathfrak{t})])$. Let us first examine what happens to $\mathbf{J}(\t,\inf\mathbb{I}(\mathfrak{t});\mathbb{I}(\mathfrak{t}))$. To this end, {we give} a preliminary observation. By Definition \ref{definition:le10}, we know $\mathbf{J}(\t,\inf\mathbb{I}(\mathfrak{t});\mathbb{I}(\mathfrak{t}))$ is the solution of an SDE on the manifold $\mathbb{S}(\N)$; see Remark \ref{remark:le6}. Its initial data at $\t=\mathfrak{t}(\mathrm{in})$ has law $\mathrm{Leb}(\mathbb{S}(\N))$, which is the invariant measure for said $\mathbb{S}(\N)$-SDE. (For this last claim, technically, by invariant measure, we mean that given any $\t\geq\mathfrak{t}(\mathrm{in})$, we know $\mathbf{J}(\t,\inf\mathbb{I}(\mathfrak{t});\mathbb{I}(\mathfrak{t}))\sim\mathrm{Leb}(\mathbb{S}(\N))$. This can be seen via Kolmogorov PDEs as discussed in Remark \ref{remark:kv9}. See also Lemma \ref{lemma:le7}.) We now use Theorem 2.1 in \cite{HP}. This says that time-reversing the $\mathbb{S}(\N)$ SDE, whose solution is $\mathbf{J}(\t,\inf\mathbb{I}(\mathfrak{t});\mathbb{I}(\mathfrak{t}))$, just requires adding a sign to each drift term in \eqref{eq:hfloc}, reversing time via $\t\mapsto\t[\ast]$, and adding another drift term given by derivatives of log of the Radon-Nikodym derivative for the law of $\mathbf{J}(\t,\inf\mathbb{I}(\mathfrak{t});\mathbb{I}(\mathfrak{t}))$ with respect to $\mathrm{Leb}(\mathbb{S}(\N))$. (Strictly speaking, \cite{HP} addresses Euclidean SDEs, though its proof is entirely based on the Ito formula, which certainly holds for SDEs on tori.) But, this density is constant as we have just argued, so the additional drift is zero. This shows that $\mathbf{J}^{\ast}(\t,\inf\mathbb{I}(\mathfrak{t});\mathbb{I}(\mathfrak{t}))$ is the solution {for} the same $\mathbb{S}(\N)$ SDE but with signs for drifts and $\t\mapsto\t[\ast]$. This gives our claim for $\mathbf{J}^{\ast}(\t,\inf\mathbb{I}(\mathfrak{t});\mathbb{I}(\mathfrak{t}))$ (and \eqref{eq:kv13IIa}-\eqref{eq:kv13IIb}). We move to $\mathbf{U}^{\t,\cdot,\ast}[\mathbb{I}(\mathfrak{t})]$. We claim $\d\mathbf{U}^{\t,\x,\ast}[\mathbb{I}(\mathfrak{t})] $ is
\begin{align}
&-\N^{2}\Delta^{\mathbb{I}(\mathfrak{t})}\mathscr{U}'(\t[\ast],\mathbf{U}^{\t,\x,\ast}[\mathbb{I}(\mathfrak{t})])\d\t-\N^{\frac32}\grad^{\mathbb{I}(\mathfrak{t}),\mathrm{a}}\mathscr{U}'(\t[\ast],\mathbf{U}^{\t,\x,\ast}[\mathbb{I}(\mathfrak{t})])\d\t+\mathscr{I}^{\dagger}(\t[\ast],\x)\d\t \label{eq:kv13I1a}\\
&-\sqrt{2}\N\grad^{\mathbb{I}(\mathfrak{t}),-}\d\mathbf{b}(\t[\ast],\x), \nonumber
\end{align}
where $\mathscr{I}^{\dagger}$ is the following modification of $\mathscr{I}$ obtained by replacing $\mathfrak{p}[\mathfrak{t}(\mathrm{in}),\t]$ with $\mathfrak{p}[\mathfrak{t}(\mathrm{in}),\t]^{\dagger}$, which we define to be the Radon-Nikodym derivative for $\mathbf{U}^{\t,\cdot}[\mathbb{I}(\mathfrak{t})]$ with respect to Lebesgue measure on $\mathbb{H}^{\sigma,\mathbb{I}(\mathfrak{t})}$:
\begin{align}
\mathscr{I}^{\dagger}(\t,\x) \ := \ 2\N^{2}\{\mathfrak{p}[\mathfrak{t}(\mathrm{in}),\t]^{\dagger}\}^{-1}\grad^{\mathbb{I}(\mathfrak{t}),-}\grad^{\mathbb{I}(\mathfrak{t}),+}\partial_{\mathbf{U}(\x)}\mathfrak{p}[\mathfrak{t}(\mathrm{in}),\t]^{\dagger}. \label{eq:kv13I1b}
\end{align}
{Indeed, we have added a sign to the drifts in \eqref{eq:glsdeloc} like Theorem 2.1 in \cite{HP} says to do. The additional drift \eqref{eq:kv13I1b} comes from the exact form of what Theorem 2.1 in \cite{HP} tells us we must add. (Indeed, it tells us to add Doob transform drift that we explain. {Consider $\log\mathfrak{p}[\mathfrak{t}(\mathrm{in}),\t]^{\dagger}$ and take its $\mathbf{U}(\x)$-partial.} {Next, consider} the operator hitting Brownian motions; this is $\sqrt{2}\N\grad^{\mathbb{I}(\mathfrak{t}),-}$, up to an unimportant sign. Take the associated covariance matrix. Now hit {\small$\partial_{\mathbf{U}(\x)}$}$\log\mathfrak{p}[\mathfrak{t}(\mathrm{in}),\t]^{\dagger}$ with this covariance matrix. This gives us \eqref{eq:kv13I1b}. Thus \eqref{eq:kv13I1a}-\eqref{eq:kv13I1b} follows.} We now claim the following, in which $\mathbf{U}:=\mathbf{U}^{\t,\cdot}[\mathbb{I}(\mathfrak{t})]$, and we use notation in Lemma \ref{lemma:kv3}:
\begin{align}
\mathscr{I}^{\dagger}(\t,\x) \ &= \ 2\N^{2}\{\mathfrak{p}[\mathfrak{t}(\mathrm{in}),\t]^{\dagger}\}^{-1}\grad^{\mathbb{I}(\mathfrak{t}),-}\grad^{\mathbb{I}(\mathfrak{t}),+}\partial_{\mathbf{U}(\x)}\{\mathfrak{p}[\mathfrak{t}(\mathrm{in}),\t]\exp[-\mathscr{HP}(\t,\mathbf{U};\sigma)]\} \label{eq:kv13I2a}\\
&= \ 2\N^{2}\{\mathfrak{p}[\mathfrak{t}(\mathrm{in}),\t]^{\dagger}\}^{-1}\grad^{\mathbb{I}(\mathfrak{t}),-}\grad^{\mathbb{I}(\mathfrak{t}),+}(\partial_{\mathbf{U}(\x)}\{\mathfrak{p}[\mathfrak{t}(\mathrm{in}),\t]\}\times\exp[-\mathscr{HP}(\t,\mathbf{U};\sigma)])\label{eq:kv13I2b}\\
&+ \ 2\N^{2}\{\mathfrak{p}[\mathfrak{t}(\mathrm{in}),\t]^{\dagger}\}^{-1}\grad^{\mathbb{I}(\mathfrak{t}),-}\grad^{\mathbb{I}(\mathfrak{t}),+}(\mathfrak{p}[\mathfrak{t}(\mathrm{in}),\t]\times\partial_{\mathbf{U}(\x)}\exp[-\mathscr{HP}(\t,\mathbf{U};\sigma)])\label{eq:kv13I2c} \\
&= \ 2\N^{2}\{\mathfrak{p}[\mathfrak{t}(\mathrm{in}),\t]^{\dagger}\}^{-1}\times\exp[-\mathscr{HP}(\t,\mathbf{U};\sigma)]\grad^{\mathbb{I}(\mathfrak{t}),-}\grad^{\mathbb{I}(\mathfrak{t}),+}\partial_{\mathbf{U}(\x)}\{\mathfrak{p}[\mathfrak{t}(\mathrm{in}),\t]\}\label{eq:kv13I2d}\\
&+ \ 2\N^{2}\{\mathfrak{p}[\mathfrak{t}(\mathrm{in}),\t]^{\dagger}\}^{-1}\mathfrak{p}[\mathfrak{t}(\mathrm{in}),\t]\grad^{\mathbb{I}(\mathfrak{t}),-}\grad^{\mathbb{I}(\mathfrak{t}),+}\partial_{\mathbf{U}(\x)}\exp[-\mathscr{HP}(\t,\mathbf{U};\sigma)]\label{eq:kv13I2e} \\
&= \ 2\N^{2}\mathfrak{p}[\mathfrak{t}(\mathrm{in}),\t]^{-1}\grad^{\mathbb{I}(\mathfrak{t}),-}\grad^{\mathbb{I}(\mathfrak{t}),+}\partial_{\mathbf{U}(\x)}\{\mathfrak{p}[\mathfrak{t}(\mathrm{in}),\t]\}\label{eq:kv13I2f}\\
&+ \ 2\N^{2}\exp[\mathscr{HP}(\t,\mathbf{U};\sigma)]\grad^{\mathbb{I}(\mathfrak{t}),-}\grad^{\mathbb{I}(\mathfrak{t}),+}\{-\mathscr{U}'(\t,\mathbf{U}(\x))\exp[-\mathscr{HP}(\t,\mathbf{U};\sigma)]\}\label{eq:kv13I2g} \\
&= \ 2\N^{2}\mathfrak{p}[\mathfrak{t}(\mathrm{in}),\t]^{-1}\grad^{\mathbb{I}(\mathfrak{t}),-}\grad^{\mathbb{I}(\mathfrak{t}),+}\partial_{\mathbf{U}(\x)}\{\mathfrak{p}[\mathfrak{t}(\mathrm{in}),\t]\}-2\N^{2}\grad^{\mathbb{I}(\mathfrak{t}),-}\grad^{\mathbb{I}(\mathfrak{t}),+}\mathscr{U}'(\t,\mathbf{U}(\x)) \label{eq:kv13I2h} \\
&= \ \mathscr{I}(\t,\x) + 2\N^{2}\Delta^{\mathbb{I}(\mathfrak{t})}\mathscr{U}'(\t,\mathbf{U}(\x)). \label{eq:kv13I2i}
\end{align}
{\eqref{eq:kv13I2a} {comes from} \eqref{eq:kv13I1b} and writing $\mathfrak{p}[\mathfrak{t}(\mathrm{in}),\t]^{\dagger}=\mathfrak{p}[\mathfrak{t}(\mathrm{in}),\t]\exp[-\mathscr{HP}(\t,\mathbf{U};\sigma)]$. Indeed, going from $\mathfrak{p}[\mathfrak{t}(\mathrm{in}),\t]^{\dagger}$ to $\mathfrak{p}[\mathfrak{t}(\mathrm{in}),\t]$ requires a change-of-measure factor \eqref{eq:kv3I}. \eqref{eq:kv13I2b}-\eqref{eq:kv13I2c} follows by the Leibniz rule for the $\mathbf{U}(\x)$-partial. \eqref{eq:kv13I2d}-\eqref{eq:kv13I2e} follows by pulling out from $\grad$-operators whatever term the $\mathbf{U}(\x)$-partial does not hit in \eqref{eq:kv13I2b}-\eqref{eq:kv13I2c}. (Indeed, the $\grad$-operators act on $\x$, which is only present in the $\mathbf{U}(\x)$-partial.) \eqref{eq:kv13I2f}-\eqref{eq:kv13I2g} holds for the following reason. In \eqref{eq:kv13I2d}, we again use the identity $\mathfrak{p}[\mathfrak{t}(\mathrm{in}),\t]^{\dagger}=\mathfrak{p}[\mathfrak{t}(\mathrm{in}),\t]\exp[-\mathscr{HP}(\t,\mathbf{U};\sigma)]$. This gives us \eqref{eq:kv13I2f}. In \eqref{eq:kv13I2e}, we again use this identity. To compute the partial therein, recall from Lemma \ref{lemma:kv3} that $\mathscr{HP}(\t,\mathbf{U};\sigma)$ is the sum over $\x\in\mathbb{I}(\mathfrak{t})$ of $\mathscr{UP}(\t,\mathbf{U}(\x);\sigma)$, which itself is just $\mathscr{U}(\t,\mathbf{U}(\x))$ plus something independent of $\mathbf{U}$. \eqref{eq:kv13I2h} follows from moving $\exp[-\mathscr{HP}]$ outside of the $\grad$-operators in \eqref{eq:kv13I2g}. This cancels the $\exp[\mathscr{HP}]$ factor therein. \eqref{eq:kv13I2i} follows by \eqref{eq:kv13I} and noting $\Delta^{\mathbb{I}(\mathfrak{t})}=-\grad^{\mathbb{I}(\mathfrak{t}),-}\grad^{\mathbb{I}(\mathfrak{t}),+}$. (This follows from construction; see Definition \ref{definition:le10}. In particular, $\Delta^{\mathbb{I}(\mathfrak{t})}$ is {the} discrete Laplacian, $\grad^{\mathbb{I}(\mathfrak{t}),+}$ is {the} discrete gradient, and $\grad^{\mathbb{I}(\mathfrak{t}),-}$ is a discrete gradient with a negative sign.) The $\mathbf{U}^{\t,\cdot,\ast}[\mathbb{I}(\mathfrak{t})]$ SDE we claimed now follows by \eqref{eq:kv13I1a} and \eqref{eq:kv13I2a}-\eqref{eq:kv13I2i}.}
\end{proof}
As discussed in Section \ref{section:sqle}, we now want to remove the $\mathscr{I}$-drift in the $\mathbf{U}^{\t,\cdot,\ast}[\mathbb{I}(\mathfrak{t})]$ SDE in Lemma \ref{lemma:kv13} by Girsanov.
\begin{definition}\label{definition:kv14}
 Take the setting of Definition \ref{definition:kv12} and $\Delta^{\mathbb{I}(\mathfrak{t})}$, $\grad^{\mathbb{I}(\mathfrak{t}),?}$ in Definition \ref{definition:le10}. Now, let $\t\mapsto(\mathbf{J}^{\sim}(\t;\mathbb{I}(\mathfrak{t})),\mathbf{U}^{\t,\cdot,\sim}[\mathbb{I}(\mathfrak{t})])$ solve the following SDE, which uses the notation $\mathbf{U}^{\t,\sim}:=\mathbf{U}^{\t,\inf\mathbb{I}(\mathfrak{t}),\sim}[\mathbb{I}(\mathfrak{t})]$ and $\mathbf{U}^{\t,\sim,+}:=\mathbf{U}^{\t,\inf\mathbb{I}(\mathfrak{t})+1,\sim}[\mathbb{I}(\mathfrak{t})]$ and $\mathbf{b}=\mathbf{b}(\t[\ast],\x)$:
{\small
\begin{align}
\d\mathbf{U}^{\t,\x,\sim}[\mathbb{I}(\mathfrak{t})] &= \N^{2}\Delta^{\mathbb{I}(\mathfrak{t})}\mathscr{U}'(\t[\ast],\mathbf{U}^{\t,\x,\sim}[\mathbb{I}(\mathfrak{t})])\d\t-\N^{\frac32}\grad^{\mathbb{I}(\mathfrak{t}),\mathrm{a}}\mathscr{U}'(\t[\ast],\mathbf{U}^{\t,\x,\sim}[\mathbb{I}(\mathfrak{t})])\d\t- \sqrt{2}\N\grad^{\mathbb{I}(\mathfrak{t}),-}\d\mathbf{b}\label{eq:kv14I}\\
\d\mathbf{J}^{\sim}(\t;\mathbb{I}(\mathfrak{t})) &=-\N^{\frac32}\grad^{\mathbb{I}(\mathfrak{t}),+}\mathscr{U}'(\t[\ast],\mathbf{U}^{\t,\sim})\d\t-\N\{\mathscr{U}'(\t[\ast],\mathbf{U}^{\t,\sim})+\mathscr{U}'(\t[\ast],\mathbf{U}^{\t,\sim,+})\}\d\t+\mathscr{R}(\t[\ast])\d\t\label{eq:kv14II}\\
&+\sqrt{2}\N^{\frac12}\d\mathbf{b}. \nonumber
\end{align}
}Finally, let us set the torus-projection $\mathbf{J}^{\sim}(\t,\inf\mathbb{I}(\mathfrak{t});\mathbb{I}(\mathfrak{t})):=\Pi^{\mathbb{S}(\N)}\mathbf{J}^{\sim}(\t;\mathbb{I}(\mathfrak{t}))$, in spirit of what is done in Definition \ref{definition:le10}.
\end{definition}
In words, the only difference between the superscript $\ast$ processes in Lemma \ref{lemma:kv13} and the $\sim$ processes above is the lack of an $\mathscr{I}$-drift in the latter. So, we can estimate the cost in removing this drift via Girsanov for joint processes. The whole point of the following result and its proof is to rigorously and carefully estimate this cost by using our energy estimate of {Lemma \ref{lemma:kv10}} to control said $\mathscr{I}$-drift \eqref{eq:kv13I}.
\begin{lemma}\label{lemma:kv15}
 Fix $\mathfrak{t}(\mathrm{in}),\mathfrak{t}\geq0$. Retain {the} notation of {Definitions \ref{definition:kv12}, \ref{definition:kv14}} and {Lemma \ref{lemma:kv13}}. For $\t\in\mathfrak{t}(\mathrm{in})+[0,\mathfrak{t}]$, we define the ``forward-backward processes" below, the first of which is for the time-reversed process with superscript $\ast$ in Lemma \ref{lemma:kv13}, and the second of which is for the auxiliary modification with superscript $\sim$ in Definition \ref{definition:kv14}:
\begin{align}
\t &\mapsto \ \mathbf{FB}^{\ast}[\t;\mathbb{I}(\mathfrak{t})] \ := \ (\mathbf{J}(\t,\inf\mathbb{I}(\mathfrak{t});\mathbb{I}(\mathfrak{t})),\mathbf{U}^{\t,\cdot}[\mathbb{I}(\mathfrak{t})],\mathbf{J}^{\ast}(\t,\inf\mathbb{I}(\mathfrak{t});\mathbb{I}(\mathfrak{t})),\mathbf{U}^{\t,\cdot,\ast}[\mathbb{I}(\mathfrak{t})]) \label{eq:kv15I1a}\\
\t &\mapsto \ \mathbf{FB}^{\sim}[\t;\mathbb{I}(\mathfrak{t})] \ := \ (\mathbf{J}(\t,\inf\mathbb{I}(\mathfrak{t});\mathbb{I}(\mathfrak{t})),\mathbf{U}^{\t,\cdot}[\mathbb{I}(\mathfrak{t})],\mathbf{J}^{\sim}(\t,\inf\mathbb{I}(\mathfrak{t});\mathbb{I}(\mathfrak{t})),\mathbf{U}^{\t,\cdot,\sim}[\mathbb{I}(\mathfrak{t})]). \label{eq:kv15I1b}
\end{align}
Suppose $(\mathbf{J}(\t,\inf\mathbb{I}(\mathfrak{t});\mathbb{I}(\mathfrak{t})),\mathbf{U}^{\t,\cdot}[\mathbb{I}(\mathfrak{t})])$ at time $\t=\mathfrak{t}(\mathrm{in})$ is distributed according to the law $\mathbb{P}^{\mathrm{Leb},\sigma,\mathfrak{t}(\mathrm{in}),\mathbb{I}(\mathfrak{t})}$, and $(\mathbf{J}^{\sim}(\t,\inf\mathbb{I}(\mathfrak{t});\mathbb{I}(\mathfrak{t})),\mathbf{U}^{\t,\cdot,\sim}[\mathbb{I}(\mathfrak{t})])$ at time $\t=\mathfrak{t}(\mathrm{in})$ is distributed according to the law of $(\mathbf{J}(\t,\inf\mathbb{I}(\mathfrak{t});\mathbb{I}(\mathfrak{t})),\mathbf{U}^{\t,\cdot}[\mathbb{I}(\mathfrak{t})])$ at time $\t=\mathfrak{t}(\mathrm{in})+\mathfrak{t}$. (This is also the law of $(\mathbf{J}^{\ast}(\t,\inf\mathbb{I}(\mathfrak{t});\mathbb{I}(\mathfrak{t})),\mathbf{U}^{\t,\cdot,\ast}[\mathbb{I}(\mathfrak{t})])$ at time $\t=\mathfrak{t}(\mathrm{in})$.) Now, let $\mathrm{RN}[\mathfrak{t}(\mathrm{in}),\mathfrak{t};\mathbb{I}(\mathfrak{t})]$ be the Radon-Nikodym derivative of $\mathbf{FB}^{\ast}[\t;\mathbb{I}(\mathfrak{t})]$ with respect to $\mathbf{FB}^{\sim}[\t;\mathbb{I}(\mathfrak{t})]$ (as probability measures on the path space for times $\t\in\mathfrak{t}(\mathrm{in})+[0,\mathfrak{t}]$). 

Now, assume $\mathfrak{t}|\mathbb{I}(\mathfrak{t})|\lesssim\N^{\gamma_{\mathrm{av}}}$; see {Proposition \ref{prop:kv1}} for $\gamma_{\mathrm{av}}$. Also assume $|\mathbb{I}(\mathfrak{t})|\gtrsim\N^{1/10}$. We let $\E^{\sim}$ be {the} path-space expectation with respect to the law of $\mathbf{FB}^{\sim}[\cdot;\mathbb{I}(\mathfrak{t})]$. {Then, we} have the following for $\gamma_{\mathrm{KV}}$ from {Proposition \ref{prop:kv1}} {and for any $\mathrm{D}>0$ large but finite:}
\begin{align}
\E^{\sim}\mathrm{RN}[\mathfrak{t}(\mathrm{in}),\mathfrak{t};\mathbb{I}(\mathfrak{t})]\log\mathrm{RN}[\mathfrak{t}(\mathrm{in}),\mathfrak{t};\mathbb{I}(\mathfrak{t})] \ \lesssim \ \N^{\gamma_{\mathrm{KV}}}\mathfrak{t}|\mathbb{I}(\mathfrak{t})|^{\frac12}+\N^{-{\mathrm{D}}}. \label{eq:kv15II}
\end{align}
\end{lemma}
\begin{proof}
Before we start the proof, we emphasize the following clarification. The last two components of the process $\mathbf{FB}^{\sim}[\t;\mathbb{I}(\mathfrak{t})]$, in principle, have nothing to do with time-reversing its first two components. (The latter two components are just SDEs driven by the time-reversed Brownian motions, and they have initial data determined using the final-time law of the first two components.) In particular, the first two components of $\mathbf{FB}^{\ast}[\t;\mathbb{I}(\mathfrak{t})]$ and $\mathbf{FB}^{\sim}[\t;\mathbb{I}(\mathfrak{t})]$ are the \emph{same} process. Now, define $\mathrm{RN}^{\mathrm{back}}[\mathfrak{t}(\mathrm{in}),\mathfrak{t};\mathbb{I}(\mathfrak{t})]$ as the Radon-Nikodym derivative for the last two components of $\mathbf{FB}^{\ast}[\t;\mathbb{I}(\mathfrak{t})]$ with respect to those of $\mathbf{FB}^{\sim}[\t;\mathbb{I}(\mathfrak{t})]$. (Again, this is as measures on the path space for times in $\mathfrak{t}(\mathrm{in})+[0,\mathfrak{t}]$.) Let $\E^{\sim,\mathrm{back}}$ be {the} expectation with respect to the path-space law of the last two components of $\mathbf{FB}^{\sim}[\t;\mathbb{I}(\mathfrak{t})]$. We claim the following holds (with explanation given afterwards):
\begin{align}
\mathrm{LHS}\eqref{eq:kv15II} \ = \ \E^{\sim,\mathrm{back}}\mathrm{RN}^{\mathrm{back}}[\mathfrak{t}(\mathrm{in}),\mathfrak{t};\mathbb{I}(\mathfrak{t})]\log\mathrm{RN}^{\mathrm{back}}[\mathfrak{t}(\mathrm{in}),\mathfrak{t};\mathbb{I}(\mathfrak{t})]. \label{eq:kv15II1}
\end{align}
By chain rule (or martingale decomposition) for relative entropy, we know \eqref{eq:kv15II1} holds if we include an additional term given by the relative entropy between the first two components of $\mathbf{FB}^{\ast}[\t;\mathbb{I}(\mathfrak{t})]$ with respect to those of $\mathbf{FB}^{\sim}[\t;\mathbb{I}(\mathfrak{t})]$, after \emph{conditioning} on the paths of the last two components for each to be the same. But, the first two components of $\mathbf{FB}^{\ast}[\t;\mathbb{I}(\mathfrak{t})]$ and $\mathbf{FB}^{\sim}[\t;\mathbb{I}(\mathfrak{t})]$ are the same process, as we noted at the beginning of this proof. Thus, their relative entropy is zero regardless of whatever values the other components of $\mathbf{FB}^{\ast}[\t;\mathbb{I}(\mathfrak{t})]$ and $\mathbf{FB}^{\sim}[\t;\mathbb{I}(\mathfrak{t})]$ take. Therefore, the additional relative entropy terms that we must add to \eqref{eq:kv15II1} to make it true are just zero. \eqref{eq:kv15II1} follows. By Girsanov, we can compute
\begin{align}
\mathrm{RN}^{\mathrm{back}}[\mathfrak{t}(\mathrm{in}),\mathfrak{t};\mathbb{I}(\mathfrak{t})] \ = \ \exp\{\mathscr{N}(\mathfrak{t}(\mathrm{in}),\mathfrak{t})+\tfrac12[\mathscr{N}(\mathfrak{t}(\mathrm{in}),\mathfrak{t})]\}, \label{eq:kv15II2a}
\end{align}
where $\mathscr{N}(\mathfrak{t}(\mathrm{in}),\mathfrak{t})$ denotes a stochastic integral, and $[\mathscr{N}(\mathfrak{t}(\mathrm{in}),\mathfrak{t})]$ is its bracket. (This stochastic integral $\mathscr{N}(\mathfrak{t}(\mathrm{in}),\mathfrak{t})$ is adapted to the time-reversed Brownian motion. Indeed, we used Girsanov to SDEs driven by time-reversed Brownian motions.) We claim
\begin{align}
[\mathscr{N}(\mathfrak{t}(\mathrm{in}),\mathfrak{t})] \ \lesssim \ {\textstyle\int_{\mathfrak{t}(\mathrm{in})}^{\mathfrak{t}(\mathrm{in})+\mathfrak{t}}}{\textstyle\sum_{\x\in\mathbb{I}(\mathfrak{t})}}\N^{-2}|[\grad^{\mathbb{I}(\mathfrak{t}),-}]^{-1}\mathscr{I}(\t,\x)|^{2}\d\t. \label{eq:kv15II2b}
\end{align}
{To get \eqref{eq:kv15II2a}-\eqref{eq:kv15II2b}, note the SDEs in Lemma \ref{lemma:kv13} and Definitions \ref{definition:kv14} are the same SDEs, except the $\mathbf{U}^{\t,\cdot,\ast}[\mathbb{I}(\mathfrak{t})]$ SDE in Lemma \ref{lemma:kv13} has additional drift $\mathscr{I}$. Girsanov implies that $\mathrm{RN}^{\mathrm{back}}$ is then the exponential martingale in \eqref{eq:kv15II2a}. It also implies that $\mathscr{N}(\mathfrak{t}(\mathrm{in}),\mathfrak{t})$ has the following bracket process. Take this additional drift $\mathscr{I}(\t,\x)$ per $\x\in\mathbb{I}(\mathfrak{t})$. Apply the inverse of the operator hitting Brownian motions in \eqref{eq:kv14I} to $\mathscr{I}$. This inverse operator is equal to $[\sqrt{2}\N]^{-1}[\grad^{\mathbb{I}(\mathfrak{t}),-}]^{-1}$. Because brackets are quadratic, square the resulting term. Then integrate on the time-domain of interest, which is $\mathfrak{t}(\mathrm{in})+[0,\mathfrak{t}]$. Sum over all $\x\in\mathbb{I}(\mathfrak{t})$. This gives $\mathrm{RHS}\eqref{eq:kv15II2b}$ up to $\mathrm{O}(1)$ factor. Now, let $\E^{\ast,\mathrm{back}}$ be {the} expectation with respect to the path-space law of the last two components of $\mathbf{FB}^{\ast}[\t;\mathbb{I}(\mathfrak{t})]$. Next, we build on \eqref{eq:kv15II1} and claim the following calculation holds, which we explain after:}
\begin{align}
\mathrm{LHS}\eqref{eq:kv15II} \ &= \ \E^{\ast,\mathrm{back}}\log\mathrm{RN}^{\mathrm{back}}[\mathfrak{t}(\mathrm{in}),\mathfrak{t};\mathbb{I}(\mathfrak{t})] \ = \ \E^{\ast,\mathrm{back}}\{\mathscr{N}(\mathfrak{t}(\mathrm{in}),\mathfrak{t})+\tfrac12[\mathscr{N}(\mathfrak{t}(\mathrm{in}),\mathfrak{t})]\} \label{eq:kv15II3a}\\
&= \ \E^{\ast,\mathrm{back}}\tfrac12[\mathscr{N}(\mathfrak{t}(\mathrm{in}),\mathfrak{t})] \ \lesssim \ \E^{\ast,\mathrm{back}} {\textstyle\int_{\mathfrak{t}(\mathrm{in})}^{\mathfrak{t}(\mathrm{in})+\mathfrak{t}}}{\textstyle\sum_{\x\in\mathbb{I}(\mathfrak{t})}}\N^{-2}|[\grad^{\mathbb{I}(\mathfrak{t}),-}]^{-1}\mathscr{I}(\t,\x)|^{2}\d\t. \label{eq:kv15II3b}
\end{align}
\eqref{eq:kv15II3a} follows by construction of $\mathrm{RN}^{\mathrm{back}}[\mathfrak{t}(\mathrm{in}),\mathfrak{t};\mathbb{I}(\mathfrak{t})]$ as the density of $\E^{\ast,\mathrm{back}}$ with respect to $\E^{\sim,\mathrm{back}}$ and then by \eqref{eq:kv15II2a}. \eqref{eq:kv15II3b} follows first by noting $\mathscr{N}(\mathfrak{t}(\mathrm{in}),\mathfrak{t})$ is a stochastic integral and therefore mean-zero. Then, we use \eqref{eq:kv15II2b}. Let us now compute the time-integrand in \eqref{eq:kv15II3b}. By \eqref{eq:kv13I}, using notation of Lemma \ref{lemma:kv13}, we claim the following estimate, in which all $\mathfrak{p}$-factors (and their derivatives) are evaluated at the reversed process $\mathbf{U}^{\t,\cdot,\ast}[\mathbb{I}(\mathfrak{t})]$:
\begin{align}
&\N^{-2}|[\grad^{\mathbb{I}(\mathfrak{t}),-}]^{-1}\mathscr{I}(\t,\x)|^{2} \ \lesssim \ \N^{2}\mathfrak{p}[\mathfrak{t}(\mathrm{in}),\t[\ast]]^{-2}|[\grad^{\mathbb{I}(\mathfrak{t}),-}]^{-1}\grad^{\mathbb{I}(\mathfrak{t}),-}\grad^{\mathbb{I}(\mathfrak{t}),+}\partial_{\mathbf{U}(\x)}\mathfrak{p}[\mathfrak{t}(\mathrm{in}),\t[\ast]]|^{2} \label{eq:kv15II4a}\\
&= \ \N^{2}\mathfrak{p}[\mathfrak{t}(\mathrm{in}),\t[\ast]]^{-2}|\grad^{\mathbb{I}(\mathfrak{t}),+}\partial_{\mathbf{U}(\x)}\mathfrak{p}[\mathfrak{t}(\mathrm{in}),\t[\ast]]|^{2} \ = \ \N^{2}\mathfrak{p}[\mathfrak{t}(\mathrm{in}),\t[\ast]]^{-2}|\mathrm{D}_{\x}\mathfrak{p}[\mathfrak{t}(\mathrm{in}),\t[\ast]]|^{2}. \label{eq:kv15II4b}
\end{align}
Indeed, \eqref{eq:kv15II4b} follows from construction of $\mathrm{D}$-operators in Definition \ref{definition:le5}. We now plug \eqref{eq:kv15II4a}-\eqref{eq:kv15II4b} into \eqref{eq:kv15II3b}. We then pull $\E^{\ast,\mathrm{back}}$ into the time-integral. Ultimately, we claim the following, which we clarify afterwards:
\begin{align}
&\E^{\ast,\mathrm{back}} {\textstyle\int_{\mathfrak{t}(\mathrm{in})}^{\mathfrak{t}(\mathrm{in})+\mathfrak{t}}}{\textstyle\sum_{\x}}\N^{-2}|[\grad^{\mathbb{I}(\mathfrak{t}),-}]^{-1}\mathscr{I}(\t,\x)|^{2}\d\t \nonumber\\
&\lesssim \ \N^{2}\E^{\ast,\mathrm{back}} {\textstyle\int_{\mathfrak{t}(\mathrm{in})}^{\mathfrak{t}(\mathrm{in})+\mathfrak{t}}}{\textstyle\sum_{\x}}\mathfrak{p}[\mathfrak{t}(\mathrm{in}),\t[\ast]]^{-2}|\mathrm{D}_{\x}\mathfrak{p}[\mathfrak{t}(\mathrm{in}),\t[\ast]]|^{2}\d\t \nonumber \\
&= \ \N^{2} {\textstyle\int_{\mathfrak{t}(\mathrm{in})}^{\mathfrak{t}(\mathrm{in})+\mathfrak{t}}}\E^{\ast,\mathrm{back}}{\textstyle\sum_{\x}}\mathfrak{p}[\mathfrak{t}(\mathrm{in}),\t[\ast],\mathbf{U}^{\t,\cdot,\ast}[\mathbb{I}(\mathfrak{t})]]^{-2}|\mathrm{D}_{\x}\mathfrak{p}[\mathfrak{t}(\mathrm{in}),\t[\ast],\mathbf{U}^{\t,\cdot,\ast}[\mathbb{I}(\mathfrak{t})]]|^{2}\d\t \label{eq:kv15II5a} \\
&= \ \N^{2} {\textstyle\int_{\mathfrak{t}(\mathrm{in})}^{\mathfrak{t}(\mathrm{in})+\mathfrak{t}}}\E^{\sigma,\t[\ast],\mathbb{I}(\mathfrak{t})}{\textstyle\sum_{\x}}\mathfrak{p}[\mathfrak{t}(\mathrm{in}),\t[\ast]]^{-1}|\mathrm{D}_{\x}\mathfrak{p}[\mathfrak{t}(\mathrm{in}),\t[\ast]]|^{2}\d\t \label{eq:kv15II5b} \\
&= \ \N^{2} {\textstyle\int_{\mathfrak{t}(\mathrm{in})}^{\mathfrak{t}(\mathrm{in})+\mathfrak{t}}}\E^{\sigma,\t,\mathbb{I}(\mathfrak{t})}{\textstyle\sum_{\x}}\mathfrak{p}[\mathfrak{t}(\mathrm{in}),\t]^{-1}|\mathrm{D}_{\x}\mathfrak{p}[\mathfrak{t}(\mathrm{in}),\t]|^{2}\d\t. \label{eq:kv15II5c}
\end{align}
\eqref{eq:kv15II5a} follows by pulling the expectation through the time-integral in the first line. (We have now emphasized the evaluation at $\mathbf{U}^{\t,\cdot,\ast}[\mathbb{I}(\mathfrak{t})]$.) \eqref{eq:kv15II5b} follows as the expectation in \eqref{eq:kv15II5a} is over $\mathbf{U}^{\t,\cdot,\ast}[\mathbb{I}(\mathfrak{t})]$ at one time. And by construction in Definition \ref{definition:kv12} and Lemma \ref{lemma:kv13}, the law of $\mathbf{U}^{\t,\cdot,\ast}[\mathbb{I}(\mathfrak{t})]$ has density $\mathfrak{p}[\mathfrak{t}(\mathrm{in}),\t[\ast]]$ with respect to $\mathbb{P}^{\sigma,\t[\ast],\mathbb{I}(\mathfrak{t})}$. \eqref{eq:kv15II5c} follows by changing-variables along the linear map $\t[\ast]\mapsto\t$. Now, let $\mathfrak{p}[\mathfrak{t}(\mathrm{in}),\t;\mathcal{E}]$ be the Radon-Nikodym derivative for the law of $\mathbf{U}^{\t,\cdot}[\mathbb{I}(\mathfrak{t}),\mathcal{E}]$ with respect to $\mathbb{P}^{\sigma,\t,\mathbb{I}(\mathfrak{t}),\mathcal{E}}$; see Definition \ref{definition:kv4}. (Here, we assume that the time $\t=\mathfrak{t}(\mathrm{in})$ data of $\mathbf{U}^{\t,\cdot}[\mathbb{I}(\mathfrak{t}),\mathcal{E}]$ equals $\mathbb{P}^{\sigma,\mathfrak{t}(\mathrm{in}),\mathbb{I}(\mathfrak{t}),\mathcal{E}}$. So, $\mathfrak{p}[\mathfrak{t}(\mathrm{in}),\mathfrak{t}(\mathrm{in});\mathcal{E}]\equiv1$. Thus, by Remark \ref{remark:kv9}, $\mathfrak{p}[\mathfrak{t}(\mathrm{in}),\t;\mathcal{E}]$ equals that in Lemma \ref{lemma:kv8}.) Write $\eqref{eq:kv15II5c}=\Phi[1]+\Phi[2]$, where
\begin{align}
\Phi[1] \ &:= \ \N^{2} {\textstyle\int_{\mathfrak{t}(\mathrm{in})}^{\mathfrak{t}(\mathrm{in})+\mathfrak{t}}}\E^{\sigma,\t,\mathbb{I}(\mathfrak{t})}{\textstyle\sum_{\x}}\mathfrak{p}[\mathfrak{t}(\mathrm{in}),\t;\mathcal{E}]^{-1}|\mathrm{D}_{\x}\mathfrak{p}[\mathfrak{t}(\mathrm{in}),\t]|^{2}\d\t \label{eq:kv15II6a}\\
\Phi[2] \ &:= \ \N^{2} {\textstyle\int_{\mathfrak{t}(\mathrm{in})}^{\mathfrak{t}(\mathrm{in})+\mathfrak{t}}}\E^{\sigma,\t,\mathbb{I}(\mathfrak{t})}{\textstyle\sum_{\x}}\{\mathfrak{p}[\mathfrak{t}(\mathrm{in}),\t]^{-1}-\mathfrak{p}[\mathfrak{t}(\mathrm{in}),\t;\mathcal{E}]^{-1}\}|\mathrm{D}_{\x}\mathfrak{p}[\mathfrak{t}(\mathrm{in}),\t]|^{2}\d\t. \label{eq:kv15II6b}
\end{align}
Let $|\mathbf{U}|$ be Euclidean length of $\mathbf{U}\in\R^{\mathbb{I}(\mathfrak{t})}$. We first compute and estimate \eqref{eq:kv15II6b} as follows (with explanation given after):
{\small
\begin{align}
&|\Phi(2)| \nonumber\\
&= \ |\eqref{eq:kv15II6b}| \nonumber\\
&\lesssim \ \N^{2}|\mathbb{I}(\mathfrak{t})|\exp[\mathrm{O}(\N^{\gamma_{\mathrm{av}}})]{\textstyle\int_{\mathfrak{t}(\mathrm{in})}^{\mathfrak{t}(\mathrm{in})+\mathfrak{t}}}\E^{\sigma,\t,\mathbb{I}(\mathfrak{t})}\{1+|\mathbf{U}|^{100}\}|\mathfrak{p}[\mathfrak{t}(\mathrm{in}),\t]^{-1}-\mathfrak{p}[\mathfrak{t}(\mathrm{in}),\t;\mathcal{E}]^{-1}|\d\t \label{eq:kv15II7a}\\
&\lesssim \ \N^{2}|\mathbb{I}(\mathfrak{t})|^{2}\exp[\mathrm{O}(\N^{\gamma_{\mathrm{av}}})]{\textstyle\int_{\mathfrak{t}(\mathrm{in})}^{\mathfrak{t}(\mathrm{in})+\mathfrak{t}}}\{\E^{\sigma,\t,\mathbb{I}(\mathfrak{t})}|\mathfrak{p}[\mathfrak{t}(\mathrm{in}),\t]^{-1}-\mathfrak{p}[\mathfrak{t}(\mathrm{in}),\t;\mathcal{E}]^{-1}|^{2}\}^{\frac12}\d\t \label{eq:kv15II7b}\\
&\lesssim \ \N^{2}|\mathbb{I}(\mathfrak{t})|^{2}\exp[\mathrm{O}(\N^{\gamma_{\mathrm{av}}})]{\textstyle\int_{\mathfrak{t}(\mathrm{in})}^{\mathfrak{t}(\mathrm{in})+\mathfrak{t}}}\{\E^{\sigma,\t,\mathbb{I}(\mathfrak{t})}|\mathfrak{p}[\mathfrak{t}(\mathrm{in}),\t]^{-1}\mathfrak{p}[\mathfrak{t}(\mathrm{in}),\t;\mathcal{E}]^{-1}|^{2}|\mathfrak{p}[\mathfrak{t}(\mathrm{in}),\t]-\mathfrak{p}[\mathfrak{t}(\mathrm{in}),\t;\mathcal{E}]|^{2}\}^{\frac12}\d\t \label{eq:kv15II7c}\\
&\lesssim \ \N^{2}|\mathbb{I}(\mathfrak{t})|^{2}\exp[\mathrm{O}(\N^{\gamma_{\mathrm{av}}})]{\textstyle\int_{\mathfrak{t}(\mathrm{in})}^{\mathfrak{t}(\mathrm{in})+\mathfrak{t}}}\{\E^{\sigma,\t,\mathbb{I}(\mathfrak{t})}|\mathfrak{p}[\mathfrak{t}(\mathrm{in}),\t]-\mathfrak{p}[\mathfrak{t}(\mathrm{in}),\t;\mathcal{E}]|\}^{\frac12}\d\t \label{eq:kv15II7d}\\
&\lesssim \ \N^{2}\mathfrak{t}|\mathbb{I}(\mathfrak{t})|^{2}\exp[\mathrm{O}(\N^{\gamma_{\mathrm{av}}})]\exp[-\N^{\frac15\gamma_{\mathrm{KV}}}] \ \lesssim \ \exp[\mathrm{O}(\N^{\gamma_{\mathrm{av}}})]\exp[-\N^{\frac15\gamma_{\mathrm{KV}}}] \ \lesssim \ \exp[-\N^{\frac16\gamma_{\mathrm{KV}}}]. \label{eq:kv15II7e}
\end{align}
}\eqref{eq:kv15II7a} follows from bounding the sum over $\x\in\mathbb{I}(\mathfrak{t})$ by $|\mathbb{I}(\mathfrak{t})|$ times the supremum, and then by using \eqref{eq:kv10II}. \eqref{eq:kv15II7b} follows by Cauchy-Schwarz. Indeed, because $|\sigma|\lesssim1$ and because $\mathbf{U}(\x)$ is sub-Gaussian for each $\x\in\mathbb{I}(\mathfrak{t})$, we know $\E^{\sigma,\t,\mathbb{I}(\mathfrak{t})}|\mathbf{U}(\x)|^{200}\lesssim1$ for all $\x\in\mathbb{I}(\mathfrak{t})$. Thus, $\E^{\sigma,\t,\mathbb{I}(\mathfrak{t})}|\mathbf{U}|^{200}\lesssim|\mathbb{I}(\mathfrak{t})|$, with the $|\mathbb{I}(\mathfrak{t})|$-factor accounting for all $\x\in\mathbb{I}(\mathfrak{t})$. \eqref{eq:kv15II7c} follows by the bound $|\mathrm{a}^{-1}-\mathrm{b}^{-1}|\lesssim|\mathrm{a}\mathrm{b}|^{-1}|\mathrm{a}-\mathrm{b}|$ with $\mathrm{a}=\mathfrak{p}[\mathfrak{t}(\mathrm{in}),\mathfrak{t}]$ and $\mathrm{b}=\mathfrak{p}[\mathfrak{t}(\mathrm{in}),\mathfrak{t};\mathcal{E}]$. \eqref{eq:kv15II7d} follows from \eqref{eq:kv8I}-\eqref{eq:kv8II}. The first bound in \eqref{eq:kv15II7e} follows via Lemma \ref{lemma:kv11}; see \eqref{eq:kv10I2c}-\eqref{eq:kv10I3}. The rest of \eqref{eq:kv15II7e} follows from $\gamma_{\mathrm{av}}\leq{\mathrm{c}}\gamma_{\mathrm{KV}}$ {for some $\mathrm{c}>0$ small}; see Proposition \ref{prop:kv1}. We now control \eqref{eq:kv15II6a}. To this end, we claim the following, again with explanation given after:
\begin{align}
|\Phi(1)| \ = \ |\eqref{eq:kv15II6a}| \ &\lesssim \ \N^{2}{\textstyle\int_{\mathfrak{t}(\mathrm{in})}^{\mathfrak{t}(\mathrm{in})+\mathfrak{t}}}\E^{\sigma,\t,\mathbb{I}(\mathfrak{t})}{\textstyle\sum_{\x}}|\mathrm{D}_{\x}\mathfrak{p}[\mathfrak{t}(\mathrm{in}),\t]|^{2}\d\t \nonumber\\
&\lesssim \ \N^{2}{\textstyle\int_{\mathfrak{t}(\mathrm{in})}^{\mathfrak{t}(\mathrm{in})+\mathfrak{t}}}\mathfrak{D}_{\mathrm{FI}}^{\sigma,\t,\mathbb{I}(\mathfrak{t})}(\mathfrak{p}[\mathfrak{t}(\mathrm{in}),\t]^{2})\d\t. \label{eq:kv15II8a}
\end{align}
\eqref{eq:kv15II8a} is by \eqref{eq:kv8II} and Definition \ref{definition:le3}. \eqref{eq:kv15II} now follows by \eqref{eq:kv15II3a}-\eqref{eq:kv15II3b}, \eqref{eq:kv15II5a}-\eqref{eq:kv15II5c}, the fact $\eqref{eq:kv15II5c}=\Phi[1]+\Phi[2]$ from right before \eqref{eq:kv15II6a}, \eqref{eq:kv15II7a}-\eqref{eq:kv15II7e}, \eqref{eq:kv15II8a}, and \eqref{eq:kv10I}. This finishes the proof.
\end{proof}
\subsection{Proof of Proposition \ref{prop:kv1}}
We now combine the ingredients in this section in exactly the way we motivated them to derive the estimates \eqref{eq:kv1II} and \eqref{eq:kv1III}. 

We first show \eqref{eq:kv1II}. The proof of \eqref{eq:kv1III} uses almost the exact same idea. For convenience of notation, we will assume $\t=\mathfrak{t}(\mathrm{in})+\mathfrak{t}$. (Of course, this is sufficient; in Proposition \ref{prop:kv1}, for any $\t$, we can always redefine $\mathfrak{t}$ while ensuring conditions of Proposition \ref{prop:kv1} are met so that $\t=\mathfrak{t}(\mathrm{in})+\mathfrak{t}$. For convenience, we also continue writing $\t$ for $\mathfrak{t}(\mathrm{in})+\mathfrak{t}$.)
\subsubsection{Changing measure}
Now, note that for all times $\s\geq\mathfrak{t}(\mathrm{in})$, the process $\mathrm{U}+\mathbf{J}(\s,\inf\mathbb{I}(\mathfrak{t});\mathbb{I}(\mathfrak{t}))$ starting from $\mathbf{J}(\mathfrak{t}(\mathrm{in}),\inf\mathbb{I}(\mathfrak{t});\mathbb{I}(\mathfrak{t}))=0$ is the same as running the $\mathbb{S}(\N)$-valued SDE $\mathbf{J}(\s,\inf\mathbb{I}(\mathfrak{t});\mathbb{I}(\mathfrak{t}))$ with initial data $\mathrm{U}$ at $\s=\mathfrak{t}(\mathrm{in})$. This is because \eqref{eq:hfloc}, to which we apply $\Pi^{\mathbb{S}(\N)}$ to get said $\mathbb{S}(\N)$-valued SDE, has RHS that depends only on $\mathbf{U}^{\t,\cdot}[\mathbb{I}(\mathfrak{t})]$. (So, the evolution of $\mathbf{J}(\s,\inf\mathbb{I}(\mathfrak{t});\mathbb{I}(\mathfrak{t}))$ commutes with additive shifts.) Using this, we claim
\begin{align}
\mathrm{LHS}\eqref{eq:kv1II} \ = \ \N^{-20\gamma_{\mathrm{reg}}}\E\{|\mathscr{A}(\t)|^{2}\mathbf{1}[|\mathscr{A}(\t)|\leq\mathscr{B}]\} \ \lesssim \ \E^{\mathrm{path},\to}\{|\mathscr{A}(\t)|^{2}\mathbf{1}[|\mathscr{A}(\t)|\leq\mathscr{B}]\}. \label{eq:kv1II1}
\end{align}
In $\mathrm{RHS}\eqref{eq:kv1II1}$, $\E^{\mathrm{path},\to}$ denotes an expectation with respect to the law of the \emph{forward} joint process $\s\mapsto(\mathbf{J}(\s,\inf\mathbb{I}(\mathfrak{t});\mathbb{I}(\mathfrak{t})),\mathbf{U}^{\s,\cdot}[\mathbb{I}(\mathfrak{t})])$ for $\s\geq\mathfrak{t}(\mathrm{in})$ with initial data at $\s=\mathfrak{t}(\mathrm{in})$ distributed as $\mathbb{P}^{\mathrm{Leb},\sigma,\mathfrak{t}(\mathrm{in}),\mathbb{I}(\mathfrak{t})}$ in Definition \ref{definition:le5}. (Also, see \eqref{eq:kv1I} for $\mathscr{A}(\t)$.) \eqref{eq:kv1II1} follows via changing measure $\mathrm{Unif}[-1,1]\otimes\mathbb{P}^{\sigma,\mathfrak{t}(\mathrm{in}),\mathbb{I}(\mathfrak{t})}\mapsto\mathrm{Leb}(\mathbb{S}(\N))\otimes\mathbb{P}^{\sigma,\mathfrak{t}(\mathrm{in}),\mathbb{I}(\mathfrak{t})}$. Indeed, the change-of-measure factor is just the Radon-Nikodym derivative of $\mathrm{Unif}[-1,1]$ with respect to $\mathrm{Leb}(\mathbb{S}(\N))=\mathrm{Unif}[-\N^{20\gamma_{\mathrm{reg}}},\N^{20\gamma_{\mathrm{reg}}}]$. So, it is $\lesssim\N^{20\gamma_{\mathrm{reg}}}$. {\emph{For the rest of this proof}}, we drop $\mathrm{U}$ from $\mathscr{A}(\t)$ in \eqref{eq:kv1I}. Indeed, $\mathrm{Unif}[-1,1]$ convolved with $\mathrm{Leb}(\mathbb{S}(\N))$ is just $\mathrm{Leb}(\mathbb{S}(\N))$.
\subsubsection{Forward-Backward decomposition}
The point of this step is to write $\mathscr{A}(\t)$ in \eqref{eq:kv1I} (but without $\mathrm{U}$ therein!) in terms of the joint process $\s\mapsto(\mathbf{J}(\s,\inf\mathbb{I}(\mathfrak{t});\mathbb{I}(\mathfrak{t})),\mathbf{U}^{\s,\cdot}[\mathbb{I}(\mathfrak{t})])$ and its time-reversal in Definition \ref{definition:kv12}. This is what is usually done when analyzing fluctuations in time; see Section 4 of \cite{CLO}, for example. Let $\mathsf{F}(\s,\cdot,\cdot):\R\times\mathbb{H}^{\sigma,\mathbb{I}(\mathfrak{t})}\to\R$ solve the resolvent equation
\begin{align}
\{[1-\mathscr{L}^{\mathrm{tot},\mathrm{sym}}(\s,\mathbb{I}(\mathfrak{t}))]\mathsf{F}\}(\s,\mathrm{a},\mathbf{U}) \ = \ \varphi(\mathrm{a})\mathfrak{a}(\s,\mathbf{U}). \label{eq:kv1II2}
\end{align}
We recall that $\mathscr{L}^{\mathrm{tot},\mathrm{sym}}(\s,\mathbb{I}(\mathfrak{t}))$ denotes the symmetric part of the generator $\mathscr{L}^{\mathrm{tot}}(\s,\mathbb{I}(\mathfrak{t}))$ for the process $\s\mapsto(\mathbf{J}(\s,\inf\mathbb{I}(\mathfrak{t});\mathbb{I}(\mathfrak{t})),\mathbf{U}^{\s,\cdot}[\mathbb{I}(\mathfrak{t})])$; see Lemma \ref{lemma:le7}. (The symmetric part is with respect to $\mathbb{P}^{\mathrm{Leb},\sigma,\s,\mathbb{I}(\mathfrak{t})}$ in Definition \ref{definition:le5}.) If $\mathscr{L}^{\mathrm{tot}}(\s,\mathbb{I}(\mathfrak{t}))^{\ast}$ denotes the adjoint of $\mathscr{L}^{\mathrm{tot}}(\s,\mathbb{I}(\mathfrak{t}))$, then
\begin{align}
[1-\mathscr{L}^{\mathrm{tot},\mathrm{sym}}(\s,\mathbb{I}(\mathfrak{t}))]\mathsf{F} \ &= \ \mathsf{F}-\tfrac12[\mathscr{L}^{\mathrm{tot}}(\s,\mathbb{I}(\mathfrak{t}))+\mathscr{L}^{\mathrm{tot}}(\s,\mathbb{I}(\mathfrak{t}))^{\ast}]\mathsf{F} \ = \ \Phi^{(1)}+\Phi^{(2)}\label{eq:kv1II3a}\\
&:= \ \{\mathsf{F}-\tfrac12[\partial_{\s}+\mathscr{L}^{\mathrm{tot}}(\s,\mathbb{I}(\mathfrak{t}))]\mathsf{F}]\}+\{-\tfrac12[-\partial_{\s}+\mathscr{L}^{\mathrm{tot}}(\s,\mathbb{I}(\mathfrak{t}))^{\ast}]\mathsf{F}\}. \label{eq:kv1II3b}
\end{align}
(Technically, we cannot differentiate in $\s$ as $\mathsf{F}$ was not assumed to be time-differentiable at a finite set of times. The workaround for this is a standard density/mollification argument for $\mathfrak{a}$ by smooth approximations that converge locally uniformly away from a finite set of points. This works since all our bounds, including for $\mathrm{LHS}\eqref{eq:kv1II}$ and $\mathrm{LHS}\eqref{eq:kv1III}$, are with respect to time-integrated norms of $\mathfrak{a}$, and these are continuous with respect to locally uniform convergence outside a measure zero set. Indeed, the upper bounds $\mathrm{RHS}\eqref{eq:kv1II}$ and $\mathrm{RHS}\eqref{eq:kv1III}$ themselves come by estimating time-integrated norms by their suprema.) Combine \eqref{eq:kv1II2} and \eqref{eq:kv1II3a}-\eqref{eq:kv1II3b}. Recall $\mathscr{A}(\t)$ in \eqref{eq:kv1I} (now without $\mathrm{U}$ therein). We have
\begin{align}
\mathscr{A}(\t) \ = \ \mathscr{A}^{(1)}(\t)+\mathscr{A}^{(2)}(\t) \ =: \ \mathscr{A}^{(1)}+\mathscr{A}^{(2)}, \label{eq:kv1II4a}
\end{align}
where $\mathscr{A}^{(\mathrm{i})}$ are integrals of $\Phi^{(\mathrm{i})}$ from \eqref{eq:kv1II3a}-\eqref{eq:kv1II3b} but with additional ``boundary terms" (in the sense of Ito formula):
\begin{align}
&\mathscr{A}^{(1)}(\t) \nonumber\\
&:= \tfrac{1}{2[\t-\mathfrak{t}(\mathrm{in})]}\{\mathsf{F}(\t,\mathbf{J}(\t,\inf\mathbb{I}(\mathfrak{t});\mathbb{I}(\mathfrak{t})),\mathbf{U}^{\t,\cdot}[\mathbb{I}(\mathfrak{t})])-\mathsf{F}(\mathfrak{t}(\mathrm{in}),\mathbf{J}(\mathfrak{t}(\mathrm{in}),\inf\mathbb{I}(\mathfrak{t});\mathbb{I}(\mathfrak{t})),\mathbf{U}^{\mathfrak{t}(\mathrm{in}),\cdot}[\mathbb{I}(\mathfrak{t})])\}\label{eq:kv1II5a}\\
&+ \ \tfrac{1}{[\t-\mathfrak{t}(\mathrm{in})]}{\textstyle\int_{\mathfrak{t}(\mathrm{in})}^{\t}}\Phi^{(1)}(\s,\mathbf{J}(\s,\inf\mathbb{I}(\mathfrak{t});\mathbb{I}(\mathfrak{t})),\mathbf{U}^{\s,\cdot}[\mathbb{I}(\mathfrak{t})])\d\s\label{eq:kv1II5b}\\
&\mathscr{A}^{(2)}(\t) \nonumber\\
&:= -\tfrac{1}{2[\t-\mathfrak{t}(\mathrm{in})]}\{\mathsf{F}(\t,\mathbf{J}(\t,\inf\mathbb{I}(\mathfrak{t});\mathbb{I}(\mathfrak{t})),\mathbf{U}^{\t,\cdot}[\mathbb{I}(\mathfrak{t})])-\mathsf{F}(\mathfrak{t}(\mathrm{in}),\mathbf{J}(\mathfrak{t}(\mathrm{in}),\inf\mathbb{I}(\mathfrak{t});\mathbb{I}(\mathfrak{t})),\mathbf{U}^{\mathfrak{t}(\mathrm{in}),\cdot}[\mathbb{I}(\mathfrak{t})])\}\label{eq:kv1II5c}\\
&+ \ \tfrac{1}{[\t-\mathfrak{t}(\mathrm{in})]}{\textstyle\int_{\mathfrak{t}(\mathrm{in})}^{\t}}\Phi^{(2)}(\s,\mathbf{J}(\s,\inf\mathbb{I}(\mathfrak{t});\mathbb{I}(\mathfrak{t})),\mathbf{U}^{\s,\cdot}[\mathbb{I}(\mathfrak{t})])\d\s.\label{eq:kv1II5d}
\end{align}
(Briefly, $\mathrm{RHS}\eqref{eq:kv1II5a}+\mathrm{RHS}\eqref{eq:kv1II5c}=0$, and $\eqref{eq:kv1II5b}+\eqref{eq:kv1II5d}=\mathscr{A}(\t)$ via \eqref{eq:kv1II3a}-\eqref{eq:kv1II3b}.) By \eqref{eq:kv1II4a}, we deduce
\begin{align}
\mathrm{RHS}\eqref{eq:kv1II1} \ = \ \E^{\mathrm{path},\to}\{|\mathscr{A}^{(1)}(\t)+\mathscr{A}^{(2)}(\t)|^{2}\mathbf{1}[|\mathscr{A}^{(1)}(\t)+\mathscr{A}^{(2)}(\t)|\lesssim\mathscr{B}]\}. \label{eq:kv1II6}
\end{align}
We conclude this step by rewriting \eqref{eq:kv1II5c}-\eqref{eq:kv1II5d}. We claim the following with notation and clarification explained after:
\begin{align}
\mathscr{A}^{(2)}(\t) \ &= \ \tfrac{1}{2[\t-\mathfrak{t}(\mathrm{in})]}\mathsf{F}([\mathfrak{t}(\mathrm{in})+\mathfrak{t}][\ast],\mathbf{J}([\mathfrak{t}(\mathrm{in})+\mathfrak{t}][\ast],\inf\mathbb{I}(\mathfrak{t});\mathbb{I}(\mathfrak{t})),\mathbf{U}^{[\mathfrak{t}(\mathrm{in})+\mathfrak{t}][\ast],\cdot}[\mathbb{I}(\mathfrak{t})])\label{eq:kv1II7a}\\
&- \ \tfrac{1}{2[\t-\mathfrak{t}(\mathrm{in})]}\mathsf{F}(\mathfrak{t}(\mathrm{in})[\ast],\mathbf{J}(\mathfrak{t}(\mathrm{in})[\ast],\inf\mathbb{I}(\mathfrak{t});\mathbb{I}(\mathfrak{t})),\mathbf{U}^{\mathfrak{t}(\mathrm{in})[\ast],\cdot}[\mathbb{I}(\mathfrak{t})])\label{eq:kv1II7b}\\
&+ \ \tfrac{1}{[\t-\mathfrak{t}(\mathrm{in})]}{\textstyle\int_{\mathfrak{t}(\mathrm{in})}^{\mathfrak{t}(\mathrm{in})+\mathfrak{t}}}\Phi^{(2)}(\s[\ast],\mathbf{J}(\s[\ast],\inf\mathbb{I}(\mathfrak{t});\mathbb{I}(\mathfrak{t})),\mathbf{U}^{\s[\ast],\cdot}[\mathbb{I}(\mathfrak{t})])\d\s\label{eq:kv1II7c}\\
&= \ \tfrac{1}{2[\t-\mathfrak{t}(\mathrm{in})]}\mathsf{F}([\mathfrak{t}(\mathrm{in})+\mathfrak{t}][\ast],\mathbf{J}^{\ast}(\mathfrak{t}(\mathrm{in})+\mathfrak{t},\inf\mathbb{I}(\mathfrak{t});\mathbb{I}(\mathfrak{t})),\mathbf{U}^{\mathfrak{t}(\mathrm{in})+\mathfrak{t},\cdot,\ast}[\mathbb{I}(\mathfrak{t})])\label{eq:kv1II7d}\\
&- \ \tfrac{1}{2[\t-\mathfrak{t}(\mathrm{in})]}\mathsf{F}(\mathfrak{t}(\mathrm{in})[\ast],\mathbf{J}^{\ast}(\mathfrak{t}(\mathrm{in}),\inf\mathbb{I}(\mathfrak{t});\mathbb{I}(\mathfrak{t})),\mathbf{U}^{\mathfrak{t}(\mathrm{in}),\cdot,\ast}[\mathbb{I}(\mathfrak{t})])\label{eq:kv1II7e}\\
&+ \ \tfrac{1}{[\t-\mathfrak{t}(\mathrm{in})]}{\textstyle\int_{\mathfrak{t}(\mathrm{in})}^{\mathfrak{t}(\mathrm{in})+\mathfrak{t}}}\Phi^{(2)}(\s[\ast],\mathbf{J}^{\ast}(\s,\inf\mathbb{I}(\mathfrak{t});\mathbb{I}(\mathfrak{t})),\mathbf{U}^{\s,\cdot,\ast}[\mathbb{I}(\mathfrak{t})])\d\s.\label{eq:kv1II7f}
\end{align}
Recall {the} $\ast$ notation from Definition \ref{definition:kv12}; it just reverses time. Let us explain \eqref{eq:kv1II7a}-\eqref{eq:kv1II7c}. Recall $\t=\mathfrak{t}(\mathrm{in})+\mathfrak{t}$. In \eqref{eq:kv1II5c}, all we do is write $\t=\mathfrak{t}(\mathrm{in})[\ast]$ and $\mathfrak{t}(\mathrm{in})=[\mathfrak{t}(\mathrm{in})+\mathfrak{t}][\ast]$. The first term in $\mathrm{RHS}\eqref{eq:kv1II5c}$
 becomes \eqref{eq:kv1II7b}, and the last term in $\mathrm{RHS}\eqref{eq:kv1II5c}$ becomes $\mathrm{RHS}\eqref{eq:kv1II7a}$. Now, in \eqref{eq:kv1II5d}, all we do is change-of-variables $\s\mapsto\s[\ast]$ in the time-integral. This is an affine bijection $[\mathfrak{t}(\mathrm{in}),\t]\to[\mathfrak{t}(\mathrm{in}),\t]$ since $\t=\mathfrak{t}(\mathrm{in})+\mathfrak{t}$; see Definition \ref{definition:kv12}. So the change-of-variables factor is 1, from which we have $\eqref{eq:kv1II5d}=\eqref{eq:kv1II7c}$. \eqref{eq:kv1II7d}-\eqref{eq:kv1II7f} follows from \eqref{eq:kv1II7a}-\eqref{eq:kv1II7c}, since putting $\ast$ in the superscript just means evaluate at $\s[\ast]$ instead of $\s$; see Definition \ref{definition:kv12}. We have now written $\mathscr{A}(\t)$ in terms of just $\s\mapsto(\mathbf{J}(\s,\inf\mathbb{I}(\mathfrak{t});\mathbb{I}(\mathfrak{t})),\mathbf{U}^{\s,\cdot}[\mathbb{I}(\mathfrak{t})])$ and its time-reversal; see \eqref{eq:kv1II5a}-\eqref{eq:kv1II5b}, \eqref{eq:kv1II6}, and \eqref{eq:kv1II7a}-\eqref{eq:kv1II7f}. (We clarify that even with this representation of $\mathscr{A}^{(2)}(\t)$ in terms of the reversed process, it is still a function of the forwards-time process, so $\E^{\mathrm{path},\to}$ of it makes sense.)
\subsubsection{Girsanov step}
As discussed in Section \ref{section:sqle} and alluded to via Lemma \ref{lemma:kv15}, working with time-reversal in Definition \ref{definition:kv12} itself is not the best idea. We use Lemma \ref{lemma:kv15} to replace it with the $\sim$ processes in Definition \ref{definition:kv14}. In order to state this precisely, let us first introduce the following version of \eqref{eq:kv1II7d}-\eqref{eq:kv1II7f} but for $\sim$ processes instead of $\ast$ processes:
\begin{align}
\mathscr{A}^{(2),\sim}(\t):= \mathscr{A}^{(2),\sim} &:=\tfrac{1}{2[\t-\mathfrak{t}(\mathrm{in})]}\mathsf{F}([\mathfrak{t}(\mathrm{in})+\mathfrak{t}][\ast],\mathbf{J}^{\sim}(\mathfrak{t}(\mathrm{in})+\mathfrak{t},\inf\mathbb{I}(\mathfrak{t});\mathbb{I}(\mathfrak{t})),\mathbf{U}^{\mathfrak{t}(\mathrm{in})+\mathfrak{t},\cdot,\sim}[\mathbb{I}(\mathfrak{t})])\label{eq:kv1II8a}\\
&- \ \tfrac{1}{2[\t-\mathfrak{t}(\mathrm{in})]}\mathsf{F}(\mathfrak{t}(\mathrm{in})[\ast],\mathbf{J}^{\sim}(\mathfrak{t}(\mathrm{in}),\inf\mathbb{I}(\mathfrak{t});\mathbb{I}(\mathfrak{t})),\mathbf{U}^{\mathfrak{t}(\mathrm{in}),\cdot,\sim}[\mathbb{I}(\mathfrak{t})])\label{eq:kv1II8b}\\
&+ \ \tfrac{1}{[\t-\mathfrak{t}(\mathrm{in})]}{\textstyle\int_{\mathfrak{t}(\mathrm{in})}^{\mathfrak{t}(\mathrm{in})+\mathfrak{t}}}\Phi^{(2)}(\s[\ast],\mathbf{J}^{\sim}(\s,\inf\mathbb{I}(\mathfrak{t});\mathbb{I}(\mathfrak{t})),\mathbf{U}^{\s,\cdot,\sim}[\mathbb{I}(\mathfrak{t})])\d\s.\label{eq:kv1II8c}
\end{align}
We now claim the following estimate {for any large but finite $\mathrm{D}>0$}, which we explain afterwards:
\begin{align}
\mathrm{RHS}\eqref{eq:kv1II6} &\lesssim \N^{\gamma_{\mathrm{KV}}}\mathfrak{t}|\mathbb{I}(\mathfrak{t})|^{\frac12}\mathscr{B}^{2}+\N^{-{\mathrm{D}}}\mathscr{B}^{2}\nonumber\\
&+\mathscr{B}^{2}\log\E^{\mathrm{path},\to,\leftarrow}\exp(\mathscr{B}^{-2}|\mathscr{A}^{(1)}+\mathscr{A}^{(2),\sim}|^{2}\mathbf{1}[|\mathscr{A}^{(1)}+\mathscr{A}^{(2),\sim}|\lesssim\mathscr{B}]). \nonumber
\end{align}
On the RHS, $\E^{\mathrm{path},\to,\leftarrow}$ is with respect to the law of $\s\mapsto\mathbf{FB}^{\sim}[\s;\mathbb{I}(\mathfrak{t})]$ from Lemma \ref{lemma:kv15} for $\s\in\mathfrak{t}(\mathrm{in})+[0,\mathfrak{t}]$. (In particular, the law of $(\mathbf{J}^{\sim}(\mathfrak{t}(\mathrm{in}),\inf\mathbb{I}(\mathfrak{t});\mathbb{I}(\mathfrak{t})),\mathbf{U}^{\mathfrak{t}(\mathrm{in}),\cdot,\sim}[\mathbb{I}(\mathfrak{t})])$ equals the time $\mathfrak{t}(\mathrm{in})+\mathfrak{t}$ law of $(\mathbf{J}(\s,\inf\mathbb{I}(\mathfrak{t});\mathbb{I}(\mathfrak{t})),\mathbf{U}^{\s,\cdot}[\mathbb{I}(\mathfrak{t})])$.) Let us now explain this bound. First, we use duality between relative entropy and exponential moments. (This is what we used to get \eqref{eq:le4II}, with $\kappa$ there equal to $\mathscr{B}^{-2}$ here. The point is that this relative entropy inequality, which can be found in Appendix 1.8 of \cite{KL}, holds for any Polish space.) This would give the previous display if we replace $\N^{\gamma_{\mathrm{KV}}}\mathfrak{t}|\mathbb{I}(\mathfrak{t})|^{1/2}\mathscr{B}^{2}+\N^{-{\mathrm{D}}}\mathscr{B}^{2}$ by $\mathscr{B}^{2}$ times relative entropy of $\mathbf{FB}^{\ast}[\cdot;\mathbb{I}(\mathfrak{t})]$ with respect to $\mathbf{FB}^{\sim}[\cdot;\mathbb{I}(\mathfrak{t})]$. (Here, we use notation of Lemma \ref{lemma:kv15}.) However, $\mathscr{B}^{2}$ times said relative entropy is $\lesssim\N^{\gamma_{\mathrm{KV}}}\mathfrak{t}|\mathbb{I}(\mathfrak{t})|^{1/2}\mathscr{B}^{2}+\N^{-{\mathrm{D}}}\mathscr{B}^{2}$; see Lemma \ref{lemma:kv15}. The previous bound follows. Now, observe that the term in $\exp$ on the RHS of the above display is $\mathrm{O}(1)$. Indeed, it has the form $\mathscr{B}^{-2}|\mathfrak{X}|^{2}\mathbf{1}[|\mathfrak{X}|\lesssim\mathscr{B}]$. We note $\exp[\mathfrak{x}]\leq1+\mathrm{O}(|\mathfrak{x}|)$ for $|\mathfrak{x}|\lesssim1$; this is by smoothness of $\exp$. Also, by concavity, we know $\log[1+|\mathfrak{x}|]\leq|\mathfrak{x}|$. So, the last term in the above display is
\begin{align}
&\lesssim \ \mathscr{B}^{2}\log\E^{\mathrm{path},\to,\leftarrow}\{1+\mathrm{O}(\mathscr{B}^{-2}|\mathscr{A}^{(1)}+\mathscr{A}^{(2),\sim}|^{2}\mathbf{1}[|\mathscr{A}^{(1)}+\mathscr{A}^{(2)}|\lesssim\mathscr{B}])\}\label{eq:kv1II9a}\\
&= \ \mathscr{B}^{2}\log\{1+\E^{\mathrm{path},\to,\leftarrow}\mathrm{O}(\mathscr{B}^{-2}|\mathscr{A}^{(1)}+\mathscr{A}^{(2),\sim}|^{2}\mathbf{1}[|\mathscr{A}^{(1)}+\mathscr{A}^{(2)}|\lesssim\mathscr{B}])\}\label{eq:kv1II9b}\\
&\lesssim \ \mathscr{B}^{2}\E^{\mathrm{path},\to,\leftarrow}\mathrm{O}(\mathscr{B}^{-2}|\mathscr{A}^{(1)}+\mathscr{A}^{(2),\sim}|^{2}\mathbf{1}[|\mathscr{A}^{(1)}+\mathscr{A}^{(2)}|\lesssim\mathscr{B}])\label{eq:kv1II9c}\\
&\lesssim \ \E^{\mathrm{path},\to,\leftarrow}|\mathscr{A}^{(1)}+\mathscr{A}^{(2),\sim}|^{2}\mathbf{1}[|\mathscr{A}^{(1)}+\mathscr{A}^{(2)}|\lesssim\mathscr{B}]. \label{eq:kv1II9d}
\end{align}
(\eqref{eq:kv1II9b} follows by linearity of expectation and $\E1=1$.) Combining the previous two displays gives
\begin{align}
\mathrm{RHS}\eqref{eq:kv1II6} \ &\lesssim \ \N^{\gamma_{\mathrm{KV}}}\mathfrak{t}|\mathbb{I}(\mathfrak{t})|^{\frac12}\mathscr{B}^{2}+\N^{-{\mathrm{D}}}\mathscr{B}^{2}+\E^{\mathrm{path},\to,\leftarrow}|\mathscr{A}^{(1)}+\mathscr{A}^{(2),\sim}|^{2}.\label{eq:kv1II10}
\end{align}
%
\subsubsection{A decoupling step}
We now separate $\mathscr{A}^{(1)}$ and $\mathscr{A}^{(2),\sim}$ in $\mathrm{RHS}\eqref{eq:kv1II10}$. By the inequality $|\mathfrak{x}+\mathfrak{y}|^{2}\lesssim|\mathfrak{x}|^{2}+|\mathfrak{y}|^{2}$, we know
\begin{align}
\E^{\mathrm{path},\to,\leftarrow}|\mathscr{A}^{(1)}+\mathscr{A}^{(2),\sim}|^{2} \ \lesssim \ {\E^{\mathrm{path},\to,\leftarrow}|\mathscr{A}^{(1)}|^{2}}+\E^{\mathrm{path},\to,\leftarrow}|\mathscr{A}^{(2),\sim}|^{2}. \label{eq:kv1II11}
\end{align}
In this step, we leave the first term in $\mathrm{RHS}\eqref{eq:kv1II11}$ alone. Observe the second term in $\mathrm{RHS}\eqref{eq:kv1II1}$ is now just expectation with respect to the law of the $\sim$ processes from Definition \ref{definition:kv14}. The (initial) time $\mathfrak{t}(\mathrm{in})$ data of said $\sim$ process is given by the law of $(\mathbf{J}(\mathfrak{t}(\mathrm{in})+\mathfrak{t},\inf\mathbb{I}(\mathfrak{t});\mathbb{I}(\mathfrak{t})),\mathbf{U}^{\mathfrak{t}(\mathrm{in})+\mathfrak{t},\cdot}[\mathbb{I}(\mathfrak{t})])$, assuming the law $(\mathbf{J}(\mathfrak{t}(\mathrm{in}),\inf\mathbb{I}(\mathfrak{t});\mathbb{I}(\mathfrak{t})),\mathbf{U}^{\mathfrak{t}(\mathrm{in}),\cdot}[\mathbb{I}(\mathfrak{t})])$ is $\sim\mathbb{P}^{\mathrm{Leb},\sigma,\mathfrak{t}(\mathrm{in}),\mathbb{I}(\mathfrak{t})}$. This is not $\mathbb{P}^{\mathrm{Leb},\sigma,\mathfrak{t}(\mathrm{in})+\mathfrak{t},\mathbb{I}(\mathfrak{t})}$. We now take care of this. With explanation given afterwards, we write
\begin{align}
\E^{\mathrm{path},\to,\leftarrow}|\mathscr{A}^{(2),\sim}|^{2} \ = \ \E^{\mathrm{path},\leftarrow,1}|\mathscr{A}^{(2),\sim}|^{2}+\E^{\mathrm{path},\to,\leftarrow}|\mathscr{A}^{(2),\sim}|^{2}- \E^{\mathrm{path},\leftarrow,1}|\mathscr{A}^{(2),\sim}|^{2}. \label{eq:kv1II12}
\end{align}
In \eqref{eq:kv1II12}, $\E^{\mathrm{path},\leftarrow,1}$ is with respect to the law of $\s\mapsto(\mathbf{J}^{\sim}(\s,\inf\mathbb{I}(\mathfrak{t});\mathbb{I}(\mathfrak{t})),\mathbf{U}^{\s,\cdot,\sim}[\mathbb{I}(\mathfrak{t})])$ given initial data at time $\s=\mathfrak{t}(\mathrm{in})$ distributed as the law of $(\mathbf{J}(\s;\mathbb{I}(\mathfrak{t}),\mathcal{E}),\mathbf{U}^{\s,\cdot}[\mathbb{I}(\mathfrak{t}),\mathcal{E}])$ from Definition \ref{definition:kv4} at $\s=\mathfrak{t}(\mathrm{in})+\mathfrak{t}$. Therefore, we are changing {the} initial data of the process $\s\mapsto(\mathbf{J}^{\sim}(\s,\inf\mathbb{I}(\mathfrak{t});\mathbb{I}(\mathfrak{t})),\mathbf{U}^{\s,\cdot,\sim}[\mathbb{I}(\mathfrak{t})])$. Recall {the} notation of Lemma \ref{lemma:kv8}. We study the first term in $\mathrm{RHS}\eqref{eq:kv1II12}$:
\begin{align}
\E^{\mathrm{path},\leftarrow,1}|\mathscr{A}^{(2),\sim}|^{2} \ &= \ \E^{\mathrm{Leb},\sigma,\mathfrak{t}(\mathrm{in})+\mathfrak{t},\mathbb{I}(\mathfrak{t}),\mathcal{E}}\mathfrak{p}[\mathfrak{t}(\mathrm{in}),\mathfrak{t}(\mathrm{in})+\mathfrak{t};\mathcal{E}]\E^{\mathrm{path},\sim}|\mathscr{A}^{(2),\sim}|^{2} \label{eq:kv1II13a}\\
&\lesssim \ \E^{\mathrm{Leb},\sigma,\mathfrak{t}(\mathrm{in})+\mathfrak{t},\mathbb{I}(\mathfrak{t}),\mathcal{E}}\E^{\mathrm{path},\sim}|\mathscr{A}^{(2),\sim}|^{2}\label{eq:kv1II13b}\\
&\lesssim \ \E^{\mathrm{Leb},\sigma,\mathfrak{t}(\mathrm{in})+\mathfrak{t},\mathbb{I}(\mathfrak{t})}\E^{\mathrm{path},\sim}|\mathscr{A}^{(2),\sim}|^{2}\label{eq:kv1II13c}\\
&=: \ \E^{\mathrm{path},\leftarrow}|\mathscr{A}^{(2),\sim}|^{2}. \label{eq:kv1II13d}
\end{align}
Let us now explain \eqref{eq:kv1II13a}. $\E^{\mathrm{path},\sim}$ denotes an expectation with respect to the law of the \emph{process} $\s\mapsto(\mathbf{J}^{\sim}(\s,\inf\mathbb{I}(\mathfrak{t});\mathbb{I}(\mathfrak{t})),\mathbf{U}^{\s,\cdot,\sim}[\mathbb{I}(\mathfrak{t})])$ for all $\s\in\mathfrak{t}(\mathrm{in})+[0,\mathfrak{t}]$. The initial data at $\s=\mathfrak{t}(\mathrm{in})$ is sampled in the outer expectation $\E^{\mathrm{Leb},\sigma,\mathfrak{t}(\mathrm{in})+\mathfrak{t},\mathbb{I}(\mathfrak{t}),\mathcal{E}}\mathfrak{p}[\mathfrak{t}(\mathrm{in}),\mathfrak{t}(\mathrm{in})+\mathfrak{t};\mathcal{E}]$ over $\mathbf{U}\in\mathbb{H}^{\sigma,\mathbb{I}(\mathfrak{t})}$. In particular, \eqref{eq:kv1II13a} is just the factorization of a path-space expectation in terms of expectation with respect to the dynamic and expectation with respect to its initial data. (It follows by law of total expectation.) \eqref{eq:kv1II13b} follows by \eqref{eq:kv8II}. \eqref{eq:kv1II13c} follows by \eqref{eq:kv1II13b} if we prove the change-of-measure factor to go from $\E^{\mathrm{Leb},\sigma,\mathfrak{t}(\mathrm{in})+\mathfrak{t},\mathbb{I}(\mathfrak{t}),\mathcal{E}}\mapsto\E^{\mathrm{Leb},\sigma,\mathfrak{t}(\mathrm{in})+\mathfrak{t},\mathbb{I}(\mathfrak{t})}$ is $\lesssim1$. To this end, note {that} the former $\E^{\mathrm{Leb},\sigma,\mathfrak{t}(\mathrm{in})+\mathfrak{t},\mathbb{I}(\mathfrak{t}),\mathcal{E}}$ is just the latter $\E^{\mathrm{Leb},\sigma,\mathfrak{t}(\mathrm{in})+\mathfrak{t},\mathbb{I}(\mathfrak{t})}$ conditioned on a set whose complement has probability $\lesssim\N^{-100}$; see Definition \ref{definition:kv4} and Lemma \ref{lemma:kv6}. Thus, \eqref{eq:kv1II13c} follows. \eqref{eq:kv1II13d} follows by definition. We note {that} \eqref{eq:kv1II13d} is just {the} expectation with respect to the law of $\s\mapsto(\mathbf{J}^{\sim}(\s,\inf\mathbb{I}(\mathfrak{t});\mathbb{I}(\mathfrak{t})),\mathbf{U}^{\s,\cdot,\sim}[\mathbb{I}(\mathfrak{t})])$ for time $\s\in\mathfrak{t}(\mathrm{in})+[0,\mathfrak{t}]$. This is assuming its initial data at $\s=\mathfrak{t}(\mathrm{in})$ is $\mathbb{P}^{\mathrm{Leb},\sigma,\mathfrak{t}(\mathrm{in})+\mathfrak{t},\mathbb{I}(\mathfrak{t})}$ from Definition \ref{definition:le5}. Let us now study the difference in $\mathrm{RHS}\eqref{eq:kv1II12}$. We proceed like \eqref{eq:kv1II13a}-\eqref{eq:kv1II13d}, by factorizing path-space expectations as expectations with respect to the dynamics and with respect to initial data. First, recall notation of Lemma \ref{lemma:kv8}. We write the following, which we explain and clarify after:
{\small
\begin{align}
&\E^{\mathrm{path},\to,\leftarrow}|\mathscr{A}^{(2),\sim}|^{2}- \E^{\mathrm{path},\leftarrow,1}|\mathscr{A}^{(2),\sim}|^{2} \label{eq:kv1II14a}\\
&= \ \E^{\mathrm{Leb},\sigma,\mathfrak{t}(\mathrm{in})+\mathfrak{t},\mathbb{I}(\mathfrak{t})}\{\mathfrak{p}[\mathfrak{t}(\mathrm{in}),\mathfrak{t}(\mathrm{in})+\mathfrak{t}]-\mathfrak{p}[\mathfrak{t}(\mathrm{in}),\mathfrak{t}(\mathrm{in})+\mathfrak{t};\mathcal{E}]\times p\}\E^{\mathrm{path},\sim}|\mathscr{A}^{(2),\sim}|^{2}\label{eq:kv1II14b}\\
&\lesssim \ [\E^{\mathrm{Leb},\sigma,\mathfrak{t}(\mathrm{in})+\mathfrak{t},\mathbb{I}(\mathfrak{t})}\{\mathfrak{p}[\mathfrak{t}(\mathrm{in}),\mathfrak{t}(\mathrm{in})+\mathfrak{t}]-\mathfrak{p}[\mathfrak{t}(\mathrm{in}),\mathfrak{t}(\mathrm{in})+\mathfrak{t};\mathcal{E}]\times p\}^{2}]^{\frac12}\label{eq:kv1II14c}\\
&\quad\quad\times[\E^{\mathrm{Leb},\sigma,\mathfrak{t}(\mathrm{in})+\mathfrak{t},\mathbb{I}(\mathfrak{t})}\E^{\mathrm{path},\sim}|\mathscr{A}^{(2),\sim}|^{4}]^{\frac12}\nonumber\\
&\lesssim \ \exp[\mathrm{O}(\N^{\gamma_{\mathrm{av}}})][\E^{\mathrm{Leb},\sigma,\mathfrak{t}(\mathrm{in})+\mathfrak{t},\mathbb{I}(\mathfrak{t})}|\mathfrak{p}[\mathfrak{t}(\mathrm{in}),\mathfrak{t}(\mathrm{in})+\mathfrak{t}]-\mathfrak{p}[\mathfrak{t}(\mathrm{in}),\mathfrak{t}(\mathrm{in})+\mathfrak{t};\mathcal{E}]\times p|]^{\frac12}[\E^{\mathrm{path},\leftarrow}|\mathscr{A}^{(2),\sim}|^{4}]^{\frac12}\label{eq:kv1II14d}\\
&\lesssim \ \exp[\mathrm{O}(\N^{\gamma_{\mathrm{av}}})]\exp[-\N^{\frac19\gamma_{\mathrm{KV}}}][\E^{\mathrm{path},\leftarrow}|\mathscr{A}^{(2),\sim}|^{4}]^{\frac12}\lesssim \ \exp[-\N^{\frac{1}{10}\gamma_{\mathrm{KV}}}][\E^{\mathrm{path},\leftarrow}|\mathscr{A}^{(2),\sim}|^{4}]^{\frac12}.\label{eq:kv1II14e}
\end{align}
}Let us explain \eqref{eq:kv1II14b}. Set $p=\{\mathbb{P}^{\mathrm{Leb},\sigma,\mathfrak{t}(\mathrm{in})+\mathfrak{t},\mathbb{I}(\mathfrak{t})}(\mathcal{E}[\mathfrak{t}(\mathrm{in}),\mathfrak{t};\mathbb{I}(\mathfrak{t})])\}^{-1}$. \eqref{eq:kv1II14b} follows by noting each term in \eqref{eq:kv1II14a} is the expectation of $\E^{\mathrm{path},\sim}|\mathscr{A}^{(2),\sim}|^{2}$ with respect to different initial data. One initial data is $\mathfrak{p}[\mathfrak{t}(\mathrm{in}),\mathfrak{t}(\mathrm{in})+\mathfrak{t}]\d\mathbb{P}^{\mathrm{Leb},\sigma,\mathfrak{t}(\mathrm{in})+\mathfrak{t},\mathbb{I}(\mathfrak{t})}$. It is the law for initial data of the first term in \eqref{eq:kv1II14a}; see the beginning of the paragraph before \eqref{eq:kv1II9a}. As noted after \eqref{eq:kv1II13d}, initial data is $\mathfrak{p}[\mathfrak{t}(\mathrm{in}),\mathfrak{t}(\mathrm{in})+\mathfrak{t};\mathcal{E}]\d\mathbb{P}^{\mathrm{Leb},\sigma,\mathfrak{t}(\mathrm{in})+\mathfrak{t},\mathbb{I}(\mathfrak{t}),\mathcal{E}}$ for the second term in \eqref{eq:kv1II14a}. But $\d\mathbb{P}^{\mathrm{Leb},\sigma,\mathfrak{t}(\mathrm{in})+\mathfrak{t},\mathbb{I}(\mathfrak{t}),\mathcal{E}}$, by Definition \ref{definition:kv4}, is just $\d\mathbb{P}^{\mathrm{Leb},\sigma,\mathfrak{t}(\mathrm{in})+\mathfrak{t},\mathbb{I}(\mathfrak{t})}$ conditioned on $\mathcal{E}[\mathfrak{t}(\mathrm{in}),\mathfrak{t};\mathbb{I}(\mathfrak{t})]$. So, $\mathfrak{p}[\mathfrak{t}(\mathrm{in}),\mathfrak{t}(\mathrm{in})+\mathfrak{t};\mathcal{E}]\d\mathbb{P}^{\mathrm{Leb},\sigma,\mathfrak{t}(\mathrm{in})+\mathfrak{t},\mathbb{I}(\mathfrak{t}),\mathcal{E}}$ is equal to the law $p\times\mathfrak{p}[\mathfrak{t}(\mathrm{in}),\mathfrak{t}(\mathrm{in})+\mathfrak{t};\mathcal{E}]\d\mathbb{P}^{\mathrm{Leb},\sigma,\mathfrak{t}(\mathrm{in})+\mathfrak{t},\mathbb{I}(\mathfrak{t})}$ (upon changing reference measures). Thus, \eqref{eq:kv1II14b} is just comparing expectations of $\E^{\mathrm{path},\sim}|\mathscr{A}^{(2),\sim}|^{2}$ with respect to two different initial data. \eqref{eq:kv1II14c} is by Cauchy-Schwarz. \eqref{eq:kv1II14d} follows from the bounds $p\lesssim1$ (see Lemma \ref{lemma:kv6}) and \eqref{eq:kv8I}-\eqref{eq:kv8II}. Now, note the first expectation in \eqref{eq:kv1II14d} is at most a total variation distance that is controlled by Lemma \ref{lemma:kv11}. This gives the first bound \eqref{eq:kv1II14e}. The rest follows from $\gamma_{\mathrm{av}}\leq{\mathrm{c}}\gamma_{\mathrm{KV}}$ {for some $\mathrm{c}>0$ small}; see Proposition \ref{prop:kv1}. We now combine every display in this step, starting with \eqref{eq:kv1II11}. This gives the following estimate:
\begin{align}
\mathrm{LHS}\eqref{eq:kv1II11} \ \lesssim \ &\E^{\mathrm{path},\to}|\mathscr{A}^{(1)}|^{2}+\E^{\mathrm{path},\leftarrow}|\mathscr{A}^{(2),\sim}|^{2}+\exp[-\N^{\frac{1}{10}\gamma_{\mathrm{KV}}}][\E^{\mathrm{path},\leftarrow}|\mathscr{A}^{(2),\sim}|^{4}]^{\frac12}. \label{eq:kv1II15}
\end{align}
Before we proceed, {we present} a clarifying point. Recall $\E^{\mathrm{path},\to}$ that we defined after \eqref{eq:kv1II1}. $\E^{\mathrm{path},\leftarrow}$ is basically the same expectation except the asymmetric parts of \eqref{eq:glsdeloc}-\eqref{eq:hfloc} have an additional sign and an additional reparameterization of $\s\in\mathfrak{t}(\mathrm{in})+[0,\mathfrak{t}]$ that does not change any of our ideas. Thus, \eqref{eq:kv1II15} is just a decoupling of forwards and $\sim$ processes from Definition \ref{definition:kv14}.
\subsubsection{Final steps: $\mathscr{A}^{(1)}$ estimates}
We now focus on the first term in $\mathrm{RHS}\eqref{eq:kv1II15}$. We first decompose
\begin{align}
\mathscr{A}^{(1)} \ &= \ \mathscr{A}^{(1),1}(\t)+\mathscr{A}^{(1),2}(\t) \nonumber\\
&:= \ \mathscr{A}^{(1),1}(\t)+[\t-\mathfrak{t}(\mathrm{in})]^{-1}{\textstyle\int_{\mathfrak{t}(\mathrm{in})}^{\t}}\mathsf{F}(\s,\mathbf{J}(\s,\inf\mathbb{I}(\mathfrak{t});\mathbb{I}(\mathfrak{t})),\mathbf{U}^{\s,\cdot}[\mathbb{I}(\mathfrak{t})])\d\s, \label{eq:kv1II16a}
\end{align}
where $\mathscr{A}^{(1),1}(\t)$ is defined to make \eqref{eq:kv1II16a} true; by \eqref{eq:kv1II3a}-\eqref{eq:kv1II3b} and \eqref{eq:kv1II5a}-\eqref{eq:kv1II5c}, it is easy to check that
\begin{align}
&\mathscr{A}^{(1),1}(\t) \nonumber\\
&= \ \tfrac{1}{2[\t-\mathfrak{t}(\mathrm{in})]}\{\mathsf{F}(\t,\mathbf{J}(\t,\inf\mathbb{I}(\mathfrak{t});\mathbb{I}(\mathfrak{t})),\mathbf{U}^{\t,\cdot}[\mathbb{I}(\mathfrak{t})])-\mathsf{F}(\mathfrak{t}(\mathrm{in}),\mathbf{J}(\mathfrak{t}(\mathrm{in}),\inf\mathbb{I}(\mathfrak{t});\mathbb{I}(\mathfrak{t})),\mathbf{U}^{\mathfrak{t}(\mathrm{in}),\cdot}[\mathbb{I}(\mathfrak{t})])\}\label{eq:kv1II16b}\\
&- \ \tfrac{1}{2[\t-\mathfrak{t}(\mathrm{in})]}{\textstyle\int_{\mathfrak{t}(\mathrm{in})}^{\t}}\{[\partial_{\s}+\mathscr{L}^{\mathrm{tot}}(\s,\mathbb{I}(\mathfrak{t}))]\mathsf{F}\}(\s,\mathbf{J}(\s,\inf\mathbb{I}(\mathfrak{t});\mathbb{I}(\mathfrak{t})),\mathbf{U}^{\s,\cdot}[\mathbb{I}(\mathfrak{t})])\d\s.\label{eq:kv1II16c}
\end{align}
We first treat $\mathscr{A}^{(1),2}(\t)$ in \eqref{eq:kv1II16a}. We claim the following set of calculations, which we explain after:
\begin{align}
\E^{\mathrm{path},\to}|\mathscr{A}^{(1),2}(\t)|^{2} \ &\lesssim \ [\t-\mathfrak{t}(\mathrm{in})]^{-1}{\textstyle\int_{\mathfrak{t}(\mathrm{in})}^{\t}}\E^{\mathrm{path},\to}|\mathsf{F}(\s,\mathbf{J}(\s,\inf\mathbb{I}(\mathfrak{t});\mathbb{I}(\mathfrak{t})),\mathbf{U}^{\s,\cdot}[\mathbb{I}(\mathfrak{t})])|^{2}\d\s\label{eq:kv1II17a}\\
&\lesssim \ [\t-\mathfrak{t}(\mathrm{in})]^{-1}{\textstyle\int_{\mathfrak{t}(\mathrm{in})}^{\t}}\E^{\mathrm{Leb},\sigma,\s,\mathbb{I}(\mathfrak{t})}\mathfrak{p}[\mathfrak{t}(\mathrm{in}),\s]|\mathsf{F}(\s,\mathrm{a},\mathbf{U})|^{2}\d\s\label{eq:kv1II17b}\\
&= \ [\t-\mathfrak{t}(\mathrm{in})]^{-1}{\textstyle\int_{\mathfrak{t}(\mathrm{in})}^{\t}}\E^{\mathrm{Leb},\sigma,\s,\mathbb{I}(\mathfrak{t})}\mathfrak{p}[\mathfrak{t}(\mathrm{in}),\s;\mathcal{E}]|\mathsf{F}(\s,\mathrm{a},\mathbf{U})|^{2}\d\s\label{eq:kv1II17c}\\
&+ \ [\t-\mathfrak{t}(\mathrm{in})]^{-1}{\textstyle\int_{\mathfrak{t}(\mathrm{in})}^{\t}}\E^{\mathrm{Leb},\sigma,\s,\mathbb{I}(\mathfrak{t})}\{\mathfrak{p}[\mathfrak{t}(\mathrm{in}),\s]-\mathfrak{p}[\mathfrak{t}(\mathrm{in}),\s;\mathcal{E}]\}|\mathsf{F}(\s,\mathrm{a},\mathbf{U})|^{2}\d\s\label{eq:kv1II17d}\\
&\lesssim \ [\t-\mathfrak{t}(\mathrm{in})]^{-1}{\textstyle\int_{\mathfrak{t}(\mathrm{in})}^{\t}}\E^{\mathrm{Leb},\sigma,\s,\mathbb{I}(\mathfrak{t})}|\mathsf{F}(\s,\mathrm{a},\mathbf{U})|^{2}\d\s\label{eq:kv1II17e}\\
&+ \ [\t-\mathfrak{t}(\mathrm{in})]^{-1}{\textstyle\int_{\mathfrak{t}(\mathrm{in})}^{\t}}\{\E^{\mathrm{Leb},\sigma,\s,\mathbb{I}(\mathfrak{t})}|\mathfrak{p}[\mathfrak{t}(\mathrm{in}),\s]-\mathfrak{p}[\mathfrak{t}(\mathrm{in}),\s;\mathcal{E}]|\}\|\mathsf{F}(\s,\cdot,\cdot)\|_{\infty}^{2}\d\s \label{eq:kv1II17f}\\
&\lesssim \ [\t-\mathfrak{t}(\mathrm{in})]^{-1}{\textstyle\int_{\mathfrak{t}(\mathrm{in})}^{\t}}\E^{\mathrm{Leb},\sigma,\s,\mathbb{I}(\mathfrak{t})}|\mathsf{F}(\s,\mathrm{a},\mathbf{U})|^{2}\d\s+\mathrm{RHS}\eqref{eq:kv1II}. \label{eq:kv1II17g}
\end{align}
\eqref{eq:kv1II17a} follows by Cauchy-Schwarz for the time-integral. \eqref{eq:kv1II17b} follows as the expectation in $\mathrm{RHS}\eqref{eq:kv1II17a}$ is of a function of the process at a single time. By {the} construction in Lemma \ref{lemma:kv8}, this law is $\mathfrak{p}[\mathfrak{t}(\mathrm{in}),\s]\d\mathbb{P}^{\mathrm{Leb},\sigma,\s,\mathbb{I}(\mathfrak{t})}$. \eqref{eq:kv1II17c}-\eqref{eq:kv1II17d} is trivial to check. \eqref{eq:kv1II17e}-\eqref{eq:kv1II17f} follows by \eqref{eq:kv8II} and otherwise elementary considerations. Let us explain \eqref{eq:kv1II17g}. Leave \eqref{eq:kv1II17e} alone. Let us now recall $\mathsf{F}$ from \eqref{eq:kv1II2} as the resolvent $\{1-\mathscr{L}^{\mathrm{tot},\mathrm{sym}}(\s,\mathbb{I}(\mathfrak{t}))\}^{-1}$ acting on $\varphi\times\mathfrak{a}$. The resolvent is uniformly bounded as an operator on $\mathscr{L}^{\infty}(\mathbb{S}(\N)\times\mathbb{H}^{\sigma,\mathbb{I}(\mathfrak{t})})$; indeed, we regularized the resolvent by the additional $1$ in the inverse. So, we deduce $\|\mathsf{F}(\s,\cdot,\cdot)\|_{\infty}$ is bounded {from} above by the sup-norm of $|\varphi|$ times that of $|\mathfrak{a}|$, which is $\lesssim$ the maximal sup-norm of $|\mathfrak{a}(\cdot;\mathrm{k})|$ by construction in Proposition \ref{prop:kv1}. Lemma \ref{lemma:kv11} shows the expectation in \eqref{eq:kv1II17f} is $\lesssim\exp[-\N^{\gamma_{\mathrm{KV}}/10}]$; see \eqref{eq:kv10I2c}-\eqref{eq:kv10I3}. So by combining the past two sentences, it is not hard to see that $|\eqref{eq:kv1II17f}|\lesssim\mathrm{RHS}\eqref{eq:kv1II}$. (Note that all $|\mathbb{I}(\mathrm{k})|$-factors in $\mathrm{RHS}\eqref{eq:kv1II}$ are $\geq1$, since they are sizes of non-empty sets.) We now claim the following for the first term in \eqref{eq:kv1II17g}. For $\s\in\mathfrak{t}(\mathrm{in})+[0,\mathfrak{t}]$,
\begin{align}
&\E^{\mathrm{Leb},\sigma,\s,\mathbb{I}(\mathfrak{t})}|\mathsf{F}(\s,\mathrm{a},\mathbf{U})|^{2} \ = \ \E^{\mathrm{Leb},\sigma,\s,\mathbb{I}(\mathfrak{t})}|\{1-\mathscr{L}^{\mathrm{tot},\mathrm{sym}}(\s,\mathbb{I}(\mathfrak{t}))\}^{-1}\varphi(\mathrm{a})\mathfrak{a}(\s,\mathbf{U})|^{2}\label{eq:kv1II18a}\\
&\lesssim \ \E^{\mathrm{Leb},\sigma,\s,\mathbb{I}(\mathfrak{t})}\varphi(\mathrm{a})\mathfrak{a}(\s,\mathbf{U})\times\{1-\mathscr{L}^{\mathrm{tot},\mathrm{sym}}(\s,\mathbb{I}(\mathfrak{t}))\}^{-2}\varphi(\mathrm{a})\mathfrak{a}(\s,\mathbf{U})\label{eq:kv1II18b}\\
&\lesssim \ \E^{\mathrm{Leb},\sigma,\s,\mathbb{I}(\mathfrak{t})}\varphi(\mathrm{a})\mathfrak{a}(\s,\mathbf{U})\times\{1-\mathscr{L}^{\mathrm{tot},\mathrm{sym}}(\s,\mathbb{I}(\mathfrak{t}))\}^{-1}\varphi(\mathrm{a})\mathfrak{a}(\s,\mathbf{U}) \ \lesssim \ \mathrm{RHS}\eqref{eq:kv1II}. \label{eq:kv1II18c}
\end{align}
\eqref{eq:kv1II18a} is by construction; see \eqref{eq:kv1II2}. \eqref{eq:kv1II18b} follows as $\{1-\mathscr{L}^{\mathrm{tot},\mathrm{sym}}(\s,\mathbb{I}(\mathfrak{t}))\}^{-1}$ is self-adjoint by definition; see Lemma \ref{lemma:le7}. The first bound in \eqref{eq:kv1II18c} follows since, in terms of spectrum and quadratic forms, we have $\{1-\mathscr{L}^{\mathrm{tot},\mathrm{sym}}(\s,\mathbb{I}(\mathfrak{t}))\}^{-2}\leq\{1-\mathscr{L}^{\mathrm{tot},\mathrm{sym}}(\s,\mathbb{I}(\mathfrak{t}))\}^{-1}$. (Indeed, since $\{1-\mathscr{L}^{\mathrm{tot},\mathrm{sym}}(\s,\mathbb{I}(\mathfrak{t}))\}$ is self-adjoint, it suffices to note $\{1-\mathscr{L}^{\mathrm{tot},\mathrm{sym}}(\s,\mathbb{I}(\mathfrak{t}))\}\geq1$; this follows because $\mathscr{L}^{\mathrm{tot},\mathrm{sym}}(\s,\mathbb{I}(\mathfrak{t}))\leq0$, as it is the symmetric part of a Markov generator.) The last bound in \eqref{eq:kv1II18c} holds by \eqref{eq:le7II}. We now study $\mathscr{A}^{(1),1}(\t)$ in \eqref{eq:kv1II16a}. By the Ito formula, we know $\mathscr{A}^{(1),1}(\t)$ is a martingale (with respect to the law in $\E^{\mathrm{path},\to}$). Its bracket process enjoys the following deterministic estimate, which we explain afterwards:
\begin{align}
[\mathscr{A}^{(1),1}(\t)] \ \lesssim \ \tfrac{1}{|\t-\mathfrak{t}(\mathrm{in})|^{2}}{\textstyle\int_{\mathfrak{t}(\mathrm{in})}^{\t}}\Gamma^{\mathsf{F}}(\s,\mathbf{J}(\s,\inf\mathbb{I}(\mathfrak{t});\mathbb{I}(\mathfrak{t})),\mathbf{U}^{\s,\cdot}[\mathbb{I}(\mathfrak{t})])\d\s, \label{eq:kv1II19a}
\end{align}
where $\Gamma^{\mathsf{F}}(\s,\cdot,\cdot):\mathbb{S}(\N)\times\mathbb{H}^{\sigma,\mathbb{I}(\mathfrak{t})}\to\R$ is the following carre-du-champ form of $\mathsf{F}$:
\begin{align}
\Gamma^{\mathsf{F}}(\s,\mathrm{a},\mathbf{U}) \ = \ \{\mathscr{L}^{\mathrm{tot}}(\s,\mathbb{I}(\mathfrak{t}))[\mathsf{F}^{2}]\}(\s,\mathrm{a},\mathbf{U})-2\mathsf{F}(\s,\mathrm{a},\mathbf{U})\mathscr{L}^{\mathrm{tot}}(\s,\mathbb{I}(\mathfrak{t}))\mathsf{F}(\s,\mathrm{a},\mathbf{U}). \label{eq:kv1II19b}
\end{align}
{The proof of} \eqref{eq:kv1II19a}-\eqref{eq:kv1II19b} can be found in Appendix 1.5 of \cite{KL}; this proof is based entirely on semigroup theory and extends to any Polish space. We now give an a priori estimate for $\mathrm{RHS}\eqref{eq:kv1II19a}$. We claim the following, the first part of which is standard:
\begin{align}
0 \ \leq \ \Gamma^{\mathsf{F}}(\s,\mathrm{a},\mathbf{U}) \ \lesssim \ \N|\partial_{\mathrm{a}}\mathsf{F}(\s,\mathrm{a},\mathbf{U})|^{2}+\N^{3}{\textstyle\sup_{\x}}|\mathrm{D}_{\x}\mathsf{F}(\s,\mathrm{a},\mathbf{U})|^{2}. \label{eq:kv1II20}
\end{align}
\eqref{eq:kv1II20} can be checked by a direct computation using \eqref{eq:kv1II19b} and Definition \ref{definition:le5}. (Indeed, when we apply Definition \ref{definition:le5} to compute \eqref{eq:kv1II19b}, all the second-order operators vanish by the Leibniz and chain rules. Alternatively, \eqref{eq:kv1II19b} is the infinitesimal quadratic variation of the stochastic integrals appearing when applying Ito to $\mathsf{F}$. But Ito says that said stochastic integrals are just first-order operators of $\mathsf{F}$ integrated against constant-coefficient Brownian motions.) Let $\mathfrak{D}$ be $\N^{1/2}\partial_{\mathrm{a}}$ or $\N\mathrm{D}_{\x}$. We now claim the following \emph{a priori estimate}, which does not reveal any square-root cancellations, since it does not use any precise information about the spectrum of $\mathscr{L}^{\mathrm{tot},\mathrm{sym}}(\s,\mathbb{I}(\mathfrak{t}))$ (besides the fact that it is supported on $(-\infty,0]$):
\begin{align}
|\mathfrak{D}\mathsf{F}(\s,\mathrm{a},\mathbf{U})| \ &= \ |\mathfrak{D}\{1-\mathscr{L}^{\mathrm{tot},\mathrm{sym}}(\s,\mathbb{I}(\mathfrak{t}))\}^{-1}\varphi(\mathrm{a})\mathfrak{a}(\s,\mathbf{U})| \label{eq:kv1II21a}\\
&= \ |\mathfrak{D}{\textstyle\int_{0}^{\infty}}\exp[-\r]\times\exp[\r\mathscr{L}^{\mathrm{tot},\mathrm{sym}}(\s,\mathbb{I}(\mathfrak{t}))]\varphi(\mathrm{a})\mathfrak{a}(\s,\mathbf{U})|\d\r\label{eq:kv1II21b}\\
&\leq \ {\textstyle\int_{0}^{\infty}}\exp[-\r]\times|\mathfrak{D}\exp[\r\mathscr{L}^{\mathrm{tot},\mathrm{sym}}(\s,\mathbb{I}(\mathfrak{t}))]\varphi(\mathrm{a})\mathfrak{a}(\s,\mathbf{U})|\d\r\label{eq:kv1II21c}\\
&\lesssim \ \|\varphi\|_{\infty}\|\mathfrak{a}(\cdot,\cdot)\|_{\infty}{\textstyle\int_{0}^{\infty}}\exp[-\r]\r^{-\frac12}\d\r \ \lesssim \ \|\varphi\|_{\infty}\|\mathfrak{a}(\cdot,\cdot)\|_{\infty}. \label{eq:kv1II21d}
\end{align}
\eqref{eq:kv1II21a} is definition; see \eqref{eq:kv1II2}. \eqref{eq:kv1II21b} is the usual spectral theorem, since $\mathscr{L}^{\mathrm{tot},\mathrm{sym}}(\s,\mathbb{I}(\mathfrak{t}))$ is a symmetric generator for a finite-dimensional Ito diffusion (see Lemma \ref{lemma:le7}). \eqref{eq:kv1II21c} is by triangle inequality. The last bound in \eqref{eq:kv1II21d} is integration. It remains to justify the first bound in \eqref{eq:kv1II21d}. To this end, we recall $\mathscr{L}^{\mathrm{tot},\mathrm{sym}}(\s,\mathbb{I}(\mathfrak{t}))$ is the generator for an Ito diffusion on a flat space $\mathbb{S}(\N)\times\mathbb{H}^{\sigma,\mathbb{I}(\mathfrak{t})}$. (Here, $\mathbb{S}(\N)$ is flat since we identify it as a quotient of $\R$, not as a subset of $\R^{2}$, for instance.) In particular, $\mathscr{L}^{\mathrm{tot},\mathrm{sym}}(\s,\mathbb{I}(\mathfrak{t}))$, which we compute in \eqref{eq:le7II0}, has two parts. {The first} is a Laplacian determined by a metric whose tangent space is spanned by mutually orthogonal $\mathfrak{D}$ operators. What remains is the first-order term in $\mathscr{L}^{\mathrm{S}}(\s,\mathbb{I}(\mathfrak{t}))$, which is a vector field given by the gradient of a convex function (see Assumption \ref{ass:intro8}). Thus, the semigroup $\exp[\r\mathscr{L}^{\mathrm{tot},\mathrm{sym}}(\s,\mathbb{I}(\mathfrak{t}))]:\mathscr{C}^{0}(\mathbb{S}(\N)\times\mathbb{H}^{\sigma,\mathbb{I}(\mathfrak{t})})\to\mathscr{C}^{1}(\mathbb{S}(\N)\times\mathbb{H}^{\sigma,\mathbb{I}(\mathfrak{t})})$ has a norm of $\lesssim\r^{-1/2}$. This can be seen from the introduction of \cite{W}. (The point is that the Laplacian has similar smoothing as discussed in the proof of \eqref{eq:kv10II}, and first-order convex perturbations do not make things worse. We note that if we only had $\mathscr{U}''=\mathrm{O}(1)$, then the proposed semigroup bound would hold with an additional factor $\lesssim\N^{\mathrm{O}(1)}$, which is fine, since we always multiply \eqref{eq:kv1II21d} by an exponentially small factor anyway. Indeed, said estimate would follow by interpolating the bound in the introduction of \cite{W} with Corollary 4.2 in \cite{W0}. Alternatively, run a Duhamel strategy as in the proof of Lemma \ref{lemma:kv10}.) This gets the first bound in \eqref{eq:kv1II21d}. We now claim
\begin{align}
&\E^{\mathrm{path},\to}|\mathscr{A}^{(1),1}(\t)|^{2} \ \lesssim \ \E^{\mathrm{path},\to}[\mathscr{A}^{(1),1}(\t)] \nonumber\\
&\lesssim \ \tfrac{1}{[\t-\mathfrak{t}(\mathrm{in})]^{2}}{\textstyle\int_{\mathfrak{t}(\mathrm{in})}^{\t}}\E^{\mathrm{Leb},\sigma,\s,\mathbb{I}(\mathfrak{t})}\mathfrak{p}[\mathfrak{t}(\mathrm{in}),\s]\Gamma^{\mathsf{F}}(\s,\mathrm{a},\mathbf{U})\d\s \label{eq:kv1II22a}\\
&= \ \tfrac{1}{[\t-\mathfrak{t}(\mathrm{in})]^{2}}{\textstyle\int_{\mathfrak{t}(\mathrm{in})}^{\t}}\E^{\mathrm{Leb},\sigma,\s,\mathbb{I}(\mathfrak{t})}\mathfrak{p}[\mathfrak{t}(\mathrm{in}),\s;\mathcal{E}]\Gamma^{\mathsf{F}}(\s,\mathrm{a},\mathbf{U})\d\s\nonumber\\
&+\tfrac{1}{[\t-\mathfrak{t}(\mathrm{in})]^{2}}{\textstyle\int_{\mathfrak{t}(\mathrm{in})}^{\t}}\E^{\mathrm{Leb},\sigma,\s,\mathbb{I}(\mathfrak{t})}\{\mathfrak{p}[\mathfrak{t}(\mathrm{in}),\s]-\mathfrak{p}[\mathfrak{t}(\mathrm{in}),\s;\mathcal{E}]\}\Gamma^{\mathsf{F}}(\s,\mathrm{a},\mathbf{U})\d\s \nonumber\\
&\lesssim \ \tfrac{1}{[\t-\mathfrak{t}(\mathrm{in})]^{2}}{\textstyle\int_{\mathfrak{t}(\mathrm{in})}^{\t}}\E^{\mathrm{Leb},\sigma,\s,\mathbb{I}(\mathfrak{t})}\Gamma^{\mathsf{F}}(\s,\mathrm{a},\mathbf{U})\d\s + \tfrac{1}{[\t-\mathfrak{t}(\mathrm{in})]}\exp[-\N^{\frac19\gamma_{\mathrm{KV}}}]\|\varphi\|_{\infty}\|\mathfrak{a}(\cdot,\cdot)\|_{\infty} \nonumber\\
&\lesssim \ \tfrac{1}{[\t-\mathfrak{t}(\mathrm{in})]^{2}}{\textstyle\int_{\mathfrak{t}(\mathrm{in})}^{\t}}\E^{\mathrm{Leb},\sigma,\s,\mathbb{I}(\mathfrak{t})}\mathsf{F}(\s,\mathrm{a},\mathbf{U})\times\{[-\mathscr{L}^{\mathrm{tot}}(\s,\mathbb{I}(\mathfrak{t}))]\mathsf{F}(\s,\mathrm{a},\mathbf{U})\}\d\s\nonumber\\
&+ \ \tfrac{1}{[\t-\mathfrak{t}(\mathrm{in})]}\exp[-\N^{\frac19\gamma_{\mathrm{KV}}}]\|\varphi\|_{\infty}\|\mathfrak{a}(\cdot,\cdot)\|_{\infty} \nonumber\\
&\lesssim \ \tfrac{1}{[\t-\mathfrak{t}(\mathrm{in})]^{2}}{\textstyle\int_{\mathfrak{t}(\mathrm{in})}^{\t}}\E^{\mathrm{Leb},\sigma,\s,\mathbb{I}(\mathfrak{t})}[\varphi(\mathrm{a})\mathfrak{a}(\s,\mathbf{U})\times\{1-\mathscr{L}^{\mathrm{tot},\mathrm{sym}}(\s,\mathbb{I}(\mathfrak{t}))\}^{-1}\varphi(\mathrm{a})\mathfrak{a}(\s,\mathbf{U})]\d\s \nonumber\\
&+ \ \tfrac{1}{[\t-\mathfrak{t}(\mathrm{in})]}\exp[-\N^{\frac19\gamma_{\mathrm{KV}}}]\|\varphi\|_{\infty}\|\mathfrak{a}(\cdot,\cdot)\|_{\infty} \nonumber\\[-0mm]
&\lesssim \ \mathrm{RHS}\eqref{eq:kv1II}. \label{eq:kv1II22b}
\end{align}
\eqref{eq:kv1II22a} follows first by the Ito isometry. Then pull $\E^{\mathrm{path},\to}$ through the time-integral in $\mathrm{RHS}\eqref{eq:kv1II19a}$. What we get is $\E^{\mathrm{path},\to}$ acting {on} something depending only on the process at time $\s$; this has law $\mathfrak{p}[\mathfrak{t}(\mathrm{in}),\s]\d\mathbb{P}^{\mathrm{Leb},\sigma,\s,\mathbb{I}(\mathfrak{t})}$ by construction in Lemma \ref{lemma:kv8}. This gives \eqref{eq:kv1II22a}. The second bound follows by writing $\mathfrak{p}[\mathfrak{t}(\mathrm{in}),\s]=\mathfrak{p}[\mathfrak{t}(\mathrm{in}),\s;\mathcal{E}]+\mathfrak{p}[\mathfrak{t}(\mathrm{in}),\s]-\mathfrak{p}[\mathfrak{t}(\mathrm{in}),\s;\mathcal{E}]$. To get the third bound, we first use \eqref{eq:kv8II} to control the first term in the second bound. We then use the deterministic estimate \eqref{eq:kv1II21a}-\eqref{eq:kv1II21d} and Lemma \ref{lemma:kv11} to control the second term in the second bound; see \eqref{eq:kv10I2c}-\eqref{eq:kv10I3}. To show the fourth bound, we apply $\E^{\mathrm{Leb},\sigma,\s,\mathbb{I}(\mathfrak{t})}$ to $\mathrm{RHS}\eqref{eq:kv1II19b}$. The first calculation in the proof of Lemma \ref{lemma:le7} shows {that} the first term in $\mathrm{RHS}\eqref{eq:kv1II19b}$ vanishes in the expectation. This explains the fourth bound above. For the fifth bound, we first (trivially) replace $\mathscr{L}^{\mathrm{tot}}$ with $\mathscr{L}^{\mathrm{tot},\mathrm{sym}}$ in the fourth bound above. We then use $-\mathscr{L}^{\mathrm{tot},\mathrm{sym}}(\s,\mathbb{I}(\mathfrak{t}))\leq1-\mathscr{L}^{\mathrm{tot},\mathrm{sym}}(\s,\mathbb{I}(\mathfrak{t}))$. Thus, it suffices to apply the resolvent equation \eqref{eq:kv1II2} to get the fifth bound above. The sixth bound above follows by \eqref{eq:le7II} and {the} reasoning for \eqref{eq:kv1II18c}. This completes our $\mathscr{A}^{(1)}$ estimates.
\subsubsection{Final steps: $\mathscr{A}^{(2),\sim}$ estimates}
Let us bound the last two terms in \eqref{eq:kv1II15}. First recall $\mathscr{A}^{(2),\sim}$ from \eqref{eq:kv1II8a}-\eqref{eq:kv1II8c}. Upon replacing $\s\mapsto(\mathbf{J}(\s,\inf\mathbb{I}(\mathfrak{t})),\mathbf{U}^{\s,\cdot}[\mathbb{I}(\mathfrak{t})])$ by the adjoint process in Definition \ref{definition:kv14} (given initial data $\s=\mathfrak{t}(\mathrm{in})$ distributed as $\mathbb{P}^{\mathrm{Leb},\sigma,\mathfrak{t}(\mathrm{in})+\mathfrak{t},\mathbb{I}(\mathfrak{t})}$), analogous versions of the bounds \eqref{eq:kv1II19a}-\eqref{eq:kv1II19b} and \eqref{eq:kv1II20} hold for $\mathscr{A}^{(2),\sim}$. Precisely, we have 
\begin{align}
[\mathscr{A}^{(2),\sim}(\t)] \ \lesssim \ \tfrac{1}{|\t-\mathfrak{t}(\mathrm{in})|^{2}}{\textstyle\int_{\mathfrak{t}(\mathrm{in})}^{\t}}\Gamma^{\mathsf{F},\sim}(\s[\ast],\mathbf{J}^{\sim}(\s,\inf\mathbb{I}(\mathfrak{t});\mathbb{I}(\mathfrak{t})),\mathbf{U}^{\s,\cdot,\sim}[\mathbb{I}(\mathfrak{t})])\d\s, \label{eq:kv1II23a}
\end{align}
where $\Gamma^{\mathsf{F},\sim}(\s,\cdot,\cdot):\mathbb{S}(\N)\times\mathbb{H}^{\sigma,\mathbb{I}(\mathfrak{t})}\to\R$ is the following ``adjoint" carre-du-champ form of $\mathsf{F}$:
\begin{align}
\Gamma^{\mathsf{F},\sim}(\s,\mathrm{a},\mathbf{U}) \ = \ \{\mathscr{L}^{\mathrm{tot}}(\s[\ast],\mathbb{I}(\mathfrak{t}))^{\ast}[\mathsf{F}^{2}]\}(\s[\ast],\mathrm{a},\mathbf{U})-2\mathsf{F}(\s[\ast],\mathrm{a},\mathbf{U})\mathscr{L}^{\mathrm{tot}}(\s[\ast],\mathbb{I}(\mathfrak{t}))^{\ast}\mathsf{F}(\s[\ast],\mathrm{a},\mathbf{U}). \label{eq:kv1II23b}
\end{align}
Moreover, we have the following a priori estimate for said adjoint form; see the reasoning for \eqref{eq:kv1II20}:
\begin{align}
0 \ \leq \ \Gamma^{\mathsf{F},\sim}(\s,\mathrm{a},\mathbf{U}) \ \lesssim \ \N|\partial_{\mathrm{a}}\mathsf{F}(\s[\ast],\mathrm{a},\mathbf{U})|^{2}+\N^{3}{\textstyle\sup_{\x}}|\mathrm{D}_{\x}\mathsf{F}(\s[\ast],\mathrm{a},\mathbf{U})|^{2}. \label{eq:kv1II23c}
\end{align}
Let $\mathfrak{p}^{\sim}[\mathfrak{t}(\mathrm{in}),\s]$ be the density for the law of $(\mathbf{J}^{\sim}(\s,\inf\mathbb{I}(\mathfrak{t});\mathbb{I}(\mathfrak{t})),\mathbf{U}^{\s,\cdot,\sim}[\mathbb{I}(\mathfrak{t})])$ with respect to $\mathbb{P}^{\mathrm{Leb},\sigma,\s[\ast],\mathbb{I}(\mathfrak{t})}$, and $\mathfrak{p}^{\sim}[\mathfrak{t}(\mathrm{in}),\s;\mathcal{E}]$ be the density, with respect to $\mathbb{P}^{\mathrm{Leb},\sigma,\s[\ast],\mathbb{I}(\mathfrak{t}),\mathcal{E}}$, of $(\mathbf{J}^{\sim}(\s,\inf\mathbb{I}(\mathfrak{t});\mathbb{I}(\mathfrak{t})),\mathbf{U}^{\s,\cdot,\sim}[\mathbb{I}(\mathfrak{t})])$ stopped at the minimal $\s\in\mathfrak{t}(\mathrm{in})+[0,\mathfrak{t}]$ that $\mathbf{U}^{\s,\cdot,\sim}[\mathbb{I}(\mathfrak{t})]$ leaves $\mathcal{E}[\mathfrak{t}(\mathrm{in}),\mathfrak{t};\mathbb{I}(\mathfrak{t})]$ from Definition \ref{definition:kv4}. We now note that \eqref{eq:kv8I}-\eqref{eq:kv8II} hold if we add a superscript $\sim$ to the $\mathfrak{p}$-terms therein. Indeed, to use the same proof, we just need to know $\mathscr{L}^{\mathrm{tot}}(\s,\mathbb{I}(\mathfrak{t}))^{\ast}$ has $\E^{\mathrm{Leb},\sigma,\s,\mathbb{I}(\mathfrak{t})}$-adjoint given by a Markovian generator. (In this case, said adjoint is just the original generator $\mathscr{L}^{\mathrm{tot}}(\s,\mathbb{I}(\mathfrak{t}))$ itself. Thus, it is true.) We also claim Lemma \ref{lemma:kv11} holds if we add $\sim$-superscripts to the processes therein. (Indeed, to use the same proof, all we need is \eqref{eq:kv8I} but for $\mathfrak{p}^{\sim}[\mathfrak{t}(\mathrm{in}),\s]$, which we just explained.) Given this paragraph and \eqref{eq:kv1II23a}-\eqref{eq:kv1II23c}, we can follow the calculation \eqref{eq:kv1II22a}-\eqref{eq:kv1II22b} to get
\begin{align}
\E^{\mathrm{path},\leftarrow}|\mathscr{A}^{(2),\sim}(\t)|^{2} \ \lesssim \ \mathrm{RHS}\eqref{eq:kv1II}. \label{eq:kv1II24}
\end{align}
(We clarify $\mathscr{A}^{(2),\sim}$ is the $\sim$ version of $\mathscr{A}^{(1),1}$, not of $\mathscr{A}^{(1)}$. Indeed, compare \eqref{eq:kv1II8a}-\eqref{eq:kv1II8c} to \eqref{eq:kv1II16b}-\eqref{eq:kv1II16c} via \eqref{eq:kv1II3a}-\eqref{eq:kv1II3b}.) We now estimate the last term in \eqref{eq:kv1II15}. To this end, we claim the following (with explanation after):
\begin{align}
\E^{\mathrm{path},\leftarrow}|\mathscr{A}^{(2),\sim}(\t)|^{4} \ &\lesssim \ \E^{\mathrm{path},\leftarrow}[\mathscr{A}^{(2),\sim}(\t)]^{2} \nonumber\\
&\lesssim \ \tfrac{1}{[\t-\mathfrak{t}(\mathrm{in})]^{3}}{\textstyle\int_{\mathfrak{t}(\mathrm{in})}^{\t}}\E^{\mathrm{path},\leftarrow}|\Gamma^{\mathsf{F},\sim}(\s[\ast],\mathbf{J}^{\sim}(\s,\inf\mathbb{I}(\mathfrak{t});\mathbb{I}(\mathfrak{t})),\mathbf{U}^{\s,\cdot,\sim}[\mathbb{I}(\mathfrak{t})])|^{2}\d\s \nonumber\\
&\lesssim \ \tfrac{1}{[\t-\mathfrak{t}(\mathrm{in})]^{2}}\N^{3}\|\varphi\|_{\infty}\|^{4}\mathfrak{a}(\cdot,\cdot)\|_{\infty}^{4}. \label{eq:kv1II25}
\end{align}
The first two lines follow by BDG, \eqref{eq:kv1II23a}, and Cauchy-Schwarz. The last line follows by \eqref{eq:kv1II23c} and \eqref{eq:kv1II21a}-\eqref{eq:kv1II21d}. Now, it is direct to verify that the following estimate for the last term in \eqref{eq:kv1II15} {holds}:
\begin{align}
\exp[-\N^{\frac{1}{10}\gamma_{\mathrm{KV}}}][\E^{\mathrm{path},\leftarrow}|\mathscr{A}^{(2),\sim}|^{4}]^{\frac12} \ \lesssim \ \mathrm{RHS}\eqref{eq:kv1II}. \label{eq:kv1II26}
\end{align}
The point is that the exponentially-decaying factor overwhelms the $\N^{3}$ factor in \eqref{eq:kv1II25}. The dependence on $\t-\mathfrak{t}(\mathrm{in})$ is correct because we take square-roots of \eqref{eq:kv1II25}. The same argument also shows the correct scaling dependence on $\varphi,\mathfrak{a}$.
\subsubsection{Putting it altogether}
\eqref{eq:kv1II} follows from \eqref{eq:kv1II1}, \eqref{eq:kv1II6}, \eqref{eq:kv1II10}, \eqref{eq:kv1II15}, \eqref{eq:kv1II16a}, \eqref{eq:kv1II17a}-\eqref{eq:kv1II17g}, \eqref{eq:kv1II18a}-\eqref{eq:kv1II18c}, \eqref{eq:kv1II22a}-\eqref{eq:kv1II22b}, \eqref{eq:kv1II24}, \eqref{eq:kv1II26}. Thus, we are left to prove the estimate \eqref{eq:kv1III}.
\subsubsection{Proof of \eqref{eq:kv1III}}
We use the reasoning for \eqref{eq:kv1II1} and then \eqref{eq:kv1II4a}. This lets us write
\begin{align}
\mathrm{LHS}\eqref{eq:kv1III} \ = \ \N^{-20\gamma_{\mathrm{reg}}}\E\{\mathscr{B}^{2}\mathbf{1}[|\mathscr{A}(\t)|\geq\mathscr{B}\} \ \lesssim \ \E^{\mathrm{path},\to}\{\mathscr{B}^{2}\mathbf{1}[|\mathscr{A}^{(1)}(\t)+\mathscr{A}^{(2)}(\t)|\geq\mathscr{B}]\}. \label{eq:kv1III1}
\end{align}
We now claim the following analog of the estimate \eqref{eq:kv1II10}, which we explain afterwards:
\begin{align}
\mathrm{RHS}\eqref{eq:kv1III1} \ \lesssim \ \N^{\gamma_{\mathrm{KV}}}\mathfrak{t}|\mathbb{I}(\mathfrak{t})|^{\frac12}\mathscr{B}^{2}+\N^{-{\mathrm{D}}}\mathscr{B}^{2}+\E^{\mathrm{path},\to,\leftarrow}\{\mathscr{B}^{2}\mathbf{1}[|\mathscr{A}^{(1)}(\t)+\mathscr{A}^{(2),\sim}(\t)|\geq\mathscr{B}]\}. \label{eq:kv1III2}
\end{align}
Indeed, the same Girsanov argument giving \eqref{eq:kv1II10} also provides \eqref{eq:kv1III2}. The only difference is that we consider a different path-space functional in which we replace $\mathscr{A}^{(2)}(\t)$ by $\mathscr{A}^{(2),\sim}(\t)$. But all we need for \eqref{eq:kv1III2} to hold is the path-space functional to be deterministically $\lesssim\mathscr{B}^{2}$, which is clear for anything of the type $\mathscr{B}^{2}\mathbf{1}[|\mathfrak{X}|\geq\mathscr{B}]$ that we consider here. By Chebyshev,
\begin{align}
\E^{\mathrm{path},\to,\leftarrow}\{\mathscr{B}^{2}\mathbf{1}[|\mathscr{A}^{(1)}(\t)+\mathscr{A}^{(2),\sim}(\t)|\geq\mathscr{B}]\} \ \lesssim \ \mathrm{LHS}\eqref{eq:kv1II11}.
\end{align}
We already have $\mathrm{LHS}\eqref{eq:kv1II11}\lesssim\mathrm{RHS}\eqref{eq:kv1II}$; see \eqref{eq:kv1II15}, \eqref{eq:kv1II16a}, \eqref{eq:kv1II17a}-\eqref{eq:kv1II17g}, \eqref{eq:kv1II18a}-\eqref{eq:kv1II18c}, \eqref{eq:kv1II22a}-\eqref{eq:kv1II22b}, \eqref{eq:kv1II24}, \eqref{eq:kv1II26}. Using this with the previous three displays shows \eqref{eq:kv1III}. This completes the proof. \qed
%
%
%
\section{Proof of Proposition \ref{prop:bg213} (modulo its main ingredient, Proposition \ref{prop:bg21310})}\label{section:bg213main}
For the purposes of making everything clear, we specify $\d=0$ in the statement of Proposition \ref{prop:bg213}. The case of $\d=1,2,3,4$ uses the same argument. Indeed, we just replace $\mathbf{H}^{\N}$ with $(\mathscr{T}^{\pm,\mathrm{j}})^{\d}\mathbf{H}^{\N}$. (In a nutshell, all we need are some pointwise bounds for the kernel of $\mathbf{H}^{\N}$. Said bounds are satisfied by the kernel of $(\mathscr{T}^{\pm,\mathrm{j}})^{\d}\mathbf{H}^{\N}$; see Proposition \ref{prop:hke}.) The key ingredient in the proof of Proposition \ref{prop:bg213} that we leave out in this section will be proved in the next section. It amounts to establishing all the technical estimates that we need in the doubly-multiscale scheme from Section \ref{section:msII}. We clarify this point when relevant. 

Again, we refer to Section \ref{section:msII} for context in comparing this section (and the next) to Section 7 of \cite{Y} (and also for an intuitive description of what this section is trying to do).
\subsection{Ingredients and constructions}
As for the organization of this section, we first list the ingredients and constructions, use them to get Proposition \ref{prop:bg213}, and then give their proofs. (Again, this is modulo proof of the main ingredient.)
\subsubsection{Multiscale decomposition}
Recall the multiscale scheme in Section \ref{section:msII}. The first construction we give is a collection of space-time scales and blocks to run that strategy. (In particular, these are the local blocks that we average on and glue.)
\begin{definition}\label{definition:bg2131}
 Recall $\mathfrak{l}(\mathrm{j})=\N^{\mathrm{j}\delta_{\mathrm{BG}}}\wedge\mathfrak{l}(\mathrm{hom})$ with $\mathfrak{l}(\mathrm{hom})=\lfloor\N^{2/3+\gamma_{\mathrm{KL}}}\rfloor$ and $\mathfrak{m}(\mathrm{j})=\N^{3/4+\alpha(\mathrm{j})}\mathfrak{l}(\mathrm{j})^{-1}$ in Definitions \ref{definition:bg24}, \ref{definition:bg29}. Here, $1\leq\mathrm{j}\leq\mathrm{j}(\infty)$, and $2\beta_{\mathrm{BG}}\leq\delta_{\mathrm{BG}}\leq3\beta_{\mathrm{BG}}$ satisfies $\N^{\mathrm{j}(\infty)\delta_{\mathrm{BG}}}=\mathfrak{l}(\mathrm{hom})$, and ${\mathrm{c}}\gamma_{\mathrm{KL}}\leq\alpha(\mathrm{j})\leq2{\mathrm{c}}\gamma_{\mathrm{KL}}$ {for a small but fixed $\mathrm{c}>0$}. Recall $\tau(\mathrm{j})=\N^{-3/2}\N^{3/4+\alpha(\mathrm{j})}=\N^{-3/2}\mathfrak{m}(\mathrm{j})\mathfrak{l}(\mathrm{j})$ in Definition \ref{definition:bg211}. Let us now construct the following for $1\leq\mathrm{j}\leq\mathrm{j}(\infty)$.
\begin{itemize}

\item Set $\mathfrak{m}(\mathrm{j},1)=\lfloor\N^{1/9}\rfloor\mathfrak{l}(\mathrm{j})^{-1}\vee1$. { For any integer $\mathrm{i}\geq1$, define} $\mathfrak{m}(\mathrm{j},\mathrm{i}+1)=\lfloor\N^{\delta_{\mathrm{ap}}}\rfloor\mathfrak{m}(\mathrm{j},\mathrm{i})$, where $\gamma_{\mathrm{ap}}\leq\delta_{\mathrm{ap}}\leq2\gamma_{\mathrm{ap}}$ and $\lfloor\N^{\delta_{\mathrm{ap}}}\rfloor$ divides $\mathfrak{m}(\mathrm{j})\mathfrak{m}(\mathrm{j},1)^{-1}$.
\item Let us now set $\mathrm{i}(\mathrm{j})>0$ as the unique positive integer for which $\mathfrak{m}(\mathrm{j},\mathrm{i}(\mathrm{j}))=\mathfrak{m}(\mathrm{j})$.
\item If $\mathfrak{m}(\mathrm{j},\mathrm{i})\mathfrak{l}(\mathrm{j})\leq\N^{1/2}$, set $\tau(\mathrm{j},\mathrm{i})=\N^{-2}\mathfrak{m}(\mathrm{j},\mathrm{i})^{2}\mathfrak{l}(\mathrm{j})^{2}$. If $\mathfrak{m}(\mathrm{j},\mathrm{i})\mathfrak{l}(\mathrm{j})\geq\N^{1/2}$, then set $\tau(\mathrm{j},\mathrm{i})=\N^{-3/2}\mathfrak{m}(\mathrm{j},\mathrm{i})\mathfrak{l}(\mathrm{j})$.
\item Now, set $\beta(\mathrm{j},\mathrm{i})$ so that $\N^{\beta(\mathrm{j},\mathrm{i})}\N^{-1+20\gamma_{\mathrm{reg}}}\tau(\mathrm{j},\mathrm{i})^{-1}\mathfrak{m}(\mathrm{j},\mathrm{i})^{-1}\mathfrak{l}(\mathrm{j})^{-1}=\N^{-90\beta_{\mathrm{BG}}}$. For later convenience, set $\beta(\mathrm{j},0)=\beta(\mathrm{j},1)$.
\end{itemize}
We conclude by defining the associated space-time blocks $\mathds{Q}[\mathrm{j},\mathrm{i}]:=(-\tau(\mathrm{j},\mathrm{i}),0]\times\llbracket0,\mathfrak{m}(\mathrm{j},\mathrm{i})-1\rrbracket$. 
\end{definition}
\begin{rem}\label{remark:bg2132}
 We first recall {the} notation of Definition \ref{definition:bg211}. The space-time blocks $\mathds{Q}[\mathrm{j},\mathrm{i}]$ will serve as the space-time domains for the averaging indices therein. (At least, for our purposes, these are the space-time domains that we take. We clarify that the addition defining $\mathds{Q}$ blocks are with respect to the group structure on $\R\times\mathbb{T}(\N)$.) We now explain the bullet points in Definition \ref{definition:bg2131}. The choice of $\beta(\mathrm{j},\mathrm{i})$ is technical. They determine upper bound cutoffs in Section \ref{section:msII}. Next, we picked $\tau(\mathrm{j},\mathrm{i})$ such that if $\mathbb{I}\subseteq\mathbb{T}(\N)$ is a discrete interval with length $\mathfrak{m}(\mathrm{j},\mathrm{i})\mathfrak{l}(\mathrm{j})$, then the relevant neighborhood for time $\tau(\mathrm{j},\mathrm{i})$ (in the sense and notation of Definition \ref{definition:le10}) is $\mathbb{I}(\tau(\mathrm{j},\mathrm{i}))$. (The relevance of discrete intervals whose length is equal to $\mathfrak{m}(\mathrm{j},\mathrm{i})\mathfrak{l}(\mathrm{j})$ comes by taking the averaging indices in \eqref{eq:bg211I} to belong to $\mathds{Q}[\mathrm{j},\mathrm{i}]$; we eventually take $\mathfrak{l}(\mathsf{F})=\mathfrak{l}(\mathrm{j})$ in \eqref{eq:bg211I}.) So, in some sense, $\mathds{Q}[\mathrm{j},\mathrm{i}]$ are ``dimensionally optimized".
\end{rem}
The sets $\mathds{Q}[\mathrm{j},\mathrm{i}]$ are each subsets of the biggest block $\mathds{Q}[\mathrm{j}(\infty),\mathrm{i}(\mathrm{j}(\infty))]$. As we have alluded to in Section \ref{section:msII} and Remark \ref{remark:bg2132}, we want to decompose averages over the biggest block into averages on the sub-blocks $\mathds{Q}[\mathrm{j},\mathrm{i}]$. The following construction introduces convenient notation for breaking averages on the biggest block into those on sub-blocks, as well as going between the sub-blocks of different space-time scales (parameterized by $1\leq\mathrm{i}\leq\mathrm{i}(\mathrm{j})$). 
\begin{definition}\label{definition:bg2133}
 {Let us fix $1\leq\mathrm{j}\leq\mathrm{j}(\infty)$ and  $1\leq\mathrm{i}\leq\mathrm{i}(\mathrm{j})$. Let us also fix $\s\geq0$ and $\y\in\mathbb{T}(\N)$.} By our choices in Definition \ref{definition:bg2131}, we may write $\mathds{Q}[\mathrm{j}(\infty),\mathrm{i}(\mathrm{j}(\infty))]=:\mathds{Q}[\mathrm{big}]$ as a disjoint union of shifts of $\mathds{Q}[\mathrm{j},\mathrm{i}]$. (This comes from the even-division-constraint in the first bullet of Definition \ref{definition:bg2131}.) Let $\mathscr{Q}[\mathrm{j},\mathrm{i}]$ denote the set of all such shifts. Now, given any $1\leq\mathrm{i}<\mathrm{i}(\mathrm{j})$ and $\mathds{Q}\in\mathscr{Q}[\mathrm{j},\mathrm{i}+1]$, write $\mathds{Q}$ as a disjoint union of shifts of $\mathds{Q}[\mathrm{j},\mathrm{i}]$. (We can do this because, by Definition \ref{definition:bg2131}, the $\mathfrak{m}(\mathrm{j},\mathrm{i})$-scales differ by positive integer factors.) Let $\mathscr{Q}[\mathds{Q}]$ be the set of all such shifts of $\mathds{Q}[\mathrm{j},\mathrm{i}]$ (appearing in our decomposition of $\mathds{Q}\in\mathscr{Q}[\mathrm{j},\mathrm{i}+1]$).

Now, fix $\mathrm{j}\geq1$. Define the following double-average at scale-indices $\mathrm{i},\mathrm{i}+1$:
\begin{align}
{\textstyle{\sum}^{\mathrm{j},\mathrm{i}+1}{\sum}^{\mathrm{i},+}}\mathds{A}^{\mathds{Q},\pm}[\s,\y(\s)] \ := \ |\mathscr{Q}[\mathrm{j},\mathrm{i}+1]|^{-1}\sum_{\mathds{Q}[+]\in\mathscr{Q}[\mathrm{j},\mathrm{i}+1]}|\mathscr{Q}[\mathds{Q}[+]]|^{-1}\sum_{\mathds{Q}\in\mathscr{Q}[\mathds{Q}[+]]}\mathds{A}^{\mathds{Q},\pm}[\s,\y(\s)], \label{eq:bg2133I}
\end{align}
where $\mathds{A}^{\mathds{Q},\pm}[\s,\y(\s)]$ is {defined to be} the RHS of \eqref{eq:bg211I} {with} $\mathsf{F}=\mathds{R}^{\chi,\mathfrak{q},\pm,\mathrm{j}}$ in Definition \ref{definition:bg26}, {but} instead of averaging over $(\r,\mathrm{k})\in(-\tau,0]\times\llbracket0,\mathfrak{l}-1\rrbracket$, {we} average over $(\r,\mathrm{k})\in\mathds{Q}$. (In this case, $\mathfrak{l}(\mathsf{F})=\mathfrak{l}(\mathrm{j})$; see the paragraph before Lemma \ref{lemma:bg27}.) More generally, given $\phi:\R\to\R$, {we define}
\begin{align}
&{\textstyle{\sum}^{\mathrm{j},\mathrm{i}+1}{\sum}^{\mathrm{i},+}}\phi(\mathds{A}^{\mathds{Q},\pm}[\s,\y(\s)]) \nonumber\\
&:= \ |\mathscr{Q}[\mathrm{j},\mathrm{i}+1]|^{-1}\sum_{\mathds{Q}[+]\in\mathscr{Q}[\mathrm{j},\mathrm{i}+1]}|\mathscr{Q}[\mathds{Q}[+]]|^{-1}\sum_{\mathds{Q}\in\mathscr{Q}[\mathds{Q}[+]]}\phi(\mathds{A}^{\mathds{Q},\pm}[\s,\y(\s)]). \label{eq:bg2133Ia}
\end{align}
If we omit the first (resp. second) sum in $\mathrm{LHS}\eqref{eq:bg2133Ia}$, we mean $\mathrm{RHS}\eqref{eq:bg2133Ia}$ but without the first (resp. second) average. We also establish the following notation for $\mathds{Q}[+]\in\mathscr{Q}[\mathrm{j},\mathrm{i}+1]$ fixed, which we explain in Remark \ref{remark:bg2134} below:
\begin{align}
\mathds{A}^{\mathds{Q}[+],\pm}[\s,\y(\s)] \ &:= \ |\mathscr{Q}[\mathds{Q}[+]]|^{-1}\sum_{\mathds{Q}\in\mathscr{Q}[\mathds{Q}[+]]}\mathds{A}^{\mathds{Q},\pm}[\s,\y(\s)] \label{eq:bg2133IIa}\\
\mathds{S}^{\mathds{Q}[+],\pm}[\s,\y(\s)] \ &:= \ \sup_{\mathds{Q}\in\mathscr{Q}[\mathds{Q}[+]]}|\mathds{A}^{\mathds{Q},\pm}[\s,\y(\s)]|. \label{eq:bg2133IIb}
\end{align}
At certain points, we include the LHS of \eqref{eq:bg2133IIa}-\eqref{eq:bg2133IIb} in the double sums in the LHS of \eqref{eq:bg2133I}-\eqref{eq:bg2133Ia}. In doing so, as suggested by our notation, the LHS of \eqref{eq:bg2133IIa}-\eqref{eq:bg2133IIb} will depend on the sum-variables from the LHS of \eqref{eq:bg2133I}-\eqref{eq:bg2133Ia} only through the outer sum-variable $\mathds{Q}[+]$ therein. (To be completely clear, this dependence on $\mathds{Q}[+]$ is exactly the depedence on $\mathds{Q}[+]$ in \eqref{eq:bg2133IIa}-\eqref{eq:bg2133IIb}.)
\end{definition}
\begin{rem}\label{remark:bg2134}
 The notation \eqref{eq:bg2133I} is meant to encode the following idea/procedure. First, let us take the big block $\mathds{Q}[\mathrm{big}]$. This is the space-time set of shift indices with which we average $\mathds{R}^{\chi,\mathfrak{q},\pm,\mathrm{j}}\mathbf{G}$ in \eqref{eq:bg213I}; see \eqref{eq:bg211I}. Tile this space-time set by mutually disjoint shifts of $\mathds{Q}[\mathrm{j},\mathrm{i}+1]$ in Definition \ref{definition:bg2131}. (This is like tiling a $10\times10$ square by using 100 mutually disjoint $1\times1$ squares.) Then take each shift of $\mathds{Q}[\mathrm{j},\mathrm{i}+1]$ in said tiling, and further tile it with mutually disjoint shifts of $\mathds{Q}[\mathrm{j},\mathrm{i}]$. Of course, averaging on the big block $\mathds{Q}[\mathrm{big}]$ is equivalent to first averaging over each copy of $\mathds{Q}[\mathrm{j},\mathrm{i}+1]$, and then averaging over all copies. By the same token, this also equals averaging the following over all copies of $\mathds{Q}[\mathrm{j},\mathrm{i}+1]$. For each copy of $\mathds{Q}[\mathrm{j},\mathrm{i}+1]$, average over each shift of $\mathds{Q}[\mathrm{j},\mathrm{i}]$ in its own tiling, and then average over all copies of $\mathds{Q}[\mathrm{j},\mathrm{i}+1]$. Thus, in particular, \eqref{eq:bg2133I} is just equal to $\mathds{A}^{\mathds{Q},\pm}[\s,\y(\s)]$ with $\mathds{Q}=\mathds{Q}[\mathrm{big}]$, which is defined with the prescription after \eqref{eq:bg2133I} with $\mathrm{j}=\mathrm{j}(\infty)$ and $\mathrm{i}=\mathrm{i}(\mathrm{j}(\infty))$. \eqref{eq:bg2133IIa} encodes the ``gluing" (from local to global scales) in Section \ref{section:msII}. (\eqref{eq:bg2133IIa} can be treated as a definition. Alternatively, if one defines $\mathrm{LHS}\eqref{eq:bg2133IIa}$ using the prescription after \eqref{eq:bg2133I} with $\mathrm{i}\mapsto\mathrm{i}+1$ and $\mathds{Q}\mapsto\mathds{Q}[+]$, then \eqref{eq:bg2133I} can be checked by reasoning in the previous few sentences about decomposing averages.) Finally, \eqref{eq:bg2133IIb} is just convenient for writing the error terms that we described in Section \ref{section:msII}.
\end{rem}
We now present two multiscale algebraic relations (that we eventually do analysis on to prove Proposition \ref{prop:bg213}). The first basically follows via our intuitive explanation of \eqref{eq:bg2133I} from Remark \ref{remark:bg2134}, but it takes it one step further by introducing cutoffs. The second is an iterative procedure, which we explain more of after its statement.
\begin{lemma}\label{lemma:bg2135}
 {We first recall the} notation of {Definition \ref{definition:bg211}}. Fix $1\leq\mathrm{j}\leq\mathrm{j}(\infty)$ and $\s\geq\tau(\mathrm{j})$ and $\y\in\mathbb{T}(\N)$. With probability 1, we have
\begin{align}
\mathds{A}^{\mathfrak{m}(\mathrm{j}),\tau(\mathrm{j}),\pm}(\mathds{R}^{\chi,\mathfrak{q},\pm,\mathrm{j}}\mathbf{Z};\s,\y(\s)) \ = \ &{\textstyle\sum^{\mathrm{j},2}\sum^{1,+}}\mathds{A}^{\mathds{Q},\pm}[\s,\y(\s)]\mathbf{1}\{|\mathds{A}^{\mathds{Q},\pm}[\s,\y(\s)]|\leq\N^{-\beta(\mathrm{j},1)}\} \label{eq:bg2135Ia}\\
+ \ &{\textstyle\sum^{\mathrm{j},2}\sum^{1,+}}\mathds{A}^{\mathds{Q},\pm}[\s,\y(\s)]\mathbf{1}\{|\mathds{A}^{\mathds{Q},\pm}[\s,\y(\s)]|>\N^{-\beta(\mathrm{j},1)}\}. \label{eq:bg2135Ib}
\end{align}
\end{lemma}
\begin{lemma}\label{lemma:bg2136}
 Fix $1\leq\mathrm{j}\leq\mathrm{j}(\infty)$ and $\s\geq\tau(\mathrm{j})$ and $\y\in\mathbb{T}(\N)$. For any $1\leq\mathrm{i}\leq\mathrm{i}(\mathrm{j})-2$, we have the following decomposition:
\begin{align}
&{\textstyle\sum^{\mathrm{j},\mathrm{i}+1}\sum^{\mathrm{i},+}}\mathds{A}^{\mathds{Q},\pm}[\s,\y(\s)]\mathbf{1}\{|\mathds{A}^{\mathds{Q},\pm}[\s,\y(\s)]|\leq\N^{-\beta(\mathrm{j},\mathrm{i})}\}\label{eq:bg2136Ia}\\
= \ &{\textstyle\sum^{\mathrm{j},\mathrm{i}+2}\sum^{\mathrm{i}+1,+}}\mathds{A}^{\mathds{Q},\pm}[\s,\y(\s)]\mathbf{1}\{|\mathds{A}^{\mathds{Q},\pm}[\s,\y(\s)]|\leq\N^{-\beta(\mathrm{j},\mathrm{i}+1)}\}\label{eq:bg2136Ib}\\
+ \ &{\textstyle\sum^{\mathrm{j},\mathrm{i}+1}\sum^{\mathrm{i},+}}\mathds{A}^{\mathds{Q},\pm}[\s,\y(\s)]\mathbf{1}\{|\mathds{A}^{\mathds{Q},\pm}[\s,\y(\s)]|\leq\N^{-\beta(\mathrm{j},\mathrm{i})}<|\mathds{S}^{\mathds{Q}[+],\pm}[\s,\y(\s)]|\}\label{eq:bg2136Ic}\\
- \ &{\textstyle\sum^{\mathrm{j},\mathrm{i}+1}\sum^{\mathrm{i},+}}\mathds{A}^{\mathds{Q},\pm}[\s,\y(\s)]\mathbf{1}\{|\mathds{A}^{\mathds{Q}[+],\pm}[\s,\y(\s)]|\leq\N^{-\beta(\mathrm{j},\mathrm{i})}<|\mathds{S}^{\mathds{Q}[+],\pm}[\s,\y(\s)]|\}\label{eq:bg2136Id}\\
+ \ &{\textstyle\sum^{\mathrm{j},\mathrm{i}+2}\sum^{\mathrm{i}+1,+}}\mathds{A}^{\mathds{Q},\pm}[\s,\y(\s)]\mathbf{1}\{\N^{-\beta(\mathrm{j},\mathrm{i}+1)}<|\mathds{A}^{\mathds{Q},\pm}[\s,\y(\s)]|\leq\N^{-\beta(\mathrm{j},\mathrm{i})}\}. \label{eq:bg2136Ie}
\end{align}
\end{lemma}
\begin{rem}\label{remark:bg2137}
 Lemma \ref{lemma:bg2135} was explained prior to its statement; again, we refer to Remark \ref{remark:bg2134} for a brief explanation. We note that the constraint $\s\geq\tau(\mathrm{j})$ in Lemmas \ref{lemma:bg2135} is just to avoid looking at \eqref{eq:hf}-\eqref{eq:glsde} at negative times. (Indeed, $\mathrm{LHS}\eqref{eq:bg2135Ia}$ averages backwards-in-time.) Let us explain Lemma \ref{lemma:bg2136}. Observe \eqref{eq:bg2136Ib} is just \eqref{eq:bg2136Ia} but we replace $\mathrm{i}\mapsto\mathrm{i}+1$ in the latter. In particular, Lemma \ref{lemma:bg2136} is a precise formulation for the multiscale upgrading in space-time scales from Section \ref{section:msII}; \eqref{eq:bg2136Ic}-\eqref{eq:bg2136Ie} are the resulting error terms at each step in this multiscale strategy. (We emphasize that $\eqref{eq:bg2136Ia}\mapsto\eqref{eq:bg2136Ib}$ also includes an upgrade in the cutoff exponent $\beta(\mathrm{j},\mathrm{i})\mapsto\beta(\mathrm{j},\mathrm{i}+1)$. As discussed in Section \ref{section:msII}, this just quantitatively encodes the idea that averaging over larger sets implies more cancellations.) Indeed (see notation in Definition \ref{definition:bg2133}), we do:
\begin{itemize}
\item In \eqref{eq:bg2136Ia}, put in the indicator that the sup over all $\mathds{A}^{\mathds{Q},\pm}[\s,\y(\s)]$-terms (which are averages over a space-time block of dimension-index $\mathrm{i}$) in a fixed/common space-time block of dimension-index $\mathrm{i}+1$ is $\leq\N^{-\beta(\mathrm{j},\mathrm{i})}$. I.e., multiply the summand in \eqref{eq:bg2136Ia} by the indicator of $|\mathds{S}^{\mathds{Q}[+],\pm}[\s,\y(\s)]|\leq\N^{-\beta(\mathrm{j},\mathrm{i})}$; the error we have to account for is \eqref{eq:bg2136Ic}.
\item Since $|\mathds{S}^{\mathds{Q}[+],\pm}[\s,\y(\s)]|$ controls $\mathds{A}^{\mathds{Q},\pm}[\s,\y(\s)]$ for all $\mathds{Q}\subseteq\mathds{Q}[+]$, we can now drop the indicator $\mathbf{1}\{|\mathds{A}^{\mathds{Q},\pm}[\s,\y(\s)]|\leq\N^{-\beta(\mathrm{j},\mathrm{i})}\}$ in \eqref{eq:bg2136Ia}. Then, average $\mathds{A}^{\mathds{Q},\pm}[\s,\y(\s)]$-terms in \eqref{eq:bg2136Ia} over all $\mathds{Q}\subseteq\mathds{Q}[+]$ to get $\mathds{A}^{\mathds{Q}[+],\pm}[\s,\y(\s)]$. In words, the average over the bigger block $\mathds{Q}[+]$ is the average of the averages over the smaller blocks $\mathds{Q}\subseteq\mathds{Q}[+]$. 
\item For each $\mathds{A}^{\mathds{Q}[+],\pm}[\s,\y(\s)]$, we have with it a factor of $\mathbf{1}\{|\mathds{S}^{\mathds{Q}[+],\pm}[\s,\y(\s)]|\leq\N^{-\beta(\mathrm{j},\mathrm{i})}\}$. Trade this in for the indicator $\mathbf{1}\{|\mathds{A}^{\mathds{Q}[+],\pm}[\s,\y(\s)]|\leq\N^{-\beta(\mathrm{j},\mathrm{i})}\}$; the error is given by \eqref{eq:bg2136Id}. (Indeed, the supremum $|\mathds{S}^{\mathds{Q}[+],\pm}[\s,\y(\s)]|$ controls the average $|\mathds{A}^{\mathds{Q}[+],\pm}[\s,\y(\s)]|$. Also, the sign in \eqref{eq:bg2136Id} is not important, since we will control everything in absolute value.) We are then left with 
\begin{align}
{\textstyle\sum^{\mathrm{j},\mathrm{i}+2}\sum^{\mathrm{i}+1,+}}\mathds{A}^{\mathds{Q},\pm}[\s,\y(\s)]\mathbf{1}\{|\mathds{A}^{\mathds{Q},\pm}[\s,\y(\s)]|\leq\N^{-\beta(\mathrm{j},\mathrm{i})}\},
\end{align}
which is just \eqref{eq:bg2136Ib} but with $\beta(\mathrm{j},\mathrm{i})$ instead of $\beta(\mathrm{j},\mathrm{i}+1)$. So, to get \eqref{eq:bg2136Ib} from the previous display, the error we have to pay is \eqref{eq:bg2136Ie}, and Lemma \ref{lemma:bg2136} follows. We write this whole argument precisely shortly.
\end{itemize}
We clarify that Lemmas \ref{lemma:bg2135}, \ref{lemma:bg2136} will only be used to derive Lemma \ref{lemma:bg2138} below.
\end{rem}
Ultimately, we apply Lemma \ref{lemma:bg2135} and then Lemma \ref{lemma:bg2136} for indices $\mathrm{i}=1,\ldots,\mathrm{i}(\mathrm{j})-2$. This gives us \eqref{eq:bg2136Ib} with $\mathrm{i}=\mathrm{i}(\mathrm{j})-2$, which is an average of terms that are deterministically $\lesssim\N^{-\beta(\mathrm{j},\mathrm{i}(\mathrm{j})-1)}\ll\N^{-1}$. This certainly beats the factor of $\N$ in $\mathrm{LHS}\eqref{eq:bg213I}$; see Definition \ref{definition:bg211}. In particular, we are left with controlling error terms \eqref{eq:bg2136Ic}-\eqref{eq:bg2136Ie} coming from each application of Lemma \ref{lemma:bg2136}. (Because $\mathrm{i}(\mathrm{j})$ is bounded by the number of steps of size $\gtrsim1$ needed to go from $0\mapsto1$, we know $\mathrm{i}(\mathrm{j})\lesssim1$. Thus, the number of times we must apply Lemma \ref{lemma:bg2136} is $\mathrm{O}(1)$. This implies that we have $\mathrm{O}(1)$-many sets of error terms \eqref{eq:bg2136Ic}-\eqref{eq:bg2136Ie}, so they do not accumulate in the large-$\N$ limit.) The following result just precisely says what we heuristically discussed in this paragraph.
\begin{lemma}\label{lemma:bg2138}
 Fix any $1\leq\mathrm{j}\leq\mathrm{j}(\infty)$ and set $\|\|=\|\|_{\t_{\mathrm{st}};\mathbb{T}(\N)}$. Recall {Definition \ref{definition:bg211}}. With notation explained after, we have
\begin{align}
\E\|\mathscr{A}^{\mathbf{X},\mathbf{T}}\mathscr{R}^{\chi,\mathfrak{q},\pm,\mathrm{j}}\| \ &\lesssim \ \E\|{\textstyle\int_{\tau(\mathrm{j})}^{\t}}\mathbf{1}[\s\leq\t_{\mathrm{st}}]\mathbf{H}^{\N}(\s,\t(\N),\x)\{\N|\Phi^{\pm,\mathrm{j}}(\s,\cdot(\s))|\}\d\s\| \label{eq:bg2138Ia}\\
&+ \ \N^{\beta_{\mathrm{BG}}}{\textstyle\int_{\tau(\mathrm{j})}^{1}}|\mathbb{T}(\N)|^{-1}{\textstyle\sum_{\y}}\E[\mathbf{1}(\s\leq\t_{\mathrm{st}})\N|\Psi^{\pm,\mathrm{j}}(\s,\y(\s))|]\d\s \label{eq:bg2138Ib}\\
&+ \ {\textstyle\sup_{\mathrm{i}}}\N^{\beta_{\mathrm{BG}}}{\textstyle\int_{\tau(\mathrm{j})}^{1}}|\mathbb{T}(\N)|^{-1}{\textstyle\sum_{\y}}\E[\mathbf{1}(\s\leq\t_{\mathrm{st}})\N|\Upsilon^{\pm,\mathrm{j},\mathrm{i}}(\s,\y(\s))|]\d\s\label{eq:bg2138Ic}\\
&+ \ {\textstyle\sup_{\mathrm{i}}}\N^{\beta_{\mathrm{BG}}}{\textstyle\int_{\tau(\mathrm{j})}^{1}}|\mathbb{T}(\N)|^{-1}{\textstyle\sum_{\y}}\E[\mathbf{1}(\s\leq\t_{\mathrm{st}})\N|\Lambda^{\pm,\mathrm{j},\mathrm{i}}(\s,\y(\s))|]\d\s. \label{eq:bg2138Id}
\end{align}
In $\mathrm{RHS}\eqref{eq:bg2138Ia}$, $\|\|$ is with respect to $(\t,\x)$ therein. See {Definition \ref{definition:method8}} for $\beta_{\mathrm{BG}}$. Suprema in \eqref{eq:bg2138Ib}-\eqref{eq:bg2138Id} are over $1\leq\mathrm{i}<\mathrm{i}(\mathrm{j})$; see {Definition \ref{definition:bg2131}} for $\mathrm{i}(\mathrm{j})$. We have also introduced the following capital Greek letters, which we explain \emph{Remark \ref{remark:bg2139}}:
\begin{align}
\Phi^{\pm,\mathrm{j}}(\s,\y(\s)) \ &:= \ {\textstyle\sum^{\mathrm{j},\mathrm{i}(\mathrm{j})}\sum^{\mathrm{i}(\mathrm{j})-1,+}}\mathds{A}^{\mathds{Q},\pm}[\s,\y(\s)]\mathbf{1}\{|\mathds{A}^{\mathds{Q},\pm}[\s,\y(\s)]|\leq\N^{-\beta(\mathrm{j},\mathrm{i}(\mathrm{j})-1)}\} \label{eq:bg2138IIa}\\
\Psi^{\pm,\mathrm{j}}(\s,\y(\s)) \ &:= \ {\textstyle\sum^{\mathrm{j},2}\sum^{1,+}}\mathds{A}^{\mathds{Q},\pm}[\s,\y(\s)]\mathbf{1}\{|\mathds{A}^{\mathds{Q},\pm}[\s,\y(\s)]|>\N^{-\beta(\mathrm{j},1)}\} \label{eq:bg2138IIb}\\
\Upsilon^{\pm,\mathrm{j},\mathrm{i}}(\s,\y(\s)) \ &:= \ {\textstyle\sum^{\mathrm{j},\mathrm{i}+1}\sum^{\mathrm{i},+}}\N^{\beta(\mathrm{j},\mathrm{i}-1)}|\mathds{A}^{\mathds{Q},\pm}[\s,\y(\s)]|^{2}\mathbf{1}\{|\mathds{A}^{\mathds{Q},\pm}[\s,\y(\s)]|\leq\N^{-\beta(\mathrm{j},\mathrm{i}-1)}\} \label{eq:bg2138IIc}\\
\Lambda^{\pm,\mathrm{j},\mathrm{i}}(\s,\y(\s)) \ &:= \ {\textstyle\sum^{\mathrm{j},\mathrm{i}+1}\sum^{\mathrm{i},+}}\N^{-\beta(\mathrm{j},\mathrm{i}-1)}\mathbf{1}\{|\mathds{A}^{\mathds{Q},\pm}[\s,\y(\s)]|>\N^{-\beta(\mathrm{j},\mathrm{i})}\}.\label{eq:bg2138IId}
\end{align}
\end{lemma}
\begin{rem}\label{remark:bg2139}
 Recall from Definition \ref{definition:bg211} that $\mathscr{A}^{\mathbf{X},\mathbf{T}}\mathscr{R}^{\chi,\mathfrak{q},\pm,\mathrm{j}}$ is a time-integrated heat operator acting on $\mathrm{LHS}\eqref{eq:bg2135Ia}$. $\mathrm{RHS}\eqref{eq:bg2138Ia}$ comes from applying Lemma \ref{lemma:bg2135} to $\mathrm{LHS}\eqref{eq:bg2135Ia}$, and iteratively applying Lemma \ref{lemma:bg2136} to $\mathrm{RHS}\eqref{eq:bg2135Ia}$ and \eqref{eq:bg2136Ia} until we hit index $\mathrm{i}=\mathrm{i}(\mathrm{j})-2$. This gives \eqref{eq:bg2138IIa}, which we then integrate against the heat operator. In doing so, the first error term we pick up is \eqref{eq:bg2135Ib}, which is just $\Psi^{\pm,\mathrm{j}}$. Integrating this against the heat operator gives \eqref{eq:bg2138Ib}. The other error terms we have are the sum over $\mathrm{i}$ of \eqref{eq:bg2136Ic}-\eqref{eq:bg2136Ie}. Using the Schwarz inequality, we will eventually show that the error terms are $\lesssim\Upsilon^{\pm,\mathrm{j},\mathrm{i}}+\Lambda^{\pm,\mathrm{j},\mathrm{i}}$. Integrating these against the heat operator gives \eqref{eq:bg2138Ic}-\eqref{eq:bg2138Id}. This basically explains Lemma \ref{lemma:bg2138} up to some cosmetic differences. Such cosmetics include the factors of $\N^{\beta_{\mathrm{BG}}}$ in \eqref{eq:bg2138Ib}-\eqref{eq:bg2138Id}. (These arise via technical calculations and are harmless.) There is also the fact that the cutoff exponents in \eqref{eq:bg2138IIc}-\eqref{eq:bg2138IId} use index $\mathrm{i}-1$, while the dimensions of $\mathds{Q}$-blocks therein use index $\mathrm{i}$. (This is also harmless, because changing the cutoff exponent indices $\mathrm{i}-1\mapsto\mathrm{i}$ introduces very small powers of $\N$; see Definition \ref{definition:bg2131}.)
\end{rem}
\subsubsection{Multiscale estimates}
To prove Proposition \ref{prop:bg213}, it certainly suffices to obtain appropriate bounds for $\mathrm{RHS}\eqref{eq:bg2138Ia}$ and \eqref{eq:bg2138Ib}-\eqref{eq:bg2138Id}. This is the goal of the following, which is the main ingredient whose proof we defer to the next section (since it is fairly complicated). We emphasize the RHS of \eqref{eq:bg21310I} below is much smaller than what we claim in Proposition \ref{prop:bg213}. This is one reason why we are very willing to discuss estimates modulo very small powers of $\N$ (like in Remark \ref{remark:bg2139}).
\begin{prop}\label{prop:bg21310}
 Retain the notation of {Lemma \ref{lemma:bg2138}}, and recall $\beta_{\mathrm{BG}}$ from {Definition \ref{definition:method8}}. We have the expectation estimates
\begin{align}
\mathrm{RHS}\eqref{eq:bg2138Ia}+\eqref{eq:bg2138Ib}+\eqref{eq:bg2138Ic}+\eqref{eq:bg2138Id} \ \lesssim \ \N^{-9\beta_{\mathrm{BG}}}. \label{eq:bg21310I}
\end{align}
\end{prop}
\subsection{Proof of Proposition \ref{prop:bg213}, assuming Lemmas \ref{lemma:bg2135}, \ref{lemma:bg2136}, \ref{lemma:bg2138}, and Proposition \ref{prop:bg21310}}
Lemma \ref{lemma:bg2138} and Proposition \ref{prop:bg21310} {lead to} $\E\mathrm{LHS}\eqref{eq:bg213I}\lesssim\N^{-8\beta_{\mathrm{BG}}}$ for $\d=0$. (This also extends to each $\d=1,2,3,4$; see the beginning of this section.) Thus, by the Markov inequality, we know $\mathrm{LHS}\eqref{eq:bg213I}\lesssim\N^{-7\beta_{\mathrm{BG}}}$ with high probability for any $\d=0,1,2,3,4$. A union bound shows that $\mathrm{LHS}\eqref{eq:bg213I}\lesssim\N^{-7\beta_{\mathrm{BG}}}$ for all $\d=0,1,2,3,4$ on the same high probability event. It now suffices to combine this bound and the deterministic lower bound $\|\mathbf{Z}\|_{\t_{\mathrm{st}};\mathbb{T}(\N)}\gtrsim\N^{-\beta_{\mathrm{BG}}}$, which follows by construction in Definition \ref{definition:method8}. \qed
\subsubsection{The rest of this section}
We will prove Lemmas \ref{lemma:bg2135}, \ref{lemma:bg2136}, \ref{lemma:bg2138}. Again, Proposition \ref{prop:bg21310} is for the next section.
\subsection{Proof of Lemma \ref{lemma:bg2135}}
We claim the following, which was justified in Remark \ref{remark:bg2134} (as we explain after):
\begin{align}
\mathds{A}^{\mathfrak{m}(\mathrm{j}),\tau(\mathrm{j}),\pm}(\mathds{R}^{\chi,\mathfrak{q},\pm,\mathrm{j}}\mathbf{Z};\s,\y(\s)) \ = \ {\textstyle\sum^{\mathrm{j},2}\sum^{1,+}}\mathds{A}^{\mathds{Q},\pm}[\s,\y(\s)]. \label{eq:bg2135I1}
\end{align}
We note that $\mathrm{LHS}\eqref{eq:bg2135I1}$ is an average of the integrand in $\mathrm{RHS}\eqref{eq:bg211I}$ over $(\r,\mathrm{j})\in(-\tau(\mathrm{j}),0]\times\llbracket0,\mathfrak{m}(\mathrm{j})-1\rrbracket=\mathds{Q}[\mathrm{j}(\infty),\mathrm{i}(\mathrm{j}(\infty))]=\mathds{Q}[\mathrm{big}]$; for the last two blocks, see Definitions \ref{definition:bg2133}, \ref{definition:bg2131}. As we explained in Remark \ref{remark:bg2134}, this is the same as the following. Tile $\mathds{Q}[\mathrm{big}]$ by mutually disjoint shifts of $\mathds{Q}[\mathrm{j},2]$. Then, tile each copy of $\mathds{Q}[\mathrm{j},2]$ with mutually disjoint shifts of $\mathds{Q}[\mathrm{j},1]$. This gives us a tiling of $\mathds{Q}[\mathrm{big}]$ by mutually disjoint shifts of $\mathds{Q}[\mathrm{j},1]$. Average the integrand in $\mathrm{RHS}\eqref{eq:bg211I}$ over indices $(\r,\mathrm{j})$ not in $\mathds{Q}[\mathrm{big}]$, but rather in a fixed shift of $\mathds{Q}[\mathrm{j},1]$. Then, average each of these $\mathds{Q}[\mathrm{j},1]$-averages over all copies of $\mathds{Q}[\mathrm{j},1]$ in our tiling of $\mathds{Q}[\mathrm{big}]$. (This is basically the same as writing an average of 10 terms as the average of two separate averages, one over terms with even index and one over terms with odd index.) But $\mathrm{RHS}\eqref{eq:bg2135I1}$ is exactly the two-scale averaging that we just described. This gives us \eqref{eq:bg2135I1}. To get \eqref{eq:bg2135Ia}-\eqref{eq:bg2135Ib} and thus finish this proof, note the indicators therein add to 1. So \eqref{eq:bg2135I1} equals $\mathrm{RHS}\eqref{eq:bg2135Ia}+\eqref{eq:bg2135Ib}$. \qed
\subsection{Proof of Lemma \ref{lemma:bg2136}}
We start by claiming the following calculation holds (we explain it after):
\begin{align}
&{\textstyle\sum^{\mathrm{j},\mathrm{i}+1}\sum^{\mathrm{i},+}}\mathds{A}^{\mathds{Q},\pm}[\s,\y(\s)]\mathbf{1}\{|\mathds{A}^{\mathds{Q},\pm}[\s,\y(\s)]|\leq\N^{-\beta(\mathrm{j},\mathrm{i})}\}\label{eq:bg2136I1a}\\
= \ &{\textstyle\sum^{\mathrm{j},\mathrm{i}+1}\sum^{\mathrm{i},+}}\mathds{A}^{\mathds{Q},\pm}[\s,\y(\s)]\mathbf{1}\{|\mathds{A}^{\mathds{Q},\pm}[\s,\y(\s)]|\leq\N^{-\beta(\mathrm{j},\mathrm{i})}\}\mathbf{1}\{|\mathds{S}^{\mathds{Q}[+],\pm}[\s,\y(\s)]|\leq\N^{-\beta(\mathrm{j},\mathrm{i})}\}\label{eq:bg2136I1b}\\
+ \ &{\textstyle\sum^{\mathrm{j},\mathrm{i}+1}\sum^{\mathrm{i},+}}\mathds{A}^{\mathds{Q},\pm}[\s,\y(\s)]\mathbf{1}\{|\mathds{A}^{\mathds{Q},\pm}[\s,\y(\s)]|\leq\N^{-\beta(\mathrm{j},\mathrm{i})}\}\mathbf{1}\{|\mathds{S}^{\mathds{Q}[+],\pm}[\s,\y(\s)]|>\N^{-\beta(\mathrm{j},\mathrm{i})}\}\label{eq:bg2136I1c} \\
= \ &{\textstyle\sum^{\mathrm{j},\mathrm{i}+1}\sum^{\mathrm{i},+}}\mathds{A}^{\mathds{Q},\pm}[\s,\y(\s)]\mathbf{1}\{|\mathds{S}^{\mathds{Q}[+],\pm}[\s,\y(\s)]|\leq\N^{-\beta(\mathrm{j},\mathrm{i})}\}\label{eq:bg2136I1d}\\
+ \ &{\textstyle\sum^{\mathrm{j},\mathrm{i}+1}\sum^{\mathrm{i},+}}\mathds{A}^{\mathds{Q},\pm}[\s,\y(\s)]\mathbf{1}\{|\mathds{A}^{\mathds{Q},\pm}[\s,\y(\s)]|\leq\N^{-\beta(\mathrm{j},\mathrm{i})}\}\mathbf{1}\{|\mathds{S}^{\mathds{Q}[+],\pm}[\s,\y(\s)]|>\N^{-\beta(\mathrm{j},\mathrm{i})}\}.\label{eq:bg2136I1e}
\end{align}
\eqref{eq:bg2136I1b}-\eqref{eq:bg2136I1c} follows since the last indicator functions in each add to 1. \eqref{eq:bg2136I1d}-\eqref{eq:bg2136I1e} follows because the first indicator in \eqref{eq:bg2136I1b} is redundant, given the second indicator therein. (Indeed, note $|\mathds{A}^{\mathds{Q},\pm}[\s,\y(\s)]|\leq|\mathds{S}^{\mathds{Q}[+],\pm}[\s,\y(\s)]|$; see Definition \ref{definition:bg2133}.) Note \eqref{eq:bg2136I1e} equals \eqref{eq:bg2136Ic}. So by the previous display, in order to prove Lemma \ref{lemma:bg2136}, it suffices to show
\begin{align}
\eqref{eq:bg2136I1d} \ = \ \eqref{eq:bg2136Ib}+\eqref{eq:bg2136Id}+\eqref{eq:bg2136Ie}. \label{eq:bg2136I2}
\end{align}
To prove \eqref{eq:bg2136I2}, the first step is the following calculation that we explain afterwards:
\begin{align}
&\eqref{eq:bg2136I1d}\nonumber\\
&= \ {\textstyle\sum^{\mathrm{j},\mathrm{i}+1}\sum^{\mathrm{i},+}}\mathds{A}^{\mathds{Q},\pm}[\s,\y(\s)]\mathbf{1}\{|\mathds{S}^{\mathds{Q}[+],\pm}[\s,\y(\s)]|\leq\N^{-\beta(\mathrm{j},\mathrm{i})}\}\mathbf{1}\{|\mathds{A}^{\mathds{Q}[+],\pm}[\s,\y(\s)]|\leq\N^{-\beta(\mathrm{j},\mathrm{i})}\}\label{eq:bg2136I3a}\\
&= \ {\textstyle\sum^{\mathrm{j},\mathrm{i}+1}\sum^{\mathrm{i},+}}\mathds{A}^{\mathds{Q},\pm}[\s,\y(\s)]\mathbf{1}\{|\mathds{A}^{\mathds{Q}[+],\pm}[\s,\y(\s)]|\leq\N^{-\beta(\mathrm{j},\mathrm{i})}\}\label{eq:bg2136I3b}\\
&- \ {\textstyle\sum^{\mathrm{j},\mathrm{i}+1}\sum^{\mathrm{i},+}}\mathds{A}^{\mathds{Q},\pm}[\s,\y(\s)]\mathbf{1}\{|\mathds{S}^{\mathds{Q}[+],\pm}[\s,\y(\s)]|>\N^{-\beta(\mathrm{j},\mathrm{i})}\}\mathbf{1}\{|\mathds{A}^{\mathds{Q}[+],\pm}[\s,\y(\s)]|\leq\N^{-\beta(\mathrm{j},\mathrm{i})}\}. \label{eq:bg2136I3c}
\end{align}
\eqref{eq:bg2136I3a} holds since $|\mathds{A}^{\mathds{Q}[+],\pm}[\s,\y(\s)]|\leq|\mathds{S}^{\mathds{Q}[+],\pm}[\s,\y(\s)]|$; see Definition \ref{definition:bg2133} for relevant notation. (Roughly, $\mathds{A}^{\mathds{Q}[+],\pm}[\s,\y(\s)]$ is an average of terms whose absolute values we take a supremum over to get $\mathds{S}^{\mathds{Q}[+],\pm}[\s,\y(\s)]$.) \eqref{eq:bg2136I3b}-\eqref{eq:bg2136I3c} follows because $\mathds{S}$-based indicators in \eqref{eq:bg2136I3a} and \eqref{eq:bg2136I3c} add to 1. Note \eqref{eq:bg2136I3c} equals \eqref{eq:bg2136Id} (signs included). Therefore, to show \eqref{eq:bg2136I2} and thereby complete this proof, by the previous display and the previous sentence, it suffices to show the following instead:
\begin{align}
\eqref{eq:bg2136I3b} \ = \ \eqref{eq:bg2136Ib}+\eqref{eq:bg2136Ie}. \label{eq:bg2136I4}
\end{align}
To show \eqref{eq:bg2136I4}, we begin with the following calculation that we explain afterwards:
\begin{align}
\eqref{eq:bg2136I3b} \ = \ &{\textstyle\sum^{\mathrm{j},\mathrm{i}+1}}\mathbf{1}\{|\mathds{A}^{\mathds{Q}[+],\pm}[\s,\y(\s)]|\leq\N^{-\beta(\mathrm{j},\mathrm{i})}\}{\textstyle\sum^{\mathrm{i},+}}\mathds{A}^{\mathds{Q},\pm}[\s,\y(\s)] \label{eq:bg2136I5a}\\
= \ &{\textstyle\sum^{\mathrm{j},\mathrm{i}+1}}\mathds{A}^{\mathds{Q}[+],\pm}[\s,\y(\s)]\mathbf{1}\{|\mathds{A}^{\mathds{Q}[+],\pm}[\s,\y(\s)]|\leq\N^{-\beta(\mathrm{j},\mathrm{i})}\}\label{eq:bg2136I5b}\\
= \ &{\textstyle\sum^{\mathrm{j},\mathrm{i}+1}}\mathds{A}^{\mathds{Q}[+],\pm}[\s,\y(\s)]\mathbf{1}\{|\mathds{A}^{\mathds{Q}[+],\pm}[\s,\y(\s)]|\leq\N^{-\beta(\mathrm{j},\mathrm{i}+1)}\}\label{eq:bg2136I5c}\\
+ \ &{\textstyle\sum^{\mathrm{j},\mathrm{i}+1}}\mathds{A}^{\mathds{Q}[+],\pm}[\s,\y(\s)]\mathbf{1}\{\N^{-\beta(\mathrm{j},\mathrm{i}+1)}<|\mathds{A}^{\mathds{Q}[+],\pm}[\s,\y(\s)]|\leq\N^{-\beta(\mathrm{j},\mathrm{i})}\}.\label{eq:bg2136I5d}
\end{align}
\eqref{eq:bg2136I5a} follows because the indicator in \eqref{eq:bg2136I3b} does not depend on the inner summation variable $\mathds{Q}$ therein. (It depends only on the outer sum variable $\mathds{Q}[+]$; see Definition \ref{definition:bg2133}.) \eqref{eq:bg2136I5b} holds by definition of $\mathds{A}^{\mathds{Q}[+],\pm}[\s,\y(\s)]$ as the inner sum in $\mathrm{RHS}\eqref{eq:bg2136I5a}$. (Again, see Definition \ref{definition:bg2133}; namely \eqref{eq:bg2133I} and \eqref{eq:bg2133IIa}.) \eqref{eq:bg2136I5c}-\eqref{eq:bg2136I5d} follows because the indicators therein sum to the indicator in \eqref{eq:bg2136I5b}. We now claim the following, which we justify afterwards:
\begin{align}
\eqref{eq:bg2136I5c} \ = \ {\textstyle\sum^{\mathrm{j},\mathrm{i}+2}\sum^{\mathrm{i}+1,+}}\mathds{A}^{\mathds{Q},\pm}[\s,\y(\s)]\mathbf{1}\{|\mathds{A}^{\mathds{Q},\pm}[\s,\y(\s)]|\leq\N^{-\beta(\mathrm{j},\mathrm{i}+1)}\} \ = \ \eqref{eq:bg2136Ib}. \label{eq:bg2136I6}
\end{align}
The last identity in \eqref{eq:bg2136I6} is easy to check. Let us explain the first identity in \eqref{eq:bg2136I6}. Note \eqref{eq:bg2136I5c} equals the following. Tile the big block $\mathds{Q}[\mathrm{big}]$ using mutually disjoint shifts of $\mathds{Q}[\mathrm{j},\mathrm{i}+1]$. Assign a value (given by the summands in \eqref{eq:bg2136I5c}) to each copy of $\mathds{Q}[\mathrm{j},\mathrm{i}+1]$, and then average these values. This is exactly the same as the following. Tile the big block $\mathds{Q}[\mathrm{big}]$ with mutually disjoint shifts of $\mathds{Q}[\mathrm{j},\mathrm{i}+2]$. Tile each copy of $\mathds{Q}[\mathrm{j},\mathrm{i}+2]$ by mutually disjoint shifts of $\mathds{Q}[\mathrm{j},\mathrm{i}+1]$. Attach the same value to each copy of $\mathds{Q}[\mathrm{j},\mathrm{i}+1]$ that we did in our unfolding of \eqref{eq:bg2136I5c} in this paragraph. Average these values over all copies of $\mathds{Q}[\mathrm{j},\mathrm{i}+1]$ per copy of $\mathds{Q}[\mathrm{j},\mathrm{i}+2]$. This assigns a value to each copy of $\mathds{Q}[\mathrm{j},\mathrm{i}+2]$. Average all of these values over all copies of $\mathds{Q}[\mathrm{j},\mathrm{i}+2]$. (This is basically the same reasoning from Remark \ref{remark:bg2134} that gave us \eqref{eq:bg2135I1}.) Recalling notation of Definition \ref{definition:bg2133}, the equivalence of these two averaging mechanisms is exactly what the first identity in \eqref{eq:bg2136I6} says. By the same token, we also have
\begin{align}
\eqref{eq:bg2136I5d} \ = \ {\textstyle\sum^{\mathrm{j},\mathrm{i}+2}\sum^{\mathrm{i}+1,+}}\mathds{A}^{\mathds{Q},\pm}[\s,\y(\s)]\mathbf{1}\{\N^{-\beta(\mathrm{j},\mathrm{i}+1)}<|\mathds{A}^{\mathds{Q},\pm}[\s,\y(\s)]|\leq\N^{-\beta(\mathrm{j},\mathrm{i})}\} \ = \ \eqref{eq:bg2136Ie} \nonumber
\end{align}
\eqref{eq:bg2136I4} follows by \eqref{eq:bg2136I5a}-\eqref{eq:bg2136I5d}, \eqref{eq:bg2136I6}, {and the previous display}. As noted before \eqref{eq:bg2136I4}, we are done. \qed
\subsection{Proof of Lemma \ref{lemma:bg2138}}
We start with the following calculation, which we explain afterwards:
\begin{align}
\mathds{A}^{\mathfrak{m}(\mathrm{j}),\tau(\mathrm{j}),\pm}&(\mathds{R}^{\chi,\mathfrak{q},\pm,\mathrm{j}}\mathbf{Z};\s,\y(\s)) \ = \ {\textstyle\sum^{\mathrm{j},2}\sum^{1,+}}\mathds{A}^{\mathds{Q},\pm}[\s,\y(\s)]\mathbf{1}\{|\mathds{A}^{\mathds{Q},\pm}[\s,\y(\s)]|\leq\N^{-\beta(\mathrm{j},1)}\} \label{eq:bg2138I1a}\\
+ \ &{\textstyle\sum^{\mathrm{j},2}\sum^{1,+}}\mathds{A}^{\mathds{Q},\pm}[\s,\y(\s)]\mathbf{1}\{|\mathds{A}^{\mathds{Q},\pm}[\s,\y(\s)]|>\N^{-\beta(\mathrm{j},1)}\} \label{eq:bg2138I1b}\\
= \ & {\textstyle\sum^{\mathrm{j},\mathrm{i}(\mathrm{j})}\sum^{\mathrm{i}(\mathrm{j})-1,+}}\mathds{A}^{\mathds{Q},\pm}[\s,\y(\s)]\mathbf{1}\{|\mathds{A}^{\mathds{Q},\pm}[\s,\y(\s)]|\leq\N^{-\beta(\mathrm{j},\mathrm{i}(\mathrm{j})-1)}\} \label{eq:bg2138I1c}\\
+ \ &{\textstyle\sum_{\mathrm{i}=1}^{\mathrm{i}(\mathrm{j})-2}}{\textstyle\sum^{\mathrm{j},\mathrm{i}+1}\sum^{\mathrm{i},+}}\mathds{A}^{\mathds{Q},\pm}[\s,\y(\s)]\mathbf{1}\{|\mathds{A}^{\mathds{Q},\pm}[\s,\y(\s)]|\leq\N^{-\beta(\mathrm{j},\mathrm{i})}<|\mathds{S}^{\mathds{Q}[+],\pm}[\s,\y(\s)]|\}\label{eq:bg2138I1e}\\
- \ &{\textstyle\sum_{\mathrm{i}=1}^{\mathrm{i}(\mathrm{j})-2}}{\textstyle\sum^{\mathrm{j},\mathrm{i}+1}\sum^{\mathrm{i},+}}\mathds{A}^{\mathds{Q},\pm}[\s,\y(\s)]\mathbf{1}\{|\mathds{A}^{\mathds{Q}[+],\pm}[\s,\y(\s)]|\leq\N^{-\beta(\mathrm{j},\mathrm{i})}<|\mathds{S}^{\mathds{Q}[+],\pm}[\s,\y(\s)]|\}\label{eq:bg2138I1d}\\
+ \ &{\textstyle\sum_{\mathrm{i}=1}^{\mathrm{i}(\mathrm{j})-2}}{\textstyle\sum^{\mathrm{j},\mathrm{i}+2}\sum^{\mathrm{i}+1,+}}\mathds{A}^{\mathds{Q},\pm}[\s,\y(\s)]\mathbf{1}\{\N^{-\beta(\mathrm{j},\mathrm{i}+1)}<|\mathds{A}^{\mathds{Q},\pm}[\s,\y(\s)]|\leq\N^{-\beta(\mathrm{j},\mathrm{i})}\}\label{eq:bg2138I1f}\\
+ \ &{\textstyle\sum^{\mathrm{j},2}\sum^{1,+}}\mathds{A}^{\mathds{Q},\pm}[\s,\y(\s)]\mathbf{1}\{|\mathds{A}^{\mathds{Q},\pm}[\s,\y(\s)]|>\N^{-\beta(\mathrm{j},1)}\}. \label{eq:bg2138I1g}
\end{align}
\eqref{eq:bg2138I1a}-\eqref{eq:bg2138I1b} holds by Lemma \ref{lemma:bg2135}. \eqref{eq:bg2138I1c}-\eqref{eq:bg2138I1g} follows by applying Lemma \ref{lemma:bg2136} for all $\mathrm{i}=1,\ldots,\mathrm{i}(\mathrm{j})-2$ and summing up all of the resulting terms. (Precisely, this replaces $\mathrm{RHS}\eqref{eq:bg2138I1a}\mapsto\eqref{eq:bg2138I1c}$. The terms we must account for are \eqref{eq:bg2138I1e}-\eqref{eq:bg2138I1f}. \eqref{eq:bg2138I1g} is just \eqref{eq:bg2138I1b}.) Note $\mathscr{A}^{\mathbf{X},\mathbf{T}}\mathscr{R}^{\chi,\mathfrak{q},\pm,\mathrm{j}}$ is just a time-integrated heat operator acting on $\mathrm{LHS}\eqref{eq:bg2138I1a}$; see Definition \ref{definition:bg211}. {We now combine this with \eqref{eq:bg2138I1a}-\eqref{eq:bg2138I1g}, the triangle inequality, linearity of integration, and linearity of heat operators}. This gives
\begin{align}
&\E\|\mathscr{A}^{\mathbf{X},\mathbf{T}}\mathscr{R}^{\chi,\mathfrak{q},\pm,\mathrm{j}}\| \nonumber\\
&\lesssim \ \E\|{\textstyle\int_{\tau(\mathrm{j})}^{\t}}\mathbf{H}^{\N}(\s,\t(\N),\x)\{\N\eqref{eq:bg2138I1c}\}\d\s\|+\E\|{\textstyle\int_{\tau(\mathrm{j})}^{\t}}\mathbf{H}^{\N}(\s,\t(\N),\x)\{\N\eqref{eq:bg2138I1g}\}\d\s\|\label{eq:bg2138I2a}\\
&+ \ \E\|{\textstyle\int_{\tau(\mathrm{j})}^{\t}}\mathbf{H}^{\N}(\s,\t(\N),\x)\{\N\eqref{eq:bg2138I1e}\}\d\s\|+\E\|{\textstyle\int_{\tau(\mathrm{j})}^{\t}}\mathbf{H}^{\N}(\s,\t(\N),\x)\{\N\eqref{eq:bg2138I1f}\}\d\s\|\label{eq:bg2138I2b}\\
&+ \ \E\|{\textstyle\int_{\tau(\mathrm{j})}^{\t}}\mathbf{H}^{\N}(\s,\t(\N),\x)\{\N\eqref{eq:bg2138I1d}\}\d\s\|. \label{eq:bg2138I2c}
\end{align}
The first term on the RHS of \eqref{eq:bg2138I2a} equals the RHS of \eqref{eq:bg2138Ia}. Indeed, $\Phi^{\pm,\mathrm{j}}$ is just \eqref{eq:bg2138I1c} by construction. (The only remaining distinction between the first term on the RHS of \eqref{eq:bg2138I2a} and the RHS of \eqref{eq:bg2138Ia} is that the latter has an indicator $\mathbf{1}[\s\leq\t_{\mathrm{st}}]$ inside its $\d\s$ integral. But we are taking $\|\|$-norms, which restricts to times before $\t_{\mathrm{st}}$. So this indicator comes for free.) We now control the last term in \eqref{eq:bg2138I2a} by \eqref{eq:bg2138Ib}. First, note \eqref{eq:bg2138I1g} equals $\Psi^{\pm,\mathrm{j}}$ in the statement of Lemma \ref{lemma:bg2138}. Second, as we just explained, we can put the indicator $\mathbf{1}[\s\leq\t_{\mathrm{st}}]$ inside the $\d\s$ integral in the last term in \eqref{eq:bg2138I2a} for free. Third, observe that $\mathbf{H}^{\N}(\s,\t(\N),\x)$ is {an} integration over $\y\in\mathbb{T}(\N)$ against a kernel that is uniformly $\lesssim|\mathbb{T}(\N)|^{-1}|\t(\N)-\s|^{-1/2}\lesssim\N^{100\gamma_{\mathrm{reg}}}|\mathbb{T}(\N)|^{-1}\lesssim\N^{\beta_{\mathrm{BG}}}|\mathbb{T}(\N)|^{-1}$. (These bounds follow by {\eqref{eq:hke2} with $\mathrm{m}=0$}, and Definitions \ref{definition:method5}, \ref{definition:reg}, \ref{definition:method8}.) Using all of this, for any $\t\leq\t_{\mathrm{st}}$ and $\x\in\mathbb{T}(\N)$, we get
\begin{align}
|{\textstyle\int_{\tau(\mathrm{j})}^{\t}}\mathbf{H}^{\N}(\s,\t(\N),\x)\{\N\eqref{eq:bg2138I1g}\}\d\s| \ &\lesssim \ \N^{100\gamma_{\mathrm{reg}}}{\textstyle\int_{\tau(\mathrm{j})}^{\t}}|\mathbb{T}(\N)|^{-1}{\textstyle\sum_{\y}}\{\mathbf{1}[\s\leq\t_{\mathrm{st}}]\N|\eqref{eq:bg2138I1g}|\}\d\s \label{eq:bg2138I3a}\\
&\leq \ \N^{\beta_{\mathrm{BG}}}{\textstyle\int_{\tau(\mathrm{j})}^{1}}|\mathbb{T}(\N)|^{-1}{\textstyle\sum_{\y}}\{\mathbf{1}[\s\leq\t_{\mathrm{st}}]\N|\eqref{eq:bg2138I1g}|\}\d\s, \label{eq:bg2138I3b}
\end{align}
where the last bound follows by extending the domain of integration (of a non-negative integrand). Note \eqref{eq:bg2138I3b} is independent of the original $(\t,\x)$-variables. Thus, the last term in \eqref{eq:bg2138I2a} is $\lesssim\E\eqref{eq:bg2138I3b}$. As $\eqref{eq:bg2138I1g}=\Psi^{\pm,\mathrm{j}}(\s,\y(\s))$, we deduce the last term in \eqref{eq:bg2138I2a} is big-Oh of \eqref{eq:bg2138Ib}. By this, the previous paragraph, and \eqref{eq:bg2138I2a}-\eqref{eq:bg2138I2c}, to show \eqref{eq:bg2138Ia}-\eqref{eq:bg2138Id}, it suffices to show
\begin{align}
\eqref{eq:bg2138I2b}+\eqref{eq:bg2138I2c} \ \lesssim \ \eqref{eq:bg2138Ic}+\eqref{eq:bg2138Id}. \label{eq:bg2138I4}
\end{align}
We first tackle \eqref{eq:bg2138I2b}, beginning with the first term therein. We claim the following calculation holds, which we explain after:
\begin{align}
|\eqref{eq:bg2138I1e}| \ &\lesssim \ {\textstyle\sum_{\mathrm{i}=1}^{\mathrm{i}(\mathrm{j})-2}}{\textstyle\sum^{\mathrm{j},\mathrm{i}+1}\sum^{\mathrm{i},+}}\N^{\beta(\mathrm{j},\mathrm{i}-1)}|\mathds{A}^{\mathds{Q},\pm}[\s,\y(\s)]|^{2}\mathbf{1}\{|\mathds{A}^{\mathds{Q},\pm}[\s,\y(\s)]|\leq\N^{-\beta(\mathrm{j},\mathrm{i})}\} \label{eq:bg2138I5a}\\
&+ \ {\textstyle\sum_{\mathrm{i}=1}^{\mathrm{i}(\mathrm{j})-2}}{\textstyle\sum^{\mathrm{j},\mathrm{i}+1}\sum^{\mathrm{i},+}}\N^{-\beta(\mathrm{j},\mathrm{i}-1)}\mathbf{1}\{|\mathds{S}^{\mathds{Q}[+],\pm}[\s,\y(\s)]|>\N^{-\beta(\mathrm{j},\mathrm{i})}\}\label{eq:bg2138I5b}\\
&\lesssim \ {\textstyle\sum_{\mathrm{i}=1}^{\mathrm{i}(\mathrm{j})-2}}{\textstyle\sum^{\mathrm{j},\mathrm{i}+1}\sum^{\mathrm{i},+}}\N^{\beta(\mathrm{j},\mathrm{i}-1)}|\mathds{A}^{\mathds{Q},\pm}[\s,\y(\s)]|^{2}\mathbf{1}\{|\mathds{A}^{\mathds{Q},\pm}[\s,\y(\s)]|\leq\N^{-\beta(\mathrm{j},\mathrm{i}-1)}\} \label{eq:bg2138I5c}\\
&+ \ {\textstyle\sum_{\mathrm{i}=1}^{\mathrm{i}(\mathrm{j})-2}}{\textstyle\sum^{\mathrm{j},\mathrm{i}+1}}\N^{-\beta(\mathrm{j},\mathrm{i}-1)}\mathbf{1}\{|\mathds{S}^{\mathds{Q}[+],\pm}[\s,\y(\s)]|>\N^{-\beta(\mathrm{j},\mathrm{i})}\}\label{eq:bg2138I5d}\\
&\lesssim \ {\textstyle\sum_{\mathrm{i}=1}^{\mathrm{i}(\mathrm{j})-2}}{\textstyle\sum^{\mathrm{j},\mathrm{i}+1}\sum^{\mathrm{i},+}}\N^{\beta(\mathrm{j},\mathrm{i}-1)}|\mathds{A}^{\mathds{Q},\pm}[\s,\y(\s)]|^{2}\mathbf{1}\{|\mathds{A}^{\mathds{Q},\pm}[\s,\y(\s)]|\leq\N^{-\beta(\mathrm{j},\mathrm{i}-1)}\} \label{eq:bg2138I5e}\\
&+ \ \N^{{\mathrm{O}(\gamma_{\mathrm{ap}})}}{\textstyle\sum_{\mathrm{i}=1}^{\mathrm{i}(\mathrm{j})-2}}{\textstyle\sum^{\mathrm{j},\mathrm{i}+1}\sum^{\mathrm{i},+}}\N^{-\beta(\mathrm{j},\mathrm{i}-1)}\mathbf{1}\{|\mathds{A}^{\mathds{Q},\pm}[\s,\y(\s)]|>\N^{-\beta(\mathrm{j},\mathrm{i})}\} \label{eq:bg2138I5f}\\
&= \ {\textstyle\sum_{\mathrm{i}=1}^{\mathrm{i}(\mathrm{j})-2}}\Upsilon^{\pm,\mathrm{j},\mathrm{i}}(\s,\y(\s))+{\textstyle\sum_{\mathrm{i}=1}^{\mathrm{i}(\mathrm{j})-2}}\N^{{\mathrm{O}(\gamma_{\mathrm{ap}})}}\Lambda^{\pm,\mathrm{j},\mathrm{i}}(\s,\y(\s)). \label{eq:bg2138I5g}
\end{align}
\eqref{eq:bg2138I5a}-\eqref{eq:bg2138I5b} follows by the Schwarz inequality applied to each summand in \eqref{eq:bg2138I1e}. \eqref{eq:bg2138I5c}-\eqref{eq:bg2138I5d} follows by relaxing the constraint in the indicator in $\mathrm{RHS}\eqref{eq:bg2138I5a}$. Then, we observe that the summands in \eqref{eq:bg2138I5b} do not depend on the inner summation variable $\mathds{Q}$, only the outer one $\mathds{Q}[+]$. Thus, we can remove the inner average in \eqref{eq:bg2138I5b} to get \eqref{eq:bg2138I5d}. \eqref{eq:bg2138I5e}-\eqref{eq:bg2138I5f} follows by first leaving \eqref{eq:bg2138I5c} alone to get \eqref{eq:bg2138I5e}. Next, recall $\mathds{S}$ from Definition \ref{definition:bg2133} as a supremum of $|\mathds{A}^{\mathds{Q},\pm}|$-terms over $\mathds{Q}\in\mathscr{Q}[\mathds{Q}[+]]$. So, if $|\mathds{S}^{\mathds{Q}[+],\pm}|>\N^{-\beta(\mathrm{j},\mathrm{i})}$, then $|\mathds{A}^{\mathds{Q},\pm}|>\N^{-\beta(\mathrm{j},\mathrm{i})}$ for at least one $\mathds{Q}\in\mathscr{Q}[\mathds{Q}[+]]$. By a union bound, we may therefore control the indicator in \eqref{eq:bg2138I5d} by a sum of indicators of $|\mathds{A}^{\mathds{Q},\pm}|>\N^{-\beta(\mathrm{j},\mathrm{i})}$ over all $\mathds{Q}\in\mathscr{Q}[\mathds{Q}[+]]$. This is exactly what the last line says to do to \eqref{eq:bg2138I5d}.  (The extra factor $\N^{{\mathrm{O}(\gamma_{\mathrm{ap}})}}$ comes from the fact that the inner-most sum in \eqref{eq:bg2138I5f} averages over $\mathds{Q}\in\mathscr{Q}[\mathds{Q}[+]]$, not sum. But we know $|\mathscr{Q}[\mathds{Q}[+]]|\lesssim\N^{{\mathrm{O}(\gamma_{\mathrm{ap}})}}$, because the dimensions of $\mathds{Q}[+]$ are small powers of $\N^{\gamma_{\mathrm{ap}}}$ bigger than dimensions of $\mathds{Q}$; see Definition \ref{definition:bg2131}. (Indeed, recall from Definition \ref{definition:bg2133} that $|\mathscr{Q}[\mathds{Q}[+]]|$ is the number of mutually disjoint shifts of $\mathds{Q}[\mathrm{j},\mathrm{i}]$ needed to cover $\mathds{Q}[\mathrm{j},\mathrm{i}+1]$. As we have just explained, by Definition \ref{definition:bg2131}, this is big-Oh of some small power of $\N^{\gamma_{\mathrm{ap}}}$.) Therefore, the extra factor of $\N^{{\mathrm{O}(\gamma_{\mathrm{ap}})}}$ is more than enough to turn the inner-most sum in \eqref{eq:bg2138I5f} from an average to an honest sum.) \eqref{eq:bg2138I5g} holds by \eqref{eq:bg2138IIc}-\eqref{eq:bg2138IId}. By \eqref{eq:bg2138I5a}-\eqref{eq:bg2138I5g} and {the} reasoning for \eqref{eq:bg2138I3a}-\eqref{eq:bg2138I3b}, we get an integrated heat-operator bound:
\begin{align}
&|{\textstyle\int_{\tau(\mathrm{j})}^{\t}}\mathbf{H}^{\N}(\s,\t(\N),\x)\{\N\eqref{eq:bg2138I1e}\}\d\s| \nonumber\\
&\lesssim \ \N^{{\mathrm{O}(\gamma_{\mathrm{reg}})}}{\textstyle\int_{\tau(\mathrm{j})}^{\t}}|\mathbb{T}(\N)|^{-1}{\textstyle\sum_{\y}}\{\mathbf{1}[\s\leq\t_{\mathrm{st}}]\N\eqref{eq:bg2138I1e}\}\d\s \label{eq:bg2138I6a}\\
&\lesssim \ {\textstyle\sum_{\mathrm{i}=1}^{\mathrm{i}(\mathrm{j})-2}}\N^{{\mathrm{O}(\gamma_{\mathrm{reg}})}}{\textstyle\int_{\tau(\mathrm{j})}^{\t}}|\mathbb{T}(\N)|^{-1}{\textstyle\sum_{\y}}\{\mathbf{1}[\s\leq\t_{\mathrm{st}}]\N|\Upsilon^{\pm,\mathrm{j},\mathrm{i}}(\s,\y(\s))|\}\d\s \label{eq:bg2138I6b}\\
&+ \ {\textstyle\sum_{\mathrm{i}=1}^{\mathrm{i}(\mathrm{j})-2}}\N^{{\mathrm{O}(\gamma_{\mathrm{reg}})}+{\mathrm{O}(\gamma_{\mathrm{ap}})}}{\textstyle\int_{\tau(\mathrm{j})}^{\t}}|\mathbb{T}(\N)|^{-1}{\textstyle\sum_{\y}}\{\mathbf{1}[\s\leq\t_{\mathrm{st}}]\N|\Lambda^{\pm,\mathrm{j},\mathrm{i}}(\s,\y(\s))|\}\d\s \label{eq:bg2138I6c} \\
&\lesssim \ {\textstyle\sum_{\mathrm{i}=1}^{\mathrm{i}(\mathrm{j})-2}}\N^{\beta_{\mathrm{BG}}}{\textstyle\int_{\tau(\mathrm{j})}^{\t}}|\mathbb{T}(\N)|^{-1}{\textstyle\sum_{\y}}\{\mathbf{1}[\s\leq\t_{\mathrm{st}}]\N|\Upsilon^{\pm,\mathrm{j},\mathrm{i}}(\s,\y(\s))|\}\d\s \label{eq:bg2138I6d}\\
&+ \ {\textstyle\sum_{\mathrm{i}=1}^{\mathrm{i}(\mathrm{j})-2}}\N^{\beta_{\mathrm{BG}}}{\textstyle\int_{\tau(\mathrm{j})}^{\t}}|\mathbb{T}(\N)|^{-1}{\textstyle\sum_{\y}}\{\mathbf{1}[\s\leq\t_{\mathrm{st}}]\N|\Lambda^{\pm,\mathrm{j},\mathrm{i}}(\s,\y(\s))|\}\d\s. \label{eq:bg2138I6e}
\end{align}
(Above, we used {that $\gamma_{\mathrm{reg}},\gamma_{\mathrm{ap}}$ are small factors times $\beta_{\mathrm{BG}}$}; see Definitions \ref{definition:reg}, \ref{definition:method8}.) As \eqref{eq:bg2138I6d}-\eqref{eq:bg2138I6e} are independent of the $(\t,\x)$-variables in $\mathrm{LHS}\eqref{eq:bg2138I6a}$, these last two lines also bound the norm inside of $\E$ in the first term in \eqref{eq:bg2138I2b}. We also note that $\mathrm{i}(\mathrm{j})\lesssim1$, since it is at most the number of steps of size $\gtrsim1$ to go from $0\mapsto1$; see Definition \ref{definition:bg2131}. Therefore, we deduce the following:
\begin{align}
&\E\|{\textstyle\int_{\tau(\mathrm{j})}^{\t}}\mathbf{H}^{\N}(\s,\t(\N),\x)\{\N\eqref{eq:bg2138I1e}\}\d\s\| \nonumber\\
&\lesssim \ {\textstyle\sum_{\mathrm{i}=1}^{\mathrm{i}(\mathrm{j})-2}}\N^{\beta_{\mathrm{BG}}}{\textstyle\int_{\tau(\mathrm{j})}^{\t}}|\mathbb{T}(\N)|^{-1}{\textstyle\sum_{\y}}\E\{\mathbf{1}[\s\leq\t_{\mathrm{st}}]\N|\Upsilon^{\pm,\mathrm{j},\mathrm{i}}(\s,\y(\s))|\}\d\s \label{eq:bg2138I7a}\\
&+ \ {\textstyle\sum_{\mathrm{i}=1}^{\mathrm{i}(\mathrm{j})-2}}\N^{\beta_{\mathrm{BG}}}{\textstyle\int_{\tau(\mathrm{j})}^{\t}}|\mathbb{T}(\N)|^{-1}{\textstyle\sum_{\y}}\E\{\mathbf{1}[\s\leq\t_{\mathrm{st}}]\N|\Lambda^{\pm,\mathrm{j},\mathrm{i}}(\s,\y(\s))|\}\d\s \label{eq:bg2138I7b}\\
&\lesssim \ \eqref{eq:bg2138Ic}+\eqref{eq:bg2138Id}. \label{eq:bg2138I7c}
\end{align}
We now treat the last term in \eqref{eq:bg2138I2b}. This uses basically the same argument. By the reasoning giving \eqref{eq:bg2138I5a}-\eqref{eq:bg2138I5g},
\begin{align}
|\eqref{eq:bg2138I1f}| \ &\lesssim \ {\textstyle\sum_{\mathrm{i}=1}^{\mathrm{i}(\mathrm{j})-2}}{\textstyle\sum^{\mathrm{j},\mathrm{i}+2}\sum^{\mathrm{i}+1,+}}\N^{\beta(\mathrm{j},\mathrm{i})}|\mathds{A}^{\mathds{Q},\pm}[\s,\y(\s)]|^{2}\mathbf{1}\{|\mathds{A}^{\mathds{Q},\pm}[\s,\y(\s)]|\leq\N^{-\beta(\mathrm{j},\mathrm{i})}\} \label{eq:bg2138I8a}\\
&+ \ {\textstyle\sum_{\mathrm{i}=1}^{\mathrm{i}(\mathrm{j})-2}}{\textstyle\sum^{\mathrm{j},\mathrm{i}+2}\sum^{\mathrm{i}+1,+}}\N^{-\beta(\mathrm{j},\mathrm{i})}\mathbf{1}\{|\mathds{A}^{\mathds{Q},\pm}[\s,\y(\s)]|>\N^{-\beta(\mathrm{j},\mathrm{i}+1)}\} \label{eq:bg2138I8b}\\
&\lesssim \ {\textstyle\sum_{\mathrm{i}=1}^{\mathrm{i}(\mathrm{j})-2}}\Upsilon^{\pm,\mathrm{j},\mathrm{i}+1}(\s,\y(\s))+{\textstyle\sum_{\mathrm{i}=1}^{\mathrm{i}(\mathrm{j})-2}}\Lambda^{\pm,\mathrm{j},\mathrm{i}+1}(\s,\y(\s)). \label{eq:bg2138I8c}
\end{align}
(We clarify {that} \eqref{eq:bg2138I8c} follows with the indices $\mathrm{i}+1$ because the summations in \eqref{eq:bg2138I8a}-\eqref{eq:bg2138I8b} have indices $\mathrm{i}+2,\mathrm{i}+1$, which are one more than the sum-indices in \eqref{eq:bg2138I5a}-\eqref{eq:bg2138I5g}.) We can now use \eqref{eq:bg2138I8a}-\eqref{eq:bg2138I8c} in the same way that we used \eqref{eq:bg2138I5a}-\eqref{eq:bg2138I5g} to get \eqref{eq:bg2138I7a}-\eqref{eq:bg2138I7c}. This gives \eqref{eq:bg2138I7a}-\eqref{eq:bg2138I7c} but with \eqref{eq:bg2138I1f} instead of \eqref{eq:bg2138I1e} in $\mathrm{LHS}\eqref{eq:bg2138I7a}$. By this and \eqref{eq:bg2138I7a}-\eqref{eq:bg2138I7c},
\begin{align}
\eqref{eq:bg2138I2b} \ \lesssim \ \eqref{eq:bg2138Ic}+\eqref{eq:bg2138Id}. \label{eq:bg2138I9}
\end{align}
We now tackle \eqref{eq:bg2138I2c}. To this end, we claim that the following calculation holds:
\begin{align}
&\eqref{eq:bg2138I1d} \nonumber\\
&= \ {\textstyle\sum_{\mathrm{i}=1}^{\mathrm{i}(\mathrm{j})-2}}{\textstyle\sum^{\mathrm{j},\mathrm{i}+1}}\mathbf{1}\{|\mathds{A}^{\mathds{Q}[+],\pm}[\s,\y(\s)]|\leq\N^{-\beta(\mathrm{j},\mathrm{i})}<|\mathds{S}^{\mathds{Q}[+],\pm}[\s,\y(\s)]|\}
{\textstyle\sum^{\mathrm{i},+}}\mathds{A}^{\mathds{Q},\pm}[\s,\y(\s)] \label{eq:bg2138I10a}\\
&= \ {\textstyle\sum_{\mathrm{i}=1}^{\mathrm{i}(\mathrm{j})-2}}{\textstyle\sum^{\mathrm{j},\mathrm{i}+1}}\mathds{A}^{\mathds{Q}[+],\pm}[\s,\y(\s)]\mathbf{1}\{|\mathds{A}^{\mathds{Q}[+],\pm}[\s,\y(\s)]|\leq\N^{-\beta(\mathrm{j},\mathrm{i})}<|\mathds{S}^{\mathds{Q}[+],\pm}[\s,\y(\s)]|\} \label{eq:bg2138I10b}\\
&\lesssim \ {\textstyle\sum_{\mathrm{i}=1}^{\mathrm{i}(\mathrm{j})-2}}{\textstyle\sum^{\mathrm{j},\mathrm{i}+1}}\N^{\beta(\mathrm{j},\mathrm{i})}|\mathds{A}^{\mathds{Q}[+],\pm}[\s,\y(\s)]|^{2}\mathbf{1}\{|\mathds{A}^{\mathds{Q}[+],\pm}[\s,\y(\s)]|\leq\N^{-\beta(\mathrm{j},\mathrm{i})}\}\label{eq:bg2138I10c}\\
&+ \ {\textstyle\sum_{\mathrm{i}=1}^{\mathrm{i}(\mathrm{j})-2}}{\textstyle\sum^{\mathrm{j},\mathrm{i}+1}}\N^{-\beta(\mathrm{j},\mathrm{i})}\mathbf{1}\{|\mathds{S}^{\mathds{Q}[+],\pm}[\s,\y(\s)]|>\N^{-\beta(\mathrm{j},\mathrm{i})}\} \label{eq:bg2138I10d} \\
&= \ {\textstyle\sum_{\mathrm{i}=1}^{\mathrm{i}(\mathrm{j})-2}}{\textstyle\sum^{\mathrm{j},\mathrm{i}+2}\sum^{\mathrm{i}+1,+}}\N^{\beta(\mathrm{j},\mathrm{i})}|\mathds{A}^{\mathds{Q},\pm}[\s,\y(\s)]|^{2}\mathbf{1}\{|\mathds{A}^{\mathds{Q},\pm}[\s,\y(\s)]|\leq\N^{-\beta(\mathrm{j},\mathrm{i})}\}\label{eq:bg2138I10e}\\
&+ \ {\textstyle\sum_{\mathrm{i}=1}^{\mathrm{i}(\mathrm{j})-2}}{\textstyle\sum^{\mathrm{j},\mathrm{i}+1}}\N^{-\beta(\mathrm{j},\mathrm{i})}\mathbf{1}\{|\mathds{S}^{\mathds{Q}[+],\pm}[\s,\y(\s)]|>\N^{-\beta(\mathrm{j},\mathrm{i})}\} \label{eq:bg2138I10f} \\
&\lesssim \ {\textstyle\sum_{\mathrm{i}=1}^{\mathrm{i}(\mathrm{j})-2}}{\textstyle\sum^{\mathrm{j},\mathrm{i}+2}\sum^{\mathrm{i}+1,+}}\N^{\beta(\mathrm{j},\mathrm{i})}|\mathds{A}^{\mathds{Q},\pm}[\s,\y(\s)]|^{2}\mathbf{1}\{|\mathds{A}^{\mathds{Q},\pm}[\s,\y(\s)]|\leq\N^{-\beta(\mathrm{j},\mathrm{i})}\}\label{eq:bg2138I10g}\\
&+ \ \N^{{\mathrm{O}(\gamma_{\mathrm{ap}})}}{\textstyle\sum_{\mathrm{i}=1}^{\mathrm{i}(\mathrm{j})-2}}{\textstyle\sum^{\mathrm{j},\mathrm{i}+1}\sum^{\mathrm{i},+}}\N^{-\beta(\mathrm{j},\mathrm{i}-1)}\mathbf{1}\{|\mathds{A}^{\mathds{Q},\pm}[\s,\y(\s)]|>\N^{-\beta(\mathrm{j},\mathrm{i})}\} \label{eq:bg2138I10h}\\
&\lesssim \ {\textstyle\sum_{\mathrm{i}=1}^{\mathrm{i}(\mathrm{j})-2}}\Upsilon^{\pm,\mathrm{j},\mathrm{i}+1}(\s,\y(\s))+{\textstyle\sum_{\mathrm{i}=1}^{\mathrm{i}(\mathrm{j})-2}}\N^{{\mathrm{O}(\gamma_{\mathrm{ap}})}}\Lambda^{\pm,\mathrm{j},\mathrm{i}}(\s,\y(\s)). \label{eq:bg2138I10i}
\end{align}
\eqref{eq:bg2138I10a} holds since the indicator in \eqref{eq:bg2138I1d} does not depend on the inner-sum variable $\mathds{Q}$. So, we pull it outside this sum. \eqref{eq:bg2138I10b} holds by gluing $\mathds{Q}$-averages into the $\mathds{Q}[+]$-average; see \eqref{eq:bg2133IIa}. \eqref{eq:bg2138I10c}-\eqref{eq:bg2138I10d} follows via Schwarz inequality. \eqref{eq:bg2138I10e}-\eqref{eq:bg2138I10f} follows by first leaving alone \eqref{eq:bg2138I10d} to obtain \eqref{eq:bg2138I10f}. Next, in \eqref{eq:bg2138I10c}, reparameterize the tiling by shifts of $\mathds{Q}[+]$ into two-scale tilings, giving the double sum in \eqref{eq:bg2138I10e}. This is the reasoning that gave us \eqref{eq:bg2136I6}. \eqref{eq:bg2138I10g}-\eqref{eq:bg2138I10h} holds by leaving alone \eqref{eq:bg2138I10e} to get \eqref{eq:bg2138I10g}. Next, bound \eqref{eq:bg2138I10f} with the calculation starting at \eqref{eq:bg2138I5d}. (This is an argument via union bound. We clarify the change in exponent $\beta(\mathrm{j},\mathrm{i})\mapsto\beta(\mathrm{j},\mathrm{i}-1)$ is harmless and only makes things bigger, so it is allowed.) \eqref{eq:bg2138I10i} holds by construction; see \eqref{eq:bg2138IIc}-\eqref{eq:bg2138IId}. Now, we use \eqref{eq:bg2138I10a}-\eqref{eq:bg2138I10i} the same way we used \eqref{eq:bg2138I5a}-\eqref{eq:bg2138I5g} to get  \eqref{eq:bg2138I7a}-\eqref{eq:bg2138I7c}. This gives \eqref{eq:bg2138I7a}-\eqref{eq:bg2138I7c} but replacing \eqref{eq:bg2138I1e} by \eqref{eq:bg2138I1d}. In particular, this implies $\eqref{eq:bg2138I2c}\lesssim\eqref{eq:bg2138Ic}+\eqref{eq:bg2138Id}$. Combining this with \eqref{eq:bg2138I9} gives \eqref{eq:bg2138I4}. As noted right before \eqref{eq:bg2138I4}, this completes the proof. \qed
%
%
%
\section{Proof of Proposition \ref{prop:bg21310}}
This section is organized as follows. We will first estimate $\mathrm{RHS}\eqref{eq:bg2138Ia}$. This is an elementary, direct bound. Next, we bound \eqref{eq:bg2138Ic}; along the way, we gather ingredients for \eqref{eq:bg2138Id}. This uses local equilibrium and Kipnis-Varadhan bounds that we spent previous sections deriving. We finish with \eqref{eq:bg2138Ib}. This requires rather minor modifications of our analysis for \eqref{eq:bg2138Ic}-\eqref{eq:bg2138Id}. Before we start, recall from the end of Section \ref{section:main} what $\inf$ and $\sup$ of a discrete interval mean. We use these notions frequently.

Before we start, however, let us give an explanation for why the bounds in Proposition \ref{prop:bg21310} hold. (The point of this section is to make these bounds rigorous. Unfortunately, this section is quite detailed and technical by nature of the argument, hence the motivation for a more-than-intuitive but still simple explanation below.)
\begin{itemize}
\item Take the RHS of \eqref{eq:bg2138Ia}. It averages $\N\Phi^{\pm,\mathrm{j}}$, and, looking at Lemma \ref{lemma:bg2138} for notation, we know $\Phi^{\pm,\mathrm{j}}$ is an average of terms that are $\lesssim\N^{-\beta(\mathrm{j},\mathrm{i}(\mathrm{j})-1)}$, which is $\ll\N^{-1}$ by construction in Definition \ref{definition:bg2131}. So, we can bound this term essentially by definition.
\item Now take \eqref{eq:bg2138Id}, namely the $\Lambda^{\pm,\mathrm{j}}$ integrand therein. By Chebyshev, we can effectively bound it by $\N^{\beta(\mathrm{j},\mathrm{i}-1)}$ times the square of a space-time average of something that is roughly a priori $\lesssim\mathfrak{l}(\mathrm{j})^{-3/2}$ on time-scale $\tau(\mathrm{j},\mathrm{i})$ and length-scale $\mathfrak{m}(\mathrm{j},\mathrm{i})$. So, by Proposition \ref{prop:kv1}, we can bound the squared space-time average by $\lesssim\N^{-2}\tau(\mathrm{j},\mathrm{i})^{-1}\mathfrak{m}(\mathrm{j},\mathrm{i})^{-1}\mathfrak{l}(\mathrm{j})^{-1}$ (the extra factor of $\mathfrak{l}(\mathrm{j})^{2}$ comes from the support length of the functional we are averaging). After multiplying by $\N^{\beta(\mathrm{j},\mathrm{i}-1)}$, by construction of the space-time scales and exponents in Definition \ref{definition:bg2131}, the resulting bound is $\ll\N^{-1}$. 

We clarify this explanation is only rigorous if we knew the law of the system (at least locally as far as the space-time average is concerned) is a canonical measure. We eventually use Lemma \ref{lemma:le9} to reduce to {the} canonical measure; for this, it is crucial that we have the a priori bound of $\N^{-\beta(\mathrm{j},\mathrm{i}-1)}$ on the $\Lambda^{\pm,\mathrm{j}}$ integrand in \eqref{eq:bg2138Id}. (Indeed, this would let us take $\kappa$ bigger in Lemma \ref{lemma:le9} since {we would be able} to control the $\mathfrak{a}$-functional therein, which will eventually be the squared space-time average from above, at sharper exponential scales, a priori.) The fact that the cost in reduction to local canonical measure (the first term on the RHS of \eqref{eq:le9I}) is $\ll\N^{-1}$ is ultimately a power-counting, but it can be intuited as follows. Yau's relative entropy method works for hyperbolic fluctuations (i.e. dropping the second-order term in \eqref{eq:glsde}) if we have initial relative entropy of order $\ll\N^{1/2}$, and it works for parabolic fluctuations (i.e. dropping the first-order term in \eqref{eq:glsde}) if we have initial relative entropy of order $\N$. Since we are somewhere in the middle, assuming initial relative entropy of order $\ll\N^{3/4}$ should be enough. (We clarify that there is no issue of time-inhomogeneity here; that was all dealt with in the proof of the Kipnis-Varadhan bound of Proposition \ref{prop:kv1}.)
\item Now take \eqref{eq:bg2138Ic}. The same argument works, except instead of using the Kipnis-Varadhan bound \eqref{eq:kv1III}, use \eqref{eq:kv1II} instead. (Indeed, the only other ingredient we needed for analyzing \eqref{eq:bg2138Id} was that its integrand $\Lambda^{\pm,\mathrm{j}}$ {has a priori estimate of order $\N^{-\beta(\mathrm{j},\mathrm{i}-1)}$}. This is clearly true for the $\Upsilon^{\pm,\mathrm{j}}$ integrand in \eqref{eq:bg2138Ic}.)
\item As for \eqref{eq:bg2138Ib}, which is all we have left, it is essentially the same thing as \eqref{eq:bg2138Id}, except we do not have the extra power-saving of $\N^{-\beta(\mathrm{j},\mathrm{i}-1)}$ in the $\Psi^{\pm,\mathrm{j}}$ integrand in \eqref{eq:bg2138Ib}. So, we are missing a helpful factor of $\N^{\beta(\mathrm{j},1)}$. But this factor is small enough, since $\beta(\mathrm{j},1)$ is the first and smallest exponent in the sequences constructed (per $\mathrm{j}$) in Definition \ref{definition:bg2131}. (We effectively chose it to be as small as we want in Definition \ref{definition:bg2131}.) Thus, this helpful factor of $\N^{\beta(\mathrm{j},1)}$ actually has no important role.
\end{itemize}
\subsection{Bound for $\mathrm{RHS}\eqref{eq:bg2138Ia}$}
By construction in \eqref{eq:bg2138IIa}, we know $\Phi^{\pm,\mathrm{j}}$ is an average of terms that are $\lesssim\N^{-\beta(\mathrm{j},\mathrm{i}(\mathrm{j})-1)}$. Also, $\t_{\mathrm{st}}\leq1$ with probability 1; see Definitions \ref{definition:method8} and \ref{definition:reg}. Finally, the $\mathbf{H}^{\N}$ heat operator is contractive in $\mathscr{L}^{\infty}(\mathbb{T}(\N))$. Thus,
\begin{align}
\mathrm{RHS}\eqref{eq:bg2138Ia} \ \lesssim \ \N\cdot\N^{-\beta(\mathrm{j},\mathrm{i}(\mathrm{j})-1)}\E\t_{\mathrm{st}} \ \lesssim \ \N\cdot\N^{-\beta(\mathrm{j},\mathrm{i}(\mathrm{j})-1)}. \label{eq:bg2138Ia1}
\end{align}
{We know $\mathfrak{m}(\mathrm{j},\mathrm{i}(\mathrm{j})-1)\mathfrak{l}(\mathrm{j})\gtrsim\N^{-\delta_{\mathrm{ap}}}\mathfrak{m}(\mathrm{j},\mathrm{i}(\mathrm{j}))\mathfrak{l}(\mathrm{j})=\N^{\frac34+\alpha(\mathrm{j})-\delta_{\mathrm{ap}}}$ for $\delta_{\mathrm{ap}}\leq2\gamma_{\mathrm{ap}}$ and ${c}\gamma_{\mathrm{KL}}\leq\alpha(\mathrm{j})\leq2{c}\gamma_{\mathrm{KL}}$; see Definition \ref{definition:bg2131} for these bounds and Definitions \ref{definition:method8}, \ref{definition:entropydata} for $\gamma_{\mathrm{ap}},\gamma_{\mathrm{KL}}$. It can be checked via these bounds that $\mathfrak{m}(\mathrm{j},\mathrm{i}(\mathrm{j})-1)\mathfrak{l}(\mathrm{j})\geq\N^{1/2}$. Thus, by Definition \ref{definition:bg2131}, we deduce $\tau(\mathrm{j},\mathrm{i}(\mathrm{j})-1)=\N^{-3/2}\mathfrak{m}(\mathrm{j},\mathrm{i}(\mathrm{j})-1)\mathfrak{l}(\mathrm{j})$. Using all this with \eqref{eq:bg2138Ia1}, {the} construction of $\beta(\mathrm{j},\mathrm{i})$ in Definition \ref{definition:bg2131}, the bounds $\delta_{\mathrm{ap}},\gamma_{\mathrm{reg}}\leq90\beta_{\mathrm{BG}}$ (see Definitions \ref{definition:reg}, \ref{definition:method8}), and {the bound} $\alpha(\mathrm{j})\geq{\mathrm{C}}\beta_{\mathrm{BG}}$ {for some large but fixed $\mathrm{C}>0$} (see Definition \ref{definition:method8}), {we obtain}}
\begin{align}
&\mathrm{RHS}\eqref{eq:bg2138Ia} \nonumber\\
&\lesssim \ \N\cdot\N^{-1+20\gamma_{\mathrm{reg}}+90\beta_{\mathrm{BG}}}\tau(\mathrm{j},\mathrm{i}(\mathrm{j})-1)^{-1}\mathfrak{m}(\mathrm{j},\mathrm{i}(\mathrm{j})-1)^{-1}\mathfrak{l}(\mathrm{j})^{-1} \label{eq:bg2138Ia2a}\\
&\lesssim \ \N^{\frac32+20\gamma_{\mathrm{reg}}+90\beta_{\mathrm{BG}}}\mathfrak{m}(\mathrm{j},\mathrm{i}(\mathrm{j})-1)^{-2}\mathfrak{l}(\mathrm{j})^{-2} \ \lesssim \ \N^{2\delta_{\mathrm{ap}}+20\gamma_{\mathrm{reg}}+90\beta_{\mathrm{BG}}}\N^{-2\alpha(\mathrm{j})} \ \lesssim \ \N^{-90\beta_{\mathrm{BG}}}, \label{eq:bg2138Ia2b}
\end{align}
The previous display \eqref{eq:bg2138Ia2a}-\eqref{eq:bg2138Ia2b} implies the desired bound for $\mathrm{RHS}\eqref{eq:bg2138Ia}$, so we are done with this subsection.
\subsection{Bounds for \eqref{eq:bg2138Ic} and \eqref{eq:bg2138Id}}
This is the technical bulk of this section. It has several steps.
\subsubsection{Technical step (to ease notation)}
Bounding \eqref{eq:bg2138Ic}-\eqref{eq:bg2138Id} amounts to estimating macroscopic-scale space-time averages of \eqref{eq:bg2138IIc}-\eqref{eq:bg2138IId}. Each of \eqref{eq:bg2138IIc}-\eqref{eq:bg2138IId} is the average of space-time shifts of some functional. Instead of keeping track of every such shift, for each shift, we change variables in the space-time integration in \eqref{eq:bg2138Ic}-\eqref{eq:bg2138Id}. This removes the shift in each summand in \eqref{eq:bg2138IIc}-\eqref{eq:bg2138IId}; the cost that we must pay is changing the integration-domain in \eqref{eq:bg2138Ic}-\eqref{eq:bg2138Id} to something which depends on the shift. Since \eqref{eq:bg2138Ic}-\eqref{eq:bg2138Id} integrate non-negative terms, to totally forget the shifts in \eqref{eq:bg2138IIc}-\eqref{eq:bg2138IId}, we just find an integration-domain that contains each shift-dependent integration-domain. (Before we state this result, however, we make one clarifying remark. Note that although the superscript $\mathds{Q}$ in $\mathds{A}^{\mathds{Q},\pm}$ shifts the index set for averaging, it does not shift the time for the coupling constant in \eqref{eq:bg211I}; see Definition \ref{definition:bg2133}. This is a minor annoyance that we get rid of in Lemma \ref{lemma:bg213101} below.) 
\begin{lemma}\label{lemma:bg213101}
 Fix $1\leq\mathrm{j}\leq\mathrm{j}(\infty)$. With notation explained after, we have the following with probability 1:
\begin{align}
\eqref{eq:bg2138Ic} \ &\lesssim \ {\textstyle\sup_{\t}\sup_{\mathrm{i}}}\N^{\beta_{\mathrm{BG}}}{\textstyle\int_{\tau(\mathrm{j},\mathrm{i})}^{1}}|\mathbb{T}(\N)|^{-1}{\textstyle\sum_{\y}}\E[\mathbf{1}(\s\leq\t_{\mathrm{st}})\N\cdot\mathrm{Cent}\Upsilon^{\pm,\mathrm{j},\mathrm{i},\t}(\s,\y(\s))]\d\s \label{eq:bg213101Ia}\\
\eqref{eq:bg2138Id} \ &\lesssim \ {\textstyle\sup_{\t}\sup_{\mathrm{i}}}\N^{\beta_{\mathrm{BG}}}{\textstyle\int_{\tau(\mathrm{j},\mathrm{i})}^{1}}|\mathbb{T}(\N)|^{-1}{\textstyle\sum_{\y}}\E[\mathbf{1}(\s\leq\t_{\mathrm{st}})\N\cdot\mathrm{Cent}\Lambda^{\pm,\mathrm{j},\mathrm{i},\t}(\s,\y(\s))]\d\s. \label{eq:bg213101Ib}
\end{align}
See {Definition \ref{definition:bg2131}} for $\tau(\mathrm{j},\mathrm{i})$. The suprema in \eqref{eq:bg213101Ia} and \eqref{eq:bg213101Ib} are over $0\leq\t\leq1$ and $1\leq\mathrm{i}<\mathrm{i}(\mathrm{j})$. We now define:
\begin{enumerate}
\item Fix $\s\geq\tau(\mathrm{j},\mathrm{i})$ and $\t\in[0,1]$. Take $\mathds{R}^{\chi,\mathfrak{q},\pm,\mathrm{j}}$ in {Definition \ref{definition:bg26}} and $\mathbf{G}^{\t}$ in {Definition \ref{definition:bg211}}. With explanation after, we set
{\small
\begin{align*}
&\mathscr{A}^{\pm,\t} \ := \ \\
&\tau(\mathrm{j},\mathrm{i})^{-1}{\textstyle\int_{0}^{\tau(\mathrm{j},\mathrm{i})}}\{\mathfrak{m}(\mathrm{j},\mathrm{i})^{-1}{\textstyle\sum_{\mathrm{k}=0}^{\mathfrak{m}(\mathrm{j},\mathrm{i})-1}}\mathds{R}^{\chi,\mathfrak{q},\pm,\mathrm{j}}(\s-\r,\y(\s-\r)\pm2\mathrm{k}\mathfrak{l}(\mathrm{j}))\mathbf{G}^{\t}(\s-\r,\y(\s-\r)\pm2\mathrm{k}\mathfrak{l}(\mathrm{j}))\}\d\r.
\end{align*}
}In words, $\mathscr{A}^{\pm,\t}$ is almost $\mathds{A}^{\mathfrak{m}(\mathrm{j},\mathrm{i}),\tau(\mathrm{j},\mathrm{i}),\pm}(\mathds{R}^{\chi,\mathfrak{q},\pm,\mathrm{j}}\mathbf{Z};\s,\y(\s))$ (using notation of {Definition \ref{definition:bg211}} with $\mathfrak{m}(\mathrm{j},\mathrm{i})$ in {Definition \ref{definition:bg2131}}). However, in $\mathscr{A}^{\pm,\t}$ above, the coupling constant in $\mathscr{A}^{\pm,\t}$ is fixed to be $\lambda(\t)$, not $\lambda(\s)$ as was the case in \eqref{eq:bg211I}. 
\item We now define the following ``centered" or ``unshifted" terms, for which we recall exponents $\beta(\mathrm{j},\cdot)$ from {Definition \ref{definition:bg2131}}:
\begin{align}
\mathrm{Cent}\Upsilon^{\pm,\mathrm{j},\mathrm{i},\t}(\s,\y(\s)) \ &:= \ \N^{\beta(\mathrm{j},\mathrm{i}-1)}|\mathscr{A}^{\pm,\t}|^{2}\mathbf{1}\{|\mathscr{A}^{\pm,\t}|\lesssim\N^{-\beta(\mathrm{j},\mathrm{i}-1)}\} \label{eq:bg213101IIa}\\
\mathrm{Cent}\Lambda^{\pm,\mathrm{j},\mathrm{i},\t}(\s,\y(\s)) \ &:= \ \N^{-\beta(\mathrm{j},\mathrm{i}-1)}\mathbf{1}\{|\mathscr{A}^{\pm,\t}|>\N^{-\beta(\mathrm{j},\mathrm{i}-1)}\}. \label{eq:bg213101IIb}
\end{align}
\end{enumerate}
\end{lemma}
\begin{proof}
For convenience, in this (short) proof let $\wt{\sum}_{\y}:=|\mathbb{T}(\N)|^{-1}\sum_{\y}$ be a normalized sum over $\mathbb{T}(\N)$. We first claim, with explanations given afterwards, that \eqref{eq:bg2138Ic} equals the supremum over indices $\mathrm{i}$ of
{\small
\begin{align}
&\N^{\beta_{\mathrm{BG}}}{\textstyle\int_{\tau(\mathrm{j})}^{1}}{\textstyle\wt{\sum}_{\y}}\E[\mathbf{1}(\s\leq\t_{\mathrm{st}}){\textstyle\sum^{\mathrm{j},\mathrm{i}+1}\sum^{\mathrm{i},+}}\N^{\beta(\mathrm{j},\mathrm{i}-1)+1}|\mathds{A}^{\mathds{Q},\pm}[\s,\y(\s)]|^{2}\mathbf{1}\{|\mathds{A}^{\mathds{Q},\pm}[\s,\y(\s)]|\leq\N^{-\beta(\mathrm{j},\mathrm{i}-1)}\}]\d\s \nonumber\\
&= {\textstyle\sum^{\mathrm{j},\mathrm{i}+1}\sum^{\mathrm{i},+}}\N^{\beta_{\mathrm{BG}}}{\textstyle\int_{\tau(\mathrm{j})}^{1}}{\textstyle\wt{\sum}_{\y}}\E[\mathbf{1}(\s\leq\t_{\mathrm{st}})\N^{\beta(\mathrm{j},\mathrm{i}-1)+1}|\mathds{A}^{\mathds{Q},\pm}[\s,\y(\s)]|^{2}\mathbf{1}\{|\mathds{A}^{\mathds{Q},\pm}[\s,\y(\s)]|\leq\N^{-\beta(\mathrm{j},\mathrm{i}-1)}\}]\d\s \nonumber\\
&\leq {\textstyle\sum^{\mathrm{j},\mathrm{i}+1}\sum^{\mathrm{i},+}\sup_{\t}}\N^{\beta_{\mathrm{BG}}}{\textstyle\int_{\tau(\mathrm{j},\mathrm{i})}^{1}}{\textstyle\wt{\sum}_{\y}}\E[\mathbf{1}(\s\leq\t_{\mathrm{st}})\N^{\beta(\mathrm{j},\mathrm{i}-1)+1}|\mathscr{A}^{\pm,\t}|^{2}\mathbf{1}\{|\mathscr{A}^{\pm,\t}|\leq\N^{-\beta(\mathrm{j},\mathrm{i}-1)}\}]\d\s \nonumber\\
&= {\textstyle\sup_{\t}}\N^{\beta_{\mathrm{BG}}}{\textstyle\int_{\tau(\mathrm{j},\mathrm{i})}^{1}}{\textstyle\wt{\sum}_{\y}}\E[\mathbf{1}(\s\leq\t_{\mathrm{st}})\N^{\beta(\mathrm{j},\mathrm{i}-1)+1}|\mathscr{A}^{\pm,\t}|^{2}\mathbf{1}\{|\mathscr{A}^{\pm,\t}|\leq\N^{-\beta(\mathrm{j},\mathrm{i}-1)}\}]\d\s. \label{eq:bg21310I1}
\end{align}
}{The first line follows by definition. (Also, note {that} the absolute value bars around $\Upsilon$ in $\mathrm{RHS}\eqref{eq:bg2138Ic}$ are redundant, since this $\Upsilon$-term is already non-negative.) The second line follows from pulling the last two sums in the first line above outside the time-integral. The third line requires some justification. For now, we recall that $\mathds{A}^{\mathds{Q},\pm}[\s,\y(\s)]$ is $\mathrm{RHS}\eqref{eq:bg211I}$ for $\mathsf{F}=\mathds{R}^{\chi,\mathfrak{q},\pm,\mathrm{j}}$ and $\mathfrak{l}(\mathsf{F})=\mathfrak{l}(\mathrm{j})$, and upon replacing the averaging set therein with $\mathds{Q}$; see Definition \ref{definition:bg2133}. In particular, $\mathds{A}^{\mathds{Q},\pm}[\s,\y(\s)]$ is $\mathscr{A}^{\pm,\t}$ from the statement of the lemma but replacing the averaging domain $[0,\tau(\mathrm{j},\mathrm{i})]\times\llbracket0,\mathfrak{m}(\mathrm{j},\mathrm{i})-1\rrbracket$ by $\mathds{Q}$ and setting $\t=\s$. The first step that we take to justify the third line is the change-of-variables for the $[\tau(\mathrm{j}),1]\times\mathbb{T}(\N)$-integration that shifts $(\s,\y)\mapsto(\wt{\s},\wt{\y})$ where, under the $(\wt{\s},\wt{\y})$-variables, $\mathds{A}^{\mathds{Q},\pm}[\s,\y(\s)]$ becomes $\mathscr{A}^{\pm,\t}$ as written, for $\t$ equal to a $\mathds{Q}$-dependent shift of $\s$. (This change-of-variables is just the shift that maps $\mathds{Q}\to[0,\tau(\mathrm{j},\mathrm{i})]\times\llbracket0,\mathfrak{m}(\mathrm{j},\mathrm{i})-1\rrbracket$.) Because shifts on $\mathbb{T}(\N)$ are bijections, the spatial summation over $\y$ does not change after this change-of-variables. The domain of time-integration after the change-of-variables equals a shift of $[\tau(\mathrm{j}),1]$. To compute it, we first make two observations. It must have infimum at least $\tau(\mathrm{j},\mathrm{i})$. Otherwise, $\mathscr{A}^{\pm,\t}$ evaluates \eqref{eq:glsde} at negative time; this does not make sense.  Moreover, the second line above restricts to $\s\leq\t_{\mathrm{st}}$. After this change-of-variables, we can still restrict to $\s\leq\t_{\mathrm{st}}$. Indeed, changing variables in the second line does not change the fact that it depends on \eqref{eq:hf}-\eqref{eq:glsde} for times $\s\leq\t_{\mathrm{st}}$. This explains $\mathbf{1}(\s\leq\t_{\mathrm{st}})$ in the third line. (Since $\t_{\mathrm{st}}\leq1$, we can push the upper limit of integration to 1 as well). So, the new integration-domain is $\subseteq[\tau(\mathrm{j},\mathrm{i}),1]\times\mathbb{T}(\N)$. As the integrand in the second line above is non-negative, the third line follows once we take a supremum over $\mathds{Q}$-dependent shifts for the time-parameter $\t$ in the coupling constant in $\mathscr{A}^{\pm,\t}$. (This explains the supremum over $0\leq\t\leq1$. We can restrict the supremum to $\t\in[0,1]$ because the coupling constants appearing in $\mathds{A}^{\mathds{Q},\pm}[\s,\y(\s)]$ for $0\leq\s\leq1$ are all evaluated at a time in $[0,1]$. Also, the supremum over $\t$ sits outside the $[\tau(\mathrm{j},\mathrm{i}),1]\times\mathbb{T}(\N)$-integral because it depends only on the $\mathds{Q}$-shift.) The first identity in \eqref{eq:bg21310I1} holds because the integral in the third line is independent of the double-average variables $\mathds{Q}[+],\mathds{Q}$, so this double-average does nothing. (We just removed dependence on $\mathds{Q}[+],\mathds{Q}$ to get the third line.) 

Now, note that the supremum over $\mathrm{i}$ of \eqref{eq:bg21310I1} is just $\mathrm{RHS}\eqref{eq:bg213101Ia}$. So, \eqref{eq:bg213101Ia} holds. To get \eqref{eq:bg213101Ib}, identical reasoning and calculation suffices.}
\end{proof}
\subsubsection{Localization of $\mathrm{Cent}\Upsilon^{\pm,\mathrm{j},\mathrm{i},\t}$ and $\mathrm{Cent}\Lambda^{\pm,\mathrm{j},\mathrm{i},\t}$}
Note $\mathrm{Cent}\Upsilon^{\pm,\mathrm{j},\mathrm{i},\t}$ and $\mathrm{Cent}\Lambda^{\pm,\mathrm{j},\mathrm{i},\t}$ from Lemma \ref{lemma:bg213101} are determined by mesoscopic space-time averages. Thus, in the sense of Lemmas \ref{lemma:le12}, \ref{lemma:le15}, both of these objects are morally local functionals. Our goal is to now make this precise, ultimately by comparing \eqref{eq:hf}-\eqref{eq:glsde} to localized versions from Definitions \ref{definition:le10} and \ref{definition:le14}. In the following result, we will introduce (necessarily) subtle and complicated constructions. We clarify what these constructions are saying in an intuitive manner in Remark \ref{remark:finalprop3}. (In doing so, we also give an intuitive explanation of how the proof of Lemma \ref{lemma:finalprop2} ultimately goes. This is meant to clarify said proof, since it has several steps.)
\begin{lemma}\label{lemma:finalprop2}
 We have the following estimates {for any large but fixed $\mathrm{D}>0$}, which use notation to be explained afterwards:
{\small
\begin{align}
&\E[\mathbf{1}(\s\leq\t_{\mathrm{st}})\mathrm{Cent}\Upsilon^{\pm,\mathrm{j},\mathrm{i},\t}(\s,\y(\s))] \nonumber\\
&\lesssim \ \N^{2\gamma_{\mathrm{ap}}}\E[\{\E^{\mathrm{loc},\s}[\mathrm{Loc}^{(1)}\Upsilon^{\pm,\mathrm{j},\mathrm{i},\t}(\s)]\}(\Pi^{\mathrm{j},\mathrm{i}}\mathbf{U}^{\s-\tau(\mathrm{j},\mathrm{i}),\y(\s-\tau(\mathrm{j},\mathrm{i}))+\cdot})]+\N^{-{\mathrm{D}}} \label{eq:finalprop2Ia}\\
&\E[\mathbf{1}(\s\leq\t_{\mathrm{st}})\mathrm{Cent}\Lambda^{\pm,\mathrm{j},\mathrm{i},\t}(\s,\y(\s))] \nonumber\\
&\lesssim \ \N^{2\gamma_{\mathrm{ap}}}\E[\{\E^{\mathrm{loc},\s}[\mathrm{Loc}^{(1)}\Lambda^{\pm,\mathrm{j},\mathrm{i},\t}(\s)]\}(\Pi^{\mathrm{j},\mathrm{i}}\mathbf{U}^{\s-\tau(\mathrm{j},\mathrm{i}),\y(\s-\tau(\mathrm{j},\mathrm{i}))+\cdot})]+\N^{-{\mathrm{D}}}. \label{eq:finalprop2Ib}
\end{align}
}We take $\s\geq\tau(\mathrm{j},\mathrm{i})$ and $\y\in\mathbb{T}(\N)$ and $0\leq\t\leq1$. We also used notation from {Lemma \ref{lemma:bg213101}} and the following:
\begin{enumerate}
\item We define $\mathbb{I}:=\llbracket-10\mathfrak{m}(\mathrm{j},\mathrm{i})\mathfrak{l}(\mathrm{j})-\N^{3/2+\gamma_{\mathrm{ap}}}\tau(\mathrm{j},\mathrm{i}),10\mathfrak{m}(\mathrm{j},\mathrm{i})\mathfrak{l}(\mathrm{j})+\N^{3/2+\gamma_{\mathrm{ap}}}\tau(\mathrm{j},\mathrm{i})\rrbracket$; see {Definition \ref{definition:bg2131}} for relevant notation. Consider $\mathbf{U}^{\s-\r,\cdot}[\mathbb{I}(\mathfrak{t})]$ and $\mathtt{J}(\s-\r,\cdot;\mathbb{I}(\mathfrak{t}))$ for $0\leq\r\leq\tau(\mathrm{j},\mathrm{i})$ with $\mathbb{I}$ from the first sentence of this bullet point and $\mathfrak{t}=\tau(\mathrm{j},\mathrm{i})$. See {Definitions \ref{definition:le10}, \ref{definition:le14}, \ref{definition:bg2131}} for relevant notation. We give initial data $\mathbf{U}^{\s-\tau(\mathrm{j},\mathrm{i}),\x}[\mathbb{I}(\mathfrak{t})]$ in a later bullet point. For now, we only specify the initial data $\mathtt{J}(\tau,\inf\mathbb{I}(\mathfrak{t});\mathbb{I}(\mathfrak{t}))=0$ at $\tau=\s-\tau(\mathrm{j},\mathrm{i})$ (which is the initial time of $\s-\r$ for $0\leq\r\leq\tau(\mathrm{j},\mathrm{i})$.)
\item For any $\mathbf{U}\in\R^{\mathbb{T}(\N)}$, we let $\Pi^{\mathrm{j},\mathrm{i}}\mathbf{U}\in\R^{\mathbb{I}(\mathfrak{t})}$ be the canonical projection. (So, for any $\x\in\mathbb{I}(\mathfrak{t})$, we have $\Pi^{\mathrm{j},\mathrm{i}}\mathbf{U}(\x)=\mathbf{U}(\x)$. We note that the discrete interval $\mathbb{I}(\mathfrak{t})$ from the previous bullet point depends only on $\mathrm{j},\mathrm{i}$. This explains the superscripts in $\Pi^{\mathrm{j},\mathrm{i}}$.)
\item We have set $\E^{\mathrm{loc},\s}\{\cdot\}(\mathbf{U})$ to be {the} expectation of $\cdot$ with respect to the law of $\mathbf{U}^{\s-\r,\cdot}[\mathbb{I}(\mathfrak{t})]$ and $\mathtt{J}(\s-\r,\cdot;\mathbb{I}(\mathfrak{t}))$ for $0\leq\r\leq\tau(\mathrm{j},\mathrm{i})$, assuming the initial data is $\mathbf{U}^{\s-\tau(\mathrm{j},\mathrm{i}),\x}[\mathbb{I}(\mathfrak{t})]=\mathbf{U}(\x)$ and $\mathtt{J}(\s-\tau(\mathrm{j},\mathrm{i}),\inf\mathbb{I}(\mathfrak{t});\mathbb{I}(\mathfrak{t}))=0$. (Here, $\mathbf{U}\in\R^{\mathbb{I}(\mathfrak{t})}$ and $\x\in\mathbb{I}(\mathfrak{t})$.) In particular, the first term in $\mathrm{RHS}\eqref{eq:finalprop2Ia}$ is {the} expectation of a path-space functional \eqref{eq:finalprop2IIa} with initial data for $\mathbf{U}^{\s-\tau(\mathrm{j},\mathrm{i}),\cdot}[\mathbb{I}(\mathfrak{t})]$ that is sampled according to the law of \eqref{eq:glsde} at time $\s-\tau(\mathrm{j},\mathrm{i})$ and shifted in space by $\y(\s-\tau(\mathrm{j},\mathrm{i}))$. (To be completely clear, we emphasize that $\E$ in $\mathrm{RHS}\eqref{eq:finalprop2Ia}$ is with respect to the $\mathbb{I}(\mathfrak{t})$-projection $\Pi^{\mathrm{j},\mathrm{i}}\mathbf{U}^{\s-\tau(\mathrm{j},\mathrm{i}),\y(\s-\tau(\mathrm{j},\mathrm{i}))+\cdot}$.)
\item We now define two auxiliary objects, which uses notation to be introduced after defining said objects. For convenience, set
\begin{align}
\mathsf{G}^{\t}(\tau,\x) \ := \ \exp[\lambda(\t)\mathtt{J}(\tau,\x;\mathbb{I}(\mathfrak{t}))+\lambda(\t)\mathrm{U}].
\end{align}
(This is a ``localized" version of what we built in {Definition \ref{definition:bg211}}.) Here, we take $\tau=\s-\r$ for $0\leq\r\leq\tau(\mathrm{j},\mathrm{i})$ as before, and we take $\x\in\mathbb{I}(\mathfrak{t})$. The random variable $\mathrm{U}$ is uniform on $[-1,1]$, and it is independent of everything else. Moreover, recall the coupling constants $\lambda(\t)$ in {Definitions \ref{definition:intro5}, \ref{definition:intro6}}. Now set {\small$\mathscr{A}^{\mathrm{loc},(1),\pm,\t}$} as the following ``localization" of $\mathscr{A}^{\pm,\t}$ in {Lemma \ref{lemma:bg213101}}:
{\small
\begin{align}
&\mathscr{A}^{\mathrm{loc},(1),\pm,\t}:=\nonumber\\
&\tau(\mathrm{j},\mathrm{i})^{-1}{\textstyle\int_{0}^{\tau(\mathrm{j},\mathrm{i})}}\{\mathfrak{m}(\mathrm{j},\mathrm{i})^{-1}{\textstyle\sum_{\mathrm{k}=0}^{\mathfrak{m}(\mathrm{j},\mathrm{i})-1}}\mathfrak{A}^{(\mathrm{k}),(1),\pm}_{\s-\r}(\mathbf{U}^{\s-\r,\cdot}[\mathbb{I}(\mathfrak{t})])\times\mathsf{G}^{\t}(\s-\r,0(\s-\r)\pm2\mathrm{k}\mathfrak{l}(\mathrm{j}))\}\d\r. \label{eq:finalprop2Ic}
\end{align}
}{Above, {\small$\mathfrak{A}^{(\mathrm{k}),(1),\pm}_{\tau}(\mathbf{U})$} is a functional whose support $\mathbb{J}(\tau,\mathrm{k})\subseteq\mathbb{I}(\mathfrak{t})$ is a $\tau$-dependent shift of $\mathbb{I}(\mathrm{k})$. It satisfies the deterministic estimate {\small$|\mathfrak{A}^{(\mathrm{k}),(1),\pm}|\lesssim\N^{25\gamma_{\mathrm{reg}}}\mathfrak{l}(\mathrm{j})^{-3/2}$}. Also, for any $\rho,\tau$, we have {\small$\E^{\rho,\tau,\mathbb{J}(\tau,\mathrm{k})}\mathfrak{A}^{(\mathrm{k}),(1),\pm}_{\tau}=0$}. As for the discrete intervals $\mathbb{J}(\tau,\mathrm{k})$, we first have the bound $|\mathbb{J}(\tau,\mathrm{k})|\lesssim\mathfrak{l}(\mathrm{j})$. Moreover, the sets $\mathbb{J}(\tau,\mathrm{k})$ are mutually disjoint over $\mathrm{k}$. We also know that, depending on the choice of sign $\pm$, given $\tau,\mathrm{k}$, the discrete interval $\mathbb{I}(\t,\mathrm{k})$ is either disjoint from the interval $\llbracket\inf\mathbb{I}(\mathfrak{t}),0(\tau)\pm2\mathrm{k}\mathfrak{l}(\mathrm{j})\rrbracket$ or contained inside it. Lastly, $0(\s-\r)$ is the image of the characteristic $\y\mapsto\y(\tau)$ at $\y=0$ and $\tau=\s-\r$; see {Definition \ref{definition:intro6}}.}
\item Let us now define the main objects on the RHS of \eqref{eq:finalprop2Ia} and \eqref{eq:finalprop2Ib}. With notation explained after, we set
\begin{align}
&\mathrm{Loc}^{(1)}\Upsilon^{\pm,\mathrm{j},\mathrm{i},\t}(\s) \nonumber\\
&:= \ \mathbf{1}[\mathcal{E}^{\mathrm{ap},\mathrm{j},\mathrm{i}}(\s)]\mathbf{1}[\mathcal{E}^{\mathrm{dens},\mathrm{j},\mathrm{i}}(\s)]\N^{\beta(\mathrm{j},\mathrm{i}-1)}|\mathscr{A}^{\mathrm{loc},(1),\pm,\t}|^{2}\mathbf{1}[|\mathscr{A}^{\mathrm{loc},(1),\pm,\t}|\lesssim\N^{-\beta(\mathrm{j},\mathrm{i}-1)+\gamma_{\mathrm{ap}}}] \label{eq:finalprop2IIa}\\
&\mathrm{Loc}^{(1)}\Lambda^{\pm,\mathrm{j},\mathrm{i},\t}(\s) \nonumber\\
&:= \ \mathbf{1}[\mathcal{E}^{\mathrm{ap},\mathrm{j},\mathrm{i}}(\s)]\mathbf{1}[\mathcal{E}^{\mathrm{dens},\mathrm{j},\mathrm{i}}(\s)]\N^{-\beta(\mathrm{j},\mathrm{i}-1)}\mathbf{1}[|\mathscr{A}^{\mathrm{loc},(1),\pm,\t}|\gtrsim\N^{-\beta(\mathrm{j},\mathrm{i}-1)-\gamma_{\mathrm{ap}}}]. \label{eq:finalprop2IIb}
\end{align}
Above, $\mathcal{E}^{\mathrm{ap},\mathrm{j},\mathrm{i}}(\s)$ denotes the event where $\mathsf{G}^{\t}$-factors in the definition of {\small$\mathscr{A}^{\mathrm{loc},(1),\pm,\t}$} are each $\gtrsim\N^{-\gamma_{\mathrm{ap}}}$ and $\lesssim\N^{\gamma_{\mathrm{ap}}}$. The event $\mathcal{E}^{\mathrm{dens},\mathrm{j},\mathrm{i}}(\s)$ is exactly when the charge density $\sigma$ of the local process $\mathbf{U}^{\t,\x}[\mathbb{I}(\mathfrak{t})]$ over $\x\in\mathbb{I}(\mathfrak{t})$ satisfies {\small$|\sigma|\lesssim\N^{\gamma_{\mathrm{reg}}}|\mathbb{I}(\mathfrak{t})|^{-1/2}$}.
\end{enumerate}
\end{lemma}
\begin{rem}\label{remark:finalprop3}
{{\small$\mathrm{Loc}^{(1)}\Upsilon^{\pm,\mathrm{j},\mathrm{i},\t}$} and $\mathrm{Cent}\Upsilon^{\pm,\mathrm{j},\mathrm{i},\t}$ are (almost) the same form. (Indeed, {\small$\mathfrak{A}^{(\mathrm{k}),(1),\pm}$} is ultimately the underlying function defining $\mathds{R}^{\chi,\mathfrak{q},\pm,\mathrm{j}}$ from Definition \ref{definition:bg26}. We have just abstracted away all information besides what we need to estimate the RHS of \eqref{eq:finalprop2Ia} and \eqref{eq:finalprop2Ib}, respectively. Of course, showing that shifts of $\mathds{R}^{\chi,\mathfrak{q},\pm,\mathrm{j}}$ have the properties claimed in Lemma \ref{lemma:finalprop2} requires an argument. But this is not so difficult; they are actually even engineered to be true in Definition \ref{definition:bg26}.) Now, we note $\mathrm{Cent}\Upsilon^{\pm,\mathrm{j},\mathrm{i},\t}$ is evaluated at \eqref{eq:hf}-\eqref{eq:glsde} but restricted to a local space-time set; see Lemma \ref{lemma:bg213101} and Definition \ref{definition:bg211}. Also, we note {\small$\mathrm{Loc}^{(1)}\Upsilon^{\pm,\mathrm{j},\mathrm{i},\t}$} is evaluated at \eqref{eq:glsdeloc} and $\mathtt{J}$ from Definition \ref{definition:le14}, restricted to the same space-time set. Therefore, Lemmas \ref{lemma:le12}, \ref{lemma:le15} let us ultimately compare {\small$\mathrm{Loc}^{(1)}\Upsilon^{\pm,\mathrm{j},\mathrm{i},\t}$} and $\mathrm{Cent}\Upsilon^{\pm,\mathrm{j},\mathrm{i},\t}$. (Indeed, the subset $\mathbb{I}(\mathfrak{t})$ is built in Lemma \ref{lemma:finalprop2} to make Lemmas \ref{lemma:le12}, \ref{lemma:le15} applicable. There is, however, the technical issue that the initial data for $\mathtt{J}$ in the setting of Lemma \ref{lemma:le15} should be $\mathbf{J}$ evaluated at some point, not zero as we assumed in Lemma \ref{lemma:finalprop2} above. But, the $\mathtt{J}$ SDE \eqref{eq:hfloc} commutes with addition; in words, KPZ-type growth is not height dependent, only slope-dependent. Therefore, the discrepancy in initial data results in adding to $\mathtt{J}$ some value of $\mathbf{J}$. At the level of their exponentials, we can then factor out the $\mathbf{J}$ piece. A similar discussion holds for the $\mathrm{U}$-term we add to $\mathtt{J}$, too.) Lastly, we will briefly explain the $(1)$-superscripts in Lemma \ref{lemma:finalprop2}. Later in this section, we technically adjust {\small$\mathrm{Loc}^{(1)}\Upsilon^{\pm,\mathrm{j},\mathrm{i},\t}$} and {\small$\mathfrak{A}^{(\mathrm{k}),(1),\pm}$} into terms that we can apply Proposition \ref{prop:kv1} to, at which point we will drop all of the $(1)$ superscripts. (Everything that we have discussed in this remark also applies, at least in an intuitive sense, to {\small$\mathrm{Loc}^{(1)}\Lambda^{\pm,\mathrm{j},\mathrm{i},\t}$} and $\mathrm{Cent}\Lambda^{\pm,\mathrm{j},\mathrm{i},\t}$ as well.)}
\end{rem}
\begin{proof}
For convenience, assume $\y(\s-\tau(\mathrm{j},\mathrm{i}))=0$ in \eqref{eq:finalprop2Ia}-\eqref{eq:finalprop2Ib}, and set the accompanying notation $\z(\s-\r):=\y(\s-\r)$ for any $\s,\r$ and this choice of $\y$. (This assumption on $\y$ is to avoid having to re-center $\mathbb{I}$ and $\mathbb{I}(\mathfrak{t})$ sets. The notation $\z(\s-\r)$ will be useful as all functionals in the definition of $\mathscr{A}^{\pm,\t}$ in Lemma \ref{lemma:bg213101} are evaluated at or near these points. Indeed, $\z(\s-\r)$ can be reconstructed in terms of the origin $0$, but this would be inconvenient to repeatedly write.) We remove this assumption at the end.
\subsubsection*{Step 1: Construction of local SDEs}
Consider $\mathbf{U}^{\s-\r,\cdot}[\mathbb{I}(\mathfrak{t})]$ in the statement of this lemma for $\r\in[0,\tau(\mathrm{j},\mathrm{i})]$, and give it the initial data $\mathbf{U}^{\s-\tau(\mathrm{j},\mathrm{i}),\cdot}[\mathbb{I}(\mathfrak{t})]=\mathbf{U}^{\s-\tau(\mathrm{j},\mathrm{i}),\cdot}$. We also take $\mathtt{J}$ in the statement of this lemma. Its initial data for all $\x\in\mathbb{I}(\mathfrak{t})$ is determined by initial data for $\mathbf{U}^{\s-\r,\cdot}[\mathbb{I}(\mathfrak{t})]$ that we just gave and the gradient relation in Definition \ref{definition:le10}. Now, set $\mathbf{J}^{!}:=\mathbf{J}^{!}(\s-\tau(\mathrm{j},\mathrm{i}),\inf\mathbb{I}(\mathfrak{t}))$, and write $\mathtt{J}(\tau,\x;\mathbb{I}(\mathfrak{t}))=\mathtt{J}(\tau,\x;\mathbb{I}(\mathfrak{t}))+\mathbf{J}^{!}-\mathbf{J}^{!}$, where $\mathbf{J}^{!}(\tau,\x)$ is the renormalized current that we defined in the proof of Lemma \ref{lemma:le15}. The flow along the SDE \eqref{eq:hfloc} which $\mathtt{J}$ solves commutes with addition. Indeed, the RHS of said SDE does not depend on the solution, just \eqref{eq:glsdeloc}. So, $\mathtt{J}(\tau,\inf\mathbb{I}(\mathfrak{t});\mathbb{I}(\mathfrak{t}))+\mathbf{J}^{!}$ is the solution to \eqref{eq:hfloc} with initial data $\mathtt{J}(\s-\r,\inf\mathbb{I}(\mathfrak{t});\mathbb{I}(\mathfrak{t}))+\mathbf{J}^{!}=\mathbf{J}^{!}$, and $\mathsf{J}(\tau,\x;\mathbb{I}(\mathfrak{t}))$ using this and the gradient relation in Definition \ref{definition:le14} (since $\mathbf{J}^{!}$ is constant in $\x$). Now, by Lemmas \ref{lemma:le12}, \ref{lemma:le15}, there is a very high probability event $\mathcal{E}^{\mathrm{vhp}}$ such that on this event, for all $\r\in[0,\tau(\mathrm{j},\mathrm{i})]$ and $\x\in\mathbb{I}$, {we have the following for any large but fixed $\mathrm{D}>0$}:
\begin{align}
|\mathbf{U}^{\s-\r,\x}[\mathbb{I}(\mathfrak{t})]-\mathbf{U}^{\s-\r,\x}|+|\mathtt{J}(\s-\r,\x;\mathbb{I}(\mathfrak{t}))+\mathbf{J}^{!}-\mathbf{J}^{!}(\s-\r,\x)| \ \lesssim \ \N^{-{\mathrm{D}}}. \label{eq:finalprop2I1}
\end{align}
%
\subsubsection*{Step 2: Modifying $\mathscr{A}^{\pm,\mathrm{t}}$}
Recall $\mathscr{A}^{\pm,\t}$  from Lemma \ref{lemma:bg213101}. In this step, we will modify the $\mathbf{G}^{\t}$ factor therein. Set
\begin{align}
\mathbf{G}^{\t,\mathrm{mod}}(\tau,\x) \ := \ \exp[\lambda(\t)\{\mathbf{J}^{!}(\tau,\x)-\mathbf{J}^{!}\}+\lambda(\t)\mathrm{U}]. \label{eq:finalprop2I2}
\end{align}
In words, $\mathbf{G}^{\t,\mathrm{mod}}$ is just $\mathbf{G}^{\t}$ from Definition \ref{definition:bg211} with the additional shift of $-\mathbf{J}^{!}+\mathrm{U}$. Observe that $\mathbf{J}^{!}$ and $\mathrm{U}$ are independent of the $(\tau,\x)$ variables in $\mathrm{LHS}\eqref{eq:finalprop2I2}$. Using this and linearity of integration/summation, with notation explained after, we have 
\begin{align}
\mathscr{A}^{\pm,\t} \ = \ \exp[\lambda(\t)\mathbf{J}^{!}-\lambda(\t)\mathrm{U}]\mathscr{A}^{\pm,\mathrm{t},\mathrm{mod}}, \label{eq:finalprop2I3a}
\end{align}
where $\mathscr{A}^{\pm,\mathrm{t},\mathrm{mod}}$ is basically just $\mathscr{A}^{\pm,\t}$ but replacing $\mathbf{G}^{\t}\mapsto\mathbf{G}^{\t,\mathrm{mod}}$. Precisely, it is equal to the following (since $\y=0$):
{\small
\begin{align}
\tau(\mathrm{j},\mathrm{i})^{-1}{\textstyle\int_{0}^{\tau(\mathrm{j},\mathrm{i})}}\{\mathfrak{m}(\mathrm{j},\mathrm{i})^{-1}{\textstyle\sum_{\mathrm{k}=0}^{\mathfrak{m}(\mathrm{j},\mathrm{i})-1}}\mathds{R}^{\chi,\mathfrak{q},\pm,\mathrm{j}}(\s-\r,\z(\s-\r)\pm2\mathrm{k}\mathfrak{l}(\mathrm{j}))\mathbf{G}^{\t,\mathrm{mod}}(\s-\r,\z(\s-\r)\pm2\mathrm{k}\mathfrak{l}(\mathrm{j}))\}\d\r. \label{eq:finalprop2I3b}
\end{align}
}
\subsubsection{Step 3: ``Localization"}
In this step, we now replace $\mathds{R}$ and $\mathbf{G}$ terms in \eqref{eq:finalprop2I3b} by ``localizations" defined with respect to $\mathbf{U}^{\tau,\cdot}[\mathbb{I}(\mathfrak{t})]$ and $\mathtt{J}(\tau,\cdot;\mathbb{I}(\mathfrak{t}))$. Let us make this precise. Observe that $\mathbf{G}^{\t,\mathrm{mod}}$ from \eqref{eq:finalprop2I2} is the same exponential as $\mathsf{G}^{\t}$ from the statement of the lemma. Except, the former is evaluated at $\mathbf{J}^{!}(\cdot,\cdot)-\mathbf{J}^{!}$ and the latter is evaluated at $\mathtt{J}(\cdot,\cdot;\mathbb{I}(\mathfrak{t}))$. For convenience, fix $\x(\s,\mathrm{r},\mathrm{k}):=\z(\s-\r)\pm2\mathrm{k}\mathfrak{l}(\mathrm{j})$. We claim the following calculation holds (with explanation given afterwards):
\begin{align}
&|\mathbf{G}^{\t,\mathrm{mod}}(\s-\r,\x(\s,\r,\mathrm{k}))-\mathsf{G}^{\t}(\s-\r,\x(\s,\r,\mathrm{k}))| \nonumber\\
&\lesssim \ |{\textstyle\int_{\mathtt{J}(\s-\r,\x(\s,\r,\mathrm{k});\mathbb{I}(\mathfrak{t}))}^{\mathbf{J}^{!}(\s-\r,\x(\s,\r,\mathrm{k}))-\mathbf{J}^{!}}}\exp[\lambda(\t)\eta+\lambda(\t)\mathrm{U}]\d\eta| \label{eq:finalprop2I4a}\\
&\lesssim \ |\mathbf{G}^{\t,\mathrm{mod}}(\s-\r,\x(\s,\r,\mathrm{k}))+\mathsf{G}^{\t}(\s-\r,\x(\s,\r,\mathrm{k}))|\label{eq:finalprop2I4b}\\
&\cdot|\mathbf{J}^{!}(\s-\r,\x(\s,\r,\mathrm{k}))-\mathbf{J}^{!}-\mathtt{J}(\s-\r,\x(\s,\r,\mathrm{k});\mathbb{I}(\mathfrak{t}))| \nonumber\\
&\lesssim \ |\mathbf{G}^{\t,\mathrm{mod}}(\s-\r,\x(\s,\r,\mathrm{k}))-\mathsf{G}^{\t}(\s-\r,\x(\s,\r,\mathrm{k}))|\label{eq:finalprop2I4c}\\
&\cdot|\mathbf{J}^{!}(\s-\r,\x(\s,\r,\mathrm{k}))-\mathbf{J}^{!}-\mathtt{J}(\s-\r,\x(\s,\r,\mathrm{k});\mathbb{I}(\mathfrak{t}))| \nonumber\\
&+ \ \mathbf{G}^{\t,\mathrm{mod}}(\s-\r,\x(\s,\r,\mathrm{k}))\cdot|\mathbf{J}^{!}(\s-\r,\x(\s,\r,\mathrm{k}))-\mathbf{J}^{!}-\mathtt{J}(\s-\r,\x(\s,\r,\mathrm{k});\mathbb{I}(\mathfrak{t}))|. \label{eq:finalprop2I4d}
\end{align}
\eqref{eq:finalprop2I4a} follows by the fundamental theorem calculus and $|\lambda(\t)|\lesssim1$. To get \eqref{eq:finalprop2I4b}, we maximize the integrand in $\mathrm{RHS}\eqref{eq:finalprop2I4a}$ over $\eta$ in the integration-domain. This means at $\eta$ equal to the two boundary points of said integration-domain, which controls the integrand in $\mathrm{RHS}\eqref{eq:finalprop2I4a}$ by the first factor in \eqref{eq:finalprop2I4b}. Bounding the length of the integration-domain in $\mathrm{RHS}\eqref{eq:finalprop2I4a}$ gives the second factor in \eqref{eq:finalprop2I4b}. To get \eqref{eq:finalprop2I4c}-\eqref{eq:finalprop2I4d}, we replace $\mathsf{G}^{\t}\mapsto\mathbf{G}^{\t,\mathrm{mod}}$ in \eqref{eq:finalprop2I4b}. (This gives \eqref{eq:finalprop2I4d}; the cost is \eqref{eq:finalprop2I4c}.) Now, on $\mathcal{E}^{\mathrm{vhp}}$, we know \eqref{eq:finalprop2I1} holds. So on this event, we have $\eqref{eq:finalprop2I4c}\lesssim\N^{-1}\mathrm{LHS}\eqref{eq:finalprop2I4a}$. This implies that we can move \eqref{eq:finalprop2I4c} to $\mathrm{LHS}\eqref{eq:finalprop2I4a}$ and deduce $\mathrm{LHS}\eqref{eq:finalprop2I4a}\lesssim\eqref{eq:finalprop2I4d}$. Now, recall that we have restricted to $\s\leq\t_{\mathrm{st}}$ in \eqref{eq:finalprop2Ia}-\eqref{eq:finalprop2Ib}. By Definition \ref{definition:method8}, we know the first factor in \eqref{eq:finalprop2I4d} is $\lesssim\N^{\gamma_{\mathrm{ap}}}$ with probability 1. (Indeed, because $\lambda(\t)$ is uniformly bounded above and {from} below away from zero, by Definition \ref{definition:method8}, we know $\mathbf{G}^{\t,\mathrm{mod}}$ at time $\s\leq\t_{\mathrm{st}}$ is the product of exponentials of terms that are all $\lesssim\log\log\N$.) Thus, on $\mathcal{E}^{\mathrm{vhp}}$, we have the upper bound $\eqref{eq:finalprop2I4d}\lesssim\N^{-{\mathrm{D}}}$ {for any large but fixed $\mathrm{D}>0$}. Therefore, we deduce the following on $\mathcal{E}^{\mathrm{vhp}}$:
\begin{align}
|\mathbf{G}^{\t,\mathrm{mod}}(\s-\r,\x(\s,\r,\mathrm{k}))-\mathsf{G}^{\t}(\s-\r,\x(\s,\r,\mathrm{k}))| \ \lesssim \ \N^{-{\mathrm{D}}}. \label{eq:finalprop2I4e}
\end{align}
We will now move to the $\mathds{R}^{\chi,\mathfrak{q},\pm,\mathrm{j}}(\s-\r,\z(\s-\r)\pm2\mathrm{k}\mathfrak{l}(\mathrm{j}))$ term in \eqref{eq:finalprop2I3b}. By Definitions \ref{definition:bg21}, \ref{definition:bg24}, \ref{definition:bg26}, this term is a functional that has the following properties. First, take a function with support $\mathbb{I}(\mathfrak{l}(\mathrm{j}),\pm)$ from Definition \ref{definition:bg21}. Then, precompose it with the shift $\mathbf{U}(\cdot)\mapsto\mathbf{U}(\cdot+\z(\s-\r)\pm2\mathrm{k}\mathfrak{l}(\mathrm{j}))$. The term $\mathds{R}^{\chi,\mathfrak{q},\pm,\mathrm{j}}(\s-\r,\z(\s-\r)\pm2\mathrm{k}\mathfrak{l}(\mathrm{j}))$ is then a functional of this form evaluated at $\mathbf{U}^{\s-\r,\cdot}$. Denote the function (that we evaluate at $\mathbf{U}^{\s-\r,\cdot}$) by {\small$\mathfrak{A}^{(\mathrm{k}),(1),\pm}_{\s-\r}$}. In this language, for later convenience, we write
\begin{align}
&\mathscr{A}^{\pm,\t,\mathrm{mod}} \nonumber\\
&:= \ \tau(\mathrm{j},\mathrm{i})^{-1}{\textstyle\int_{0}^{\tau(\mathrm{j},\mathrm{i})}}\{\mathfrak{m}(\mathrm{j},\mathrm{i})^{-1}{\textstyle\sum_{\mathrm{k}=0}^{\mathfrak{m}(\mathrm{j},\mathrm{i})-1}}\mathfrak{A}^{(\mathrm{k}),(1),\pm}_{\s-\r}(\mathbf{U}^{\s-\r,\cdot})\mathbf{G}^{\t,\mathrm{mod}}(\s-\r,\z(\s-\r)\pm2\mathrm{k}\mathfrak{l}(\mathrm{j}))\}\d\r. \label{eq:finalprop2I5}
\end{align}
We now list more properties of {\small$\mathfrak{A}^{(\mathrm{k}),(1),\pm}$}. First, by Definition \ref{definition:bg26}, we know {\small$|\mathfrak{A}^{(\mathrm{k}),(1),\pm}|\lesssim\N^{20\gamma_{\mathrm{reg}}}\mathfrak{l}(\mathrm{j}-1)^{-3/2}\lesssim\N^{25\gamma_{\mathrm{reg}}}\mathfrak{l}(\mathrm{j})^{-3/2}$}. We claim that the support of {\small$\mathfrak{A}^{(\mathrm{k}),(1),\pm}_{\tau}$} is a $\tau$-dependent shift $\mathbb{J}(\tau,\mathrm{k})$ of an interval $|\mathbb{I}(\mathrm{k})|\leq\mathfrak{l}(\mathrm{j})$. This was justified in the previous paragraph. Next, we claim that $\mathbb{J}(\tau,\mathrm{k})$ are mutually disjoint. Indeed, these are discrete intervals of length $\leq\mathfrak{l}(\mathrm{j})$ that are shifted by distinct multiples of $2\mathfrak{l}(\mathrm{j})$, as justified in the previous paragraph. Now, we claim {\small$\E^{\rho,\tau,\mathbb{J}(\tau,\mathrm{k})}\mathfrak{A}^{(\mathrm{k}),(1),\pm}_{\tau}=0$} for any $\tau,\mathrm{k}$. This is explained before Lemma \ref{lemma:bg27}. We now claim that the support $\mathbb{J}(\tau,\mathrm{k})$ is either contained in or disjoint from $\llbracket\inf\mathbb{I}(\mathfrak{t}),\z(\tau)\pm2\mathrm{k}\mathfrak{l}(\mathrm{j})\rrbracket$. To see this, first recall from the previous paragraph that {\small$\mathfrak{A}^{(\mathrm{k}),(1),\pm}_{\tau}$} has support given by the shift $\z(\tau)\pm2\mathrm{k}\mathfrak{l}(\mathrm{j})+\mathbb{I}(\mathfrak{l}(\mathrm{j}),\pm)$. (There are only two choices of combinations of signs here.) Next, recall from Definition \ref{definition:bg21} that $\mathbb{I}(\mathfrak{l}(\mathrm{j}),\pm)$ either has infimum $\geq1$ or supremum $\leq0$. So, the only way our claim may be false is if $\z(\tau)\pm2\mathrm{k}\mathfrak{l}(\mathrm{j})+\mathbb{I}(\mathfrak{l}(\mathrm{j}),\pm)$ goes below $\inf\mathbb{I}(\mathfrak{t})$. To rule this out, we recall $\z(\tau)\pm2\mathrm{k}\mathfrak{l}(\mathrm{j})+\mathbb{I}(\mathfrak{l}(\mathrm{j}),\pm)\subseteq\mathbb{I}$, as it is a shift of $0$ by length $\leq\N^{3/2+\gamma_{\mathrm{ap}}}\tau(\mathrm{j},\mathrm{i})+3\mathfrak{m}(\mathrm{j},\mathrm{i})\mathfrak{l}(\mathrm{j})\leq|\mathbb{I}|$; refer to the statement of {Lemma \ref{lemma:finalprop2}} for this last bound. It now suffices to note that $|\mathbb{I}(\mathfrak{l}(\mathrm{j}),\pm)|\lesssim\mathfrak{l}(\mathrm{j})$, which is much smaller than the distance between $\mathbb{I}$ and the boundary of $\mathbb{I}(\mathfrak{t})$. (Indeed, $\mathbb{I}(\mathfrak{t})$ is a neighborhood of $\mathbb{I}$ of radius $\gtrsim\N^{\gamma_{\mathrm{ap}}}\mathfrak{m}(\mathrm{j},\mathrm{i})\mathfrak{l}(\mathrm{j})\gg\mathfrak{l}(\mathrm{j})$; see Definitions \ref{definition:le10} and the statement of {Lemma \ref{lemma:finalprop2}}.) Ultimately, by this paragraph, we deduce that {\small$\mathfrak{A}^{(\mathrm{k}),(1),\pm}$} satisfies the properties that we claimed in {Lemma \ref{lemma:finalprop2}}. We have shown and mentioned this for now because of convenience; we do not use it until later.

Now, we give one last property of {\small$\mathfrak{A}^{(\mathrm{k}),(1),\pm}$}. First, recall that $\s\leq\t_{\mathrm{st}}$. For any $\x\in\mathbb{T}(\N)$ and $\mathbf{U}\in\R^{\mathbb{T}(\N)}$, we claim
\begin{align}
|\partial_{\mathbf{U}(\x)}\mathfrak{A}^{(\mathrm{k}),(1),\pm}_{\s-\r}(\mathbf{U})| \ \lesssim \ \N^{\mathrm{O}(1)}. \label{eq:finalprop2I6}
\end{align}
\eqref{eq:finalprop2I6} is fairly elementary but uninteresting to show, so we defer {it} to the end of this proof and assume that it is true for now. Next, we recall that on the very high probability event $\mathcal{E}^{\mathrm{vhp}}$, the estimate \eqref{eq:finalprop2I1} holds. We also recall from the previous paragraph that the support of {\small$\mathfrak{A}^{(\mathrm{k}),(1),\pm}$} lives in $\mathbb{I}$. (In particular, it makes sense to evaluate it at $\mathbf{U}\subseteq\R^{\mathbb{I}(\mathfrak{t})}$.) By this, \eqref{eq:finalprop2I6}, and calculus, we get the following on $\mathcal{E}^{\mathrm{vhp}}$ {(where $\mathrm{D}>0$ denotes a possibly different arbitrarily large but fixed constant)}:
\begin{align}
|\mathfrak{A}^{(\mathrm{k}),(1),\pm}_{\s-\r}(\mathbf{U}^{\s-\r,\cdot})-\mathfrak{A}^{(\mathrm{k}),(1),\pm}_{\s-\r}(\mathbf{U}^{\s-\r,\cdot}[\mathbb{I}(\mathfrak{t})])| \ \lesssim \ \N^{\mathrm{O}(1)}\mathrm{LHS}\eqref{eq:finalprop2I1} \ \lesssim \ \N^{-{\mathrm{D}}}. \label{eq:finalprop2I7}
\end{align}
%
\subsubsection*{Step 4: Turning $\mathbf{1}(\mathcal{E}^{\mathrm{vhp}})\mathbf{1}(\s\leq\t_{\mathrm{st}})$ into $\mathbf{1}[\mathcal{E}^{\mathrm{ap},\mathrm{j},\mathrm{i}}(\s)]\mathbf{1}[\mathcal{E}^{\mathrm{dens},\mathrm{j},\mathrm{i}}(\s)]$}
Because $\s\leq\t_{\mathrm{st}}$, we know that $\s-\tau(\mathrm{j},\mathrm{i})\leq\t_{\mathrm{st}}\leq\t_{\mathrm{reg}}$; see Definition \ref{definition:method8}. By Definition \ref{definition:reg} and Remark \ref{remark:intro14}, we deduce that the average of $\mathbf{U}^{\s-\tau(\mathrm{j},\mathrm{i}),\x}$ over $\x\in\mathbb{I}(\mathfrak{t})$ is $\lesssim\N^{\gamma_{\mathrm{reg}}}|\mathbb{I}(\mathfrak{t})|^{-1/2}$. Thus, the average (or charge density) of $\mathbf{U}^{\s-\tau(\mathrm{j},\mathrm{i}),\x}[\mathbb{I}(\mathfrak{t})]$ over $\x\in\mathbb{I}(\mathfrak{t})$ is $\lesssim\N^{\gamma_{\mathrm{reg}}}|\mathbb{I}(\mathfrak{t})|^{-1/2}$. Since charge density is conserved in time, we have $\mathbf{1}(\s\leq\t_{\mathrm{st}})\leq\mathbf{1}[\mathcal{E}^{\mathrm{dens},\mathrm{j},\mathrm{i}}(\s)]$. On the other hand, given $\s\leq\t_{\mathrm{st}}$, the $\mathbf{G}^{\t,\mathrm{mod}}$-factors in \eqref{eq:finalprop2I3b} are $\lesssim\N^{\gamma_{\mathrm{ap}}}$; this was explained in the paragraph before \eqref{eq:finalprop2I4e}. Also, on $\mathcal{E}^{\mathrm{vhp}}$, \eqref{eq:finalprop2I4e} holds. Combining the previous two sentences with the triangle inequality shows that the $\mathsf{G}^{\t}$-factors in \eqref{eq:finalprop2Ic} are $\lesssim\N^{\gamma_{\mathrm{ap}}}$ on the event where $\mathbf{1}(\mathcal{E}^{\mathrm{vhp}})\mathbf{1}(\s\leq\t_{\mathrm{st}})=1$. In particular, we have the inequality $\mathbf{1}(\mathcal{E}^{\mathrm{vhp}})\mathbf{1}(\s\leq\t_{\mathrm{st}})\leq\mathbf{1}[\mathcal{E}^{\mathrm{ap},\mathrm{j},\mathrm{i}}(\s)]$. For convenience, let us now summarize the conclusion of this step:
\begin{align}
\mathbf{1}(\mathcal{E}^{\mathrm{vhp}})\mathbf{1}(\s\leq\t_{\mathrm{st}}) \ \leq \ \mathbf{1}[\mathcal{E}^{\mathrm{ap},\mathrm{j},\mathrm{i}}(\s)]\mathbf{1}[\mathcal{E}^{\mathrm{dens},\mathrm{j},\mathrm{i}}(\s)]. \label{eq:finalprop2I7b}
\end{align}
%
\subsubsection*{Step 5: Comparing $\mathscr{A}^{\pm,\t,\mathrm{mod}}$ and $\mathscr{A}^{\mathrm{loc},(1),\pm,\t}$ on $\mathcal{E}^{\mathrm{vhp}}$}
In \eqref{eq:finalprop2I5}, let us replace {\small$\mathfrak{A}^{(\mathrm{k}),(1),\pm}(\mathbf{U}^{\s-\r,\cdot})\mapsto\mathfrak{A}^{(\mathrm{k}),(1),\pm}(\mathbf{U}^{\s-\r,\cdot}[\mathbb{I}(\mathfrak{t})])$}. Recall that $|\mathbf{G}^{\t,\mathrm{mod}}|$-factors are $\lesssim\N$ for $\s\leq\t_{\mathrm{st}}$ on $\mathcal{E}^{\mathrm{vhp}}$ by \eqref{eq:finalprop2I7b}. So, the cost in this replacement, by this $\mathbf{G}^{\t,\mathrm{mod}}$ bound and \eqref{eq:finalprop2I7}, is $\lesssim\N^{-{\mathrm{D}}}$. After this, in \eqref{eq:finalprop2I5}, we will further replace $\mathbf{G}^{\t,\mathrm{mod}}(\s-\r,\z(\s-\r)\pm2\mathrm{k}\mathfrak{l}(\mathrm{j}))\mapsto\mathsf{G}^{\t}(\s-\r,\z(\s-\r)\pm2\mathrm{k}\mathfrak{l}(\mathrm{j}))$. Recall {\small$|\mathfrak{A}^{(\mathrm{k}),(1),\pm}|\lesssim\N^{25\gamma_{\mathrm{reg}}}\mathfrak{l}(\mathrm{j})^{-3/2}\lesssim\N$} from right after \eqref{eq:finalprop2I5}. Thus, by \eqref{eq:finalprop2I4e}, the cost in this replacement is $\lesssim\N^{-{\mathrm{D}}}$. (We clarify that this paragraph is on $\mathcal{E}^{\mathrm{vhp}}$.) Therefore, on this event, we deduce the estimate below {for a possibly different arbitrarily large but fixed constant $\mathrm{D}>0$}:
\begin{align}
|\mathscr{A}^{\pm,\t,\mathrm{mod}}-\mathscr{A}^{\mathrm{loc},(1),\pm,\t}| \ \lesssim \ \N^{-{\mathrm{D}}}, \label{eq:finalprop2I8}
\end{align}
in which $\mathscr{A}^{\mathrm{loc},(1),\pm,\t}$ is from the statement of this lemma. We now use \eqref{eq:finalprop2I8} to show two preliminary bounds. First, we claim the following calculations in which we restrict to $\s\leq\t_{\mathrm{st}}$ and we condition on $\mathcal{E}^{\mathrm{vhp}}$: 
\begin{align}
&\N^{\beta(\mathrm{j},\mathrm{i}-1)}|\mathscr{A}^{\pm,\t}|^{2}\mathbf{1}[|\mathscr{A}^{\pm,\t}|\lesssim\N^{-\beta(\mathrm{j},\mathrm{i}-1)}]\nonumber\\
&= \ \N^{\beta(\mathrm{j},\mathrm{i}-1)}\exp[2\lambda(\t)\mathbf{J}^{!}-2\lambda(\t)\mathrm{U}]|\mathscr{A}^{\pm,\t,\mathrm{mod}}|^{2}\mathbf{1}[|\mathscr{A}^{\pm,\t}|\lesssim\N^{-\beta(\mathrm{j},\mathrm{i}-1)}] \label{eq:finalprop2I9a}\\
&\lesssim \ \N^{\beta(\mathrm{j},\mathrm{i}-1)+2\gamma_{\mathrm{ap}}}|\mathscr{A}^{\pm,\t,\mathrm{mod}}|^{2}\mathbf{1}[|\mathscr{A}^{\pm,\t}|\lesssim\N^{-\beta(\mathrm{j},\mathrm{i}-1)}] \label{eq:finalprop2I9b}\\
&\lesssim \ \N^{\beta(\mathrm{j},\mathrm{i}-1)+2\gamma_{\mathrm{ap}}}|\mathscr{A}^{\pm,\t,\mathrm{mod}}|^{2}\mathbf{1}[|\mathscr{A}^{\pm,\t,\mathrm{mod}}|\lesssim\N^{-\beta(\mathrm{j},\mathrm{i}-1)+\gamma_{\mathrm{ap}}}] \label{eq:finalprop2I9c}\\
&\lesssim \ \N^{\beta(\mathrm{j},\mathrm{i}-1)+2\gamma_{\mathrm{ap}}}|\mathscr{A}^{\mathrm{loc},(1),\pm,\t}|^{2}\mathbf{1}[|\mathscr{A}^{\mathrm{loc},(1),\pm,\t}|\lesssim\N^{-\beta(\mathrm{j},\mathrm{i}-1)+\gamma_{\mathrm{ap}}}]+\N^{-{\mathrm{D}}}.\label{eq:finalprop2I9d}
\end{align}
\eqref{eq:finalprop2I9a} follows by \eqref{eq:finalprop2I3a}. To show {\eqref{eq:finalprop2I9b}}, we first note $|\lambda(\t)\mathrm{U}|\lesssim1$. Moreover, $\exp[2\lambda(\t)\mathbf{J}^{!}]\lesssim\N^{\gamma_{\mathrm{ap}}}$. Indeed, $\mathbf{J}^{!}$ is evaluated before time $\s\leq\t_{\mathrm{st}}$. By Definition \ref{definition:method8}, this implies $|\mathbf{J}^{!}(\s-\tau(\mathrm{j},\mathrm{i}),\x)|\lesssim\log\log\N$ because $\lambda(\tau)$ is bounded uniformly {from} above and {from} below away from zero; for this, see Assumption \ref{ass:intro8}. Exponentiating $\mathrm{O}(\log\log\N)$ gives the $\exp[2\lambda(\t)\mathbf{J}^{!}]\lesssim\N^{\gamma_{\mathrm{ap}}}$ upper bound. \eqref{eq:finalprop2I9b} therefore holds. By the same token, we also know that $\exp[\lambda(\t)\mathbf{J}^{!}-\lambda(\t)\mathrm{U}]\gtrsim\N^{\gamma_{\mathrm{ap}}}$. Thus, by \eqref{eq:finalprop2I3a}, we deduce that {\small$|\mathscr{A}^{\pm,\t,\mathrm{mod}}|\lesssim\N^{\gamma_{\mathrm{ap}}}|\mathscr{A}^{\pm,\t}|$}. This derives \eqref{eq:finalprop2I9c} from \eqref{eq:finalprop2I9b}. To get \eqref{eq:finalprop2I9d}, we first use the estimate \eqref{eq:finalprop2I8}. This lets us control {\small$|\mathscr{A}^{\pm,\t,\mathrm{mod}}|^{2}\lesssim|\mathscr{A}^{\mathrm{loc},(1),\pm,\t}|^{2}+\N^{-999}\{|\mathscr{A}^{\pm,\t,\mathrm{mod}}|+|\mathscr{A}^{\mathrm{loc},(1),\pm,\t}|\}\lesssim|\mathscr{A}^{\mathrm{loc},(1),\pm,\t}|^{2}+\N^{-999}|\mathscr{A}^{\pm,\t,\mathrm{mod}}|+\N^{-999}$}; the last estimate follows by another application of \eqref{eq:finalprop2I8}. After multiplying by the indicator of {\small$|\mathscr{A}^{\pm,\t,\mathrm{mod}}|\lesssim\N^{-\beta(\mathrm{j},\mathrm{i}-1)+\gamma_{\mathrm{ap}}}$}, the previous bound now provides for us {\small$|\mathscr{A}^{\pm,\t,\mathrm{mod}}|^{2}\lesssim|\mathscr{A}^{\mathrm{loc},(1),\pm,\t}|^{2}+\N^{-{\mathrm{D}}}\N^{-\beta(\mathrm{j},\mathrm{i}-1)+\gamma_{\mathrm{ap}}}$}. We again use \eqref{eq:finalprop2I8} to replace the indicator in \eqref{eq:finalprop2I9c} with \eqref{eq:finalprop2I9d}. (For this, we note $\N^{-{\mathrm{D}}}\ll\N^{-\beta(\mathrm{j},\mathrm{i}-1)+\gamma_{\mathrm{ap}}}$.) This derives \eqref{eq:finalprop2I9d} from \eqref{eq:finalprop2I9c}. By an almost identical argument, we also deduce the following estimate, again for $\s\leq\t_{\mathrm{st}}$ and upon conditioning on $\mathcal{E}^{\mathrm{vhp}}$:
\begin{align}
\N^{-\beta(\mathrm{j},\mathrm{i}-1)}\mathbf{1}[|\mathscr{A}^{\pm,\t}|>\N^{-\beta(\mathrm{j},\mathrm{i}-1)}] \ \lesssim \ \N^{-\beta(\mathrm{j},\mathrm{i}-1)}\mathbf{1}[|\mathscr{A}^{\mathrm{loc},(1),\pm,\t}|\gtrsim\N^{-\beta(\mathrm{j},\mathrm{i}-1)-\gamma_{\mathrm{ap}}}]. \label{eq:finalprop2I9e} 
\end{align}
We use \eqref{eq:finalprop2I9a}-\eqref{eq:finalprop2I9e} with \eqref{eq:finalprop2I7b} to get the following, in which we view the $\mathrm{Loc}^{(1)}$-terms in \eqref{eq:finalprop2I10a}-\eqref{eq:finalprop2I10b} from the statement of the lemma as functions of the local process $\tau\mapsto(\mathbf{U}^{\tau,\cdot}[\mathbb{I}(\mathfrak{t})],\mathtt{J}(\tau,\cdot;\mathbb{I}(\mathfrak{t})))$ with initial data given in Step 1 of this proof:
\begin{align}
\E[\mathbf{1}(\mathcal{E}^{\mathrm{vhp}})\mathbf{1}(\s\leq\t_{\mathrm{st}})\mathrm{Cent}\Upsilon^{\pm,\mathrm{j},\mathrm{i},\t}(\s,\z(\s))] \ &\lesssim \ \N^{2\gamma_{\mathrm{ap}}}\E[\mathrm{Loc}^{(1)}\Upsilon^{\pm,\mathrm{j},\mathrm{i},\t}(\s)]+\N^{-{\mathrm{D}}} \label{eq:finalprop2I10a} \\
\E[\mathbf{1}(\mathcal{E}^{\mathrm{vhp}})\mathbf{1}(\s\leq\t_{\mathrm{st}})\mathrm{Cent}\Lambda^{\pm,\mathrm{j},\mathrm{i},\t}(\s,\z(\s))] \ &\lesssim \ \N^{2\gamma_{\mathrm{ap}}}\E[\mathrm{Loc}^{(1)}\Lambda^{\pm,\mathrm{j},\mathrm{i},\t}(\s)]+\N^{-{\mathrm{D}}}. \label{eq:finalprop2I10b}
\end{align}
On the other hand, because $\mathcal{E}^{\mathrm{vhp}}$ is very high probability, and because $|\mathrm{Cent}\Upsilon^{\pm,\mathrm{j},\mathrm{i},\t}|+|\mathrm{Cent}\Lambda^{\pm,\mathrm{j},\mathrm{i},\t}|\lesssim1$ (which can be directly verified by \eqref{eq:bg213101IIa}-\eqref{eq:bg213101IIb}) we have the following estimate on the complement of $\mathcal{E}^{\mathrm{vhp}}$:
\begin{align}
&\E[\{1-\mathbf{1}(\mathcal{E}^{\mathrm{vhp}})\}\mathbf{1}(\s\leq\t_{\mathrm{st}})\mathrm{Cent}\Upsilon^{\pm,\mathrm{j},\mathrm{i},\t}(\s,\z(\s))]\nonumber\\
&+\E[\{1-\mathbf{1}(\mathcal{E}^{\mathrm{vhp}})\}\mathbf{1}(\s\leq\t_{\mathrm{st}})\mathrm{Cent}\Lambda^{\pm,\mathrm{j},\mathrm{i},\t}(\s,\z(\s))] \nonumber\\
&\lesssim \ \N^{-9999}. \label{eq:finalprop2I10c}
\end{align}
To conclude this step, we combine \eqref{eq:finalprop2I10a}, \eqref{eq:finalprop2I10b}, and \eqref{eq:finalprop2I10c}. This gives
\begin{align}
\mathrm{LHS}\eqref{eq:finalprop2Ia} \ &\lesssim \ \N^{2\gamma_{\mathrm{ap}}}\E[\mathrm{Loc}^{(1)}\Upsilon^{\pm,\mathrm{j},\mathrm{i},\t}(\s)]+\N^{-{\mathrm{D}}} \label{eq:finalprop2I11a}\\
\mathrm{LHS}\eqref{eq:finalprop2Ib} \ &\lesssim \ \N^{2\gamma_{\mathrm{ap}}}\E[\mathrm{Loc}^{(1)}\Lambda^{\pm,\mathrm{j},\mathrm{i},\t}(\s)]+\N^{-{\mathrm{D}}}. \label{eq:finalprop2I11b}
\end{align}
%
\subsubsection*{Step 6: Factorizing the $\E$-expectation}
Observe that $\E$ in $\mathrm{RHS}\eqref{eq:finalprop2I11a}$ and $\mathrm{RHS}\eqref{eq:finalprop2I11b}$ are with respect to the law of the joint process $\tau\mapsto(\mathbf{U}^{\tau,\cdot}[\mathbb{I}(\mathfrak{t})],\mathtt{J}(\tau,\cdot;\mathbb{I}(\mathfrak{t})))$. Recall the initial data in step 1 of this argument. We can rewrite expectations in $\mathrm{RHS}\eqref{eq:finalprop2I11a}$ and $\mathrm{RHS}\eqref{eq:finalprop2I11b}$ by conditioning on the initial data of the joint process, evaluating expectation of $\mathrm{Loc}$-terms therein with respect to the law of this joint process, and then taking expectation over the law of the initial data. Recalling that said initial data is just the projection $\Pi^{\mathrm{j},\mathrm{i}}$ applied to $\mathbf{U}^{\s-\tau(\mathrm{j},\mathrm{i}),\cdot}$, we therefore deduce $\mathrm{RHS}\eqref{eq:finalprop2I11a}=\mathrm{RHS}\eqref{eq:finalprop2Ia}$ and $\mathrm{RHS}\eqref{eq:finalprop2I11b}=\mathrm{RHS}\eqref{eq:finalprop2Ib}$. By combining these two identities with the estimates \eqref{eq:finalprop2I11a} and \eqref{eq:finalprop2I11b}, we get \eqref{eq:finalprop2Ia}-\eqref{eq:finalprop2Ib}.
\subsubsection*{Step 7: Proving \eqref{eq:finalprop2I6}}
Recall {\small$\mathfrak{A}^{(\mathrm{k}),(1),\pm}$} is the underlying function defining $\mathds{R}^{\chi,\mathfrak{q},\pm,\mathrm{j}}$ from Definition \ref{definition:bg26}. By Definitions \ref{definition:bg21}, \ref{definition:bg24}, \ref{definition:bg26}, we deduce {\small$\mathfrak{A}^{(\mathrm{k}),(1),\pm}$} is given by $\leq10$-many compositions, products, and sums of functionals with the following form. 
\begin{enumerate}
\item The function {\small$\mathbf{U}\mapsto|\mathbb{K}|^{-1}\sum_{\x\in\mathbb{K}}\mathbf{U}(\x)$} on $\mathbf{U}\in\R^{\mathbb{T}(\N)}$ for a subset $\mathbb{K}\subseteq\mathbb{T}(\N)$. This satisfies the gradient bound \eqref{eq:finalprop2I6}.
\item The function $\sigma\mapsto\E^{\sigma,\tau,\mathbb{K}}\mathsf{F}$, where $\tau\lesssim1$ and {\small$\mathsf{F}(\mathbf{U})\lesssim\N+\N\sum_{\x\in\mathbb{K}}|\mathbf{U}(\x)|^{100}$}. (This polynomial estimate for $\mathsf{F}$ follows from our assumption on $\mathfrak{q}$ in Proposition \ref{prop:bg22}.) Here, $\sigma$ is always given by average of $\mathbf{U}^{\tau,\x}$ for $\tau\leq\t_{\mathrm{st}}$ and over $\x\in\mathbb{K}$ for some subset $|\mathbb{K}|\gtrsim\N^{1/999}$. Thus, we deduce $|\sigma|\lesssim1$. Following the last paragraph in the proofs of Lemmas \ref{lemma:bg23} and \ref{lemma:bg1hl2} (see after Lemma \ref{lemma:ee1}), we deduce that the $\sigma$-derivative of this function is uniformly $\lesssim\N$.
\item The function $\eta\mapsto\chi(\eta)$, where $\chi$ is smooth and compactly supported with derivative $|\partial_{\eta}\chi|\lesssim\N^{2}$. (See Definition \ref{definition:bg26}.)
\end{enumerate}
With the Leibniz and chain rules from calculus and the previous bullet points, \eqref{eq:finalprop2I6} becomes elementary to check.
\subsubsection*{Step 8: Removing $\y(\s-\tau(\mathrm{j},\mathrm{i}))=0$}
When removing this assumption, all that we have to do is shift everything in space. For example, we must replace $\mathbb{I}$ by {a} $(\y,\s)$-dependent shift, and thus replace $\mathbb{I}(\mathfrak{t})$ by the same shift. The process $\mathbf{U}^{\tau,\cdot}$ is also shifted by the same. This explains the shift of $\y(\s-\tau(\mathrm{j},\mathrm{i}))$ in $\mathrm{RHS}\eqref{eq:finalprop2Ia}$ and $\mathrm{RHS}\eqref{eq:finalprop2Ib}$. We must now explain why the local processes $\mathbf{U}^{\tau,\cdot}[\mathbb{I}(\mathfrak{t})]$ and $\mathtt{J}(\tau,\cdot;\mathbb{I}(\mathfrak{t}))$ involve no $(\y,\s)$-dependent shift in $\mathrm{RHS}\eqref{eq:finalprop2Ia}$ and $\mathrm{RHS}\eqref{eq:finalprop2Ib}$. Indeed, the law of these processes is homogeneous in space. Thus, so are $\E^{\mathrm{loc},\s}$-expectations in \eqref{eq:finalprop2Ia} and \eqref{eq:finalprop2Ib}. This gives the proof for general $\y(\s-\tau(\mathrm{j},\mathrm{i}))$.
\end{proof}
\subsubsection{Reduction to local equilibrium}
The first term in the RHS of \eqref{eq:finalprop2Ia} and \eqref{eq:finalprop2Ib}, respectively, are local statistics of \eqref{eq:glsde}. Thus, we can reduce their estimation to local equilibrium estimates via Lemma \ref{lemma:le9}.
\begin{lemma}\label{lemma:finalprop4}
 Fix $1\leq\mathrm{j}\leq\mathrm{j}(\infty)$ and $1\leq\mathrm{i}<\mathrm{i}(\mathrm{j})$. Retain notation of {Lemma \ref{lemma:finalprop2}}. With notation explained after, we have
\begin{align}
&{\textstyle\int_{\tau(\mathrm{j},\mathrm{i})}^{1}}|\mathbb{T}(\N)|^{-1}{\textstyle\sum_{\y}}\N\cdot\E[\{\E^{\mathrm{loc},\s}[\mathrm{Loc}^{(1)}\Upsilon^{\pm,\mathrm{j},\mathrm{i},\t}(\s)]\}(\Pi^{\mathrm{j},\mathrm{i}}\mathbf{U}^{\s-\tau(\mathrm{j},\mathrm{i}),\y(\s-\tau(\mathrm{j},\mathrm{i}))+\cdot})]\d\s \nonumber\\
&\lesssim \ \N^{-30\beta_{\mathrm{BG}}}+\mathrm{LE}^{\t}(\Upsilon), \label{eq:finalprop4Ia} \\
&{\textstyle\int_{\tau(\mathrm{j},\mathrm{i})}^{1}}|\mathbb{T}(\N)|^{-1}{\textstyle\sum_{\y}}\N\cdot\E[\{\E^{\mathrm{loc},\s}[\mathrm{Loc}^{(1)}\Lambda^{\pm,\mathrm{j},\mathrm{i},\t}(\s)]\}(\Pi^{\mathrm{j},\mathrm{i}}\mathbf{U}^{\s-\tau(\mathrm{j},\mathrm{i}),\y(\s-\tau(\mathrm{j},\mathrm{i}))+\cdot})]\d\s \nonumber\\
&\lesssim \ \N^{-30\beta_{\mathrm{BG}}}+\mathrm{LE}^{\t}(\Lambda). \label{eq:finalprop4Ib}
\end{align}
Above, the local equilibrium error terms $\mathrm{LE}^{\t}(\Upsilon)$ and $\mathrm{LE}^{\t}(\Lambda)$ are defined below (with more notation explained after):
\begin{align}
\mathrm{LE}^{\t}(\Upsilon) \ &:= \ \sup_{\sigma\in\R}\sup_{\tau(\mathrm{j},\mathrm{i})\leq\s\leq1}\N\times\E^{\sigma,\mathrm{loc},\s}\mathrm{Loc}^{(1)}\Upsilon^{\pm,\mathrm{j},\mathrm{i},\t}(\s) \label{eq:finalprop4IIa}\\
\mathrm{LE}^{\t}(\Lambda) \ &:= \ \sup_{\sigma\in\R}\sup_{\tau(\mathrm{j},\mathrm{i})\leq\s\leq1}\N\times\E^{\sigma,\mathrm{loc},\s}\mathrm{Loc}^{(1)}\Lambda^{\pm,\mathrm{j},\mathrm{i},\t}(\s). \label{eq:finalprop4IIb}
\end{align}
Above, $\E^{\sigma,\mathrm{loc},\s}$ is $\E^{\mathrm{loc},\s}$ in {Lemma \ref{lemma:finalprop2}} but with canonical initial data $\mathbf{U}^{\s-\tau(\mathrm{j},\mathrm{i}),\cdot}[\mathbb{I}(\mathfrak{t})]\sim\mathbb{P}^{\sigma,\s-\tau(\mathrm{j},\mathrm{i}),\mathbb{I}(\mathfrak{t})}$.
\end{lemma}
\begin{proof}
We start with \eqref{eq:finalprop4Ia}. Define $\mathfrak{a}(\s,\mathbf{U}):=\E^{\mathrm{loc},\s}[\mathrm{Loc}^{(1)}\Upsilon^{\pm,\mathrm{j},\mathrm{i},\t}(\s)]\}(\Pi^{\mathrm{j},\mathrm{i}}\mathbf{U})$ for any $\mathbf{U}\in\R^{\mathbb{T}(\N)}$. It is a functional whose support equals the discrete interval $\mathbb{I}(\mathfrak{t})$ with $\mathfrak{t}=\tau(\mathrm{j},\mathrm{i})$. (Here, $\mathbb{I}$ is constructed in Lemma \ref{lemma:finalprop2}. The support claim follows since $\Pi^{\mathrm{j},\mathrm{i}}$ is projection $\R^{\mathbb{T}(\N)}\to\R^{\mathbb{I}(\mathfrak{t})}$; see Lemma \ref{lemma:finalprop2}.) We now apply Lemma \ref{lemma:le9} with $\kappa>0$ to be chosen shortly. This gives
{
\begin{align}
\mathrm{LHS}\eqref{eq:finalprop4Ia} \ &\lesssim \ \tfrac{1}{\kappa}\N^{-\frac54-\gamma_{\mathrm{KL}}}|\mathbb{I}(\mathfrak{t})|^{3} \nonumber\\
&+ \tfrac{\N}{\kappa}\sup_{\substack{\sigma\in\R\\ \tau(\mathrm{j},\mathrm{i})\leq\s\leq1}}\log\E^{\sigma,\s-\tau(\mathrm{j},\mathrm{i}),\mathbb{I}(\mathfrak{t})}\exp\{\kappa\E^{\mathrm{loc},\s}[\mathrm{Loc}^{(1)}\Upsilon^{\pm,\mathrm{j},\mathrm{i},\t}(\s)]\}(\Pi^{\mathrm{j},\mathrm{i}}\mathbf{U})\}. \label{eq:finalprop4IIa1}
\end{align}
}(Indeed, because we multiply by $\N$ in $\mathrm{LHS}\eqref{eq:finalprop4IIa}$, we must also multiply the RHS of the bound in Lemma \ref{lemma:le9} by $\N$. Moreover, the supremum in $\mathrm{RHS}\eqref{eq:finalprop4IIa1}$ is over $\tau(\mathrm{j},\mathrm{i})\leq\s\leq1$, not $0\leq\s\leq1$ as in Lemma \ref{lemma:le9}. However, this difference is completely cosmetic, since in $\mathrm{LHS}\eqref{eq:finalprop4Ia}$, the process \eqref{eq:glsde} is shifted backwards. Thus, a simple change-of-variables gives us \eqref{eq:finalprop4IIa1}.) We now estimate the double supremum in \eqref{eq:finalprop4IIa1}. By construction in \eqref{eq:finalprop2IIa}, we deduce that {\small$\E^{\mathrm{loc},\s}[\mathrm{Loc}^{(1)}\Upsilon^{\pm,\mathrm{j},\mathrm{i},\t}(\s)]\}(\Pi^{\mathrm{j},\mathrm{i}}\mathbf{U})$} has a deterministic upper bound of $\N^{-\beta(\mathrm{j},\mathrm{i}-1)+2\gamma_{\mathrm{ap}}}$. (Indeed, it is an expectation of \eqref{eq:finalprop2IIa}, which has the form $\N^{\beta(\mathrm{j},\mathrm{i}-1)}\mathrm{X}^{2}\mathbf{1}[|\mathrm{X}|\lesssim\N^{-\beta(\mathrm{j},\mathrm{i}-1)+\gamma_{\mathrm{ap}}}]$.) This motivates us to take $\kappa=\N^{\beta(\mathrm{j},\mathrm{i}-1)-2\gamma_{\mathrm{ap}}}$ to be the inverse of this deterministic bound. Here is the upshot. With this choice of $\kappa$, we know the term inside the exponential in the last term in \eqref{eq:finalprop4IIa1} is $\mathrm{O}(1)$ with probability 1. With this, we claim the following estimate, which we explain afterwards, in which $\mathrm{LD}$ denotes the last term in \eqref{eq:finalprop4IIa1} (just for convenience):
\begin{align}
\mathrm{LD} \ &\lesssim \ \tfrac{\N}{\kappa}{\textstyle\sup_{\sigma,\s}}\log\E^{\sigma,\s-\tau(\mathrm{j},\mathrm{i}),\mathbb{I}(\mathfrak{t})}\{1+\mathrm{O}(\kappa\E^{\mathrm{loc},\s}[\mathrm{Loc}^{(1)}\Upsilon^{\pm,\mathrm{j},\mathrm{i},\t}(\s)]\}(\Pi^{\mathrm{j},\mathrm{i}}\mathbf{U}))\} \label{eq:finalprop4IIa2a}\\
&\lesssim \ \tfrac{\N}{\kappa}{\textstyle\sup_{\sigma,\s}}\log\{1+\E^{\sigma,\s-\tau(\mathrm{j},\mathrm{i}),\mathbb{I}(\mathfrak{t})}\mathrm{O}(\kappa\E^{\mathrm{loc},\s}[\mathrm{Loc}^{(1)}\Upsilon^{\pm,\mathrm{j},\mathrm{i},\t}(\s)]\}(\Pi^{\mathrm{j},\mathrm{i}}\mathbf{U}))\}\label{eq:finalprop4IIa2b}\\
&\lesssim \ \N\times{\textstyle\sup_{\sigma,\s}}\E^{\sigma,\s-\tau(\mathrm{j},\mathrm{i}),\mathbb{I}(\mathfrak{t})}\E^{\mathrm{loc},\s}[\mathrm{Loc}^{(1)}\Upsilon^{\pm,\mathrm{j},\mathrm{i},\t}(\s)]\}(\Pi^{\mathrm{j},\mathrm{i}}\mathbf{U}) \ = \ \mathrm{LE}^{\t}(\Upsilon). \label{eq:finalprop4IIa2c}
\end{align}
\eqref{eq:finalprop4IIa2a} follows from $\exp[\mathrm{a}]\leq1+\mathrm{O}(|\mathrm{a}|)$ for any $|\mathrm{a}|\lesssim1$, which follows from smoothness of the exponential and $\exp[0]=1$. (The implied constant in the big-Oh depends on our a priori upper bound on $\mathrm{a}$ itself.) \eqref{eq:finalprop4IIa2b} follows via linearity of expectation. \eqref{eq:finalprop4IIa2c} follows first by the estimate $\log[1+\mathrm{a}]\leq\mathrm{a}$ for any $\mathrm{a}\geq0$. Then, we note that $\E^{\sigma,\mathrm{loc},\s}$ is exactly the double expectation in \eqref{eq:finalprop4IIa2c}. This controls the last term in \eqref{eq:finalprop4IIa1}. We now control the first term in $\mathrm{RHS}\eqref{eq:finalprop4IIa1}$. Recall from Definition \ref{definition:le10} that $|\mathbb{I}(\mathfrak{t})|\lesssim\N^{\gamma_{\mathrm{ap}}}[\N\tau(\mathrm{j},\mathrm{i})^{1/2}+\N^{3/2}\tau(\mathrm{j},\mathrm{i})+|\mathbb{I}|]$ since $\mathfrak{t}=\tau(\mathrm{j},\mathrm{i})$. By Lemma \ref{lemma:finalprop2}, $|\mathbb{I}|\lesssim\mathfrak{m}(\mathrm{j},\mathrm{i})\mathfrak{l}(\mathrm{j})+\N^{3/2+\gamma_{\mathrm{ap}}}\tau(\mathrm{j},\mathrm{i})$. Also, recall our choice of $\kappa$. This lets us compute the first term in $\mathrm{RHS}\eqref{eq:finalprop4IIa1}$, which we denote by $\mathrm{Cost}$ for convenience:
\begin{align}
&\mathrm{Cost} \nonumber\\
&\lesssim \ \N^{-\frac54-\gamma_{\mathrm{KL}}+2\gamma_{\mathrm{ap}}}\N^{-\beta(\mathrm{j},\mathrm{i}-1)}\N^{3\gamma_{\mathrm{ap}}}[\N^{3}\tau(\mathrm{j},\mathrm{i})^{\frac32}+\N^{\frac92+3\gamma_{\mathrm{ap}}}\tau(\mathrm{j},\mathrm{i})^{3}+\mathfrak{m}(\mathrm{j},\mathrm{i})^{3}\mathfrak{l}(\mathrm{j})^{3}] \label{eq:finalprop4IIa3a}\\
&\lesssim \ \N^{-\frac94-\frac12\gamma_{\mathrm{KL}}}[\N^{3}\tau(\mathrm{j},\mathrm{i})^{\frac12}\mathfrak{m}(\mathrm{j},\mathrm{i})^{-1}\mathfrak{l}(\mathrm{j})^{-1}+\N^{\frac92}\tau(\mathrm{j},\mathrm{i})^{2}\mathfrak{m}(\mathrm{j},\mathrm{i})^{-1}\mathfrak{l}(\mathrm{j})^{-1}+\tau(\mathrm{j},\mathrm{i})^{-1}\mathfrak{m}(\mathrm{j},\mathrm{i})^{2}\mathfrak{l}(\mathrm{j})^{2}]. \label{eq:finalprop4IIa3b}
\end{align}
\eqref{eq:finalprop4IIa3a} follows by plugging our bounds from the previous paragraph into the first term in $\mathrm{RHS}\eqref{eq:finalprop4IIa1}$. To show \eqref{eq:finalprop4IIa3b}, we first factor out all $\gamma_{\mathrm{ap}}$-exponents from square-brackets in \eqref{eq:finalprop4IIa3a}. This raises the $\N$-exponent in the first factor in $\eqref{eq:finalprop4IIa3a}$ by $\leq10\gamma_{\mathrm{ap}}$. (We can drop all $\gamma_{\mathrm{ap}}$ exponents if we lower $\gamma_{\mathrm{KL}}\mapsto9\gamma_{\mathrm{KL}}/10$; see Definition \ref{definition:method8}.) Next, we factor $\tau(\mathrm{j},\mathrm{i})\mathfrak{m}(\mathrm{j},\mathrm{i})\mathfrak{l}(\mathrm{j})$ from square brackets in \eqref{eq:finalprop4IIa3a}. This turns the square brackets in \eqref{eq:finalprop4IIa3a} into the square brackets in \eqref{eq:finalprop4IIa3b}. The cost in this second factoring forces us to replace $\N^{-\beta(\mathrm{j},\mathrm{i}-1)}$ in \eqref{eq:finalprop4IIa3a} by $\N^{-\beta(\mathrm{j},\mathrm{i}-1)}\tau(\mathrm{j},\mathrm{i})\mathfrak{m}(\mathrm{j},\mathrm{i})\mathfrak{l}(\mathrm{j})$. But, modulo a factor $\lesssim\N^{\gamma_{\mathrm{KL}}/10}$, which we account for by further lowering $9\gamma_{\mathrm{KL}}/10\mapsto2\gamma_{\mathrm{KL}}/3$, we get $\N^{-\beta(\mathrm{j},\mathrm{i}-1)}\tau(\mathrm{j},\mathrm{i})\mathfrak{m}(\mathrm{j},\mathrm{i})\mathfrak{l}(\mathrm{j})\lesssim\N^{-1}$ by Definition \ref{definition:bg2131}, which gives \eqref{eq:finalprop4IIa3b}. (Note we paired $\beta(\mathrm{j},\mathrm{i}-1)$ with $\tau(\mathrm{j},\mathrm{i})\mathfrak{m}(\mathrm{j},\mathrm{i})\mathfrak{l}(\mathrm{j})$, not $\tau(\mathrm{j},\mathrm{i}-1)\mathfrak{m}(\mathrm{j},\mathrm{i}-1)\mathfrak{l}(\mathrm{j})$. But replacing $\tau(\mathrm{j},\mathrm{i})\mathfrak{m}(\mathrm{j},\mathrm{i})\mathfrak{l}(\mathrm{j})\mapsto\tau(\mathrm{j},\mathrm{i}-1)\mathfrak{m}(\mathrm{j},\mathrm{i}-1)\mathfrak{l}(\mathrm{j})$ has a multiplicative cost of $\lesssim\N^{\gamma_{\mathrm{KV}}/99}$, which can be ignored by further lowering $2\gamma_{\mathrm{KL}}/3\mapsto\gamma_{\mathrm{KL}}/2$.) 

We now bound \eqref{eq:finalprop4IIa3b}. First, assume $\mathfrak{m}(\mathrm{j},\mathrm{i})\mathfrak{l}(\mathrm{j})\leq\N^{1/2}$. In this case, by Definition \ref{definition:bg2131}, $\tau(\mathrm{j},\mathrm{i})=\N^{-2}\mathfrak{m}(\mathrm{j},\mathrm{i})^{2}\mathfrak{l}(\mathrm{j})^{2}$. So
\begin{align}
\eqref{eq:finalprop4IIa3b} \ &\lesssim \ \N^{-\frac94-\frac12\gamma_{\mathrm{KL}}}[\N^{3}\N^{-1}+\N^{\frac92-4+\frac32}+\N^{2}] \ = \ \N^{-\frac14-\frac12\gamma_{\mathrm{KL}}}. \label{eq:finalprop4IIa4}
\end{align}
(We emphasize that \eqref{eq:finalprop4IIa4} uses the bound $\mathfrak{m}(\mathrm{j},\mathrm{i})\mathfrak{l}(\mathrm{j})\leq\N^{1/2}$.) Now, if $\mathfrak{m}(\mathrm{j},\mathrm{i})\mathfrak{l}(\mathrm{j})\geq\N^{1/2}$, we know $\tau(\mathrm{j},\mathrm{i})=\N^{-3/2}\mathfrak{m}(\mathrm{j},\mathrm{i})\mathfrak{l}(\mathrm{j})$; see Definition \ref{definition:bg2131}. Definition \ref{definition:bg2131} also tells us $\mathfrak{m}(\mathrm{j},\mathrm{i})\mathfrak{l}(\mathrm{j})\leq\N^{3/4+\gamma_{\mathrm{KL}}/10}$. Using all of this, we deduce the following estimate:
\begin{align}
\eqref{eq:finalprop4IIa3b} \ &\lesssim \ \N^{-\frac94-\frac12\gamma_{\mathrm{KL}}}[\N^{\frac94}\mathfrak{m}(\mathrm{j},\mathrm{i})^{-\frac12}\mathfrak{l}(\mathrm{j})^{-\frac12}+\N^{\frac32}\mathfrak{m}(\mathrm{j},\mathrm{i})\mathfrak{l}(\mathrm{j})+\N^{\frac32}\mathfrak{m}(\mathrm{j},\mathrm{i})\mathfrak{l}(\mathrm{j})] \label{eq:finalprop4IIa5a}\\
&\lesssim \ \N^{-\frac94-\frac12\gamma_{\mathrm{KL}}}[\N^{\frac94}\N^{-\frac14}+\N^{\frac32}\N^{\frac34+\frac{1}{10}\gamma_{\mathrm{KL}}}] \ \lesssim \ \N^{-30\beta_{\mathrm{BG}}}. \label{eq:finalprop4IIa5b}
\end{align}
(The last bound follows since $\beta_{\mathrm{BG}}$ is equal to a small factor times $\gamma_{\mathrm{KL}}$; see Definition \ref{definition:method8}.) In any case, we deduce $\eqref{eq:finalprop4IIa3b}\lesssim\N^{-30\beta_{\mathrm{BG}}}$. Combining this with \eqref{eq:finalprop4IIa1}, \eqref{eq:finalprop4IIa2a}-\eqref{eq:finalprop4IIa2c}, and \eqref{eq:finalprop4IIa3a}-\eqref{eq:finalprop4IIa3b} completes the proof of \eqref{eq:finalprop4Ia}. To get \eqref{eq:finalprop4Ib}, the same argument works by replacing $\Upsilon\mapsto\Lambda$ formally. (Indeed, the only thing we used about $\Upsilon$ was the deterministic bound $\eqref{eq:finalprop2IIa}\lesssim\N^{-\beta(\mathrm{j},\mathrm{i}-1)+2\gamma_{\mathrm{ap}}}$, which is true for $\Lambda$ in place of $\Upsilon$; see \eqref{eq:finalprop2IIb}. Everything else in this argument was about $\mathfrak{m}(\mathrm{j},\mathrm{i}),\mathfrak{l}(\mathrm{j}),\tau(\mathrm{j},\mathrm{i}),\beta(\mathrm{j},\mathrm{i})$.)
\end{proof}
\subsubsection{The technical adjustments discussed in \emph{Remark \ref{remark:finalprop3}}}
We now use a priori bounds coming from $\mathcal{E}^{\mathrm{ap},\mathrm{j},\mathrm{i}}(\s)$ and $\mathcal{E}^{\mathrm{dens},\mathrm{j},\mathrm{i}}(\s)$ to technically adjust {\small$\mathscr{A}^{\mathrm{loc},(1),\pm,\t}$}, thereby adjusting {\small$\mathrm{Loc}^{(1)}\Upsilon^{\pm,\mathrm{j},\mathrm{i},\t}$} and {\small$\mathrm{Loc}^{(1)}\Lambda^{\pm,\mathrm{j},\mathrm{i},\t}$}; see Lemma \ref{lemma:finalprop2} for notation in this sentence. In a nutshell, these two events place us in the setting of Lemma \ref{lemma:le16}, so we can swap $\mathtt{J}$ in the constructions of Lemma \ref{lemma:finalprop2} with $\mathbf{J}$ from Definition \ref{definition:le10}. The event $\mathcal{E}^{\mathrm{ap},\mathrm{j},\mathrm{i}}(\s)$, in particular, then lets us introduce a priori estimates for the exponential of $\mathbf{J}$ from Definition \ref{definition:le10}. We will clarify more what the following lemma is saying once we have stated it precisely.
\begin{lemma}\label{lemma:finalprop5}
 Fix $1\leq\mathrm{j}\leq\mathrm{j}(\infty)$ and $1\leq\mathrm{i}<\mathrm{i}(\mathrm{j})$. Retain notation of {Lemmas \ref{lemma:finalprop2}, \ref{lemma:finalprop4}}. With notation explained after, {we have the following in which $\mathrm{D}>0$ is any large but fixed constant:}
\begin{align}
\E^{\sigma,\mathrm{loc},\s}\mathrm{Loc}^{(1)}\Upsilon^{\pm,\mathrm{j},\mathrm{i},\t}(\s) \ &\lesssim \ \mathbf{1}[|\sigma|\lesssim1]\E^{\sigma,\mathrm{loc},\s}\mathrm{Loc}\Upsilon^{\pm,\mathrm{j},\mathrm{i},\t}(\s)+\N^{-{\mathrm{D}}} \label{eq:finalprop5Ia} \\
\E^{\sigma,\mathrm{loc},\s}\mathrm{Loc}^{(1)}\Lambda^{\pm,\mathrm{j},\mathrm{i},\t}(\s) \ &\lesssim \ \mathbf{1}[|\sigma|\lesssim1]\E^{\sigma,\mathrm{loc},\s}\mathrm{Loc}\Lambda^{\pm,\mathrm{j},\mathrm{i},\t}(\s)+\N^{-{\mathrm{D}}}. \label{eq:finalprop5Ib}
\end{align}
For clarity, we emphasize that $\E^{\sigma,\mathrm{loc},\s}$ is defined in {Lemma \ref{lemma:finalprop4}}. Let us now define the RHS of \eqref{eq:finalprop5Ia} and \eqref{eq:finalprop5Ib} as follows.
\begin{enumerate}
\item Let $\mathbf{J}(\s-\r,\cdot;\mathbb{I}(\mathfrak{t}))$ be the $\mathbf{J}$-process from {Definition \ref{definition:le10}} with initial data $\mathbf{J}(\s-\tau(\mathrm{j},\mathrm{i}),\inf\mathbb{I}(\mathfrak{t});\mathbb{I}(\mathfrak{t}))=0$. (Initial data for any other $\x\in\mathbb{I}(\mathfrak{t})$ is then determined by initial data $\mathbf{U}^{\s-\tau(\mathrm{j},\mathrm{i}),\cdot}[\mathbb{I}(\mathfrak{t})]$, which, under $\E^{\sigma,\mathrm{loc},\s}$ is distributed as $\mathbb{P}^{\sigma,\s-\tau(\mathrm{j},\mathrm{i}),\mathbb{I}(\mathfrak{t})}$.
\item Next, we set the following version of $\mathsf{G}^{\t}$ in {Lemma \ref{lemma:finalprop2}} but replacing $\mathtt{J}$ by $\mathbf{J}$ from the previous bullet point and with a cutoff:
\begin{align}
&\mathbf{G}^{\t}(\tau,\x;\mathbb{I}(\mathfrak{t})) \nonumber\\
&:= \ \exp[\lambda(\t)\mathbf{J}(\tau,\x;\mathbb{I}(\mathfrak{t}))+\lambda(\t)\mathrm{U}]\mathbf{1}\{\N^{-2\gamma_{\mathrm{ap}}}\lesssim\exp[\lambda(\t)\mathbf{J}(\tau,\x;\mathbb{I}(\mathfrak{t}))+\lambda(\t)\mathrm{U}]\lesssim\N^{2\gamma_{\mathrm{ap}}}\}.
\end{align}
Again, we clarify that $\mathrm{U}$ is uniform on $[-1,1]$ and independent of everything else. We now set, with notation explained after,
\begin{align}
&\mathscr{A}^{\mathrm{loc},\pm,\t} \nonumber\\
&:= \ \tau(\mathrm{j},\mathrm{i})^{-1}{\textstyle\int_{0}^{\tau(\mathrm{j},\mathrm{i})}}\{\mathbf{G}^{\t}(\s-\r,\inf\mathbb{I}(\mathfrak{t});\mathbb{I}(\mathfrak{t}))\times\mathfrak{m}(\mathrm{j},\mathrm{i})^{-1}{\textstyle\sum_{\mathrm{k}=0}^{\mathfrak{m}(\mathrm{j},\mathrm{i})-1}}\mathfrak{A}^{(\mathrm{k}),\pm}_{\s-\r}(\mathbf{U}^{\s-\r,\cdot}[\mathbb{I}(\mathfrak{t})])\}\d\r.
\end{align}
Above, {\small$\mathfrak{A}^{(\mathrm{k}),\pm}_{\tau}(\mathbf{U})$} is a collection of functionals that satisfies the constraints of {Lemma \ref{lemma:le2}} with respect to discrete intervals $\mathbb{J}(\tau,\mathrm{k})$ that we built in {Lemma \ref{lemma:finalprop2}}. We also know the deterministic bounds {\small$|\mathfrak{A}^{(\mathrm{k}),\pm}_{\tau}(\mathbf{U})|\lesssim\N^{30\gamma_{\mathrm{reg}}}\mathfrak{l}(\mathrm{j})^{-3/2}$} for all $\tau,\mathrm{k},\mathbf{U}$.
\item Let us now define the main objects on the RHS of \eqref{eq:finalprop5Ia} and \eqref{eq:finalprop5Ib}, respectively:
\begin{align}
\mathrm{Loc}\Upsilon^{\pm,\mathrm{j},\mathrm{i},\t}(\s) \ &:= \ \N^{\beta(\mathrm{j},\mathrm{i}-1)}|\mathscr{A}^{\mathrm{loc},\pm,\t}|^{2}\mathbf{1}[|\mathscr{A}^{\mathrm{loc},\pm,\t}|\lesssim\N^{-\beta(\mathrm{j},\mathrm{i}-1)+\gamma_{\mathrm{ap}}}] \label{eq:finalprop5IIa} \\
\mathrm{Loc}\Lambda^{\pm,\mathrm{j},\mathrm{i},\t}(\s) \ &:= \ \N^{-\beta(\mathrm{j},\mathrm{i}-1)}\mathbf{1}[|\mathscr{A}^{\mathrm{loc},\pm,\t}|\gtrsim\N^{-\beta(\mathrm{j},\mathrm{i}-1)-\gamma_{\mathrm{ap}}}]. \label{eq:finalprop5IIb}
\end{align}
\end{enumerate}
\end{lemma}
\begin{rem}\label{remark:finalprop6}
 The cutoff in $\mathbf{G}^{\tau}$ comes (basically) for free because of Lemma \ref{lemma:le16} and the indicator of $\mathcal{E}^{\mathrm{ap},\mathrm{j},\mathrm{i}}(\s)\cap\mathcal{E}^{\mathrm{dens},\mathrm{j},\mathrm{i}}(\s)$. The a priori estimate $|\sigma|\lesssim1$ comes for free due to the indicator of $\mathcal{E}^{\mathrm{dens},\mathrm{j},\mathrm{i}}(\s)$ (modulo a couple {of} uninteresting technical details). Instead of discussing those in this remark, let us address the functions {\small$\mathfrak{A}^{(\mathrm{k}),\pm}$} above. The $\mathbf{G}^{\t}$ factor in the definition of $\mathscr{A}^{\mathrm{loc},\pm,\t}$ is evaluated at the spatial point $\inf\mathbb{I}(\mathfrak{t})$, not the $\mathrm{k}$-dependent points $0(\s-\r)\pm2\mathrm{k}\mathfrak{l}(\mathrm{j})$ like in the definition of {\small$\mathscr{A}^{\mathrm{loc},(1),\pm,\t}$} in Lemma \ref{lemma:finalprop2}. If we, for now, forget about the cutoff for $\mathbf{G}^{\tau}$ in Lemma \ref{lemma:finalprop5}, the cost behind replacing the spatial point is the exponential of a $\mathbf{J}$-increment, which is a (weighted) sum of $\mathbf{U}[\mathbb{I}(\mathfrak{t})]$ charges. Indeed, we ultimately let {\small$\mathfrak{A}^{(\mathrm{k}),\pm}$} be the product of this exponential with {\small$\mathfrak{A}^{(\mathrm{k}),(1),\pm}$} in Lemma \ref{lemma:finalprop2}. Verifying its properties in Lemma \ref{lemma:finalprop5} amounts to fairly direct reasoning; see the proof below.
\end{rem}
\begin{proof}
We split this argument into several steps. In what follows, we let $\mathcal{E}^{?,\mathrm{j},\mathrm{i}}:=\mathcal{E}^{?,\mathrm{j},\mathrm{i}}(\s)$ for $?=\mathrm{dens},\mathrm{ap}$ out of convenience.
\subsubsection*{Step 0: The a priori estimate for $\sigma$}
Recall $\mathcal{E}^{\mathrm{dens},\mathrm{j},\mathrm{i}}$ from Lemma \ref{lemma:finalprop2}. On this event, the average of $\mathbf{U}^{\tau,\x}[\mathbb{I}(\mathfrak{t})]$ over $\x\in\mathbb{I}(\mathfrak{t})$ is $\lesssim\N^{\gamma_{\mathrm{reg}}}|\mathbb{I}(\mathfrak{t})|^{-1/2}$ for all $\tau=\s-\r$ with $\r\in[0,\tau(\mathrm{j},\mathrm{i})]$. By definition of charge density, this is exactly $\sigma$ in \eqref{eq:finalprop5Ia}. Now, observe $|\mathbb{I}(\mathfrak{t})|\geq|\mathbb{I}|$, since $\mathbb{I}(\mathfrak{t})$ is a neighborhood of $\mathbb{I}$. Also, $|\mathbb{I}|\geq\mathfrak{m}(\mathrm{j},\mathrm{i})\mathfrak{l}(\mathrm{j})\geq\N^{1/9}$. This first bound follows by construction in Lemma \ref{lemma:finalprop2}; the second follows by Definition \ref{definition:bg2131}. Thus, because $\gamma_{\mathrm{reg}}$ {is small} (see Definition \ref{definition:reg}), we deduce $|\sigma|\lesssim1$ in \eqref{eq:finalprop5Ia}, so
\begin{align}
\E^{\sigma,\mathrm{loc},\s}\mathrm{Loc}^{(1)}\Upsilon^{\pm,\mathrm{j},\mathrm{i},\t}(\s) \ &\leq \ \mathbf{1}[|\sigma|\lesssim1]\E^{\sigma,\mathrm{loc},\s}\mathrm{Loc}^{(1)}\Upsilon^{\pm,\mathrm{j},\mathrm{i},\t}(\s) \label{eq:finalprop5Ia0} \\
\E^{\sigma,\mathrm{loc},\s}\mathrm{Loc}^{(1)}\Lambda^{\pm,\mathrm{j},\mathrm{i},\t}(\s) \ &\leq \ \mathbf{1}[|\sigma|\lesssim1]\E^{\sigma,\mathrm{loc},\s}\mathrm{Loc}^{(1)}\Lambda^{\pm,\mathrm{j},\mathrm{i},\t}(\s). \label{eq:finalprop5Ia0b} 
\end{align}
%
\subsubsection*{Step 1: Replace $\mathtt{J}(\t,\cdot;\mathbb{I}(\mathfrak{t}))$ by $\mathbf{J}(\t,\cdot;\mathbb{I}(\mathfrak{t}))$}
Recall {\small$\mathrm{Loc}^{(1)}\Upsilon^{\pm,\mathrm{j},\mathrm{i},\t}$} in Lemma \ref{lemma:finalprop2}. It is a function of $\tau\mapsto(\mathtt{J}(\tau,\cdot;\mathbb{I}(\mathfrak{t})),\mathbf{U}^{\tau,\cdot}[\mathbb{I}(\mathfrak{t})])$. Now, set the following version of $\mathbf{G}^{\tau}$ from the statement of {Lemma \ref{lemma:finalprop5}} without the cutoff:
\begin{align}
\mathbf{G}^{\t,\mathrm{uncut}}(\tau,\x;\mathbb{I}(\mathfrak{t})) \ := \ \exp[\lambda(\t)\mathbf{J}(\tau,\x;\mathbb{I}(\mathfrak{t}))+\lambda(\t)\mathrm{U}]. \label{eq:finalprop5Ia1}
\end{align}
Equivalently, \eqref{eq:finalprop5Ia1} is $\mathsf{G}^{\t}$ in Lemma \ref{lemma:finalprop2} but for $\mathbf{J}$ from Definition \ref{definition:le10} instead of $\mathtt{J}$. Restrict to $\mathcal{E}^{\mathrm{ap},\mathrm{j},\mathrm{i}}\cap\mathcal{E}^{\mathrm{dens},\mathrm{j},\mathrm{i}}$, which we defined in Lemma \ref{lemma:finalprop2}. Fix $\mathrm{k}\in\llbracket0,\mathfrak{m}(\mathrm{j},\mathrm{i})-1\rrbracket$. Given $\r\in[0,\tau(\mathrm{j},\mathrm{i})]$, for convenience let us set $\x(\s,\r,\mathrm{k}):=0(\s-\r)\pm2\mathrm{k}\mathfrak{l}(\mathrm{j})$. We claim $\x(\s,\r,\mathrm{k})\in\mathbb{I}$ from Lemma \ref{lemma:finalprop2}. (Indeed, it is a shift by $\leq2\mathfrak{m}(\mathrm{j},\mathrm{i})\mathfrak{l}(\mathrm{j})$ of something in the image of the characteristic map $\z\mapsto\z(\tau)$ of speed {\small$\lesssim\N^{3/2}$} for $\tau\leq\tau(\mathrm{j},\mathrm{i})$.) Thus, by Lemma \ref{lemma:le16} (and justification as for why Lemma \ref{lemma:le16} applies to be given after), we get the following on a very high probability event that we denote by $\mathcal{E}^{\mathrm{vhp},1}$ (see Lemma \ref{lemma:finalprop2} for $\mathsf{G}^{\t}$ below):
\begin{align}
\mathsf{G}^{\t}(\s-\r,\x(\s,\r,\mathrm{k})) \ = \ \mathbf{G}^{\t,\mathrm{uncut}}(\s-\r,\x(\s,\r,\mathrm{k});\mathbb{I}(\mathfrak{t})). \label{eq:finalprop5Ia2}
\end{align}
Indeed, the LHS and RHS of \eqref{eq:finalprop5Ia2} are the same exponential evaluated at $\mathtt{J}(\s-\r,\x;\mathbb{I}(\mathfrak{t}))$ and $\mathbf{J}(\s-\r,\x;\mathbb{I}(\mathfrak{t}))$, respectively, for $\x\in\mathbb{I}$. Note that Lemma \ref{lemma:le16} requires an a priori estimate for $\mathtt{J}$; said a priori bound holds because we have restricted {ourselves to the} event $\mathcal{E}^{\mathrm{ap},\mathrm{j},\mathrm{i}}$ on which we have two sided bounds for the exponential of $\mathtt{J}$. We also require the bound $|\mathbb{I}(\mathfrak{t})|\lesssim\N^{5/6}$. This can be checked by our upper bounds $|\mathbb{I}(\mathfrak{t})|\lesssim\N^{\gamma_{\mathrm{ap}}}[\N\tau(\mathrm{j},\mathrm{i})^{1/2}+\N^{3/2}\tau(\mathrm{j},\mathrm{i})+|\mathbb{I}|]$ and $|\mathbb{I}|\lesssim\mathfrak{m}(\mathrm{j},\mathrm{i})\mathfrak{l}(\mathrm{j})+\N^{3/2+\gamma_{\mathrm{ap}}}\tau(\mathrm{j},\mathrm{i})$ given prior to \eqref{eq:finalprop4IIa3a}, along with Definition \ref{definition:bg2131} and calculations that boil down to $\mathfrak{m}(\mathrm{j},\mathrm{i})\mathfrak{l}(\mathrm{j})\leq\N^{3/4+\gamma_{\mathrm{KL}}}$ and $\tau(\mathrm{j},\mathrm{i})\lesssim\N^{-3/4+\gamma_{\mathrm{KL}}}$. Last, to use Lemma \ref{lemma:le16}, we need $\sigma$ in \eqref{eq:finalprop5Ia} to satisfy $|\sigma|\lesssim\N^{\gamma_{\mathrm{reg}}}|\mathbb{I}(\mathfrak{t})|^{-1/2}$. This holds since we work on $\mathcal{E}^{\mathrm{dens},\mathrm{j},\mathrm{i}}$. 
\subsubsection*{Step 2: Gymnastics for $\mathbf{G}^{\t,\mathrm{uncut}}$}
We start this step with the following trivial decomposition of the RHS of \eqref{eq:finalprop5Ia2}:
\begin{align}
\mathbf{G}^{\t,\mathrm{uncut}}(\s-\r,\x(\s,\r,\mathrm{k});\mathbb{I}(\mathfrak{t})) &= \mathbf{G}^{\t,\mathrm{uncut}}(\s-\r,\inf\mathbb{I}(\mathfrak{t});\mathbb{I}(\mathfrak{t}))\nonumber\\
&\times \ \exp\{\lambda(\t)[\mathbf{J}(\t,\x(\s,\r,\mathrm{k});\mathbb{I}(\mathfrak{t}))-\mathbf{J}(\t,\inf\mathbb{I}(\mathfrak{t});\mathbb{I}(\mathfrak{t}))]\}. \label{eq:finalprop5Ia3}
\end{align}
On the intersection $\mathcal{E}^{\mathrm{vhp},1}\cap\mathcal{E}^{\mathrm{ap},\mathrm{j},\mathrm{i}}$, where $\mathcal{E}^{\mathrm{vhp},1}$ is when \eqref{eq:finalprop5Ia2} holds, we know $\N^{-2\gamma_{\mathrm{ap}}}\lesssim\mathrm{LHS}\eqref{eq:finalprop5Ia3}\lesssim\N^{2\gamma_{\mathrm{ap}}}$. Moreover, as in the proof of Lemma \ref{lemma:le16}, there exists a very high probability $\mathcal{E}^{\mathrm{vhp},2}$ such that on {$\mathcal{E}^{\mathrm{vhp}}\cap\mathcal{E}^{\mathrm{ap},\mathrm{j},\mathrm{i}}\cap\mathcal{E}^{\mathrm{dens},\mathrm{j},\mathrm{i}}$, where $\mathcal{E}^{\mathrm{vhp}}:=\mathcal{E}^{\mathrm{vhp},1}\cap\mathcal{E}^{\mathrm{vhp},2}$,} the $\exp$-factor in $\mathrm{RHS}\eqref{eq:finalprop5Ia3}$ is both $\gtrsim1$ and $\lesssim1$. (See \eqref{eq:le16I5} and the definition of $\mathcal{E}$ therein from right before \eqref{eq:le16I3}. Now, let us also clarify that this argument works because the setting of Lemma \ref{lemma:le16} holds in this proof as justified in the previous step.) So, on the intersection $\mathcal{E}^{\mathrm{vhp},1}\cap\mathcal{E}^{\mathrm{vhp},2}\cap\mathcal{E}^{\mathrm{ap},\mathrm{j},\mathrm{i}}\cap\mathcal{E}^{\mathrm{dens},\mathrm{j},\mathrm{i}}$, the first factor in $\mathrm{RHS}\eqref{eq:finalprop5Ia3}$ is unchanged if $\mathbf{G}^{\t,\mathrm{uncut}}\mapsto\mathbf{G}^{\t}$. Thus, upon recalling {\small$\mathscr{A}^{\mathrm{loc},(1),\pm,\t}$} from Lemma \ref{lemma:finalprop2} and $\mathscr{A}^{\mathrm{loc},\pm,\t}$ from {Lemma \ref{lemma:finalprop5}}, on $\mathcal{E}^{\mathrm{vhp},1}\cap\mathcal{E}^{\mathrm{vhp},2}\cap\mathcal{E}^{\mathrm{ap},\mathrm{j},\mathrm{i}}\cap\mathcal{E}^{\mathrm{dens},\mathrm{j},\mathrm{i}}$, we know that
\begin{align}
\mathscr{A}^{\mathrm{loc},(1),\pm,\t} \ = \ \mathscr{A}^{\mathrm{loc},\pm,\t}, \label{eq:finalprop5Ia4a}
\end{align}
if $\mathfrak{A}^{(\mathrm{k}),\pm}$ in the formula for $\mathscr{A}^{\mathrm{loc},\pm,\t}$ equals $\mathfrak{A}^{(\mathrm{k}),(1),\pm}$ in the formula for $\mathscr{A}^{\mathrm{loc},(1),\pm,\t}$ multiplied by the $\exp$-factor in \eqref{eq:finalprop5Ia3}:
\begin{align}
\mathfrak{A}^{(\mathrm{k}),\pm}_{\tau}(\mathbf{U}) \ = \ \mathfrak{A}^{(\mathrm{k}),(1),\pm}_{\tau}(\mathbf{U})\mathrm{CutExp}\{\lambda(\t)\N^{-\frac12}{\textstyle\sum_{\w=\inf\mathbb{I}(\mathfrak{t})+1}^{0(\tau)\pm2\mathrm{k}\mathfrak{l}(\mathrm{j})}}\mathbf{U}(\w)\} \quad\mathrm{for}\quad\mathbf{U}\in\R^{\mathbb{I}(\mathfrak{t})}. \label{eq:finalprop5Ia4b}
\end{align}
Above, $\mathrm{CutExp}(\mathrm{a})=\mathrm{Exp}(\mathrm{a})\mathbf{1}[\mathrm{Exp}(\mathrm{a})\lesssim1]$. Indeed, as argued after \eqref{eq:finalprop5Ia3}, we know the $\exp$-factor in $\mathrm{RHS}\eqref{eq:finalprop5Ia3}$ is $\lesssim1$ on the intersection $\mathcal{E}^{\mathrm{vhp},1}\cap\mathcal{E}^{\mathrm{vhp},2}\cap\mathcal{E}^{\mathrm{ap},\mathrm{j},\mathrm{i}}\cap\mathcal{E}^{\mathrm{dens},\mathrm{j},\mathrm{i}}$. Thus, upgrading to $\mathrm{CutExp}$ from $\mathrm{Exp}$ comes for free.
\subsubsection*{Step 3: Finishing the proof of \eqref{eq:finalprop5Ia} assuming that $\mathfrak{A}^{(\mathrm{k}),\pm}$ satisfies the necessary properties}  
{We first recall the notation} $\mathcal{E}^{\mathrm{vhp}}:=\mathcal{E}^{\mathrm{vhp},1}\cap\mathcal{E}^{\mathrm{vhp},2}$. By {the} union bound (for complements of $\mathcal{E}^{\mathrm{vhp},1},\mathcal{E}^{\mathrm{vhp},2}$), we know $\mathcal{E}^{\mathrm{vhp}}$ is very high probability. Since \eqref{eq:finalprop5Ia4a} holds on $\mathcal{E}^{\mathrm{vhp}}\cap\mathcal{E}^{\mathrm{ap},\mathrm{j},\mathrm{i}}\cap\mathcal{E}^{\mathrm{dens},\mathrm{j},\mathrm{i}}$, we know {\small$\mathrm{Loc}^{(1)}\Upsilon^{\pm,\mathrm{j},\mathrm{i},\t}$} equals $\mathrm{Loc}\Upsilon^{\pm,\mathrm{j},\mathrm{i},\t}$ on this intersection. Indeed, these two terms are given by the same function, but the former is evaluated at $\mathrm{LHS}\eqref{eq:finalprop5Ia4a}$, and the latter is evaluated at $\mathrm{RHS}\eqref{eq:finalprop5Ia4b}$; see \eqref{eq:finalprop2IIa} and \eqref{eq:finalprop5IIa}. Thus, we have {the following for any large but fixed $\mathrm{D}>0$}:
\begin{align}
&\E^{\sigma,\mathrm{loc},\s}\mathrm{Loc}^{(1)}\Upsilon^{\pm,\mathrm{j},\mathrm{i},\t}(\s) \nonumber\\
&= \ \E^{\sigma,\mathrm{loc},\s}\mathbf{1}[\mathcal{E}^{\mathrm{vhp}}]\mathrm{Loc}^{(1)}\Upsilon^{\pm,\mathrm{j},\mathrm{i},\t}(\s)+\E^{\sigma,\mathrm{loc},\s}\{1-\mathbf{1}[\mathcal{E}^{\mathrm{vhp}}]\}\mathrm{Loc}^{(1)}\Upsilon^{\pm,\mathrm{j},\mathrm{i},\t}(\s) \label{eq:finalprop5IIa5a}\\
&\leq \ \E^{\sigma,\mathrm{loc},\s}\mathrm{Loc}\Upsilon^{\pm,\mathrm{j},\mathrm{i},\t}(\s)+\E^{\sigma,\mathrm{loc},\s}\{1-\mathbf{1}[\mathcal{E}^{\mathrm{vhp}}]\}\mathrm{Loc}^{(1)}\Upsilon^{\pm,\mathrm{j},\mathrm{i},\t}(\s) \label{eq:finalprop5IIa5b}\\
&\lesssim \ \E^{\sigma,\mathrm{loc},\s}\mathrm{Loc}\Upsilon^{\pm,\mathrm{j},\mathrm{i},\t}(\s)+\N^{-{\mathrm{D}}}. \label{eq:finalprop5IIa5c}
\end{align}
(We clarify \eqref{eq:finalprop5IIa5c} follows since $\mathcal{E}^{\mathrm{vhp}}$ is very high probability, and {\small$\mathrm{Loc}^{(1)}\Upsilon^{\pm,\mathrm{j},\mathrm{i},\t}$} is polynomially bounded in $\N$ by construction in \eqref{eq:finalprop2IIa}.) Combining \eqref{eq:finalprop5IIa5a}-\eqref{eq:finalprop5IIa5c}, which is uniform in $\sigma,\t,\s$, with \eqref{eq:finalprop5Ia0}, gives \eqref{eq:finalprop5Ia}.
\subsubsection*{Step 4: Using the same argument to derive \eqref{eq:finalprop5Ib}}
Again, because \eqref{eq:finalprop5Ia4a} is true on $\mathcal{E}^{\mathrm{vhp}}\cap\mathcal{E}^{\mathrm{ap},\mathrm{j},\mathrm{i}}\cap\mathcal{E}^{\mathrm{dens},\mathrm{j},\mathrm{i}}$, by construction in \eqref{eq:finalprop2IIb} and \eqref{eq:finalprop5IIb}, we know {\small$\mathrm{Loc}^{(1)}\Lambda^{\pm,\mathrm{j},\mathrm{i},\t}$} equals $\mathrm{Loc}\Lambda^{\pm,\mathrm{j},\mathrm{i},\t}$ on this intersection for the same reason this was true with $\Upsilon$ in place of $\Lambda$. (See right before \eqref{eq:finalprop5IIa5a}.) We also know that {\small$\mathrm{Loc}^{(1)}\Lambda^{\pm,\mathrm{j},\mathrm{i},\t}$} is polynomially bounded in $\N$; see \eqref{eq:finalprop2IIb}. Therefore, the reasoning for \eqref{eq:finalprop5IIa5a}-\eqref{eq:finalprop5IIa5c} also provides us the following estimate {for any large but fixed $\mathrm{D}>0$}, which has {an} implied constant that is independent of $\sigma,\t,\s$:
\begin{align}
\E^{\sigma,\mathrm{loc},\s}\mathrm{Loc}^{(1)}\Lambda^{\pm,\mathrm{j},\mathrm{i},\t}(\s) \ \lesssim \ \E^{\sigma,\mathrm{loc},\s}\mathrm{Loc}\Lambda^{\pm,\mathrm{j},\mathrm{i},\t}(\s)+\N^{-{\mathrm{D}}}. 
\end{align}
Combining this with \eqref{eq:finalprop5Ia0b} gives us \eqref{eq:finalprop5Ib}, thereby completing this step.
\subsubsection*{Step 5: Showing that $\mathfrak{A}^{(\mathrm{k}),\pm}$ satisfies the necessary properties}
We are left to explain the paragraph after the definition of $\mathscr{A}^{\mathrm{loc},\pm,\t}$ in the statement of {Lemma \ref{lemma:finalprop5}}. Look at $\mathrm{RHS}\eqref{eq:finalprop5Ia4b}$ for {a} fixed $\tau$. Recall that the discrete intervals $\mathbb{J}(\tau,\mathrm{k})$ that we built in Lemma \ref{lemma:finalprop2} are mutually disjoint. Also recall $\mathbb{J}(\tau,\mathrm{k})$ is either disjoint from or contained in the set of $\w$ appearing in the sum in $\mathrm{CutExp}$ in \eqref{eq:finalprop5Ia4b}; this comes from Lemma \ref{lemma:finalprop2}. Now, take $\mathrm{k}\in\llbracket0,\mathfrak{m}(\mathrm{j},\mathrm{i})-1\rrbracket$. {Condition on $\mathbf{U}(\w)$ for all $\w\not\in\mathbb{J}(\tau,\mathrm{k})$, condition on the average of $\mathbf{U}(\w)$ over $\w\in\mathbb{J}(\tau,\mathrm{k})$, and denote the average by $\rho$.} We claim that after this conditioning, the expectation $\E^{\sigma,\tau,\mathbb{I}(\mathfrak{t})}$ becomes $\E^{\rho,\tau,\mathbb{J}(\mathfrak{t},\mathrm{k})}$. Indeed, this is the statement that if one takes a random walk bridge, looks at an increment of a given length, and conditions on its average drift, one gets a random walk bridge for the law of said increment. We additionally claim that after the conditioning, the $\mathrm{CutExp}$-factor in \eqref{eq:finalprop5Ia4b} is constant. To see this, suppose first that the sum in said $\mathrm{CutExp}$-factor is over $\mathbf{U}(\w)$ for some subset of $\w\not\in\mathbb{J}(\tau,\mathrm{k})$. We have conditioned on all such $\mathbf{U}(\w)$, so the $\mathrm{CutExp}$-factor is indeed constant. Assume that the sum in the $\mathrm{CutExp}$-factor in \eqref{eq:finalprop5Ia4b} has $\mathbf{U}(\w)$ for all $\w\in\mathbb{J}(\tau,\mathrm{k})$. The dependence on such $\mathbf{U}(\w)$ is through their sum (or equivalently, their average). Therefore, the $\mathrm{CutExp}$-factor is still constant in this case. Ultimately, we know, in the language of Lemma \ref{lemma:le2}, that {\small$\E^{\sigma,\tau,\mathbb{I}(\mathfrak{t})}(\mathfrak{A}^{(\mathrm{k}),\pm}_{\tau}|\mathscr{F}(\mathrm{k}))$} is proportional to {\small$\E^{\rho,\tau,\mathbb{J}(\tau,\mathrm{k})}\mathfrak{A}^{(\mathrm{k}),(1),\pm}_{\tau}$}, which vanishes by Lemma \ref{lemma:finalprop2}. This proves the first property that we claimed after the formula for $\mathscr{A}^{\mathrm{loc},\pm,\t}$ (namely, satisfaction of the condition in Lemma \ref{lemma:le2}). (Technically, we have shown more, because we have shown vanishing under conditional expectation after more conditioning than is required in Lemma \ref{lemma:le2}. We never use this stronger statement, however.) It is left to get the deterministic bound $|\mathfrak{A}^{(\mathrm{k}),\pm}|\lesssim\N^{30\gamma_{\mathrm{reg}}}\mathfrak{l}(\mathrm{j})^{-3/2}$. This follows by the deterministic bound for {\small$\mathfrak{A}^{(\mathrm{k}),(1),\pm}$} in Lemma \ref{lemma:finalprop2} and {the fact that} $\mathrm{CutExp}$ in \eqref{eq:finalprop5Ia4b} is $\mathrm{O}(1)$ by construction.
\end{proof}
\subsubsection{Applying the Kipnis-Varadhan estimate ({Proposition \ref{prop:kv1}})}
Observe that the RHS of \eqref{eq:finalprop5Ia} and \eqref{eq:finalprop5Ib}, respectively, can be directly treated by Proposition \ref{prop:kv1}. Doing so and computing the resulting bounds is the purpose of the following.
\begin{lemma}\label{lemma:finalprop7}
 Fix $1\leq\mathrm{j}\leq\mathrm{j}(\infty)$ and $1\leq\mathrm{i}<\mathrm{i}(\mathrm{j})$. Retain {the} notation of {Lemmas \ref{lemma:finalprop2}, \ref{lemma:finalprop4}, \ref{lemma:finalprop5}}. Uniformly in $\s,\t,\sigma$, we have 
\begin{align}
\N\cdot\mathbf{1}[|\sigma|\lesssim1]\{\E^{\sigma,\mathrm{loc},\s}\mathrm{Loc}\Upsilon^{\pm,\mathrm{j},\mathrm{i},\t}(\s)+\E^{\sigma,\mathrm{loc},\s}\mathrm{Loc}\Lambda^{\pm,\mathrm{j},\mathrm{i},\t}(\s)\} \ &\lesssim \ \N^{-30\beta_{\mathrm{BG}}}. \label{eq:finalprop7I}
\end{align}
\end{lemma}
\begin{proof}
We bound the $\Upsilon$-term. To bound the $\Lambda$-term in \eqref{eq:finalprop7I}, the same argument (and calculation) works; the only difference is that we must use \eqref{eq:kv1III} instead of \eqref{eq:kv1II}. To this end, we apply \eqref{eq:kv1II} with the following choices. First take $\mathscr{B}=\N^{-\beta(\mathrm{j},\mathrm{i}-1)+\gamma_{\mathrm{ap}}}$. Take $\varphi(\cdot)=\exp[\lambda(\t)\cdot]\mathbf{1}\{\N^{-2\gamma_{\mathrm{ap}}}\lesssim\exp[\lambda(\t)\cdot]\lesssim\N^{2\gamma_{\mathrm{ap}}}\}$. (Here, we emphasize that $\t$ is not a time-variable which we integrate or average; it is a fixed parameter for the coupling constant $\lambda(\t)$.) We choose {\small$\mathfrak{a}(\tau,\cdot;\mathrm{k}):=\mathfrak{A}^{(\mathrm{k}),\pm}_{\tau}(\cdot)$}, which satisfies {the} constraints of Lemma \ref{lemma:le2} (and thus of Proposition \ref{prop:kv1}) with sets $\mathbb{J}(\tau,\mathrm{k})$ given by $\tau$-dependent shifts of $\mathbb{I}(\mathrm{k})$, which themselves have length $|\mathbb{I}(\mathrm{k})|\lesssim\mathfrak{l}(\mathrm{j})$. We take $\mathrm{m}$ (the number of $\mathfrak{a}(\cdot,\cdot;\mathrm{k})$ that we average) to be $\mathfrak{m}(\mathrm{j},\mathrm{i})$. Next, we take $\mathfrak{t}=\tau(\mathrm{j},\mathrm{i})\lesssim1$. (This upper bound follows by construction in Definition \ref{definition:bg2131}.) We emphasize the bound $\mathfrak{t}|\mathbb{I}(\mathfrak{t})|\lesssim\N^{\gamma_{\mathrm{KV}}}$ can be directly verified via Definition \ref{definition:bg2131}. This lets us use \eqref{eq:kv1II} to deduce the following estimate (with explanation given after), where $\|\varphi\|$ is defined in Proposition \ref{prop:kv1}:
\begin{align}
\mathbf{1}[|\sigma|\lesssim1]\E^{\sigma,\mathrm{loc},\s}\mathrm{Loc}\Upsilon^{\pm,\mathrm{j},\mathrm{i},\t}(\s) \ \lesssim \ &\N^{\beta(\mathrm{j},\mathrm{i}-1)+20\gamma_{\mathrm{reg}}}\N^{\gamma_{\mathrm{KV}}}\tau(\mathrm{j},\mathrm{i})|\mathbb{I}(\tau(\mathrm{j},\mathrm{i}))|^{\frac12}\N^{-2\beta(\mathrm{j},\mathrm{i}-1)+2\gamma_{\mathrm{ap}}} \label{eq:finalprop7I1a}\\
+ \ &\N^{\beta(\mathrm{j},\mathrm{i}-1)+100\gamma_{\mathrm{reg}}}\N^{-2}\mathfrak{m}(\mathrm{j},\mathrm{i})^{-1}\tau(\mathrm{j},\mathrm{i})^{-1}\|\varphi\|^{2}\mathfrak{l}(\mathrm{j})^{2}\mathfrak{l}(\mathrm{j})^{-3}. \label{eq:finalprop7I1b}
\end{align}
Indeed, we just used \eqref{eq:kv1II} with the above choices and then the a priori bounds for $\mathfrak{a}(\cdot,\cdot;\mathrm{k})$ from Lemma \ref{lemma:finalprop5} with $|\mathbb{I}(\mathrm{k})|\lesssim\mathfrak{l}(\mathrm{j})$ (which we noted in the previous paragraph). In particular, this is where the second line comes from. We now estimate the upper bound in \eqref{eq:finalprop7I1a}-\eqref{eq:finalprop7I1b}. First, recall from Definition \ref{definition:le10} that $|\mathbb{I}(\tau(\mathrm{j},\mathrm{i}))|\lesssim\N^{\gamma_{\mathrm{ap}}}[\N\tau(\mathrm{j},\mathrm{i})^{1/2}+\N^{3/2}\tau(\mathrm{j},\mathrm{i})+|\mathbb{I}|]$. Also recall that $|\mathbb{I}|\lesssim\mathfrak{m}(\mathrm{j},\mathrm{i})\mathfrak{l}(\mathrm{j})+\N^{3/2+\gamma_{\mathrm{ap}}}\tau(\mathrm{j},\mathrm{i})$ from Lemma \ref{lemma:finalprop2}. So, $|\mathbb{I}(\tau(\mathrm{j},\mathrm{i}))|\lesssim\N^{1+\gamma_{\mathrm{ap}}}\tau(\mathrm{j},\mathrm{i})^{1/2}+\N^{3/2+2\gamma_{\mathrm{ap}}}\tau(\mathrm{j},\mathrm{i})+\N^{\gamma_{\mathrm{ap}}}\mathfrak{m}(\mathrm{j},\mathrm{i})\mathfrak{l}(\mathrm{j})$. Thus, everything in $\mathrm{RHS}\eqref{eq:finalprop7I1a}$ can be expressed in terms of $\N$ and constants from Definition \ref{definition:bg2131}. We now claim this means
{
\begin{align}
&\mathrm{RHS}\eqref{eq:finalprop7I1a} \nonumber\\
&\lesssim \ \N^{-\beta(\mathrm{j},\mathrm{i}-1)+20\gamma_{\mathrm{reg}}+2\gamma_{\mathrm{ap}}+\gamma_{\mathrm{KV}}}\label{eq:finalprop7I2a}\\
&\times[\N^{\frac12+\frac12\gamma_{\mathrm{ap}}}\tau(\mathrm{j},\mathrm{i})^{\frac54}+\N^{\frac34+\gamma_{\mathrm{ap}}}\tau(\mathrm{j},\mathrm{i})^{\frac32}+\N^{\frac12\gamma_{\mathrm{ap}}}\tau(\mathrm{j},\mathrm{i})\mathfrak{m}(\mathrm{j},\mathrm{i})^{\frac12}\mathfrak{l}(\mathrm{j})^{\frac12}] \nonumber\\
&\lesssim \ \N^{3\gamma_{\mathrm{KV}}}\N^{-\beta(\mathrm{j},\mathrm{i}-1)}\tau(\mathrm{j},\mathrm{i})\mathfrak{m}(\mathrm{j},\mathrm{i})\mathfrak{l}(\mathrm{j})\label{eq:finalprop7I2b}\\
&\times[\N^{\frac12}\tau(\mathrm{j},\mathrm{i})^{\frac14}\mathfrak{m}(\mathrm{j},\mathrm{i})^{-1}\mathfrak{l}(\mathrm{j})^{-1}+\N^{\frac34}\tau(\mathrm{j},\mathrm{i})^{\frac12}\mathfrak{m}(\mathrm{j},\mathrm{i})^{-1}\mathfrak{l}(\mathrm{j})^{-1}+\mathfrak{m}(\mathrm{j},\mathrm{i})^{-\frac12}\mathfrak{l}(\mathrm{j})^{-\frac12}]\nonumber\\
&\lesssim \ \N^{-1+20\gamma_{\mathrm{reg}}+90\beta_{\mathrm{BG}}+4\gamma_{\mathrm{KV}}}\label{eq:finalprop7I2c}\\
&\times[\N^{\frac12}\tau(\mathrm{j},\mathrm{i})^{\frac14}\mathfrak{m}(\mathrm{j},\mathrm{i})^{-1}\mathfrak{l}(\mathrm{j})^{-1}+\N^{\frac34}\tau(\mathrm{j},\mathrm{i})^{\frac12}\mathfrak{m}(\mathrm{j},\mathrm{i})^{-1}\mathfrak{l}(\mathrm{j})^{-1}+\mathfrak{m}(\mathrm{j},\mathrm{i})^{-\frac12}\mathfrak{l}(\mathrm{j})^{-\frac12}]\nonumber\\
&\lesssim \ \N^{-1+5\gamma_{\mathrm{KV}}}[\N^{\frac12}\tau(\mathrm{j},\mathrm{i})^{\frac14}\mathfrak{m}(\mathrm{j},\mathrm{i})^{-1}\mathfrak{l}(\mathrm{j})^{-1}+\N^{\frac34}\tau(\mathrm{j},\mathrm{i})^{\frac12}\mathfrak{m}(\mathrm{j},\mathrm{i})^{-1}\mathfrak{l}(\mathrm{j})^{-1}+\mathfrak{m}(\mathrm{j},\mathrm{i})^{-\frac12}\mathfrak{l}(\mathrm{j})^{-\frac12}].\label{eq:finalprop7I2d}
\end{align}
}\eqref{eq:finalprop7I2a} follows from plugging our estimate for $|\mathbb{I}(\tau(\mathrm{j},\mathrm{i}))|$ from the previous paragraph into $\mathrm{RHS}\eqref{eq:finalprop7I1a}$. To get \eqref{eq:finalprop7I2b}, we first pull out all $\gamma_{\mathrm{ap}}$-exponents inside the square brackets in $\mathrm{RHS}\eqref{eq:finalprop7I2a}$. For the sake of an upper bound, this means we can change the exponent for $\N$ outside the square bracket in $\mathrm{RHS}\eqref{eq:finalprop7I2a}$ to $-\beta(\mathrm{j},\mathrm{i}-1)+20\gamma_{\mathrm{reg}}+3\gamma_{\mathrm{ap}}+\gamma_{\mathrm{KV}}$. Now, by Definitions \ref{definition:reg}, \ref{definition:method8}, and Proposition \ref{prop:kv1}, we know that $\gamma_{\mathrm{KV}}$ is large  compared to $\gamma_{\mathrm{reg}},\gamma_{\mathrm{ap}}$, which means said exponent is $\leq-\beta(\mathrm{j},\mathrm{i}-1)+3\gamma_{\mathrm{KV}}$. This explains the first factor in \eqref{eq:finalprop7I2b}. The rest of \eqref{eq:finalprop7I2b} follows by factoring out $\tau(\mathrm{j},\mathrm{i})\mathfrak{m}(\mathrm{j},\mathrm{i})\mathfrak{l}(\mathrm{j})$ from each term inside the square brackets in $\mathrm{RHS}\eqref{eq:finalprop7I2a}$. \eqref{eq:finalprop7I2c} follows by construction of $\beta(\mathrm{j},\mathrm{i}-1)$ in Definition \ref{definition:bg2131}. (Technically, in \eqref{eq:finalprop7I2b}, we paired $\beta(\mathrm{j},\mathrm{i}-1)$ with $\tau(\mathrm{j},\mathrm{i})\mathfrak{m}(\mathrm{j},\mathrm{i})\mathfrak{l}(\mathrm{j})$, not $\tau(\mathrm{j},\mathrm{i}-1)\mathfrak{m}(\mathrm{j},\mathrm{i}-1)\mathfrak{l}(\mathrm{j})$. But, this only introduces the multiplicative cost of $\lesssim\N^{\gamma_{\mathrm{KV}}}$.) \eqref{eq:finalprop7I2d} follows by exponent-counting as in the derivation of \eqref{eq:finalprop7I2b}. We now bound \eqref{eq:finalprop7I2d}. In the following reasoning, we invite the reader to refer to Definition \ref{definition:bg2131} for details. First, note that if $\mathfrak{m}(\mathrm{j},\mathrm{i})\mathfrak{l}(\mathrm{j})\geq\N^{1/2}$, then $\tau(\mathrm{j},\mathrm{i})\lesssim\N^{-2/3}$. So, in this case, the term in square brackets in \eqref{eq:finalprop7I2d} is $\lesssim\N^{-1/99}$. Since $\gamma_{\mathrm{KV}}$ is small (see Proposition \ref{prop:kv1}), $\mathfrak{m}(\mathrm{j},\mathrm{i})\geq\N^{1/2}$ implies $\eqref{eq:finalprop7I2d}\lesssim\N^{-1-1/99}$. This is certainly $\lesssim\N^{-1-30\beta_{\mathrm{BG}}}$ (see Definition \ref{definition:method8}). Thus, it suffices to assume $\mathfrak{m}(\mathrm{j},\mathrm{i})\mathfrak{l}(\mathrm{j})\leq\N^{1/2}$. In this case, we claim 
\begin{align}
\eqref{eq:finalprop7I2d} \ &\lesssim \ \N^{-1+5\gamma_{\mathrm{KV}}}[\N^{\frac12}\N^{-\frac12}\mathfrak{m}(\mathrm{j},\mathrm{i})^{-\frac12}\mathfrak{l}(\mathrm{j})^{-\frac12}+\N^{\frac34}\N^{-1}+\mathfrak{m}(\mathrm{j},\mathrm{i})^{-\frac12}\mathfrak{l}(\mathrm{j})^{-\frac12}] \nonumber\\
&\lesssim \ \N^{-1+5\gamma_{\mathrm{KV}}-\frac{1}{20}} \nonumber\\
&\lesssim \ \N^{-1-30\beta_{\mathrm{BG}}}. \label{eq:finalprop7I3}
\end{align}
Indeed, we first use $\tau(\mathrm{j},\mathrm{i})=\N^{-2}\mathfrak{m}(\mathrm{j},\mathrm{i})\mathfrak{l}(\mathrm{j})$. Then, we use $\mathfrak{m}(\mathrm{j},\mathrm{i})\mathfrak{l}(\mathrm{j})\gtrsim\N^{1/9}$. Again, for both of these, see Definition \ref{definition:bg2131}. The last bound follows because $\gamma_{\mathrm{KV}},\beta_{\mathrm{BG}}$ are small; see Definition \ref{definition:method8} and Proposition \ref{prop:kv1}. We now control \eqref{eq:finalprop7I1b}. To this end, using the definition of $\beta(\mathrm{j},\mathrm{i}-1)$ in Definition \ref{definition:bg2131} gives the following (like in the derivation of \eqref{eq:finalprop7I2c}):
\begin{align}
\eqref{eq:finalprop7I1b} \ &\lesssim \ \N^{-1-90\beta_{\mathrm{BG}}+120\gamma_{\mathrm{reg}}}\|\varphi\|^{2}. \label{eq:finalprop7I4}
\end{align}
Recall $\varphi$ from the first paragraph of this proof. It is direct to check by calculus that $\|\varphi\|^{2}\lesssim\N^{10\gamma_{\mathrm{ap}}}$. (In words, $\varphi$ has support in an interval of length $\lesssim\log\N$ since it is defined by the two-sided cutoff of an exponential. It is also bounded uniformly above by $\N^{2\gamma_{\mathrm{ap}}}$.) Thus, $\mathrm{RHS}\eqref{eq:finalprop7I4}\lesssim\N^{-1-90\beta_{\mathrm{BG}}+120\gamma_{\mathrm{reg}}+10\gamma_{\mathrm{ap}}}\lesssim\N^{-1-30\beta_{\mathrm{BG}}}$, where this final bound follows from Definition \ref{definition:method8}. We deduce $\N\cdot\eqref{eq:finalprop7I1b}\lesssim\mathrm{RHS}\eqref{eq:finalprop7I}$. Using this with \eqref{eq:finalprop7I1a}-\eqref{eq:finalprop7I1b}, \eqref{eq:finalprop7I2a}-\eqref{eq:finalprop7I2d}, and \eqref{eq:finalprop7I3} completes the proof.
\end{proof}
\subsection{Bound for \eqref{eq:bg2138Ib}}
We ultimately use Lemmas \ref{lemma:bg213101}, \ref{lemma:finalprop2}, \ref{lemma:finalprop4}, \ref{lemma:finalprop5}, \ref{lemma:finalprop7} to bound \eqref{eq:bg2138Ic}-\eqref{eq:bg2138Id}. To estimate \eqref{eq:bg2138Ib}, we roughly follow the same strategy. (There are only some cosmetic adjustments to make; we explain them below.) First:
\begin{lemma}\label{lemma:finalprop8}
 Adopt {the} notation of {Lemma \ref{lemma:bg213101}} and specialize to $\mathrm{i}=1$. We have the following estimate, with notation explained afterwards, with probability 1:
\begin{align}
\eqref{eq:bg2138Ib} \ \lesssim \ {\textstyle\sup_{\t}}\N^{\beta_{\mathrm{BG}}}{\textstyle\int_{\tau(\mathrm{j},1)}^{1}}|\mathbb{T}(\N)|^{-1}{\textstyle\sum_{\y}}\E[\mathbf{1}(\s\leq\t_{\mathrm{st}})\N\cdot\mathrm{Cent}\Psi^{\pm,\mathrm{j},\t}(\s,\y(\s))]\d\s. \label{eq:finalprop8I}
\end{align}
Here, $\mathrm{Cent}\Psi^{\pm,\mathrm{j},\t}(\s,\y(\s))$ is defined to be the following, in which we recall $\mathscr{A}^{\pm,\t}$ from {Lemma \ref{lemma:bg213101}} for $\mathrm{i}=1$:
\begin{align}
\mathrm{Cent}\Psi^{\pm,\mathrm{j},\t}(\s,\y(\s)) \ := \ |\mathscr{A}^{\pm,\t}|\mathbf{1}[|\mathscr{A}^{\pm,\t}|>\N^{-\beta(\mathrm{j},1)}]\mathbf{1}[|\mathscr{A}^{\pm,\t}|\lesssim\N^{30\gamma_{\mathrm{reg}}}\mathfrak{l}(\mathrm{j})^{-\frac32}]. \label{eq:finalprop8II}
\end{align}
\end{lemma}
\begin{proof}
Refer to {the} notation in \eqref{eq:bg2138Ib} and \eqref{eq:bg2138IIb}. By the calculation in the proof of Lemma \ref{lemma:bg213101}, we first deduce
\begin{align}
\eqref{eq:bg2138Ib} \ &\leq \ {\textstyle\sup_{\t}}\N^{\beta_{\mathrm{BG}}}{\textstyle\int_{\tau(\mathrm{j},1)}^{1}}|\mathbb{T}(\N)|^{-1}{\textstyle\sum_{\y}}\N\cdot\E[\mathbf{1}(\s\leq\t_{\mathrm{st}})|\mathscr{A}^{\pm,\t}|\mathbf{1}(|\mathscr{A}^{\pm,\t}|>\N^{-\beta(\mathrm{j},1)})]\d\s. \label{eq:finalprop8I1}
\end{align}
(Indeed, the calculation in the proof of Lemma \ref{lemma:bg213101} is just about changing variables; the underlying functional of $\mathscr{A}^{\pm,\t}$ plays no role in its validity, except that it is non-negative to make triangle inequalities work.) We are left to show $|\mathscr{A}^{\pm,\t}|\lesssim\N^{30\gamma_{\mathrm{reg}}}\mathfrak{l}(\mathrm{j})^{-3/2}$ with probability 1. This lets us put an indicator that turns the expectation in $\mathrm{RHS}\eqref{eq:finalprop8I1}$ into expectation of $\mathrm{Cent}\Psi^{\pm,\mathrm{j},\t}(\s,\y(\s))$. See the definition of $\mathscr{A}^{\pm,\t}$ from the statement of Lemma \ref{lemma:bg213101}. By construction in Definition \ref{definition:bg26}, $\mathds{R}^{\chi,\mathfrak{q},\pm,\mathrm{j}}$ is $\mathrm{O}(\N^{20\gamma_{\mathrm{reg}}}\mathfrak{l}(\mathrm{j})^{-3/2})$. Also, because $\s\leq\t_{\mathrm{st}}$, by Definition \ref{definition:method8}, we know $\mathbf{G}^{\t}$ is exponential of something $\lesssim\log\log\N$. Thus, $|\mathsf{G}^{\t}|\lesssim\N^{\gamma_{\mathrm{ap}}}$. This implies that $|\mathscr{A}^{\pm,\t}|\lesssim\N^{20\gamma_{\mathrm{reg}}+\gamma_{\mathrm{ap}}}\mathfrak{l}(\mathrm{j})^{-3/2}$. By Definitions \ref{definition:reg} and \ref{definition:method8}, we have $\gamma_{\mathrm{ap}}\leq10\gamma_{\mathrm{reg}}$, so we are done.
\end{proof}
The next ingredient is a localization of $\mathrm{Cent}\Psi^{\pm,\mathrm{j},\t}$. This is an analog of Lemma \ref{lemma:finalprop2}, and its proof is basically identical.
\begin{lemma}\label{lemma:finalprop9}
 Adopt {the} notation of {Lemma \ref{lemma:finalprop2}} and specialize to $\mathrm{i}=1$. With more notation explained after, we have {the following for any fixed $\mathrm{D}>0$}:
\begin{align}
&\E[\mathbf{1}(\s\leq\t_{\mathrm{st}})\mathrm{Cent}\Psi^{\pm,\mathrm{j},\t}(\s,\y(\s))] \nonumber\\
&\lesssim \ \N^{2\gamma_{\mathrm{ap}}}\E[\{\E^{\mathrm{loc},\s}[\mathrm{Loc}^{(1)}\Psi^{\pm,\mathrm{j},\t}(\s)]\}(\Pi^{\mathrm{j},1}\mathbf{U}^{\s-\tau(\mathrm{j},1),\y(\s-\tau(\mathrm{j},1))+\cdot})]+\N^{-{\mathrm{D}}}. \label{eq:finalprop9I}
\end{align}
Recall $\mathscr{A}^{\mathrm{loc},(1),\pm,\t}$ for $\mathrm{i}=1$ and events $\mathcal{E}^{\mathrm{ap},\mathrm{j},1}(\s),\mathcal{E}^{\mathrm{dens},\mathrm{j},1}(\s)$ from {Lemma \ref{lemma:finalprop2}}. We define $\mathrm{Loc}^{(1)}\Psi^{\pm,\mathrm{j},\t}(\s)$ to be the following product of four indicator functions:
\begin{align}
&\mathbf{1}[\mathcal{E}^{\mathrm{ap},\mathrm{j},1}(\s)]\mathbf{1}[\mathcal{E}^{\mathrm{dens},\mathrm{j},1}(\s)]|\mathscr{A}^{\mathrm{loc},(1),\pm,\t}|\nonumber\\
&\times\mathbf{1}[|\mathscr{A}^{\mathrm{loc},(1),\pm,\t}|\gtrsim\N^{-\beta(\mathrm{j},1)-\gamma_{\mathrm{ap}}}]\nonumber\\
&\times\mathbf{1}[|\mathscr{A}^{\mathrm{loc},(1),\pm,\t}|\lesssim\N^{30\gamma_{\mathrm{reg}}+\gamma_{\mathrm{ap}}}\mathfrak{l}(\mathrm{j})^{-\frac32}]. \label{eq:finalprop9II}
\end{align}
\end{lemma}
\begin{proof}
See the proof of Lemma \ref{lemma:finalprop2}. The first four steps have nothing to do with the function of {\small$\mathscr{A}^{\mathrm{loc},(1),\pm,\t}$} that we are interested in for the current lemma. In Step 5, we just show $\mathscr{A}^{\pm,\t}$ is equal to {\small$\mathrm{O}(\N^{\gamma_{\mathrm{ap}}}\mathscr{A}^{\mathrm{loc},(1),\pm,\t})$} plus $\mathrm{O}(\N^{-{\mathrm{D}}})$. This, combined with the factorization of $\E$ in Step 6 of the proof of Lemma \ref{lemma:finalprop2}, gives \eqref{eq:finalprop9I} up to an additional error given by the expectation of $|\mathrm{Cent}\Psi^{\pm,\mathrm{j},\t}|$ on a very low probability event. But by construction in Lemma \ref{lemma:finalprop8}, we know $|\mathrm{Cent}\Psi^{\pm,\mathrm{j},\t}|\lesssim\N$ with probability 1, so this expectation is $\lesssim_{\mathrm{D}}\N^{-\mathrm{D}}$. This finishes the proof.
\end{proof}
The next step is a local equilibrium reduction like Lemma \ref{lemma:finalprop4}. The only difference in their proofs is power-counting.
\begin{lemma}\label{lemma:finalprop10}
 Adopt {the} notation of {Lemmas \ref{lemma:finalprop8}, \ref{lemma:finalprop9}}. With explanation given after, we have
\begin{align}
&{\textstyle\int_{\tau(\mathrm{j},1)}^{1}}|\mathbb{T}(\N)|^{-1}{\textstyle\sum_{\y}}\N\cdot\E[\{\E^{\mathrm{loc},\s}[\mathrm{Loc}^{(1)}\Psi^{\pm,\mathrm{j},\t}(\s)]\}(\Pi^{\mathrm{j},1}\mathbf{U}^{\s-\tau(\mathrm{j},1),\y(\s-\tau(\mathrm{j},1))+\cdot})]\d\s \nonumber\\
&\lesssim \ \N^{-30\beta_{\mathrm{BG}}}+\mathrm{LE}^{\t}(\Psi). \label{eq:finalprop10I}
\end{align}
Above, the local equilibrium error $\mathrm{LE}^{\t}(\Psi)$ is the following analog of \eqref{eq:finalprop4IIa}-\eqref{eq:finalprop4IIb}, where $\E^{\sigma,\mathrm{loc},\s}$ is from {Lemma \ref{lemma:finalprop4}}:
\begin{align}
\mathrm{LE}^{\t}(\Psi) \ := \ \sup_{\sigma\in\R}\sup_{\tau(\mathrm{j},1)\leq\s\leq1}\N\times\E^{\sigma,\mathrm{loc},\s}\mathrm{Loc}^{(1)}\Psi^{\pm,\mathrm{j},\t}(\s). \label{eq:finalprop10II}
\end{align}
\end{lemma}
\begin{proof}
We use Lemma \ref{lemma:le9} with the following choices. We choose the functional $\mathfrak{a}(\s,\mathbf{U})=\E^{\mathrm{loc},\s}[\mathrm{Loc}^{(1)}\Psi^{\pm,\mathrm{j},\t}(\s)](\Pi^{\mathrm{j},1}\mathbf{U})$ with $\mathbf{U}\in\R^{\mathbb{T}(\N)}$. It is a functional whose support is $\mathbb{I}(\mathfrak{t})$ for $\mathfrak{t}=\tau(\mathrm{j},1)$; see Lemma \ref{lemma:finalprop2} for this discrete interval. (The support claim holds since {\small$\Pi^{\mathrm{j},1}$} is the projection onto {\small$\R^{\mathbb{I}(\mathfrak{t})}$}; see Lemma \ref{lemma:finalprop2}. For $\kappa>0$ to be chosen shortly, this implies the following analog of \eqref{eq:finalprop4IIa1}:
\begin{align}
&\mathrm{LHS}\eqref{eq:finalprop10I} \nonumber\\
&\lesssim \ \tfrac{1}{\kappa}\N^{-\frac54-\gamma_{\mathrm{KL}}}|\mathbb{I}(\mathfrak{t})|^{3} + \tfrac{\N}{\kappa}\sup_{\sigma\in\R}\sup_{\tau(\mathrm{j},1)\leq\s\leq1}\log\E^{\sigma,\s-\tau(\mathrm{j},1),\mathbb{I}(\mathfrak{t})}\exp\{\kappa\E^{\mathrm{loc},\s}\mathrm{Loc}^{(1)}\Psi^{\pm,\mathrm{j},\t}(\s)\}. \label{eq:finalprop10I1}
\end{align}
Let us now choose $\kappa=\N^{-30\gamma_{\mathrm{reg}}-\gamma_{\mathrm{ap}}}\mathfrak{l}(\mathrm{j})^{3/2}$. By construction in \eqref{eq:finalprop9II}, for this choice of $\kappa$, the term inside the exponential in $\mathrm{RHS}\eqref{eq:finalprop10I1}$ is $\mathrm{O}(1)$. So, by the proof of \eqref{eq:finalprop4IIa2a}-\eqref{eq:finalprop4IIa2c}, we get the following in which $\mathrm{LD}$ is now the second term in \eqref{eq:finalprop10I1}:
\begin{align}
\mathrm{LD} \ \lesssim \ \N\times{\textstyle\sup_{\sigma,\s}}\E^{\sigma,\s-\tau(\mathrm{j},1),\mathbb{I}(\mathfrak{t})}\E^{\mathrm{loc},\s}\mathrm{Loc}^{(1)}\Psi^{\pm,\mathrm{j},\t}(\s) \ = \ \mathrm{LE}^{\t}(\Psi). \label{eq:finalprop10I2}
\end{align}
Let us now estimate the first term on the RHS of \eqref{eq:finalprop10I1}, which we denote by $\mathrm{Cost}$ in this proof. Recall from the paragraph after \eqref{eq:finalprop4IIa2a}-\eqref{eq:finalprop4IIa2c} that $|\mathbb{I}(\mathfrak{t})|\lesssim\N^{\gamma_{\mathrm{ap}}}[\N\tau(\mathrm{j},1)^{1/2}+\N^{3/2}\tau(\mathrm{j},1)+|\mathbb{I}|]$ for $\mathfrak{t}=\tau(\mathrm{j},1)$, and $|\mathbb{I}|\lesssim\mathfrak{m}(\mathrm{j},1)\mathfrak{l}(\mathrm{j})+\N^{3/2+\gamma_{\mathrm{ap}}}\tau(\mathrm{j},1)$. So, we deduce the following estimate, in which we use that $\gamma_{\mathrm{ap}},\gamma_{\mathrm{reg}}$ are small compared to $\gamma_{\mathrm{KL}}$ (see Definitions \ref{definition:reg}, \ref{definition:method8}):
\begin{align}
\mathrm{Cost} \ &\lesssim \ \N^{30\gamma_{\mathrm{reg}}+\gamma_{\mathrm{ap}}}\mathfrak{l}(\mathrm{j})^{-\frac32}\N^{-\frac54-\gamma_{\mathrm{KL}}}\N^{3\gamma_{\mathrm{ap}}}[\N^{3}\tau(\mathrm{j},1)^{\frac32}+\N^{\frac92+3\gamma_{\mathrm{ap}}}\tau(\mathrm{j},1)^{3}+\mathfrak{m}(\mathrm{j},1)^{3}\mathfrak{l}(\mathrm{j})^{3}] \label{eq:finalprop10I3a}\\
&\lesssim \ \N^{\frac74-\frac23\gamma_{\mathrm{KL}}}\tau(\mathrm{j},1)^{\frac32}\mathfrak{l}(\mathrm{j})^{-\frac32}+\N^{\frac{13}{4}-\frac23\gamma_{\mathrm{KL}}}\tau(\mathrm{j},1)^{3}\mathfrak{l}(\mathrm{j})^{-\frac32}+\N^{-\frac54-\frac23\gamma_{\mathrm{KL}}}\mathfrak{l}(\mathrm{j})^{-\frac32}\mathfrak{m}(\mathrm{j},1)^{3}\mathfrak{l}(\mathrm{j})^{3}.\label{eq:finalprop10I3b}
\end{align}
Suppose $\mathfrak{l}(\mathrm{j})\leq\mathfrak{m}(\mathrm{j},1)\mathfrak{l}(\mathrm{j})\leq\N^{1/9}$. In this case, by Definition \ref{definition:bg2131}. we know $\tau(\mathrm{j},1)=\N^{-2}\mathfrak{m}(\mathrm{j},1)^{2}\mathfrak{l}(\mathrm{j})^{2}\lesssim\N^{-16/9}$. Thus, 
\begin{align}
\eqref{eq:finalprop10I3b} \ &\lesssim \ \N^{\frac74-\frac83-\frac23\gamma_{\mathrm{KL}}}\mathfrak{l}(\mathrm{j})^{-\frac32}+\N^{\frac{13}{4}-\frac{16}{3}-\frac23\gamma_{\mathrm{KL}}}\mathfrak{l}(\mathrm{j})^{-\frac32}+\N^{-\frac54-\frac23\gamma_{\mathrm{KL}}+\frac13}\mathfrak{l}(\mathrm{j})^{-\frac32} \ \lesssim \ \N^{-\gamma_{\mathrm{KL}}}. \label{eq:finalprop10I4a}
\end{align}
We now suppose $\mathfrak{m}(\mathrm{j},1)\mathfrak{l}(\mathrm{j})\geq\N^{1/9}$. In this case, we know $\mathfrak{m}(\mathrm{j},1)=1$ by construction in Definition \ref{definition:bg2131}. Additionally, suppose that $\mathfrak{l}(\mathrm{j})\leq\N^{1/2}$. In this case, by construction in Definition \ref{definition:bg2131}, we know $\tau(\mathrm{j},1)=\N^{-2}\mathfrak{l}(\mathrm{j})^{2}$. Therefore, we deduce
\begin{align}
\eqref{eq:finalprop10I3b} \ &\lesssim \ \N^{\frac74-\frac23\gamma_{\mathrm{KL}}-3}\mathfrak{l}(\mathrm{j})^{\frac32}+\N^{\frac{13}{4}-\frac23\gamma_{\mathrm{KL}}-6}\mathfrak{l}(\mathrm{j})^{\frac92}+\N^{-\frac54-\frac23\gamma_{\mathrm{KL}}}\mathfrak{l}(\mathrm{j})^{\frac32} \label{eq:finalprop10I4b}\\
&\lesssim \ \N^{-\frac54-\frac23\gamma_{\mathrm{KL}}}\mathfrak{l}(\mathrm{j})^{\frac32}+\N^{-\frac{11}{4}-\frac23\gamma_{\mathrm{KL}}}\mathfrak{l}(\mathrm{j})^{\frac92} \nonumber\\
&\lesssim \ \N^{-\frac54-\frac23\gamma_{\mathrm{KL}}+\frac34}+\N^{-\frac{11}{4}-\frac23\gamma_{\mathrm{KL}}+\frac94} \ \lesssim \ \N^{-\frac12-\frac23\gamma_{\mathrm{KL}}}. \label{eq:finalprop10I4c}
\end{align}
(In order to derive \eqref{eq:finalprop10I4c}, we use the assumption $\mathfrak{l}(\mathrm{j})\leq\N^{1/2}$.) Now, suppose $\mathfrak{m}(\mathrm{j},1)\mathfrak{l}(\mathrm{j})\geq\N^{1/9}$ (so that $\mathfrak{m}(\mathrm{j},1)=1$) and that $\mathfrak{l}(\mathrm{j})\geq\N^{1/2}$. This is the last case. By Definition \ref{definition:bg2131}, we get $\tau(\mathrm{j},1)=\N^{-3/2}\mathfrak{l}(\mathrm{j})$ and $\mathfrak{l}(\mathrm{j})\leq\N^{3/4+\gamma_{\mathrm{KL}}}$. Therefore, we have 
\begin{align}
\eqref{eq:finalprop10I3b} \ &\lesssim \ \N^{\frac74-\frac23\gamma_{\mathrm{KL}}-\frac94}+\N^{\frac{13}{4}-\frac23\gamma_{\mathrm{KL}}-\frac92}\mathfrak{l}(\mathrm{j})^{\frac32}+\N^{-\frac54-\frac23\gamma_{\mathrm{KL}}}\mathfrak{l}(\mathrm{j})^{\frac32} \nonumber\\
&\lesssim \ \N^{-\frac12-\frac23\gamma_{\mathrm{KL}}}+\N^{-\frac54-\frac23\gamma_{\mathrm{KL}}}\N^{\frac98+\frac32\gamma_{\mathrm{KL}}} \nonumber\\
&\lesssim \ \N^{-\gamma_{\mathrm{KL}}}. \nonumber
\end{align}
(The last bound follows as $\gamma_{\mathrm{KL}}$ is small; see Definition \ref{definition:entropydata}.) Using every display in this proof gives \eqref{eq:finalprop10I}.
\end{proof}
We now give a technical adjustment analog to Lemma \ref{lemma:finalprop5}. Its proof is basically identical.
\begin{lemma}\label{lemma:finalprop11}
 Fix $1\leq\mathrm{j}\leq\mathrm{j}(\infty)$. Adopt the notation of {Lemmas \ref{lemma:finalprop9}, \ref{lemma:finalprop10}}. With notation explained after, we have {the following for any fixed $\mathrm{D}>0$}:
\begin{align}
\E^{\sigma,\mathrm{loc},\s}\mathrm{Loc}^{(1)}\Psi^{\pm,\mathrm{j},\t}(\s) \ \lesssim \ \mathbf{1}[|\sigma|\lesssim1]\E^{\sigma,\mathrm{loc},\s}\mathrm{Loc}\Psi^{\pm,\mathrm{j},\t}(\s)+\N^{-{\mathrm{D}}}. \label{eq:finalprop11I}
\end{align}
To explain $\mathrm{RHS}\eqref{eq:finalprop11I}$, take the notation from {Lemma \ref{lemma:finalprop5}} for $\mathrm{i}=1$. We then define
\begin{align}
\mathrm{Loc}\Psi^{\pm,\mathrm{j},\t}(\s) \ := \ |\mathscr{A}^{\mathrm{loc},\pm,\t}|\mathbf{1}[|\mathscr{A}^{\mathrm{loc},\pm,\t}|\gtrsim\N^{-\beta(\mathrm{j},1)-\gamma_{\mathrm{ap}}}]\mathbf{1}[|\mathscr{A}^{\mathrm{loc},\pm,\t}|\lesssim\N^{30\gamma_{\mathrm{reg}}+\gamma_{\mathrm{ap}}}\mathfrak{l}(\mathrm{j})^{-\frac32}]. \label{eq:finalprop11II}
\end{align}
\end{lemma}
\begin{proof}
First, we claim that $\mathrm{LHS}\eqref{eq:finalprop11I}\leq\mathbf{1}[|\sigma|\lesssim1]\mathrm{LHS}\eqref{eq:finalprop11I}$. This follows by the reasoning in Step 0 in the proof of Lemma \ref{lemma:finalprop5}. (Indeed, said argument holds for $\mathrm{i}=1$ as allowed in Lemma \ref{lemma:finalprop5}.) Next, we note that Steps 1 and 2 in the proof of Lemma \ref{lemma:finalprop5} hold for $\mathrm{i}=1$. We now claim $\mathrm{LHS}\eqref{eq:finalprop5IIa5a}\lesssim\eqref{eq:finalprop5IIa5c}$ holds upon replacing $\Upsilon^{\pm,\mathrm{j},\mathrm{i},\t}\mapsto\Psi^{\pm,\mathrm{j},\t}$, as all we need in the argument are the ingredients from Steps 1 and 2 in the proof of Lemma \ref{lemma:finalprop5} and the polynomial bound {\small$|\mathrm{Loc}^{(1)}\Psi^{\pm,\mathrm{j},\t}|\lesssim\N$}. (This estimate follows by construction in Lemma \ref{lemma:finalprop10}.) Combining this with the first sentence in this paragraph gives \eqref{eq:finalprop11I}.
\end{proof}
We conclude with a Kipnis-Varadhan estimate, like Lemma \ref{lemma:finalprop7}. The only difference in proofs is some power-counting.
\begin{lemma}\label{lemma:finalprop12}
 Fix $1\leq\mathrm{j}\leq\mathrm{j}(\infty)$. Retain the notation of {Lemma \ref{lemma:finalprop11}}. We have the following estimate uniformly in $\s,\t,\sigma$:
\begin{align}
\N\cdot\mathbf{1}[|\sigma|\lesssim1]\E^{\sigma,\mathrm{loc},\s}\mathrm{Loc}\Psi^{\pm,\mathrm{j},\t}(\s) \ \lesssim \ \N^{-30\beta_{\mathrm{BG}}}. \label{eq:finalprop12I}
\end{align}
\end{lemma}
\begin{proof}
We first apply the following Schwarz inequality, which decouples the indicators in \eqref{eq:finalprop11II}:
\begin{align}
\mathrm{Loc}\Psi^{\pm,\mathrm{j},\t}(\s) \ &\lesssim \ \N^{-\beta(\mathrm{j},1)}\mathbf{1}[|\mathscr{A}^{\mathrm{loc},\pm,\t}|\gtrsim\N^{-\beta(\mathrm{j},1)-\gamma_{\mathrm{ap}}}] \nonumber\\
&+ \ \N^{\beta(\mathrm{j},1)}|\mathscr{A}^{\mathrm{loc},\pm,\t}|^{2}\mathbf{1}[|\mathscr{A}^{\mathrm{loc},\pm,\t}|\lesssim\N^{30\gamma_{\mathrm{reg}}+\gamma_{\mathrm{ap}}}\mathfrak{l}(\mathrm{j})^{-\frac32}]. \label{eq:finalprop12I1}
\end{align}
Let $\Psi[1]$ be the first term in $\mathrm{RHS}\eqref{eq:finalprop12I1}$, and let $\Psi[2]$ be the second term therein. Note $\Psi[1]$ equals $\mathrm{Loc}\Lambda^{\pm,\mathrm{j},2,\t}(\s)$ from \eqref{eq:finalprop5IIb}. So, by Lemma \ref{lemma:finalprop7}, we get $\N\cdot\mathbf{1}[|\sigma|\lesssim1]\E^{\sigma,\mathrm{loc},\s}\Psi[1]\lesssim\N^{-30\beta_{\mathrm{BG}}}$. Thus, to prove \eqref{eq:finalprop12I}, it suffices to show
\begin{align}
\N\cdot\mathbf{1}[|\sigma|\lesssim1]\E^{\sigma,\mathrm{loc},\s}\Psi[2] \ \lesssim \ \N^{-30\beta_{\mathrm{BG}}}. \label{eq:finalprop12I2}
\end{align}
To this end, we use Proposition \ref{prop:kv1} with the following choices. First, take $\mathscr{B}=\N^{30\gamma_{\mathrm{reg}}+\gamma_{\mathrm{ap}}}\mathfrak{l}(\mathrm{j})^{-3/2}$. We then take every other choice in Proposition \ref{prop:kv1} to be what we chose to get \eqref{eq:finalprop7I1a}-\eqref{eq:finalprop7I1b} but specialized to $\mathrm{i}=1$. Ultimately, we deduce the bound below, which is an analog to \eqref{eq:finalprop7I1a}-\eqref{eq:finalprop7I1b} but with a different first line (coming from the different choice of $\mathscr{B}$) and basically same second line (with the harmless discrepancy in $\mathrm{i}$-indices by 1):
\begin{align}
\mathbf{1}[|\sigma|\lesssim1]\E^{\sigma,\mathrm{loc},\s}\Psi[2] \ \lesssim \ &\N^{\beta(\mathrm{j},1)+20\gamma_{\mathrm{reg}}}\N^{\gamma_{\mathrm{KV}}}\tau(\mathrm{j},1)|\mathbb{I}(\tau(\mathrm{j},1))|^{\frac12}\N^{60\gamma_{\mathrm{reg}}+2\gamma_{\mathrm{ap}}}\mathfrak{l}(\mathrm{j})^{-3} \label{eq:finalprop12I3a}\\
+ \ &\N^{\beta(\mathrm{j},1)+100\gamma_{\mathrm{reg}}}\N^{-2}\mathfrak{m}(\mathrm{j},1)^{-1}\tau(\mathrm{j},1)^{-1}\|\varphi\|^{2}\mathfrak{l}(\mathrm{j})^{2}\mathfrak{l}(\mathrm{j})^{-3}. \label{eq:finalprop12I3b}
\end{align}
(We clarify $\varphi$ was one of the choices made for Proposition \ref{prop:kv1} that we took from \eqref{eq:finalprop7I1a}-\eqref{eq:finalprop7I1b}.) We already have an estimate for \eqref{eq:finalprop12I3b}; see \eqref{eq:finalprop7I4} with $\mathrm{i}=1$ and the $\|\varphi\|$-estimate given right after \eqref{eq:finalprop7I4}. Including the extra factor of $\N^{\beta_{\mathrm{BG}}}$ to account for the discrepancy in $\mathrm{i}$-indices in the $\beta$-exponents in \eqref{eq:finalprop7I4} and \eqref{eq:finalprop12I3b}, we deduce
\begin{align}
\eqref{eq:finalprop12I3b} \ \lesssim \ \N^{-1-89\beta_{\mathrm{BG}}+120\gamma_{\mathrm{reg}}+10\gamma_{\mathrm{ap}}} \ \lesssim \ \N^{-1-88\beta_{\mathrm{BG}}}, \label{eq:finalprop12I4}
\end{align}
where the last estimate above follows since $\gamma_{\mathrm{reg}},\gamma_{\mathrm{ap}}$ are small compared to $\beta_{\mathrm{BG}}$ (see Definitions \ref{definition:reg}, \ref{definition:method8}). We now bound \eqref{eq:finalprop12I3a}. Recall from before \eqref{eq:finalprop10I3a} that $|\mathbb{I}(\mathfrak{t})|\lesssim\N^{\gamma_{\mathrm{ap}}}[\N\tau(\mathrm{j},1)^{1/2}+\N^{3/2}\tau(\mathrm{j},1)+|\mathbb{I}|]$ for $\mathfrak{t}=\tau(\mathrm{j},1)$, and $|\mathbb{I}|\lesssim\mathfrak{m}(\mathrm{j},1)\mathfrak{l}(\mathrm{j})+\N^{3/2+\gamma_{\mathrm{ap}}}\tau(\mathrm{j},1)$. Also, we recall from Definitions \ref{definition:reg} and \ref{definition:method8} that $\gamma_{\mathrm{reg}},\gamma_{\mathrm{ap}}$ are small compared to $\gamma_{\mathrm{KL}}$. Ultimately, we deduce the following, where the second line is a multiplication-by-1 (that can be verified directly), and the third line follows by construction of $\beta(\mathrm{j},1)$ in Definition \ref{definition:bg2131}:
\begin{align}
&\mathrm{RHS}\eqref{eq:finalprop12I3a} \nonumber\\
&\lesssim \ \N^{\beta(\mathrm{j},\mathrm{i})+2\gamma_{\mathrm{KV}}}\mathfrak{l}(\mathrm{j})^{-3}\tau(\mathrm{j},1)[\N^{\frac12}\tau(\mathrm{j},1)^{\frac14}+\N^{\frac34}\tau(\mathrm{j},1)^{\frac12}+\mathfrak{m}(\mathrm{j},1)^{\frac12}\mathfrak{l}(\mathrm{j})^{\frac12}] \nonumber\\
&\lesssim \ \N^{\beta(\mathrm{j},\mathrm{i})-1}\tau(\mathrm{j},1)^{-1}\mathfrak{m}(\mathrm{j},1)^{-1}\mathfrak{l}(\mathrm{j})^{-1}\N^{1+2\gamma_{\mathrm{KV}}}\mathfrak{l}(\mathrm{j})^{-3}\nonumber\\
&\times \ [\N^{\frac12}\tau(\mathrm{j},1)^{\frac94}\mathfrak{m}(\mathrm{j},1)\mathfrak{l}(\mathrm{j})+\N^{\frac34}\tau(\mathrm{j},1)^{\frac52}\mathfrak{m}(\mathrm{j},1)\mathfrak{l}(\mathrm{j})+\tau(\mathrm{j},1)^{2}\mathfrak{m}(\mathrm{j},1)^{\frac32}\mathfrak{l}(\mathrm{j})^{\frac32}] \nonumber\\
&\lesssim \ \N^{1+2\gamma_{\mathrm{KV}}-20\gamma_{\mathrm{reg}}-90\beta_{\mathrm{BG}}}\mathfrak{l}(\mathrm{j})^{-3}\nonumber\\
&\times[\N^{\frac12}\tau(\mathrm{j},1)^{\frac94}\mathfrak{m}(\mathrm{j},1)\mathfrak{l}(\mathrm{j})+\N^{\frac34}\tau(\mathrm{j},1)^{\frac52}\mathfrak{m}(\mathrm{j},1)\mathfrak{l}(\mathrm{j})+\tau(\mathrm{j},1)^{2}\mathfrak{m}(\mathrm{j},1)^{\frac32}\mathfrak{l}(\mathrm{j})^{\frac32}].\label{eq:finalprop12I5}
\end{align}
For now, assume $\mathfrak{l}(\mathrm{j})\leq\mathfrak{m}(\mathrm{j},1)\mathfrak{l}(\mathrm{j})\leq\N^{1/9}\leq\N^{1/2}$. By Definition \ref{definition:bg2131}, this means $\tau(\mathrm{j},1)=\N^{-2}\mathfrak{m}(\mathrm{j},1)^{2}\mathfrak{l}(\mathrm{j})^{2}$. Therefore,
\begin{align}
&\eqref{eq:finalprop12I5} \nonumber\\
&\lesssim \ \N^{1+2\gamma_{\mathrm{KV}}-20\gamma_{\mathrm{reg}}-90\beta_{\mathrm{BG}}}\mathfrak{l}(\mathrm{j})^{-3}[\N^{\frac12-\frac92}\mathfrak{m}(\mathrm{j},1)^{\frac92}\mathfrak{l}(\mathrm{j})^{\frac92}+\N^{\frac34-5}\mathfrak{m}(\mathrm{j},1)^{5}\mathfrak{l}(\mathrm{j})^{5}+\N^{-4}\mathfrak{m}(\mathrm{j},1)^{\frac92}\mathfrak{l}(\mathrm{j})^{\frac92}] \label{eq:finalprop12I6a}\\
&\lesssim \ \N^{-3+2\gamma_{\mathrm{KV}}}\mathfrak{l}(\mathrm{j})^{-3}\mathfrak{m}(\mathrm{j},1)^{\frac92}\mathfrak{l}(\mathrm{j})^{\frac92}+\N^{-\frac{13}{4}+2\gamma_{\mathrm{KV}}}\mathfrak{l}(\mathrm{j})^{-3}\mathfrak{m}(\mathrm{j},1)^{5}\mathfrak{l}(\mathrm{j})^{5} \nonumber\\
&\lesssim \ \N^{-\frac52+2\gamma_{\mathrm{KV}}}+\N^{-\frac94+2\gamma_{\mathrm{KV}}}\lesssim\N^{-1-100\beta_{\mathrm{BG}}}, \label{eq:finalprop12I6b}
\end{align}
{where the last two bounds use $\mathfrak{l}(\mathrm{j})\geq1$, $\mathfrak{m}(\mathrm{j},1)\mathfrak{l}(\mathrm{j})\leq\N^{1/9}$, and $\gamma_{\mathrm{KV}},\beta_{\mathrm{BG}}\leq c$ for some small $c>0$} (see Definitions \ref{definition:entropydata}, \ref{definition:method8}, and Proposition \ref{prop:kv1}). Now assume $\N^{1/9}\leq\mathfrak{l}(\mathrm{j})\leq\N^{1/2}$. By Definition \ref{definition:bg2131}, this means $\mathfrak{m}(\mathrm{j},1)=1$ and $\tau(\mathrm{j},1)=\N^{-2}\mathfrak{l}(\mathrm{j})^{2}$. Therefore, we get
\begin{align}
\eqref{eq:finalprop12I5} \ &\lesssim \ \N^{1+2\gamma_{\mathrm{KV}}}\mathfrak{l}(\mathrm{j})^{-3}[\N^{\frac12-\frac92}\mathfrak{l}(\mathrm{j})^{\frac92}+\N^{\frac34-5}\mathfrak{l}(\mathrm{j})^{5}+\N^{-4}\mathfrak{l}(\mathrm{j})^{\frac92}] \nonumber\\
&\lesssim \ \N^{-3+2\gamma_{\mathrm{KV}}}\mathfrak{l}(\mathrm{j})^{\frac32}+\N^{-\frac{13}{4}+2\gamma_{\mathrm{KV}}}\mathfrak{l}(\mathrm{j})^{2} \ \lesssim \ \N^{-1-100\beta_{\mathrm{BG}}}; \nonumber
\end{align}
the last estimates use $\mathfrak{l}(\mathrm{j})\leq\N^{1/2}$ and $\gamma_{\mathrm{KV}},\beta_{\mathrm{BG}}\leq c$ for some small $c>0$. Now, we assume $\mathfrak{l}(\mathrm{j})\geq\N^{1/2}$. This is the last case left. By Definition \ref{definition:bg2131}, we have $\mathfrak{m}(\mathrm{j},1)=1$ as before. We also have $\tau(\mathrm{j},1)=\N^{-3/2}\mathfrak{l}(\mathrm{j})$ and $\mathfrak{l}(\mathrm{j})\leq\N^{3/4+\gamma_{\mathrm{KL}}}$; again, see Definition \ref{definition:bg2131}. So,
\begin{align}
\eqref{eq:finalprop12I5} \ &\lesssim \ \N^{1+2\gamma_{\mathrm{KV}}}\mathfrak{l}(\mathrm{j})^{-3}[\N^{\frac12-\frac{27}{8}}\mathfrak{l}(\mathrm{j})^{\frac{13}{4}}+\N^{\frac34-\frac{15}{4}}\mathfrak{l}(\mathrm{j})^{\frac72}+\N^{-3}\mathfrak{l}(\mathrm{j})^{\frac72}] \nonumber\\
&\lesssim \ \N^{-\frac{15}{8}+2\gamma_{\mathrm{KV}}}\mathfrak{l}(\mathrm{j})^{\frac14}+\N^{-2+2\gamma_{\mathrm{KV}}}\mathfrak{l}(\mathrm{j})^{\frac12} \ \lesssim \ \N^{-1-100\beta_{\mathrm{BG}}}. \nonumber
\end{align}
Combining the previous four displays shows that $\mathrm{RHS}\eqref{eq:finalprop12I3a}\lesssim\N^{-1-100\beta_{\mathrm{BG}}}$. Combine this with \eqref{eq:finalprop12I3a}-\eqref{eq:finalprop12I3b} and \eqref{eq:finalprop12I4}. This gives \eqref{eq:finalprop12I2}. As we explained right before \eqref{eq:finalprop12I2}, this gives the desired estimate \eqref{eq:finalprop12I}, so the proof is complete.
\end{proof}
\subsection{Proof of Proposition \ref{prop:bg21310}}
To bound $\mathrm{RHS}\eqref{eq:bg2138Ia}$, use \eqref{eq:bg2138Ia2a}-\eqref{eq:bg2138Ia2b}. For \eqref{eq:bg2138Ic} and \eqref{eq:bg2138Id}, use \eqref{eq:bg213101Ia}-\eqref{eq:bg213101Ib}, \eqref{eq:finalprop2Ia}-\eqref{eq:finalprop2Ib}, \eqref{eq:finalprop4Ia}-\eqref{eq:finalprop4Ib}, \eqref{eq:finalprop5Ia}-\eqref{eq:finalprop5Ib}, and \eqref{eq:finalprop7I}. For \eqref{eq:bg2138Ib}, use \eqref{eq:finalprop8I}, \eqref{eq:finalprop9I}, \eqref{eq:finalprop10I}, \eqref{eq:finalprop11I}, and \eqref{eq:finalprop12I}. \qed
%
%
%
\section{Proof of Proposition \ref{prop:bg212}}\label{section:propbg212proof}
This is basically a time-version of Lemma \ref{lemma:bg210}. The main difference is the need of a technical stochastic bound (Lemma \ref{lemma:bg2121}), which takes advantage of the fact that we are not just trying to replace a general function by its local time-average inside the heat operator, but rather a function given by the average of fluctuating terms. See the paragraphs after the statement of Proposition \ref{prop:bg212} for a more precise explanation of this heuristic. (Ultimately, what we do below is essentially an easier version of the analysis in the previous section.)

Like we did in the proof of Proposition \ref{prop:bg213}, we restrict to $\d=0$ in \eqref{eq:bg212I} solely out of convenience. (The proof of \eqref{eq:bg212I} for $\d=1,2,3,4$ follows by the same argument. In particular, formally replace the heat operator $\mathbf{H}^{\N}$ by its composition $(\mathscr{T}^{\pm,\mathrm{j}})^{\d}\mathbf{H}^{\N}$ with powers of the differential operator $\mathscr{T}^{\pm,\mathrm{j}}$. As for why this works, see the beginning of Section \ref{section:bg213main}.)
\subsection{The ingredients}
This argument needs a few preliminaries that we now present. (We then use these to prove Proposition \ref{prop:bg212}. We finish the section by presenting proofs of each ingredient used.) The first is a stochastic estimate that we explain after.
\begin{lemma}\label{lemma:bg2121}
{Consider a collection $\{\mathfrak{A}^{(\mathrm{k})}(\tau,\mathbf{U})\}_{\mathrm{k}\in\mathds{B}}$ (for some fixed countable index set $\mathds{B}$) of functionals of $\tau\geq0$ and $\mathbf{U}\in\R^{\mathbb{T}(\N)}$, which satisfy the following assumptions.} First, they satisfy the constraints in Lemma \ref{lemma:le2} with respect to discrete intervals $\mathbb{J}(\mathrm{k})\subseteq\mathbb{T}(\N)$. We also assume that $|\mathfrak{A}^{(\mathrm{k})}(\tau,\mathbf{U})|\lesssim\N^{30\gamma_{\mathrm{reg}}}\mathfrak{l}(\mathrm{j})^{-3/2}$ for all $\tau\geq0$ and $\mathbf{U}\in\R^{\mathbb{T}(\N)}$. Let us now define the average
\begin{align}
\mathfrak{A}^{\mathds{B}}(\s,\mathbf{U}^{\s,\y+\cdot}) \ := \ \mathfrak{A}^{\mathds{B}}(\s,\y) \ := \ |\mathds{B}|^{-1}\sum_{\mathrm{k}\in\mathds{B}}\mathfrak{A}^{(\mathrm{k})}(\s,\mathbf{U}^{\s,\y+\cdot}){.}\label{eq:bg21210new}
\end{align}
Suppose $|\mathds{B}|\mathfrak{l}(\mathrm{j})=\N^{1/2+\delta_{\mathrm{KL}}}$, for fixed $20^{-1}\gamma_{\mathrm{KL}}\leq\delta_{\mathrm{KL}}\leq10^{-1}\gamma_{\mathrm{KL}}$. Moreover, suppose that there exists a discrete interval of length $\lesssim|\mathds{B}|\mathfrak{l}(\mathrm{j})$ that contains both $\mathbb{J}(\mathrm{k})$ and the support of {\small$\mathfrak{A}^{(\mathrm{k})}(\tau,\mathbf{U})$} for all $\mathrm{k}\in\mathds{B}$ (for fixed $\tau\geq0$). {Then, we} have the estimate
\begin{align}
\E{\textstyle\int_{0}^{1}}|\mathbb{T}(\N)|^{-1}{\textstyle\sum_{\y}}|\mathfrak{A}^{\mathds{B}}(\s,\y(\s))|\mathbf{1}[|\mathfrak{A}^{\mathds{B}}(\s,\y(\s))|\gtrsim\N^{-\frac14-\frac{\gamma_{\mathrm{KL}}}{100}}]\d\s \ \lesssim \ \N^{-1-100\beta_{\mathrm{BG}}}. \label{eq:bg2121I}
\end{align}
\end{lemma}
Lemma \ref{lemma:bg2121} is simply saying that at canonical measures, Lemma \ref{lemma:le2} gives square-root cancellation at large-deviations scale. By the local equilibrium reduction estimate of Lemma \ref{lemma:le9}, we can also show quantitative estimates without the assumption of local equilibrium. In particular, the proof of Lemma \ref{lemma:bg2121} is similar to proofs of Lemmas \ref{lemma:finalprop4}, \ref{lemma:finalprop7}. {However, it is much easier here because the length-scales on which we need to reduce to local equilibrium are much smaller. Moreover, the local equilibrium estimate itself is not dynamical.} (Indeed, we only need to apply Lemma \ref{lemma:le2}, not the dynamical estimates from Proposition \ref{prop:kv1}.)

The second ingredient is totally elementary; it bounds the error of introducing time-average inside a space-time integration by moving the resulting ``time-gradients" onto the heat kernel (in the spirit of integration-by-parts).
\begin{lemma}\label{lemma:bg2122}
 Take any $\phi:\R\times\mathbb{T}(\N)\to\R$ and $\tau>\N^{-{\mathrm{C}}}$, {where $\mathrm{C}>0$ is any fixed constant}. For any $0\leq\t\leq\t_{\mathrm{st}}$ and $\x\in\mathbb{T}(\N)$, we have the following estimate:
\begin{align}
|{\textstyle\int_{0}^{\t}}\mathbf{H}^{\N}(\s,\t(\N),\x)[\phi_{\s,\cdot(\s)}]\d\s-{\mathbf{1}_{\t\geq\tau}}{\textstyle\int_{\tau}^{\t}}\mathbf{H}^{\N}(\s,\t(\N),\x)[\tau^{-1}{\textstyle\int_{0}^{\tau}}\phi_{\s-\r,\cdot(\s-\r)}\d\r]\d\s| \ \lesssim \ \N^{100\gamma_{\mathrm{reg}}}\tau\|\phi\|_{\t_{\mathrm{st}};\mathbb{T}(\N)}. \label{eq:bg2122I}
\end{align}
\end{lemma}
\subsection{Proof of Proposition \ref{prop:bg212}}
Recall that we want to show \eqref{eq:bg212I}. This argument is basically using Lemma \ref{lemma:bg2121} to introduce an a priori cutoff for the spatial-average in $\mathscr{A}^{\mathbf{X}}\mathscr{R}^{\chi,\mathfrak{q},\pm,\mathrm{j}}$ in \eqref{eq:bg212I} (see Definition \ref{definition:bg29}), using Lemma \ref{lemma:bg2122} to introduce the time-average that turns $\mathscr{A}^{\mathbf{X}}\mathscr{R}^{\chi,\mathfrak{q},\pm,\mathrm{j}}$ into $\mathscr{A}^{\mathbf{X},\mathbf{T}}\mathscr{R}^{\chi,\mathfrak{q},\pm,\mathrm{j}}$ in \eqref{eq:bg212I}, and removing the cutoff by again using Lemma \ref{lemma:bg2121}. (We need to additionally freeze the coupling constant $\lambda(\t)$ in the time-average in $\mathscr{A}^{\mathbf{X},\mathbf{T}}\mathscr{R}^{\chi,\mathfrak{q},\pm,\mathrm{j}}$. This will be done by smoothness of $\lambda(\t)$ in $\t$ as noted after Definition \ref{definition:bg211}, as well as another technical and perhaps uninteresting cutoff argument based on Lemma \ref{lemma:bg2121}.) To make the presentation clearer, we will consider each step separately. First, however, it will be convenient to set $\|\|:=\|\|_{\t_{\mathrm{st}};\mathbb{T}(\N)}$ throughout this proof.
\subsubsection{The case $\mathfrak{l}(\mathrm{j})\geq\N^{1/2}$}
We focus on this case first as the $\mathds{B}$-set in Lemma \ref{lemma:bg2121} does not quite make sense if $\mathfrak{l}(\mathrm{j})$ is slightly bigger (in powers of $\N$) than $\N^{1/2}$. It is notationally annoying to single out this case when making the previous paragraph precise, so we deal with this case directly. (The only difference between this case and $\mathfrak{l}(\mathrm{j})\leq\N^{1/2}$ is a technical need for a priori estimates via Lemma \ref{lemma:bg2121} that come for free if $\mathfrak{l}(\mathrm{j})\geq\N^{1/2}$ by construction of the $\chi$-cutoffs in Definition \ref{definition:bg26}.) By Lemma \ref{lemma:bg2122}, we deduce the following estimates, in which we take the choices $\phi_{\s,\y}=\N{\mathds{A}_{\mathbf{X}}^{\mathfrak{m}(\mathrm{j}),\pm}}[\mathds{R}^{\chi,\mathfrak{q},\pm,\mathrm{j}}\mathbf{Z};\s,\y]$ and $\tau=\tau(\mathrm{j})$ from Definition \ref{definition:bg211}:
\begin{align}
\mathrm{LHS}\eqref{eq:bg2122I} \ &\lesssim \ \N^{\gamma_{\mathrm{ap}}}\tau(\mathrm{j})\|\N{\mathds{A}_{\mathbf{X}}^{\mathfrak{m}(\mathrm{j}),\pm}}[\mathds{R}^{\chi,\mathfrak{q},\pm,\mathrm{j}}\mathbf{Z};\t,{\y(\x,\t)})]\| \nonumber\\
&\lesssim \ \N^{\frac14+c\gamma_{\mathrm{KL}}}\|{\mathds{A}_{\mathbf{X}}^{\mathfrak{m}(\mathrm{j}),\pm}}[\mathds{R}^{\chi,\mathfrak{q},\pm,\mathrm{j}}\mathbf{Z};\t,{\y(\x,\t)})]\|.\label{eq:bg2121}
\end{align}
(In \eqref{eq:bg2121}, the norms are with respect to $(\t,\x)$-variables, {and the constant $c>0$ is small}.) We clarify {that} the last bound in \eqref{eq:bg2121} follows by construction of $\tau(\mathrm{j})$ from Definition \ref{definition:bg211}. We now bound the last norm in \eqref{eq:bg2121}. Recall from Definition \ref{definition:bg29} that what is inside the norm is an average of products between $\mathds{R}^{\chi,\mathfrak{q},\pm,\mathrm{j}}$ and $\mathbf{Z}$ evaluated at time before $\t_{\mathrm{st}}$. By Definitions \ref{definition:method8} and \ref{definition:bg26}, these are $\lesssim\N^{\gamma_{\mathrm{ap}}+20\gamma_{\mathrm{reg}}}\mathfrak{l}(\mathrm{j}-1)^{-3/2}\lesssim\N^{-1/4-\gamma_{\mathrm{KL}}/100}$, since $\mathfrak{l}(\mathrm{j}-1)\gtrsim\N^{-\gamma_{\mathrm{KL}}}\mathfrak{l}(\mathrm{j})$ by Definition \ref{definition:bg24} and by smallness of $\gamma_{\mathrm{ap}},\gamma_{\mathrm{reg}},\gamma_{\mathrm{KL}}$. {We emphasize that it is this upper bound $\N^{\gamma_{\mathrm{ap}}+20\gamma_{\mathrm{reg}}}\mathfrak{l}(\mathrm{j}-1)^{-3/2}\lesssim\N^{-1/4-\gamma_{\mathrm{KL}}/100}$ which uses the assumption $\mathfrak{l}(\mathrm{j})\geq\N^{1/2}$.} Thus, from \eqref{eq:bg2121}, we deduce {the following for some large constant $\mathrm{C}>0$}:
\begin{align}
\mathrm{LHS}\eqref{eq:bg2122I} \ \lesssim \ \N^{100\gamma_{\mathrm{reg}}+c\gamma_{\mathrm{KL}}-\frac{\gamma_{\mathrm{KL}}}{100}} \ \lesssim \ \N^{-{\mathrm{D}}\beta_{\mathrm{BG}}}, \label{eq:bg2121b}
\end{align}
where the last bound follows because $\gamma_{\mathrm{reg}},\beta_{\mathrm{BG}}$ are small compared to $\gamma_{\mathrm{KL}}$. We now compute as follows (recall $\mathbf{G}^{\s}$ from Definition \ref{definition:bg211}):
{\small
\begin{align}
&\N^{-1}\tau(\mathrm{j})^{-1}{\textstyle\int_{0}^{\tau(\mathrm{j})}}\phi_{\s-\r,\y(\s-\r)}\d\r \nonumber\\
&= \ \tau(\mathrm{j})^{-1}{\textstyle\int_{0}^{\tau(\mathrm{j})}}\{\mathfrak{m}(\mathrm{j})^{-1}{\textstyle\sum_{\mathrm{k}=0}^{\mathfrak{m}(\mathrm{j})-1}}\mathds{R}^{\chi,\mathfrak{q},\pm,\mathrm{j}}(\s-\r,\y(\s-\r)\pm2\mathrm{k}\mathfrak{l}(\mathrm{j}))\mathbf{G}(\s-\r,\y(\s-\r)\pm2\mathrm{k}\mathfrak{l}(\mathrm{j}))\}\d\r \nonumber\\
&= \ \tau(\mathrm{j})^{-1}{\textstyle\int_{0}^{\tau(\mathrm{j})}}\d\r\{\mathfrak{m}(\mathrm{j})^{-1}{\textstyle\sum_{\mathrm{k}=0}^{\mathfrak{m}(\mathrm{j})-1}}\mathds{R}^{\chi,\mathfrak{q},\pm,\mathrm{j}}(\s-\r,\y(\s-\r)\pm2\mathrm{k}\mathfrak{l}(\mathrm{j}))\mathbf{G}^{\s}(\s-\r,\y(\s-\r)\pm2\mathrm{k}\mathfrak{l}(\mathrm{j}))\nonumber\\
&\quad\quad\quad\quad\quad\quad\quad\quad\quad\quad\quad\quad\quad\quad\times\frac{\mathbf{G}(\s-\r,\y(\s-\r)\pm2\mathrm{k}\mathfrak{l}(\mathrm{j}))}{\mathbf{G}^{\s}(\s-\r,\y(\s-\r)\pm2\mathrm{k}\mathfrak{l}(\mathrm{j}))}\}\nonumber\\
&= \ \tau(\mathrm{j})^{-1}{\textstyle\int_{0}^{\tau(\mathrm{j})}}\{\mathfrak{m}(\mathrm{j})^{-1}{\textstyle\sum_{\mathrm{k}=0}^{\mathfrak{m}(\mathrm{j})-1}}\mathds{R}^{\chi,\mathfrak{q},\pm,\mathrm{j}}(\s-\r,\y(\s-\r)\pm2\mathrm{k}\mathfrak{l}(\mathrm{j}))\mathbf{G}^{\s}(\s-\r,\y(\s-\r)\pm2\mathrm{k}\mathfrak{l}(\mathrm{j}))\}\d\r \label{eq:bg2122}\\
&+ \ \tau(\mathrm{j})^{-1}{\textstyle\int_{0}^{\tau(\mathrm{j})}}\d\r\{\mathfrak{m}(\mathrm{j})^{-1}{\textstyle\sum_{\mathrm{k}=0}^{\mathfrak{m}(\mathrm{j})-1}}\mathds{R}^{\chi,\mathfrak{q},\pm,\mathrm{j}}(\s-\r,\y(\s-\r)\pm2\mathrm{k}\mathfrak{l}(\mathrm{j}))\mathbf{G}^{\s}(\s-\r,\y(\s-\r)\pm2\mathrm{k}\mathfrak{l}(\mathrm{j}))\nonumber\\
&\quad\quad\quad\quad\quad\quad\quad\quad\quad\quad\quad\quad\quad\quad\times[1-\frac{\mathbf{G}(\s-\r,\y(\s-\r)\pm2\mathrm{k}\mathfrak{l}(\mathrm{j}))}{\mathbf{G}^{\s}(\s-\r,\y(\s-\r)\pm2\mathrm{k}\mathfrak{l}(\mathrm{j}))}]\}. \nonumber
\end{align}
}We assume $\s\leq\t_{\mathrm{st}}$. The last factor in the last line has the form $1-\exp[\lambda(\s-\r)\mathrm{a}]\exp[-\lambda(\s)\mathrm{a}]$ for $|\mathrm{a}|\lesssim\log\log\N$ by Definition \ref{definition:method8}. By smoothness of $\lambda(\t)$ in $\t$ (see Assumption \ref{ass:intro8}), we know $|\lambda(\s-\r)-\lambda(\s)|\lesssim|\r|\lesssim\tau(\mathrm{j})$. Therefore, an elementary Taylor expansion estimate shows that the last factor in the last line is $\lesssim\N^{\gamma_{\mathrm{ap}}}\tau(\mathrm{j})$. By Definition \ref{definition:method8}, we know the $\mathbf{G}^{\s}$-factor in the last line is $\lesssim\N^{\gamma_{\mathrm{ap}}}$ for $\s\leq\t_{\mathrm{st}}$, because it is the exponential of something $\lesssim\log\log\N$.  Thus, by the a priori estimates from Definition \ref{definition:bg26} for the $\mathds{R}^{\chi,\mathfrak{q},\pm,\mathrm{j}}$-factor in the last line, we deduce the last line is $\lesssim\N^{10\gamma_{\mathrm{ap}}+20\gamma_{\mathrm{reg}}}\mathfrak{l}(\mathrm{j}-1)^{-3/2}\tau(\mathrm{j})\lesssim\N^{-1-\gamma_{\mathrm{KL}}/100}$, similar to the reasoning given immediately before \eqref{eq:bg2121b}. (This estimate is deterministic and uniform in $\s\leq\t_{\mathrm{st}}$ and $\y$.) As $\mathbf{H}^{\N}$-operators are contractive, we deduce from this paragraph, the previous display, and \eqref{eq:bg2121b} (with $\phi_{\s,\y}=\N{\mathds{A}_{\mathbf{X}}^{\mathfrak{m}(\mathrm{j}),\pm}}[\mathds{R}^{\chi,\mathfrak{q},\pm,\mathrm{j}}\mathbf{Z};\s,\y]$) that
\begin{align}
|\mathscr{A}^{\mathbf{X},\pm}\mathscr{R}^{\chi,\mathfrak{q},\pm,\mathrm{j}}(\t,\x)-{\textstyle\int_{0}^{\t}}\mathbf{H}^{\N}(\s,\t(\N),\x)[\eqref{eq:bg2122}]\d\s| \ \lesssim \ \N^{-{\mathrm{C}}\beta_{\mathrm{BG}}}+\N^{-100\gamma_{\mathrm{KL}}} \ \lesssim \ \N^{-{\mathrm{C}}\beta_{\mathrm{BG}}}, \label{eq:bg2123}
\end{align}
where the last bound uses the $\beta_{\mathrm{BG}}$ is small compared to $\gamma_{\mathrm{KL}}$. The second term in the absolute value above is $\mathscr{A}^{\mathbf{X},\mathbf{T}}\mathscr{R}^{\chi,\mathfrak{q},\pm,\mathrm{j}}$; see Definition \ref{definition:bg211}. As the previous bound is deterministic and holds uniformly over $\t\leq\t_{\mathrm{st}}$ and $\x\in\mathbb{T}(\N)$, the desired bound \eqref{eq:bg212I} (for $\d=0$) follows. (Again, bounds for $\d=1,2,3,4$ follow by the same argument upon replacing the heat kernel by $(\mathscr{T}^{\pm,\mathrm{j}})^{\d}$ acting on it.)
\subsubsection{A priori cutoff via {Lemma \ref{lemma:bg2121}}}
We first set some notation. Take the discrete interval $\llbracket0,\mathfrak{m}(\mathrm{j})-1\rrbracket$. Write it as a disjoint union of shifts of $\mathds{X}$, where {$\mathds{X}\subseteq\mathbb{T}(\N)$ is a discrete interval that satisfies $|\mathds{X}|\lesssim|\mathds{B}|\mathfrak{l}(\mathrm{j})$}. Let us denote the set of such shifts of $\mathds{X}$ by $\mathscr{X}$. Now, we use this notation to rewrite $\phi_{\s,\y}=\N{\mathds{A}_{\mathbf{X}}^{\mathfrak{m}(\mathrm{j}),\pm}}[\mathds{R}^{\chi,\mathfrak{q},\pm,\mathrm{j}}\mathbf{Z};\s,\y]$ below, which we justify afterwards:
\begin{align}
&{\mathds{A}_{\mathbf{X}}^{\mathfrak{m}(\mathrm{j}),\pm}}[\mathds{R}^{\chi,\mathfrak{q},\pm,\mathrm{j}}\mathbf{Z};\s,\y(\s)] \nonumber\\
&= \ |\mathscr{X}|^{-1}\sum_{\mathds{B}\in\mathscr{X}}|\mathds{B}|^{-1}\sum_{\mathrm{k}\in\mathds{B}}\mathds{R}^{\chi,\mathfrak{q},\pm,\mathrm{j}}(\s,\y(\s)\pm2\mathrm{k}\mathfrak{l}(\mathrm{j}))\mathbf{G}(\s,\y(\s)\pm2\mathrm{k}\mathfrak{l}(\mathrm{j})). \label{eq:bg2124a}
\end{align}
(Indeed, by Definition \ref{definition:bg28}, $\mathrm{LHS}\eqref{eq:bg2124a}$ is the average of the summands in $\mathrm{RHS}\eqref{eq:bg2124a}$ over $\mathrm{k}\in\llbracket0,\mathfrak{m}(\mathrm{j})-1\rrbracket$. This is the same as averaging over the set of mutually disjoint shifts $\mathds{B}$ of $\mathds{X}$ whose union equals $\llbracket0,\mathfrak{m}(\mathrm{j})-1\rrbracket$ and then averaging over $\mathrm{k}$ in each $\mathds{B}$.) We now do some gymnastics with $\mathrm{RHS}\eqref{eq:bg2124a}$ similar to those in Lemma \ref{lemma:finalprop5}. (Except, it is much easier here, because we never need to use the ``local" SDE \eqref{eq:glsdeloc}, so we can reason with \eqref{eq:glsde} itself and the stopping time $\t_{\mathrm{st}}$. Also, there is no time-averaging here.) Fix $\mathds{B}\in\mathscr{X}$. Recall $\mathbf{G}$ in Definition \ref{definition:intro6}, and recall the gradient relation $\mathbf{U}^{\t,\x}=\N^{1/2}[\mathbf{J}(\t,\x)-\mathbf{J}(\t,\x-1)]$; see Definition \ref{definition:intro4}. Now, restrict to $\s\leq\t_{\mathrm{st}}$. We also set $\mathfrak{k}[\mathds{B}]=\sup\mathds{B}$. (The role of this choice is to satisfy two conditions. First, $\pm2\mathrm{k}\mathfrak{l}(\mathrm{j})$ is always ``to the right" of $-2\mathfrak{k}[\mathds{B}]\mathfrak{l}(\mathrm{j})$ for any $\mathrm{k}\in\mathds{B}$. Second, the geodesic distance between $\pm2\mathrm{k}\mathfrak{l}(\mathrm{j})$ and $-2\mathfrak{k}[\mathds{B}]\mathfrak{l}(\mathrm{j})$ is $\lesssim|\mathds{B}|\mathfrak{l}(\mathrm{j})$. In particular, the smallest discrete interval containing $\pm2\mathrm{k}\mathfrak{l}(\mathrm{j})$ and $-2\mathfrak{k}[\mathds{B}]\mathfrak{l}(\mathrm{j})$ is $\lesssim|\mathds{B}|\mathfrak{l}(\mathrm{j})$ {in} size. Also, its infimum is $-2\mathfrak{k}[\mathds{B}]\mathfrak{l}(\mathrm{j})$ while its supremum is $\pm2\mathrm{k}\mathfrak{l}(\mathrm{j})$.) With explanation after, we now claim the $(\mathds{B},\mathrm{k})$-summand in $\mathrm{RHS}\eqref{eq:bg2124a}$ equals 
\begin{align}
&\mathbf{G}(\s,\y(\s)-2\mathfrak{k}[\mathds{B}]\mathfrak{l}(\mathrm{j}))\times\mathds{R}^{\chi,\mathfrak{q},\pm,\mathrm{j}}(\s,\y(\s)\pm2\mathrm{k}\mathfrak{l}(\mathrm{j}))\exp[\lambda(\s)\N^{-\frac12}{\textstyle\sum_{\w=-2\mathfrak{k}[\mathds{B}]\mathfrak{l}(\mathrm{j})+1}^{\pm2\mathrm{k}\mathfrak{l}(\mathrm{j})}}\mathbf{U}^{\s,\y(\s)+\w}]\label{eq:bg2125a}\\
&= \ \mathbf{G}(\s,\y(\s)-2\mathfrak{k}[\mathds{B}]\mathfrak{l}(\mathrm{j}))\times\mathds{R}^{\chi,\mathfrak{q},\pm,\mathrm{j}}(\s,\y(\s)\pm2\mathrm{k}\mathfrak{l}(\mathrm{j}))\label{eq:bg2125b}\\
&\times\mathrm{CutExp}[\lambda(\s)\N^{-\frac12}{\textstyle\sum_{\w=-2\mathfrak{k}[\mathds{B}]\mathfrak{l}(\mathrm{j})+1}^{\pm2\mathrm{k}\mathfrak{l}(\mathrm{j})}}\mathbf{U}^{\s,\y(\s)+\w}].\nonumber
\end{align}
\eqref{eq:bg2125a} follows by the definition of $\mathbf{G}$ and the gradient relation for \eqref{eq:hf}-\eqref{eq:glsde} from the previous paragraph. In \eqref{eq:bg2125b}, we used $\mathrm{CutExp}(\mathrm{a})=\exp(\mathrm{a})\mathbf{1}[\exp(\mathrm{a})\lesssim1]$ from \eqref{eq:finalprop5Ia4b}. To get \eqref{eq:bg2125b}, it suffices to show $\mathrm{CutExp}(\mathrm{a})=\exp(\mathrm{a})$ for appropriate $\mathrm{a}$ above. Equivalently, it suffices to show that $\lambda(\s)\N^{-1/2}$ times the sum in \eqref{eq:bg2125a}-\eqref{eq:bg2125b} is $\lesssim1$. To this end, we recall $\s\leq\t_{\mathrm{reg}}$. By Remark \ref{remark:intro14}, this means the sum in \eqref{eq:bg2125a}-\eqref{eq:bg2125b} is $\lesssim\N^{\gamma_{\mathrm{reg}}}|\pm2\mathrm{k}\mathfrak{l}(\mathrm{j})-2\mathfrak{k}[\mathds{B}]\mathfrak{l}(\mathrm{j})|^{1/2}+\N^{\gamma_{\mathrm{reg}}}\lesssim\N^{\gamma_{\mathrm{reg}}}|\mathds{B}|^{1/2}\mathfrak{l}(\mathrm{j})^{1/2}+\N^{\gamma_{\mathrm{reg}}}\lesssim\N^{1/4+\delta_{\mathrm{KL}}+\gamma_{\mathrm{reg}}}$, where $\delta_{\mathrm{KL}}\leq10^{-1}\gamma_{\mathrm{KL}}$ is from Lemma \ref{lemma:bg2121}. Because $\lambda(\t)$ is smooth in $\t$, we deduce the sum in \eqref{eq:bg2125a}-\eqref{eq:bg2125b} is $\lesssim1$, thereby giving \eqref{eq:bg2125b}. Now let {\small$\mathfrak{A}^{(\mathrm{k})}(\s,\mathbf{U})$} denote the functional of $\mathbf{U}\in\R^{\mathbb{T}(\N)}$ such that {\small$\mathfrak{A}^{(\mathrm{k})}(\s,\mathbf{U}^{\s,\y(\s)+\cdot})$} equals everything in \eqref{eq:bg2125b} after the $\times$ symbol. The above display lets us write \eqref{eq:bg2124a} as follows (where sums are still over $\mathds{B}\in\mathscr{X}$ and $\mathrm{k}\in\mathds{B}$):
\begin{align}
{\mathds{A}_{\mathbf{X}}^{\mathfrak{m}(\mathrm{j}),\pm}}[\mathds{R}^{\chi,\mathfrak{q},\pm,\mathrm{j}}\mathbf{Z};\s,\y(\s)] \ = \ |\mathscr{X}|^{-1}{\textstyle\sum_{\mathds{B}}}\mathbf{G}(\s,\y(\s)-2\mathfrak{k}[\mathds{B}]\mathfrak{l}(\mathrm{j}))\times|\mathds{B}|^{-1}{\textstyle\sum_{\mathrm{k}}}\mathfrak{A}^{(\mathrm{k})}(\s,\mathbf{U}^{\y(\s)+\cdot}). \label{eq:bg2124b}
\end{align}
Let $\mathds{A}^{\mathds{B},\pm}[\s,\y(\s)]$ denote the average over $\mathrm{k}\in\mathds{B}$ in $\mathrm{RHS}\eqref{eq:bg2124b}$. We now define the following ``cut-off version" of \eqref{eq:bg2124b}:
\begin{align}
&{\mathds{C}_{\mathbf{X}}^{\mathfrak{m}(\mathrm{j}),\pm}}[\mathds{R}^{\chi,\mathfrak{q},\pm,\mathrm{j}}\mathbf{Z};\s,\y(\s)] \nonumber\\
&{:=} \ |\mathscr{X}|^{-1}{\textstyle\sum_{\mathds{B}}}\mathbf{G}(\s,\y(\s)-2\mathfrak{k}[\mathds{B}]\mathfrak{l}(\mathrm{j}))\cdot\mathds{A}^{\mathds{B},\pm}[\s,\y(\s)]\mathbf{1}\{|\mathds{A}^{\mathds{B},\pm}[\s,\y(\s)]|\lesssim\N^{-\frac14-\frac{\gamma_{\mathrm{KL}}}{100}}\}. \label{eq:bg2124c}
\end{align}
For convenience, we set $\Phi[\s,\y(\s)]:=\eqref{eq:bg2124b}-\eqref{eq:bg2124c}$. By construction, we know $\Phi[\s,\y(\s)]$ is just $\mathrm{RHS}\eqref{eq:bg2124c}$, but $\lesssim$ therein turns into its opposite, namely $\gtrsim$. We now fix $\t\leq\t_{\mathrm{st}}\leq1$ and use this to compute as follows (with explanation given after):
\begin{align}
&|{\textstyle\int_{0}^{\t}}\mathbf{H}^{\N}(\s,\t(\N),\x)\{\N\Phi[\s,\cdot(\s)]\}\d\s| \nonumber\\
&\lesssim \ |\mathscr{X}|^{-1}{\textstyle\sum_{\mathds{B}}}\N^{\gamma_{\mathrm{ap}}}|{\textstyle\int_{0}^{\t}}\mathbf{H}^{\N}(\s,\t(\N),\x)\{\N\mathds{A}^{\mathds{B},\pm}[\s,\cdot(\s)]\mathbf{1}\{|\mathds{A}^{\mathds{B},\pm}[\s,\cdot(\s)]|\gtrsim\N^{-\frac14-\frac{\gamma_{\mathrm{KL}}}{100}}\}\}\d\s| \nonumber\\
&\lesssim \ |\mathscr{X}|^{-1}{\textstyle\sum_{\mathds{B}}}\N^{\beta_{\mathrm{BG}}}{\textstyle\int_{0}^{1}}|\mathbb{T}(\N)|^{-1}{\textstyle\sum_{\y}}\N|\mathds{A}^{\mathds{B},\pm}[\s,\y(\s)]|\mathbf{1}\{|\mathds{A}^{\mathds{B},\pm}[\s,\y(\s)]|\gtrsim\N^{-\frac14-\frac{\gamma_{\mathrm{KL}}}{100}}\}\d\s. \nonumber
\end{align}
The first line follows from our calculation of $\Phi[\s,\y(\s)]$ in the previous paragraph, as well as bounding $\mathbf{G}\lesssim\N^{\gamma_{\mathrm{ap}}}$ since we work before time $\t_{\mathrm{st}}$ (see Definition \ref{definition:method8}). The second line follows from {\eqref{eq:hke2} with $\mathrm{m}=0$}. More precisely, this estimates the heat kernel in $\mathbf{H}^{\N}(\s,\t(\N),\x)$ by $|\mathbb{T}(\N)|^{-1}|\t(\N)-\s|^{-1/2}\lesssim\N^{100\gamma_{\mathrm{reg}}}|\mathbb{T}(\N)|^{-1}$; see Definition \ref{definition:method5}. We then use $100\gamma_{\mathrm{reg}}+\gamma_{\mathrm{ap}}\leq\beta_{\mathrm{BG}}$; see Definitions \ref{definition:reg}, \ref{definition:method8}. Observe the previous estimates are deterministic and uniform in $\t\leq\t_{\mathrm{st}}$ and $\x\in\mathbb{T}(\N)$. Moreover, the final bound does not depend on these variables. Thus, the $\|\|_{\t_{\mathrm{st}};\mathbb{T}(\N)}$-norm of the LHS of the first line is $\lesssim$ the second line. Assuming that our choices of {\small$\mathfrak{A}^{(\mathrm{k})}$} (given by everything after the $\times$ symbol in \eqref{eq:bg2125b}) satisfy the constraints of Lemma \ref{lemma:bg2121}, we deduce
\begin{align}
&{\E\|{\textstyle\int_{0}^{\t}}\mathbf{H}^{\N}(\s,\t(\N),\x)\{\N\Phi[\s,\cdot(\s)]\}\d\s\|_{\t_{\mathrm{st}};\mathbb{T}_{N}}} \nonumber\\
&\lesssim \ {\textstyle\sup_{\mathds{B}}}\N^{\beta_{\mathrm{BG}}}\E{\textstyle\int_{0}^{1}}|\mathbb{T}(\N)|^{-1}{\textstyle\sum_{\y}}\N|\mathds{A}^{\mathds{B},\pm}[\s,\y(\s)]|\mathbf{1}\{|\mathds{A}^{\mathds{B},\pm}[\s,\y(\s)]|\gtrsim\N^{-\frac14-\frac{\gamma_{\mathrm{KL}}}{100}}\}\d\s \nonumber\\
&\lesssim \ \N^{\beta_{\mathrm{BG}}}\N\N^{-1-100\beta_{\mathrm{BG}}} \ \lesssim \ \N^{-99\beta_{\mathrm{BG}}}. \label{eq:bg2126}
\end{align}
Let us now verify that our choices of {\small$\mathfrak{A}^{(\mathrm{k})}$}, in fact, satisfy constraints of Lemma \ref{lemma:bg2121} with respect to $\mathbb{J}(\mathrm{k})$ given by the support of the $\mathds{R}^{\chi,\mathfrak{q},\pm,\mathrm{j}}(\s,\pm2\mathrm{k}\mathfrak{l}(\mathrm{j}))$-factor in \eqref{eq:bg2125b} (centered so that $\y(\s)=0$). First, $\mathbb{J}(\mathrm{k})$ are all contained in a common discrete interval with length $\lesssim|\mathds{B}|\mathfrak{l}(\mathrm{j})$. Indeed, the support of $\mathds{R}^{\chi,\mathfrak{q},\pm,\mathrm{j}}(\s,\pm2\mathrm{k}\mathfrak{l}(\mathrm{j}))$ has length $\lesssim\mathfrak{l}(\mathrm{j})$. Ranging over all $\mathrm{k}\in\mathds{B}$ gives supports that are shifts of this length-$\mathfrak{l}(\mathrm{j})$ interval by $\lesssim|\mathds{B}|\mathfrak{l}(\mathrm{j})$, so the claim follows. Next, observe the $\mathrm{CutExp}$-terms in \eqref{eq:bg2125b} (upon setting $\y(\s)=0$) have support in $\llbracket-2\mathfrak{k}[\mathds{B}]\mathfrak{l}(\mathrm{j}),\pm2\mathrm{k}\mathfrak{l}(\mathrm{j})\rrbracket$, which are all contained in some common discrete interval of length $\lesssim|\mathds{B}|\mathfrak{l}(\mathrm{j})$ by construction. This paragraph shows that $\mathbb{J}(\mathrm{k})$ and the support of {\small$\mathfrak{A}^{(\mathrm{k})}$} are all in a common discrete interval of length $\lesssim|\mathds{B}|\mathfrak{l}(\mathrm{j})$. The other constraints in Lemma \ref{lemma:bg2121} follow from Step 5 in the proof of Lemma \ref{lemma:finalprop5}. (Step 5 discusses the functional underlying {\small$\mathfrak{A}^{(\mathrm{k})}$}; whether we evaluate it at \eqref{eq:glsde} or \eqref{eq:glsdeloc} is irrelevant. As the functional {\small$\mathfrak{A}^{(\mathrm{k}),\pm}$} is the same as the functional {\small$\mathfrak{A}^{(\mathrm{k})}$} here, modulo differences in $\mathfrak{m}(\mathrm{j})$-scales that are irrelevant to the constraints of Lemma \ref{lemma:bg2121}, the reasoning therein applies.)
\subsubsection{Introducing a time-average}
Take $\s\leq\t_{\mathrm{st}}$. Observe \eqref{eq:bg2124c} is an average of terms of the form $\mathfrak{X}\mathbf{1}[|\mathfrak{X}|\lesssim\N^{-1/4-\gamma_{\mathrm{KL}}/100}]$ times $\mathbf{G}(\s,\cdot)$-terms that are $\lesssim\N^{\gamma_{\mathrm{ap}}}$ by Definition \ref{definition:method8}. Thus, we know $|\eqref{eq:bg2124c}|\lesssim\N^{-1/4-\gamma_{\mathrm{KL}}/100+\gamma_{\mathrm{ap}}}$. Thus, upon setting {\small$\phi_{\s,\y[\s]}$} equal to \eqref{eq:bg2124c}, we get the following in exactly the same fashion as \eqref{eq:bg2121}-\eqref{eq:bg2121b}, in which $\mathrm{LHS}\eqref{eq:bg2122I}$ uses $\tau=\tau(\mathrm{j})$. (Indeed, all we used about the choice of $\phi$ therein is the a priori estimate $|\phi|\lesssim\N^{-1/4-\gamma_{\mathrm{KL}}/100}$. In a nutshell, the following estimates the cost in replacing ${\mathds{C}_{\mathbf{X}}^{\mathfrak{m}(\mathrm{j}),\pm}}[\mathds{R}^{\chi,\mathfrak{q},\pm,\mathrm{j}}\mathbf{Z};\s,\y(\s)]$ by the average of ${\mathds{C}_{\mathbf{X}}^{\mathfrak{m}(\mathrm{j}),\pm}}[\mathds{R}^{\chi,\mathfrak{q},\pm,\mathrm{j}}\mathbf{Z};\s-\r,\y(\s-\r)]$ over $\r\in[0,\tau(\mathrm{j})]$ inside a time-integrated $\mathbf{H}^{\N}$-heat operator.)
\begin{align}
\mathrm{LHS}\eqref{eq:bg2122I} \ \lesssim \ \N^{\gamma_{\mathrm{ap}}+100\gamma_{\mathrm{reg}}+{\mathrm{c}}\gamma_{\mathrm{KL}}-\frac{\gamma_{\mathrm{KL}}}{100}} \ \lesssim \ \N^{-{\mathrm{C}}\beta_{\mathrm{BG}}}. \label{eq:bg2127}
\end{align}
We used that $\gamma_{\mathrm{reg}},\gamma_{\mathrm{ap}},\beta_{\mathrm{BG}}$ are small compared to $\gamma_{\mathrm{KL}}$ to derive the last bound. (There is a discrepancy of $\N^{\gamma_{\mathrm{ap}}}$ in \eqref{eq:bg2127} compared to \eqref{eq:bg2121b}, but this comes just from the fact that our estimate for $\phi$ is now worse by a factor of $\lesssim\N^{\gamma_{\mathrm{ap}}}$, and $\mathrm{LHS}\eqref{eq:bg2122I}$ scales linearly in $\phi$.) 

Let us recap what we have done. First, \eqref{eq:bg2126} replaces ${\mathds{A}_{\mathbf{X}}^{\mathfrak{m}(\mathrm{j}),\pm}}[\mathds{R}^{\chi,\mathfrak{q},\pm,\mathrm{j}}\mathbf{Z};\s,\y(\s)]$ in $\mathscr{A}^{\mathbf{X}}\mathscr{R}^{\chi,\mathfrak{q},\pm,\mathrm{j}}$ by ${\mathds{C}^{\mathfrak{m}(\mathrm{j}),\pm}_{\mathbf{X}}}[\mathds{R}^{\chi,\mathfrak{q},\pm,\mathrm{j}}\mathbf{Z};\s,\y(\s)]$ with error $\lesssim\N^{-99\beta_{\mathrm{BG}}}$. \eqref{eq:bg2127} further replaces ${\mathds{C}^{\mathfrak{m}(\mathrm{j}),\pm}_{\mathbf{X}}}[\mathds{R}^{\chi,\mathfrak{q},\pm,\mathrm{j}}\mathbf{Z};\s,\y(\s)]$ by its time-average. We finish this step by replacing the time-average of ${\mathds{C}^{\mathfrak{m}(\mathrm{j}),\pm}_{\mathbf{X}}}[\mathds{R}^{\chi,\mathfrak{q},\pm,\mathrm{j}}\mathbf{Z};\s,\y(\s)]$ with that (over the same set of times) of ${\mathds{A}^{\mathfrak{m}(\mathrm{j}),\pm}_{\mathbf{X}}}[\mathds{R}^{\chi,\mathfrak{q},\pm,\mathrm{j}}\mathbf{Z};\s,\y(\s)]$. To put this precisely, recall from right after \eqref{eq:bg2124c} that $\Phi[\s,\y(\s)]:=\eqref{eq:bg2124b}-\eqref{eq:bg2124c}$. We want to prove the estimate
\begin{align}
\E\|{\textstyle\int_{\tau(\mathrm{j})}^{\t}}\mathbf{H}^{\N}(\s,\t(\N),\x)\{\N\times\tau(\mathrm{j})^{-1}{\textstyle\int_{0}^{\tau(\mathrm{j})}}\Phi[\s-\r,\cdot(\s-\r)]\d\r\}\d\s\|_{\t_{\mathrm{st}};\mathbb{T}(\N)} \ \lesssim \ \N^{-99\beta_{\mathrm{BG}}}. \label{eq:bg2128}
\end{align}
To this end, by using the triangle inequality, we can move the $\d\r$-time-average outside of the norm and expectation. At this point, we can follow the proof of the display after \eqref{eq:bg2124c}. This shows that $\mathrm{LHS}\eqref{eq:bg2128}$ is big-Oh of the RHS of the first line of \eqref{eq:bg2126}, upon replacing the time-integration domain therein by $[\tau(\mathrm{j}),1]$ and introducing a shift $\s\mapsto\s-\r$ in the integrand (and then finally taking a supremum over $\r\in[0,\tau(\mathrm{j})]$ to account for the average in the previous sentence). But, this is bounded {from} above by the RHS of the first line of \eqref{eq:bg2126} itself; to see this, use a change-of-variables to remove the aforementioned shift $\s\mapsto\s-\r$ and extend the integration domain to all of $[0,1]$.  \eqref{eq:bg2128} now follows by applying the last bound in \eqref{eq:bg2126}. 
\subsubsection{Freezing the coupling constant $\lambda(\t)$}
The previous two steps estimate the cost in replacing $\mathscr{A}^{\mathbf{X}}\mathscr{R}^{\chi,\mathfrak{q},\pm,\mathrm{j}}$ (see Definition \ref{definition:bg29}) by something similar to $\mathscr{A}^{\mathbf{X},\mathbf{T}}\mathscr{R}^{\chi,\mathfrak{q},\pm,\mathrm{j}}$ (see Definition \ref{definition:bg211}), which is $\mathscr{A}^{\mathbf{X},\mathbf{T}}\mathscr{R}^{\chi,\mathfrak{q},\pm,\mathrm{j}}$, except $\mathds{A}^{\mathfrak{m}(\mathrm{j}),\tau(\mathrm{j}),\pm}[\mathds{R}^{\chi,\mathfrak{q},\pm,\mathrm{j}}\mathbf{Z};\s,\cdot(\s)]$ therein must be modified by replacing $\mathbf{G}^{\s}$ by $\mathbf{G}$ (or equivalently, by $\mathbf{G}^{\s-\r}$). We now bound the cost in ``undoing" this replacement, which would therefore finish {the} estimates for $\mathrm{LHS}\eqref{eq:bg212I}$. (We summarize how to combine our estimates so far to get \eqref{eq:bg212I} in the next step.) We start by recalling notation of Definition \ref{definition:bg211}. Now, for any $\mathrm{k}$, let us set $\mathds{R}^{\mathrm{k}}(\tau,\w):=\mathds{R}^{\chi,\mathfrak{q},\pm,\mathrm{j}}(\tau,\w\pm2\mathrm{k}\mathfrak{l}(\mathrm{j}))$ purely for convenience. Next, recall {the} notation introduced right before \eqref{eq:bg2124a}. For any $0\leq\r\leq\tau(\mathrm{j})$ (where $\tau(\mathrm{j})$ is from Definition \ref{definition:bg211}), we can write the following that we justify and explain afterwards, in which all sums but the first one {and the last one} are over $\mathds{B}\in\mathscr{X}$ and $\mathrm{k}\in\mathds{B}$:
\begin{align}
&\mathfrak{m}(\mathrm{j})^{-1}{\textstyle\sum_{\mathrm{k}=0}^{\mathfrak{m}(\mathrm{j})-1}}\mathds{R}^{\mathrm{k}}(\tau,\w)\mathbf{G}^{\s}(\tau,\w\pm2\mathrm{k}\mathfrak{l}(\mathrm{j})) \nonumber\\
&= \ |\mathscr{X}|^{-1}{\textstyle\sum_{\mathds{B}}}|\mathds{B}|^{-1}{\textstyle\sum_{\mathrm{k}}}\mathds{R}^{\mathrm{k}}(\tau,\w)\mathbf{G}^{\s}(\tau,\w\pm2\mathrm{k}\mathfrak{l}(\mathrm{j})) \label{eq:bg2129a}\\
&= \ |\mathscr{X}|^{-1}{\textstyle\sum_{\mathds{B}}}\mathbf{G}^{\s}(\tau,\w-2\mathfrak{l}(\mathrm{j})\sup\mathds{B})\times|\mathds{B}|^{-1}{\textstyle\sum_{\mathrm{k}}}\mathds{R}^{\mathrm{k}}(\tau,\w)\tfrac{\mathbf{G}^{\s}(\tau,\w\pm2\mathrm{k}\mathfrak{l}(\mathrm{j})) }{\mathbf{G}^{\s}(\tau,\w-2\mathfrak{l}(\mathrm{j})\sup\mathds{B})}\label{eq:bg2129b}\\
&= \ |\mathscr{X}|^{-1}{\textstyle\sum_{\mathds{B}}}\mathbf{G}^{\s}(\tau,\w-2\mathfrak{l}(\mathrm{j})\sup\mathds{B})\times|\mathds{B}|^{-1}{\textstyle\sum_{\mathrm{k}}}\mathds{R}^{\mathrm{k}}(\tau,\w)\exp[\lambda(\s)\N^{-\frac12}{\sum_{\z=-2\mathfrak{l}(\mathrm{j})\sup\mathds{B}+1}^{\pm2\mathrm{k}\mathfrak{l}(\mathrm{j})}}\mathbf{U}^{\tau,\w+\z}].\label{eq:bg2129c}
\end{align}
\eqref{eq:bg2129a} follows from the same reasoning as \eqref{eq:bg2124a}. (In fact, these are the same statement, except \eqref{eq:bg2124a} uses $\mathbf{G}$ and \eqref{eq:bg2129a} uses $\mathbf{G}^{\s}$.) \eqref{eq:bg2129b} is obvious once we note the first $\mathbf{G}^{\s}$-factor is independent of the $\mathrm{k}$-variable and can therefore be moved in or out of the $\mathrm{k}$-sum. \eqref{eq:bg2129c} follows by construction of $\mathbf{G}^{\s}$ in Definition \ref{definition:bg211} (as basically the Gartner transform but with frozen coupling constant $\lambda(\s)$). By the same token, we have the following (which is \eqref{eq:bg2129a}-\eqref{eq:bg2129c} but without freezing the coupling constant):
\begin{align}
&\mathfrak{m}(\mathrm{j})^{-1}{\textstyle\sum_{\mathrm{k}=0}^{\mathfrak{m}(\mathrm{j})-1}}\mathds{R}^{\mathrm{k}}(\tau,\w)\mathbf{G}(\tau,\w\pm2\mathrm{k}\mathfrak{l}(\mathrm{j})) \nonumber\\
&= \ |\mathscr{X}|^{-1}{\textstyle\sum_{\mathds{B}}}|\mathds{B}|^{-1}{\textstyle\sum_{\mathrm{k}}}\mathds{R}^{\mathrm{k}}(\tau,\w)\mathbf{G}(\tau,\w\pm2\mathrm{k}\mathfrak{l}(\mathrm{j})) \label{eq:bg2129d}\\
&= \ |\mathscr{X}|^{-1}{\textstyle\sum_{\mathds{B}}}\mathbf{G}(\tau,\w-2\mathfrak{l}(\mathrm{j})\sup\mathds{B})\times|\mathds{B}|^{-1}{\textstyle\sum_{\mathrm{k}}}\mathds{R}^{\mathrm{k}}(\tau,\w)\tfrac{\mathbf{G}(\tau,\w\pm2\mathrm{k}\mathfrak{l}(\mathrm{j})) }{\mathbf{G}(\tau,\w-2\mathfrak{l}(\mathrm{j})\sup\mathds{B})}\label{eq:bg2129e}\\
&= \ |\mathscr{X}|^{-1}{\textstyle\sum_{\mathds{B}}}\mathbf{G}(\tau,\w-2\mathfrak{l}(\mathrm{j})\sup\mathds{B})\times|\mathds{B}|^{-1}{\textstyle\sum_{\mathrm{k}}}\mathds{R}^{\mathrm{k}}(\tau,\w)\exp[\lambda(\tau)\N^{-\frac12}{\sum_{\z=-2\mathfrak{l}(\mathrm{j})\sup\mathds{B}+1}^{\pm2\mathrm{k}\mathfrak{l}(\mathrm{j})}}\mathbf{U}^{\tau,\w+\z}].\label{eq:bg2129f}
\end{align}
Now, {we present} some more useful notation. Given any $\mathds{B},\mathrm{k}$, let $\mathfrak{Y}^{\mathds{B},\mathrm{k}}(\tau,\w)$ be the inner-most summand in \eqref{eq:bg2129f}, and let $\mathfrak{Y}^{\mathds{B},\mathrm{k},\s}(\tau,\w)$ be the inner-most summand in \eqref{eq:bg2129c}. Let us also set $\mathbf{G}^{\s,\mathds{B}}(\tau,\w)$ as the $\mathbf{G}^{\s}$-factor in \eqref{eq:bg2129c}, and let $\mathbf{G}^{\mathds{B}}(\tau,\w)$ be the $\mathbf{G}$-factor in \eqref{eq:bg2129f}. A simple calculation shows that $\eqref{eq:bg2129a}-\eqref{eq:bg2129d}=\Gamma(\tau,\w;1)+\Gamma(\tau,\w;2)$, where
\begin{align}
\Gamma(\tau,\w;1) &:= \ |\mathscr{X}|^{-1}{\textstyle\sum_{\mathds{B}}}\mathbf{G}^{\s,\mathds{B}}(\tau,\w)\times|\mathds{B}|^{-1}{\textstyle\sum_{\mathrm{k}}}[\mathfrak{Y}^{\mathds{B},\mathrm{k},\s}(\tau,\w)-\mathfrak{Y}^{\mathds{B},\mathrm{k}}(\tau,\w)] \label{eq:bg2129g}\\
\Gamma(\tau,\w;2) \ &:= \ |\mathscr{X}|^{-1}{\textstyle\sum_{\mathds{B}}}[\mathbf{G}^{\s,\mathds{B}}(\tau,\w)-\mathbf{G}^{\mathds{B}}(\tau,\w)]\times|\mathds{B}|^{-1}{\textstyle\sum_{\mathrm{k}}}\mathfrak{Y}^{\mathds{B},\mathrm{k}}(\tau,\w). \label{eq:bg2129h}
\end{align}
For the rest of this step, assume $\tau=\s-\r$ for $0\leq\r\lesssim\tau(\mathrm{j})$, where $\tau(\mathrm{j})$ is from Definition \ref{definition:bg211}. Now, additionally assume that $\s\leq\t_{\mathrm{st}}$. In this case, we claim $|\mathbf{G}^{\s,\mathds{B}}(\tau,\w)-\mathbf{G}^{\mathds{B}}(\tau,\w)|\lesssim\N^{\gamma_{\mathrm{ap}}}\tau(\mathrm{j})$ with probability 1. Indeed, the LHS of this proposed bound is (the absolute value of) a difference of the exponentials $\mathrm{a}\mapsto\exp[\lambda(\s)\mathrm{a}]$ and $\mathrm{a}\mapsto\exp[\lambda(\tau)\mathrm{a}]$ evaluated at the same $\mathrm{a}$, which, by Definition \ref{definition:method8}, satisfies $|\mathrm{a}|\lesssim\log\log\N$. (It now suffices to use $|\lambda(\s)-\lambda(\tau)|\lesssim|\s-\tau|$ by smoothness of $\lambda(\cdot)$, which follows by Assumption \ref{ass:intro8}.) Thus, treating $\tau=\tau(\s)$ as a shift of $\s$, we deduce the following almost sure estimate for $\t\leq\t_{\mathrm{st}}$:
\begin{align}
&|{\textstyle\int_{\tau(\mathrm{j})}^{\t}}\mathbf{H}^{\N}(\s,\t(\N),\x)\{\N\Gamma(\tau,\cdot(\tau);2)\}\d\s| \nonumber\\
&\lesssim \ |\mathscr{X}|^{-1}{\textstyle\sum_{\mathds{B}}}\N^{\gamma_{\mathrm{ap}}}\tau(\mathrm{j}){\textstyle\int_{\tau(\mathrm{j})}^{\t}}\mathbf{H}^{\N}(\s,\t(\N),\x)\{\N||\mathds{B}|^{-1}{\textstyle\sum_{\mathrm{k}}}\mathfrak{Y}^{\mathds{B},\mathrm{k}}(\tau,\cdot(\tau))|\}\d\s. \label{eq:bg21210}
\end{align}
We now claim {$\E\|\mathrm{RHS}\eqref{eq:bg21210}\|_{\t_{\mathrm{st}};\mathbb{T}_{N}}\lesssim\N^{-99\beta_{\mathrm{BG}}}$}. To justify this, we use \eqref{eq:bg2126} to prove that the cost (in $\E\|\|$-norm) for introducing into the $\mathrm{k}\in\mathds{B}$-average in $\mathrm{RHS}\eqref{eq:bg21210}$ a cutoff of $\lesssim\N^{-1/4-\gamma_{\mathrm{KL}}/100}$ is $\lesssim\N^{-99\beta_{\mathrm{BG}}}$. (See the paragraph after \eqref{eq:bg2128} for why it does not matter that $\tau$ is not $\s$ but a time-shift of $\s$.) After introducing cutoff into $\mathrm{RHS}\eqref{eq:bg21210}$, we deduce via contractivity of the $\mathbf{H}^{\N}$ semigroup that $\mathrm{RHS}\eqref{eq:bg21210}\lesssim\N^{1+\gamma_{\mathrm{ap}}}\tau(\mathrm{j})\N^{-1/4-\gamma_{\mathrm{KL}}/100}$, which is $\lesssim\N^{-99\beta_{\mathrm{BG}}}$ because $\tau(\mathrm{j})\lesssim\N^{-3/4+{\mathrm{c}}\gamma_{\mathrm{KL}}}$ {for some small but fixed $c>0$} (see Definition \ref{definition:bg211}) and because $\gamma_{\mathrm{ap}},\beta_{\mathrm{BG}}$ are small compared to $\gamma_{\mathrm{KL}}$ (see Definition \ref{definition:method8}). This and \eqref{eq:bg21210} gives
\begin{align}
\E\|{\textstyle\int_{\tau(\mathrm{j})}^{\t}}\mathbf{H}^{\N}(\s,\t(\N),\x)\{\N\Gamma(\tau,\cdot(\tau);2)\}\d\s\|_{\t_{\mathrm{st}};\mathbb{T}(\N)} \ \lesssim \ \N^{-99\beta_{\mathrm{BG}}}. \label{eq:bg21210b}
\end{align}
Next, we claim the following, which is just the estimate \eqref{eq:bg21210b} but replacing $\Gamma(\cdot,\cdot;2)\mapsto\Gamma(\cdot,\cdot;1)$:
\begin{align}
\E\|{\textstyle\int_{\tau(\mathrm{j})}^{\t}}\mathbf{H}^{\N}(\s,\t(\N),\x)\{\N\Gamma(\tau,\cdot(\tau);1)\}\d\s\|_{\t_{\mathrm{st}};\mathbb{T}(\N)} \ \lesssim \ \N^{-99\beta_{\mathrm{BG}}}. \label{eq:bg21210c}
\end{align}
To show \eqref{eq:bg21210c}, we first note that in \eqref{eq:bg2129g}, we can remove the $\mathbf{G}^{\s,\mathds{B}}$-factor if we include a factor $\N^{\gamma_{\mathrm{ap}}}$ and if we put absolute values around the $\mathrm{k}\in\mathds{B}$-average. This is because before time $\t_{\mathrm{st}}$, we have $|\mathbf{G}^{\s,\mathds{B}}|\lesssim\N^{\gamma_{\mathrm{ap}}}$ by Definition \ref{definition:method8}. We now claim that $\N^{-\gamma_{\mathrm{ap}}}\tau(\mathrm{j})^{-1}[\mathfrak{Y}^{\mathds{B},\mathrm{k},\s}(\tau,\w)-\mathfrak{Y}^{\mathds{B},\mathrm{k}}(\tau,\w)]$ has the form of {\small$\mathfrak{A}^{(\mathrm{k})}(\tau,\mathbf{U}^{\tau,\w+\cdot})$}, where {\small$\mathfrak{A}^{(\mathrm{k})}$} is a function that satisfies the assumptions of Lemma \ref{lemma:bg2121}. Indeed, recall $\mathfrak{Y}^{\mathds{B},\mathrm{k},\s}(\tau,\w)-\mathfrak{Y}^{\mathds{B},\mathrm{k}}(\tau,\w)$ is the difference between the $\mathrm{k}\in\mathds{B}$-averages in \eqref{eq:bg2129c} and \eqref{eq:bg2129f}, respectively. As we explained after \eqref{eq:bg2126}, the $\mathrm{k}\in\mathds{B}$-average in \eqref{eq:bg2129c} satisfies {the} constraints of Lemma \ref{lemma:le2} (which are part of the constraints of Lemma \ref{lemma:bg2121}). This argument does not depend on the value of the coupling constant in \eqref{eq:bg2129c} to be $\lambda(\s)$. Thus the $\mathrm{k}\in\mathds{B}$-average in \eqref{eq:bg2129f} also satisfies {the} constraints of Lemma \ref{lemma:le2}. This implies the constraints in Lemma \ref{lemma:le2} are satisfied (since these constraints are defined by vanishing of linear functionals) by $\N^{-\gamma_{\mathrm{ap}}}\tau(\mathrm{j})^{-1}[\mathfrak{Y}^{\mathds{B},\mathrm{k},\s}(\tau,\w)-\mathfrak{Y}^{\mathds{B},\mathrm{k}}(\tau,\w)]$. The upper bound constraint in Lemma \ref{lemma:bg2121} follows since $|\mathds{R}^{\mathrm{k}}|\lesssim\N^{20\gamma_{\mathrm{reg}}}\mathfrak{l}(\mathrm{j})^{-3/2}$, and since $\mathfrak{Y}^{\mathds{B},\mathrm{k},\s}(\tau,\w)-\mathfrak{Y}^{\mathds{B},\mathrm{k}}(\tau,\w)$ equals $\mathds{R}^{\mathrm{k}}(\tau,\w)$ times something that, as explained before \eqref{eq:bg21210}, is $\lesssim\N^{\gamma_{\mathrm{ap}}}\tau(\mathrm{j})$. Thus, we can use Lemma \ref{lemma:bg2121} to show that after multiplying $\mathrm{LHS}\eqref{eq:bg21210c}$ by $\N^{\gamma_{\mathrm{ap}}}\tau(\mathrm{j})\ll1$, we can replace $\Gamma(\tau,\cdot(\tau);1)$ in \eqref{eq:bg21210c} by something that is $\lesssim\N^{-1/4-\gamma_{\mathrm{KL}}/100}$ with an error that is $\lesssim\N^{-99\beta_{\mathrm{BG}}}$. We are then left with $\N^{\gamma_{\mathrm{ap}}}\tau(\mathrm{j})$ times $\lesssim\N^{-1/4-\gamma_{\mathrm{KL}}/100}$, which, as explained before \eqref{eq:bg21210b}, produces $\lesssim\N^{-99\beta_{\mathrm{BG}}}$. Therefore, \eqref{eq:bg21210c} follows. Now, recall $\eqref{eq:bg2129a}-\eqref{eq:bg2129d}=\Gamma(\tau,\w;1)+\Gamma(\tau,\w;2)$ from right before \eqref{eq:bg2129e}. Combining this with \eqref{eq:bg21210b}-\eqref{eq:bg21210c}, linearity of integration and $\mathbf{H}^{\N}$ operators, and the triangle inequality, we ultimately deduce
\begin{align}
\E\|{\textstyle\int_{\tau(\mathrm{j})}^{\t}}\mathbf{H}^{\N}(\s,\t(\N),\x)\{\N\times\tau(\mathrm{j})^{-1}{\textstyle\int_{0}^{\tau(\mathrm{j})}}[\eqref{eq:bg2129a}-\eqref{eq:bg2129d}]_{\tau=\s-\r}\d\r\}\d\s\|_{\t_{\mathrm{st}};\mathbb{T}(\N)} \ \lesssim \ \N^{-99\beta_{\mathrm{BG}}}. \label{eq:bg21211}
\end{align}
%
\subsubsection{Putting it together}
The desired \eqref{eq:bg212I} asks to replace the spatial-average ${\mathds{A}^{\mathfrak{m}(\mathrm{j}),\pm}_{\mathbf{X}}}[\mathds{R}^{\chi,\mathfrak{q},\pm,\mathrm{j}}\mathbf{Z};\s,\cdot(\s)]$ with the space-time average $\mathds{A}^{\mathfrak{m}(\mathrm{j}),\tau(\mathrm{j}),\pm}[\mathds{R}^{\chi,\mathfrak{q},\pm,\mathrm{j}}\mathbf{Z};\s,\cdot(\s)]$ (see Definitions \ref{definition:bg29}, \ref{definition:bg211}), after integrating against a time-integrated $\mathbf{H}^{\N}$-operator. \eqref{eq:bg2126}, \eqref{eq:bg2127}, and \eqref{eq:bg2128} together show that the cost in replacing said average ${\mathds{A}^{\mathfrak{m}(\mathrm{j}),\pm}_{\mathbf{X}}}[\mathds{R}^{\chi,\mathfrak{q},\pm,\mathrm{j}}\mathbf{Z};\s,\cdot(\s)]$ by time-average on time-scale $\tau(\mathrm{j})$ is $\lesssim\N^{-99\beta_{\mathrm{BG}}}$ ({in the $\E\|\|_{\t_{\mathrm{st}};\mathbb{T}_{N}}$-norm}). \eqref{eq:bg21211} shows the cost in freezing the coupling constant (to make it independent of the time-averaging variable $\r$) is $\lesssim\N^{-99\beta_{\mathrm{BG}}}$ ({in the $\E\|\|_{\t_{\mathrm{st}};\mathbb{T}_{N}}$-norm}). Now note that freezing the coupling constant in said time-average gives $\mathds{A}^{\mathfrak{m}(\mathrm{j}),\tau(\mathrm{j}),\pm}[\mathds{R}^{\chi,\mathfrak{q},\pm,\mathrm{j}}\mathbf{Z};\s,\cdot(\s)]$. So, $\E\mathrm{LHS}\eqref{eq:bg212I}\lesssim\N^{-99\beta_{\mathrm{BG}}}$. By Markov inequality, $\mathrm{LHS}\eqref{eq:bg212I}\lesssim\N^{-98\beta_{\mathrm{BG}}}\lesssim\mathrm{RHS}\eqref{eq:bg212I}$ with high probability (where the last bound follows since {$\|\mathbf{Z}\|_{\t_{\mathrm{st}};\mathbb{T}_{N}}\gtrsim\N^{-\gamma_{\mathrm{ap}}}$} by Definition \ref{definition:method8}), so for any $\mathrm{j}$ and $\d=0$, the bound \eqref{eq:bg212I} holds with high probability. (Again, the proof for $\d=1,2,3,4$ is the same after replacing $\mathbf{H}^{\N}\mapsto(\mathscr{T}^{\pm,\mathrm{j}})^{\d}\mathbf{H}^{\N}$.) To extend this to \eqref{eq:bg212I} holding simultaneously for all $\d=0,1,2,3,4$ and $1\leq\mathrm{j}\leq\mathrm{j}(\infty)$ with high probability, it suffices to use $\mathrm{j}(\infty)\lesssim1$ (by Definition \ref{definition:bg24}, it is $\lesssim$ the number of steps of size $\gtrsim1$ needed to go from 0 to 1) and the fact that the intersection of $\mathrm{O}(1)$-many high probability events is high probability (by union bound for their complements).
\qed
\subsection{Proof of Lemma \ref{lemma:bg2121}}
Use Lemma \ref{lemma:le9} with the following choices (in the language of Lemma \ref{lemma:le8a}). Choose $\mathfrak{a}(\tau,\mathbf{U})=\mathfrak{A}^{\mathds{B}}(\tau,\mathbf{U})$ (so we can take $\mathbb{I}$ of length $\lesssim|\mathds{B}|\mathfrak{l}(\mathrm{j})$ by assumption). This gives the following (for $\kappa>0$ {to be} determined shortly):
\begin{align}
&\mathrm{LHS}\eqref{eq:bg2121I} \nonumber\\
&\lesssim \ \tfrac{1}{\kappa}\N^{-\frac94-\gamma_{\mathrm{KL}}}|\mathds{B}|^{3}\mathfrak{l}(\mathrm{j})^{3} + \tfrac{1}{\kappa}\sup_{\sigma\in\R}\sup_{0\leq\s\leq\t}\log\E^{\sigma,\s,\mathbb{I}}\exp[\kappa|\mathfrak{A}^{\mathds{B}}(\s,\mathbf{U})|\mathbf{1}\{|\mathfrak{A}^{\mathds{B}}(\s,\mathbf{U})|\gtrsim\N^{-\frac14-\frac{\gamma_{\mathrm{KL}}}{100}}\}]. \label{eq:bg2121I1a}
\end{align}
Before we proceed, we first use $\exp[0]=1$ to deduce that the second term in $\mathrm{RHS}\eqref{eq:bg2121I1a}$ equals
\begin{align}
\kappa^{-1}\textstyle{\sup_{\sigma,\s}}\log\{1+\E^{\sigma,\s,\mathbb{I}}\mathbf{1}\{|\mathfrak{A}^{\mathds{B}}(\s,\mathbf{U})|\gtrsim\N^{-\frac14-\frac{\gamma_{\mathrm{KL}}}{100}}\}\exp[\kappa|\mathfrak{A}^{\mathds{B}}(\s,\mathbf{U})|]\}. \label{eq:bg2121I1b}
\end{align}
{Let us now choose $\kappa=\N^{1/4}$; here is the motivation for this choice.} Recall {\small$\mathfrak{A}^{\mathds{B}}$} is the average of $|\mathds{B}|$-many {\small$\mathfrak{A}^{(\mathrm{k})}$}-terms that satisfy the upper bounds {\small$|\mathfrak{A}^{(\mathrm{k})}|\lesssim\N^{30\gamma_{\mathrm{reg}}}\mathfrak{l}(\mathrm{j})^{-3/2}$} and the constraints of Lemma \ref{lemma:le2}. Thus, Lemma \ref{lemma:le2} says {\small$\mathfrak{A}^{\mathds{B}}$} is sub-Gaussian with variance parameter $\lesssim\N^{60\gamma_{\mathrm{reg}}}\mathfrak{l}(\mathrm{j})^{-3}|\mathds{B}|^{-1}\lesssim\N^{-1/2-\delta_{\mathrm{KL}}+60\gamma_{\mathrm{reg}}}\lesssim\N^{-1/2-\gamma_{\mathrm{KL}}/30}$, where the last two bounds follow from {the assumptions} in the statement of {Lemma \ref{lemma:bg2121}} and the bounds $\delta_{\mathrm{KL}}\geq\gamma_{\mathrm{KL}}/20$ and $\gamma_{\mathrm{reg}}\leq{\mathrm{c}}\gamma_{\mathrm{KL}}$ {for some small $\mathrm{c}>0$}. Therefore, the probability of {\small$|\mathfrak{A}^{\mathds{B}}(\s,\mathbf{U})|\gtrsim\N^{-1/4-\gamma_{\mathrm{KL}}/100}$} is $\lesssim\N^{-{\mathrm{D}}}$ {for any large but fixed $\mathrm{D}>0$} by standard Gaussian concentration. Also, the product {\small$2\kappa|\mathfrak{A}^{\mathds{B}}|$} is sub-Gaussian of variance-parameter $\lesssim1$, so its exponential moment is $\lesssim1$. Using the previous two sentences with Cauchy-Schwarz shows
\begin{align}
&\E^{\sigma,\s,\mathbb{I}}\mathbf{1}\{|\mathfrak{A}^{\mathds{B}}(\s,\mathbf{U})|\gtrsim\N^{-\frac14-\frac{\gamma_{\mathrm{KL}}}{100}}\}\exp[\kappa|\mathfrak{A}^{\mathds{B}}(\s,\mathbf{U})|] \nonumber\\
&\lesssim \ (\E^{\sigma,\s,\mathbb{I}}\mathbf{1}\{|\mathfrak{A}^{\mathds{B}}(\s,\mathbf{U})|\gtrsim\N^{-\frac14-\frac{\gamma_{\mathrm{KL}}}{100}}\})^{\frac12}(\E^{\sigma,\s,\mathbb{I}}\exp[2\kappa|\mathfrak{A}^{\mathds{B}}(\s,\mathbf{U})|])^{\frac12} \nonumber\\
&\lesssim \ \N^{-99}. \label{eq:bg2121I2}
\end{align}
{Combining \eqref{eq:bg2121I2} with $\kappa\geq1$, \eqref{eq:bg2121I1b}, and $\log[1+\mathrm{a}]\leq\mathrm{a}$} (for all $\mathrm{a}\geq0$) then shows that the second term in $\mathrm{RHS}\eqref{eq:bg2121I1a}$ is $\lesssim\N^{-99}$. Thus, to complete the proof, it suffices to show that the first term in $\mathrm{RHS}\eqref{eq:bg2121I1a}$ is $\lesssim\N^{-1-{\mathrm{D}}\beta_{\mathrm{BG}}}$. For this, recall $|\mathds{B}|\mathfrak{l}(\mathrm{j})\lesssim\N^{1/2+\gamma_{\mathrm{KL}}/10}$ by assumption in {Lemma \ref{lemma:bg2121}}, and recall our choice $\kappa=\N^{1/4}$. Therefore, the first term in $\mathrm{RHS}\eqref{eq:bg2121I1a}$ is 
\begin{align}
\lesssim \N^{-\frac14}\N^{-\frac94-\gamma_{\mathrm{KL}}}\N^{\frac32+\frac{3}{10}\gamma_{\mathrm{KL}}} \ \lesssim \ \N^{-\frac52+\frac32-\frac{7}{10}\gamma_{\mathrm{KL}}} \ \lesssim \ \N^{-1-\frac{7}{10}\gamma_{\mathrm{KL}}},
\end{align}
which is $\lesssim\N^{-1-{\mathrm{D}}\beta_{\mathrm{BG}}}$ since $\beta_{\mathrm{BG}}$ is small compared to $\gamma_{\mathrm{KL}}$ (see Definition \ref{definition:method8}). This finishes the proof. \qed
\subsection{Proof of Lemma \ref{lemma:bg2122}}
If $\tau\geq\t$, then the desired bound follows by the contractivity of $\mathbf{H}^{\N}$ (see {\eqref{eq:hke3} with $\mathrm{m}=0$}). Thus, it suffices to assume $\tau\leq\t$. We first claim the following calculation holds (with explanation to be given afterwards):
\begin{align}
{\textstyle\int_{\tau}^{\t}}\mathbf{H}^{\N}(\s,\t(\N),\x)\{\tau^{-1}{\textstyle\int_{0}^{\tau}}\phi_{\s-\r,\cdot(\s-\r)}\d\r\}\d\s \ &= \ \tau^{-1}{\textstyle\int_{0}^{\tau}}\d\r{\textstyle\int_{\tau-\r}^{\t-\r}}\mathbf{H}^{\N}(\s+\r,\t(\N),\x)\{\phi_{\s,\cdot(\s)}\}\d\s. \label{eq:bg2122I1}
\end{align}
Indeed, by using linearity of integration and heat operators, we can move the $\d\r$-average outside {both the $\d\s$-integral and the $\mathbf{H}^{\N}$ operator}. Then, we change variables $\s\mapsto\s+\r$, which removes the time-shift in $\phi$ at the cost of a time-shift in the integration-domain and heat operator. We now decompose $\mathrm{RHS}\eqref{eq:bg2122I1}$ in the following fashion (with explanation given afterwards):
{
\begin{align}
&\mathrm{RHS}\eqref{eq:bg2122I1} \nonumber\\
&= \ \tau^{-1}{\textstyle\int_{0}^{\tau}}\d\r{\textstyle\int_{\tau-\r}^{\t-\r}}\mathbf{H}^{\N}(\s,\t(\N),\x)\{\phi_{\s,\cdot(\s)}\}\d\s\nonumber\\
&+ \ \tau^{-1}{\textstyle\int_{0}^{\tau}}\d\r{\textstyle\int_{\tau-\r}^{\t-\r}}[\mathbf{H}^{\N}(\s+\r,\t(\N),\x)-\mathbf{H}^{\N}(\s,\t(\N),\x)]\{\phi_{\s,\cdot(\s)}\}\d\s\nonumber\\
&= \ \tau^{-1}{\textstyle\int_{0}^{\tau}}\d\r{\textstyle\int_{0}^{\t}}\mathbf{H}^{\N}(\s,\t(\N),\x)\{\phi_{\s,\cdot(\s)}\}\d\s\nonumber\\
&+ \ \tau^{-1}{\textstyle\int_{0}^{\tau}}\d\r{\textstyle\int_{\tau-\r}^{\t-\r}}[\mathbf{H}^{\N}(\s+\r,\t(\N),\x)-\mathbf{H}^{\N}(\s,\t(\N),\x)]\{\phi_{\s,\cdot(\s)}\}\d\s \nonumber\\
&- \ \tau^{-1}{\textstyle\int_{0}^{\tau}}\d\r{\textstyle\int_{0}^{\tau-\r}}\mathbf{H}^{\N}(\s,\t(\N),\x)\{\phi_{\s,\cdot(\s)}\}\d\s-\tau^{-1}{\textstyle\int_{0}^{\tau}}\d\r{\textstyle\int_{\t-\r}^{\t}}\mathbf{H}^{\N}(\s,\t(\N),\x)\{\phi_{\s,\cdot(\s)}\}\d\s. \label{eq:bg2122I2}
\end{align}
}The first identity just removes the time-shift in $\mathbf{H}^{\N}$ in $\mathrm{RHS}\eqref{eq:bg2122I1}$ (with the appropriate cost). The second and third identities follow by applying the disjoint union $[0,\t]=[0,\tau-\r]\cup(\tau-\r,\t-\r]\cup(\t-\r,\t]$ to the first integral in the RHS of the first identity. Because the operator $\mathbf{H}^{\N}$ is contractive (see {\eqref{eq:hke3} with $\mathrm{m}=0$}), we deduce $\|\eqref{eq:bg2122I2}\|_{\t_{\mathrm{st}};\mathbb{T}(\N)}\lesssim\tau\|\phi\|_{\t_{\mathrm{st}};\mathbb{T}(\N)}$. Proposition \ref{prop:hke}, {namely \eqref{eq:hke4}}, also implies the operator $\mathbf{H}^{\N}(\s+\r,\t(\N),\cdot)-\mathbf{H}^{\N}(\s,\t(\N),\cdot)$ has operator norm (on $\mathscr{L}^{\infty}(\mathbb{T}(\N))$) $\lesssim[|\t(\N)-\s|^{-1}+|\t(\N)-\s-\r|^{-1}][\r+\N^{-1}]$. But $\s,\s+\r\leq\t$ and $\t(\N)-\t=\N^{-100\gamma_{\mathrm{reg}}}$ by Definition \ref{definition:method5}, so the operator norm is $\lesssim\N^{100\gamma_{\mathrm{reg}}}[\r+\N^{-1}]$. As $\r+\N^{-1}\lesssim\tau$, the last term in the second identity has $\|\|_{\t_{\mathrm{st}};\mathbb{T}(\N)}$ that is $\lesssim\N^{100\gamma_{\mathrm{reg}}}\tau\|\phi\|_{\t_{\mathrm{st}};\mathbb{T}(\N)}$. Now, observe {that}, up to a sign that is irrelevant because we take norms, $\mathrm{LHS}\eqref{eq:bg2122I}$ is just $\mathrm{RHS}\eqref{eq:bg2122I1}$ minus the first term in the second identity of the previous display, so we are done. \qed
%
%
%
\section{Proof of Corollary \ref{corollary:kpz}}
By Theorem \ref{theorem:kpz}, it suffices to prove $\t_{\mathrm{reg}}=1$ with high probability. Define the discrete time-set $\mathbb{X}:=[0,1]\cap\N^{-{\mathrm{D}}}\Z$ {for some large but finite $\mathrm{D}>0$}, and assume that with high probability, we have $|\mathbf{h}(\t,\x)-\mathbf{h}(\t,\y)|\lesssim\N^{\gamma_{\mathrm{reg}}/2}\N^{-1/2}|\x-\y|^{1/2}$ for $\t\in\mathbb{X}$ and $\x,\y\in\mathbb{T}(\N)$. Now, by Lemma \ref{lemma:ste}, this implies that with high probability, we can extend this upper bound from the very fine discretization $\mathbb{X}$ to all $\t\in[0,1]$. In particular, with high probability, we get $|\mathbf{h}(\t,\x)-\mathbf{h}(\t,\y)|\lesssim\N^{\gamma_{\mathrm{reg}}/2}\N^{-1/2}|\x-\y|^{1/2}+\N^{-99}\llangle\mathbf{U}\rrangle$ for $\t\in[0,1]$, where $\llangle\mathbf{U}\rrangle$ is the supremum of $|\mathbf{U}^{\s,\x}|$ over all $\s\in\mathbb{X}$ and $\x\in\mathbb{T}(\N)$. (The relevance of $\llangle\mathbf{U}\rrangle$ is just the fact that $\mathbf{h}(\t,\x)-\mathbf{h}(\t,\y)$ is a sum of $\N^{1/2}\mathbf{U}^{\t,\cdot}$ terms, and values of $\mathbf{U}^{\t,\cdot}$ for $\t\in[0,1]$ are controlled by its values for $\t\in\mathbb{X}$, as the SDE \eqref{eq:glsde} that it solves has uniformly Lipschitz coefficients; see Lemma \ref{lemma:ste}.) But $\llangle\mathbf{U}\rrangle$ is big-Oh of $\N^{1/2}$ times the max of $|\mathbf{h}(\s,\z)-\mathbf{h}(\s,\w)|$ over all $\s\in\mathbb{X}$ and $\z,\w\in\mathbb{T}(\N)$ because of the gradient relation from Definition \ref{definition:intro4}. So, we get $\llangle\mathbf{U}\rrangle\lesssim\N^{\gamma_{\mathrm{reg}}}\max_{\z,\w}|\z-\w|^{1/2}\lesssim\N^{2}$ (since we assumed a priori control on $|\mathbf{h}(\s,\z)-\mathbf{h}(\s,\w)|$ over all $\s\in\mathbb{X}$ and $\z,\w\in\mathbb{T}(\N)${)}. Thus, with high probability, we get that $|\mathbf{h}(\t,\x)-\mathbf{h}(\t,\y)|\lesssim\N^{\gamma_{\mathrm{reg}}/2}\N^{-1/2}|\x-\y|^{1/2}+\N^{-99}$, which is $\leq\N^{\gamma_{\mathrm{reg}}}\N^{-1/2}|\x-\y|^{1/2}$ if $\N\gtrsim1$ is large enough, since we only take $\x\neq\y$ in the formula of $\t_{\mathrm{reg}}$ (see Definition \ref{definition:reg}), and $|\x-\y|$ are integers. This gives $\t_{\mathrm{reg}}=1$ with high probability. 

It is left to show that with high probability, we get $|\mathbf{h}(\t,\x)-\mathbf{h}(\t,\y)|\lesssim\N^{\gamma_{\mathrm{reg}}/2}\N^{-1/2}|\x-\y|^{1/2}$ for $\t\in\mathbb{X}$ and $\x,\y\in\mathbb{T}(\N)$. Suppose $\mathbf{U}^{0,\cdot}\sim\mathbb{P}^{0,0,\mathbb{T}(\N)}$. As shown in Section 2 of \cite{DT} (combined with conservation of the average of $\mathbf{U}^{\t,\x}$ over $\x\in\mathbb{T}(\N)$), we know $\mathbb{P}^{0,0,\mathbb{T}(\N)}$ is an invariant measure since we assume $\mathscr{U}(\t,\cdot)=\mathscr{U}(0,\cdot)$. So $\mathbf{U}^{\t,\cdot}\sim\mathbb{P}^{0,0,\mathbb{T}(\N)}$ for all $\t\geq0$. Now, fix $\t\in\mathbb{X}$ and $\x,\y\in\mathbb{T}(\N)$. As $\mathbf{U}^{\t,\cdot}\sim\mathbb{P}^{0,0,\mathbb{T}(\N)}$, $\mathbf{h}(\t,\x)-\mathbf{h}(\t,\y)$ is a length-$|\x-\y|$ increment of a random walk bridge with zero average drift. The increments of said random walk bridge are $\N^{-1/2}\mathbf{U}^{\t,\z}$-terms, which are sub-Gaussian by the convexity in Assumption \ref{ass:intro8} (or the log-Sobolev inequality that it implies). Thus, by standard Gaussian random walk bridge concentration (see the end of the proof of Lemma \ref{lemma:bg27}), we get {that} $|\mathbf{h}(\t,\x)-\mathbf{h}(\t,\y)|\lesssim\N^{\gamma_{\mathrm{reg}}/2}\N^{-1/2}|\x-\y|^{1/2}$ \emph{fails} with probability $\lesssim\exp[-\Omega(1)\N^{\gamma_{\mathrm{reg}}}]$ for $\Omega(1)\gtrsim1$. As the size of $\mathbb{X}\times\mathbb{T}(\N)\times\mathbb{T}(\N)$ is $\lesssim\N^{\mathrm{O}(1)}$, {a} union bound says {that} $|\mathbf{h}(\t,\x)-\mathbf{h}(\t,\y)|\lesssim\N^{\gamma_{\mathrm{reg}}/2}\N^{-1/2}|\x-\y|^{1/2}$ \emph{fails} for some $(\t,\x,\y)\in\mathbb{X}\times\mathbb{T}(\N)\times\mathbb{T}(\N)$ with probability $\lesssim\exp[-\Omega(1)\N^{\gamma_{\mathrm{reg}}}]$ (for a possibly different $\Omega(1)\gtrsim1$). Denote this event (whose probability we have just now shown is $\lesssim\exp[-\Omega(1)\N^{\gamma_{\mathrm{reg}}}]$)  by $\mathcal{E}$. (Note $\mathcal{E}$ is a path-space event.) Let us also return the general setting of initial law for \eqref{eq:glsde} having density $\mathfrak{p}$ with respect to $\mathbb{P}^{0,0,\mathbb{T}(\N)}$. Let $\mathbb{P}^{\mathfrak{p}}$ be {the} probability with respect to the law of \eqref{eq:glsde} with initial measure $\mathfrak{p}\d\mathbb{P}^{0,0,\mathbb{T}(\N)}$, and let $\mathbb{P}^{\mathrm{stat}}$ be the same but replacing $\mathfrak{p}\mapsto1$. By the entropy inequality (see before the estimate prior to (5.28) in \cite{CYau}), we get the following estimate:
\begin{align}
\mathbb{P}^{\mathfrak{p}}[\mathcal{E}] \ \lesssim \ \{1+\mathfrak{D}[\mathfrak{p}]\}\times\{\log(1+\mathbb{P}^{\mathrm{stat}}[\mathcal{E}]^{-1})\}^{-1}, \label{eq:corollary1} 
\end{align}
where $\mathfrak{D}[\mathfrak{p}]$ is {the} relative entropy of the initial measure $\mathfrak{p}\d\mathbb{P}^{0,0,\mathbb{T}(\N)}$ with respect to the stationary one $\mathbb{P}^{0,0,\mathbb{T}(\N)}$. So, $\mathfrak{D}[\mathfrak{p}]\lesssim\N^{\alpha_{\mathrm{KL}}}$ by assumption. In the previous paragraph, we showed $\mathbb{P}^{\mathrm{stat}}[\mathcal{E}]\lesssim\exp[-\Omega(1)\N^{\gamma_{\mathrm{reg}}}]$. Since $\gamma_{\mathrm{reg}}$ is an arbitrarily large (but uniformly-bounded-in-$\N$) multiple of $\alpha_{\mathrm{KL}}$ (see Definition \ref{definition:entropydata}, \ref{definition:reg}, and Corollary \ref{corollary:kpz}), an elementary calculation then shows $\mathrm{RHS}\eqref{eq:corollary1}\lesssim\N^{-\beta}$ for $\beta\gtrsim1$, so we are done. \qed
\appendix
\section{Proofs of Theorem \ref{theorem:she}, Proposition \ref{prop:method11} modulo technical, elementary steps}
\subsection{Proof sketch of Theorem \ref{theorem:she}}
Fix any $\delta\in(0,1]$. Observe {that} the function $\mathrm{a}\mapsto\mathrm{a}\log\mathrm{a}$ is uniformly Lipschitz on $\mathrm{a}\in[\delta,\delta^{-1}]$ (with $\delta$-dependent Lipschitz norm). By {the assumptions} in Theorem \ref{theorem:she}, we know $\mathbf{Z}^{\infty,\mathrm{in}}(0,\cdot)$ is continuous and strictly positive. So, Theorem \ref{theorem:she} holds by standard theory for one-dimensional stochastic heat equations if we replace $[0,\infty)$ by $[0,\tau)$, where $\tau$ is a random but almost surely positive stopping time. (This requires that $\lambda(\t)\gtrsim1$ is uniformly smooth and the heat kernel estimates for $\mathbf{H}$ in Proposition \ref{prop:hkecont}. As for $\tau$, if $\tau(\delta)$ denotes the first time that $\mathbf{Z}^{\infty}$ either goes below $\delta$ or above $\delta^{-1}$, then $\tau$ is the limit of $\tau(\delta)$ as $\delta\to0$.) We are then tasked with showing that for $\mathfrak{t}\geq0$ and $\e>0$, there exists $\delta=\delta(\mathfrak{t},\e)\in(0,1]$ such that $\tau(\delta)\leq\mathfrak{t}$ has probability $\leq\e$. To this end, we proceed formally. (To make this rigorous, one just needs to mollify the noise in \eqref{eq:she2a}-\eqref{eq:she2b} in a reasonable way to get a legitimate {SPDE}.) Recall $\mathbf{h}^{\infty}=\lambda(\t)^{-1}\log\mathbf{Z}^{\infty}(\t,\x)$ from Definition \ref{definition:intro2}. Chain/Ito rule shows
\begin{align}
\d\mathbf{h}^{\infty}(\t,\x) \ = \ \bar{\alpha}(\t)\partial_{\x}^{2}\mathbf{h}^{\infty}(\t,\x)\d\t+\{\bar{\alpha}(\t;\wedge)|\partial_{\x}\mathbf{h}^{\infty}(\t,\x)|^{2}-\lambda(\t)"\infty"\}\d\t+\xi(\t,\x)\d\t. \label{eq:sheI1}
\end{align}
(This is just \eqref{eq:kpz} but with additional infinite renormalization $\lambda(\t)"\infty"$, which, again, is interpreted as something that diverges as we take the implicit mollification away.) Now, define $\mathbf{o}^{\infty}(\t,\x):=\lambda(\t)\mathbf{h}^{\infty}(\t,\x)$. Because $\lambda(\t)$ is smooth in $\t$, when we use Ito to compute $\d\mathbf{o}^{\infty}$, there are no cross-variations. It is therefore not too hard to see from \eqref{eq:sheI1} that
\begin{align}
\d\mathbf{o}^{\infty}(\t,\x) \ &= \ \bar{\alpha}(\t)\partial_{\x}^{2}\mathbf{o}^{\infty}(\t,\x)\d\t+\{\tfrac{\bar{\alpha}(\t;\wedge)}{\lambda(\t)}|\partial_{\x}\mathbf{o}^{\infty}(\t,\x)|^{2}-\lambda(\t)^{2}"\infty"\}\d\t\nonumber\\
&+ \ \lambda(\t)\xi(\t,\x)\d\t+\tfrac{\lambda'(\t)}{\lambda(\t)}\mathbf{o}^{\infty}(\t,\x)\d\t. \label{eq:sheI2}
\end{align}
Now, let $\mathbf{i}^{\infty}(\t,\x)$ solve \eqref{eq:sheI2} with $\mathbf{i}^{\infty}(0,\cdot)=\mathbf{o}^{\infty}(0,\cdot)$ but after replacing all $\mathbf{o}^{\infty}\mapsto\mathbf{i}^{\infty}$ and removing the last term on the RHS:
\begin{align}
\d\mathbf{i}^{\infty}(\t,\x) \ = \ \bar{\alpha}(\t)\partial_{\x}^{2}\mathbf{i}^{\infty}(\t,\x)\d\t+\{\tfrac{\bar{\alpha}(\t;\wedge)}{\lambda(\t)}|\partial_{\x}\mathbf{i}^{\infty}(\t,\x)|^{2}-\lambda(\t)^{2}"\infty"\}\d\t+\lambda(\t)\xi(\t,\x)\d\t. \label{eq:sheI3}
\end{align}
Set $\mathbf{d}^{\infty}:=\mathbf{o}^{\infty}-\mathbf{i}^{\infty}$, so $\mathbf{d}^{\infty}(0,\cdot)\equiv0$. An easy calculation shows {that the SPDE} for $\mathbf{d}^{\infty}$ is the following linearization of \eqref{eq:sheI2}:
\begin{align}
\d\mathbf{d}^{\infty}(\t,\x) \ &= \ \bar{\alpha}(\t)\partial_{\x}^{2}\mathbf{d}^{\infty}(\t,\x)\d\t + \tfrac{\bar{\alpha}(\t;\wedge)}{\lambda(\t)}[\partial_{\x}\mathbf{o}^{\infty}(\t,\x)+\partial_{\x}\mathbf{i}^{\infty}(\t,\x)]\partial_{\x}\mathbf{d}^{\infty}(\t,\x)\d\t \nonumber\\
&+ \  \tfrac{\lambda'(\t)}{\lambda(\t)}\mathbf{d}^{\infty}(\t,\x)\d\t+\tfrac{\lambda'(\t)}{\lambda(\t)}\mathbf{i}^{\infty}(\t,\x)\d\t. \nonumber
\end{align}
Thus, $\mathbf{d}^{\infty}$ solves a linear {parabolic} equation with zero-order term $\lambda'(\t)\lambda(\t)^{-1}\mathbf{i}^{\infty}(\t,\x)\d\t$. This, uniform smoothness, boundedness, positivity of $\lambda(\t)$, and $\mathbf{d}^{\infty}(0,\cdot)\equiv0$ then give a maximum principle bound (with continuous $\mathfrak{t}$-dependence in $\lesssim_{\mathfrak{t}}$):
\begin{align}
\sup_{0\leq\t\leq\mathfrak{t}}\sup_{\x\in\mathbb{T}}|\mathbf{d}^{\infty}(\t,\x)| \ \lesssim_{\mathfrak{t}} \ \sup_{0\leq\t\leq\mathfrak{t}}\sup_{\x\in\mathbb{T}}|\mathbf{i}^{\infty}(\t,\x)|. \label{eq:sheI4}
\end{align}
(Although the drift coefficient $\partial_{\x}\mathbf{o}^{\infty}(\t,\x)+\partial_{\x}\mathbf{i}^{\infty}(\t,\x)$ blows up as we remove mollification, we still get \eqref{eq:sheI4} after removing mollification, because the maximum principle does not depend on quantitative bounds on first-order coefficients.) Now, observe that $\tau(\delta)$ is at least the first time the LHS of \eqref{eq:sheI4} is $\gtrsim\log\delta^{-1}$ by calculus with the exponential. Thus, it suffices to show that for any $\mathfrak{t}\geq0$ and $\e>0$, the probability of $\mathrm{RHS}\eqref{eq:sheI4}\gtrsim\log\delta^{-1}$ is $\leq\e$ for some $\delta=\delta(\mathfrak{t},\e)$. We claim $\mathfrak{Z}:=\exp[\mathbf{i}^{\infty}]$ solves
\begin{align}
\d\mathfrak{Z}(\t,\x) \ = \ \bar{\alpha}(\t)\partial_{\x}^{2}\mathfrak{Z}(\t,\x)\d\t + \lambda(\t)\mathfrak{Z}(\t,\x)\xi(\t,\x)\d\t. \label{eq:sheI5}
\end{align}
\eqref{eq:sheI5} is just saying that KPZ can be linearized by Cole-Hopf. The more precise point is that the coupling constant one needs to multiply \eqref{eq:sheI3} by in order to linearize via Cole-Hopf is the ratio between the quadratic coefficient $\bar{\alpha}(\t;\wedge)\lambda(\t)^{-1}$ and $\bar{\alpha}(\t)$. But this is just 1! (The whole point of multiplying $\mathbf{h}^{\infty}$ by $\lambda(\t)$ to get $\mathbf{o}^{\infty}$ is to normalize the coupling constant to be 1. The price we must pay is the last term in \eqref{eq:sheI2}. This is a potential term that we ignore in \eqref{eq:sheI3} anyway and deal with afterwards to get \eqref{eq:sheI4}.) Now, global-in-time upper bounds for $\mathrm{RHS}\eqref{eq:sheI4}$ follow by standard comparison principle methods for \eqref{eq:sheI5}. For example, see \cite{Mu} (which does not care if $\bar{\alpha}(\t),\lambda(\t)$ are time-dependent; uniformly {bounded from above and from} below away from zero is enough). \qed
\subsection{Proof sketch of Proposition \ref{prop:method11}}
We first show the claim that $|\mathbf{W}(\t,\N\x)-\mathbf{Z}^{\infty}(\t,\x)|\to0$ uniformly in $(\t,\x)\in[0,1]\times\mathbb{T}$ in probability, where $\mathbf{W}(\t,\N\x)$ extends $\N^{-1}\mathbb{T}(\N)\subseteq\mathbb{T}$ via linear interpolation. First replace $\t(\N)\mapsto\t$ in \eqref{eq:method8IIa}-\eqref{eq:method8IIb}. (The cost in $\t(\N)\mapsto\t$ is $\mathrm{o}(1)$ as $\mathbf{H}^{\N}(\s,\t(\N),\x)\approx\mathbf{H}^{\N}(\s,\t,\x)$. Indeed, $\mathbf{H}^{\N}(\s,\t(\N),\x)$ is just a smoothing via $\mathbf{H}^{\N}(\t,\t(\N),\x)$, which is a short-time heat kernel and thus approximately a delta function, of $\mathbf{H}^{\N}(\s,\t,\x)$, which satisfies spatial regularity estimates that are uniform in $\N$ by Proposition \ref{prop:hke}.) Next, we again observe that $\mathrm{a}\mapsto\mathrm{a}\log\mathrm{a}$ is uniformly Lipschitz on $\mathrm{a}\in[\delta,\delta^{-1}]$ for any $\delta\in(0,1]$. So by standard methods for stochastic heat equations \cite{G}, regularity for $\mathbf{H}^{\N}$ in Proposition \ref{prop:hke}, and convergence of discrete-to-continuum kernels in Proposition \ref{prop:hkecont}, we get $|\mathbf{W}(\t,\N\x)-\mathbf{Z}^{\infty}(\t,\x)|\to0$ uniformly in $(\t,\x)\in[0,\tau(\delta)]\times\mathbb{T}$ in probability for any $\delta\in(0,1]$, where $\tau(\delta)$ is a stopping time given by the first time that $\mathbf{W}$ or $\mathbf{Z}^{\infty}$ either go below $\delta$ or above $\delta^{-1}$. It now suffices to show that for $\e>0$, there exists $\delta=\delta(\e)$ such that $\tau(\delta)\wedge1=1$ with probability $\geq1-\e$. To this end, let $\mathcal{E}[\mathbf{W};\delta]$ be the event that $\mathbf{W}$ either goes below $\delta$ or above $\delta^{-1}$ at time $\tau(\delta)\wedge1$, and let $\mathcal{E}[\mathbf{Z}^{\infty};\delta]$ be the same but for $\mathbf{Z}^{\infty}$ instead of $\mathbf{W}$. We know that for some $\mathcal{E}$ with probability $\mathrm{o}(1)$, we have $\mathcal{E}[\mathbf{W};\delta]\subseteq\mathcal{E}\cup\mathcal{E}[\mathbf{Z}^{\infty};\delta/2]$. (This comes from $|\mathbf{W}(\t,\N\x)-\mathbf{Z}^{\infty}(\t,\x)|\to0$ uniformly in $(\t,\x)\in[0,\tau(\delta)]\times\mathbb{T}$ in probability.) The claim now follows by noting that $\mathbb{P}(\mathcal{E}[\mathbf{Z}^{\infty};\delta])\to0$ as $\delta\to0$, which holds by Theorem \ref{theorem:she}. (In words, this is a bootstrapping argument. Locally in time, we know $\mathbf{W}(\t,\N\x)\approx\mathbf{Z}^{\infty}(\t,\x)$. So the upper and lower bounds from Theorem \ref{theorem:she} transfer to $\mathbf{W}$, which then lets us propagate $\mathbf{W}(\t,\N\x)\approx\mathbf{Z}^{\infty}(\t,\x)$. The point is that the bounds we get via Theorem \ref{theorem:she} are uniform throughout the propagation.)

We are left to estimate $\mathbf{Y}-\mathbf{W}$. This {reasoning} has a similar flavor. First, note the $\mathbf{Y}$ and $\mathbf{W}$ equation in Definition \ref{definition:method8} are the same, but the $\mathbf{Y}$ PDE has additional terms that are uniformly $\lesssim\N^{-\gamma_{\mathrm{reg}}/10}$ with probability 1 by definition. If the stochastic heat equation $\mathbf{W}$ solves had uniformly Lipschitz coefficients, standard practice for one-dimensional stochastic heat equations would bound $\mathbf{Y}-\mathbf{W}$. The issue is the $\mathbf{Y}\log\mathbf{Y}-\mathbf{W}\log\mathbf{W}$ difference of nonlinearities that we must control. By calculus, we deduce $|\mathbf{Y}\log\mathbf{Y}-\mathbf{W}\log\mathbf{W}|\lesssim\{|\mathbf{W}|+|\mathbf{W}|^{-1}+|\mathbf{Y}|+|\mathbf{Y}|^{-1}\}|\mathbf{W}-\mathbf{Y}|$. Thus, instead of the proposed estimate on $\mathbf{Y}-\mathbf{W}$, we would have $\|\mathbf{Y}-\mathbf{W}\|_{\mathfrak{t}}\lesssim\N^{-\gamma_{\mathrm{reg}}/90}\{1+\|\mathbf{Y}\|_{\mathfrak{t}}+\|\mathbf{Y}^{-1}\|_{\mathfrak{t}}+\|\mathbf{W}\|_{\mathfrak{t}}+\|\mathbf{W}^{-1}\|_{\mathfrak{t}}\}$, where $\|\|_{\mathfrak{t}}:=\|\|_{\mathfrak{t};\mathbb{T}(\N)}$ for any stopping time $\mathfrak{t}\in[0,1]$. It is left to get $\|\mathbf{Y}\|+\|\mathbf{Y}^{-1}\|\lesssim\N^{\gamma_{\mathrm{ap}}}[1+\|\mathbf{W}\|+\|\mathbf{W}^{-1}\|]$, where $\|\|$ is $\|\|_{\mathfrak{t}}$ for $\mathfrak{t}\equiv1$, with high probability. Let $\mathcal{F}$ be the event where $\|\mathbf{W}\|+\|\mathbf{W}^{-1}\|\lesssim\log\N$. By the first paragraph of this proof, $\mathcal{F}$ is high probability. Next, let $\tau\in[0,1]$ be 1 or the first time $\mathfrak{t}\in[0,1]$ that $\|\mathbf{Y}\|_{\mathfrak{t}}+\|\mathbf{Y}^{-1}\|_{\mathfrak{t}}\lesssim\N^{\gamma_{\mathrm{ap}}}[1+\|\mathbf{W}\|_{\mathfrak{t}}+\|\mathbf{W}^{-1}\|_{\mathfrak{t}}]$. It suffices to show $\tau=1$ with high probability, conditioning on $\mathcal{F}$ and $\|\mathbf{Y}-\mathbf{W}\|_{\tau}\lesssim\N^{-\gamma_{\mathrm{reg}}/90}\{1+\|\mathbf{Y}\|_{\tau}+\|\mathbf{Y}^{-1}\|_{\tau}+\|\mathbf{W}\|_{\tau}+\|\mathbf{W}^{-1}\|_{\tau}\}$ (since both of these events are high probability). Suppose $\tau\neq1$. Then we know $\|\mathbf{Y}\|_{\mathfrak{t}}+\|\mathbf{Y}^{-1}\|_{\mathfrak{t}}\lesssim\N^{\gamma_{\mathrm{ap}}}[1+\|\mathbf{W}\|_{\mathfrak{t}}+\|\mathbf{W}^{-1}\|_{\mathfrak{t}}]$ for all $\mathfrak{t}<\tau$. By the events we are conditioning on and calculus, we then know $\|\mathbf{Y}-\mathbf{W}\|_{\mathfrak{t}}=\mathrm{o}(1)$ and thus $\|\mathbf{Y}\|_{\mathfrak{t}}+\|\mathbf{Y}^{-1}\|_{\mathfrak{t}}\lesssim\log^{10}\N[1+\|\mathbf{W}\|_{\mathfrak{t}}+\|\mathbf{W}^{-1}\|_{\mathfrak{t}}]$. Since this is true for all $\mathfrak{t}<\tau$, by almost sure continuity of $\mathbf{Y},\mathbf{W}$, we get $\|\mathbf{Y}\|_{\tau}+\|\mathbf{Y}^{-1}\|_{\tau}\lesssim\log^{10}\N[1+\|\mathbf{W}\|_{\tau}+\|\mathbf{W}^{-1}\|_{\tau}]$. But we assume $\tau\neq1$, so $\|\mathbf{Y}\|_{\tau}+\|\mathbf{Y}^{-1}\|_{\tau}\gtrsim\N^{\gamma_{\mathrm{ap}}}[1+\|\mathbf{W}\|_{\tau}+\|\mathbf{W}^{-1}\|_{\tau}]$. (We did not claim any quantitative continuity of $\mathbf{Y}$ and $\mathbf{W}$; this is {a} qualitative statement.) This is a contradiction, so after conditioning on high probability events, we deduce that $\tau\neq1$ is impossible. In particular, $\tau=1$ with high probability. As noted earlier in this paragraph, this finishes the proof. \qed
%
%
%
\section{Proofs of Proposition \ref{prop:method2} and Lemma \ref{lemma:method4}: the $\mathbf{Z}$ SDE}\label{section:derivesde}
\subsection{Ito formula}
Throughout this section, {we} refer to Definition \ref{definition:intro6} for what $\mathbf{Z}$ is. Computing its SDE via Ito and \eqref{eq:hf}, we get:
\begin{lemma}\label{lemma:method21}
 Set $\mathbf{U}(\t):=\mathbf{U}^{\t,{\y(\x,\t)}}$ and $\mathbf{V}(\t):=\mathbf{U}^{\t,{\y(\x,\t)}+1}$ for convenience. Given any $\t\geq0$ and $\x\in\mathbb{T}(\N)$, we compute $\d\mathbf{Z}(\t,\x)$ to be equal to
{\small
\begin{align}
&-\tfrac12\lambda(\t)\N^{\frac32}\{\grad^{-}\mathscr{U}'(\t,\mathbf{V}(\t))\}\mathbf{Z}(\t,\x)\d\t + \tfrac12\lambda(\t)\N^{\frac32}\{\grad^{+}\mathscr{U}'(\t,\mathbf{U}(\t))\}\mathbf{Z}(\t,\x)\d\t + \delta(\t\in\mathbb{J})\grad^{-}\mathbf{Z}(\t,\x)\d\t \nonumber \\
&+ \lambda(\t)\N\mathscr{U}'(\t,\mathbf{U}(\t))\mathbf{Z}(\t,\x)\d\t + \lambda(\t)\N\mathscr{U}'(\t,\mathbf{V}(\t))\mathbf{Z}(\t,\x)\d\t + \lambda(\t)^{2}\N\mathbf{Z}(\t,\x)\d\t \nonumber \\
&- \lambda(\t)\mathscr{R}(\t)\mathbf{Z}(\t,\x)\d\t + \{\partial_{\t}\log|\lambda(\t)|\}\times\mathbf{Z}(\t,\x)\log\mathbf{Z}(\t,\x)\d\t + \sqrt{2}\lambda(\t)\N^{\frac12}\mathbf{Z}(\t,\x)\d\mathbf{b}(\t,{\y(\x,\t)}). \label{eq:method21I}
\end{align}
}\end{lemma}
\begin{proof}
We note that $\mathbf{Z}(\t,\x)$ has a jump if and only if $\t\in\mathbb{J}$. Indeed, $\mathbf{h}$ does (see Definition \ref{definition:intro6}), and $\mathbf{Z}(\t,\x)=\exp[\lambda(\t)\mathbf{h}(\t,\x)]$. Since said jump is in space to the left by 1, this gives the last term in the first line of \eqref{eq:method21I}. It now remains to compute $\d\mathbf{Z}(\t,\x)$ if $\t\not\in\mathbb{J}$. In this case, there is no jump, and we can use Ito on $\mathbf{Z}$ as a function of $\t$ and $\mathbf{h}$. We claim the following:
\begin{align}
&\d\mathbf{Z}(\t,\x) \ = \ \mathbf{Z}(\t,\x)\d\{\lambda(\t)\mathbf{h}(\t,\x)\}+\tfrac12\mathbf{Z}(\t,\x)\d[\lambda(\t)\mathbf{h}(\t,\x),\lambda(\t)\mathbf{h}(\t,\x)] \label{eq:method21I1a}\\
&= \ \tfrac{\d\lambda(\t)}{\lambda(\t)}\lambda(\t)\mathbf{h}(\t,\x)\mathbf{Z}(\t,\x)\d\t+\lambda(\t)\mathbf{Z}(\t,\x)\d\mathbf{h}(\t,\x)+\tfrac12\lambda(\t)^{2}\mathbf{Z}(\t,\x)\d[\mathbf{h}(\t,\x),\mathbf{h}(\t,\x)] \label{eq:method21I1b}\\
&= \ \tfrac{\d\lambda(\t)}{\lambda(\t)}\mathbf{Z}(\t,\x)\log\mathbf{Z}(\t,\x)\d\t+\lambda(\t)\mathbf{Z}(\t,\x)\d\mathbf{J}(\t,{\y(\x,\t)})-\lambda(\t)\mathscr{R}(\t)\mathbf{Z}(\t,\x)\d\t \label{eq:method21I1c}\\
&+ \ \tfrac12\lambda(\t)^{2}\mathbf{Z}(\t,\x)\d[\mathbf{J}(\t,{\y(\x,\t)}),\mathbf{J}(\t,{\y(\x,\t)})]. \nonumber
\end{align}
\eqref{eq:method21I1a} follows by Ito. \eqref{eq:method21I1b} follows by Leibniz rule for $\d\{\lambda(\t)\mathbf{h}(\t,\x)\}$ and $\d[\lambda(\t)\mathbf{h}(\t,\x),\lambda(\t)\mathbf{h}(\t,\x)]=\lambda(\t)^{2}\d[\mathbf{h}(\t,\x),\mathbf{h}(\t,\x)]$. (The last claim follows because $\lambda(\t)$ is smooth in $\t$, so there is no cross-variation.) \eqref{eq:method21I1c} follows from $\log\mathbf{Z}(\t,\x)=\lambda(\t)\mathbf{h}(\t,\x)$ and $\d\mathbf{h}(\t,\x)=\d\mathbf{J}(\t,{\y(\x,\t)})-\mathscr{R}(\t)\d\t$; see Definition \ref{definition:intro6}. (Again, since $\mathscr{R}(\t)$ is smooth in $\t$, the $[\mathbf{h},\mathbf{h}]$ bracket in \eqref{eq:method21I1b} turns into the $[\mathbf{J},\mathbf{J}]$ bracket in \eqref{eq:method21I1c}; there is no $\mathscr{R}(\t)$-cross variation to account for.) Now, \eqref{eq:hf} says $\lambda(\t)\mathbf{Z}(\t,\x)\d\mathbf{J}(\t,{\y(\x,\t)})$ is equal to
\begin{align}
&\lambda(\t)\N^{\frac32}\{\grad^{+}\mathscr{U}'(\t,\mathbf{U}(\t))\}\mathbf{Z}(\t,\x)\d\t + \lambda(\t)\N\{\mathscr{U}'(\t,\mathbf{U}(\t))+\mathscr{U}'(\t,\mathbf{V}(\t))\}\mathbf{Z}(\t,\x)\d\t\nonumber\\
&+\sqrt{2}\lambda(\t)\N^{\frac12}\mathbf{Z}(\t,\x)\d\mathbf{b}(\t,{\y(\x,\t)}). \nonumber
\end{align}
Observe {that} $\grad^{+}\mathscr{U}'(\t,\mathbf{U}(\t))=-\grad^{-}\mathscr{U}'(\t,\mathbf{V}(\t))$. This implies that the first term in the above display is equal to the first two terms on the RHS of the first line in \eqref{eq:method21I}. Next, observe the second term in the previous display is {equal to} the first two terms in the second line of \eqref{eq:method21I}. Third, the last term in the above display, plus the first and third terms in \eqref{eq:method21I1c}, equals the last line of \eqref{eq:method21I}. Therefore, to deduce \eqref{eq:method21I}, we are left to show {that} the last (bracket) term in \eqref{eq:method21I1c} is equal to the remaining last term in the second line of \eqref{eq:method21I}:
\begin{align}
\tfrac12\lambda(\t)^{2}\mathbf{Z}(\t,\x)\d[\mathbf{J}(\t,{\y(\x,\t)}),\mathbf{J}(\t,{\y(\x,\t)})] \ = \ \lambda(\t)^{2}\N\mathbf{Z}(\t,\x)\d\t. 
\end{align}
But this follows{,} since the martingale in the $\mathbf{J}$-SDE is $\sqrt{2}\N^{1/2}$ times a standard Brownian motion; see \eqref{eq:hf}. So we are done.
\end{proof}
\subsection{"Symmetric terms": the first line of \eqref{eq:method21I}}
Let us expand the RHS of the first line in \eqref{eq:method21I}. This comes from the symmetric part of the \eqref{eq:glsde} dynamics. At leading-order, we match it to a nonlinear second-order operator acting on $\mathbf{Z}$. We then linearize this nonlinear operator into $\N^{2}\bar{\alpha}(\t)\mathscr{T}(\t)$ (while collecting error terms); see Definition \ref{definition:method1}. \emph{In this calculation, what also falls out is a ``non-Gaussian degree-2 Hermite polynomial"} in \eqref{eq:method22Ic}. This will be grouped with the second line of \eqref{eq:method21I}.
\begin{lemma}\label{lemma:method22}
 Recall ${\mathfrak{a}}^{\pm}$ and $\mathfrak{b}$ in {Proposition \ref{prop:method2}}. Set $\mathbf{U}(\t)=\mathbf{U}^{\t,{\y(\x,\t)}}$ and $\mathbf{V}(\t)=\mathbf{U}^{\t,{\y(\x,\t)}+1}$. For any $\t\geq0$ and $\x\in\mathbb{T}(\N)$,
\begin{align}
&-\tfrac12\lambda(\t)\N^{\frac32}\{\grad^{-}\mathscr{U}'(\t,\mathbf{V}(\t))\}\mathbf{Z}(\t,\x)\d\t + \tfrac12\lambda(\t)\N^{\frac32}\{\grad^{+}\mathscr{U}'(\t,\mathbf{U}(\t))\}\mathbf{Z}(\t,\x)\d\t \label{eq:method22Ia}\\
&= \ \N^{2}\bar{\alpha}(\t)\Delta\mathbf{Z}(\t,\x)\d\t + \tfrac14\N\lambda(\t)^{2}\bar{\alpha}(\t)\Delta\mathbf{Z}(\t,\x)\d\t \label{eq:method22Ib}\\
&- \ \tfrac12\lambda(\t)^{2}\N\{\mathscr{U}'(\t,\mathbf{U}(\t))\mathbf{U}(\t)+\mathscr{U}'(\t,\mathbf{V}(\t))\mathbf{V}(\t)\}\mathbf{Z}(\t,\x)\d\t\label{eq:method22Ic} \\
&+ \ \tfrac12\lambda(\t)\N^{\frac32}\grad^{+}\{\mathscr{W}'(\t,\mathbf{U}(\t))\mathbf{Z}(\t,\x)\}\d\t - \tfrac12\lambda(\t)\N^{\frac32}\grad^{-}\{\mathscr{W}'(\t,\mathbf{V}(\t))\mathbf{Z}(\t,\x)\}\d\t \label{eq:method22Id}\\
&+ \ \tfrac14\lambda(\t)^{2}\N\grad^{+}\{[\bar{\alpha}(\t)\mathbf{U}(\t)^{2}-1]\mathbf{Z}(\t,\x)\}\d\t + \tfrac14\lambda(\t)^{2}\N\grad^{-}\{[\bar{\alpha}(\t)\mathbf{V}(\t)^{2}-1]\mathbf{Z}(\t,\x)\}\d\t \label{eq:method22Ie}\\
&+ \ \tfrac{1}{12}\lambda(\t)^{4}\E^{0,\t}[\mathscr{U}'(\t,\mathbf{u})\mathbf{u}^{3}]\mathbf{Z}(\t,\x)\d\t + \tfrac{1}{12}\lambda(\t)^{4}\{\mathscr{U}'(\t,\mathbf{V}(\t))\mathbf{V}(\t)^{3}-\E^{0,\t}[\mathscr{U}'(\t,\mathbf{u})\mathbf{u}^{3}]\}\mathbf{Z}(\t,\x)\d\t \label{eq:method22If}\\
&+ \ \N^{\frac12}\grad^{-}\{{\mathfrak{a}}^{-}(\t,{\y(\x,\t)})\mathbf{Z}(\t,\x)\}\d\t + \N^{\frac12}\grad^{+}\{{\mathfrak{a}}^{+}(\t,{\y(\x,\t)})\mathbf{Z}(\t,\x)\}\d\t + \N^{-\frac12}{\mathfrak{b}}(\t,{\y(\x,\t)})\mathbf{Z}(\t,\x)\d\t.\label{eq:method22Ig}
\end{align}
\end{lemma}
\begin{proof}
We start by computing discrete gradients of $\mathbf{Z}$ in terms of $\mathbf{U}(\t),\mathbf{V}(\t)$. By Taylor expansion (and Definition \ref{definition:intro6}), we have
\begin{align}
&\grad^{+}\mathbf{Z}(\t,\x) \ = \ \mathbf{Z}(\t,\x)\{\exp[\lambda(\t)\{\mathbf{J}(\t,{\y(\x,\t)}+1)-\mathbf{J}(\t,{\y(\x,\t)})\}]-1\} \label{eq:method22I1a}\\
&= \ \mathbf{Z}(\t,\x){\textstyle\sum_{\mathrm{j}=1}^{\infty}}\tfrac{\lambda(\t)^{\mathrm{j}}}{\mathrm{j}!}\{\mathbf{J}(\t,{\y(\x,\t)}+1)-\mathbf{J}(\t,{\y(\x,\t)})\}^{\mathrm{j}} \ = \ \mathbf{Z}(\t,\x){\textstyle\sum_{\mathrm{j}=1}^{\infty}}\tfrac{\lambda(\t)^{\mathrm{j}}}{\mathrm{j}!}\N^{-\frac12\mathrm{j}}\mathbf{V}(\t)^{\mathrm{j}}, \label{eq:method22I1b}
\end{align}
where the second identity in \eqref{eq:method22I1b} follows as $\mathbf{V}(\t)=\mathbf{U}^{\t,{\y(\x,\t)}+1}=\N^{1/2}\{\mathbf{J}(\t,{\y(\x,\t)}+1)-\mathbf{J}(\t,{\y(\x,\t)})\}$; see Definition \ref{definition:intro6}. (The calculation \eqref{eq:method22I1a}-\eqref{eq:method22I1b} is the usual one in \cite{BG}, for example.) Similarly, we have the following formula for negative gradients:
\begin{align}
\grad^{-}\mathbf{Z}(\t,\x) \ = \ \mathbf{Z}(\t,\x){\textstyle\sum_{\mathrm{j}=1}^{\infty}}\tfrac{(-1)^{\mathrm{j}}\lambda(\t)^{\mathrm{j}}}{\mathrm{j}!}\N^{-\frac12\mathrm{j}}\mathbf{U}(\t)^{\mathrm{j}}. \label{eq:method22I1c}
\end{align}
Note {that} $\Delta=-\grad^{+}\grad^{-}=-\grad^{-}\grad^{+}$; see Definition \ref{definition:intro4}. {(Indeed, $\Delta$ is a discrete Laplacian, and $\grad^{+}$ is a discrete gradient. Moreover, $\grad^{-}$ is a negative discrete gradient, because it has the opposite orientation.)} Therefore, by using \eqref{eq:method22I1a}-\eqref{eq:method22I1c}, we can compute $\Delta\mathbf{Z}$ in two ways:
\begin{align}
\Delta\mathbf{Z}(\t,\x) \ = \ -{\textstyle\sum_{\mathrm{j}=1}^{\infty}}\tfrac{\lambda(\t)^{\mathrm{j}}}{\mathrm{j}!}\N^{-\frac12\mathrm{j}}\grad^{-}\{\mathbf{Z}(\t,\x)\mathbf{V}(\t)^{\mathrm{j}}\} \ = \ -{\textstyle\sum_{\mathrm{j}=1}^{\infty}}\tfrac{(-1)^{\mathrm{j}}\lambda(\t)^{\mathrm{j}}}{\mathrm{j}!}\N^{-\frac12\mathrm{j}}\grad^{+}\{\mathbf{Z}(\t,\x)\mathbf{U}(\t)^{\mathrm{j}}\}. \label{eq:method22I1d}
\end{align}
{Now, we state a final preliminary.} For any $\phi,\psi:\mathbb{T}(\N)\to\R$ and $\mathfrak{l}>0$, we have the discrete Leibniz rule $\grad^{\mathfrak{l}}(\phi\psi)(\x)=\psi(\x)\grad^{\mathfrak{l}}\phi(\x)+\phi(\x+\mathfrak{l})\grad^{\mathfrak{l}}\psi(\x)$. (This differs from the continuum Leibniz rule in the spatial shift, which vanishes when $\mathfrak{l}\to0$ to recover continuum derivatives.) It will also be convenient to note that for any $\t\leq\t_{\mathrm{reg}}$ and $\x\in\mathbb{T}(\N)$, we have $|\mathbf{U}^{\t,\x}|\lesssim\N^{\gamma_{\mathrm{reg}}}$. This is deterministic. (See Remark \ref{remark:intro14}.) Since $|\mathscr{U}''|\lesssim1$ by Assumption \ref{ass:intro8}, this gives $|\mathscr{U}'(\t,\mathbf{U}^{\t,\x})[\mathbf{U}^{\t,\x}]^{\d}|\lesssim\N^{(\d+1)\gamma_{\mathrm{reg}}}$ {for $\t\leq\t_{\mathrm{reg}}$, $\x\in\mathbb{T}(\N)$, and $\d\geq0$}. (This bound is also deterministic.) We will take this paragraph for granted and use it possibly without reference.

We now inspect the first term in \eqref{eq:method22Ia}. By the discrete Leibniz rule and then \eqref{eq:method22I1c}, we first compute 
\begin{align}
&-\tfrac12\lambda(\t)\N^{\frac32}\{\grad^{-}\mathscr{U}'(\t,\mathbf{V}(\t))\}\mathbf{Z}(\t,\x) \nonumber\\
&= \ -\tfrac12\lambda(\t)\N^{\frac32}\grad^{-}\{\mathscr{U}'(\t,\mathbf{V}(\t))\mathbf{Z}(\t,\x)\}+\tfrac12\lambda(\t)\N^{\frac32}\mathscr{U}'(\t,\mathbf{U}(\t))\grad^{-}\mathbf{Z}(\t,\x) \label{eq:method22I2a}\\
&= \ -\tfrac12\lambda(\t)\N^{\frac32}\grad^{-}\{\mathscr{U}'(\t,\mathbf{V}(\t))\mathbf{Z}(\t,\x)\}-\tfrac12\lambda(\t)^{2}\N\mathscr{U}'(\t,\mathbf{U}(\t))\mathbf{U}(\t)\mathbf{Z}(\t,\x)\label{eq:method22I2b}\\
&+ \ \tfrac14\lambda(\t)^{3}\N^{\frac12}\mathscr{U}'(\t,\mathbf{U}(\t))\mathbf{U}(\t)^{2}\mathbf{Z}(\t,\x) - \tfrac{1}{12}\lambda(\t)^{4}\mathscr{U}'(\t,\mathbf{U}(\t))\mathbf{U}(\t)^{3}\mathbf{Z}(\t,\x)\label{eq:method22I2c}\\
&+ \ \tfrac12\lambda(\t)\N^{\frac32}\{{\textstyle\sum_{\mathrm{j}=4}^{\infty}}\tfrac{(-1)^{\mathrm{j}}\lambda(\t)^{\mathrm{j}}}{\mathrm{j}!}\N^{-\frac12\mathrm{j}}\mathscr{U}'(\t,\mathbf{U}(\t))\mathbf{U}(\t)^{\mathrm{j}}\}\mathbf{Z}(\t,\x). \label{eq:method22I2d}
\end{align}
(\eqref{eq:method22I2b}-\eqref{eq:method22I2d} {follows by extracting the $1\leq\mathrm{j}\leq3$ terms} in $\mathrm{RHS}\eqref{eq:method22I1c}$.) By \eqref{eq:method22I1a}-\eqref{eq:method22I1b} and the same token, {in order to evaluate the second term in \eqref{eq:method22Ia}, we observe that},
\begin{align}
&\tfrac12\lambda(\t)\N^{\frac32}\{\grad^{+}\mathscr{U}'(\t,\mathbf{U}(\t))\}\mathbf{Z}(\t,\x) \nonumber\\
&= \ \tfrac12\lambda(\t)\N^{\frac32}\grad^{+}\{\mathscr{U}'(\t,\mathbf{U}(\t))\mathbf{Z}(\t,\x)\}-\tfrac12\lambda(\t)\N^{\frac32}\mathscr{U}'(\t,\mathbf{V}(\t))\grad^{+}\mathbf{Z}(\t,\x) \label{eq:method22I2e}\\
&= \ \tfrac12\lambda(\t)\N^{\frac32}\grad^{+}\{\mathscr{U}'(\t,\mathbf{V}(\t))\mathbf{Z}(\t,\x)\}-\tfrac12\lambda(\t)^{2}\N\mathscr{U}'(\t,\mathbf{V}(\t))\mathbf{V}(\t)\mathbf{Z}(\t,\x)\label{eq:method22I2f}\\
&-\tfrac14\lambda(\t)^{3}\N^{\frac12}\mathscr{U}'(\t,\mathbf{V}(\t))\mathbf{V}(\t)^{2}\mathbf{Z}(\t,\x)-\tfrac{1}{12}\lambda(\t)^{4}\mathscr{U}'(\t,\mathbf{V}(\t))\mathbf{V}(\t)^{3}\mathbf{Z}(\t,\x)\label{eq:method22I2g}\\
&-\tfrac12\lambda(\t)\N^{\frac32}\{{\textstyle\sum_{\mathrm{j}=4}^{\infty}}\tfrac{\lambda(\t)^{\mathrm{j}}}{\mathrm{j}!}\N^{-\frac12\mathrm{j}}\mathscr{U}'(\t,\mathbf{V}(\t))\mathbf{V}(\t)^{\mathrm{j}}\}\mathbf{Z}(\t,\x). \label{eq:method22I2h}
\end{align}
From \eqref{eq:method22I2a}-\eqref{eq:method22I2h}, we deduce the following expansion for \eqref{eq:method22Ia}:
\begin{align}
&\tfrac{\eqref{eq:method22Ia}}{{\d\t}} \nonumber\\
&= \ -\tfrac12\lambda(\t)\N^{\frac32}\grad^{-}\{\mathscr{U}'(\t,\mathbf{V}(\t))\mathbf{Z}(\t,\x)\}+\tfrac12\lambda(\t)\N^{\frac32}\grad^{+}\{\mathscr{U}'(\t,\mathbf{U}(\t))\mathbf{Z}(\t,\x)\} \label{eq:method22I2i}\\
&- \ \tfrac12\lambda(\t)^{2}\N\mathscr{U}'(\t,\mathbf{U}(\t))\mathbf{U}(\t)\mathbf{Z}(\t,\x)-\tfrac12\lambda(\t)^{2}\N\mathscr{U}'(\t,\mathbf{V}(\t))\mathbf{V}(\t)\mathbf{Z}(\t,\x)\label{eq:method22I2j}\\
&+ \ \tfrac14\lambda(\t)^{3}\N^{\frac12}\mathscr{U}'(\t,\mathbf{U}(\t))\mathbf{U}(\t)^{2}\mathbf{Z}(\t,\x)-\tfrac14\lambda(\t)^{3}\N^{\frac12}\mathscr{U}'(\t,\mathbf{V}(\t))\mathbf{V}(\t)^{2}\mathbf{Z}(\t,\x) \label{eq:method22I2k}\\
&- \ \tfrac{1}{12}\lambda(\t)^{4}\mathscr{U}'(\t,\mathbf{U}(\t))\mathbf{U}(\t)^{3}\mathbf{Z}(\t,\x)-\tfrac{1}{12}\lambda(\t)^{4}\mathscr{U}'(\t,\mathbf{V}(\t))\mathbf{V}(\t)^{3}\mathbf{Z}(\t,\x)+\{\eqref{eq:method22I2d}+\eqref{eq:method22I2h}\}.\label{eq:method22I2l}
\end{align}
We now unfold each line in the previous display \eqref{eq:method22I2i}-\eqref{eq:method22I2l}. We first claim {that}
\begin{align}
&\mathrm{RHS}\eqref{eq:method22I2i}\nonumber\\
&= \ -\tfrac12\lambda(\t)\N^{\frac32}\grad^{-}\{\bar{\alpha}(\t)\mathbf{V}(\t)\mathbf{Z}(\t,\x)\}+\tfrac12\lambda(\t)\N^{\frac32}\grad^{+}\{\bar{\alpha}(\t)\mathbf{U}(\t)\mathbf{Z}(\t,\x)\}\label{eq:method22I3a}\\
&- \ \tfrac12\lambda(\t)\N^{\frac32}\grad^{-}\{\mathscr{W}'(\t,\mathbf{V}(\t))\mathbf{Z}(\t,\x)\}+\tfrac12\lambda(\t)\N^{\frac32}\grad^{+}\{\mathscr{W}'(\t,\mathbf{U}(\t))\mathbf{Z}(\t,\x)\}\label{eq:method22I3b}\\
&= \ \tfrac12\N^{2}\bar{\alpha}(\t)\Delta\mathbf{Z}(\t,\x)+\tfrac12\N^{2}\bar{\alpha}(\t)\{{\textstyle\sum_{\mathrm{j}=2}^{\infty}}\tfrac{\lambda(\t)^{\mathrm{j}}}{\mathrm{j}!}\N^{-\frac12\mathrm{j}}\grad^{-}[\mathbf{V}(\t)^{\mathrm{j}}\mathbf{Z}(\t,\x)]\}\label{eq:method22I3c}\\
&+ \ \tfrac12\N^{2}\bar{\alpha}(\t)\Delta\mathbf{Z}(\t,\x)+\tfrac12\N^{2}\bar{\alpha}(\t)\{{\textstyle\sum_{\mathrm{j}=2}^{\infty}}\tfrac{(-1)^{\mathrm{j}}\lambda(\t)^{\mathrm{j}}}{\mathrm{j}!}\N^{-\frac12\mathrm{j}}\grad^{+}[\mathbf{U}(\t)^{\mathrm{j}}\mathbf{Z}(\t,\x)]\}\label{eq:method22I3d}\\
&- \ \tfrac12\lambda(\t)\N^{\frac32}\grad^{-}\{\mathscr{W}'(\t,\mathbf{V}(\t))\mathbf{Z}(\t,\x)\}+\tfrac12\lambda(\t)\N^{\frac32}\grad^{+}\{\mathscr{W}'(\t,\mathbf{U}(\t))\mathbf{Z}(\t,\x)\} \label{eq:method22I3e}\\
&= \ \N^{2}\bar{\alpha}(\t)\Delta\mathbf{Z}(\t,\x)+\tfrac14\N\lambda(\t)^{2}\bar{\alpha}(\t)\grad^{-}\{\mathbf{V}(\t)^{2}\mathbf{Z}(\t,\x)\}+\tfrac14\N\lambda(\t)^{2}\bar{\alpha}(\t)\grad^{+}\{\mathbf{U}(\t)^{2}\mathbf{Z}(\t,\x)\} \label{eq:method22I3f}\\
&+ \ \tfrac12\N^{2}\bar{\alpha}(\t)\{{\textstyle\sum_{\mathrm{j}=3}^{\infty}}\tfrac{\lambda(\t)^{\mathrm{j}}}{\mathrm{j}!}\N^{-\frac12\mathrm{j}}\grad^{-}[\mathbf{V}(\t)^{\mathrm{j}}\mathbf{Z}(\t,\x)]\}+\tfrac12\N^{2}\bar{\alpha}(\t)\{{\textstyle\sum_{\mathrm{j}=3}^{\infty}}\tfrac{(-1)^{\mathrm{j}}\lambda(\t)^{\mathrm{j}}}{\mathrm{j}!}\N^{-\frac12\mathrm{j}}\grad^{+}[\mathbf{U}(\t)^{\mathrm{j}}\mathbf{Z}(\t,\x)]\}\label{eq:method22I3g}\\
&- \ \tfrac12\lambda(\t)\N^{\frac32}\grad^{-}\{\mathscr{W}'(\t,\mathbf{V}(\t))\mathbf{Z}(\t,\x)\}+\tfrac12\lambda(\t)\N^{\frac32}\grad^{+}\{\mathscr{W}'(\t,\mathbf{U}(\t))\mathbf{Z}(\t,\x)\} \label{eq:method22I3h}\\
&= \ \N^{2}\bar{\alpha}(\t)\Delta\mathbf{Z}(\t,\x)+\tfrac14\N\lambda(\t)^{2}\bar{\alpha}(\t)\grad^{-}\mathbf{Z}(\t,\x)+\tfrac14\N\lambda(\t)^{2}\bar{\alpha}(\t)\grad^{+}\mathbf{Z}(\t,\x) \label{eq:method22I3i}\\
&+ \ \tfrac14\N\lambda(\t)^{2}\grad^{-}\{[\bar{\alpha}(\t)\mathbf{V}(\t)^{2}-1]\mathbf{Z}(\t,\x)\}+\tfrac14\N\lambda(\t)^{2}\grad^{+}\{[\bar{\alpha}(\t)\mathbf{U}(\t)^{2}-1]\mathbf{Z}(\t,\x)\} \label{eq:method22I3j} \\
&+ \ \tfrac12\N^{2}\bar{\alpha}(\t)\{{\textstyle\sum_{\mathrm{j}=3}^{\infty}}\tfrac{\lambda(\t)^{\mathrm{j}}}{\mathrm{j}!}\N^{-\frac12\mathrm{j}}\grad^{-}[\mathbf{V}(\t)^{\mathrm{j}}\mathbf{Z}(\t,\x)]\}+\tfrac12\N^{2}\bar{\alpha}(\t)\{{\textstyle\sum_{\mathrm{j}=3}^{\infty}}\tfrac{(-1)^{\mathrm{j}}\lambda(\t)^{\mathrm{j}}}{\mathrm{j}!}\N^{-\frac12\mathrm{j}}\grad^{+}[\mathbf{U}(\t)^{\mathrm{j}}\mathbf{Z}(\t,\x)]\}\label{eq:method22I3k}\\
&- \ \tfrac12\lambda(\t)\N^{\frac32}\grad^{-}\{\mathscr{W}'(\t,\mathbf{V}(\t))\mathbf{Z}(\t,\x)\}+\tfrac12\lambda(\t)\N^{\frac32}\grad^{+}\{\mathscr{W}'(\t,\mathbf{U}(\t))\mathbf{Z}(\t,\x)\}. \label{eq:method22I3l}
\end{align}
\eqref{eq:method22I3a}-\eqref{eq:method22I3b} follows because $\mathscr{W}'(\t,\mathrm{a})=\mathscr{U}'(\t,\mathrm{a})-\bar{\alpha}(\t)\mathrm{a}$; see Proposition \ref{prop:method2}. \eqref{eq:method22I3c}-\eqref{eq:method22I3e} follows by a couple of steps. First, leave \eqref{eq:method22I3b} alone. Next, for terms on $\mathrm{RHS}\eqref{eq:method22I3a}$, we use \eqref{eq:method22I1c} and \eqref{eq:method22I1a}-\eqref{eq:method22I1b}, respectively. While doing so, we separate the $\mathrm{j}=1$ and $\mathrm{j}\geq2$ summands. (The $\mathrm{j}=1$ summand corresponds to $-\lambda(\t)\mathbf{V}(\t)\mathbf{Z}(\t,\x)$ and $\lambda(\t)\mathbf{U}(\t)\mathbf{Z}(\t,\x)$ on $\mathrm{RHS}\eqref{eq:method22I3a}$, while the $\mathrm{j}\geq2$-summands give the second terms in \eqref{eq:method22I3c} and \eqref{eq:method22I3d}, respectively.) \eqref{eq:method22I3f}-\eqref{eq:method22I3h} follows just from further separating $\mathrm{j}=2$ and $\mathrm{j}\geq3$ summands in \eqref{eq:method22I3c}-\eqref{eq:method22I3d}. \eqref{eq:method22I3i}-\eqref{eq:method22I3l} follows by writing $\bar{\alpha}(\t)\mathrm{a}^{2}=1+\{\bar{\alpha}(\t)\mathrm{a}^{2}-1\}$ for $\mathrm{a}=\mathbf{U}(\t),\mathbf{V}(\t)$ inside $\grad^{-}$ and $\grad^{+}$ in \eqref{eq:method22I3f}. We now inspect \eqref{eq:method22I2k}. For this term, we claim the following:
\begin{align}
\eqref{eq:method22I2k} \ &= \ \tfrac14\lambda(\t)^{3}\N^{\frac12}\{\mathscr{U}'(\t,\mathbf{U}(\t))\mathbf{U}(\t)^{2}-\mathscr{U}'(\t,\mathbf{V}(\t))\mathbf{V}(\t)^{2}\}\mathbf{Z}(\t,\x) \nonumber\\
&= \ -\tfrac14\lambda(\t)^{3}\N^{\frac12}\{\grad^{+}[\mathscr{U}'(\t,\mathbf{U}(\t))\mathbf{U}(\t)^{2}]\}\mathbf{Z}(\t,\x) \nonumber \\
&= \ -\tfrac14\lambda(\t)^{3}\N^{\frac12}\grad^{+}\{\mathscr{U}'(\t,\mathbf{U}(\t))\mathbf{U}(\t)^{2}\mathbf{Z}(\t,\x)\}+\tfrac14\lambda(\t)^{3}\N^{\frac12}\mathscr{U}'(\t,\mathbf{V}(\t))\mathbf{V}(\t)^{2}\grad^{+}\mathbf{Z}(\t,\x) \label{eq:method22I4a}\\
&= \ -\tfrac14\lambda(\t)^{3}\N^{\frac12}\grad^{+}\{\mathscr{U}'(\t,\mathbf{U}(\t))\mathbf{U}(\t)^{2}\mathbf{Z}(\t,\x)\}+\tfrac14\lambda(\t)^{4}\mathscr{U}'(\t,\mathbf{V}(\t))\mathbf{V}(\t)^{3}\mathbf{Z}(\t,\x) \label{eq:method22I4b}\\
&+ \ \tfrac14\lambda(\t)^{3}\N^{\frac12}\{{\textstyle\sum_{\mathrm{j}=2}^{\infty}}\tfrac{\lambda(\t)^{\mathrm{j}}}{\mathrm{j}!}\N^{-\frac12\mathrm{j}}\mathscr{U}'(\t,\mathbf{V}(\t))\mathbf{V}(\t)^{\mathrm{j}}\}\mathbf{Z}(\t,\x). \label{eq:method22I4c}
\end{align}
The first two lines follow by definition (and by $\mathbf{V}(\t)$ being {a} spatial shift of $\mathbf{U}(\t)$ by $+1$). \eqref{eq:method22I4a} follows from the discrete Leibniz rule. \eqref{eq:method22I4b}-\eqref{eq:method22I4c} follows by \eqref{eq:method22I1a}-\eqref{eq:method22I1b}. We now inspect \eqref{eq:method22I2l} (without touching the last two terms therein). We claim
\begin{align}
&\eqref{eq:method22I2l} \nonumber\\
&= \ -\tfrac16\lambda(\t)^{4}\mathscr{U}'(\t,\mathbf{V}(\t))\mathbf{V}(\t)^{3}\mathbf{Z}(\t,\x)+\tfrac{1}{12}\lambda(\t)^{4}\{\grad^{+}[\mathscr{U}'(\t,\mathbf{U}(\t))\mathbf{U}(\t)^{3}]\}\mathbf{Z}(\t,\x)+\{\eqref{eq:method22I2d}+\eqref{eq:method22I2h}\}\label{eq:method22I5b}\\
&= \ -\tfrac16\lambda(\t)^{4}\mathscr{U}'(\t,\mathbf{V}(\t))\mathbf{V}(\t)^{3}\mathbf{Z}(\t,\x)+\tfrac{1}{12}\lambda(\t)^{4}\grad^{+}\{\mathscr{U}'(\t,\mathbf{U}(\t))\mathbf{U}(\t)^{3}\mathbf{Z}(\t,\x)\} \label{eq:method22I5c}\\
&- \ \tfrac{1}{12}\lambda(\t)^{4}\mathscr{U}'(\t,\mathbf{V}(\t))\mathbf{V}(\t)^{3}\grad^{+}\mathbf{Z}(\t,\x)+\{\eqref{eq:method22I2d}+\eqref{eq:method22I2h}\}\label{eq:method22I5d}\\
&= \ -\tfrac16\lambda(\t)^{4}\mathscr{U}'(\t,\mathbf{V}(\t))\mathbf{V}(\t)^{3}\mathbf{Z}(\t,\x)+\tfrac{1}{12}\lambda(\t)^{4}\grad^{+}\{\mathscr{U}'(\t,\mathbf{U}(\t))\mathbf{U}(\t)^{3}\mathbf{Z}(\t,\x)\} \label{eq:method22I5e}\\
&- \ \tfrac{1}{12}\lambda(\t)^{4}\{{\textstyle\sum_{\mathrm{j}=1}^{\infty}}\tfrac{\lambda(\t)^{\mathrm{j}}}{\mathrm{j}!}\N^{-\frac12\mathrm{j}}\mathscr{U}'(\t,\mathbf{V}(\t))\mathbf{V}(\t)^{3+\frac12\mathrm{j}}\}\mathbf{Z}(\t,\x)+\{\eqref{eq:method22I2d}+\eqref{eq:method22I2h}\}. \label{eq:method22I5f}
\end{align}
\eqref{eq:method22I5b} follows{,} since $\mathbf{V}(\t)$ is a spatial shift of $\mathbf{U}(\t)$ by $+1$. \eqref{eq:method22I5c}-\eqref{eq:method22I5d} follows from {the} discrete Leibniz rule. \eqref{eq:method22I5e}-\eqref{eq:method22I5f} follows by \eqref{eq:method22I1a}-\eqref{eq:method22I1b}. We now use \eqref{eq:method22I2i}-\eqref{eq:method22I2l}, \eqref{eq:method22I3a}-\eqref{eq:method22I3l}, \eqref{eq:method22I4a}-\eqref{eq:method22I4c}, and \eqref{eq:method22I5b}-\eqref{eq:method22I5f}. All we do below is {to} copy the results of these calculations, except {for} two points. First, for \eqref{eq:method22I3i}, recall $\grad^{+}+\grad^{-}=\Delta$ from Definition \ref{definition:intro4}. Next, the first term in \eqref{eq:method22I5e} is the last term in \eqref{eq:method22I4b}, just with different coefficients; we combine these coefficients below. Ultimately, we get
\begin{align}
&\tfrac{\eqref{eq:method22Ia}}{{\d\t}} \nonumber\\
&= \ \N^{2}\bar{\alpha}(\t)\Delta\mathbf{Z}(\t,\x)+\tfrac14\N\lambda(\t)^{2}\bar{\alpha}(\t)\Delta\mathbf{Z}(\t,\x) \label{eq:method22I6a}\\
&+ \ \tfrac14\N\lambda(\t)^{2}\grad^{-}\{[\bar{\alpha}(\t)\mathbf{V}(\t)^{2}-1]\mathbf{Z}(\t,\x)\}+\tfrac14\N\lambda(\t)^{2}\grad^{+}\{[\bar{\alpha}(\t)\mathbf{U}(\t)^{2}-1]\mathbf{Z}(\t,\x)\} \label{eq:method22I6b} \\
&+ \ \tfrac12\N^{2}\bar{\alpha}(\t)\{{\textstyle\sum_{\mathrm{j}=3}^{\infty}}\tfrac{\lambda(\t)^{\mathrm{j}}}{\mathrm{j}!}\N^{-\frac12\mathrm{j}}\grad^{-}[\mathbf{V}(\t)^{\mathrm{j}}\mathbf{Z}(\t,\x)]\}+\tfrac12\N^{2}\bar{\alpha}(\t)\{{\textstyle\sum_{\mathrm{j}=3}^{\infty}}\tfrac{(-1)^{\mathrm{j}}\lambda(\t)^{\mathrm{j}}}{\mathrm{j}!}\N^{-\frac12\mathrm{j}}\grad^{+}[\mathbf{U}(\t)^{\mathrm{j}}\mathbf{Z}(\t,\x)]\}\label{eq:method22I6c}\\
&- \ \tfrac12\lambda(\t)\N^{\frac32}\grad^{-}\{\mathscr{W}'(\t,\mathbf{V}(\t))\mathbf{Z}(\t,\x)\}+\tfrac12\lambda(\t)\N^{\frac32}\grad^{+}\{\mathscr{W}'(\t,\mathbf{U}(\t))\mathbf{Z}(\t,\x)\} \label{eq:method22I6d}\\
&- \ \tfrac12\lambda(\t)^{2}\N\mathscr{U}'(\t,\mathbf{U}(\t))\mathbf{U}(\t)\mathbf{Z}(\t,\x)-\tfrac12\lambda(\t)^{2}\N\mathscr{U}'(\t,\mathbf{V}(\t))\mathbf{V}(\t)\mathbf{Z}(\t,\x)\label{eq:method22I6e}\\
&- \ \tfrac14\lambda(\t)^{3}\N^{\frac12}\grad^{+}\{\mathscr{U}'(\t,\mathbf{U}(\t))\mathbf{U}(\t)^{2}\mathbf{Z}(\t,\x)\}+\tfrac{1}{12}\lambda(\t)^{4}\mathscr{U}'(\t,\mathbf{V}(\t))\mathbf{V}(\t)^{3}\mathbf{Z}(\t,\x) \label{eq:method22I6f}\\
&+ \ \tfrac14\lambda(\t)^{3}\N^{\frac12}\{{\textstyle\sum_{\mathrm{j}=2}^{\infty}}\tfrac{\lambda(\t)^{\mathrm{j}}}{\mathrm{j}!}\N^{-\frac12\mathrm{j}}\mathscr{U}'(\t,\mathbf{V}(\t))\mathbf{V}(\t)^{\mathrm{j}}\}\mathbf{Z}(\t,\x)+\tfrac{1}{12}\lambda(\t)^{4}\grad^{+}\{\mathscr{U}'(\t,\mathbf{U}(\t))\mathbf{U}(\t)^{3}\mathbf{Z}(\t,\x)\} \label{eq:method22I6g}\\
&- \ \tfrac{1}{12}\lambda(\t)^{4}\{{\textstyle\sum_{\mathrm{j}=1}^{\infty}}\tfrac{\lambda(\t)^{\mathrm{j}}}{\mathrm{j}!}\N^{-\frac12\mathrm{j}}\mathscr{U}'(\t,\mathbf{V}(\t))\mathbf{V}(\t)^{3+\frac12\mathrm{j}}\}\mathbf{Z}(\t,\x)+\{\eqref{eq:method22I2d}+\eqref{eq:method22I2h}\}. \label{eq:method22I6h}
\end{align}
It suffices to show \eqref{eq:method22I6a}-\eqref{eq:method22I6h} is equal to \eqref{eq:method22Ib}-\eqref{eq:method22Ig} with appropriate choices of ${\mathfrak{a}}^{\pm}$ and ${\mathfrak{b}}$ (after multiplying by $\d\t$). The RHS of \eqref{eq:method22I6a} equals \eqref{eq:method22Ib}. \eqref{eq:method22I6b} equals \eqref{eq:method22Ie}. \eqref{eq:method22I6d} equals \eqref{eq:method22Id}. \eqref{eq:method22I6e} equals \eqref{eq:method22Ic}. The second term in \eqref{eq:method22I6f} equals \eqref{eq:method22If}. We clarify {that} all that is left in \eqref{eq:method22I6a}-\eqref{eq:method22I6h} are \eqref{eq:method22I6c}, the first term in \eqref{eq:method22I6f}, and \eqref{eq:method22I6g}-\eqref{eq:method22I6h}. All that is left in \eqref{eq:method22Ib}-\eqref{eq:method22Ig} is the last line \eqref{eq:method22Ig}. Therefore, we get the proposed identity $\eqref{eq:method22Ia}=\eqref{eq:method22Ib}+\ldots+\eqref{eq:method22Ig}$ if we choose
{\small
\begin{align}
&{\mathfrak{a}}^{-}(\t,{\y(\x,\t)}) \nonumber\\
&:= \ \tfrac12\N^{\frac32}\bar{\alpha}(\t){\textstyle\sum_{\mathrm{j}=3}^{\infty}}\tfrac{\lambda(\t)^{\mathrm{j}}}{\mathrm{j}!}\N^{-\frac12\mathrm{j}}\mathbf{V}(\t)^{\mathrm{j}} \label{eq:method22I7a}\\
&{\mathfrak{a}}^{+}(\t,{\y(\x,\t)}) \nonumber\\
 &:= \ \tfrac12\N^{\frac32}\bar{\alpha}(\t){\textstyle\sum_{\mathrm{j}=3}^{\infty}}\tfrac{(-1)^{\mathrm{j}}\lambda(\t)^{\mathrm{j}}}{\mathrm{j}!}\N^{-\frac12\mathrm{j}}\mathbf{U}(\t)^{\mathrm{j}} -\tfrac14\lambda(\t)^{3}\mathscr{U}'(\t,\mathbf{U}(\t))\mathbf{U}(\t)^{2}+\tfrac{1}{12}\lambda(\t)^{4}\mathscr{U}'(\t,\mathbf{U}(\t))\mathbf{U}(\t)^{3} \label{eq:method22I7b}\\
&{\mathfrak{b}}(\t,{\y(\x,\t)}) \nonumber\\
&:= \ \tfrac14\N\lambda(\t)^{3}{\textstyle\sum_{\mathrm{j}=2}^{\infty}}\tfrac{\lambda(\t)^{\mathrm{j}}}{\mathrm{j}!}\N^{-\frac12\mathrm{j}}\mathbf{V}(\t)^{\mathrm{j}}-\tfrac{1}{12}\N^{\frac12}\lambda(\t)^{4}{\textstyle\sum_{\mathrm{j}=1}^{\infty}}\tfrac{\lambda(\t)^{\mathrm{j}}}{\mathrm{j}}\N^{-\frac12\mathrm{j}}\mathscr{U}'(\t,\mathbf{V}(\t))\mathbf{V}(\t)^{3+\frac12\mathrm{j}} \label{eq:method22I7c}\\
&+ \ \tfrac12\N^{2}\lambda(\t){\textstyle\sum_{\mathrm{j}=4}^{\infty}}\tfrac{(-1)^{\mathrm{j}}\lambda(\t)^{\mathrm{j}}}{\mathrm{j}!}\N^{-\frac12\mathrm{j}}\mathscr{U}'(\t,\mathbf{U}(\t))\mathbf{U}(\t)^{\mathrm{j}}-\tfrac12\N^{2}\lambda(\t){\textstyle\sum_{\mathrm{j}=4}^{\infty}}\tfrac{\lambda(\t)^{\mathrm{j}}}{\mathrm{j}!}\N^{-\frac12\mathrm{j}}\mathscr{U}'(\t,\mathbf{V}(\t))\mathbf{V}(\t)^{\mathrm{j}}. \label{eq:method22I7d}
\end{align}
}It now remains to show ${\mathfrak{a}}^{\pm}$ and ${\mathfrak{b}}$ satisfy the estimates we claim in Proposition \ref{prop:method2}. Recall from after \eqref{eq:method22I1d} the bound $|\mathbf{V}(\t)|\lesssim\N^{\gamma_{\mathrm{reg}}}$ for $\t\leq\t_{\mathrm{reg}}$. Elementary geometric series bounds show {that} the infinite sum in \eqref{eq:method22I7a} is $\lesssim\N^{-3/2+3\gamma_{\mathrm{reg}}}$ for $\t_{\mathrm{reg}}$. (This requires $\gamma_{\mathrm{reg}}<1/2$; see Definitions \ref{definition:entropydata} and \ref{definition:reg}.) Therefore, $|{\mathfrak{a}}^{\pm}(\t,\cdot)|\lesssim\N^{10\gamma_{\mathrm{reg}}}$ if $\t\leq\t_{\mathrm{reg}}$. A similar argument also shows the same for ${\mathfrak{a}}^{+}$ and ${\mathfrak{b}}$. The (intuitive) point here is that for every geometric series, the factor $\N^{-\mathrm{j}/2}$ (with the minimal index $\mathrm{j}$) cancels the $\N$-dependent prefactor. As noted at the beginning of this paragraph, this finishes the proof.
\end{proof}
\subsection{"Asymmetric terms": second line of \eqref{eq:method21I}}
The asymmetric drift contribution in \eqref{eq:method21I} homogenizes into three (non-error) terms. {The first} and second are the diverging counter-term and characteristic shift in Definition \ref{definition:intro6}. The third is an a priori mysterious KPZ quadratic. But, by deep (and universal!) algebraic structure of Cole-Hopf, this quadratic is already accounted for in \eqref{eq:method22Ic}.
\begin{lemma}\label{lemma:method23}
 Recall ${\mathfrak{c}}$ from {Proposition \ref{prop:method2}}. Set $\mathbf{U}(\t):=\mathbf{U}^{\t,{\y(\x,\t)}}$ and $\mathbf{V}(\t):=\mathbf{U}^{\t,{\y(\x,\t)}+1}$. For any $\t\geq0$ and $\x\in\mathbb{T}(\N)$, we have 
{\small
\begin{align}
&\lambda(\t)\N\mathscr{U}'(\t,\mathbf{U}(\t))\mathbf{Z}(\t,\x)\d\t + \lambda(\t)\N\mathscr{U}'(\t,\mathbf{V}(\t))\mathbf{Z}(\t,\x)\d\t \label{eq:method23Ia}\\
&= \ \N^{\frac32}\bar{\alpha}(\t)\grad^{+}\mathbf{Z}(\t,\x)\d\t - \N^{\frac32}\bar{\alpha}(\t)\grad^{-}\mathbf{Z}(\t,\x)\d\t + \lambda(\t)\N\{\mathscr{U}'(\t,\mathbf{U}(\t))-\bar{\alpha}(\t)\mathbf{U}(\t)\}\mathbf{Z}(\t,\x)\d\t \label{eq:method23Ib}\\
&+ \ \lambda(\t)\N\{\mathscr{U}'(\t,\mathbf{V}(\t))-\bar{\alpha}(\t)\mathbf{V}(\t)\}\mathbf{Z}(\t,\x)\d\t + \tfrac16\lambda(\t)^{3}\bar{\alpha}(\t)\mathbf{V}(\t)^{3}\mathbf{Z}(\t,\x)\d\t \label{eq:method23Ic} \\
&+ \ \tfrac16\lambda(\t)^{3}\bar{\alpha}(\t)\grad^{+}\{\mathbf{U}(\t)^{3}\mathbf{Z}(\t,\x)\}\d\t + \N^{-\frac12}{\mathfrak{c}}(\t,{\y(\x,\t)})\mathbf{Z}(\t,\x)\d\t- \tfrac12\lambda(\t)^{2}\N^{\frac12}\bar{\alpha}(\t)\grad^{+}\{\mathbf{U}(\t)^{2}\mathbf{Z}(\t,\x)\}\d\t. \label{eq:method23Id}
\end{align}
}\end{lemma}
\begin{proof}
We first claim {that the following calculation holds}:
\begin{align}
&\tfrac{\eqref{eq:method23Ia}}{\d\t} \nonumber\\
&= \ \lambda(\t)\N\{\mathscr{U}'(\t,\mathbf{U}(\t))-\bar{\alpha}(\t)\mathbf{U}(\t)\}\mathbf{Z}(\t,\x)+\lambda(\t)\N\{\mathscr{U}'(\t,\mathbf{V}(\t))-\bar{\alpha}(\t)\mathbf{V}(\t)\}\mathbf{Z}(\t,\x) \label{eq:method23I1a}\\
&+ \ \N\lambda(\t)\bar{\alpha}(\t)\mathbf{U}(\t)\mathbf{Z}(\t,\x)+\N\lambda(\t)\bar{\alpha}(\t)\mathbf{V}(\t)\mathbf{Z}(\t,\x)\label{eq:method23I1b}\\
&= \ \lambda(\t)\N\{\mathscr{U}'(\t,\mathbf{U}(\t))-\bar{\alpha}(\t)\mathbf{U}(\t)\}\mathbf{Z}(\t,\x)+\lambda(\t)\N\{\mathscr{U}'(\t,\mathbf{V}(\t))-\bar{\alpha}(\t)\mathbf{V}(\t)\}\mathbf{Z}(\t,\x) \label{eq:method23I1c}\\
&- \ \N^{\frac32}\bar{\alpha}(\t)\grad^{-}\mathbf{Z}(\t,\x)+\N^{\frac32}\bar{\alpha}(\t)\grad^{+}\mathbf{Z}(\t,\x) + \N^{\frac32}\bar{\alpha}(\t)\{{\textstyle\sum_{\mathrm{j}=2}^{\infty}}\tfrac{\lambda(\t)^{\mathrm{j}}}{\mathrm{j}!}\N^{-\frac12\mathrm{j}}[(-1)^{\mathrm{j}}\mathbf{U}(\t)^{\mathrm{j}}-\mathbf{V}(\t)^{\mathrm{j}}]\}\mathbf{Z}(\t,\x). \label{eq:method23I1d}
\end{align}
\eqref{eq:method23I1a}-\eqref{eq:method23I1b} follows by adding and subtracting $\lambda(\t)\bar{\alpha}(\t)\mathbf{U}(\t)\mathbf{Z}(\t,\x)+\lambda(\t)\bar{\alpha}(\t)\mathbf{V}(\t)\mathbf{Z}(\t,\x)$ to \eqref{eq:method23Ia}. \eqref{eq:method23I1c}-\eqref{eq:method23I1d} follows by using \eqref{eq:method22I1a}-\eqref{eq:method22I1c} to \eqref{eq:method23I1b}. (While doing so, we separate {the $\mathrm{j}=1$ term} from $\mathrm{j}\geq2$ terms in the infinite sums in \eqref{eq:method22I1a}-\eqref{eq:method22I1c}.) Note {that} $\mathrm{RHS}\eqref{eq:method23I1c}$ plus the first two terms in \eqref{eq:method23I1d} gives \eqref{eq:method23Ib} plus the first term in \eqref{eq:method23Ic}. So, by the above display, it suffices to show
\begin{align}
&\N^{\frac32}\bar{\alpha}(\t)\{{\textstyle\sum_{\mathrm{j}=2}^{\infty}}\tfrac{\lambda(\t)^{\mathrm{j}}}{\mathrm{j}!}\N^{-\frac12\mathrm{j}}[(-1)^{\mathrm{j}}\mathbf{U}(\t)^{\mathrm{j}}-\mathbf{V}(\t)^{\mathrm{j}}]\}\mathbf{Z}(\t,\x) \nonumber\\
&= \ \tfrac16\lambda(\t)^{3}\bar{\alpha}(\t)\mathbf{V}(\t)^{3}\mathbf{Z}(\t,\x)- \tfrac12\lambda(\t)^{2}\N^{\frac12}\bar{\alpha}(\t)\grad^{+}\{\mathbf{U}(\t)^{2}\mathbf{Z}(\t,\x)\} \nonumber\\
&+ \ \tfrac16\lambda(\t)^{3}\bar{\alpha}(\t)\grad^{+}\{\mathbf{U}(\t)^{3}\mathbf{Z}(\t,\x)\} + \N^{-\frac12}{\mathfrak{c}}(\t,{\y(\x,\t)})\mathbf{Z}(\t,\x). \label{eq:method23I2}  
\end{align}
To this end, we extract from the infinite series in \eqref{eq:method23I2} the terms with indices $\mathrm{j}=2,3$. (Higher indices have high-enough powers of $\N^{-1}$ to beat $\N^{3/2}$.) We then leave all $\mathrm{j}\geq4$ grouped together. We claim that this ultimately gives
\begin{align}
&\mathrm{LHS}\eqref{eq:method23I2} \nonumber\\
&= \ \tfrac{\lambda(\t)^{2}}{2}\N^{\frac12}\bar{\alpha}(\t)\{\mathbf{U}(\t)^{2}-\mathbf{V}(\t)^{2}\}\mathbf{Z}(\t,\x)-\tfrac{\lambda(\t)^{3}}{6}\bar{\alpha}(\t)\{\mathbf{U}(\t)^{3}+\mathbf{V}(\t)^{3}\}\mathbf{Z}(\t,\x)\label{eq:method23I3a}\\
&+ \ \N^{\frac32}\bar{\alpha}(\t)\{{\textstyle\sum_{\mathrm{j}=4}^{\infty}}\tfrac{\lambda(\t)^{\mathrm{j}}}{\mathrm{j}!}\N^{-\frac12\mathrm{j}}[(-1)^{\mathrm{j}}\mathbf{U}(\t)^{\mathrm{j}}-\mathbf{V}(\t)^{\mathrm{j}}]\}\mathbf{Z}(\t,\x) \label{eq:method23I3b}\\
&= \ \tfrac{\lambda(\t)^{2}}{2}\N^{\frac12}\bar{\alpha}(\t)\{\mathbf{U}(\t)^{2}-\mathbf{V}(\t)^{2}\}\mathbf{Z}(\t,\x) - \tfrac{\lambda(\t)^{3}}{3}\bar{\alpha}(\t)\mathbf{V}(\t)^{3}\mathbf{Z}(\t,\x)-\tfrac{\lambda(\t)^{3}}{6}\bar{\alpha}(\t)\{\mathbf{U}(\t)^{3}-\mathbf{V}(\t)^{3}\}\mathbf{Z}(\t,\x)\label{eq:method23I3c}\\
&+ \ \N^{\frac32}\bar{\alpha}(\t)\{{\textstyle\sum_{\mathrm{j}=4}^{\infty}}\tfrac{\lambda(\t)^{\mathrm{j}}}{\mathrm{j}!}\N^{-\frac12\mathrm{j}}[(-1)^{\mathrm{j}}\mathbf{U}(\t)^{\mathrm{j}}-\mathbf{V}(\t)^{\mathrm{j}}]\}\mathbf{Z}(\t,\x) \label{eq:method23I3d}\\
&= \ -\tfrac{\lambda(\t)^{2}}{2}\N^{\frac12}\bar{\alpha}(\t)\{\grad^{+}\mathbf{U}(\t)^{2}\}\mathbf{Z}(\t,\x) - \tfrac{\lambda(\t)^{3}}{3}\bar{\alpha}(\t)\mathbf{V}(\t)^{3}\mathbf{Z}(\t,\x)+\tfrac{\lambda(\t)^{3}}{6}\bar{\alpha}(\t)\{\grad^{+}\mathbf{U}(\t)^{3}\}\mathbf{Z}(\t,\x)\label{eq:method23I3e}\\
&+ \ \N^{\frac32}\bar{\alpha}(\t)\{{\textstyle\sum_{\mathrm{j}=4}^{\infty}}\tfrac{\lambda(\t)^{\mathrm{j}}}{\mathrm{j}!}\N^{-\frac12\mathrm{j}}[(-1)^{\mathrm{j}}\mathbf{U}(\t)^{\mathrm{j}}-\mathbf{V}(\t)^{\mathrm{j}}]\}\mathbf{Z}(\t,\x)\label{eq:method23I3f} \\
&= \ -\tfrac{\lambda(\t)^{2}}{2}\N^{\frac12}\bar{\alpha}(\t)\grad^{+}\{\mathbf{U}(\t)^{2}\mathbf{Z}(\t,\x)\}+\tfrac{\lambda(\t)^{2}}{2}\N^{\frac12}\bar{\alpha}(\t)\mathbf{V}(\t)^{2}\grad^{+}\mathbf{Z}(\t,\x)-\tfrac{\lambda(\t)^{3}}{3}\bar{\alpha}(\t)\mathbf{V}(\t)^{3}\mathbf{Z}(\t,\x)\label{eq:method23I3g}\\
&+ \ \tfrac{\lambda(\t)^{3}}{6}\bar{\alpha}(\t)\grad^{+}\{\mathbf{U}(\t)^{3}\mathbf{Z}(\t,\x)\} - \tfrac{\lambda(\t)^{3}}{6}\bar{\alpha}(\t)\mathbf{V}(\t)^{3}\grad^{+}\mathbf{Z}(\t,\x)+\eqref{eq:method23I3f} \label{eq:method23I3h}\\
&= \ -\tfrac{\lambda(\t)^{2}}{2}\N^{\frac12}\bar{\alpha}(\t)\grad^{+}\{\mathbf{U}(\t)^{2}\mathbf{Z}(\t,\x)\}+\tfrac{\lambda(\t)^{3}}{2}\bar{\alpha}(\t)\mathbf{V}(\t)^{3}\mathbf{Z}(\t,\x)-\tfrac{\lambda(\t)^{3}}{3}\bar{\alpha}(\t)\mathbf{V}(\t)^{3}\mathbf{Z}(\t,\x) \label{eq:method23I3i}\\
&+ \ \tfrac{\lambda(\t)^{2}}{2}\N^{\frac12}\bar{\alpha}(\t)\mathbf{V}(\t)^{2}\{{\textstyle\sum_{\mathrm{j}=2}^{\infty}}\tfrac{\lambda(\t)^{\mathrm{j}}}{\mathrm{j}!}\N^{-\frac12\mathrm{j}}\mathbf{V}(\t)^{\mathrm{j}}\}\mathbf{Z}(\t,\x)+\tfrac{\lambda(\t)^{3}}{6}\bar{\alpha}(\t)\grad^{+}\{\mathbf{U}(\t)^{3}\mathbf{Z}(\t,\x)\}\label{eq:method23I3j}\\
&- \ \tfrac{\lambda(\t)^{3}}{6}\bar{\alpha}(\t)\mathbf{V}(\t)^{3}\{{\textstyle\sum_{\mathrm{j}=1}^{\infty}}\tfrac{\lambda(\t)^{\mathrm{j}}}{\mathrm{j}!}\N^{-\frac12\mathrm{j}}\mathbf{V}(\t)^{\mathrm{j}}\}\mathbf{Z}(\t,\x)+\eqref{eq:method23I3f}\label{eq:method23I3k}\\
&= \ -\tfrac{\lambda(\t)^{2}}{2}\N^{\frac12}\bar{\alpha}(\t)\grad^{+}\{\mathbf{U}(\t)^{2}\mathbf{Z}(\t,\x)\}+\tfrac{\lambda(\t)^{3}}{6}\bar{\alpha}(\t)\mathbf{V}(\t)^{3}\mathbf{Z}(\t,\x)+\eqref{eq:method23I3j}+\eqref{eq:method23I3k}. \label{eq:method23I3l}
\end{align}
\eqref{eq:method23I3a}-\eqref{eq:method23I3b} follows by extracting {the} $\mathrm{j}=2,3$ summands from the infinite series in \eqref{eq:method23I2}. \eqref{eq:method23I3c}-\eqref{eq:method23I3d} follows from, on the RHS of \eqref{eq:method23I3a}, writing $\mathbf{U}(\t)^{3}+\mathbf{V}(\t)^{3}=2\mathbf{V}(\t)^{3}+[\mathbf{U}(\t)^{3}-\mathbf{V}(\t)^{3}]$. \eqref{eq:method23I3e}-\eqref{eq:method23I3f} follows by $\mathbf{U}(\t)^{2}-\mathbf{V}(\t)^{2}=-\grad^{+}\mathbf{U}(\t)^{2}$ (since $\mathbf{V}(\t)$ shifts $\mathbf{U}(\t)$ in space by $+1$). \eqref{eq:method23I3g}-\eqref{eq:method23I3h} by the discrete Leibniz rule for the first and third terms in \eqref{eq:method23I3e}. (We refer to the paragraph following \eqref{eq:method22I1d} for the discrete Leibniz rule.) \eqref{eq:method23I3i}-\eqref{eq:method23I3k} follows from \eqref{eq:method22I1a}-\eqref{eq:method22I1b} applied to the second term in \eqref{eq:method23I3g} and {to the} second term in \eqref{eq:method23I3h}. (For the former application, we separate the $\mathrm{j}=1$ term from {the} $\mathrm{j}\geq2$ terms.) \eqref{eq:method23I3l} follows by combining the last two terms in \eqref{eq:method23I3i}. The desired identity \eqref{eq:method23I2} follows from $\mathrm{LHS}\eqref{eq:method23I2}=\eqref{eq:method23I3l}$, which is a consequence of the previous display, if we make the following choice for ${\mathfrak{c}}$:
\begin{align}
&{\mathfrak{c}}(\t,{\y(\x,\t)}) \nonumber\\
&= \ \tfrac{\lambda(\t)^{2}}{2}\N\bar{\alpha}(\t)\mathbf{V}(\t)^{2}\{{\textstyle\sum_{\mathrm{j}=2}^{\infty}}\tfrac{\lambda(\t)^{\mathrm{j}}}{\mathrm{j}!}\N^{-\frac12\mathrm{j}}\mathbf{V}(\t)^{\mathrm{j}}\}-\N^{\frac12}\tfrac{\lambda(\t)^{3}}{6}\bar{\alpha}(\t)\mathbf{V}(\t)^{3}\{{\textstyle\sum_{\mathrm{j}=1}^{\infty}}\tfrac{\lambda(\t)^{\mathrm{j}}}{\mathrm{j}!}\N^{-\frac12\mathrm{j}}\mathbf{V}(\t)^{\mathrm{j}}\} \\
&+ \ \N^{2}\bar{\alpha}(\t)\{{\textstyle\sum_{\mathrm{j}=4}^{\infty}}\tfrac{\lambda(\t)^{\mathrm{j}}}{\mathrm{j}!}\N^{-\frac12\mathrm{j}}[(-1)^{\mathrm{j}}\mathbf{U}(\t)^{\mathrm{j}}-\mathbf{V}(\t)^{\mathrm{j}}]\}.
\end{align}
The last paragraph in the proof of Lemma \ref{lemma:method22} shows $|{\mathfrak{c}}(\t,\cdot)|\lesssim\N^{10\gamma_{\mathrm{reg}}}$ deterministically for $\t\leq\t_{\mathrm{reg}}$, so we are done.
\end{proof}
\subsection{Proof of Proposition \ref{prop:method2}}
Combine \eqref{eq:method21I}, \eqref{eq:method22Ia}-\eqref{eq:method22Ig}, and \eqref{eq:method23Ia}-\eqref{eq:method23Id}. This gives
\begin{align}
&\d\mathbf{Z}(\t,\x) \nonumber\\
&= \ \N^{2}\bar{\alpha}(\t)\Delta\mathbf{Z}(\t,\x)\d\t + \tfrac14\N\lambda(\t)^{2}\bar{\alpha}(\t)\Delta\mathbf{Z}(\t,\x)\d\t+\delta(\t\in\mathbb{J})\grad^{-}\mathbf{Z}(\t,\x)\d\t \label{eq:method2I1a}\\
&- \ \tfrac12\lambda(\t)^{2}\N\{\mathscr{U}'(\t,\mathbf{U}(\t))\mathbf{U}(\t)+\mathscr{U}'(\t,\mathbf{V}(\t))\mathbf{V}(\t)\}\mathbf{Z}(\t,\x)\d\t\label{eq:method2I1b} \\
&+ \ \tfrac12\lambda(\t)\N^{\frac32}\grad^{+}\{\mathscr{W}'(\t,\mathbf{U}(\t))\mathbf{Z}(\t,\x)\}\d\t - \tfrac12\lambda(\t)\N^{\frac32}\grad^{-}\{\mathscr{W}'(\t,\mathbf{V}(\t))\mathbf{Z}(\t,\x)\}\d\t \label{eq:method2I1c}\\
&+ \ \tfrac14\lambda(\t)^{2}\N\grad^{+}\{[\bar{\alpha}(\t)\mathbf{U}(\t)^{2}-1]\mathbf{Z}(\t,\x)\}\d\t + \tfrac14\lambda(\t)^{2}\N\grad^{-}\{[\bar{\alpha}(\t)\mathbf{V}(\t)^{2}-1]\mathbf{Z}(\t,\x)\}\d\t \label{eq:method2I1d}\\
&+ \ \tfrac{1}{12}\lambda(\t)^{4}\E^{0,\t}[\mathscr{U}'(\t,\mathbf{u})\mathbf{u}^{3}]\mathbf{Z}(\t,\x)\d\t + \tfrac{1}{12}\lambda(\t)^{4}\{\mathscr{U}'(\t,\mathbf{V}(\t))\mathbf{V}(\t)^{3}-\E^{0,\t}[\mathscr{U}'(\t,\mathbf{u})\mathbf{u}^{3}]\}\mathbf{Z}(\t,\x)\d\t \label{eq:method2I1e}\\
&+ \ \N^{\frac12}\grad^{-}\{{\mathfrak{a}}^{-}(\t,{\y(\x,\t)})\mathbf{Z}(\t,\x)\}\d\t + \N^{\frac12}\grad^{+}\{{\mathfrak{a}}^{+}(\t,{\y(\x,\t)})\mathbf{Z}(\t,\x)\}\d\t + \N^{-\frac12}{\mathfrak{b}}(\t,{\y(\x,\t)})\mathbf{Z}(\t,\x)\d\t\label{eq:method2I1f}\\
&+ \ \N^{\frac32}\bar{\alpha}(\t)\grad^{+}\mathbf{Z}(\t,\x)\d\t - \N^{\frac32}\bar{\alpha}(\t)\grad^{-}\mathbf{Z}(\t,\x)\d\t+\lambda(\t)^{2}\N\mathbf{Z}(\t,\x)\d\t\label{eq:method2I1g}\\
&+ \ \lambda(\t)\N\{\mathscr{U}'(\t,\mathbf{U}(\t))-\bar{\alpha}(\t)\mathbf{U}(\t)\}\mathbf{Z}(\t,\x)\d\t+\lambda(\t)\N\{\mathscr{U}'(\t,\mathbf{V}(\t))-\bar{\alpha}(\t)\mathbf{V}(\t)\}\mathbf{Z}(\t,\x)\d\t \label{eq:method2I1h}\\
&+ \ \tfrac16\lambda(\t)^{3}\bar{\alpha}(\t)\mathbf{V}(\t)^{3}\mathbf{Z}(\t,\x)\d\t - \tfrac12\lambda(\t)^{2}\N^{\frac12}\bar{\alpha}(\t)\grad^{+}\{\mathbf{U}(\t)^{2}\mathbf{Z}(\t,\x)\}\d\t \label{eq:method2I1i} \\
&+ \ \tfrac16\lambda(\t)^{3}\bar{\alpha}(\t)\grad^{+}\{\mathbf{U}(\t)^{3}\mathbf{Z}(\t,\x)\}\d\t + \N^{-\frac12}{\mathfrak{c}}(\t,{\y(\x,\t)})\mathbf{Z}(\t,\x)\d\t \label{eq:method2I1j} \\
&- \ \lambda(\t)\mathscr{R}(\t)\mathbf{Z}(\t,\x)\d\t + \{\partial_{\t}\log|\lambda(\t)|\}\times\mathbf{Z}(\t,\x)\log\mathbf{Z}(\t,\x)\d\t + \sqrt{2}\lambda(\t)\N^{\frac12}\mathbf{Z}(\t,\x)\d\mathbf{b}(\t,{\y(\x,\t)}). \label{eq:method2I1k}
\end{align}
$\mathscr{T}(\t)\mathbf{Z}(\t,\x)\d\t$ on the RHS of \eqref{eq:method2I} equals $\mathrm{RHS}\eqref{eq:method2I1a}$ plus the first two terms in \eqref{eq:method2I1g}. The second plus third term on $\mathrm{RHS}\eqref{eq:method2I}$ equal the last two terms in \eqref{eq:method2I1k}. We now identify the rest of \eqref{eq:method2I1a}-\eqref{eq:method2I1k} with the remaining term ${\mathfrak{z}}(\t,{\y(\x,\t)})\mathbf{Z}(\t,\x)\d\t$ in \eqref{eq:method2I}. The RHS of \eqref{eq:method2IIb} plus \eqref{eq:method2IIc} equal \eqref{eq:method2I1b} plus \eqref{eq:method2I1h} and the last term in \eqref{eq:method2I1g}. \eqref{eq:method2IId} equals \eqref{eq:method2I1c}. \eqref{eq:method2IIe} equals \eqref{eq:method2I1d}. \eqref{eq:method2IIf} plus \eqref{eq:method2IIg} equals \eqref{eq:method2I1e} plus the first term in \eqref{eq:method2I1i} and the first term in \eqref{eq:method2I1k}. \eqref{eq:method2IIi} plus \eqref{eq:method2IIj} equals the remaining terms, namely \eqref{eq:method2I1f} plus the last term in \eqref{eq:method2I1i} and \eqref{eq:method2I1j}. (For this, we set ${\mathfrak{a}}^{+}$ in \eqref{eq:method2IIi} to be ${\mathfrak{a}}^{+}$ in \eqref{eq:method2I1f} plus $-2^{-1}\lambda(\t)^{2}\bar{\alpha}(\t)\mathbf{U}(\t)^{2}$ in \eqref{eq:method2I1i} plus $6^{-1}\N^{-1/2}\lambda(\t)^{3}\bar{\alpha}(\t)\mathbf{U}(\t)^{3}$ in \eqref{eq:method2I1j}. Lemma \ref{lemma:method22} and the reasoning in the last paragraph of its proof imply {that} this choice of ${\mathfrak{a}}^{+}$ satisfies the bound claimed in Proposition \ref{prop:method2}.) This finishes the proof. \qed
\subsection{Proof of Lemma \ref{lemma:method4}}
We start with the following preliminary calculation for any generic $\mathsf{F}$ with suitable growth at infinity:
\begin{align}
\E^{\sigma,\t}[\mathsf{F}(\mathbf{u})\mathscr{U}'(\t,\mathbf{u})] \ &= \ {\textstyle\int_{\R}}\mathsf{F}(\mathbf{u})\mathscr{U}'(\t,\mathbf{u})\exp\{\lambda(\sigma,\t)\mathbf{u}-\mathscr{U}(\t,\mathbf{u})+\mathscr{Z}(\sigma,\t)\}\d\mathbf{u} \label{eq:method41a} \\
&= \ \lambda(\sigma,\t){\textstyle\int_{\R}}\mathsf{F}(\mathbf{u})\exp\{\lambda(\sigma,\t)\mathbf{u}-\mathscr{U}(\t,\mathbf{u})+\mathscr{Z}(\sigma,\t)\}\d\mathbf{u} \label{eq:method41b}\\
&- {\textstyle\int_{\R}}\mathsf{F}(\mathbf{u})\left(\lambda(\sigma,\t)-\mathscr{U}'(\t,\mathbf{u})\right)\exp\{\lambda(\sigma,\t)\mathbf{u}-\mathscr{U}(\t,\mathbf{u})+\mathscr{Z}(\sigma,\t)\}\d\mathbf{u}\label{eq:method41c} \\
&= \ \lambda(\sigma,\t)\E^{\sigma,\t}\mathsf{F}(\mathbf{u}) - {\textstyle\int_{\R}}\mathsf{F}(\mathbf{u})\partial_{\mathbf{u}}\exp\{\lambda(\sigma,\t)\mathbf{u}-\mathscr{U}(\t,\mathbf{u})+\mathscr{Z}(\sigma,\t)\}\d\mathbf{u} \label{eq:method41d}\\
&= \  \lambda(\sigma,\t)\E^{\sigma,\t}\mathsf{F}(\mathbf{u}) + {\textstyle\int_{\R}}\partial_{\mathbf{u}}\mathsf{F}(\mathbf{u})\exp\{\lambda(\sigma,\t)\mathbf{u}-\mathscr{U}(\t,\mathbf{u})+\mathscr{Z}(\sigma,\t)\}\d\mathbf{u}  \label{eq:method41e}\\
&= \ \lambda(\sigma,\t)\E^{\sigma,\t}\mathsf{F}(\mathbf{u}) + \E^{\sigma,\t}\mathsf{F}'(\mathbf{u}).\label{eq:method41f}
\end{align}
\eqref{eq:method41a} is by definition. \eqref{eq:method41b}-\eqref{eq:method41c} is by adding and subtracting $\lambda(\sigma,\t)\E^{\sigma,\t}\mathsf{F}$. \eqref{eq:method41d} is calculus. \eqref{eq:method41e} is integration-by-parts. \eqref{eq:method41f} is definition. By \eqref{eq:method41a}-\eqref{eq:method41f} for $\mathsf{F}(\mathbf{u})=1$, we also get
\begin{align}
\lambda(\sigma,\t) \ = \ \lambda(\sigma,\t)\E^{\sigma,\t}1 \ = \ \E^{\sigma,\t}\mathscr{U}'(\t,\mathbf{u}) \quad{\mathrm{and}}\quad \partial_{\sigma}\lambda(\sigma,\t) \ = \ \partial_{\sigma}\E^{\sigma,\t}\mathscr{U}'(\t,\mathbf{u}). \label{eq:method41g}
\end{align}
Let us now show the first claim, that \eqref{eq:method2IIb} divided by $\mathbf{Z}(\t,\x)\d\t$ is in $\mathrm{QCT}$. In particular, we want to show
\begin{align}
\partial_{\sigma}^{\d}\E^{\sigma,\t}\{\mathscr{U}'(\t,\mathbf{u})-\bar{\alpha}(\t)\mathbf{u}-\tfrac12\lambda(\t)[\mathscr{U}'(\t,\mathbf{u})\mathbf{u}-1]\}|_{\sigma=0} \ = \ 0 \quad\mathrm{for}\quad \d=0,1,2. \label{eq:method42}
\end{align}
We first note {that} $\E^{\sigma,\t}\mathbf{u}=\sigma$ by definition of $\E^{\sigma,\t}$. We use this {together} with \eqref{eq:method41g} and \eqref{eq:method41a}-\eqref{eq:method41f} for $\mathsf{F}(\mathbf{u})=\mathbf{u}$ to get
\begin{align}
&\E^{\sigma,\t}\{\mathscr{U}'(\t,\mathbf{u})-\bar{\alpha}(\t)\mathbf{u}-\tfrac12\lambda(\t)[\mathscr{U}'(\t,\mathbf{u})\mathbf{u}-1]\} \nonumber\\
&= \ \lambda(\sigma,\t) - \bar{\alpha}(\t)\sigma - \tfrac12\lambda(\t)\lambda(\sigma,\t)\E^{\sigma,\t}\mathbf{u} - \tfrac12\lambda(\t)\E^{\sigma,\t}1 + \tfrac12\lambda(\t) \nonumber \\
&= \ \lambda(\sigma,\t)-\bar{\alpha}(\t)\sigma - \tfrac12\lambda(\t)\lambda(\sigma,\t)\sigma. \nonumber
\end{align}
The last line vanishes for $\sigma=0$ (see Assumption \ref{ass:intro8} for $\lambda(0,\t)=0$). As for its derivatives, we get
\begin{align}
&\partial_{\sigma}\E^{\sigma,\t}\{\mathscr{U}'(\t,\mathbf{u})-\bar{\alpha}(\t)\mathbf{u}-\tfrac12\lambda(\t)[\mathscr{U}'(\t,\mathbf{u})\mathbf{u}-1]\} \nonumber\\
&= \ \partial_{\sigma}\lambda(\sigma,\t) - \bar{\alpha}(\t) - \tfrac12\lambda(\t)\partial_{\sigma}\lambda(\sigma,\t)\sigma - \tfrac12\lambda(\t)\lambda(\sigma,\t){;} \label{eq:method43a}\\
&\partial_{\sigma}^{2}\E^{\sigma,\t}\{\mathscr{U}'(\t,\mathbf{u})-\bar{\alpha}(\t)\mathbf{u}-\tfrac12\lambda(\t)[\mathscr{U}'(\t,\mathbf{u})\mathbf{u}-1]\} \nonumber\\
&= \ \partial_{\sigma}^{2}\lambda(\sigma,\t) - \lambda(\t)\partial_{\sigma}\lambda(\sigma,\t) - \tfrac12\lambda(\t)\partial_{\sigma}^{2}\lambda(\sigma,\t)\sigma. \label{eq:method43b}
\end{align}
\eqref{eq:method43a} vanishes at $\sigma=0$ because $\bar{\alpha}(\t)=\partial_{\sigma}\E^{\sigma,\t}\mathscr{U}'(\t,\mathbf{u})|_{\sigma=0}$ (see Definition \ref{definition:intro5}) and $\lambda(0,\t)=0$ (again, see Assumption \ref{ass:intro8}). \eqref{eq:method43b} vanishes at $\sigma=0$ because, by \eqref{eq:method41g} and construction (see Definition \ref{definition:intro5}), we have $\lambda(\t)=\partial_{\sigma}^{2}\lambda(\sigma,\t)/\partial_{\sigma}\lambda(\sigma,\t)|_{\sigma=0}$. Thus, we get \eqref{eq:method42}. Let us now show the second claim, that $\bar{\mathscr{U}}'(\t,\mathbf{u})=\mathscr{U}'(\t,\mathbf{u})-\bar{\alpha}(\t)\mathbf{u}$ belongs to $\mathrm{LCT}$:
\begin{align}
\partial_{\sigma}^{\mathrm{d}}\E^{\sigma,\t}\{\mathscr{U}'(\t,\mathbf{u})-\bar{\alpha}(\t)\mathbf{u}\}|_{\sigma=0} \ = \ 0 \quad\mathrm{for}\quad\d=0,1. \label{eq:method44}
\end{align}
By \eqref{eq:method41g} and $\E^{\sigma,\t}\mathbf{u}=\sigma$, we have {that} $\mathrm{LHS}\eqref{eq:method44}=\lambda(\sigma,\t)-\bar{\alpha}(\t)\sigma$. Note {that} $\lambda(\sigma,\t)-\bar{\alpha}(\t)\sigma$ vanishes at $\sigma=0$ as $\lambda(0,\t)=0$ (see Assumption \ref{ass:intro8}). This proves \eqref{eq:method44} for $\d=0$. Its $\sigma$-derivative is $\partial_{\sigma}\lambda(\sigma,\t)-\bar{\alpha}(\t)$. By \eqref{eq:method41g}, this equals $\partial_{\sigma}\E^{\sigma,\t}\mathscr{U}'(\t,\mathbf{u})-\bar{\alpha}(\t)$. By Definition \ref{definition:intro5}, we know {that} $\partial_{\sigma}\E^{\sigma,\t}\mathscr{U}'(\t,\mathbf{u})-\bar{\alpha}(\t)$ at $\sigma=0$. This gives \eqref{eq:method44} for $\mathrm{d}=1$. We now show the last claim:
\begin{align}
&\E^{0,\t}\{\bar{\alpha}(\t)\mathbf{u}^{2}-1\} \ = \ \E^{0,\t}\{\tfrac{1}{12}\lambda(\t)^{4}[\mathscr{U}'(\t,\mathbf{u})\mathbf{u}^{3}-\E^{0,\t}(\mathscr{U}'(\t,\mathbf{u})\mathbf{u}^{3})]\} \label{eq:method45a} \\
&= \ \E^{0,\t}\{\tfrac{1}{12}\lambda(\t)^{4}\E^{0,\t}(\mathscr{U}'(\t,\mathbf{u})\mathbf{u}^{3})+\tfrac16\lambda(\t)^{3}\bar{\alpha}(\t)\mathbf{u}^{3}-\lambda(\t)\mathscr{R}(\t)\} \ = \ 0. \label{eq:method45b}
\end{align}
The second expectation in \eqref{eq:method45a} is zero because it is the expectation of {a} centered function. The same is true about the expectation in \eqref{eq:method45b} (see Definition \ref{definition:intro6} for $\mathscr{R}(\t)$). Next, apply $\partial_{\sigma}$ to both sides of the tautological identity $\sigma=\E^{\sigma,\t}\mathbf{u}$. This gives
{\small
\begin{align}
1 \ &= \ \partial_{\sigma}{\textstyle\int_{\R}}\mathbf{u}\exp\{\lambda(\sigma,\t)\mathbf{u}-\mathscr{U}(\t,\mathbf{u})+\mathscr{Z}(\sigma,\t)\}\d\mathbf{u} \nonumber\\
&= \ \partial_{\sigma}\lambda(\sigma,\t){\textstyle\int_{\R}}\mathbf{u}^{2}\exp\{\lambda(\sigma,\t)\mathbf{u}-\mathscr{U}(\t,\mathbf{u})+\mathscr{Z}(\sigma,\t)\}\d\mathbf{u} \ = \ \partial_{\sigma}\lambda(\sigma,\t)\E^{\sigma,\t}\mathbf{u}^{2}. \nonumber
\end{align}
}where the second-to-last identity follows from differentiating under the integral sign (allowed because the exponential has sub-Gaussian decay; see Assumption \ref{ass:intro8}). We get $\partial_{\sigma}\lambda(\sigma,\t)\E^{\sigma,\t}\mathbf{u}^{2}=1$. If $\sigma=0$, then $\bar{\alpha}(\t)\E^{0,\t}\mathbf{u}^{2}=1$, since $\bar{\alpha}(\t)=\partial_{\sigma}\lambda(\sigma,\t)|_{\sigma=0}$ (see Definition \ref{definition:intro5}). Thus, the first expectation in \eqref{eq:method45a} is zero. Since \eqref{eq:method45a}-\eqref{eq:method45b} is the last claim, we are done. \qed
%
%
%
\section{Proofs of Lemmas \ref{lemma:bg23}, \ref{lemma:bg1hl2}}\label{section:bg1hl2proof}
As we mentioned before, this section is a simple extension of the equivalence of ensembles estimates in \cite{DGP} (see Appendix B therein, for example). We start with a more general bound, whose sole purpose is to make precise the heuristics from after Lemma \ref{lemma:method4}. We then use it to deduce Lemmas \ref{lemma:bg23}, \ref{lemma:bg1hl2}.
\begin{lemma}\label{lemma:ee1}
 Suppose $\mathfrak{a}_{1}\in\mathrm{CT}$, $\mathfrak{a}_{2}\in\mathrm{LCT}$, and $\mathfrak{a}_{3}\in\mathrm{QCT}$. Fix $\mathfrak{l}\geq \N^{\beta_{\mathrm{BG}}}$ like in {Lemmas \ref{lemma:bg23}, \ref{lemma:bg1hl2}}. Recall $\gamma_{\mathrm{reg}}$ from {Definition \ref{definition:reg}}. We have the following with probability 1 for $\mathrm{j}=1,2,3$ and $\t\leq\t_{\mathrm{reg}}$ and $\x\in\mathbb{T}(\N)$:
\begin{align}
|\E^{\mathfrak{l},\pm}(\mathfrak{a}_{\mathrm{j}}(\t,\cdot);\t,\x)| \ \lesssim \  \N^{10\gamma_{\mathrm{reg}}}|\mathfrak{l}|^{-\frac12\mathrm{j}}\sup_{0\leq\t\leq1}\sup_{\d=0,1,2,3}\sup_{|\sigma|\lesssim1}|\partial_{\sigma}^{\d}\E^{\sigma,\t}\mathfrak{a}_{\mathrm{j}}(\t,\cdot)| \ =: \ \N^{10\gamma_{\mathrm{reg}}}|\mathfrak{l}|^{-\frac12\mathrm{j}}\|\mathfrak{a}_{\mathrm{j}}\|.\label{eq:ee1I}
\end{align}
\end{lemma}
\begin{proof}
We first recall the following observation; see Remark \ref{remark:intro14}. With probability 1, for any $\t\leq\t_{\mathrm{reg}}$ and $\x\in\mathbb{T}(\N)$, we have the following estimate on the charge density for {the} canonical measure expectation on the LHS of \eqref{eq:ee1I}, where $?\in\{+,-\}$:
\begin{align}
|\sigma(\t,\x;\mathfrak{l},?)| \ = \ |\mathfrak{l}|^{-1}| \N^{\frac12}\grad^{?\mathfrak{l}}\mathbf{J}(\t,\x)| \ \lesssim \  \N^{\gamma_{\mathrm{reg}}}|\mathfrak{l}|^{-\frac12}. \label{eq:ee1I1}
\end{align}
We now prove \eqref{eq:ee1I}. We use Corollary B.3 in \cite{DGP} with the following choices of objects therein. First, we set $F=\mathfrak{a}_{\mathrm{j}}$, which has support length, denoted by $\ell$ in Corollary B.3 of \cite{DGP}, equal to 1. Next, set charge density $\rho=\sigma(\t,\x;\mathfrak{l},\pm)$. Now, we clarify that $N$ in Corollary B.3 of \cite{DGP} is the length-scale for {the} canonical measure expectation, so it is $\mathfrak{l}$ (\emph{not} our scaling parameter $\N$). Lastly, $\sigma^{2}$ in Corollary B.3 of \cite{DGP} is the variance $\E^{\rho,\t}[\mathbf{u}-\E^{\rho,\t}\mathbf{u}]^{2}$. By Assumption \ref{ass:intro8}, $\E^{\rho,\t}$ is an expectation with respect to a measure on $\R$ whose Lebesgue density is sub-Gaussian locally uniformly in $(\rho,\t)$. Thus if $|\rho|\lesssim1$, we know said variance is $\lesssim1$. Indeed, by \eqref{eq:ee1I1}, we have $|\rho|\lesssim1$ for $\t\leq\t_{\mathrm{reg}}$ with probability 1 for our choice of $\rho$. By this application of Corollary B.3 of \cite{DGP}, we get
\begin{align}
|\E^{\mathfrak{l},\pm}(\mathfrak{a}_{\mathrm{j}}(\t,\cdot);\t,\x) - \E^{\sigma(\t,\x;\mathfrak{l},\pm),\t}\mathfrak{a}_{\mathrm{j}} - \mathrm{O}(|\mathfrak{l}|^{-1})[\partial_{\sigma}^{2}\E^{\sigma,\t}\mathfrak{a}_{\mathrm{j}}]|_{\sigma=\sigma(\t,\x;\mathfrak{l},\pm)}| \ \lesssim \ |\mathfrak{l}|^{-\frac32}\|\mathfrak{a}_{\mathrm{j}}\| \ \leq \ |\mathfrak{l}|^{-\frac12\mathrm{j}}\|\mathfrak{a}_{\mathrm{j}}\|. \label{eq:ee1I2}
\end{align}
Let us clarify \eqref{eq:ee1I2}. On the far LHS, the first term is the canonical measure expectation in \eqref{eq:ee1I}. The second expectation is of $\mathfrak{a}_{\mathrm{j}}$ with respect to the grand-canonical product measure with the same charge density $\sigma(\t,\x;\mathfrak{l},\pm)$ and time-$\t$ potential $\mathscr{U}(\t,\cdot)$; see Definition \ref{definition:intro5}. Finally, the last bound in \eqref{eq:ee1I2} follows because $|\mathfrak{l}|\geq1$ and $\mathrm{j}\leq3$. Now, by Taylor expansion,
\begin{align}
\E^{\sigma(\t,\x;\mathfrak{l},\pm),\t}\mathfrak{a}_{\mathrm{j}} \ &= \ \E^{0,\t}\mathfrak{a}_{\mathrm{j}} + [\partial_{\sigma}\E^{\sigma,\t}\mathfrak{a}_{\mathrm{j}}]|_{\sigma=0}\sigma(\t,\x;\mathfrak{l},\pm) + \tfrac12[\partial_{\sigma}^{2}\E^{\sigma,\t}\mathfrak{a}_{\mathrm{j}}]|_{\sigma=0}|\sigma(\t,\x;\mathfrak{l},\pm)|^{2}\nonumber\\
&+ \mathrm{O}(|\sigma(\t,\x;\mathfrak{l},\pm)|^{3}\|\mathfrak{a}_{\mathrm{j}}\|). \label{eq:ee1I3}
\end{align}
The big-Oh on the RHS of \eqref{eq:ee1I3} comes from estimating the third derivative $|\partial_{\sigma}^{3}\E^{\sigma,\t}\mathfrak{a}_{\mathrm{j}}|$ by $\|\mathfrak{a}_{\mathrm{j}}\|$ for $|\sigma|\leq|\sigma(\t,\x;\mathfrak{l},\pm)|\lesssim1$. (See \eqref{eq:ee1I1} and recall {that} $\mathfrak{l}\geq\N^{\beta_{\mathrm{BG}}}$ and $\beta_{\mathrm{BG}}\geq{\mathrm{C}}\gamma_{\mathrm{reg}}$ {for some large but fixed $\mathrm{C}>0$} in Definitions \ref{definition:reg}, \ref{definition:method8}.) By Taylor and a similar derivative bound for $\E^{\sigma,\t}$,
\begin{align}
[\partial_{\sigma}^{2}\E^{\sigma,\t}\mathfrak{a}_{\mathrm{j}}]|_{\sigma=\sigma(\t,\x;\mathfrak{l},\pm)} \ = \ [\partial_{\sigma}^{2}\E^{\sigma,\t}\mathfrak{a}_{\mathrm{j}}]|_{\sigma=0} + \mathrm{O}\left(|\sigma(\t,\x;\mathfrak{l},\pm)|\|\mathfrak{a}_{\mathrm{j}}\|\right). \label{eq:ee1I4}
\end{align}
Combining \eqref{eq:ee1I2}, \eqref{eq:ee1I3}, and \eqref{eq:ee1I4} and then using \eqref{eq:ee1I1} now gives the following estimate {for $\mathrm{j}=1,2,3$, $\t\leq1$ and $\x\in\mathbb{T}(\N)$}:
\begin{align}
&|\E^{\mathfrak{l},\pm}(\mathfrak{a}_{\mathrm{j}}(\t,\cdot);\t,\x)| \nonumber\\
&\lesssim \ |\E^{0,\t}\mathfrak{a}_{\mathrm{j}}| + |[\partial_{\sigma}\E^{\sigma,\t}\mathfrak{a}_{\mathrm{j}}]|_{\sigma=0}\sigma(\t,\x;\mathfrak{l},\pm)| + |[\partial_{\sigma}^{2}\E^{\sigma,\t}\mathfrak{a}_{\mathrm{j}}]|_{\sigma=0}|\sigma(\t,\x;\mathfrak{l},\pm)|^{2}| + |\sigma(\t,\x;\mathfrak{l},\pm)|^{3}\|\mathfrak{a}_{\mathrm{j}}\| \label{eq:ee1I5a} \\
&+ \ |\mathfrak{l}|^{-1}|[\partial_{\sigma}^{2}\E^{\sigma,\t}\mathfrak{a}_{\mathrm{j}}]|_{\sigma=0}| + |\mathfrak{l}|^{-1}|\sigma(\t,\x;\mathfrak{l},\pm)|\|\mathfrak{a}_{\mathrm{j}}\| + |\mathfrak{l}|^{-\frac12\mathrm{j}}\|\mathfrak{a}_{\mathrm{j}}\| \label{eq:ee1I5b} \\
&\lesssim \ |\E^{0,\t}\mathfrak{a}_{\mathrm{j}}| +  \N^{\gamma_{\mathrm{reg}}}|[\partial_{\sigma}\E^{\sigma,\t}\mathfrak{a}_{\mathrm{j}}]|_{\sigma=0}||\mathfrak{l}|^{-\frac12}+ \N^{2\gamma_{\mathrm{reg}}}|[\partial_{\sigma}^{2}\E^{\sigma,\t}\mathfrak{a}_{\mathrm{j}}]|_{\sigma=0}||\mathfrak{l}|^{-1} +  \N^{3\gamma_{\mathrm{reg}}}|\mathfrak{l}|^{-\frac32}\|\mathfrak{a}_{\mathrm{j}}\| \label{eq:ee1I5c}\\
&+ \ |\mathfrak{l}|^{-1}|[\partial_{\sigma}^{2}\E^{\sigma,\t}\mathfrak{a}_{\mathrm{j}}]|_{\sigma=0}| +  \N^{\gamma_{\mathrm{reg}}}|\mathfrak{l}|^{-\frac32}\|\mathfrak{a}_{\mathrm{j}}\| + |\mathfrak{l}|^{-\frac12\mathrm{j}}\|\mathfrak{a}_{\mathrm{j}}\|. \label{eq:ee1I5d}
\end{align}
By the reasoning in the paragraph before \eqref{eq:ee1I4}, we can bound all derivatives in \eqref{eq:ee1I5c}-\eqref{eq:ee1I5d} by $\|\mathfrak{a}_{\mathrm{j}}\|$. Now, for $\mathrm{j}=1$, the first term in \eqref{eq:ee1I5c} vanishes. The rest of \eqref{eq:ee1I5c}-\eqref{eq:ee1I5d} is $\mathrm{O}(\N^{10\gamma_{\mathrm{reg}}}|\mathfrak{l}|^{-1/2}\|\mathfrak{a}_{\mathrm{j}}\|)$. This gives \eqref{eq:ee1I} if $\mathrm{j}=1$. If $\mathrm{j}=2$, the first two terms in \eqref{eq:ee1I5c} vanish; the rest of \eqref{eq:ee1I5c}-\eqref{eq:ee1I5d} is $\mathrm{O}(\N^{10\gamma_{\mathrm{reg}}}|\mathfrak{l}|^{-1}\|\mathfrak{a}_{\mathrm{j}}\|)$. If $\mathrm{j}=3$, the first three terms in \eqref{eq:ee1I5c} and first term in \eqref{eq:ee1I5d} vanish; the rest of \eqref{eq:ee1I5c}-\eqref{eq:ee1I5d} is $\mathrm{O}(\N^{10\gamma_{\mathrm{reg}}}|\mathfrak{l}|^{-3/2}\|\mathfrak{a}_{\mathrm{j}}\|)$. The last two sentences give \eqref{eq:ee1I} for $\mathrm{j}=2,3$, thereby completing the proof.
\end{proof}
\begin{proof}[Proofs of {Lemmas \ref{lemma:bg23}, \ref{lemma:bg1hl2}}]
Use Lemma \ref{lemma:ee1} for $\mathfrak{a}_{1}=\mathfrak{w}$ and $\mathfrak{a}_{2}=\mathfrak{d}$ in Proposition \ref{prop:bg1hl1}, and for $\mathfrak{a}_{3}=\mathfrak{q}$ in Proposition \ref{prop:bg22}. This finishes the proofs, if we can control the $\|\|$-norm in \eqref{eq:ee1I} for any $\mathfrak{a}(\t,\cdot):\R\to\R$ satisfying $|\partial_{\mathbf{u}}^{\d}\mathfrak{a}(\t,\mathbf{u})|\lesssim_{\d}1+\mathbf{u}^{10}$:
\begin{align}
|\partial_{\sigma}^{\d}\E^{\sigma,\t}\mathfrak{a}(\t,\cdot)| \ \lesssim \ 1. \label{eq:eep1}
\end{align}
(This bound must be uniform in $\d=0,1,2,3$ and $0\leq\t\leq1$ and $|\sigma|\lesssim1$.) By construction (see Definition \ref{definition:intro5}) and Assumption \ref{ass:intro8}, $\E^{\sigma,\t}$ is {the} expectation with respect to a measure on $\R$ that is sub-Gaussian uniformly in $0\leq\t\leq1$ and $|\sigma|\lesssim1$. Its density with respect to Lebesgue on $\R$ is smooth in $\sigma$ and also has sub-Gaussian decay (uniformly in $0\leq\t\leq1$ and $|\sigma|\lesssim1$). Because we assumed $\mathfrak{a}$ has polynomial growth uniformly in $0\leq\t\leq1$, standard measure theory and calculus implies {that} $\E^{\sigma,\t}\mathfrak{a}(\t,\cdot)$ is smooth in $\sigma$ (uniformly in $0\leq\t\leq1$). Thus, \eqref{eq:eep1} follows (as smooth functions are locally uniformly bounded), and we are done.
\end{proof}
%
%
%
\section{Proofs for Lemmas \ref{lemma:method6}, \ref{lemma:method7}, \ref{lemma:method10}, \ref{lemma:method13}}
\subsection{Proof of Lemma \ref{lemma:method6}}
Set $\|\|=\|\|_{\t_{\mathrm{reg}};\mathbb{T}(\N)}$ for this proof. We claim {that} it suffices to prove $\|\mathbf{Z}-\mathbf{S}\|\lesssim\N^{-20\gamma_{\mathrm{reg}}}\|\mathbf{Z}\|$. Indeed, by this and triangle inequality, we get $\|\mathbf{Z}-\mathbf{S}\|\lesssim\N^{-20\gamma_{\mathrm{reg}}}\|\mathbf{Z}-\mathbf{S}\|+\N^{-20\gamma_{\mathrm{reg}}}\|\mathbf{S}\|$, so $\|\mathbf{Z}-\mathbf{S}\|\lesssim\N^{-20\gamma_{\mathrm{reg}}}\|\mathbf{S}\|$, which would finish the proof. To prove $\|\mathbf{Z}-\mathbf{S}\|\lesssim\N^{-20\gamma_{\mathrm{reg}}}\|\mathbf{Z}\|$, first note {that} $\mathbf{H}^{\N}(\s,\t,\x)[1]=1$, i.e. constants are eigenfunctions of $\mathbf{H}^{\N}$. Thus,
\begin{align}
\mathbf{Z}(\t,\x)-\mathbf{S}(\t,\x) \ = \ \mathbf{H}^{\N}(\t,\t(\N),\x)[\mathbf{Z}(\t,\x)-\mathbf{Z}(\t,\cdot)]. \label{eq:method61}
\end{align}
Fix $\t\leq\t_{\mathrm{reg}}$ and $\x,\y\in\mathbb{T}(\N)$. Set $\mathbf{h}(\vee):=\mathbf{h}(\t,\x)\vee\mathbf{h}(\t,\y)$ and $\mathbf{h}(\wedge):=\mathbf{h}(\t,\x)\wedge\mathbf{h}(\t,\y)$. (We are abusing notation by not writing dependence on $\t,\x,\y$ here, but we shortly stop using this notation anyway.) By calculus and $|\lambda(\t)|\lesssim1$, we have
\begin{align}
|\mathbf{Z}(\t,\x)-\mathbf{Z}(\t,\y)| \ \lesssim \ {\textstyle\int_{\mathbf{h}(\wedge)}^{\mathbf{h}(\vee)}}\exp[|\lambda(\t)|\upsilon]\d\upsilon \ &\lesssim \ |\mathbf{h}(\vee)-\mathbf{h}(\wedge)||\mathbf{Z}(\t,\x)+\mathbf{Z}(\t,\y)| \nonumber\\
&\lesssim \ \N^{-\frac12+\gamma_{\mathrm{reg}}}|\x-\y|^{\frac12}|\mathbf{Z}\|_{\t_{\mathrm{st}};\mathbb{T}(\N)}. \label{eq:method62}
\end{align}
(The second bound in \eqref{eq:method62} is by bounding the integrand at the two limits of integration. The last bound holds by construction of $\t_{\mathrm{reg}}$ in Definition \ref{definition:reg}.) By \eqref{eq:method61}-\eqref{eq:method62}, we deduce the following deterministic estimate for any $\t\leq\t_{\mathrm{st}}$ and $\x\in\mathbb{T}(\N)$:
\begin{align}
|\mathbf{Z}(\t,\x)-\mathbf{S}(\t,\x)| \ &\lesssim \ \N^{\gamma_{\mathrm{reg}}}\|\mathbf{Z}\|\sum_{\y\in\mathbb{T}(\N)}\mathbf{H}^{\N}(\t,\t(\N),\x,\y)\times\N^{-\frac12}|\x-\y|^{\frac12} \nonumber\\
&\lesssim \ \N^{\gamma_{\mathrm{reg}}}\|\mathbf{Z}\|[|\t(\N)-\t|^{\frac14}+\N^{-\frac12}]. \label{eq:method63}
\end{align}
(The last bound follows by a standard moment bound for random walks; see {\eqref{eq:hke6}}.) {Since the previous bound \eqref{eq:method63} is deterministic and uniform in $\t\leq\t_{\mathrm{st}}$, $\x\in\mathbb{T}(\N)$}, we can replace the far LHS of \eqref{eq:method63} by its $\|\|_{\t_{\mathrm{st}};\mathbb{T}(\N)}$-norm. It now suffices to recall from Definition \ref{definition:method5} that $\t(\N)-\t=\N^{-100\gamma_{\mathrm{reg}}}$ in order to prove the desired bound in Lemma \ref{lemma:method6}. \qed
\subsection{Proof of Lemma \ref{lemma:method7}}
Note {that} $\d_{\t}\mathbf{H}^{\N}(\s,\t,\x,\y)=\mathscr{T}(\t)\mathbf{H}^{\N}(\s,\t,\x,\y)\d\t$ as measures, in which $\mathscr{T}(\t)$ acts on $\x$; see Definition \ref{definition:method1}. We also have $\d_{\s}\mathbf{H}^{\N}(\s,\t,\x,\y)=-\mathscr{T}(\s)^{\ast}\mathbf{H}^{\N}(\s,\t,\x,\y)\d\s$ as measures, where $\mathscr{T}(\s)^{\ast}$ is the adjoint of $\mathscr{T}(\s)$ with respect to {the} uniform measure on $\mathbb{T}(\N)$, and it now acts on the $\y$-variable. (This is the so-called ``adjoint equation"; it is the Kolmogorov PDE for the heat kernel of the adjoint random walk. The sign in front of $\mathscr{T}(\s)^{\ast}$ is there because $\partial_{\s}$ acts forward in $\s$, while the adjoint random walk moves ``backwards" in $\s$. See {\eqref{eq:hke00}}.) In what follows, all $\mathscr{T}$-operators without adjoints always act on the heat kernel in $\x$; the adjoints always act on $\y$. Moreover, all differential equations are shorthand for the corresponding (rigorous) integrated equations, in which terms like $\mathscr{T}(\t)\mathbf{H}^{\N}(\s,\t,\x,\y)$ and $-\mathscr{T}(\s)^{\ast}\mathbf{H}^{\N}(\s,\t,\x,\y)$ have a clear meaning.

We now compute the stochastic differential $\d\mathbf{S}(\t,\x)=\d\mathbf{H}^{\N}(\t,\t(\N),\x)[\mathbf{Z}(\t,\cdot)]$. We first claim {that} the (Ito) chain rule gives
\begin{align}
\d\mathbf{S}(\t,\x) \ = \ {\textstyle\sum_{\y}}\{\d\mathbf{H}^{\N}(\t,\t(\N),\x,\y)\}\mathbf{Z}(\t,\y) + {\textstyle\sum_{\y}}\mathbf{H}^{\N}(\t,\t(\N),\x,\y)\d\mathbf{Z}(\t,\y) \ =: \ \Xi[1]+\Xi[2], \label{eq:method71}
\end{align}
where sums are over $\y\in\mathbb{T}(\N)$. Indeed, there is no cross-variation term because $\mathbf{H}^{\N}$ is smooth in both time-variables. Using the $\mathbf{H}^{\N}$ PDEs that we explained in the first paragraph, we claim that the first term on the RHS is computed as
\begin{align}
\d\Xi[1] \ &= \ {\textstyle\sum_{\y}}\mathscr{T}(\t(\N))\mathbf{H}(\t,\t(\N),\x,\y)\mathbf{Z}(\t,\y)\d\t-{\textstyle\sum_{\y}}\mathscr{T}(\t)^{\ast}\mathbf{H}^{\N}(\t,\t(\N),\x,\y)\mathbf{Z}(\t,\y)\d\t \nonumber\\
&= \ \mathscr{T}(\t(\N))\mathbf{S}(\t,\x)\d\t-\mathbf{H}^{\N}(\t,\t(\N),\x)[\mathscr{T}(\t)\mathbf{Z}(\t,\cdot)]\d\t. \label{eq:method72}
\end{align}
(To get the last identity in \eqref{eq:method72}, pull $\mathscr{T}(\t(\N))$ outside the $\y$-sum, since it acts on $\x$, and use $\mathbf{S}(\t,\x)=\mathbf{H}^{\N}(\t,\t(\N),\x)[\mathbf{Z}(\t,\cdot)]$. For the second sum in the first line, swap $\mathscr{T}(\t)^{\ast}$ for $\mathscr{T}(\t)$ and move it on $\mathbf{Z}(\t,\y)$.) We compute $\Xi[2]$ by letting $\mathbf{H}^{\N}(\t,\t(\N),\x)$ act on $\mathrm{RHS}\eqref{eq:method2I}$. What falls out is every term in $\mathrm{RHS}\eqref{eq:method7Ia}$ and $\eqref{eq:method7Ib}$, except we get $\mathbf{H}^{\N}(\t,\t(\N),\x)[\mathscr{T}(\t)\mathbf{Z}(\t,\cdot)]\d\t$ instead of $\mathscr{T}(\t(\N))\mathbf{S}(\t,\x)\d\t$ in $\mathrm{RHS}\eqref{eq:method7Ia}$. But once we add $\Xi[1]$ (see \eqref{eq:method72}), the $\mathbf{H}^{\N}(\t,\t(\N),\x)[\mathscr{T}(\t)\mathbf{Z}(\t,\cdot)]\d\t$ that we got is cancelled, and the $\mathscr{T}(\t(\N))\mathbf{S}(\t,\x)\d\t$ that we were missing is now there. This establishes the first identity $\mathrm{LHS}\eqref{eq:method7Ia}=\mathrm{RHS}\eqref{eq:method7Ia}+\eqref{eq:method7Ib}$. The rest of that display follows by elementary manipulations. The display after follows by the Duhamel principle. In particular, Duhamel gives the following, in which we use notation to be explained afterwards:
\begin{align}
\mathbf{S}(\t,\x) \ &= \ \mathbf{H}^{\N}(0(\N),\t(\N),\x)(\mathbf{S}(0,\cdot))\nonumber\\
&+ \ {\textstyle\int_{0}^{\t}}[\mathbf{H}^{\N}(\s(\N),\t(\N),\cdot)\circ\mathbf{H}^{\N}(\s,\s(\N),\cdot)](\x)(\sqrt{2}\lambda(\s) \N^{\frac12}\mathbf{Z}(\s,\cdot)\d\mathbf{b}(\s,\cdot(\s))) \nonumber \\
&+ \ {\textstyle\int_{0}^{\t}}[\mathbf{H}^{\N}(\s(\N),\t(\N),\cdot)\circ\mathbf{H}^{\N}(\s,\s(\N),\cdot)](\x)(\partial_{\s}\log|\lambda(\s)|\cdot\mathbf{Z}(\s,\cdot)\log\mathbf{Z}(\s,\cdot))\d\s \\
&+ \ {\textstyle\int_{0}^{\t}}[\mathbf{H}^{\N}(\s(\N),\t(\N),\cdot)\circ\mathbf{H}^{\N}(\s,\s(\N),\cdot)](\x)({\mathfrak{z}}(\s,\cdot(\s))\mathbf{Z}(\s,\cdot))\d\s.
\end{align}
Here, $[\mathbf{H}^{\N}(\r,\t,\cdot)\circ\mathbf{H}^{\N}(\s,\r,\cdot)](\x)$ is the composition of first acting by $\mathbf{H}^{\N}(\s,\r,\cdot)$ and then by $\mathbf{H}^{\N}(\r,\t,\cdot)$. But by the semigroup property (see {\eqref{eq:hke0}}), we know that $[\mathbf{H}^{\N}(\r,\t,\cdot)\circ\mathbf{H}^{\N}(\s,\r,\cdot)](\x)=\mathbf{H}^{\N}(\s,\t,\x)$. In particular, we deduce 
\begin{align}
\mathbf{H}^{\N}(0(\N),\t(\N),\x)(\mathbf{S}(0,\cdot))&=[\mathbf{H}^{\N}(0(\N),\t(\N),\cdot)\circ\mathbf{H}^{\N}(0,0(\N),\cdot)](\x)(\mathbf{Z}(0,\cdot))\\
&=\mathbf{H}^{\N}(0,\t(\N),\x)(\mathbf{Z}(0,\cdot)).
\end{align}
Plugging this into the previous display gives \eqref{eq:method7IIa}-\eqref{eq:method7IId}. Finally, \eqref{eq:method7III} follows by construction of ${\mathfrak{z}}$ from Proposition \ref{prop:method2}. ($\mathrm{QCT}$ collects the RHS of \eqref{eq:method2IIb} and \eqref{eq:method2IIc}. $\mathrm{LCT}$ collects \eqref{eq:method2IId}. $\mathrm{CT}$ collects \eqref{eq:method2IIe}-\eqref{eq:method2IIg}. $\mathrm{An}$ collects \eqref{eq:method2IIi}-\eqref{eq:method2IIj}.) This completes the proof. \qed
\subsection{Proof of Lemma \ref{lemma:method10}}
We let $\tau$ be the minimum of $\t_{\mathrm{st}}$ and the explosion times for $\mathbf{Y},\mathbf{S},\log\mathbf{Y},\log\mathbf{S}$. We emphasize that $\tau$ is a stopping time, since all processes are adapted to the same Brownian filtration (and $\t_{\mathrm{st}}$ is a stopping time; see Definition \ref{definition:reg}). Note that the SDEs for $\mathbf{Y}$ and $\mathbf{S}$ (see Lemma \ref{lemma:method7} and Remark \ref{remark:method9}) are the same for times in $[0,\tau)$. So, standard Ito theory implies $\mathbf{Y}(\t,\x)=\mathbf{S}(\t,\x)$ for all $\t\leq\tau$. We are now left to show $\tau=\t_{\mathrm{st}}$ with probability 1. We first note {that} the explosion time of $\mathbf{S}$ is infinite. Indeed, it is a bounded operator acting on $\mathbf{Z}$, which is a continuous function of \eqref{eq:hf}-\eqref{eq:glsde}, which has infinite life-time because its coefficients are uniformly Lipschitz in the solution by Assumption \ref{ass:intro8}. It now suffices to show the explosion time of $\mathbf{Y},\log\mathbf{Y}$ are both at least $\t_{\mathrm{st}}$ with probability 1. Indeed, if one of them is strictly less than $\t_{\mathrm{st}}$, then $\mathbf{Y}=\mathbf{S}$ at the explosion time of $\mathbf{Y}$ or $\log\mathbf{Y}$. This would imply finiteness of explosion time of $\mathbf{S}$ or $\log\mathbf{S}$, thus giving a contradiction, so we are done. \qed
\subsection{Proof of Lemma \ref{lemma:method13}}
We bound $\mathrm{An}$. Fix $\t\leq\t_{\mathrm{st}}$ and $\x\in\mathbb{T}(\N)$. Again, the operator $\mathbf{H}^{\N}(\s,\t(\N),\x)$ is a convolution operator, {since it is the semigroup for a spatially-homogeneous infinitesimal generator}. Thus, it commutes {with the} constant-coefficient gradients $\grad^{+}$ and $\grad^{-}$, so we can move $\grad^{+}$ and $\grad^{-}$ to the $\mathbf{H}^{\N}$ heat kernel in \eqref{eq:method7IVi}. Now, for $\t\leq\t_{\mathrm{st}}$, recall from Proposition \ref{prop:method2} that the $\mathfrak{a}^{\pm},\mathfrak{b},\mathfrak{c}$-terms in \eqref{eq:method7IVh}-\eqref{eq:method7IVi} are $\lesssim\N^{100\gamma_{\mathrm{reg}}}$ with probability 1. Moreover, the operator norms of $\mathbf{H}^{\N}(\s,\t(\N),\x)$ and $\N\grad^{?}\mathbf{H}^{\N}(\s,\t(\N),\x)$ (for {$?\in\{\pm\}$}) on $\mathscr{L}^{\infty}(\mathbb{T}(\N))$ are $\lesssim1+|\t(\N)-\s|^{-1/2}\leq1+|\t-\s|^{-1/2}$ (see {\eqref{eq:hke3}}). We ultimately deduce
\begin{align}
|\mathrm{An}(\t,\x)| \ \lesssim \ \N^{-\frac12+100\gamma_{\mathrm{reg}}}{\textstyle\int_{0}^{\t}}[1+|\t-\s|^{-\frac12}]\d\s \ \lesssim \ \N^{-\frac13}{,} \label{eq:method131}
\end{align}
for $\t\leq\t_{\mathrm{st}}$, where the last bound follows since $\t_{\mathrm{st}}\leq1$ (see Definition \ref{definition:method8}) and since $\gamma_{\mathrm{reg}}$ is small (see Definitions \ref{definition:entropydata}, \ref{definition:reg}). The previous estimate \eqref{eq:method131} is deterministic for $\t\leq\t_{\mathrm{st}}$ and $\x\in\mathbb{T}(\N)$. Thus, the first bound in \eqref{eq:method13I} follows. We now control \eqref{eq:method7IIc}. Note the function $\mathrm{a}\mapsto\mathrm{a}\log\mathrm{a}$ has derivative $1+\log\mathrm{a}\lesssim\mathrm{a}+\mathrm{a}^{-1}$ (where we restrict to $\mathrm{a}>0$). Thus, given any $\s,\y$, we have the upper bound $|\mathbf{Z}(\s,\y)\log\mathbf{Z}(\s,\y)-\mathbf{S}(\s,\y)\log\mathbf{S}(\s,\y)|\lesssim\{|\mathbf{Z}(\s,\y)|+|\mathbf{Z}(\s,\y)|^{-1}+|\mathbf{Y}(\s,\y)|+|\mathbf{Y}(\s,\y)|^{-1}\}|\mathbf{Z}(\s,\y)-\mathbf{S}(\s,\y)|$. Now, if $\s\leq\t_{\mathrm{st}}$, then by Lemma \ref{lemma:method6}, we have that $|\mathbf{Z}(\s,\y)-\mathbf{S}(\s,\y)|\lesssim\N^{-20\gamma_{\mathrm{reg}}}\|\mathbf{Z}\|_{\t_{\mathrm{st}};\mathbb{T}(\N)}$. Also, by Definition \ref{definition:method8}, we have $\|\mathbf{Z}\|+\|\mathbf{Z}^{-1}\|+\|\mathbf{S}\|+\|\mathbf{S}^{-1}\|\lesssim\N^{\gamma_{\mathrm{ap}}}$, where $\|\|=\|\|_{\t_{\mathrm{st}};\mathbb{T}(\N)}$. So, for $\s\leq\t_{\mathrm{st}}$, we get $|\mathbf{Z}(\s,\y)\log\mathbf{Z}(\s,\y)-\mathbf{S}(\s,\y)\log\mathbf{S}(\s,\y)|\lesssim\N^{2\gamma_{\mathrm{ap}}-20\gamma_{\mathrm{reg}}}\lesssim\N^{-10\gamma_{\mathrm{reg}}}$, since $\gamma_{\mathrm{ap}}$ is small compared to $\gamma_{\mathrm{reg}}$ (see Definitions \ref{definition:reg}, \ref{definition:method8}). It now suffices to use {the} contractivity of $\mathbf{H}^{\N}$ (namely, what we used to get \eqref{eq:method131}) to deduce the \eqref{eq:method7IIc} estimate in \eqref{eq:method13I}. (We controlled what is inside the heat operator in \eqref{eq:method7IIc} by $\lesssim\N^{-10\gamma_{\mathrm{reg}}}$ with probability 1, so we actually deduce the \eqref{eq:method7IIc} bound in \eqref{eq:method13I} with an improved upper bound of $\lesssim\N^{-10\gamma_{\mathrm{reg}}}$.) It now suffices to control \eqref{eq:method7IId}, which is the only stochastic bound. (In particular, the previous estimates are all deterministic, since all the randomness was put into the stopping time $\t_{\mathrm{st}}$.) We claim {that} the following holds for any $\t\leq\t_{\mathrm{st}}$ and $\x\in\mathbb{T}(\N)$:
\begin{align}
\eqref{eq:method7IId} \ &= \ {\textstyle\int_{0}^{\t}}\mathbf{H}(\s,\t(\N),\x)\{\sqrt{2}\lambda(\s)\N^{\frac12}\Theta(\s,\cdot)\d\mathbf{b}(\s,\cdot(\s))\}, \nonumber\\
\Theta(\s,\cdot) \ &:= \ [\mathbf{Z}(\s,\cdot)-\mathbf{S}(\s,\cdot)]\mathbf{1}\{|\mathbf{Z}(\s,\cdot)-\mathbf{S}(\s,\cdot)|\lesssim\N^{-10\gamma_{\mathrm{reg}}}\}. \nonumber
\end{align}
Indeed, this identity is just saying we can put the indicator of $|\mathbf{Z}(\s,\cdot)-\mathbf{S}(\s,\cdot)|\lesssim\N^{-10\gamma_{\mathrm{reg}}}$ into the stochastic integral in \eqref{eq:method7IId} for free. Indeed, as shown in the previous paragraph, we have the deterministic estimate $|\mathbf{Z}(\s,\cdot)-\mathbf{S}(\s,\cdot)|\lesssim\N^{-10\gamma_{\mathrm{reg}}}$ for any $\s\leq\t_{\mathrm{st}}$. The upshot of this identity is that it now suffices to bound a stochastic integral of the heat operator acting on an adapted process that is uniformly $\lesssim\N^{-10\gamma_{\mathrm{reg}}}$. The \eqref{eq:method7IId} bound in \eqref{eq:method13I} now holds by standard methods for one-dimensional stochastic heat equations, i.e. BDG and Kolmogorov continuity estimates for space-time Holder norms. (We refer to Section 3 of \cite{DT}, for example, which actually deals with a more complicated noise than independent standard Brownian motions. We also clarify that the upper bound we proposed in \eqref{eq:method13I} has a less negative exponent than the uniform bound of $\lesssim\N^{-10\gamma_{\mathrm{reg}}}$ we have for $\Theta$ above. This is just the statement that if we want \eqref{eq:method7IId} to exceed its natural a priori scale of $\lesssim\N^{-10\gamma_{\mathrm{reg}}}$ by a strictly positive power of $\N$, this should happen with low probability.) Having now shown the \eqref{eq:method7IId} bound in \eqref{eq:method13I}, we are done. \qed
%
%
%
\section{Proof of Proposition \ref{prop:bg1hl1}}
As we mentioned before Section \ref{section:bg22proofoutline}, this section amounts to modifying (very straightforwardly) the proof of Proposition \ref{prop:bg22}. We explain this as follows.

We first bound $\mathscr{R}^{\mathfrak{d},?}$ for $?\in\{+,-\}$. We start by introducing the notation $\t\mapsto\t(\N;2)=\t+2^{-1}\N^{-100\gamma_{\mathrm{reg}}}$ as the midpoint between $\t$ and $\t(\N)$ from Definition \ref{definition:method5}. With explanation given after, we claim the following:
\begin{align}
\mathscr{R}^{\mathfrak{d},?}(\t,\x) \ &= \ {\textstyle\int_{0}^{\t}}[\mathbf{H}^{\N}(\t(\N;2),\t(\N),\cdot)\circ\mathbf{H}^{\N}(\s,\t(\N;2),\cdot)](\x)\{\N^{\frac32}\grad^{?}[\mathds{R}^{\mathfrak{d},?}(\s,\cdot(\s))\mathbf{Z}(\s,\cdot)]\}\d\s \label{eq:bg1hl11a}\\
&= \ {\textstyle\int_{0}^{\t}}[\N\grad^{?}\mathbf{H}^{\N}(\t(\N;2),\t(\N),\cdot)\circ\mathbf{H}^{\N}(\s,\t(\N;2),\cdot)](\x)\{\N^{\frac12}\mathds{R}^{\mathfrak{d},?}(\s,\cdot(\s))\mathbf{Z}(\s,\cdot)\}\d\s \label{eq:bg1hl11b}\\
&= \ \N\grad^{?}\mathbf{H}^{\N}(\t(\N;2),\t(\N),\x)\{{\textstyle\int_{0}^{\t}}\mathbf{H}^{\N}(\s,\t(\N;2),\star)[\N^{\frac12}\mathds{R}^{\mathfrak{d},?}(\s,\cdot(\s))\mathbf{Z}(\s,\cdot)]\d\s\}. \label{eq:bg1hl11c}
\end{align}
\eqref{eq:bg1hl11a} follows by the Chapman-Kolmogorov equation (or semigroup property) from {\eqref{eq:hke0}} and {the} definition of $\mathscr{R}^{\mathfrak{d},?}$ from Proposition \ref{prop:bg1hl1}. (We clarify that the $\circ$ in $\mathrm{RHS}\eqref{eq:bg1hl11a}$ means first acting by the heat operator $\mathbf{H}^{\N}(\s,\t(\N;2),\cdot)$, composing with $\mathbf{H}^{\N}(\t(\N;2),\t(\N),\cdot)$, and evaluating the composition at $\x$, i.e. dot product with its $\x$-th row.) \eqref{eq:bg1hl11b} holds because heat operators are convolution operators. Thus we can move the constant-coefficient gradient $\N\grad^{?}$ past convolutions such that it hits the outermost convolution. Let us now explain \eqref{eq:bg1hl11c}. It applies $\N\grad^{?}$ to the heat operator $\mathbf{H}^{\N}(\t(\N;2),\t(\N),\x)$ acting on a space-time function $\mathds{F}(\t,\star)$. This function is what is inside the curly brackets in \eqref{eq:bg1hl11c}. (In particular,  the outermost $\mathbf{H}^{\N}$ in \eqref{eq:bg1hl11c} is a convolution operator with integration-variable $\star$. The innermost $\mathbf{H}^{\N}$ is convolution with integration-variable $\cdot$.) Thus,  \eqref{eq:bg1hl11c} holds because the map $\N\grad^{?}\mathbf{H}^{\N}(\t(\N;2),\t(\N),\cdot)$ is linear, so it commutes with $\d\s$-integration. (Note that $\N\grad^{?}\mathbf{H}^{\N}(\t(\N;2),\t(\N),\cdot)$ is independent of $\s$.) By {\eqref{eq:hke3}}, the $\mathscr{L}^{\infty}$-operator norm of $\N\grad^{?}\mathbf{H}^{\N}(\t(\N;2),\t(\N),\cdot)$ is $\lesssim|\t(\N)-\t(\N;2)|^{-1/2}\lesssim\N^{50\gamma_{\mathrm{reg}}}$. Using this with \eqref{eq:bg1hl11a}-\eqref{eq:bg1hl11c} and linearity of the map $\mathfrak{d}\mapsto\mathds{R}^{\mathfrak{d},?}$ (see Proposition \ref{prop:bg1hl1} for notation), we deduce the following (that we clarify afterwards), where $\|\|:=\|\|_{\t_{\mathrm{st}};\mathbb{T}(\N)}$ for the rest of this proof, and where $\mathfrak{d}[\wedge]:=\N^{-1/2}\mathfrak{d}$:
\begin{align}
\|\mathscr{R}^{\mathfrak{d},?}\| \ &\lesssim \ \N^{50\gamma_{\mathrm{reg}}}\|{\textstyle\int_{0}^{\t}}\mathbf{H}^{\N}(\s,\t(\N;2),\x)[\N^{\frac12}\mathds{R}^{\mathfrak{d},?}(\s,\cdot(\s))\mathbf{Z}(\s,\cdot)]\d\s\| \label{eq:bg1hl12a}\\
&= \ \N^{50\gamma_{\mathrm{reg}}}\|{\textstyle\int_{0}^{\t}}\mathbf{H}^{\N}(\s,\t(\N;2),\x)[\N\mathds{R}^{\mathfrak{d}[\wedge],?}(\s,\cdot(\s))\mathbf{Z}(\s,\cdot)]\d\s\| \label{eq:bg1hl12b}.
\end{align}
(In \eqref{eq:bg1hl12b}, $\mathds{R}^{\mathfrak{d}[\wedge],?}$ is $\mathds{R}^{\mathfrak{d},?}$ from Proposition \ref{prop:bg1hl1} but replacing $\mathfrak{d}$ therein by $\mathfrak{d}[\wedge]$.) By Lemma \ref{lemma:bg1hl2} and the fact that the length-scale $\mathfrak{l}$ always satisfies $\mathfrak{l}\lesssim\N$, we deduce that $\mathfrak{d}[\wedge]$ satisfies the bounds in Lemma \ref{lemma:bg23}. As noted in Section \ref{subsection:bg2word}, this means {that} Proposition \ref{prop:bg22} applies to \eqref{eq:bg1hl12b}. (Section \ref{subsection:bg2word} also notes that taking $\t(\N;2)$ in \eqref{eq:bg1hl12b} instead of $\t(\N)$ in Proposition \ref{prop:bg22} is unimportant.) Using this and \eqref{eq:bg1hl12a}-\eqref{eq:bg1hl12b}, we know that $\|\mathscr{R}^{\mathfrak{d},?}\|\lesssim\N^{50\gamma_{\mathrm{reg}}-3\beta_{\mathrm{BG}}}\|\mathbf{Z}\|\lesssim\N^{-2\beta_{\mathrm{BG}}}\|\mathbf{Z}\|$, where the last bound holds because $\gamma_{\mathrm{reg}}$ is small compared to $\beta_{\mathrm{BG}}$ (see Definitions \ref{definition:reg}, \ref{definition:method8}), with high probability for each $?\in\{+,-\}$. We now control $\mathscr{R}^{\mathfrak{w},2,?}$ for $?\in\{+,-\}$. The argument is the same with the following adjustments. First, $\mathscr{R}^{\mathfrak{w},2,?}$ is $\mathscr{R}^{\mathfrak{d},?}$ but {multiplied} by $\N^{-1/2}$ and {by replacing} $\mathfrak{d}$ by $\mathfrak{w}$. So, like \eqref{eq:bg1hl12a}-\eqref{eq:bg1hl12b},
\begin{align}
\|\mathscr{R}^{\mathfrak{w},2,?}\| \ &\lesssim \ \N^{50\gamma_{\mathrm{reg}}}\|{\textstyle\int_{0}^{\t}}\mathbf{H}^{\N}(\s,\t(\N;2),\x)[\N\mathds{R}^{\mathfrak{w}[\wedge],?}(\s,\cdot(\s))\mathbf{Z}(\s,\cdot)]\d\s\|, \label{eq:bg1hl13}
\end{align}
in which $\mathfrak{w}[\wedge]:=\N^{-1}\mathfrak{w}$. (The point is that $\mathfrak{w}[\wedge]$ multiplies by an additional $\N^{-1/2}$ compared to $\mathfrak{d}[\wedge]$. This $\N^{-1/2}$-factor comes for free, however, as $\mathscr{R}^{\mathfrak{w},2,?}$ has an additional $\N^{-1/2}$ compared to $\mathscr{R}^{\mathfrak{d},?}$; see before \eqref{eq:bg1hl13}.) By Lemma \ref{lemma:bg1hl2} and the fact that the length $\mathfrak{l}$ satisfies $\mathfrak{l}\lesssim\N$, we know {that} $\mathfrak{w}[\wedge]$ satisfies the bound in Lemma \ref{lemma:bg23}. Thus, $\|\mathscr{R}^{\mathfrak{w},2,?}\|\lesssim\N^{50\gamma_{\mathrm{reg}}-3\beta_{\mathrm{BG}}}\|\mathbf{Z}\|\lesssim\N^{-2\beta_{\mathrm{BG}}}\|\mathbf{Z}\|$ with high probability as before. We now control $\mathscr{R}^{\mathfrak{w},1,?}$ for $?\in\{+,-\}$. As $\mathscr{R}^{\mathfrak{w},1,?}$ has no $\N\grad$-operator, we can more easily get
\begin{align}
\|\mathscr{R}^{\mathfrak{w},1,?}\| \ &\lesssim \ \|{\textstyle\int_{0}^{\t}}\mathbf{H}^{\N}(\s,\t(\N),\x)[\N\mathds{R}^{\mathfrak{w}[\wedge],?}(\s,\cdot(\s))\mathbf{Z}(\s,\cdot)]\d\s\|. \label{eq:bg1hl14}
\end{align}
Now, use the paragraph between \eqref{eq:bg1hl13}-\eqref{eq:bg1hl14} to deduce {that} $\|\mathscr{R}^{\mathfrak{w},1,?}\|\lesssim\N^{-2\beta_{\mathrm{BG}}}\|\mathbf{Z}\|$ with high probability from \eqref{eq:bg1hl14}. {The union} bound now shows that with high probability, we get $\|\Gamma\|\lesssim\N^{-2\beta_{\mathrm{BG}}}\|\mathbf{Z}\|$ for $\Gamma=\mathscr{R}^{\mathfrak{d},?},\mathscr{R}^{\mathfrak{w},\mathrm{i},?}$ and $\mathrm{i}=1,2$ and $?=+,-$. On this high probability event, the triangle inequality gives the desired estimate \eqref{eq:bg1hl1I}, so the proof is complete. \qed
%
%
%
\section{Heat kernel estimates}
The results in this section basically make precise the idea that the $\mathbf{H}^{\N}$-kernel is very close to the heat kernel for the continuum $\mathbf{H}$-semigroup from the introduction, and that ``very close" means in a fairly strong topology.

The first result in this section concerns the semi-discrete heat kernel $\mathbf{H}^{\N}$. In a nutshell, we establish its semigroup property and adjoint equation. We then remove the delta-function in \eqref{eq:method1IT} by following a constant-speed characteristic. We conclude with pointwise and summed (regularity) bounds for the $\mathbf{H}^{\N}$ kernel, as well as a moment bound for the underlying random walk. (The reader is invited to compare Proposition \ref{prop:hke} to what happens for Brownian motion heat kernels on $\mathbb{T}$ of speed $\N^{2}$.)
\begin{prop}\label{prop:hke}
 Take $\s\leq\r\leq\t$. We {have} the following Chapman-Kolmogorov equation (or semigroup property), which is an identity of operators, where $[\mathbf{H}^{\N}(\r,\t,\cdot)\circ\mathbf{H}^{\N}(\s,\r,\cdot)](\x)$ is the composition of first acting by $\mathbf{H}^{\N}(\s,\r,\cdot)$ and then by $\mathbf{H}^{\N}(\r,\t,\cdot)$:
\begin{align}
[\mathbf{H}^{\N}(\r,\t,\cdot)\circ\mathbf{H}^{\N}(\s,\r,\cdot)](\x) \ = \ \mathbf{H}^{\N}(\s,\t,\x). \label{eq:hke0}
\end{align}
Fix $\s\in\R$. Let $\mathscr{T}(\s)^{\ast}$ be the adjoint of $\mathscr{T}(\s)$ with respect to {the} uniform measure on $\mathbb{T}(\N)$. (In particular, replace all $\grad$-operators in $\mathrm{RHS}\eqref{eq:method1IT}$, which we evaluate at time $\s$ in this context, by their negatives.) We have the following, where $\mathscr{T}(\s)^{\ast}$ acts on $\y$:
\begin{align}
\partial_{\s}\mathbf{H}^{\N}(\s,\t,\x,\y) \ = \ -\mathscr{T}(\s)^{\ast}\mathbf{H}^{\N}(\s,\t,\x,\y). \label{eq:hke00}
\end{align}
Now, for any $\t\in\R$, set $\mathscr{T}^{!}(\t)$ as $\mathscr{T}(\t)$ from {Definition \ref{definition:method1}} but without the delta-function term:
\begin{align}
\mathscr{T}^{!}(\t) \ := \ \N^{2}\bar{\alpha}(\t)\Delta + \tfrac14\N\lambda(\t)^{2}\bar{\alpha}(\t)\Delta+\N^{\frac32}\bar{\alpha}(\t)\grad^{+}-\N^{\frac32}\bar{\alpha}(\t)\grad^{-}. 
\end{align}
Let $\mathbf{H}^{\N,!}(\s,\t,\x,\y)$ solve $\partial_{\t}\mathbf{H}^{\N,!}(\s,\t,\x,\y)=\mathscr{T}^{!}(\t)\mathbf{H}^{\N,!}(\s,\t,\x,\y)$ and $\mathbf{H}^{\N,!}(\s,\s,\x,\y)=\mathbf{1}[\x=\y]$, so $\mathbf{H}^{\N,!}$ is the heat kernel for $\partial_{\t}-\mathscr{T}^{!}(\t)$. (Here, $\mathscr{T}^{!}(\t)$ acts on the $\x$-variable.) We first claim the following relation between $\mathbf{H}^{\N}$ and $\mathbf{H}^{\N,!}$ heat kernels:
\begin{align}
\mathbf{H}^{\N}(\s,\t,\x,\y) \ = \ \mathbf{H}^{\N,!}(\s,\t,\x^{\s,\t},\y), \quad\mathrm{where}\quad \x^{\s,\t} \ := \ \x-\lfloor2\N^{\frac32}{\textstyle\int_{\s}^{\t}}\bar{\alpha}(\tau)\d\tau\rfloor. \label{eq:hke1}
\end{align}
In particular, the $\mathbf{H}^{\N}(\s,\t,\x,\y)$ depends only on $\x-\y$, and thus $\mathbf{H}^{\N}(\s,\t,\x)$ is a convolution operator. We now give some bounds for the $\mathbf{H}^{\N,!}$ heat kernel. Take any $\s\leq\t$, and fix a positive integer $\mathrm{m}\geq1$. Recall the length-$\mathfrak{l}$ gradient $\grad^{\mathfrak{l}}$ (for $\mathfrak{l}\in\Z$), and take $\mathfrak{l}_{1},\ldots,\mathfrak{l}_{\mathrm{m}}\in\Z$. We have the following, in which the product of $\grad^{\mathfrak{l}_{\mathrm{k}}}$-operators means their composition:
\begin{align}
|[{\textstyle\prod_{\mathrm{k}=1}^{\mathrm{m}}}\N\grad^{\mathfrak{l}_{\mathrm{k}}}]\mathbf{H}^{\N,!}(\s,\t,\x,\y)| \ \ \lesssim \ \N^{-1}|\t-\s|^{-\frac12-\frac{\mathrm{m}}{2}}{\textstyle\prod_{\mathrm{k}=1}^{\mathrm{m}}}|\mathfrak{l}_{\mathrm{k}}|. \label{eq:hke2}
\end{align}
The same holds for $\mathrm{m}=0$, upon removing all products over $\mathrm{k}=1,\ldots,\mathrm{m}$ in \eqref{eq:hke2}. Next, take any $\phi:\mathbb{T}(\N)\to\R$. In the same setting as \eqref{eq:hke2} (including the case $\mathrm{m}=0$), we have the following, where the sum and sup are both over $\y\in\mathbb{T}(\N)$:
\begin{align}
{\textstyle\sum_{\y}}|[{\textstyle\prod_{\mathrm{k}=1}^{\mathrm{m}}}\N\grad^{\mathfrak{l}_{\mathrm{k}}}]\mathbf{H}^{\N,!}(\s,\t,\x,\y)||\phi_{\y}| \ \lesssim \ |\t-\s|^{-\frac{\mathrm{m}}{2}}{\textstyle\prod_{\mathrm{k}=1}^{\mathrm{m}}}|\mathfrak{l}_{\mathrm{k}}|\times{\textstyle\sup_{\y}}|\phi_{\y}|. \label{eq:hke3}
\end{align}
We claim {that} \eqref{eq:hke2}-\eqref{eq:hke3} are true as written for $\mathbf{H}^{\N}$ in place of $\mathbf{H}^{\N,!}$. In particular, the operator norm of $[{\textstyle\prod_{\mathrm{k}=1}^{\mathrm{m}}}\N\grad^{\mathfrak{l}_{\mathrm{i}}}]\mathbf{H}^{\N}(\s,\t,\x)$ (on $\mathscr{L}^{\infty}(\mathbb{T}(\N))\to\mathscr{L}^{\infty}(\mathbb{T}(\N))$) is $\lesssim\mathrm{RHS}\eqref{eq:hke3}$. We now give time-regularity for $\mathbf{H}^{\N}$. Fix $\tau\geq0$ and $0\leq\s\leq\t\lesssim1$. We claim
\begin{align}
{\textstyle\sum_{\y}}|\mathbf{H}^{\N}(\s,\t+\tau,\x,\y)-\mathbf{H}^{\N}(\s,\t,\x,\y)||\phi_{\y}| \ &\lesssim \ |\t-\s|^{-1}[\tau+\N^{-1}]\times{\textstyle\sup_{\y}}|\phi_{\y}|\label{eq:hke4}\\
|\mathbf{H}^{\N}(\s,\t+\tau,\x,\y)-\mathbf{H}^{\N}(\s,\t,\x,\y)| \ &\lesssim \ \N^{-1}|\t-\s|^{-\frac32}[\tau+\N^{-1}]. \label{eq:hke5}
\end{align}
Next, take any $0\leq p\leq2$. We now claim the following moment estimate:
\begin{align}
{\textstyle\sum_{\y}}\mathbf{H}^{\N}(\s,\t,\x,\y)\times\N^{-p}|\x-\y|^{p} \ \lesssim_{p} \ |\t-\s|^{\frac{p}{2}}+\N^{-p}. \label{eq:hke6}
\end{align}
\end{prop}
\begin{proof}
We first show \eqref{eq:hke0}. Fix $\s\leq\r\leq\t$, and fix $\x,\y\in\mathbb{T}(\N)$. It suffices to show that the kernels of $\mathrm{LHS}\eqref{eq:hke0}$ and $\mathrm{RHS}\eqref{eq:hke0}$ are equal. In particular, it suffices to prove the following ``matrix multiplication" identity, where the sum is over all $\z\in\mathbb{T}(\N)$:
\begin{align}
{\textstyle\sum_{\z}}\mathbf{H}^{\N}(\r,\t,\x,\z)\mathbf{H}^{\N}(\s,\r,\z,\y) \ = \ \mathbf{H}^{\N}(\s,\t,\x,\y). \label{eq:hke01}
\end{align}
Both sides of the proposed identity clearly vanish under $\partial_{\t}-\mathscr{T}(\t)$, because, by construction, the LHS is a linear combination of terms that do, and the RHS does. Both sides are also equal at $\t=\r$. So, by standard uniqueness for linear PDEs, both sides are equal for all $\t\geq\r$. Next, let us show \eqref{eq:hke00}. Differentiate \eqref{eq:hke01} in $\r$. As $\mathrm{RHS}\eqref{eq:hke01}$ is independent of $\r$, the Leibniz rule gives
\begin{align}
{\textstyle\sum_{\z}}\partial_{\r}\mathbf{H}^{\N}(\r,\t,\x,\z)\times\mathbf{H}^{\N}(\s,\r,\z,\y) + {\textstyle\sum_{\z}}\mathbf{H}^{\N}(\r,\t,\x,\z)\times\mathscr{T}(\r)\mathbf{H}^{\N}(\s,\r,\z,\y) \ = \ 0. \label{eq:hke001}
\end{align}
We now replace $\mathscr{T}(\r)$ in $\mathrm{LHS}\eqref{eq:hke001}$ by $\mathscr{T}(\r)^{\ast}$ and let it instead act on $\mathbf{H}^{\N}(\r,\t,\x,\z)$ in $\z$. Then, take $\r\to\s$ from above. Since $\mathbf{H}^{\N}(\s,\s,\z,\y)=\mathbf{1}[\z=\y]$, we deduce $\partial_{\s}\mathbf{H}^{\N}(\s,\t,\x,\y)+\mathscr{T}(\s)^{\ast}\mathbf{H}^{\N}(\s,\t,\x,\y)=0$, which implies \eqref{eq:hke00}. (The only subtlety here is if $\s\in\mathbb{J}$, in which case $\mathscr{T}(\s)^{\ast}$ has a discrete gradient that $\mathscr{T}(\r)^{\ast}$ does not see as we take $\r\to\s$ from above. This is compensated for{,} since the difference between $\partial_{\s}\mathbf{H}^{\N}(\s,\t,\x,\y)$ and $\partial_{\r}\mathbf{H}^{\N}(\r,\t,\x,\y)$ as $\r\to\s$ from above is the same discrete gradient, coming from the fact that $\partial_{\s}$ in $\partial_{\s}\mathbf{H}^{\N}(\s,\t,\x,\y)$ must act on $\x^{\s,\t}$ in \eqref{eq:hke1}, which is shown shortly and does not need \eqref{eq:hke00}. In particular, letting $\partial_{\s}$ act on $\x^{\s,\t}$ in $\mathrm{RHS}\eqref{eq:hke1}$ produces $\delta(\s\in\mathbb{J})\grad^{-}$, which matches the additional $\delta(\s\in\mathbb{J})\grad^{-}$ we get when we replace the limit of $\mathscr{T}(\r)^{\ast}$ as $\r\to\s$ from above by $\mathscr{T}(\s)^{\ast}$ itself.) This concludes the proof of \eqref{eq:hke00}.

We now prove \eqref{eq:hke1}. Take the following random walk $\tau\mapsto\mathfrak{X}(\tau)$. It jumps according to Poisson clocks that have generator given by the first four terms in $\mathrm{RHS}\eqref{eq:method1IT}$. It also jumps, independently, to the left by 1 unit at the deterministic set of times $\mathbb{J}$. As noted after Definition \ref{definition:method1}, $\mathbf{H}^{\N}(\s,\t,\x,\y)$ is the probability that $\mathfrak{X}(\tau)$ goes from $\x$ at $\tau=\s$ to $\y$ at $\tau=\t$. Because the speed of the jumps in $\mathfrak{X}$ are independent of the position of $\mathfrak{X}$, we deduce that $\mathfrak{X}$ is the same as following a random walk whose generator equals the first four terms in $\mathrm{RHS}\eqref{eq:method1IT}$, and then adding deterministic jumps to the left by 1 unit for every time in $\mathbb{J}$. Using this, we deduce {that} the probability that $\mathfrak{X}(\tau)$ goes from $\x$ at $\tau=\s$ to $\y$ at $\tau=\t$ is $\mathrm{RHS}\eqref{eq:hke1}$, so \eqref{eq:hke1} follows.  For the sentence after \eqref{eq:hke1}, it suffices to note {that} $\mathbf{H}^{\N,!}(\s,\t,\w,\z)$ depends only on $\w-\z$ (its generator is space-homogeneous), and $\x^{\s,\t}-\x$ is independent of $\x$. We move to \eqref{eq:hke2}. Set $\mathscr{T}^{!}(\t)=\mathscr{T}^{!,1}(\t)+\mathscr{T}^{!,2}(\t)$ for $\mathscr{T}^{!,2}(\t)=2\N^{3/2}\bar{\alpha}(\t)\grad^{+}$ and $\mathscr{T}^{!,1}(\t):=[\N^{2}-\N^{3/2}+4^{-1}\N\lambda(\t)^{2}]\bar{\alpha}(\t)\Delta$. Since $\mathscr{T}^{!,2}(\t)$ commutes with $\mathscr{T}^{!,1}(\t)$ (as they are both in the algebra generated by commuting operators $\grad^{+}$ and $\grad^{-}$), we get
\begin{align}
\mathbf{H}^{\N,!}(\s,\t,\cdot) \ = \ \exp[{\textstyle\int_{\s}^{\t}}\mathscr{T}^{!}(\tau)\d\tau] \ = \ \exp[{\textstyle\int_{\s}^{\t}}\mathscr{T}^{!,2}(\tau)\d\tau]\exp[{\textstyle\int_{\s}^{\t}}\mathscr{T}^{!,1}(\tau)\d\tau]. \label{eq:hkebound}
\end{align}
Since $\mathscr{T}^{!,1}(\tau)$ and $\mathscr{T}^{!,2}(\tau)$ are spatially homogeneous discrete differentials, the kernel $\mathbf{H}^{\N,!}(\s,\t,\x,\y)$ for the LHS is the spatial convolution on $\mathbb{T}(\N)$ between the kernels for the exponentials on the far RHS. But the $\mathscr{T}^{!,2}(\tau)$-semigroup is uniformly bounded as an operator $\mathscr{L}^{p}(\mathbb{T}(\N))\to\mathscr{L}^{p}(\mathbb{T}(\N))$ for any $p\in[1,\infty]$, since it is the semigroup for a totally asymmetric random walk. Since convolution on $\mathbb{T}(\N)$ commutes with $\grad$-operators, it suffices to assume {that} $\mathbf{H}^{\N,!}$ is instead the kernel for $\mathscr{T}^{!,1}(\tau)$ when showing the bounds \eqref{eq:hke2}-\eqref{eq:hke3}. To this end, now let $\mathbf{H}^{\N,\mathrm{line}}$ be the heat kernel satisfying $\partial_{\t}\mathbf{H}^{\N,\mathrm{line}}(\s,\t,\z,\w)=\mathscr{T}^{!,\mathrm{line}}(\t)\mathbf{H}^{\N,\mathrm{line}}(\s,\t,\z,\w)$ and $\mathbf{H}^{\N,\mathrm{line}}(\s,\s,\z,\w)=\mathbf{1}[\z=\w]$ for $\s,\t\in\R$ and $\z,\w\in\Z$, where $\mathscr{T}^{!,\mathrm{line}}(\t)$ is just $\mathscr{T}^{!,1}(\t)$ but replacing $\Delta$ {by} the discrete Laplacian on the full line $\Z$. It is standard that $\mathbf{H}^{\N,!}(\s,\t,\x,\y)$ equals the sum over all $\mathrm{k}\in\Z$ of $\mathbf{H}^{\N,\mathrm{line}}(\s,\t,\x,\y+\mathrm{k}|\mathbb{T}(\N)|)$. At this point, \eqref{eq:hke2}-\eqref{eq:hke3} now follow by sub-exponentially decaying regularity bounds in Proposition A.1 and Corollary A.2 of \cite{DT}, whose extension to higher-order derivatives follows by taking more differentials and then doing the exact same analysis. (Technically, \cite{DT} deals with the case of $\bar{\alpha}(\t)=\bar{\alpha}(0)$. However, the formulas and analysis therein for the time $\s\mapsto\t$ heat kernel hold if we replace the Laplacian coefficient with the integral between $\s$ and $\t$ of $\bar{\alpha}(\tau)\d\tau$.) The fact that \eqref{eq:hke2}-\eqref{eq:hke3} hold for $\mathbf{H}^{\N}$ in place of $\mathbf{H}^{\N,!}$ follows immediately by \eqref{eq:hke1} and noting {that} \eqref{eq:hke2}-\eqref{eq:hke3} are uniform in $\x$. 

Next, instead of proving \eqref{eq:hke4}-\eqref{eq:hke5}, which is a somewhat involved argument, we now prove \eqref{eq:hke6}, which is much simpler. (We defer \eqref{eq:hke4}-\eqref{eq:hke5} to the end of this proof.) \eqref{eq:hke3} for $\phi\equiv1$ and $\mathrm{m}=0$ gives \eqref{eq:hke6} for $p=0$. So by interpolation, it suffices to assume $p=2$ in \eqref{eq:hke6}. By \eqref{eq:hke1}, $\mathbf{H}^{\N}(\s,\t,\x,\y)$ is the transition probability of $\x\rightsquigarrow\y$ from times $\s$ to $\t$ of the random walk $\mathfrak{X}(\tau)$ given by the random walk $\mathfrak{X}^{!}$ defined by $\mathbf{H}^{\N,!}$, but for every $\tau$ in the jump set $\mathbb{J}$, shift $\mathfrak{X}^{!}\mapsto\mathfrak{X}^{!}-1$. So $\mathfrak{X}(\tau)=\mathfrak{X}(\tau;1)+\mathfrak{X}(\tau;2)$, where $\mathfrak{X}(\tau;1)$ is a symmetric simple random walk on $\mathbb{T}(\N)$ of speed $\lesssim\N^{2}$, and $\mathfrak{X}(\tau;2)$ is a totally asymmetric random walk (to the right) of speed $2\N^{3/2}\bar{\alpha}(\tau)$ that is then pushed to the left by 1 unit at every jump time in $\mathbb{J}$. (The generator of $\mathfrak{X}(\tau;1)$ equals $\mathscr{T}^{!,1}(\tau)$. The generator of $\mathfrak{X}(\tau;2)$ is $\mathscr{T}^{!,2}(\tau)$ plus the delta function in \eqref{eq:method1IT} at time $\tau$.) To get \eqref{eq:hke6}, it suffices to show {that} the second moment of $\mathfrak{X}(\t)-\mathfrak{X}(\s)$ is $\lesssim\N^{2}|\t-\s|+1$. The second moment of $\mathfrak{X}(\t;1)-\mathfrak{X}(\s;1)$ is $\lesssim\N^{2}|\t-\s|$ by standard martingale bounds ($\mathfrak{X}(\tau;1)$ is the symmetric walk). Next, note {that} $\mathfrak{X}(\t;2)-\mathfrak{X}(\s;2)$ is a difference of a Poisson variable of speed $\lesssim\N^{3/2}|\t-\s|$ and its mean (plus $\mathrm{O}(1)$, as $[\s,\t]$ may not ``intersect exactly" with $\mathbb{J}$). Therefore, its second moment satisfies $\lesssim\N^{3/2}|\t-\s|+1$. As $|\mathfrak{X}(\t)-\mathfrak{X}(\s)|^{2}\lesssim|\mathfrak{X}(\t;1)-\mathfrak{X}(\s;1)|^{2}+|\mathfrak{X}(\t;2)-\mathfrak{X}(\s;2)|^{2}$, \eqref{eq:hke6} holds. We move to \eqref{eq:hke4}-\eqref{eq:hke5}. Assume {that \eqref{eq:hke4} holds} for $\phi\equiv1$. Set $\r=2^{-1}(\t+\s)$. By \eqref{eq:hke0} and the pointwise estimate \eqref{eq:hke2} for $\mathrm{m}=0$ and $\mathbf{H}^{\N}$ instead of $\mathbf{H}^{\N,!}$, we have
\begin{align}
\mathrm{LHS}\eqref{eq:hke5} \ &\leq \ {\textstyle\sum_{\z}}|\mathbf{H}^{\N}(\r,\t+\tau,\x,\z)-\mathbf{H}^{\N}(\r,\t,\x,\z)|\mathbf{H}^{\N}(\s,\r,\z,\y) \nonumber\\
&\lesssim \ \N^{-1}|\r-\s|^{-\frac12}|\t-\r|^{-1}[\tau+\N^{-1}] \ \lesssim \ \mathrm{RHS}\eqref{eq:hke5}, \nonumber
\end{align}
so \eqref{eq:hke5} holds. Also, for any general $\phi:\mathbb{T}(\N)\to\R$, we again use \eqref{eq:hke0} and claim the following, where sums are all over $\mathbb{T}(\N)$:
\begin{align}
\mathrm{LHS}\eqref{eq:hke4} \ \lesssim \ {\textstyle\sum_{\y}}{\textstyle\sum_{\z}}|\mathbf{H}^{\N}(\r,\t+\tau,\x,\z)-\mathbf{H}^{\N}(\r,\t,\x,\z)|\mathbf{H}^{\N}(\s,\r,\z,\y)|\phi_{\y}| \ \lesssim \ \mathrm{RHS}\eqref{eq:hke4},
\end{align}
where $\r=2^{-1}(\t+\s)$ as before. (The second bound follows from bounding the $\y$-sum via {the} contractivity of $\mathbf{H}^{\N}$, namely \eqref{eq:hke3} for $\mathrm{m}=0$ for $\mathbf{H}^{\N}$ in place of $\mathbf{H}^{\N,!}$, and then bounding the $\z$-sum by our assumption that \eqref{eq:hke4} holds for $\phi\equiv1$.) So, to complete the proof of this proposition, we are left to show \eqref{eq:hke4} for $\phi\equiv1$. First, {we give} a few preliminaries. Recall $\mathfrak{X}(\tau;2)$; it is the sum of a totally asymmetric simple random walk of speed $2\N^{3/2}\bar{\alpha}(\tau)$ to the right with jumps at times in $\mathbb{J}$ to the left by one unit, and $\mathfrak{X}(\tau;1)$ is a symmetric simple random walk of speed $\lesssim\N^{2}$ and $\gtrsim\N^{2}$ (with different implied constants) with generator $\mathscr{T}^{!,1}(\tau)$. Let $\mathbf{H}^{\N,\Delta}$ solve the PDE $\partial_{\t}\mathbf{H}^{\N,\Delta}(\s,\t,\x,\y)=\mathscr{T}^{!,1}(\t)\mathbf{H}^{\N,\Delta}(\s,\t,\x,\y)$ and $\mathbf{H}^{\N,\Delta}(\s,\s,\x,\y)=\mathbf{1}[\x=\y]$. Again, $\mathscr{T}^{!,1}$ is from \eqref{eq:hkebound} and it acts on $\x$. Next, let $\E^{\x,\s}$ denote {the} expectation with respect to the walk $\mathfrak{X}(\tau;2)$ starting at time $\s$ and position $\x$. We now claim {that}
\begin{align}
\mathbf{H}^{\N}(\s,\t,\x,\y) \ = \ \E^{\x,\s}\mathbf{H}^{\N,\Delta}(\s,\t,\mathfrak{X}(\t;2),\y). \label{eq:walkhk}
\end{align}
Indeed, \eqref{eq:walkhk} follows from first noting that \eqref{eq:hkebound} holds if we drop $!$ on the LHS and replace the first exponential on the far RHS by the semigroup for $\mathfrak{X}(\tau;2)$, which gives the $\E^{\x,\s}$ operator in $\mathrm{RHS}\eqref{eq:walkhk}$. We then match $(\x,\y)$-entries of the LHS and RHS of the resulting equation to get \eqref{eq:walkhk}. Via \eqref{eq:walkhk}, for any $\tau\geq0$, we have $\mathrm{LHS}\eqref{eq:hke4}|_{\phi\equiv1}\lesssim\Upsilon^{\tau,1}[\s,\t,\x]+\Upsilon^{\tau,2}[\s,\t,\x]$, where
\begin{align}
\Upsilon^{\tau,1}[\s,\t,\x] \ &:= \ {\textstyle\sum_{\y}}\E^{\x,\s}|\mathbf{H}^{\N,\Delta}[\s,\t+\tau,\mathfrak{X}(\t+\tau;2),\y]-\mathbf{H}^{\N,\Delta}[\s,\t,\mathfrak{X}(\t+\tau;2),\y]| \nonumber\\
\Upsilon^{\tau,2}[\s,\t,\x] \ &:= \ {\textstyle\sum_{\y}}|\E^{\x,\s}\mathbf{H}^{\N,\Delta}[\s,\t,\mathfrak{X}(\t+\tau;2),\y]-\E^{\x,\s}\mathbf{H}^{\N,\Delta}[\s,\t,\mathfrak{X}(\t;2),\y]|. \nonumber
\end{align}
We first control $\Upsilon^{\tau,1}$ above. As argued after \eqref{eq:hkebound}, we know {that} {\small$\N^{2}\sum_{\y}|\Delta\mathbf{H}^{\N,\Delta}(\s,\t,\x,\y)|\lesssim|\t-\s|^{-1}$}. By the fundamental theorem of calculus and the PDE for {\small$\mathbf{H}^{\N,\Delta}$}, we deduce {that} {\small$|\Upsilon^{\tau,1}[\s,\t,\x]|\lesssim|\t-\s|^{-1}\tau\lesssim\mathrm{RHS}\eqref{eq:hke4}|_{\phi\equiv1}$}. Thus, we are left to show {that}
\begin{align}
|\Upsilon^{\tau,2}[\s,\t,\x]| \ \lesssim \ \mathrm{RHS}\eqref{eq:hke4}|_{\phi\equiv1}. \label{eq:whatislefthke}
\end{align}
Let us first couple $\mathfrak{X}(\t+\tau;2)$ and $\mathfrak{X}(\t;2)$ in the definition of $\Upsilon^{\tau,2}$. These are the same $\mathfrak{X}(\cdot;2)$ random walk, which starts at time $\s$ and position $\x$, evaluated at time $\t+\tau$ and $\t$, respectively. Thus, we know $\mathfrak{X}(\t+\tau;2)=\mathfrak{X}(\t;2)+\mathfrak{Y}(\t,\tau;2)$, where $\mathfrak{Y}(\t,\tau;2)$ is a random walk starting from position 0 and $\tau=0$ with the same (space homogeneous) dynamics as $\mathfrak{X}(\t+\cdot;2)$. (We emphasize the increment $\mathfrak{Y}(\t,\tau;2)$ is independent of $\mathfrak{X}(\t;2)$.) As noted two paragraphs after \eqref{eq:hkebound}, $\mathfrak{Y}(\t,\tau;2)$ is $\mathrm{O}(1)$ plus the difference of a Poisson variable of speed $\lesssim\N^{3/2}\tau$ and its mean. So, $\E^{\x,\s,\t}\mathfrak{Y}(\t,\tau;2)=\mathrm{O}(1)$ and $\E^{\x,\s}|\mathfrak{Y}(\t,\tau;2)|^{2}\lesssim\N^{3/2}\tau+\mathrm{O}(1)$, where $\E^{\x,\s,\t}$ is $\E^{\x,\s}$ but further conditioning on $\mathfrak{X}(\t;2)$. We now claim (with explanation given afterwards) {that}
{\small
\begin{align}
&\mathbf{H}^{\N,\Delta}[\s,\t,\mathfrak{X}(\t+\tau;2),\y]-\mathbf{H}^{\N,\Delta}[\s,\t,\mathfrak{X}(\t;2),\y] \nonumber\\
&= \ \mathbf{H}^{\N,\Delta}[\s,\t,\mathfrak{X}(\t;2)+\mathfrak{Y}(\t,\tau;2),\y]-\mathbf{H}^{\N,\Delta}[\s,\t,\mathfrak{X}(\t;2),\y]\nonumber\\
&= \ \mathfrak{Y}(\t,\tau;2)\grad^{+}\mathbf{H}^{\N,\Delta}[\s,\t,\mathfrak{X}(\t;2),\y] \label{eq:whatislefthke1}\\
&+ {\textstyle\sum_{\w=0}^{\mathfrak{Y}(\t,\tau;2)-1}}\{\grad^{+}\mathbf{H}^{\N,\Delta}[\s,\t,\w+\mathfrak{X}(\t;2),\y]-\grad^{+}\mathbf{H}^{\N,\Delta}[\s,\t,\mathfrak{X}(\t;2),\y]\}.\nonumber
\end{align}
}The first identity follows from $\mathfrak{X}(\t+\tau;2)=\mathfrak{X}(\t;2)+\mathfrak{Y}(\t,\tau;2)$. The second identity can be checked as a ``discrete Taylor expansion". (Rewrite the difference in the RHS of the first line as a telescoping sum of discrete gradients $\grad^{+}\mathbf{H}^{\N,\Delta}[\s,\t,\w+\mathfrak{X}(\t;2),\y]$ from $\w=0$ to $\w=\mathfrak{Y}(\t,\tau;2)-1$. {Afterwards}, replace $\grad^{+}\mathbf{H}^{\N,\Delta}[\s,\t,\w+\mathfrak{X}(\t;2),\y]$ by $\grad^{+}\mathbf{H}^{\N,\Delta}[\s,\t,\mathfrak{X}(\t;2),\y]$ to obtain the first term in \eqref{eq:whatislefthke1} with an error given by the sum in \eqref{eq:whatislefthke1}.) We now claim the following (with explanation afterwards):
\begin{align}
&{\textstyle\sum_{\y}}|\E^{\x,\s}\mathfrak{Y}(\t,\tau;2)\grad^{+}\mathbf{H}^{\N,\Delta}[\s,\t,\mathfrak{X}(\t;2),\y]| \nonumber\\
&= \ {\textstyle\sum_{\y}}|\E^{\x,\s}\grad^{+}\mathbf{H}^{\N,\Delta}[\s,\t,\mathfrak{X}(\t;2),\y]\E^{\x,\s,\t}\mathfrak{Y}(\t,\tau;2)|\label{eq:whatislefthke2} \\
&\lesssim \ {\textstyle\sum_{\y}}\E^{\x,\s}|\grad^{+}\mathbf{H}^{\N,\Delta}[\s,\t,\mathfrak{X}(\t;2),\y]| \ \lesssim \ \N^{-1}|\t-\s|^{-\frac12} \ \lesssim \ \N^{-1}|\t-\s|^{-1}.\label{eq:whatislefthke3}
\end{align}
The first line follows by law of total expectation, and noting that $\grad^{+}\mathbf{H}^{\N,\Delta}[\s,\t,\mathfrak{X}(\t;2),\y]$ is deterministic once we condition on $\mathfrak{X}(\t;2)$ in $\E^{\x,\s,\t}$. {In order to get the second line, first observe that} $|\E^{\x,\s,\t}\mathfrak{Y}(\t,\tau;2)|\lesssim1$. {Afterwards, we use} \eqref{eq:hke3} for $\mathbf{H}^{\N,\Delta}$ instead of $\mathbf{H}^{\N,!}$. (Indeed, as argued after \eqref{eq:hkebound}, such a bound is true.) Finally, we use the fact that $\s\leq\t\lesssim1$, which lets us adjust the exponent from $-1/2$ to $-1$ up to a $\mathrm{O}(1)$ factor. We now claim the following estimate:
\begin{align}
&{\textstyle\sum_{\y}}{\textstyle\sum_{\w=0}^{\mathfrak{Y}(\t,\tau;2)-1}}|\grad^{+}\mathbf{H}^{\N,\Delta}[\s,\t,\w+\mathfrak{X}(\t;2),\y]-\grad^{+}\mathbf{H}^{\N,\Delta}[\s,\t,\mathfrak{X}(\t;2),\y]| \nonumber\\
&\lesssim \ \N^{-2}|\t-\s|^{-1}\times|\mathfrak{Y}(\t,\tau;2)|^{2}. \label{eq:whatislefthke4}
\end{align}
Indeed, pull the $\y$-sum in and, again, use \eqref{eq:hke3} for $\mathbf{H}^{\N,\Delta}$ instead of $\mathbf{H}^{\N,!}$. Now, by \eqref{eq:whatislefthke4} and $\E^{\x,\s}|\mathfrak{Y}(\t,\tau;2)|^{2}\lesssim\N^{3/2}\tau+\mathrm{O}(1)$,
\begin{align}
|\Upsilon^{\tau,2}[\s,\t,\x]| \ \lesssim \ \N^{-2}|\t-\s|^{-1}\E^{\x,\s}|\mathfrak{Y}(\t,\tau;2)|^{2} \ &\lesssim \ \N^{-\frac12}|\t-\s|^{-1}\tau+\N^{-2}|\t-\s|^{-1} \\
&\lesssim \ \mathrm{RHS}\eqref{eq:hke4}|_{\phi\equiv1}.
\end{align}
Therefore, \eqref{eq:whatislefthke} follows. As we noted right before \eqref{eq:whatislefthke}, the proof is now complete.
\end{proof}
The last result of this section (Proposition \ref{prop:hkecont}) presents similar regularity bounds for the continuum $\mathbf{H}$ kernel. It also gives estimates that show $\mathbf{H}^{\N}\approx\mathbf{H}$ (in some sense after suitably rescaling). These bounds are \eqref{eq:hkecont2}-\eqref{eq:hkecont3} below, which, according to Lemma 3.2 in \cite{G}, are the types of bounds that we need to combine with the $\mathbf{H}^{\N}$ estimates in Proposition \ref{prop:hke} in order to show convergence of stochastic heat equations in the proof of Proposition \ref{prop:method11}.
\begin{prop}\label{prop:hkecont}
 Recall $\mathbf{H}$ from {Definition \ref{definition:intro2}} and {the} $\mathbf{H}^{\N}$-terms from {Definition \ref{definition:method1}}. {For any $\s\leq\t$, $\x,\y\in\mathbb{T}$, $\d\geq0$, it holds that}
\begin{align}
|\t-\s|^{\frac12+\frac{\d}{2}}|\partial_{\x}^{\d}\mathbf{H}(\s,\t,\x,\y)|+|\t-\s|^{\frac12+\d}|\partial_{\t}^{\d}\mathbf{H}(\s,\t,\x,\y)| \ \lesssim \ 1. \label{eq:hkecont1}
\end{align}
The same is true if we integrate $\mathrm{LHS}\eqref{eq:hkecont1}$ over $\y\in\mathbb{T}$. Now, fix $0\leq\s\leq\t\lesssim1$ and $\x\in\mathbb{T}$. {Then, there} exist $\gamma,\beta\in(0,1)$ such that
\begin{align}
&{\textstyle\int_{\s}^{\t}\int_{\mathbb{T}}}|\mathbf{H}(\s,\t,\x,\y)-\N\mathbf{H}^{\N}(\s,\t,\N\x,\N\y)|^{2}\d\y\d\t \ \lesssim_{|\t-\s|} \ \N^{-\gamma}\label{eq:hkecont2}\\
&{\textstyle\int_{\mathbb{T}}}|\mathbf{H}(\s,\t,\x,\y)-\N\mathbf{H}^{\N}(\s,\t,\N\x,\N\y)|^{2}\d\y \ \lesssim \ \N^{-\gamma}|\t-\s|^{-\beta}. \label{eq:hkecont3}
\end{align}
($\mathbf{H}^{\N}$ extends from $\mathbb{T}(\N)\times\mathbb{T}(\N)$ to $\N\mathbb{T}\times\N\mathbb{T}$ as a piecewise constant. The implied constant in \eqref{eq:hkecont2} is continuous in $|\t-\s|$.)
\end{prop}
\begin{proof}
Let $\mathbf{H}^{\mathrm{line}}(\s,\t,\x,\y)$ be the Gaussian kernel for $\x,\y\in\R$ and $\s\leq\t$ with mean zero and variance {\small$\int_{\s}^{\t}\bar{\alpha}(\tau)\d\tau$.} Because the torus $\mathbb{T}$ is the quotient $\R/|\mathbb{T}|\Z$, the method of images gives the representation  
\begin{align}
\mathbf{H}(\s,\t,\x,\y) \ = \ {\textstyle\sum_{\mathrm{k}}}\mathbf{H}^{\mathrm{line}}(\s,\t,\x,\y+\mathrm{k}|\mathbb{T}|),
\end{align}
where the sum is over $\mathrm{k}\in\Z$. At this point, \eqref{eq:hkecont1} now follows from standard regularity calculations and estimates for Gaussian kernels. \eqref{eq:hkecont1}, but integrating $\mathrm{LHS}\eqref{eq:hkecont1}$ over $\y\in\mathbb{T}$, follows because $|\mathbb{T}|\lesssim1$. We now show \eqref{eq:hkecont2}-\eqref{eq:hkecont3}. Recall the $\mathbf{H}^{\N,\Delta}$ kernel from before \eqref{eq:walkhk}. If we replace $\mathbf{H}^{\N}$ in \eqref{eq:hkecont2}-\eqref{eq:hkecont3} by $\mathbf{H}^{\N,\Delta}$, then the resulting bound would hold by the exact same argument as the proof of Lemma 3.2 in \cite{G}. (This argument is based on spectral theory for the Laplacian and discrete Laplacian on an interval with Dirichlet boundary conditions. But, eigenvalues are also exactly computable and have the same asymptotics in the current case of periodic boundary conditions.) Thus, it suffices to control $\mathbf{H}^{\N}-\mathbf{H}^{\N,\Delta}$. In particular, it is left to show
\begin{align}
&{\textstyle\int_{\s}^{\t}\int_{\mathbb{T}}}|\N\mathbf{H}^{\N}(\s,\t,\N\x,\N\y)-\N\mathbf{H}^{\N,\Delta}(\s,\t,\N\x,\N\y)|^{2}\d\y\d\t \ \lesssim_{|\t-\s|} \ \N^{-\gamma}\label{eq:hkecontp1}\\
&{\textstyle\int_{\mathbb{T}}}|\N\mathbf{H}^{\N}(\s,\t,\N\x,\N\y)-\N\mathbf{H}^{\N,\Delta}(\s,\t,\N\x,\N\y)|^{2}\d\y \ \lesssim \ \N^{-\gamma}|\t-\s|^{-\beta}. \label{eq:hkecontp2}
\end{align}
\eqref{eq:hkecontp1} follows by \eqref{eq:hkecontp2} and elementary integration (since $\beta\in(0,1)$ is independent of $\N$), so we focus on \eqref{eq:hkecontp2}. Because we have chosen piecewise constant extensions of $\mathbf{H}^{\N}$ and $\mathbf{H}^{\N,\Delta}$ from $\mathbb{T}(\N)\subseteq\N\mathbb{T}$, we deduce the following in $\w=\w(\x)\in\mathbb{T}(\N)$:
\begin{align}
\mathrm{LHS}\eqref{eq:hkecontp2} \ &= \ \N^{-1}{\textstyle\sum_{\z}}|\N\mathbf{H}^{\N}(\s,\t,\w,\z)-\N\mathbf{H}^{\N,\Delta}(\s,\t,\w,\z)|^{2}. \label{eq:hkecontp3}
\end{align}
First assume {that} $|\t-\s|\lesssim\N^{-3/2}$. Both $\mathbf{H}^{\N}$ and $\mathbf{H}^{\N,\Delta}$ are probability measures in $\y\in\mathbb{T}(\N)$. We also have the pointwise estimate \eqref{eq:hke2} for $\mathbf{H}^{\N}$ and $\mathbf{H}^{\N,\Delta}$ instead of $\mathbf{H}^{\N,!}$ (as noted right after \eqref{eq:hke3} and after \eqref{eq:hkebound}, respectively). Therefore, we have 
\begin{align}
\mathrm{RHS}\eqref{eq:hkecontp3} \ \lesssim \ \N\times{\textstyle\sup_{\z}}|\mathbf{H}^{\N}(\s,\t,\w,\z)-\mathbf{H}^{\N,\Delta}(\s,\t,\w,\z)| \ \lesssim \ |\t-\s|^{-\frac12} \ \lesssim \ \N^{-\gamma}|\t-\s|^{-\beta}{,}
\end{align}
given any $\beta\in(1/2,1)$, where $\gamma=\gamma(\beta)\gtrsim_{\beta}1$. The previous two displays give \eqref{eq:hkecontp2} in the case $|\t-\s|\lesssim\N^{-3/2}$, so it suffices to assume {that} $|\t-\s|\gtrsim\N^{-3/2}$. We use \eqref{eq:walkhk} and refer back to the paragraph preceding it for relevant notation. We claim this gives
\begin{align}
|\N\mathbf{H}^{\N}(\s,\t,\w,\z)-\N\mathbf{H}^{\N,\Delta}(\s,\t,\w,\z)| \ &\lesssim \ \N\E^{\w,\s}|\mathbf{H}^{\N,\Delta}[\s,\t,\mathfrak{X}(\t;2),\z]-\mathbf{H}^{\N,\Delta}[\s,\t,\w,\z]| \label{eq:hkecontp4a}\\
&\lesssim \ \N^{-1}|\t-\s|^{-1}\E^{\w,\s}|\mathfrak{X}(\t;2)-\w|\label{eq:hkecontp4b}\\
&\lesssim \ \N^{-\frac14}|\t-\s|^{-\frac12}+\N^{-1}|\t-\s|^{-1} \ \lesssim \ \N^{-\frac14}|\t-\s|^{-\frac12}. \label{eq:hkecontp4c}
\end{align}
\eqref{eq:hkecontp4a} is by \eqref{eq:walkhk} and then putting $\N\mathbf{H}^{\N,\Delta}(\s,\t,\w,\z)$ inside the expectation. \eqref{eq:hkecontp4b} is by the heat kernel gradient estimate \eqref{eq:hke2} for $\mathbf{H}^{\N,\Delta}$ in place of $\mathbf{H}^{\N,!}$ (which is a valid estimate as explained after \eqref{eq:hkebound}). \eqref{eq:hkecontp4c} follows because, as noted before \eqref{eq:whatislefthke1}, we know {that} $\mathfrak{X}(\t;2)-\w$ is $\mathrm{O}(1)$ plus a centered Poisson random variable of intensity $\N^{3/2}|\t-\s|$. This gives us the first bound in \eqref{eq:hkecontp4c}. The last bound follows by the assumption $|\t-\s|\gtrsim\N^{-3/2}$. We now interpolate the trivial bounds $\mathbf{H}^{\N}+\mathbf{H}^{\N,\Delta}\lesssim1$ with \eqref{eq:hkecontp4a}-\eqref{eq:hkecontp4c} to get the following estimate with $\e\in(0,1)$ small (again, in the case where $|\t-\s|\gtrsim\N^{-3/2}$):
\begin{align}
\mathrm{RHS}\eqref{eq:hkecontp3} \ \lesssim \ \N^{2\e}\N^{-\frac12+2\e}|\t-\s|^{-1+2\e} \ \lesssim \ \N^{-\frac14+4\e}|\t-\s|^{-1+2\e},
\end{align}
which equals $\mathrm{RHS}\eqref{eq:hkecontp2}$ for some $\gamma,\beta\in(0,1)$ if $\e\in(0,1/8)$. So, the previous bound and \eqref{eq:hkecontp3} give \eqref{eq:hkecontp2} in the remaining case $|\t-\s|\gtrsim\N^{-3/2}$. This proves \eqref{eq:hkecontp2} in general, so we are done as explained both prior to and after \eqref{eq:hkecontp1}-\eqref{eq:hkecontp2}.
\end{proof}
%
%
%
\section{Other technical results}
\subsection{Short-time continuity}
Throughout this paper, we often need to bootstrap from control of a space-time function on a very fine discretization of space-time to the entire continuous space-time. Lemma \ref{lemma:ste} gives a stochastic result based on controlling the \emph{very} short-time behavior of \eqref{eq:glsde} and \eqref{eq:glsdeloc}. Lemma \ref{lemma:steeasy} is a much simpler result for short-time continuity of smooth functions that happen to be stochastic. (In particular, Lemma \ref{lemma:steeasy} is completely separate from {the} continuity of stochastic processes.) Because Lemmas \ref{lemma:ste} and \ref{lemma:steeasy} are more or less intuitive exercises in (stochastic) calculus, we will not write out every detail. Also, before we start, let us adopt the following notation. We write $\mathrm{a}\asymp\mathrm{b}$ for $\mathrm{a},\mathrm{b}\geq0$ if $\mathrm{b}\lesssim\mathrm{a}\lesssim\mathrm{b}$ with different implied constants.
\begin{lemma}\label{lemma:ste}
 Fix {a} sufficiently large $\mathrm{D}$ that is independent of $\N$. Fix times $0\leq\mathfrak{t}_{1}\leq\mathfrak{t}_{2}\lesssim\N^{\mathrm{D}}$. Take any mesh scale $0\leq\mathfrak{n}\asymp\N^{-\mathrm{D}}$, and set $\mathbb{X}^{\mathfrak{n}}:=[\mathfrak{t}_{1},\mathfrak{t}_{2}]\cap\mathfrak{n}\Z$ to be a discretization of $[\mathfrak{t}_{1},\mathfrak{t}_{2}]$. Given any $\t\in[\mathfrak{t}_{1},\mathfrak{t}_{2}]$, we define $\t_{\circ}:=\max(\mathbb{X}^{\mathfrak{n}},(-\infty,\t])$ as the biggest point in $\mathbb{X}^{\mathfrak{n}}$ bounded {from} above by $\t$. {Then, we have} the following short-time continuity for \eqref{eq:glsde}, in which $\|\|$ is {the} supremum over variables $(\t,\x)\in[\mathfrak{t}_{1},\mathfrak{t}_{2}]\times\mathbb{T}(\N)$:
\begin{align}
\mathbb{P}[\|\mathbf{U}^{\t,\x}-\mathbf{U}^{\t_{\circ},\x}\| \ \gtrsim \ \N^{-\frac{\mathrm{D}}{2}}\{1+\|\mathbf{U}^{\t_{\circ},\x}\|\}] \ \lesssim \ \exp[-\N^{99}]. \label{eq:steI}
\end{align}
Consider any jointly smooth function $\phi:\R\times\R^{\mathbb{T}(\N)}\to\R$. Let $\llangle\rrangle$ be {the} sup-norm over $(\s,\mathbf{U})\in\R\times\R^{\mathbb{T}(\N)}$. Next, for any $\mathscr{B}>0$, define $\mathcal{E}[\phi,\mathscr{B}]$ to be the event where $|\phi|\gtrsim\mathscr{B}$ at $(\s,\mathbf{U}^{\t,\cdot})$ for some $\s,\t\in[\mathfrak{t}_{1},\mathfrak{t}_{2}]$. For any $\gamma,\mathscr{B}>0$, we have the following bound, in which $\grad_{\mathbf{U}}$ means gradient with respect to the $\mathbf{U}$-variable:
\begin{align}
\mathbb{P}\{\mathcal{E}[\phi,\mathscr{B}]\} \ \lesssim_{\gamma} \ &\exp[\N^{\gamma}]\times\sup_{\s,\t\in[\mathfrak{t}_{1},\mathfrak{t}_{2}]}\mathbb{P}[|\phi_{\s,\mathbf{U}^{\t,\cdot}}|\gtrsim\mathscr{B}-\{\llangle\partial_{\s}\phi\rrangle+\llangle\grad_{\mathbf{U}}\phi\rrangle\}\N^{-999}]\label{eq:steIIa}\\
+ \ &\exp[\N^{\gamma}]\times\sup_{\t\in[\mathfrak{t}_{1},\mathfrak{t}_{2}]}\sup_{\x\in\mathbb{T}(\N)}\mathbb{P}[|\mathbf{U}^{\t,\x}|\gtrsim\N^{999}]+\exp[-\N^{99}]. \label{eq:steIIb}
\end{align}
(In words, \eqref{eq:steIIa}-\eqref{eq:steIIb} uses calculus to extend $\phi$ from $(\s,\mathbf{U}^{\t,\cdot})$ for $\s,\t\in\mathbb{X}^{\mathfrak{n}}$ to $\s,\t\in[\mathfrak{t}_{1},\mathfrak{t}_{2}]$. The cost is control on $\|\mathbf{U}^{\t,\x}\|$, or by \eqref{eq:steI}, control on $\|\mathbf{U}^{\t_{\circ},\x}\|$.) We now specialize to \eqref{eq:hf}. Fix any integer $|\mathfrak{l}|\lesssim\N$, and recall $\|\|$ from right before \eqref{eq:steI}. We have
\begin{align}
\mathbb{P}[\|[\mathbf{h}(\t,\x)-\mathbf{h}(\t,\x+\mathfrak{l})]-[\mathbf{h}(\t_{\circ},\x)-\mathbf{h}(\t_{\circ},\x+\mathfrak{l})]\|\gtrsim\N^{-999}\{1+\|\mathbf{U}^{\t_{\circ},\x}\|\}] \ \lesssim \ \exp[-\N^{99}]. \label{eq:steIII}
\end{align}
Set $\mathcal{E}[\mathfrak{l}]$ as the event inside the probability in $\mathrm{LHS}\eqref{eq:steIII}$. By {a} union bound, the probability of the intersection of $\mathcal{E}[\mathfrak{l}]$ over $|\mathfrak{l}|\lesssim\N$ is $\lesssim\exp[-\N^{98}]$. Lastly, everything in this lemma holds if we replace $\mathbb{T}(\N)$ by any discrete interval $\mathbb{K}$ and if we replace \eqref{eq:hf}-\eqref{eq:glsde} by \eqref{eq:glsdeloc} and the process $\t\mapsto\mathbf{J}(\t,\cdot;\mathbb{K})$ from {Definition \ref{definition:le10}} (namely, with $\mathbb{I}(\mathfrak{t})$ therein equal to $\mathbb{K}$), respectively.
\end{lemma}
\begin{proof}
\eqref{eq:steI} is a standard estimate of Ito calculus along with the following observations. First, the diffusion coefficient in \eqref{eq:glsde} is constant, and the Brownian motions therein {have} speed $\lesssim\N^{2}$. Second, the drifts in \eqref{eq:glsde} are uniformly Lipschitz in the solution of Lipschitz norm $\lesssim\N^{2}$ by Assumption \ref{ass:intro8}. We now move to \eqref{eq:steIIa}-\eqref{eq:steIIb}. Just by calculus, we get {that} $\mathbb{P}\{\mathcal{E}[\phi,\mathscr{B}]\}\leq\mathbb{P}\{\mathcal{E}[\phi,\mathscr{B};2]\}$, where $\mathcal{E}[\phi,\mathscr{B};2]$ is the event that $|\phi|\gtrsim\mathscr{B}-[\llangle\partial_{\s}\phi\rrangle+\llangle\grad_{\mathbf{U}}\phi\rrangle]\mathfrak{n}\|\mathbf{U}\|$ at $(\s,\mathbf{U}^{\t,\cdot})$ for some $\s,\t\in\mathbb{X}^{\mathfrak{n}}$. (In particular, we can discretize time with an error controlled by derivatives of $\phi$ and the mesh length $\mathfrak{n}$; there is no randomness here.) By {a} union bound over all $\s,\t\in\mathbb{X}^{\mathfrak{n}}$, we know {that} $\mathbb{P}\{\mathcal{E}[\phi,\mathscr{B};2]\}$ is $\lesssim_{\gamma}$ than $\mathrm{RHS}\eqref{eq:steIIa}$ \emph{if we replace $\N^{-999}$ therein by $\mathfrak{n}\|\mathbf{U}\|$}. (For this, we implicitly use $|\mathbb{X}^{\mathfrak{n}}|^{2}\lesssim_{\gamma}\exp[\N^{\gamma}]$ for any $\gamma>0$, since $|\mathbb{X}^{\mathfrak{n}}|$ is polynomial in $\N$.) By another union bound, we can undo the replacement of $\N^{-999}$ by $\mathfrak{n}\|\mathbf{U}\|$ if we add the probability that $\|\mathbf{U}\|\gtrsim\N^{999}$, since $\mathfrak{n}\lesssim\N^{-\mathrm{D}}$ by assumption. By \eqref{eq:steI}, the probability that $\|\mathbf{U}\|\gtrsim\N^{999}$ is $\exp[-\N^{99}]$ plus the probability that $|\mathbf{U}^{\t,\x}|\gtrsim\N^{999}$ for some $(\t,\x)\in\mathbb{X}^{\mathfrak{n}}\times\mathbb{T}(\N)$. {A union bound} over $\mathbb{X}^{\mathfrak{n}}\times\mathbb{T}(\N)$, whose size is polynomial in $\N$ and thus $\lesssim_{\gamma}\exp[\N^{\gamma}]$, implies that this last probability is $\lesssim_{\gamma}\eqref{eq:steIIb}$. So \eqref{eq:steIIa}-\eqref{eq:steIIb} follows. Let us now move to \eqref{eq:steIII}. By the gradient relation for $\mathbf{h}$ and $\mathbf{U}$ (see Definitions \ref{definition:intro4}, \ref{definition:intro6}), we know $\mathbf{h}(\t,\x)-\mathbf{h}(\t,\x+\mathfrak{l})$ is a sum of $\lesssim|\mathfrak{l}|\lesssim\N$ many terms of the form $\mathrm{O}(1)\mathbf{U}^{\t,\z}$ (where $\z\in\mathbb{T}(\N)$). Thus, the difference in the first norm in \eqref{eq:steIII} is $\mathrm{O}(\N)$ times the first norm in \eqref{eq:steI}. (Everything so far in this proof of \eqref{eq:steIII} is deterministic.) It now suffices to use \eqref{eq:steI}. To justify the final sentence in the statement of Lemma \ref{lemma:ste}, it suffices to note that all we used is $|\mathbb{T}(\N)|\lesssim\N$ and the gradient relation between \eqref{eq:glsde} and \eqref{eq:hf}. These are true (either by construction or an easy argument) for the objects that we replace with, so we are done.
\end{proof}
\begin{lemma}\label{lemma:steeasy}
 Fix any discrete interval $\mathbb{K}\subseteq\mathbb{T}(\N)$ and any jointly smooth $\phi:\R\times\R^{\mathbb{K}}\to\R$. Fix a probability measure $\mathbb{P}$ on $\R^{\mathbb{K}}$. Take $\mathbb{X}^{\mathfrak{n}}$ and take $\mathfrak{t}_{1},\mathfrak{t}_{2}$, all from {Lemma \ref{lemma:ste}}. We have the following estimate for any $\gamma,\mathscr{B}>0$, in which both suprema are over $\s\in[\mathfrak{t}_{1},\mathfrak{t}_{2}]$, and $\llangle\rrangle$ is {the} sup-norm over $(\s,\mathbf{U})\in\R\times\R^{\mathbb{K}}$ (as in {Lemma \ref{lemma:ste}} but for $\mathbb{K}$ instead of $\mathbb{T}(\N)$):
\begin{align}
\mathbb{P}[{\textstyle\sup_{\s}}|\phi_{\s,\mathbf{U}}|\gtrsim\mathscr{B}] \ \lesssim_{\gamma} \ \exp[\N^{\gamma}]\times{\textstyle\sup_{\s}}\mathbb{P}[|\phi_{\s,\mathbf{U}}|\gtrsim\mathscr{B}-\llangle\partial_{\s}\phi\rrangle\N^{-999}]. \label{eq:steeasyI}
\end{align}
\end{lemma}
\begin{proof}
Fix $\mathbf{U}\in\R^{\mathbb{K}}$ and recall {the} $\t_{\circ}$ notation from Lemma \ref{lemma:ste}. By the fundamental theorem of calculus, we know {that}
\begin{align}
{\textstyle\sup_{\s}}|\phi_{\s,\mathbf{U}}-\phi_{\s_{\circ},\mathbf{U}}| \ \lesssim \ \llangle\partial_{\s}\phi\rrangle{\textstyle\sup_{\s}}|\s-\s_{\circ}| \ \lesssim \ \N^{-999}\llangle\partial_{\s}\phi\rrangle,
\end{align}
as $\s_{\circ}$ is in a discretization of $[\mathfrak{t}_{1},\mathfrak{t}_{2}]$ of mesh $\lesssim\N^{-999}$. So, $\mathrm{LHS}\eqref{eq:steeasyI}$ is at most the probability of $|\phi_{\s,\mathbf{U}}|\gtrsim\mathscr{B}-\llangle\partial_{\s}\phi\rrangle\N^{-999}$ for some $\s\in\mathbb{X}^{\mathfrak{n}}$. By union bound over $\mathbb{X}^{\mathfrak{n}}$, which has size polynomial in $\N$ and thus $\lesssim_{\gamma}\exp[\N^{\gamma}]$, \eqref{eq:steeasyI} follows.
\end{proof}
\subsection{Fluctuation property for canonical ensembles}
In a nutshell, Lemma \ref{lemma:vanishcanonical} states that if we take a functional and subtract its expectation with respect to an appropriate canonical measure expectation on its support, we get something that vanishes with respect to any canonical ensemble expectation on any superset of its support. It may sound tautological, but there are subtleties. (E.g., we must center with respect to the right canonical measure to make vanishing true {for \emph{all} the} canonical measure expectations. There is also an issue of scales, namely vanishing must hold for all supersets of the support.) Ultimately, these subtleties are easy to deal with and work out very cleanly; see the proof of Lemma 2 in \cite{GJ15}. Lemma \ref{lemma:vanishcanonical} lays this issue to rest.
\begin{lemma}\label{lemma:vanishcanonical}
 Fix a discrete interval $\mathbb{K}\subseteq\mathbb{T}(\N)$. Take any function $\mathfrak{A}:\R^{\mathbb{K}}\to\R$ for which $\E^{\s,\sigma,\mathbb{K}}|\mathfrak{A}|<\infty$ for all $\s,\sigma$. Define the functional $\mathbf{U}\mapsto\sigma[\mathbf{U}]$ for $\mathbf{U}\in\R^{\mathbb{K}}$ given by the average of $\mathbf{U}(\x)$ over $\x\in\mathbb{K}$. Lastly, define $\mathfrak{A}^{\mathrm{cent},\s}(\mathbf{U}):=\mathfrak{A}(\mathbf{U})-\E^{\s,\sigma[\mathbf{U}],\mathbb{K}}\mathfrak{A}$. In particular, $\E^{\s,\sigma[\mathbf{U}],\mathbb{K}}\mathfrak{A}$ is a canonical measure expectation where the charge density $\sigma[\mathbf{U}]$ is a functional of $\mathbf{U}$. (Its dependence on $\mathbf{U}\in\R^{\mathbb{K}}$ is therefore through $\sigma[\mathbf{U}]$.) For any $\varrho\in\R$ and $\s$ and $\mathbb{K}[+]\supseteq\mathbb{K}$, we have $\E^{\s,\varrho,\mathbb{K}[+]}\mathfrak{A}^{\mathrm{cent},\s}=0$.
\end{lemma}
\begin{proof}
This fact is (basically, namely in a different guise) used in the proof of Lemma 2 in \cite{GJ15}. We give a proof here. Suppose $\mathbb{K}[+]=\mathbb{K}$. In this case, under the measure $\E^{\s,\varrho,\mathbb{K}}$, we know {that} $\sigma[\mathbf{U}]=\varrho$ with probability 1. (In words, $\sigma[\mathbf{U}]$ is the charge density on $\mathbb{K}$, which we condition to equal $\varrho$ in $\E^{\s,\varrho,\mathbb{K}}$.) Thus, under $\E^{\s,\varrho,\mathbb{K}}$, we know {that} $\mathfrak{A}^{\mathrm{cent},\s}$ is just $\mathfrak{A}$ minus its $\E^{\s,\varrho,\mathbb{K}}$-expectation, and therefore it vanishes under $\E^{\s,\varrho,\mathbb{K}}$. Now, for a general $\mathbb{K}[+]\supseteq\mathbb{K}$, we claim (with explanation given afterwards) that
\begin{align}
\E^{\s,\varrho,\mathbb{K}[+]}\mathfrak{A}^{\mathrm{cent},\s} \ = \ {\textstyle\int_{\R}}\{\E^{\s,\eta,\mathbb{K}}\mathfrak{A}^{\mathrm{cent},\s}\}\d\mathds{Q}(\eta), \label{eq:vanishcanonical1}
\end{align}
where $\mathds{Q}(\eta)$ is a probability measure on $\R$ associated to the distribution of the charge density on $\mathbb{K}$ with respect to the measure in $\E^{\s,\varrho,\mathbb{K}[+]}$. To prove \eqref{eq:vanishcanonical1}, use the law of total expectation to condition on the charge density $\eta$ on $\mathbb{K}$, and take an expectation via $\d\mathds{Q}(\eta)$. (We then project the measure in $\E^{\s,\varrho,\mathbb{K}[+]}$ after conditioning on $\eta$ onto its $\R^{\mathbb{K}}$ marginal, which is allowed because $\mathfrak{A}^{\mathrm{cent},\s}$ has support $\mathbb{K}$.) \eqref{eq:vanishcanonical1} would then follow if we knew that under $\E^{\s,\varrho,\mathbb{K}[+]}$, conditioning on the charge density on $\mathbb{K}$ and projecting onto the $\R^{\mathbb{K}}$-marginal gives the canonical ensemble expectation $\E^{\s,\eta,\mathbb{K}}$. But this is just the fact that if one takes an increment of a random walk bridge and conditions on its average drift to be $\eta$, one gets a random walk bridge of average drift $\eta$ for the law of this increment. Thus, \eqref{eq:vanishcanonical1} holds. Since we showed {that} the integrand in $\mathrm{RHS}\eqref{eq:vanishcanonical1}$ is zero for all $\eta$, the lemma follows.
\end{proof}
%
%
%
{
\section{Glossary and explanation for notation}\label{section:glossary}
We now provide a glossary to streamline and ``categorize" some of the notation used in this paper, with the hopes of easing the reading of this paper. In any case, we emphasize that the notation in this paper is almost always referred to, recalled, and explained explicitly whenever it is used outside of its initial introduction or whenever it is not standard.
\begin{enumerate}
\item The terms $\mathbf{Z},\mathbf{U},\mathbf{H},\mathbf{J},\mathbf{h},\mathbf{G}$ and related objects are given bold font to indicate that they are determined by (stochastic) differential equations. The terms $\mathbf{S},\mathbf{W},\mathbf{Y}$ are related to $\mathbf{Z}$. Moreover, bold-objects, such as $\mathbf{V},\mathbf{A},\mathbf{D},\mathbf{I}$, whose input variables are superscripts and which are used in Sections \ref{section:proofoutline} and \ref{section:le}, are modifications of the $\mathbf{U}$-process.
\item The object $\mathscr{U}$ is the potential, and $\mathscr{W}(\t,\mathbf{u}):=\mathscr{U}(\t,\mathbf{u})-\frac12\bar{\alpha}(\t)\mathbf{u}^{2}$ as a function of $\mathbf{u}\in\R$. The objects $\mathscr{UP},\mathscr{HP},\mathscr{CP},\mathscr{AP}$, which are used in Section \ref{section:kv}, are technical modifications of $\mathscr{U}$.
\item Fraktur font (e.g. $\mathfrak{q},\mathfrak{z},\mathfrak{d},\mathfrak{w}$), whenever it is used for objects whose inputs are $\t$ and $\mathbf{u}$, is reserved for functionals of the $\mathbf{U}$-process (unless otherwise explicitly mentioned). Similarly, $\mathfrak{D}_{\mathrm{KL}},\mathfrak{D}_{\mathrm{FI}}$ (see Section \ref{section:le}) are the relative entropy and Fisher information, respectively; they are technically functions of the $(\mathbf{U},\mathbf{J})$-process. However, fraktur font for objects like $\mathfrak{l},\mathfrak{m}$, which do not take $\t,\mathbf{U}$ as input variables, refer instead of length-scales (i.e. positive integers) for averaging purposes. 
\item Similarly, sans-serif is often used for functions of $\mathbf{U}$ when $\mathbf{U}$ is a dummy variable. 
\item The ``ds-font", which resembles blackboard font, is often used for expectations with respect to (grand) canonical measures as introduced in Section \ref{subsubsection:hommeasures}, as well as for related objects. Said related objects include the $\mathds{R}$-objects from Section \ref{section:bg2}, which are determined by certain differences of $\mathds{E}$-expectations. (We distinguish these objects from the real numbers $\R$ by including superscripts.) Similarly, the notation $\mathds{S}(\N)$ is used to define the Lebesgue measure $\mathrm{Leb}[\mathds{S}(\N)]$ and the resulting measures $\mathbb{P}^{\mathrm{Leb},\sigma,\t,\mathds{I}}$; see Definition \ref{definition:le5}. (We choose the notation $\mathds{R}$ for these differences of $\E$-expectations because it stands for ``replacement" or ``renormalizaton".) Relatedly, $\mathscr{R}$, when used as a function of space and time, is an integration of $\mathds{R}$-functions against the heat kernel $\mathbf{H}$.
\item Another use of ``ds-font" is for objects like $\mathds{A},\mathds{C}$, the latter of which is distinguished from the complex numbers by including superscripts and inputs. These objects refer to space-time and space averages of functions of the $\mathds{R}$-objects from the previous bullet point; see Definition \ref{definition:bg2133}. Similarly, $\mathscr{A}$-terms denote space-time integration of $\mathds{A}$-terms against the heat kernel.
\item Yet another use of ``ds-font" is for $\mathds{I}$ and $\mathds{B}$ objects. The $\mathds{I}$ objects are generally sub-intervals in $\mathbb{T}(\N)$, and $\mathds{B}$ objects are generic countable index sets.
\item Subscripts and superscripts denoted by $?$ and $!$ indicate ``dummy placeholders". For example, $\grad^{?}$ for $?\in\{\pm\}$ means the collection $\{\grad^{+},\grad^{-}\}$. The superscript $\infty$, as in \eqref{eq:kpz}, indicates a (limit) SPDE that does not depend on $\N$.
\item We often use $\sigma$, which is sometimes a function of space-time, to denote the charge density in the grand-canonical and canonical measures from Section \ref{subsubsection:hommeasures}.
\item Script font, except for what has been mentioned above, is often used for operators of some sort. For example, $\mathscr{L}$-operators are infinitesimal generators (see Remark \ref{remark:le6}).
\item The notation $\mathrm{Cent}$ denotes a ``centering map" for random space-time functions to make them mean-zero. The notation $\mathrm{Loc}$ is a map which turns functions of $\mathbf{J},\mathbf{U}$ processes into functions of localized versions of these processes. In particular, we instead evaluate said functions at the localized processes. See Lemmas \ref{lemma:bg213101} and \ref{lemma:finalprop2} for the introduction of these maps.
\item The objects $\Lambda^{\pm,\mathrm{j}},\Phi^{\pm,\mathrm{j}},\Upsilon^{\pm,\mathrm{j},\mathrm{i}}$, and $\Lambda^{\pm,\mathrm{j},\mathrm{i}}$ are the basic building blocks relating $\mathds{A}^{\mathds{Q},\pm}$-averages with respect to ``neighboring" space-time scales; see Lemma \ref{lemma:bg2138}. These building blocks are basically averages of $\mathds{A}^{\mathds{Q},\pm}$-terms, which themselves are averages of $\mathrm{QCT}$-functionals, with additional technical cutoffs. These are the main error terms in Sections \ref{section:bg213main}-\ref{section:propbg212proof}.
\item The objects $\tau(\mathrm{j},\mathrm{i})$ and $\mathfrak{m}(\mathrm{j},\mathrm{i})$ are time-scales and length-scales, respectively, on which we average $\mathds{R}^{\mathfrak{q},\pm,\mathrm{j}}$-terms from Definition \ref{definition:bg24}. These time and length scales are introduced and used in Sections \ref{section:bg213main}-\ref{section:propbg212proof}. These scales get bigger (by a small power of $\N$) as we increase the index $\mathrm{i}$ (for any fixed $\mathrm{j}$ index). On the other hand, $\mathrm{j}$ indexes the renormalization terms $\mathds{R}^{\mathfrak{q},\pm,\mathrm{j}}$. (The point is that the maximal time and length scales on which we average $\mathds{R}^{\mathfrak{q},\pm,\mathrm{j}}$ depends on $\mathrm{j}$ itself. This is a purely technical point.)
\end{enumerate}
}








%
%
%
\end{document}